\documentclass{amsbook}

\newtheorem{thm}{Theorem}[section]
\newtheorem*{thm*}{Theorem}

\newtheorem{claim}[thm]{Claim}

\newtheorem{cor}[thm]{Corollary}

\newtheorem{lem}[thm]{Lemma}
\newtheorem*{lem*}{Lemma}
\newtheorem{mainthm}{Theorem}
\newtheorem*{mainthm*}{Theorem}
\newtheorem{maincor}[mainthm]{Corollary}

\newtheorem{prop}[thm]{Proposition}

\theoremstyle{definition}

\newtheorem*{case*}{Case}

\newtheorem{defn}[thm]{Definition}
\newtheorem*{defn*}{Definition}
\newtheorem{exmp}[thm]{Example}
\newtheorem*{exmp*}{Example}

\newtheorem{maindefn}[mainthm]{Definition}

\renewcommand{\thestep}{}

\theoremstyle{remark}

\renewcommand{\thecase}{}

\numberwithin{subcase}{case}

\numberwithin{subsubcase}{subcase}

\newtheorem{rmk}[thm]{Remark}
\newtheorem*{rmk*}{Remark}

\makeatletter
\def\alphenumi{
  \def\theenumi{\alph{enumi}}
  \def\p@enumi{\theenumi}
  \def\labelenumi{(\@alph\c@enumi)}}
\makeatother

\makeatletter
\def\thecase{\@arabic\c@case}
\makeatother

\makeatletter
\def\thestep{\@arabic\c@step}
\makeatother

\newcommand{\transv}{\mathrel{\text{\tpitchfork}}}
\makeatletter
\newcommand{\tpitchfork}{%
  \vbox{
    \baselineskip\z@skip
    \lineskip-.52ex
    \lineskiplimit\maxdimen
    \m@th
    \ialign{##\crcr\hidewidth\smash{$-$}\hidewidth\crcr$\pitchfork$\crcr}
  }%
}
\makeatother

\DeclareFontFamily{U}{mathx}{\hyphenchar\font45}
\DeclareFontShape{U}{mathx}{m}{n}{
      <5> <6> <7> <8> <9> <10>
      <10.95> <12> <14.4> <17.28> <20.74> <24.88>
      mathx10
      }{}
\DeclareSymbolFont{mathx}{U}{mathx}{m}{n}
\DeclareFontSubstitution{U}{mathx}{m}{n}
\DeclareMathAccent{\widecheck}{0}{mathx}{"71}
\DeclareMathAccent{\wideparen}{0}{mathx}{"75}

\newcount\hh
\newcount\mm
\mm=\time
\hh=\time
\divide\hh by 60
\divide\mm by 60
\multiply\mm by 60
\mm=-\mm
\advance\mm by \time
\def\hhmm{\number\hh:\ifnum\mm<10{}0\fi\number\mm}

\let\oldmarginpar\marginpar
\renewcommand\marginpar[1]{\-\oldmarginpar[\raggedleft\footnotesize #1]%
{\raggedright\footnotesize #1}}

\newcommand\AAA{\mathbb{A}}

\newcommand\CC{\mathbb{C}}

\newcommand\FF{\mathbb{F}}
\newcommand\GG{\mathbb{G}}

\newcommand\KK{\mathbb{K}}

\newcommand\NN{\mathbb{N}}

\newcommand\PP{\mathbb{P}}
\newcommand\QQ{\mathbb{Q}}
\newcommand\RR{\mathbb{R}}

\newcommand\TT{\mathbb{T}}

\newcommand\VV{\mathbb{V}}

\newcommand\ZZ{\mathbb{Z}}

\newcommand\cF{{\mathcal{F}}}
\newcommand\cG{{\mathcal{G}}}

\newcommand\cM{{\mathcal{M}}}

\newcommand\cX{{\mathcal{X}}}
\newcommand\cY{{\mathcal{Y}}}

\newcommand\fa{{\mathfrak{a}}}

\newcommand\fb{{\mathfrak{b}}}

\newcommand\fm{{\mathfrak{m}}}

\newcommand\fp{{\mathfrak{p}}}

\newcommand\sA{{\mathscr{A}}}
\newcommand\sB{{\mathscr{B}}}
\newcommand\sC{{\mathscr{C}}}

\newcommand\sF{{\mathscr{F}}}
\newcommand\sG{{\mathscr{G}}}

\newcommand\sI{{\mathscr{I}}}
\newcommand\sJ{{\mathscr{J}}}

\newcommand\sL{{\mathscr{L}}}

\newcommand\sN{{\mathscr{N}}}
\newcommand\sO{{\mathscr{O}}}

\newcommand\sR{{\mathscr{R}}}
\newcommand\sS{{\mathscr{S}}}
\newcommand\sT{{\mathscr{T}}}

\newcommand\bg{{\mathbf{g}}}

\newcommand\eps{\varepsilon}

\newcommand\om{\omega}

\newcommand\GL{\operatorname{GL}}
\newcommand\Or{\operatorname{O}}
\DeclareMathOperator{\ord}{ord}

\newcommand\SL{\operatorname{SL}}
\newcommand\SO{\operatorname{SO}}

\newcommand\U{\operatorname{U}}

\newcommand\less{\setminus}

\newcommand\Ad{{\operatorname{Ad}}}

\DeclareMathOperator{\An}{An}

\DeclareMathOperator{\Aut}{Aut}

\newcommand\codim{\operatorname{codim}}

\DeclareMathOperator{\Corank}{Corank}
\newcommand\cosupp{\operatorname{cosupp}}

\DeclareMathOperator{\Crit}{Crit}

\newcommand\diag{\operatorname{diag}}

\newcommand\Diff{\operatorname{Diff}}

\newcommand\End{\operatorname{End}}

\DeclareMathOperator{\Euler}{Euler}

\newcommand\Gl{\operatorname{Gl}}

\newcommand\Gr{\operatorname{Gr}}
\newcommand\grad{\operatorname{grad}}

\newcommand\hess{\operatorname{hess}}
\newcommand\Hess{\operatorname{Hess}}
\newcommand\Hilb{\operatorname{Hilb}}

\newcommand\Hom{\operatorname{Hom}}

\DeclareMathOperator{\Imag}{Im}

\newcommand\Ker{\operatorname{Ker}}
\DeclareMathOperator{\length}{length}

\DeclareMathOperator{\Mat}{Mat}

\DeclareMathOperator{\pr}{pr}

\DeclareMathOperator{\bProj}{\textbf{Proj}}

\newcommand\Ran{\operatorname{Ran}}
\newcommand\rank{\operatorname{rank}}
\newcommand\Rank{\operatorname{Rank}}
\newcommand\Real{\operatorname{Re}}

\DeclareMathOperator{\Spec}{Spec}

\DeclareMathOperator{\spm}{spm}
\DeclareMathOperator{\Spm}{Spm}
\newcommand\Stab{\operatorname{Stab}}

\newcommand\supp{\operatorname{supp}}

\DeclareMathOperator{\Sym}{Sym}

\DeclareMathOperator{\Zero}{Zero}

\newcommand\Bl{{\mathrm{Bl}}}

\DeclareMathOperator{\Cone}{Cone}

\DeclareMathOperator{\embdim}{\mathrm{emb\ dim}}

\DeclareMathOperator{\ev}{ev}

\DeclareMathOperator{\expdim}{\mathrm{exp\ dim}}

\newcommand\id{{\mathrm{id}}}
\DeclareMathOperator{\inn}{in}
\newcommand\inv{{\mathrm{inv}}}

\newcommand\loc{{\mathrm{loc}}}

\newcommand\mutatis{{\emph{mutatis mutandis }}}

\newcommand\red{{\mathrm{red}}}

\DeclareMathOperator{\rel}{rel}

\newcommand\set{{\mathrm{set}}}
\newcommand\sing{\mathrm{sing}}

\newcommand\sm{{\mathrm{sm}}}

\numberwithin{equation}{section}
\numberwithin{section}{chapter}
\numberwithin{figure}{section}

\usepackage{amscd}
\usepackage{amssymb}
\usepackage{bm}
\usepackage{graphicx}
\usepackage{hyperref}
\usepackage{mathrsfs}
\usepackage[usenames]{color}
\usepackage{paralist}
\usepackage{slashed}
\usepackage{tikz-cd}
\usepackage{url}
\usepackage{verbatim}
\usepackage{xr}
\usepackage[all]{xy}
\usepackage{stmaryrd}

\renewcommand{\d}[1]{\ensuremath{\operatorname{d}\!{#1}}}

\hypersetup{pdftitle={Bia{\l}ynicki--Birula Theory, Morse-Bott Theory, and Resolution of Singularities for Analytic Spaces}}
\hypersetup{pdfauthor={Paul M. N. Feehan}}

\begin{document}

\frontmatter

\title[Bia{\l}ynicki--Birula and Morse--Bott Theories for Analytic Spaces]{Bia{\l}ynicki--Birula Theory, Morse--Bott Theory, and Resolution of Singularities for Analytic Spaces}

\author[Paul M. N. Feehan]{Paul M. N. Feehan}







\maketitle


\tableofcontents
\listoffigures

\chapter*{Preface}
\label{chap:Preface}
Our goal in this work is to develop aspects of Bia{\l}ynicki--Birula and Morse--Bott theory that can be extended from the classical setting of smooth manifolds to that of complex analytic spaces with a holomorphic $\CC^*$ action. We extend prior results on existence of Bia{\l}ynicki--Birula decompositions for compact, complex K\"ahler manifolds to non-compact complex manifolds and develop functorial properties of the Bia{\l}ynicki--Birula decomposition, in particular with respect to blowup along a $\CC^*$-invariant, embedded complex submanifold. We deduce the existence of a Bia{\l}ynicki--Birula decomposition for a $\CC^*$-invariant, closed, complex analytic subspace of complex manifold with a $\CC^*$ action; derive geometric consequences for the positivity of the Bia{\l}ynicki--Birula index and co-index at a fixed point; and we develop stronger versions of these results by applying resolution of singularities for analytic spaces.

\chapter*{Acknowledgments}
\label{chap:Acknowledgments}
For their helpful communications and patient explanations in response to my many questions, I warmly thank Dan Abramovich, Francesca Acquistapace, Denis Auroux, David Bayer, Lev Borisov, Michel Brion, Anders Buch, Ian Coley, Tedi Dragheci, Mariano Echeverria, Santiago Encinas, Luis Fern\'andez, Ron Fintushel, Eberhard Freitag, Robert Friedman, Mark Goresky, Ian Hambleton, Herwig Hauser, Daniel Huybrechts, Johan de Jong, Veronica Kalicki, J\'anos Koll\'ar, Thomas Leness, Tian-Jun Li, Aleksandar Milivojevic, Konstantin Mischaikow, Paul Nelson, Keith Pardue, Jacob Sturm, Loring Tu, Ravi Vakil, Scott Wilson, Sven-Ake Wegner, Richard Wentworth, and Jaros{\l}aw W{\l}odarczyk. We thank the National Science Foundation for its support of our research via the grant DMS-2104865.
This monograph is also based in part on work supported by the National Science Foundation under Grant No. 1440140, while the author was in residence at the Mathematical Sciences Research Institute in Berkeley, California, during Fall 2022 as a Research Professor in the program \emph{Analytic and Geometric Aspects of Gauge Theory}.
\bigskip
\bigskip

\leftline{Paul M. N. Feehan}
\leftline{Department of Mathematics}
\leftline{Rutgers, The State University of New Jersey}
\leftline{Piscataway, NJ 08854-8019}
\leftline{United States}
\medskip

\leftline{\texttt{feehan@math.rutgers.edu}}
\leftline{\url{math.rutgers.edu/~feehan}}
\bigskip

\leftline{This version: December 22, 2022}

\mainmatter

\chapter{Introduction}
\label{chap:Introduction}
Our goal in this work is to develop aspects of Morse--Bott theory that can be extended from the classical case of smooth manifolds to the case of analytic spaces that may be singular. For this purpose, we primarily restrict our attention to complex analytic spaces that are $\CC^*$-invariant, closed subspaces of Hermitian manifolds endowed with a holomorphic $\CC^*$ action such that the induced action of the circle $S^1 \subset \CC^*$ is Hamiltonian with respect to the fundamental two-form defined by the Hermitian metric and almost complex structure. For smooth Hermitian manifolds, it follows from results initiated by Frankel \cite{Frankel_1959} that Morse--Bott theory is equivalent to a parallel theory due to Bia{\l}ynicki--Birula \cite{Bialynicki_1973}. While Bia{\l}ynicki--Birula theory was originally developed for smooth algebraic varieties, it has since been extended to algebraic varieties in far greater generality, as well as to compact, complex K\"ahler manifolds and certain complex analytic spaces.

Our work is motivated in part by a joint project with Leness \cite{Feehan_Leness_introduction_virtual_morse_theory_so3_monopoles} and Wentworth more recently, where we consider aspects of Bia{\l}ynicki--Birula and Morse--Bott theories on a complex analytic space arising as the moduli space of non-Abelian monopoles over a compact, complex, K\"ahler surface. As explained in \cite{Feehan_Leness_introduction_virtual_morse_theory_so3_monopoles}, one of the aims of that project is to reprove the Bogomolov--Miyaoka--Yau inequality (see Miyaoka \cite{Miyaoka_1977} and Yau \cite{YauPNAS,Yau}) for compact, complex surfaces of general type using Seiberg--Witten gauge theory and ultimately extend that argument to the case of closed, connected, oriented, smooth four-dimensional manifolds with $b_1=0$, odd $b^+ \geq 3$, and Seiberg--Witten simple type. Our approach is inspired by the extraordinarily influential\footnote{Hitchin's articles \cite{Hitchin_1987, Hitchin_1987duke} on Higgs pairs over Riemann surfaces have well over 900 citations.} study due to Hitchin of the moduli space of Higgs pairs over Riemann surfaces \cite{Hitchin_1987} (see also \cite{Hitchin_1987duke}), some features of which were recently extended by Hausel and Hitchin \cite{Hausel_Hitchin_2022}. To prove two of his key results \cite[Proposition 7.1, p. 92 and Theorem 7.6, p. 96]{Hitchin_1987}, Hitchin applies the Hirzebruch--Riemann--Roch Theorem to compute the Morse--Bott index, co-index, and nullity of smoothly embedded critical submanifolds of a Hamiltonian function for the canonical circle action on the moduli space of Higgs pairs, assumed to be smooth by virtue of a hypothesis that the degree of the rank two Hermitian vector bundle over the Riemann surface is odd. We generalize this paradigm in \cite{Feehan_Leness_introduction_virtual_morse_theory_so3_monopoles} with Leness and apply the Hirzebruch--Riemann--Roch Theorem to compute the \emph{virtual} Morse--Bott index, co-index, and nullity of critical points of a Hamiltonian function for the canonical circle action on the moduli space of non-Abelian monopoles over a compact, complex, K\"ahler surface, without any assumption of smoothness. Exploring the geometric significance of the virtual Morse--Bott index, co-index, and nullity of critical points in the general context of a complex analytic space with a holomorphic $\CC^*$ action is one of the themes of our present work.

In Section \ref{sec:Functorial_properties_BB_decompositions_complex_manifolds}, we begin by providing a general definition of a Bia{\l}ynicki--Birula decomposition for a complex manifold endowed with a holomorphic $\CC^*$ action (Definition \ref{maindefn:BB_decomposition_complex_manifold}) and of the Bia{\l}ynicki--Birula nullity, co-index, and index at a fixed point (Definition \ref{maindefn:Stable_unstable_submanifolds_BB_index_co-index_nullity}). We extend the previous Bia{\l}ynicki--Birula decompositions for smooth algebraic varieties and compact, complex K\"ahler manifolds due to Bia{\l}ynicki--Birula \cite{Bialynicki_1973}, Carrell and Sommese \cite{Carrell_Sommese_1978ms, Carrell_Sommese_1979cmh, Carrell_Sommese_1983}, Fujiki \cite{Fujiki_1979}, and Yang \cite{Yang_2008} to a broader class of non-compact almost Hermitian manifolds, with the strongest results for complex Hermitian manifolds (see Theorem \ref{mainthm:BB_decomposition_complex_manifold_C*_action}). We also describe certain functorial properties of Bia{\l}ynicki--Birula decompositions, including induced decompositions on embedded complex submanifolds (Theorem \ref{mainthm:BB_decomposition_C*_invariant_complex_submanifold}) and induced decompositions on blowups of complex manifolds along embedded complex submanifolds (Theorem \ref{mainthm:BB_decomposition_blowup_complex_manifold_C*_action_along_submanifold}).

Section \ref{sec:BB_decomposition_complex_analytic_spaces} highlights our results on Bia{\l}ynicki--Birula decompositions for complex analytic spaces in the sense of Grauert and Remmert \cite{Grauert_Remmert_coherent_analytic_sheaves}. Definition \ref{maindefn:BB_decomposition_complex_analytic_space} describes a weaker notion of Bia{\l}ynicki--Birula decomposition than that of Definition \ref{maindefn:BB_decomposition_complex_manifold}, but one which is still adequate for our applications, while Definition \ref{maindefn:Stable_unstable_subspaces_BB_index_co-index_nullity} introduces the Bia{\l}ynicki--Birula nullity, co-index, and index at a fixed point as the Krull dimensions of the local rings of the fixed point subspace, stable subspace, and unstable subspace. Theorem \ref{mainthm:BB_decomposition_C*_invariant_complex_analytic_subspace} provides a Bia{\l}ynicki--Birula decomposition for a $\CC^*$-invariant, closed, complex analytic subspace of complex manifold with a holomorphic $\CC^*$ action and geometric consequences for the positivity of the Bia{\l}ynicki--Birula nullity, co-index, and index at a fixed point, respectively. In practice, the Krull dimensions of the local rings at a fixed point are difficult to compute but they may be bounded from below by more easily computed expected dimensions corresponding to a convenient choice of local model spaces (Theorem \ref{mainthm:Upper_and_lower_bounds_Krull_dimensions_X0_Xpm_at_p}) and this motivates our Definition \ref{maindefn:Virtual_BB_nullity_co-index_index} of \emph{virtual} Bia{\l}ynicki--Birula nullity, co-index, and index at a fixed point.

Section \ref{sec:BB_decomposition_embedded_resolution_singularities_complex_analytic_spaces} provides important strengthenings of Theorem \ref{mainthm:BB_decomposition_C*_invariant_complex_analytic_subspace} for Bia{\l}ynicki--Birula decompositions for complex analytic spaces that can be achieved by appealing to the resolution of singularities for analytic spaces due to Hironaka \cite{Hironaka_1964-I-II}. Theorem \ref{mainthm:BB_decomposition_strict_transform_complex_analytic_subspace} provides a Bia{\l}ynicki--Birula decomposition for the strict transform, relative to a resolution morphism, of a $\CC^*$-invariant, complex analytic subspace of a complex manifold with a holomorphic $\CC^*$ action. Corollary \ref{maincor:BB_decomposition_strict_transform_complex_analytic_subspace} yields a Bia{\l}ynicki--Birula decomposition for the strict transform of a $\CC^*$-invariant, closed complex analytic subspace of a complex manifold with a holomorphic $\CC^*$ action and several interpretations of the Bia{\l}ynicki--Birula nullity, co-index, and index at a fixed point in the complex analytic space. The fixed point subspace, stable subspace, and unstable subspace need not, of course, lie in the smooth locus of the complex analytic space. However, if the circle action $S^1\times X \to X$ induced by restricting the $\CC^*$ action to $S^1 \subset \CC^*$ has a Hamiltonian function $f:X\to\RR$, then Corollary \ref{maincor:BB_decomposition_strict_transform_complex_analytic_subspace} asserts that such a fixed point is not a local minimum\footnote{We say that $p$ is not a local minimum (respectively, maximum) of $f:Y_\sm\to \RR$ if $p$ is not a local minimum (respectively, maximum) of $f:Y_\sm\cup\{p\}\to \RR$.} (respectively, maximum) of the restriction $f:Y_\sm\to \RR$ of $f$ to the smooth locus $Y_\sm \subset Y$ if the Bia{\l}ynicki--Birula co-index (respectively, index) at a fixed point is positive. 

In Section \ref{sec:Frankel_theorem_circle_actions_almost_Hermitian_manifolds}, we review an extension, due to the author and Leness \cite{Feehan_Leness_introduction_virtual_morse_theory_so3_monopoles}, of Frankel's Theorem \ref{mainthm:Frankel_almost_Hermitian} for a Hamiltonian function of a circle action on a smooth almost Hermitian manifold \cite{Frankel_1959}. Frankel's Theorem provided the original basis for the link between Bia{\l}ynicki--Birula and Morse--Bott theories. In Section \ref{sec:Morse-Bott_decompositions_complex_analytic_subspaces_Hermitian_manifolds}, we briefly describe Morse--Bott decompositions for complex analytic subspaces of Hermitian manifolds. In Section \ref{sec:Interpretation_virtual_Morse-Bott_nullity_coindex_index}, we outline how the virtual Morse--Bott nullity, co-index, and index typically arise in the context of gauge theory. In Section \ref{sec:Outline}, we summarize the remainder of our current work.

\section{Bia{\l}ynicki--Birula decompositions for complex manifolds}
\label{sec:Functorial_properties_BB_decompositions_complex_manifolds}
In this section, we highlight our main results on the existence of Bia{\l}ynicki--Birula decompositions for complex manifolds equipped with holomorphic $\CC^*$ actions and some of their functorial properties of these decompositions, including their behavior with respect to holomorphic embeddings and blowups along complex submanifolds.

Bia{\l}ynicki--Birula decompositions for complex manifolds equipped with holomorphic $\CC^*$ actions were provided by Bia{\l}ynicki--Birula \cite{Bialynicki_1973, Bialynicki_1974}, Carrell and Sommese \cite{Carrell_Sommese_1978ms, Carrell_Sommese_1979cmh, Carrell_Sommese_1983}, Fujiki \cite{Fujiki_1979}, and Yang \cite{Yang_2008} under various additional hypotheses, such as that
\begin{inparaenum}[\itshape a\upshape)]
\item the manifold be a smooth algebraic variety\footnote{Bia{\l}ynicki--Birula allows locally affine $\GG_m$ actions on algebraic varieties over an algebraically closed field $k$ with group of units $\GG_m$.} as in Bia{\l}ynicki--Birula \cite{Bialynicki_1973}, or  
\item the manifold be compact and K\"ahler as in Carrell and Sommese \cite{Carrell_Sommese_1978ms, Carrell_Sommese_1979cmh, Carrell_Sommese_1983} or Fujiki \cite{Fujiki_1979}, or
\item the circle action induced by the restriction of the $\CC^*$ action to $S^1$ has a Hamiltonian function that is proper and bounded below and the number of connected components of the subset of fixed points of the $\CC^*$ action is finite, as in Yang \cite{Yang_2008}.
\end{inparaenum}
However, elementary examples --- such as a linear $\CC^*$ action on finite-dimensional complex vector space, as in Example \ref{exmp:Bialynicki-Birula_decomposition_C^2} --- indicate that while additional hypotheses such as those indicated above may be sufficient, they are not necessary. Sommese \cite{Sommese_1982} provides more sophisticated examples in that vein. With that in mind, we provide the following definition based in part on the conclusions of results due to Carrell and Sommese \cite{Carrell_Sommese_1978ms} and Fujiki \cite{Fujiki_1979} (see Remark \ref{rmk:BB_decomposition_complex_manifold_variants} for a more detailed explanation).

\begin{maindefn}[Bia{\l}ynicki--Birula decomposition for a holomorphic $\CC^*$ action on a complex manifold]
\label{maindefn:BB_decomposition_complex_manifold}
Let $X$ be a complex manifold and $\CC^*\times X \to X$ be a holomorphic $\CC^*$ action such that the subset $X^0 := X^{\CC^*} \subset X$ of fixed points of the $\CC^*$ action is non-empty with at most countably many connected components, $X_\alpha^0$ for $\alpha\in\sA$, that are embedded complex submanifolds of $X$. For each $\alpha\in\sA$, define
\begin{equation}
  \label{eq:Complex_manifold_Xalpha_plus_minus_submanifolds}
  X_\alpha^+ := \left\{z:\lim_{\lambda\to 0} \lambda\cdot z \in X_\alpha^0 \right\}
  \quad\text{and}\quad
  X_\alpha^- := \left\{z:\lim_{\lambda\to \infty} \lambda\cdot z \in X_\alpha^0 \right\},
\end{equation}
so the subsets $X_\alpha^+\subset X$ are $\CC^*$-invariant and mutually disjoint for all $\alpha\in\sA$ and similarly for the the subsets $X_\alpha^-\subset X$ for all $\alpha\in\sA$, and
\begin{equation}
  \label{eq:pi_Xalpha_pm_to_Xalpha_0_bundle_maps}
  \pi_\alpha^+(z) := \lim_{\lambda\to 0} \lambda\cdot z, \quad\text{for all } z \in X_\alpha^+,
  \quad\text{and}\quad
  \pi_\alpha^-(z) := \lim_{\lambda\to \infty} \lambda\cdot z, \quad\text{for all } z \in X_\alpha^-.
\end{equation}
Then $X$ has a \emph{(mixed, plus, or minus) fundamental Bia{\l}ynicki--Birula decomposition} if the following hold:
\begin{enumerate}
\item\label{item:Xalpha_pm_submanifold_X}
  Each $X_\alpha^+$ is an embedded complex submanifold of $X$;
\item\label{item:pi_Xalpha_pm_to_Xalpha_0_bundle}
  The natural map $\pi_\alpha^+:X_\alpha^+\to X_\alpha^0$ is a $\CC^*$-equivariant, holomorphic, maximal-rank surjection
  whose fibers are vector spaces;
\item\label{item:Xalpha_0_section_Xalpha_pm}
  $X_\alpha^0$ is a section of $X_\alpha^+$;
\end{enumerate}  
and the analogous properties hold for the subsets $X_\alpha^-$ and for the natural maps $\pi_\alpha^-:X_\alpha^-\to X_\alpha^0$. Furthermore, we require that:
\begin{enumerate}
\setcounter{enumi}{4}  
\item\label{item:Normal_bundle_Xalpha_0_Xalpha_pm}
  The normal bundles $N_{X_\alpha^0/X_\alpha^+}$ of $X_\alpha^0$ in $X_\alpha^+$ and $N_{X_\alpha^0/X_\alpha^-}$ of $X_\alpha^0$ in $X_\alpha^-$ are subbundles of the normal bundle $N_{X_\alpha^0/X}$ of $X_\alpha^0$ in $X$. In particular,
\begin{subequations}
\label{eq:Normal_bundle_Xalpha_0_Xalpha_pm_equals_Normal_bundle_pm_Xalpha_0_Xalpha}    
\begin{gather}
  \label{eq:Normal_bundle_Xalpha_0_Xalpha_+_equals_Normal_bundle_+_Xalpha_0_Xalpha}    
  N_{X_\alpha^0/X_\alpha^+} = N_{X_\alpha^0/X}^+
  \quad\text{and}\quad
  TX_\alpha^+\restriction X_\alpha^0 = TX_\alpha^0\oplus N_{X_\alpha^0/X}^+,
  \\
  \label{eq:Normal_bundle_Xalpha_0_Xalpha_-_equals_Normal_bundle_-_Xalpha_0_Xalpha}
  N_{X_\alpha^0/X_\alpha^-} = N_{X_\alpha^0/X}^-
  \quad\text{and}\quad
  TX_\alpha^-\restriction X_\alpha^0 = TX_\alpha^0\oplus N_{X_\alpha^0/X}^-,
\end{gather}    
\end{subequations}
where
\begin{equation}
  \label{eq:TX_Xalpha_0_bundle_weight-sign_decomposition}
  TX\restriction X_\alpha^0 = T^0X_\alpha \oplus N_{X_\alpha^0/X}^+ \oplus N_{X_\alpha^0/X}^-
\end{equation}
is the weight-sign decomposition defined by the $S^1$ action on $X$ induced by the $\CC^*$ action;
\item\label{item:BB_decomposition_complex_manifold_mixed_plus_minus} $X$ may be expressed as a disjoint union,
\begin{equation}
  \label{eq:Complex_manifold_mixed_decomposition}
   X = \bigsqcup_{(\alpha,j)\in A} X_\alpha^j,
\end{equation}
for some subset $A \subset \sA\times \{+,-\}$.
\item\label{item:Closure_Xalpha_pm_complex_analytic_subvariety_X}
If in addition to the preceding properties the topological closure in $X$ of $\bar X_\alpha^+$ or $\bar X_\alpha^-$ for each subset $X_\alpha^+$ or $X_\alpha^-$ in \eqref{eq:Complex_manifold_mixed_decomposition}, respectively, is a complex analytic subvariety of $X$, then we say that $X$ has an \emph{analytic Bia{\l}ynicki--Birula decomposition}.
\item\label{item:Xalpha_pm_Zariski_open_in_barXalpha_pm}
If in addition to the preceding properties each subset $X_\alpha^+$ or $X_\alpha^-$ in \eqref{eq:Complex_manifold_mixed_decomposition} is Zariski-open in $\bar X_\alpha^+$ or $\bar X_\alpha^-$, respectively, then we say that $X$ has a \emph{meromorphic Bia{\l}ynicki--Birula decomposition}.
\end{enumerate}
When we omit the qualifier fundamental, analytic, or meromorphic and say that $X$ has a \emph{Bia{\l}ynicki--Birula decomposition}, we mean that it is meromorphic in the sense of Item \eqref{item:Xalpha_pm_Zariski_open_in_barXalpha_pm}. If one can express $X$ as
\begin{equation}
  \label{eq:Complex_manifold_pure_plus_or_minus_decomposition}
  X = \bigsqcup_{\alpha\in \sA} X_\alpha^+
  \quad\text{or}\quad
  X = \bigsqcup_{\alpha\in \sA} X_\alpha^-,
\end{equation}
where the union is disjoint, then we say that $X$ has a \emph{plus} or \emph{minus decomposition}, respectively and, otherwise, if $X$ is expressed as in \eqref{eq:Complex_manifold_mixed_decomposition}, that it has a \emph{mixed decomposition}.
\qed
\end{maindefn}

\begin{rmk}[Fundamental, analytic, and meromorphic Bia{\l}ynicki--Birula decompositions in Definition \ref{maindefn:BB_decomposition_complex_manifold}]
\label{rmk:BB_decomposition_complex_manifold_variants}
The properties in Definition \ref{maindefn:BB_decomposition_complex_manifold} up through Item \eqref{item:BB_decomposition_complex_manifold_mixed_plus_minus} that define a \emph{fundamental} Bia{\l}ynicki--Birula decomposition are based on the conclusions of \cite[Lemma 2.1, p. 804]{Fujiki_1979} due to Fujiki and of \cite[Proposition II, p. 55, excluding Items (A), (D), and (E)]{Carrell_Sommese_1978ms} due to Carrell and Sommese.

The properties in Definition \ref{maindefn:BB_decomposition_complex_manifold} up to and including Item \eqref{item:Closure_Xalpha_pm_complex_analytic_subvariety_X} that define an \emph{analytic} Bia{\l}ynicki--Birula decomposition are based on the conclusions of \cite[Lemma 2.1, p. 804 and Proposition 2.8, p. 812]{Fujiki_1979} due to Fujiki.

The properties in Definition \ref{maindefn:BB_decomposition_complex_manifold} up to and including Item \eqref{item:Xalpha_pm_Zariski_open_in_barXalpha_pm} that define a \emph{meromorphic} Bia{\l}ynicki--Birula decomposition are based on the conclusions of \cite[Lemma 2.1, p. 804 and Theorem 2.2, p. 805]{Fujiki_1979} due to Fujiki and \cite[Proposition II, p. 55, excluding Items (D) and (E)]{Carrell_Sommese_1978ms} due to Carrell and Sommese. For $\KK=\RR$ or $\CC$ (see Fujiki \cite[p. 226]{Fujiki_1978im} for the case $\KK=\CC$), an open subset of a $\KK$-analytic space $(X,\sO_X)$ (as in Definition \ref{defn:Analytic_space} is called \emph{Zariski open} if it is the complement of a analytic subset of $X$ (as in Definition \ref{defn:Analytic_set}).
  
In \cite[Section 2, p. 805]{Fujiki_1979}, Fujiki defines a plus decomposition of $X$ in \eqref{eq:Complex_manifold_pure_plus_or_minus_decomposition} to be \emph{meromorphic} if each inclusion $X_\alpha^+\subset X$ has the properties described in Items \eqref{item:Closure_Xalpha_pm_complex_analytic_subvariety_X} and \eqref{item:Xalpha_pm_Zariski_open_in_barXalpha_pm} of Definition \ref{maindefn:BB_decomposition_complex_manifold}. Carrell and Sommese apply the term meromorphic in the same context in \cite[Proof of Proposition II, Item (A), p. 55]{Carrell_Sommese_1978ms}, citing extension theorems due to Siu \cite{Siu_1975}, \cite[Siu's Extension Theorem, p. 108]{Sommese_1975} and Sommese \cite[Lemma I-A, p. 108, and Lemma II-A, p. 109]{Sommese_1975} to verify the properties in Item \eqref{item:Closure_Xalpha_pm_complex_analytic_subvariety_X} of Definition \ref{maindefn:BB_decomposition_complex_manifold} in their setting. Siu \cite[p. 423]{Siu_1975} also uses this terminology, though his use of the term is appears different: he attributes his definition to Remmert \cite[Definition 3.15, p. 367]{Remmert_1957} and in \cite[Proposition, p. 107]{Sommese_1975}, Sommese attributes his use of the term meromorphic to Remmert as well.
\end{rmk}

\begin{rmk}[On the weight-sign decomposition of the restriction of the tangent bundle to a component of the fixed-point subset]
\label{rmk:Weight-sign_decomposition_restriction_tangent_bundle}  
Continue the notation of Definition \ref{maindefn:BB_decomposition_complex_manifold}. If $p \in X_\alpha^0$, so $p$ is a fixed point of the $\CC^*$ action on $X$, then the splitting \eqref{eq:TX_Xalpha_0_bundle_weight-sign_decomposition} simplifies to give
\begin{equation}
  \label{eq:T_pX_weight-sign_decomposition}
  T_pX = T_p^0X \oplus T_p^+X  \oplus T_p^-X,
\end{equation}
where the subspaces $T_p^0X$, $T_p^+X$, and $T_p^-X$ are the zero, positive, and negative weight subspaces of $T_pX$, respectively, for the $S^1$ action. The positive and negative weight subspaces obey the identities
\begin{equation}
  \label{eq:T_pX_pm_expressions}
  T_pX^\pm = T_p^\pm X = N_{X_\alpha^0/X}^\pm|_p = N_{X_\alpha^0/X_\alpha^\pm}|_p
\end{equation}
and the identity
\begin{equation}
  \label{eq:T_pX_0_expression}
    T_pX_\alpha^0 = T_p^0X
\end{equation}
is obeyed by the zero weight subspace.
\end{rmk}  

By analogy with the corresponding definitions in Morse--Bott theory (see Bott \cite[Definition, p. 248]{Bott_1954}, \cite{Bott_1959} or Nicolaescu \cite[Definition 2.41]{Nicolaescu_morse_theory} for a modern exposition), we have the

\begin{maindefn}[Stable and unstable submanifolds of a complex manifold and Bia{\l}ynicki--Birula index, co-index, and nullity]
\label{maindefn:Stable_unstable_submanifolds_BB_index_co-index_nullity}
Continue the notation of Definition \ref{maindefn:BB_decomposition_complex_manifold}. The subset $X_\alpha^+$ (respectively, $X_\alpha^-$) is called the \emph{stable} (respectively, \emph{unstable}) \emph{submanifold} for the fixed-point submanifold $X_\alpha^0$. For each point $p\in X_\alpha^0$, the complex dimensions $\beta_X^0(p)$, $\beta_X^+(p)$, and $\beta_X^-(p)$, of the subspaces $T_p^0X$, $T_p^+X$, and $T_p^-X$ of the tangent space $T_pX$ are called the \emph{Bia{\l}ynicki--Birula nullity}, \emph{co-index}, and \emph{index}, respectively, of the point $p$ in $X$ defined by the $\CC^*$ action and they obey
\begin{equation}
  \label{eq:Dim_Tp_X_equals_BB_nullity_plus_coindex_plus_index}
  \dim T_pX = \beta_X^0(p) + \beta_X^+(p) + \beta_X^-(p).
\end{equation}  
\end{maindefn}

By Remark \ref{rmk:Weight-sign_decomposition_restriction_tangent_bundle}, the Bia{\l}ynicki--Birula nullity, co-index, and index in Definition \ref{maindefn:Stable_unstable_submanifolds_BB_index_co-index_nullity} are equivalently given by
\begin{subequations}
  \label{eq:BB_nullity_co-index_index_complex_manifold}
  \begin{align}
     \label{eq:BB_nullity_complex_manifold}
    \beta_X^0(p) &= \dim T_p^0X = \dim X_\alpha^0,
    \\
     \label{eq:BB_co-index_complex_manifold}
    \beta_X^+(p) &= \dim T_p^+X = \dim X_p^+,
    \\
     \label{eq:BB_index_complex_manifold}
    \beta_X^-(p) &= \dim T_p^-X = \dim X_p^-,
  \end{align}
\end{subequations}  
where, for brevity, we write $X_p^\pm = X_\alpha^\pm|_p$ for the fibers of $X_\alpha^\pm$ over a point $p \in X_\alpha^0$. We have the following generalization of results due to Bia{\l}ynicki--Birula \cite{Bialynicki_1973}, Carrell and Sommese \cite{Carrell_Sommese_1978ms, Carrell_Sommese_1979cmh, Carrell_Sommese_1983}, Fujiki \cite{Fujiki_1979}, and Yang \cite{Yang_2008} (see the forthcoming Remark \ref{rmk:BB_decomposition_complex_manifold_C*_action_previous_versions}). 

\begin{mainthm}[Bia{\l}ynicki--Birula decomposition for a complex manifold with a $\CC^*$ action]
\label{mainthm:BB_decomposition_complex_manifold_C*_action}
Let $(X,g,J)$ be a real analytic, almost Hermitian manifold with fundamental two-form $\omega = g(\cdot,J\cdot)$ as in \eqref{eq:Fundamental_two-form} and $\CC^* \times X \to X$ be a real analytic action of $\CC^*$ by real analytic diffeomorphisms of $X$ with orbits
\begin{equation}
  \label{eq:Orbit_point_under_C*_action_complex_manifold}
  A_z:\CC^* \ni \lambda \mapsto \lambda\cdot z \in X, \quad\text{for } z \in X.
\end{equation}
Assume furthermore that
\begin{itemize}
\item \emph{(Hamiltonian circle action.)} The real analytic circle action $S^1\times X \to X$ obtained by restricting the $\CC^*$ action to $S^1\subset \CC^*$ is Hamiltonian in the sense of the forthcoming \eqref{eq:MomentMap}, so $\iota_\Theta\omega = df$ on $X$ for a real analytic function $f:X\to\RR$, where the real analytic vector field $\Theta$ on $X$ given by the
forthcoming \eqref{eq:Vector_field_generator_circle_action} is the generator of the circle action, and
\item \emph{(Precompact orbits.)} For each point $z\in X$, at least one of the orbits $A_z(D)$ or $A_z(\CC\less D)$ is precompact in $X$, where $D\subset \CC$ is the open unit disk centered at the origin.
\end{itemize}
Then the following hold:
\begin{enumerate}
\item\label{item:BB_complex_manifold_limits_z0_zinfty_exist_fixed_points}
  For each point $z \in X$, at least one of the limits
\[
    z_0 := \lim_{\lambda\to 0} \lambda\cdot z
    \quad\text{or}\quad
    z_\infty := \lim_{\lambda\to \infty} \lambda\cdot z
\]
exists in $X$, is a fixed point of the $\CC^*$ action on $X$, and the corresponding extended map $A_z:\CC\to X$ or $A_z:\CC\cup\{\infty\}\to X$ is $\CC^*$ equivariant and continuous. In particular, the subset $X^0 \subset X$ of fixed points of the $\CC^*$ action on $X$ is non-empty.
\item\label{item:BB_complex_manifold_Az_holomorphic_BB_decomposition_exists}
If in addition the almost complex structure $J$ is integrable, then the map $A_z:\CC\to X$ or $A_z:\CC\cup\{\infty\}\to X$ is holomorphic for each $z\in X$ and the complex manifold $X$ admits a meromorphic Bia{\l}ynicki--Birula decomposition in the sense of Definition \ref{maindefn:BB_decomposition_complex_manifold}.
\end{enumerate}
\end{mainthm}


\begin{rmk}[Previous versions of Item \eqref{item:BB_complex_manifold_Az_holomorphic_BB_decomposition_exists} in Theorem \ref{mainthm:BB_decomposition_complex_manifold_C*_action}]
\label{rmk:BB_decomposition_complex_manifold_C*_action_previous_versions}
Bia{\l}ynicki--Birula obtained Item \eqref{item:BB_complex_manifold_Az_holomorphic_BB_decomposition_exists} in Theorem \ref{mainthm:BB_decomposition_complex_manifold_C*_action} via his \cite[Theorem 4.1, p. 492]{Bialynicki_1973} under the hypothesis that $X$ is a smooth, complex, projective algebraic variety (see Theorem \ref{thm:Milne_13-47} for a slightly more general statement).

When $X$ is a compact, complex manifold with a holomorphic action $\CC^*\times X\to X$, Fujiki proved the existence of an \emph{analytic} Bia{\l}ynicki--Birula decomposition in the sense of Definition \ref{maindefn:BB_decomposition_complex_manifold} via his \cite[Lemma 2.1, p. 804 and Proposition 2.8, p. 812]{Fujiki_1979}. While Fujiki assumes that $X$ is compact, he only uses that hypothesis in \cite[Proposition 2.2, p. 231]{Fujiki_1978im} to show that a holomorphic action $\CC^*\times X\to X$ extends to a map $\PP^1\times X\to X$ that is \emph{meromorphic} in the sense of \cite[Definition 2.1, p. 230]{Fujiki_1978im}. This in turn implies that for each $z\in X$, the orbit map $A_z:\CC^* \to X$ in \eqref{eq:Orbit_point_under_C*_action_complex_manifold} extends uniquely to a holomorphic map $A_z:\PP^1 \to X$ (see Fujiki \cite[Section 2, p. 803]{Fujiki_1979}). In his \cite[Lemma II-A, p. 109]{Sommese_1975}, Sommese proved the latter result under the additional hypothesis that $X$ is K\"ahler. Thus if the conclusion in Item \eqref{item:BB_complex_manifold_Az_holomorphic_BB_decomposition_exists} is relaxed to that of existence of an \emph{analytic} Bia{\l}ynicki--Birula decomposition, Theorem \ref{mainthm:BB_decomposition_complex_manifold_C*_action} shows that a hypothesis that $X$ is K\"ahler is unnecessary and that a hypothesis that $X$ is compact can be weakened to our hypothesis on precompactness of orbits.

To obtain a \emph{meromorphic} Bia{\l}ynicki--Birula decomposition decomposition, Fujiki used his additional hypothesis that $X$ is K\"ahler only in the statement and proof of his \cite[Proposition 2.10, p. 815]{Fujiki_1979}. The statements and proofs of his \cite[Lemma 2.11, p. 816 and Theorem 2.2, p. 805]{Fujiki_1979} only rely on his hypothesis that $X$ is K\"ahler through his application of \cite[Proposition 2.10, p. 815]{Fujiki_1979}. We essentially replace his hypothesis that $X$ is K\"ahler by our hypothesis that the induced circle action $S^1\times X\to X$ is Hamiltonian. Thus Item \eqref{item:BB_complex_manifold_Az_holomorphic_BB_decomposition_exists} in Theorem \ref{mainthm:BB_decomposition_complex_manifold_C*_action} shows that a hypothesis that $X$ is K\"ahler is unnecessary and that an assumption that $X$ is compact can be weakened to the condition on precompactness of orbits. In \cite[Proposition II, p. 55]{Carrell_Sommese_1978ms}, Carrell and Sommese also prove Item \eqref{item:BB_complex_manifold_Az_holomorphic_BB_decomposition_exists} in Theorem \ref{mainthm:BB_decomposition_complex_manifold_C*_action} but, like Fujiki, only under the stronger hypotheses that $X$ is compact and K\"ahler.

Yang \cite[Theorem 4.12, p. 92]{Yang_2008} proved Item \eqref{item:BB_complex_manifold_Az_holomorphic_BB_decomposition_exists} in Theorem \ref{mainthm:BB_decomposition_complex_manifold_C*_action} for holomorphic $\CC^*$ actions on complex K\"ahler manifolds that need not be compact but where the induced circle action $S^1\times X\to X$ has a Hamiltonian function that is proper and bounded below and the number of connected components of the fixed-point subset $X^{\CC^*}$ is finite.
\end{rmk}  

The Bia{\l}ynicki--Birula decomposition enjoys useful functorial properties and we shall describe the principal ones below, beginning with the elementary

\begin{mainthm}[Bia{\l}ynicki--Birula decomposition for a properly embedded complex submanifold of a complex manifold with a $\CC^*$ action]
\label{mainthm:BB_decomposition_C*_invariant_complex_submanifold}
Let $X$ be a finite-dimensional complex manifold, $Y \subset X$ be a properly embedded complex submanifold, and $\CC^*\times X \to X$ be a holomorphic $\CC^*$ action such that $Y$ is $\CC^*$-invariant with at least one fixed point. If $X$ admits a plus (respectively, minus or mixed) Bia{\l}ynicki--Birula decomposition as in Definition \ref{maindefn:BB_decomposition_complex_manifold} with subsets $X^0$, $X^\pm$, $X_\alpha^\pm$, then $Y$ inherits a plus (respectively, minus or mixed) Bia{\l}ynicki--Birula decomposition with corresponding subsets
\begin{equation}
\label{eq:BB_decomposition_C*_invariant_complex_submanifold}  
  Y^0 = Y\cap X^0, \quad Y^\pm = Y\cap X^\pm, \quad\text{and}\quad Y_p^\pm = Y\cap X_p^\pm,
  \quad\text{for all } p \in Y^0,
\end{equation}
and the Bia{\l}ynicki--Birula nullity, co-index, and index of $p$ in $Y$ obey the identity
\begin{equation}
  \label{eq:Dim_Tp_Y_equals_BB_nullity_plus_coindex_plus_index}
  \dim T_pY = \beta_Y^0(p) + \beta_Y^+(p) + \beta_Y^-(p), \quad\text{for all } p \in Y^0.
\end{equation}
\end{mainthm}

The following result plays a crucial role in our application of Resolution of Singularities.

\begin{mainthm}[Bia{\l}ynicki--Birula decomposition for the blowup of a complex manifold with a $\CC^*$ action along an invariant complex submanifold]
\label{mainthm:BB_decomposition_blowup_complex_manifold_C*_action_along_submanifold}  
Let $X$ be a finite-dimensional complex manifold and $Z \subset X$ be an embedded complex submanifold. If $\Phi:\CC^* \to \Aut(X)$ is a homomorphism from the group $\CC^*$ onto a subgroup of the group $\Aut(X)$ of biholomorphic automorphisms of $X$ such that the corresponding action $\CC^*\times X\to X$ is holomorphic and leaves $Z$ invariant, then the following hold for the blowup $\Bl_Z(X)$ of $X$ along $Z$:
\begin{enumerate}
\item\label{item:BlZX_complex_manifold}
$\Bl_Z(X)$ is a complex manifold of dimension equal to that of $X$.
\item\label{item:BlZX_C*_action_lifts_to_BlZX_from_X_complex_manifold}
The homomorphism $\Phi:\CC^* \to \Aut(X)$ lifts uniquely to a homomorphism $\Bl_Z(\Phi):\CC^* \to \Aut(\Bl_Z(X))$ such that the blowup morphism $\pi:\Bl_Z(X)\to X$ is $\CC^*$-equivariant, the following diagram commutes (for $\pi_{\Aut}$ defined by the forthcoming Proposition \ref{prop:Functorial_property_blowup_complex_manifold}), and the corresponding action $\CC^*\times \Bl_Z(X)\to \Bl_Z(X)$ is holomorphic:
\[
  \begin{tikzcd}
    &\Aut(\Bl_Z(X)) \arrow[d, "\pi_{\Aut}"]
    \\
    \CC^* \arrow[r, "\Phi"] \arrow[ru, "\Bl_Z(\Phi)"] &\Aut(X)
  \end{tikzcd}
\]  
\end{enumerate}  
If in addition the subset $X^0 \subset X$ of fixed points of the $\CC^*$ action on $X$ is non-empty, then the following hold:
\begin{enumerate}
\setcounter{enumi}{2}  
\item\label{item:BlZX0_non-empty_complex_manifold}  
The subset $\Bl_Z(X)^0 \subset \Bl_Z(X)$ of fixed points of the $\CC^*$ action on $\Bl_Z(X)$ is non-empty.
\end{enumerate}  
If in addition $X$ admits a plus (respectively, minus or mixed) Bia{\l}ynicki--Birula decomposition in the sense of Definition \ref{maindefn:BB_decomposition_complex_manifold}, then the following hold:
\begin{enumerate}
\setcounter{enumi}{3}  
\item\label{item:BB_decomposition_BlZX_complex_manifold}
  The blowup $\Bl_Z(X)$ inherits a plus (respectively, minus or mixed) Bia{\l}ynicki--Birula decomposition in the sense of Definition \ref{maindefn:BB_decomposition_complex_manifold}.
\item\label{item:piBlZXpm0_is_Xpm0_complex_manifold}
  The following identities hold:
\begin{subequations}
  \label{eq:piBlZXpm_is_Xpm_complex_manifold}
  \begin{align}
    \label{eq:piBlZX0_is_X0_complex_manifold}
    \pi\left(\Bl_Z(X)^0\right) &= X^0,
    \\
    \label{eq:piBlZX+_is_X+_complex_manifold}
    \pi\left(\Bl_Z(X)^+\right) &= X^+,
    \\
    \label{eq:piBlZX-_is_X-_complex_manifold}
  \pi\left(\Bl_Z(X)^-\right) &= X^-.
  \end{align}
\end{subequations}
Moreover, for all points $\tilde p \in \Bl_Z(X)^0$ and $p = \pi(\tilde p) \in X^0$, the following identities hold:
\begin{subequations}
  \label{eq:piBlZXtildep_pm_is_Xp_pm_complex_manifold}
  \begin{align}
    \label{eq:piBlZXtildep+_is_Xp+_complex_manifold}
   \pi\left(\Bl_Z(X)_{\tilde p}^+\right) &= X_p^+,
    \\
    \label{eq:piBlZXtildep-_is_Xp-_complex_manifold}
   \pi\left(\Bl_Z(X)_{\tilde p}^-\right) &= X_p^-.
  \end{align}
\end{subequations}  
\item\label{item:Dimension_BlZXalphapm0_complex_manifold}  
The complex dimensions of a component $\Bl_Z(X)_\alpha^0$ of the fixed-point subset $\Bl_Z(X)^0$ and the fibers $\Bl_Z(X)_{\alpha,\tilde p}^\pm$ of the natural projections $\pi_\alpha^\pm:\Bl_Z(X)_\alpha^\pm \to \Bl_Z(X)_\alpha^0$ over points $\tilde p \in \Bl_Z(X)_\alpha^0$ with $\pi(\tilde p) = p \in X^0$ are given by
\begin{subequations}
  \label{eq:Complex_dimension_BlZX0_fibers_pm_complex_manifold}
  \begin{align}
    \label{eq:Complex_dimension_BlZX0_complex_manifold}
    \dim\Bl_Z(X)_\alpha^0 &= \dim T_pX^0,
    \\
    \label{eq:Complex_dimension_BlZXfibers_pm_complex_manifold}
    \dim\Bl_Z(X)_{\alpha,\tilde p}^\pm &= \dim X_p^\pm.
  \end{align}
\end{subequations}
In particular, the Bia{\l}ynicki--Birula signature is preserved by blowup:
\begin{equation}
  \label{eq:BBsignature_preserved_by_blowup}
  \beta_{\Bl_Z(X)}^0(\tilde p) = \beta_X^0(p)
  \quad\text{and}\quad
  \beta_{\Bl_Z(X)}^\pm(\tilde p) = \beta_X^\pm(p).
\end{equation}
\end{enumerate}
\end{mainthm}

\begin{rmk}[Bia{\l}ynicki--Birula decompositions need not commute with strict transforms]
\label{rmk:BB_decompositions_blowup_complex_manifold_need_not_commute_with_strict_transforms}  
As explained in Remark \ref{rmk:Submanifolds_blowup_strict_transforms} in a simpler setting, the fixed-point submanifolds $\Bl_Z(X)_\alpha^0$ or stable or unstable fibers $\Bl_Z(X)_{\alpha,\tilde p}^\pm$ in $\Bl_Z(X)$ provided by Theorem
\ref{mainthm:BB_decomposition_blowup_complex_manifold_C*_action_along_submanifold} need \emph{not} coincide with the strict transforms of (a connected component of) $X^0$ or $X_p^\pm$, where $\tilde p \in \pi^{-1}(p)$ and $\pi:\Bl_Z(X)\to X$ is the blowup morphism.
\end{rmk}

\section{Bia{\l}ynicki--Birula decompositions for complex analytic spaces}
\label{sec:BB_decomposition_complex_analytic_spaces}
The results that we described in Section \ref{sec:Functorial_properties_BB_decompositions_complex_manifolds} provide Bia{\l}ynicki--Birula decompositions for complex manifolds. We now proceed to the main focus of our work, which concerns aspects of Bia{\l}ynicki--Birula theory that can be recovered for complex analytic spaces that are not necessarily complex manifolds. 
We begin with the following weaker version of Definition \ref{maindefn:BB_decomposition_complex_manifold}, where the assumption that $X$ is a complex manifold is relaxed to that of a complex analytic space.

\begin{maindefn}[Bia{\l}ynicki--Birula decomposition for a holomorphic $\CC^*$ action on a complex analytic space]
\label{maindefn:BB_decomposition_complex_analytic_space}
Let $(X,\sO_X)$ be a complex analytic space as in Definition \ref{defn:Analytic_space} and $\CC^*\times X\to X$ be a holomorphic action as in Section \ref{sec:Analytic_group_actions_analytic_spaces} such that the subset $X^0\subset X$ of fixed points of the $\CC^*$ action is non-empty with at most countably many connected components, $X_\alpha^0$ for $\alpha\in\sA$, that are locally closed complex analytic subspaces of $X$ as in the forthcoming Definition \ref{defn:Analytic_subspace}. For each $\alpha\in\sA$, define $X_\alpha^\pm$ as in \eqref{eq:Complex_manifold_Xalpha_plus_minus_submanifolds} and the natural maps $\pi_\alpha^\pm$ as in \eqref{eq:pi_Xalpha_pm_to_Xalpha_0_bundle_maps}. Then $X$ has a \emph{(mixed, plus, or minus) Bia{\l}ynicki--Birula decomposition} if the following hold:
\begin{enumerate}
\item\label{item:Xalpha_pm_complex_analytic_subspace_X}
  Each $X_\alpha^+$ is a locally closed, complex analytic subspace of $X$;
\item\label{item:pi_Xalpha_pm_to_Xalpha_0_holomorphic_map}
  The map $\pi_\alpha^+:X_\alpha^+\to X_\alpha^0$ is a $\CC^*$-equivariant epimorphism of complex analytic spaces;
\item\label{item:Xalpha_0_section_Xalpha_pm_complex_analytic_subspace}
  $X_\alpha^0$ is a section of $X_\alpha^+$;
\end{enumerate}  
and the analogous properties hold for the subsets $X_\alpha^-$ and for the natural maps $\pi_\alpha^-:X_\alpha^-\to X_\alpha^0$.
Furthermore, we require that:
\begin{enumerate}
\setcounter{enumi}{3}  
\item $X$ may be expressed as a disjoint union as in \eqref{eq:Complex_manifold_mixed_decomposition}, for some subset $A \subset \sA\times \{+,-\}$.
\end{enumerate}
If one can express $X$ as a disjoint union as in \eqref{eq:Complex_manifold_pure_plus_or_minus_decomposition}, then we say that $X$ has a \emph{plus} or \emph{minus decomposition}, respectively and, otherwise, if $X$ is expressed as in \eqref{eq:Complex_manifold_mixed_decomposition}, that it has a \emph{mixed decomposition}.
\qed
\end{maindefn}


\begin{maindefn}[Stable and unstable subspaces of a complex analytic space and Bia{\l}ynicki--Birula index, co-index, and nullity]
\label{maindefn:Stable_unstable_subspaces_BB_index_co-index_nullity}
Continue the notation of Definition \ref{maindefn:BB_decomposition_complex_analytic_space}. The locally closed complex analytic subspace $X_\alpha^+$ (respectively, $X_\alpha^-$) is called the \emph{stable} (respectively, \emph{unstable}) \emph{subspace} for the fixed-point subspace $X_\alpha^0$. For each point $p\in X_\alpha^0$, the Krull dimensions $\beta_X^0(p)$, $\beta_X^+(p)$, and $\beta_X^-(p)$, of the local rings $\sO_{X_\alpha^0,p}$, $\sO_{X_p^+,p}$, and $\sO_{X_p^-,p}$ are called the \emph{Bia{\l}ynicki--Birula nullity}, \emph{co-index}, and \emph{index}, respectively, of the point $p$ in $X$ defined by the $\CC^*$ action, where we write $X_p^\pm = X_\alpha^\pm|_p$ for the fibers of $X_\alpha^\pm$ over a point $p \in X_\alpha^0$.
\end{maindefn}

Definition \ref{maindefn:Stable_unstable_subspaces_BB_index_co-index_nullity} is a generalization of Definition \ref{maindefn:Stable_unstable_submanifolds_BB_index_co-index_nullity} in the sense that, if $p$ is a smooth point of $X$ as in the forthcoming Definition \ref{defn:Analytic_space}, then the forthcoming Theorem \ref{thm:Comparison_smoothness_regularity_analytic_spaces} yields an open neighborhood $U\subset X$ of $p$ such that $X\cap U$ is a complex manifold of dimension equal to the Krull dimension of the local ring $\sO_{X,p}$ and the subspaces $X_\alpha^0\cap U$ and $X_p^\pm\cap U$ are embedded complex submanifolds of $U$ with dimensions
\begin{align*}
  \beta_X^0(p) &= \dim T_p^0X = \dim (X_\alpha^0\cap U),
  \\
  \beta_X^+(p) &= \dim T_p^+X = \dim (X_p^+\cap U),
  \\
  \beta_X^-(p) &= \dim T_p^-X = \dim (X_p^+\cap U),
\end{align*}
just as in Definition \ref{maindefn:Stable_unstable_submanifolds_BB_index_co-index_nullity}.

\begin{rmk}[Krull dimension and the restriction to complex analytic spaces in Definition \ref{maindefn:Stable_unstable_submanifolds_BB_index_co-index_nullity}]
\label{rmk:Krull_dimension_and_restriction_to_complex_analytic_spaces}  
While a version of a Bia{\l}ynicki--Birula decomposition can be obtained for $\CC^*$ actions on almost complex rather than complex manifolds (for example, see Proposition \ref{prop:Extension_C*_actions_real_analytic_almost_Hermitian_manifolds}), the primary reason for restricting to complex manifolds or complex analytic spaces in the statements of our main results is that positivity of the Krull dimension for a local ring corresponding to a point in an algebraic variety or analytic space over a field $\KK$ is not useful in our applications when the field $\KK$ is not algebraically closed. In particular, while Theorem \ref{thm:Narasimhan_section_3-1_theorem_1_p_41} and Corollary \ref{cor:Generic_smoothness_complex_analytic_space} hold for $\KK=\CC$, they do not hold for $\KK=\RR$, as mentioned in Remark \ref{rmk::Narasimhan_section_3-1_theorem_1_p_41_restriction_to_C} and as explained in Example \ref{exmp:Hypersurface}.
\end{rmk}  

\begin{mainthm}[Bia{\l}ynicki--Birula decomposition for a closed complex analytic subspace of a complex manifold with a $\CC^*$ action]
\label{mainthm:BB_decomposition_C*_invariant_complex_analytic_subspace}
Let $X$ be a finite-dimensional complex manifold, $(Y,\sO_Y)$ be a closed complex analytic subspace of $X$ as in the forthcoming Definition \ref{defn:Analytic_subspace}, and $\CC^*\times X\to X$ be a holomorphic action on $X$ that leaves $Y$ invariant with at least one fixed point in $Y$.
\begin{enumerate}
\item\label{item:BB_decomposition_C*_invariant_complex_analytic_subspace} 
If $X$ admits a plus (respectively, minus or mixed) Bia{\l}ynicki--Birula decomposition in the sense of Definition \ref{maindefn:BB_decomposition_complex_manifold} with subsets $X^0$, $X^\pm$, $X_\alpha^\pm$, then $Y$ inherits a plus (respectively, minus or mixed) Bia{\l}ynicki--Birula decomposition as in Definition \ref{maindefn:BB_decomposition_complex_analytic_space} with corresponding locally closed complex analytic subspaces:
\begin{equation}
  \label{eq:BB_decomposition_C*_invariant_complex_analytic_subspace}
  Y^0 = Y\cap X^0, \quad Y^\pm = Y\cap X^\pm, \quad\text{and}\quad Y_p^\pm = Y\cap X_p^\pm,
  \quad\text{for all } p \in Y^0.
\end{equation}
\item\label{item:BB_decomposition_C*_invariant_complex_analytic_subspace_smooth_point}
If $p \in Y$, then there is an open neighborhood $U \subset Y$ of $p$ such that $\dim (Y_\sm\cap U) = \dim\sO_{Y,p}$, where $Y_\sm \subset Y$ denotes the subset of smooth points as in Definition \ref{defn:Smooth_point_analytic_space}.
\item\label{item:BB_dimensions_nullity_co-index_index_complex_analytic_subspace}
If $p \in Y^0$ then, after possibly shrinking $U$, the Bia{\l}ynicki--Birula nullity, co-index, and index of $p$ in $Y$ in the sense of Definition \ref{maindefn:Stable_unstable_subspaces_BB_index_co-index_nullity} are given by
\begin{subequations}
\label{eq:BB_dimensions_nullity_co-index_index_complex_analytic_subspace}  
\begin{align}
  \label{eq:BB_nullity_complex_analytic_subspace}
  \beta_Y^0(p) &= \dim\sO_{Y^0,p} = \dim(Y^0)_\sm\cap U,
  \\
  \label{eq:BB_co-index_complex_analytic_subspace}
  \beta_Y^+(p) &= \dim\sO_{Y_p^+,p} = \dim(Y_p^+)_\sm\cap U,
  \\
  \label{eq:BB_index_complex_analytic_subspace}
  \beta_Y^-(p) &= \dim\sO_{Y_p^-,p} = \dim(Y_p^+)_\sm\cap U,
\end{align}
\end{subequations}
where $(Y^0)_\sm \subset Y^0$ and $(Y_p^\pm)_\sm \subset Y_p^\pm$ denote subsets of smooth points as in the forthcoming Definition \ref{defn:Smooth_point_analytic_space}.
\item\label{item:BB_decomposition_C*_invariant_complex_analytic_subspace_nonempty_BB_subspaces}
If $\beta_Y^0(p) > 0$ (respectively, $\beta_Y^+(p) > 0$ or $\beta_Y^-(p) > 0$), then $(Y^0)_\sm\cap U$ (respectively, $(Y_p^+)_\sm\cap U$ or $(Y_p^-)_\sm\cap U$) is non-empty.
\item\label{item:BB_index_coindex_at_p_positive_basic_implies_not_local_min_max}
If the induced circle action $S^1\times X \to X$ has a Hamiltonian function $f:X\to\RR$ in the sense of the forthcoming \eqref{eq:MomentMap} and $\beta_Y^-(p) > 0$ (respectively, $\beta_Y^+(p) > 0$), then $p$ is not a local minimum (respectively, maximum) of the restriction $f:Y_\sm\cup\{p\}\to \RR$. 
\end{enumerate}
\end{mainthm}

\begin{rmk}[On the hypotheses of Theorem \ref{mainthm:BB_decomposition_C*_invariant_complex_analytic_subspace}]
\label{rmk:BB_decomposition_C*_invariant_complex_analytic_subspace_hypotheses}
A result of Sumihiro \cite[Theorem 1, p. 5]{Sumihiro_1974} asserts that, over an algebraically closed field, any quasi-projective, normal algebraic variety with a regular action of a linear algebraic group can be embedded equivariantly in a projective space with a linear action. Konarski \cite{Konarski_1978} generalized the results of Bia{\l}ynicki--Birula \cite{Bialynicki_1973} to the case of complete, normal algebraic varieties over an algebraically closed field $k$ equipped with an action by the multiplicative group $\GG_m$ of units in $k$. Weber \cite[Section 2, p. 539]{Weber_2017} studied the Bia{\l}ynicki--Birula decomposition for singular complex algebraic varieties with $\CC^*$ actions, but under the assumption that they are $\CC^*$-equivariantly embedded in complete, complex, smooth algebraic varieties. Drinfeld \cite{Drinfeld_2013arxiv_v2} considered arbitrary algebraic spaces over an arbitrary field $k$ equipped with a $\GG_m$ action. Related results in the category of algebraic varieties were obtained by Carrell and Goresky \cite{Carrell_Goresky_1983}, Gonzales \cite{Gonzales_2014, Gonzales_2016}, Hausel and Hitchin \cite{Hausel_Hitchin_2022}, Jelisiejew and Sienkiewicz \cite{Jelisiejew_Sienkiewicz_2019}, and Kirwan \cite{Kirwan_1988}.

Bia{\l}ynicki--Birula and Sommese \cite{Bialynicki-Birula_Sommese_1983} considered irreducible, normal complex analytic spaces $X$ but under the additional assumptions that the holomorphic action $\CC^*\times X \to X$ is meromorphic (that is, extends to a meromorphic map $\PP_\CC^1\times X \to X$) and locally linearizable around fixed points (that is, for each $p \in X^{\CC^*}$ there exist an integer $n$ and $\CC^*$-equivariant, holomorphic, proper embedding from an open neighborhood of $p$ into $\CC^n$ equipped with a linear $\CC^*$ action).
\end{rmk}  

\begin{rmk}[Estimating Krull dimensions]
\label{rmk:Estimating_Krull_dimensions_expected_dimensions} 
Suppose $(U,\sO_U)$ is a local model space (see the forthcoming Definition \ref{defn:Analytic_model_space}) for an open neighborhood of a point $p$ in a complex analytic space $(X,\sO_X)$ (see the forthcoming Definition \ref{defn:Analytic_space}), so $U$ is the topological support of $\sO_D/\sI$ with a domain\footnote{By a \emph{domain} in $\KK^n$, where $\KK=\RR$ or $\CC$, we as usual mean a connected open subset.} $D \subset \CC^n$ and ideal $\sI \subset \sO_D$ with generators $f_1,\ldots,f_r$ and structure sheaf $\sO_U := (\sO_D/\sI)\restriction U$. According to the forthcoming Lemma \ref{lem:Algebraic_dimension_pX_geq_expdim_pX}, we have the inequality
\[
  \dim\sO_{X,p} \geq n - r,
\]
where $\expdim_p X := n - r$ is the \emph{expected dimension} of $X$ at $p$ as defined in the forthcoming Remark \ref{rmk:Expected_dimension}. When $r$ is equal to the \emph{minimal number} of generators of $\sI_p \subset \sO_{U,p}$, then the Krull dimension \eqref{eq:Krull_dimension} of $X$ at $p$ is given by
\[
  \dim\sO_{X,p} = n - r.
\]
See the forthcoming Remark \ref{rmk::Narasimhan_section_3-1_theorem_1_p_41_meaning_of_d} for an explanation of the preceding equality.
\end{rmk}

\begin{rmk}[Local complete intersections]
\label{rmk:Local_complete_intersections}
Continue the notation of Remark \ref{rmk:Estimating_Krull_dimensions_expected_dimensions}. If $U$ is \emph{pure dimensional} at $p$ (as in Grauert and Remmert \cite[Section 5.4.2, p. 106]{Grauert_Remmert_coherent_analytic_sheaves}) so that, after possibly shrinking $U$, we have $\dim\sO_{U,x} = d$ for all $x \in U$, where $d := \dim\sO_{X,p}$, then $(U,\sO_U)$ is a \emph{local complete intersection} (as in Hartshorne \cite[Chapter II, Section 8, Definition, p. 185]{Hartshorne_algebraic_geometry}, De Jong and Pfister \cite[Definition 4.1.18 (1), p. 134]{DeJong_Pfister_local_analytic_geometry}, or Ebelin \cite[Section 2.10, Definition, p. 106]{Ebeling_functions_several_complex_variables_singularities}) since $\sI$ has $r = \codim(U,D) = n - d$ generators. In general, local complete intersections need not be smooth (see Fulton and Lang \cite[Section IV.3, p. 86]{Fulton_Lang_riemann-roch_algebra}).

For a well-known example related to gauge theory, we note that Donaldson \cite{DonPoly} (quoted by Friedman and Morgan in \cite[Section 4.4.3, Theorem 4.10, p. 339]{FrM}) has shown that the moduli space of rank two, stable holomorphic bundles of fixed determinant is smooth as a scheme and of the expected dimension at a generic point in every component when the second Chern number of the bundles is sufficiently large. See also Huybrechts and Lehn \cite[Proposition 2.A.11, p. 57]{Huybrechts_Lehn_geometry_moduli_spaces_sheaves}.
\end{rmk}

While the Krull dimension at a point of a complex analytic space or a scheme may be difficult to compute, it can be estimated via the inequalities
\[
  \dim T_pX \geq \dim\sO_{X,p} \geq \expdim_pX,
\]
where $T_pX$ is the Zariski tangent space to $X$ at $p$ and the expected dimension $\expdim_pX$ is computed as in Remark \ref{rmk:Estimating_Krull_dimensions_expected_dimensions} with respect to some (convenient) local model for $X$ near $p$. In particular, while the Bia{\l}ynicki--Birula nullity, co-index, and index may be difficult to compute, one can define useful lower bounds for them that may be readily calculated in practice (see Remark \ref{rmk:Calculation_virtual_BB_nullity_co-index_index}).

\begin{maindefn}[Virtual Bia{\l}ynicki--Birula nullity, co-index, and index of a fixed point in a complex analytic space with a holomorphic $\CC^*$ action]
\label{maindefn:Virtual_BB_nullity_co-index_index}  
Let $(X,\sO_X)$ be a complex analytic space as in the forthcoming Definition \ref{defn:Analytic_space} with a holomorphic $\CC^*$ action such that $X$ admits a plus (respectively, minus or mixed) Bia{\l}ynicki--Birula decomposition as in Definition \ref{maindefn:BB_decomposition_complex_analytic_space}. Let $p \in X^0$, so $p$ is a fixed point of the $\CC^*$ action on $X$, and $(U,\sO_U)$ be a local model for an open neighborhood of $p$ in $X$, defined by
\begin{inparaenum}[\itshape i\upshape)]
\item a domain $D \subset \CC^n$ around the origin where $n$ is the embedding dimension $\embdim_pX$ of $X$ at $p$ as in the forthcoming Definition \ref{defn:Embedding_dimension} and $n = \dim T_pX$ by the forthcoming Proposition \ref{prop:Equality_dimension_zariski_tangent_space_and_embedding_dimension}, and
\item an ideal $\sI \subset \sO_D$ with generators $f_1,\ldots,f_r$, and
\item topological support $U = \supp(\sO_D/\sI) \subset D$, and
\item structure sheaf $\sO_U := (\sO_D/\sI)\restriction U$.
\end{inparaenum}
Let $T_pX = \CC^n$ have the linear isotropy $S^1$ action, with corresponding weight-sign decomposition $T_pX = T_p^0X \oplus T_p^+X \oplus T_p^-X$. Let $\Xi_p = \CC^r$ have the linear $S^1$ action induced by the holomorphic map $F = (f_1,\ldots,f_r):D \to \Xi$ as in the forthcoming Lemma \ref{lem:Blanchard_holomorphic_map_from_domain_around_fixed_point_into_vector_space}, with corresponding weight-sign decomposition, $\Xi_p = \Xi_p^0\oplus\Xi_p^+\oplus\Xi_p^-$. Then the \emph{virtual Bia{\l}ynicki--Birula nullity, co-index}, and \emph{index}, respectively, of the point $p$ in $X$ \emph{relative to the local model space} are defined by
\begin{subequations}
  \label{eq:Virtual_BB_nullity_co-index_index}
  \begin{align}
    \label{eq:Virtual_BB_nullity}
    \expdim_p X^0 := \dim T_p^0X - \dim \Xi_p^0,
    \\
    \label{eq:Virtual_BB_co-index}
    \expdim_p X_p^+ := \dim T_p^+X - \dim \Xi_p^+,
    \\
    \label{eq:Virtual_BB_index}
    \expdim_p X_p^- := \dim T_p^-X - \dim \Xi_p^-.                  
  \end{align}
\end{subequations}
\end{maindefn}

We have the important

\begin{mainthm}[Upper and lower bounds for the Krull dimensions of the local rings of the fixed point subset and stable and unstable subspaces of a complex analytic space with a holomorphic $\CC^*$ action]
\label{mainthm:Upper_and_lower_bounds_Krull_dimensions_X0_Xpm_at_p}
Continue the notation and assumptions of Definition \ref{maindefn:Virtual_BB_nullity_co-index_index}. If $p \in X^0$, then the Krull dimensions of the local rings at $p$ for $X^0$ and $X_p^\pm$ obey the following bounds:
\begin{subequations}
  \label{eq:Upper_and_lower_bounds_Krull_dimensions_X0_Xpm_at_p}
  \begin{align}
    \label{eq:Upper_and_lower_bounds_Krull_dimension_X0}
    \dim T_p^0X &\geq \dim\sO_{X^0,p} \geq \dim T_p^0X - \dim \Xi_p^0,
    \\
    \label{eq:Upper_and_lower_bounds_Krull_dimension_X+_at_p}
    \dim T_p^+X &\geq \dim\sO_{X_p^+,p} \geq \dim T_p^+X - \dim \Xi_p^+,
    \\
    \label{eq:Upper_and_lower_bounds_Krull_dimension_X-_at_p}
    \dim T_p^-X &\geq \dim\sO_{X_p^-,p} \geq \dim T_p^-X - \dim \Xi_p^-.
  \end{align}
\end{subequations}
Moreover, the Bia{\l}ynicki--Birula nullity, co-index, and index obey the following generalization of \eqref{eq:Dim_Tp_X_equals_BB_nullity_plus_coindex_plus_index}:
\begin{equation}
  \label{eq:Krulldim_X_at_p_equals_sum_BBnullity_and_BBcoindex_and_BBindex_at_p}
  \dim\sO_{X,p} \leq \beta_X^0(p) + \beta_X^+(p) + \beta_X^-(p).
\end{equation}
\end{mainthm}

Theorem \ref{mainthm:Upper_and_lower_bounds_Krull_dimensions_X0_Xpm_at_p} shows that the Bia{\l}ynicki--Birula nullity, co-index, and index of a point $p$ as in Definition \ref{maindefn:Stable_unstable_subspaces_BB_index_co-index_nullity} are bounded above by the dimensions of the subspaces $T_p^0X$ and $T_p^\pm X$ of the Zariski tangent space $T_pX$ to $X$ at $p$ defined by the weight-sign decomposition and bounded below by the virtual Bia{\l}ynicki--Birula nullity, co-index, and index, respectively, as in Definition \ref{maindefn:Virtual_BB_nullity_co-index_index}:
\begin{subequations}
\label{eq:BB_nullity_co-index_index_geq_expected_dimensions}    
\begin{align}
\label{eq:BB_nullity_geq_expected_dimension}
  \dim T_p^0X \geq \beta_X^0(p) &\geq \expdim_p X^0,
  \\
  \label{eq:BB_co-index_geq_expected_dimension}
  \dim T_p^+X \geq \beta_X^+(p) &\geq \expdim_p X^+,
  \\
  \label{eq:BB_index_geq_expected_dimension}
  \dim T_p^-X \geq \beta_X^-(p) &\geq \expdim_p X^-.                 
\end{align}  
\end{subequations}
The expected dimensions clearly always obey the following equality:
\begin{equation}
  \label{eq:Expdim_X_at_p_equals_sum_virtual_BB_nullity_index_coindex_index_X_at_p}
  \expdim_p X = \expdim_p X^0 + \expdim_p X^+ + \expdim_p X^-.
\end{equation}
The significance of the virtual Bia{\l}ynicki--Birula nullity, co-index, and index for applications is illustrated by the following

\begin{rmk}[Calculation of virtual Bia{\l}ynicki--Birula nullity, co-index, and index of a fixed point in a complex analytic moduli space with a holomorphic $\CC^*$ action]
\label{rmk:Calculation_virtual_BB_nullity_co-index_index} 
In Feehan and Leness \cite[Chapter 13]{Feehan_Leness_introduction_virtual_morse_theory_so3_monopoles}, we proved that the expected dimensions in \eqref{eq:BB_nullity_co-index_index_geq_expected_dimensions} may be readily computed using the Hirzebruch--Riemann--Roch Theorem (see Hirzebruch \cite{Hirzebruch_topological_methods_algebraic_geometry}) to calculate the Euler characteristic of the elliptic complex for the holomorphic pair equations, when $X$ is the moduli space of stable pairs of holomorphic, rank two vector bundles and sections over a compact, complex, K\"ahler surface. Our monograph \cite{Feehan_Leness_introduction_virtual_morse_theory_so3_monopoles} with Leness is the first step in our program to use gauge theory to prove that
\begin{inparaenum}[\itshape a\upshape)]
\item all compact, complex surfaces of general type obey the Bogomolov--Miyaoka--Yau inequality (see Barth, Hulek, Peters, and Van de Ven \cite[Theorem 4.1, p. 275]{Barth_Hulek_Peters_Van_de_Ven_compact_complex_surfaces}) and, more generally, that
\item all closed, oriented, smooth four-dimensional manifolds with $b_1=0$, odd $b^+ \geq 3$, and of Seiberg--Witten simple type obey the Bogomolov--Miyaoka--Yau inequality.
\end{inparaenum}  
Our present work represents the second step in our program.

When $p$ is a smooth point of $X$ in the sense of the forthcoming Definition \ref{defn:Smooth_point_analytic_space}, then $r = 0$ and $\Xi_p = (0)$, so $\dim T_pX = \dim\sO_{X,p}$, and the inequalities \eqref{eq:BB_nullity_co-index_index_geq_expected_dimensions} reduce to the equalities \eqref{eq:BB_nullity_co-index_index_complex_manifold} in Definition \ref{maindefn:Stable_unstable_submanifolds_BB_index_co-index_nullity}. In this context, Hitchin \cite[Section 7]{Hitchin_1987} computed the Bia{\l}ynicki--Birula nullity, co-index, and index from \eqref{eq:BB_nullity_co-index_index_complex_manifold} in Definition \ref{maindefn:Stable_unstable_submanifolds_BB_index_co-index_nullity} using the Hirzebruch--Riemann--Roch Theorem to calculate the Euler characteristic of the elliptic complex (with only one non-zero term) for the Higgs pair equations, when $X$ is the moduli space of Higgs pairs (for a rank two vector bundle of odd degree) over a closed Riemann surface.

In these examples, the expected dimensions are computed relative to local Kuranishi models (see Kuranishi \cite{Kuranishi} and Atiyah, Hitchin, and Singer \cite{AHS}) for open neighborhoods of fixed points in these moduli spaces. 
\end{rmk}

We return to the general setting of Definition \ref{maindefn:BB_decomposition_complex_analytic_space}. We shall see in Section \ref{sec:Morse-Bott_decompositions_complex_analytic_subspaces_Hermitian_manifolds} that the set of connected components $\{X_\alpha^0\}_{\alpha\in\sA}$ of the fixed-point subset $X^0 \subset X$ may be totally ordered by the values of a Hamiltonian function $f$ for the induced $S^1$ action, when $X$ is a closed, $\CC^*$-invariant, complex analytic subspace of a complex Hermitian manifold equipped with a holomorphic $\CC^*$ action such that the induced action of $S^1\subset\CC^*$ is isometric and Hamiltonian. However, Bia{\l}ynicki--Birula and Sommese \cite{Bialynicki-Birula_Sommese_1983} observed that a partial ordering on $\{X_\alpha^0\}_{\alpha\in\sA}$ may be defined without appealing to values of a Hamiltonian function.

\begin{maindefn}[Partial and total ordering of components of fixed-point subsets]
\label{maindefn:Bialynicki-Birula_and_Sommese_definition_page_776}  
Continue the notation of Definition \ref{maindefn:BB_decomposition_complex_analytic_space}. A component $X_\alpha^0$ is \emph{directly less than} a component $X_\beta^0$ if
\[
  \left(X_\alpha^+ \less X_\alpha^0\right)\cap \left(X_\beta^- \less X_\beta^0\right) \neq \varnothing.
\]
The component $X_\alpha^0$ is \emph{less than} $X_\beta^0$ if there exists a sequence $\alpha_0 = \alpha,\alpha_1,\ldots,\alpha_k = \beta$ such that $X_{\alpha_j}$ is directly less than $X_{\alpha_{j+1}}$ for $j = 0, \ldots,k-1$, in which case one writes $X_\alpha^0 < X_\beta^0$. One writes $X_\alpha^0 \leq X_\beta^0$ if either $X_\alpha^0 < X_\beta^0$ or $\alpha = \beta$. The set $\{X_\alpha^0\}_{\alpha\in\sA}$ is \emph{totally ordered} if every pair of connected components, $X_\alpha^0$ and $X_\beta^0$, is comparable and \emph{partially ordered} otherwise.
\end{maindefn}

\section[Bia{\l}ynicki--Birula decompositions and resolution of singularities]{Bia{\l}ynicki--Birula decompositions for complex analytic spaces and resolution of singularities}
\label{sec:BB_decomposition_embedded_resolution_singularities_complex_analytic_spaces}
The results that we described in Section \ref{sec:BB_decomposition_complex_analytic_spaces} provide Bia{\l}ynicki--Birula decompositions for complex analytic spaces. We now consider additional properties of those decompositions that can be achieved via Resolution of Singularities for complex analytic spaces. We begin with the following generalization of Theorem \ref{mainthm:BB_decomposition_blowup_complex_manifold_C*_action_along_submanifold}. (We may recover the statement of Theorem \ref{mainthm:BB_decomposition_blowup_complex_manifold_C*_action_along_submanifold} from Theorem \ref{mainthm:BB_decomposition_strict_transform_complex_analytic_subspace} by replacing the composition of blowup maps, $\Pi:X'\to X$, by a single blowup map, $\pi:\Bl_Z(X)\to X$, along a center given by a $\CC^*$-invariant, embedded complex submanifold $Z \subset X$.)

\begin{mainthm}[Bia{\l}ynicki--Birula decomposition for the strict transform of a $\CC^*$-invariant, complex analytic subspace of a complex manifold with a $\CC^*$ action]
\label{mainthm:BB_decomposition_strict_transform_complex_analytic_subspace}  
Let $X$ be a finite-dimensional complex manifold and $(Y,\sO_Y)$ be a closed, complex analytic subspace. If $\Phi:\CC^* \to \Aut(X)$ is a homomorphism from the group $\CC^*$ onto a subgroup of the group $\Aut(X)$ of biholomorphic automorphisms of $X$ such that the corresponding action $\CC^*\times X\to X$ is holomorphic and leaves $Y$ invariant, then there is a $\CC^*$-equivariant resolution morphism $\Pi:X' \to X$ in the sense of the forthcoming Theorem \ref{thm:Embedded_resolution_of_singularities_analytic_space} (see also Theorem \ref{thm:Embedded_resolution_of_singularities_algebraic_scheme}) with exceptional divisor $E' \subset X'$ and strict transform $Y' \subset X'$, so that the following hold:
\begin{enumerate}
\item\label{item:X'_and_Y'_complex_manifolds}
$X'$ is a complex manifold with dimension equal to that of $X$ and $Y'$ is an embedded complex submanifold of $X'$ with dimension equal to that of the top stratum of smooth points $Y_\sm \subset Y$. 
\item\label{item:C*_action_lifts_to_X'_from_X}
The homomorphism $\Phi:\CC^* \to \Aut(X)$ lifts uniquely to a homomorphism $\Pi(\Phi):\CC^* \to \Aut(X')$ such that the resolution morphism $\Pi:X'\to X$ is $\CC^*$-equivariant, the following diagram commutes (for $\pi_{\Aut}$ defined by repeated application of the forthcoming Proposition \ref{prop:Functorial_property_blowup_complex_manifold}), the corresponding action $\CC^*\times X'\to X'$ is holomorphic, and $Y'$ is $\CC^*$-invariant:
\end{enumerate}
\[
  \begin{tikzcd}
    &\Aut(X') \arrow[d, "\pi_{\Aut}"]
    \\
    \CC^* \arrow[r, "\Phi"] \arrow[ru, "\Pi(\Phi)"] &\Aut(X)
  \end{tikzcd}
\]  
If in addition the subset $X^0 \subset X$ of fixed points of the $\CC^*$ action on $X$ is non-empty, then the following hold:
\begin{enumerate}
\setcounter{enumi}{2}  
\item\label{item:X'0_non-empty_complex_manifold}  
The subset $X^{\prime,0} \subset X'$ of fixed points of the $\CC^*$ action on $X'$ is non-empty.
\end{enumerate}  
If in addition $X$ admits a plus (respectively, minus or mixed) Bia{\l}ynicki--Birula decomposition in the sense of Definition \ref{maindefn:BB_decomposition_complex_manifold}, then the following hold:
\begin{enumerate}
\setcounter{enumi}{3}  
\item\label{item:BB_decomposition_X'_complex_manifold}
The strict transform $X'$ inherits a plus (respectively, minus or mixed) Bia{\l}ynicki--Birula decomposition in the sense of Definition \ref{maindefn:BB_decomposition_complex_manifold}.
\item\label{item:PiX'0pm_in_Xpm_and_PiX'p'_pm_is_Xp_pm_complex_manifold}
The following identities hold:
\begin{subequations}
  \label{eq:PiX'0pm_in_Xpm_complex_manifold}
  \begin{align}
    \label{eq:PiX'0_is_X0_complex_manifold}
    \Pi\left(X^{\prime,0}\right) &= X^0,
    \\
    \label{eq:PiX'+_is_X+_complex_manifold}
    \Pi\left(X^{\prime,+}\right) &= X^+,
    \\
    \label{eq:PiX'-_is_X-_complex_manifold}
  \Pi\left(X^{\prime,-}\right) &= X^-.
  \end{align}
\end{subequations}
Moreover, for all $p' \in X^{\prime,0}$ and $p = \Pi(p') \in X^0$, the following identities hold:
\begin{subequations}
  \label{eq:PiX'p'_pm_is_Xp_pm_complex_manifold}
  \begin{align}
    \label{eq:PiX'p'+_is_Xp+_complex_manifold}
   \Pi\left(X_{p'}^{\prime,+}\right) &= X_p^+,
    \\
    \label{eq:PiX'p'-_is_Xp-_complex_manifold}
   \Pi\left(X_{p'}^{\prime,-}\right) &= X_p^-.
  \end{align}
\end{subequations}
\item\label{item:Dimension_X'alphapm_X'alpha0_complex_manifold}  
The complex dimensions of a component $X_\alpha^{\prime,0}$ of the fixed-point subset $X^{\prime,0}$ and the fibers $X_{\alpha,p'}^{\prime,\pm}$ of the natural projections $\Pi_\alpha^\pm:X_\alpha^{\prime,\pm} \to X_\alpha^{\prime,0}$ over points $p' \in X_\alpha^{\prime,0}$ with $\Pi(p') = p \in X^0$ obey
\begin{subequations}
  \label{eq:Complex_dimension_X'0_fibers_pm_complex_manifold}
  \begin{align}
    \label{eq:Complex_dimension_X'0_complex_manifold}
    \dim X_\alpha^{\prime,0} &= \dim T_pX^0,
    \\
    \label{eq:Complex_dimension_X'fibers_pm_complex_manifold}
    \dim X_{\alpha,p'}^{\prime,\pm} &= \dim X_p^\pm.
  \end{align}
\end{subequations}
\end{enumerate}
\end{mainthm}

\begin{rmk}[Smoothness of the fixed-point subset and the stable and unstable subspaces in Theorems \ref{mainthm:BB_decomposition_blowup_complex_manifold_C*_action_along_submanifold} and \ref{mainthm:BB_decomposition_strict_transform_complex_analytic_subspace}]
\label{rmk:Smoothness_fixed-point_subset_and_stable_and_unstablesubspaces}  
When the resolution morphism $\Pi:X'\to X$ in Theorem \ref{mainthm:BB_decomposition_strict_transform_complex_analytic_subspace} comprises a single blowup morphism $\pi:X'\to X$ with blowup center $Z \subset X$, as in Theorem \ref{mainthm:BB_decomposition_blowup_complex_manifold_C*_action_along_submanifold}, and the intersections $X^0\cap Z$, $X^\pm\cap Z$, and $X_p^\pm\cap Z$ are smooth, then Corollary \ref{cor:Smoothness_strict_transform_submanifold} would imply that the strict transforms $X^{0,\prime}$, $X^{\pm,\prime}$, and $X_p^{\pm,\prime}$ of $X^0$, $X^\pm$, and $X_p^\pm$ (for $p'\in X^{\prime,0}$ with $p = \Pi(p') \in X^0$) are also smooth. More generally, if the resolution morphism $\Pi:X'\to X$ involves iterated blowups as in Theorem \ref{mainthm:BB_decomposition_strict_transform_complex_analytic_subspace} then, by the same argument, the strict transforms $X^{0,\prime}$, $X^{\pm,\prime}$, and $X_p^{\pm,\prime}$ would be smooth provided the intersections with blowup centers in each intermediate blowup were smooth. One of the strengths of the conclusions of Theorems \ref{mainthm:BB_decomposition_blowup_complex_manifold_C*_action_along_submanifold} and \ref{mainthm:BB_decomposition_strict_transform_complex_analytic_subspace} is that the subvarieties $X^{\prime,0}$, $X^{\prime,\pm}$, and $X_{p'}^{\prime,\pm}$ are smooth \emph{without} any such hypotheses on smoothness of intersections with blowup centers. However, as we noted in the simpler setting of Remark \ref{rmk:Submanifolds_blowup_strict_transforms}, the strict transforms (even after restricting to connected components) of $X^0$, $X^\pm$, and $X_p^\pm$ with respect to the resolution morphism $\Pi:X'\to X$ need \emph{not} be equal to the embedded complex submanifolds $X^{\prime,0}$, $X^{\prime,\pm}$, and $X_{p'}^{\prime,\pm}$ provided by the Bia{\l}ynicki--Birula decomposition for the complex manifold $X'$.
\end{rmk}

We may combine Theorems \ref{mainthm:BB_decomposition_C*_invariant_complex_submanifold} and \ref{mainthm:BB_decomposition_strict_transform_complex_analytic_subspace} to give the

\begin{maincor}[Bia{\l}ynicki--Birula decomposition for the strict transform of a $\CC^*$-invariant, closed complex analytic subspace of a complex manifold with a $\CC^*$ action]
\label{maincor:BB_decomposition_strict_transform_complex_analytic_subspace}
Continue all of the hypotheses and notation of Theorem \ref{mainthm:BB_decomposition_strict_transform_complex_analytic_subspace} and assume in addition that the subset $Y^0 \subset Y$ of fixed points of the $\CC^*$ action on $Y$ is non-empty. Then the following hold:
\begin{enumerate} 
\item\label{item:Y'0_non-empty_strict_transform_complex_analytic_subspace}  
The subset $Y^{\prime,0} \subset Y'$ of fixed points of the $\CC^*$ action on $Y'$ is non-empty.  
\item\label{item:BB_decomposition_strict_transform_complex_analytic_subspace}
The strict transform $Y'$ inherits a plus (respectively, minus or mixed) Bia{\l}ynicki--Birula decomposition in the sense of Definition \ref{maindefn:BB_decomposition_complex_manifold}, with corresponding subsets
\begin{equation}
\label{eq:BB_decomposition_strict_transform_complex_analytic_subspace}  
  Y^{\prime,0} = Y'\cap X^{\prime,0}, \quad Y^{\prime,\pm} = Y'\cap X^{\prime,\pm}, \quad\text{and}\quad Y_{p'}^{\prime,\pm} = Y'\cap X_{p'}^{\prime,\pm},
  \quad\text{for all } p' \in Y^{\prime,0}.
\end{equation}
\item\label{item:PiY'0pm_in_Ypm_strict_transform_complex_analytic_subspace}
The following identities hold:
\begin{subequations}
  \label{eq:PiY'0pm_in_Ypm_strict_transform_complex_analytic_subspace}
  \begin{align}
    \label{eq:PiY'0_is_Y0_strict_transform_complex_analytic_subspace}
    \Pi\left(Y^{\prime,0}\right) &= Y^0,
    \\
    \label{eq:PiY'+_is_Y+_strict_transform_complex_analytic_subspace}
    \Pi\left(Y^{\prime,+}\right) &= Y^+,
    \\
    \label{eq:PiY'-_is_Y-_strict_transform_complex_analytic_subspace}
  \Pi\left(Y^{\prime,-}\right) &= Y^-.
  \end{align}
\end{subequations}
Moreover, for all $p' \in Y^{\prime,0}$ and $p = \Pi(p') \in Y^0$, the following identities hold:
\begin{subequations}
  \label{eq:PiY'p'_pm_is_Yp_pm_strict_transform_complex_analytic_subspace}
  \begin{align}
    \label{eq:PiY'p'+_is_Yp+_strict_transform_complex_analytic_subspace}
   \Pi\left(Y_{p'}^{\prime,+}\right) &= Y_p^+,
    \\
    \label{eq:PiY'p'-_is_Yp-_strict_transform_complex_analytic_subspace}
   \Pi\left(Y_{p'}^{\prime,-}\right) &= Y_p^-.
  \end{align}
\end{subequations}
\item\label{item:BB_dimensions_nullity_co-index_index_strict_transform_complex_analytic_subspace}
For all $p' \in Y^{\prime,0}$ and $p = \Pi(p') \in Y^0$, the Bia{\l}ynicki--Birula nullity, co-index, and index obey
\begin{subequations}
\label{eq:dim_p_Yprime_and_Y_equals_BB_nullity_plus_coindex_plus_index}  
\begin{align}
  \label{eq:dim_p_Yprime_equals_BB_nullity_plus_coindex_plus_index}
  \dim_{p'}Y' &= \beta_{Y'}^0(p') + \beta_{Y'}^+(p') + \beta_{Y'}^-(p'),
  \\
  \label{eq:dim_p_Y_equals_BB_nullity_plus_coindex_plus_index}
  \dim_pY &\leq \beta_Y^0(p) + \beta_Y^+(p) + \beta_Y^-(p),
\end{align}
\end{subequations}
together with
\begin{equation}
  \label{eq:BB_nullity_coindex_index_preserved_resolution_singularities}
  \beta_{Y'}^0(p') \leq \beta_Y^0(p),
  \quad \beta_{Y'}^+(p') \leq \beta_Y^+(p),
  \quad\text{and}\quad \beta_{Y'}^-(p') \leq \beta_Y^-(p).
\end{equation}
\item\label{eq:BB_index_coindex_at_p_positive_implies_not_local_min_max}
If the induced circle action $S^1\times X \to X$ has a Hamiltonian function $f:X\to\RR$ in the sense of the forthcoming \eqref{eq:MomentMap} and $\beta_Y^-(p) > 0$ (respectively, $\beta_Y^+(p) > 0$), then $p$ is not a local minimum (respectively, maximum) of the restriction $f:Y_\sm\cup\{p\}\to \RR$.
\end{enumerate}
\end{maincor}

\begin{rmk}[Bia{\l}ynicki--Birula decompositions need not commute with strict transforms]
\label{rmk:BB_decompositions_blowup_complex_analytic_space_need_not_commute_with_strict_transforms_main}  
Again, as we noted in the simpler setting of Remark \ref{rmk:Submanifolds_blowup_strict_transforms}, the embedded complex submanifolds $Y^{\prime,0}$, $Y^{\prime,\pm}$, and $Y_{p'}^{\prime,\pm}$ in Corollary \ref{maincor:BB_decomposition_strict_transform_complex_analytic_subspace} are \emph{not} necessarily equal to the strict transforms $Y^{0,\prime}$, $Y^{\pm,\prime}$, and $Y_p^{\pm,\prime}$ (as in the forthcoming Definition \ref{defn:Strict_transform_and_total_transform_analytic_spaces}) of $Y^0$, $Y^\pm$, and $Y_p^\pm$, respectively.
\end{rmk}

The important conclusion in Item \eqref{eq:BB_index_coindex_at_p_positive_implies_not_local_min_max} of Corollary \ref{maincor:BB_decomposition_strict_transform_complex_analytic_subspace} repeats that of Item \eqref{item:BB_index_coindex_at_p_positive_basic_implies_not_local_min_max} of Theorem \ref{mainthm:BB_decomposition_C*_invariant_complex_analytic_subspace}, although the local geometry is perhaps more easily understood in the context of Corollary \ref{maincor:BB_decomposition_strict_transform_complex_analytic_subspace}.

\begin{rmk}[Resolution of singularities and perfect Morse--Bott functions]
\label{rmk:Resolution_singularities_perfect_Morse-Bott_functions}  
Theorem \ref{mainthm:BB_decomposition_strict_transform_complex_analytic_subspace} asserts that the strict transform $Y'$ is a complex manifold. Suppose now that $X$ is compact, in which case $X'$ is compact and thus $Y'\subset X'$ is an embedded compact, complex submanifold. Assume in addition that $Y'$ is connected, so $Y'$ is a closed, connected, oriented manifold. Consequently (see Hatcher \cite[Theorem 3.26 (a), p. 236]{Hatcher}), one has that $H_d(Y';\ZZ) = \ZZ$ and thus $b_d(Y')=1$, where $d=\dim Y'$, while $H_0(Y';\ZZ) = \ZZ$ and thus $b_0(Y')=1$ (see Hatcher \cite[Proposition 2.7, p. 109]{Hatcher}), where $b_j(Y')$ denotes the $j$-th Betti number of $Y'$, computed as the rank of the singular homology group $H_j(Y';\ZZ)$, for $j = 1,\ldots,d$. 

If $X$ is endowed with K\"ahler metric, then its resolution of singularities, $X'$, also has a K\"ahler metric by repeated application of Proposition \ref{prop:Voisin_3-24}. If the induced circle action on $X'$ has a Hamiltonian function $f:X'\to\RR$ in the sense of the forthcoming \eqref{eq:MomentMap}, then the induced circle action on $Y'$ inherits a Hamiltonian function $f:Y'\to\RR$ by restriction.
Frankel's Theorem \ref{mainthm:Frankel_almost_Hermitian} implies that the Morse--Bott signature of the Hamiltonian $f:Y'\to\RR$ at a critical point $p'\in Y'$ is equal to the Bia{\l}ynicki--Birula signature of the fixed point $p'\in Y'$ given by the $\CC^*$ action on $Y'$. From the discussion in Section \ref{subsec:Morse_inequalities}, we obtain that $f:Y'\to\RR$ is a perfect Morse--Bott function and so there are exactly one connected, critical submanifold of $Y'$ with Morse--Bott index equal to zero, where $f$ attains its \emph{absolute minimum}, and exactly one connected, critical submanifold of $Y'$ with Morse--Bott coindex equal to zero, where $f$ attains its \emph{absolute maximum}.
\end{rmk}

\section[Frankel's Theorem for Hamiltonian circle actions on almost Hermitian manifolds]{Frankel's Theorem for a Hamiltonian function of a circle action on a smooth  almost Hermitian manifold}
\label{sec:Frankel_theorem_circle_actions_almost_Hermitian_manifolds}
The forthcoming version, Theorem \ref{mainthm:Frankel_almost_Hermitian}, of Frankel's Theorem \cite[Section 3]{Frankel_1959} that we proved with Leness in \cite{Feehan_Leness_introduction_virtual_morse_theory_so3_monopoles} is more general than that stated in \cite{Frankel_1959} because we allow for circle actions on closed, smooth manifolds $(M,g,J)$ that are only assumed to be \emph{almost Hermitian}\label{page:Almost_Hermitian_manifold}, so the $g$-orthogonal almost complex structure $J$ need not be integrable and the fundamental two-form
\begin{equation}
  \label{eq:Fundamental_two-form}
  \omega = g(\cdot,J\cdot)
\end{equation}
is non-degenerate but not required to be closed, whereas Frankel assumed in \cite[Section 3]{Frankel_1959} that $\omega$ was closed. (Our convention in \eqref{eq:Fundamental_two-form} agrees with that of Kobayashi \cite[Equation (7.6.8), p. 251]{Kobayashi_differential_geometry_complex_vector_bundles} but is opposite to that used elsewhere, for example, Huybrechts \cite[Definition 1.2.13, p. 29]{Huybrechts_2005}.) Frankel notes \cite[p. 1]{Frankel_1959} that the main results of his article hold when $\omega$ is a $g$-harmonic, symplectic form. (If $\omega$ is symplectic, then $d\omega=0$ while if $g$ is \emph{adapted} to $\omega$, it is well-known that $d^{*_g}\omega=0$ --- see Delano\"e \cite{Delanoe_2002} --- and so $\omega$ is harmonic.) Recall that $J \in C^\infty(\End(TM))$ is an \emph{almost complex structure} on $M$ if $J^2 = -\id_{TM}$ and $J$ is \emph{orthogonal with respect to} or \emph{compatible with} a Riemannian metric $g$ on $M$ if
\begin{equation}
  \label{eq:g_compatible_J}
  g(J\xi,J\eta) = g(\xi,\eta),
\end{equation}
for all vector fields $\xi, \eta \in C^\infty(TM)$. Recall that a smooth manifold is called \emph{almost symplectic} \label{page:almost_symplectic_manifold} if it admits a non-degenerate two-form and \emph{symplectic} \label{page:symplectic_manifold} if that two-form is closed (see Libermann and Marle \cite[Definition 12.4]{Libermann_Marle_symplectic_geometry_analytical_mechanics}).

If a smooth manifold $M$ admits a smooth circle action,
\begin{equation}
\label{eq:Circle_action_smooth_manifold}
\rho:S^1\times M\to M,
\end{equation}
we denote the induced circle action on the tangent bundle $TM$ by
\begin{equation}
\label{eq:Circle_action_tangent_bundle}
\rho_*:S^1\times TM\to TM,
\end{equation}
where, using $D_2\rho$ to denote the differential of $\rho$ in directions tangent to $M$,
\begin{equation}
  \label{eq:rho_*_definition}
  \rho_*(e^{i\theta},v) := D_2\rho(e^{i\theta},x)v, \quad\text{for all } x \in M, v\in T_xM, \text{ and } e^{i\theta}\in S^1.
\end{equation}
A covariant two-tensor field $\varpi \in C^\infty(T^*M\otimes T^*M)$ is called \emph{circle invariant}
(with respect to \eqref{eq:Circle_action_smooth_manifold}) if it obeys
\begin{equation}
\label{eq:Circle_invariant_covariant_2-tensor}
  \varpi\left(\rho_*(e^{i\theta})v, \rho_*(e^{i\theta})w\right) = \varpi(v,w),
  \quad\text{for all } p \in M, v, w \in T_pM, \text{ and } e^{i\theta} \in S^1.
\end{equation}
Similarly, a tensor field $\tau \in C^\infty(TM\otimes T^*M) = C^\infty(\End(TM))$ is called \emph{circle invariant} (with respect to \eqref{eq:Circle_action_smooth_manifold}) if it obeys
\begin{equation}
\label{eq:Circle_invariant_(1,1)-tensor}
  \tau\left(\rho_*(e^{i\theta})v\right) = \rho_*(e^{i\theta})\tau v,
  \quad\text{for all } p \in M, v \in T_pM, \text{ and } e^{i\theta} \in S^1.
\end{equation}
(One could, of course, absorb definitions \eqref{eq:Circle_invariant_covariant_2-tensor} and \eqref{eq:Circle_invariant_(1,1)-tensor} into a general definition of a tensor field invariant under the flow of a vector field, as in Lee \cite[Equation (12.11), p. 323]{Lee_john_smooth_manifolds} and characterize covariant tensor fields that are invariant under the flow of a vector field in terms of vanishing Lie derivatives as in  Lee \cite[Theorem 12.37, p. 324]{Lee_john_smooth_manifolds}.) A smooth two-form $\omega$ or smooth Riemannian metric $g$ on $M$ is called \emph{circle invariant} (with respect to \eqref{eq:Circle_action_smooth_manifold}) if it obeys \eqref{eq:Circle_invariant_covariant_2-tensor} with $\varpi$ replaced by $\omega$ or $g$, respectively. A circle action is called \emph{Hamiltonian} with respect to a smooth two-form $\omega$ on $M$ if \footnote{By analogy with the usual meaning \cite[Definition 2.1]{Dwivedi_Herman_Jeffrey_van_den_Hurk} of a Hamiltonian vector field and Hamiltonian function.} there exists a smooth function $f:M\to\RR$ such that
\begin{equation}
  \label{eq:MomentMap}
  df = \iota_\Theta\omega,
\end{equation}
where $\Theta \in C^\infty(TM)$ is the vector field generated by the circle action, so $\Theta_p = D_1\rho(1,p) \in T_pM$ for all $p\in M$, with $D_1\rho$ denoting the differential of $\rho$ in directions tangent to $S^1$ (see the forthcoming \eqref{eq:Vector_field_generator_circle_action}). Adapting Atiyah and Bott \cite[Section 1, p. 530]{Atiyah_Bott_1983}, Bott \cite[Definition, p. 248]{Bott_1954}, \cite{Bott_1959} (see also Nicolaescu \cite[Definition 2.41]{Nicolaescu_morse_theory}) and specializing the definitions of Feehan \cite[Definition 1.5, p. 86]{Feehan_lojasiewicz_inequality_all_dimensions_morse-bott} and Feehan and Maridakis \cite[Definition 1.10, p. 9]{Feehan_Maridakis_Lojasiewicz-Simon_Banach} from the setting of functions on Banach spaces over a field $\KK$ to Euclidean space, we have the

\begin{maindefn}[Morse--Bott function]
\label{maindefn:Morse-Bott_function}
Let $(M,g)$ be a finite-dimensional, smooth Riemannian manifold and $f:M\to\RR$ be a smooth function. We let
\[
  \Crit f := \{p\in M: df(p)=0\}
\]
denote the \emph{critical set} of $f$. The function $f$ is \emph{Morse--Bott at $p$} if there exists an open neighborhood $U\subset M$ of $p$ such that $U\cap \Crit f$ is a smooth submanifold with tangent space
\[
  T_p\Crit f = \Ker\Hess_gf(p),
\]
where $\Hess_g f \in \End(TM)$ is the Hessian operator (see the forthcoming equation \eqref{eq:DefineHessian}). The function $f$ is \emph{Morse--Bott} if it is Morse--Bott at every point of $\Crit f$. The dimensions of the maximal positive and negative subspaces of $\Hess_gf(p)$ are the \emph{Morse--Bott index} and \emph{co-index} of $f$ at $p$, respectively. 
\end{maindefn}  

\begin{mainthm}[Frankel's theorem for circle actions on almost Hermitian manifolds]
\label{mainthm:Frankel_almost_Hermitian}
(See Feehan and Leness \cite[Theorem 1]{Feehan_Leness_introduction_virtual_morse_theory_so3_monopoles} and compare Frankel \cite[Section 3]{Frankel_1959}.)
Let $M$ be a smooth manifold endowed with a smooth circle action $\rho$ as in \eqref{eq:Circle_action_smooth_manifold} and a non-degenerate two-form $\om$ that is circle invariant in the sense of \eqref{eq:Circle_invariant_covariant_2-tensor}. Then
\begin{enumerate}
\item
\label{item:Frankel_almost_Hermitian_FixedPointsAreZerosOfVField}
A point $p\in M$ is a fixed point of the action \eqref{eq:Circle_action_smooth_manifold} if and only if $\Theta_p=0$;
\item
\label{item:Frankel_almost_Hermitian_ComponentsOfFixedPointsAreSmoothSubmanifolds}
Each connected component of the fixed-point set of the circle action \eqref{eq:Circle_action_smooth_manifold} is a smooth submanifold of even dimension and even codimension in $M$.
\end{enumerate}
In addition, let $f:M\to\RR$ be a smooth function that is Hamiltonian in the sense of \eqref{eq:MomentMap}. Then
\begin{enumerate}
\setcounter{enumi}{2}
\item
\label{item:Frankel_almost_Hermitian_FixedPointsAreCriticalPoints}
A point $p\in M$ is a critical point of $f$ if and only if $p$ is a fixed point of the circle action \eqref{eq:Circle_action_smooth_manifold}.
\end{enumerate}
Furthermore, assume that $g$ is a smooth Riemannian metric on $M$ that is circle invariant in the sense of \eqref{eq:Circle_invariant_covariant_2-tensor}. Then
\begin{enumerate}
\setcounter{enumi}{3}
\item
  \label{item:Frankel_almost_Hermitian_f_is_MB}
  The function $f$ is Morse--Bott at $p$ in the sense of Definition \ref{maindefn:Morse-Bott_function}.
\end{enumerate}
Let $J$ be the almost complex structure on $TM$ determined by $\omega$ and $g$ such that $(\omega,g,J)$ is a compatible triple in the sense that they obey \eqref{eq:Fundamental_two-form} and \eqref{eq:g_compatible_J}. Then
\begin{enumerate}
\setcounter{enumi}{4}
\item
\label{item:Frankel_almost_Hermitian_ACisS1Invariant}
The almost complex structure $J$ is circle invariant in the sense of \eqref{eq:Circle_invariant_(1,1)-tensor};
\item
\label{item:Frankel_almost_Hermitian_WeightsAreEigenvalues}
The eigenvalues of $\Hess_g f \in \End(T_pM)$ are given by the weights of the circle action on $(T_pM,J)$ if the signs of the weights are chosen to be compatible\footnote{See, for example, Feehan and Leness \cite[Lemma 3.3.2 and Definition 3.3.3]{Feehan_Leness_introduction_virtual_morse_theory_so3_monopoles}.} with $J$.
\end{enumerate}
Finally, one has that
\begin{enumerate}
\setcounter{enumi}{6}
\item
\label{item:Frankel_almost_Hermitian_Sylvester}
The signature $(\lambda_p^+(f),\lambda_p^-(f),\lambda_p^0(f))$ of the Hessian operator $\Hess_g f(p)$ is independent of the choice of Riemannian metric $g$, where $\lambda_p^\pm(f)$ denotes the number of positive (negative) eigenvalues of $\Hess_g f(p)$ and $\lambda_p^0(f)$ denotes the nullity of $\Hess_g f(p)$.
\end{enumerate}
\end{mainthm}

\begin{rmk}[Fixed-point sets of isometric circle actions on Riemannian manifolds are totally geodesic submanifolds]
\label{rmk:Fixed-point_sets_isometric_circle_actions_Riemannian_manifolds_totally_geodesic} 
We recall that if $g$ is an circle invariant Riemannian metric on $M$, then the vector field $X \in C^\infty(TM)$ generated by the circle action is a Killing field of $g$. A theorem of Kobayashi \cite[p. 63]{Kobayashi_1958} implies that the zero set of a Killing field is a closed, smooth \emph{totally geodesic} submanifold of even codimension and so Item \eqref{item:Frankel_almost_Hermitian_FixedPointsAreZerosOfVField} implies that the fixed-point subset is a closed, smooth \emph{totally geodesic} submanifold of even codimension, strengthening Item  \eqref{item:Frankel_almost_Hermitian_ComponentsOfFixedPointsAreSmoothSubmanifolds}. See \cite[Remark 2.3.2]{Feehan_Leness_introduction_virtual_morse_theory_so3_monopoles} for an explanation of the significance of the totally geodesic property in the context of Morse--Bott theory.
\end{rmk}

Recall that the \emph{gradient vector field} $\grad_g f \in C^\infty(TM)$ associated to a function $f \in C^\infty(M,\RR)$ is defined by the relation
\begin{equation}
\label{eq:DefineGradient}
  g(\grad_g f,\eta) := df(\eta), \quad\text{for all } \eta \in C^\infty(TM).
\end{equation}  
If $\nabla^g$ is the covariant derivative for the Levi--Civita connection on $TM$ defined by the Riemannian metric $g$, then one can define the Hessian of $f \in C^\infty(M,\RR)$ by (see Petersen \cite[Proposition 2.2.6]{Petersen_2006})
\begin{equation}
\label{eq:DefineHessian}
  \Hess_g f := \nabla^g\grad_g f \in C^\infty(\End(TM)).
\end{equation}
See \cite[Equation (2.1.2)]{Feehan_Leness_introduction_virtual_morse_theory_so3_monopoles} for an alternative expression for $\Hess_g f$. If $p\in M$ is a critical point of $f$, then Theorem \ref{mainthm:Frankel_almost_Hermitian} implies that subspace $T_p^-M \subset T_pM$ on which the Hessian $\Hess_g f(p) \in \End(T_pM)$ is \emph{negative definite} is equal to the subspace of $T_pM$ on which the circle acts with \emph{negative weight}. Hence, the \emph{Morse--Bott index} of $f$ at a critical point $p$,
\[
  \lambda_p^-(f) := \dim_\RR T_p^-M,
\]
is equal to the dimension of the subspace of $T_pM$ on which the circle acts with \emph{negative weight}; one calls $\lambda_p^+(f)$ and $\lambda_p^0(f)$ the \emph{Morse--Bott co-index} and \emph{nullity}, respectively. In  \cite[Section 3]{Frankel_1959}, Frankel says that the critical set of $f$ in Theorem \ref{mainthm:Frankel_almost_Hermitian} is \emph{non-degenerate at $p$} whereas we use the terminology that $f$ is \emph{Morse--Bott at $p$} following \cite[Definition 1.5]{Feehan_lojasiewicz_inequality_all_dimensions_morse-bott}.

\section{Morse--Bott decompositions for complex analytic subspaces of Hermitian manifolds}
\label{sec:Morse-Bott_decompositions_complex_analytic_subspaces_Hermitian_manifolds}
Suppose now that the complex manifold $X$ in Theorem \ref{mainthm:BB_decomposition_C*_invariant_complex_analytic_subspace} is equipped with a smooth Hermitian metric $h$ such that the corresponding Riemannian metric $g = 2\Real h$ is $S^1$-invariant in the sense of \eqref{eq:Circle_invariant_covariant_2-tensor}. Let $\omega = g(\cdot,J\cdot)$ be the corresponding fundamental two-form as in \eqref{eq:Fundamental_two-form}, where $J$ is the standard (and necessarily $S^1$-invariant) almost complex structure on $TX$ and $(g,J,\omega)$ forms a compatible triple. We assume furthermore that the real analytic circle action $S^1\times X \to X$ obtained by restricting the holomorphic $\CC^*$ action to $S^1\subset \CC^*$ is Hamiltonian in the sense of \eqref{eq:MomentMap}, so $\iota_\Theta\omega = df$ on $X$ for a real analytic function $f:X\to\RR$, where the real analytic vector field $\Theta$ on $X$ given by the forthcoming \eqref{eq:Vector_field_generator_circle_action} is the generator of the $S^1$ action. We have the following extension of Definition \ref{maindefn:Virtual_BB_nullity_co-index_index} and analogue of the the classical definitions of Morse--Bott nullity, co-index, and index of a Morse--Bott function at a critical point.

\begin{maindefn}[Virtual Morse--Bott nullity, co-index, and index of a fixed point in a real analytic, almost Hermitian space with a real analytic circle action]
\label{maindefn:Virtual_Morse-Bott_nullity_co-index_index}
Continue the notation of the preceding paragraph, let $(Y,\sO_Y)$ be a $\CC^*$-invariant, closed, complex analytic subspace of $X$, and let $p \in Y^0$, so $p$ is a fixed point of the $S^1$ actions on $X$ and $Y$. Then the \emph{virtual Morse--Bott nullity, co-index}, and \emph{index}, respectively, of the restriction of the Hamiltonian function $f$ to $Y$ at $p$ \emph{relative to a local model space} are defined by
\begin{subequations}
  \label{eq:Virtual_Morse-Bott_nullity_co-index_index}
  \begin{align}
    \label{eq:Virtual_Morse-Bott_nullity}
    \lambda_Y^0(p) := 2\expdim_p Y^0,
    \\
    \label{eq:Virtual_Morse-Bott_co-index}
    \lambda_Y^+(p) := 2\expdim_p Y^+,
    \\
    \label{eq:Virtual_Morse-Bott_index}
    \lambda_Y^-(p) := 2\expdim_p Y^-,                  
  \end{align}
\end{subequations}
where the complex expected dimensions on the right-hand side of \eqref{eq:Virtual_Morse-Bott_nullity_co-index_index} are as in \eqref{eq:Virtual_BB_nullity_co-index_index}.
\end{maindefn}

We recall from equation \eqref{eq:BB_decomposition_C*_invariant_complex_analytic_subspace} in Theorem \ref{mainthm:BB_decomposition_C*_invariant_complex_analytic_subspace} that $Y^{\CC^*} = Y\cap X^{\CC^*}$ and so
\[
  Y^{S^1} = Y\cap X^{S^1}.
\]
We now apply Item \eqref{item:Frankel_almost_Hermitian_FixedPointsAreCriticalPoints} in Theorem \ref{mainthm:Frankel_almost_Hermitian} to interpret the preceding equality of fixed points of a circle action as an equality of critical points of the corresponding Hamiltonian function. Let $T_qY$ denote the Zariski tangent space\footnote{See \eqref{eq:Zariski_cotangent_and_tangent_space} for the definition of the \emph{intrinsic Zariski tangent space} and Lemma \ref{lem:Identification_extrinsic_cotangent_space_with_Zariski_cotangent_space} for the isomorphism with the \emph{extrinsic Zariski tangent space}.} to $Y$ at a point $q\in Y$. Recall from the forthcoming Definition \ref{defn:Critical_point} that we define
\[
  \Crit (f|_Y) := \left\{q\in Y: T_qY \subseteq \Ker \d(q)\right\}.
\]
Suppose $Y$ has embedding dimension, $\embdim_pY = n$, at a point $p$. By the forthcoming Proposition \ref{prop:Equality_dimension_zariski_tangent_space_and_embedding_dimension}, we know that $\embdim_pY = \dim T_pY$, so $\dim T_pY = n$.
One can show that there is an $S^1$-invariant, embedded complex submanifold $U \subset X$ of dimension $n$ such that $(Y\cap U, \sO_Y\restriction U)$ is (isomorphic to) an $S^1$-equivariant local model space for an open neighborhood in $Y$ of $p$ and $T_pU=T_pY$. If $p \in \Crit (f|_Y)$, then $p \in \Crit(f|_U)$ by definition of $\Crit (f|_Y)$ and thus $p$ is a fixed point of the $S^1$ action on $U$ by Item \eqref{item:Frankel_almost_Hermitian_FixedPointsAreCriticalPoints} in Theorem \ref{mainthm:Frankel_almost_Hermitian} and hence a fixed point of the $S^1$ action on $X$, so $p \in \Crit f$. Conversely, if $p \in Y\cap \Crit f$, then $p$ is a fixed point of the $S^1$ action on $X$ by Item \eqref{item:Frankel_almost_Hermitian_FixedPointsAreCriticalPoints} in Theorem \ref{mainthm:Frankel_almost_Hermitian} and so $p$ is a fixed point of the $S^1$ action on $U$ and thus $p \in \Crit(f|_U)$ and so $p \in \Crit(f|_Y)$ by Definition \ref{defn:Critical_point}. In other words, we have shown that
\[
  \Crit(f|_Y) = Y\cap \Crit f.
\]
The preceding identity illustrates that it is natural to interpret the restriction of the Hamiltonian $f$ to $Y$ as a type of Morse--Bott function on $Y$, as implied by Definition \ref{maindefn:Virtual_Morse-Bott_nullity_co-index_index}.

In Section \ref{sec:Morse-Bott_decomposition_smooth_manifold} we briefly recall the Morse--Bott decomposition of a smooth manifold with respect to a Morse--Bott function, analogous to the Bia{\l}ynicki--Birula decomposition of a complex manifold equipped with a holomorphic $\CC^*$ action as in Definition \ref{maindefn:BB_decomposition_complex_manifold}. The equivalence of the two decompositions for a Hermitian manifold with a holomorphic $\CC^*$ action and Hamiltonian induced $S^1$ action is explained in Section \ref{sec:BB_decomposition_complex_manifold}; see also Chriss and Ginzburg \cite[Proposition 2.4.22 and Corollary 2.4.24, p. 89]{Chriss_Ginzburg_representation_theory_complex_geometry}.

\section{Interpretation of the virtual Morse--Bott nullity, co-index, and index}
\label{sec:Interpretation_virtual_Morse-Bott_nullity_coindex_index}
The integers $\lambda_Y^0(p)$, $\lambda_Y^+(p)$, and $\lambda_Y^-(p)$ in Definition \ref{maindefn:Virtual_Morse-Bott_nullity_co-index_index} may be interpreted as dimensions of cohomology groups for elliptic complexes (see Gilkey \cite[Section 1.5, Definition, p. 43]{Gilkey2}) for moduli spaces of solutions to non-linear partial differential equations modulo groups of gauge transformations: see
\begin{inparaenum}
\item Atiyah, Hitchin, and Singer \cite{AHS} or Donaldson and Kronheimer \cite{DK} for the anti-self-dual Yang--Mills equation,
\item Bradlow \cite{Bradlow_1991} and Bradlow and Garcia--Prada \cite{BradlowGP} for the projective vortex equations,
\item Feehan and Leness \cite{FL2b} for the non-Abelian monopole equations over smooth Riemannian four-manifolds,
\item Gromov \cite{Gromov} or McDuff and Salamon \cite{McDuffSalamon2} for the pseudoholomorphic curve equation,
\item Hitchin \cite{Hitchin_1987} for the Higgs pair equations over Riemann surfaces,
\item Kobayashi \cite{Kobayashi} for the Hermitian--Einstein equations over complex K\"ahler manifolds, and \item Simpson \cite{Simpson_1988} for the Higgs pair equations over complex K\"ahler manifolds.
\end{inparaenum}
In these examples, one is typically given a Banach affine space $\cX$ over a field $\KK=\RR$ or $\CC$, a Banach space $\cY$ over $\KK$, analytic actions $\cG\times\cX\to\cX$ and $\cG\times\cY\to\cY$ of a Banach Lie group $\cG$, and a $\cG$-equivariant analytic map $\cF:\cX\to\cY$. One defines a \emph{moduli space of solutions},
\[
  \cM := \cF^{-1}(0)/\cG,
\]
that may have singularities if the action of $\cG$ on $\cF^{-1}(0)$ is not free or if there are points in $\cF^{-1}(0)$ that are not regular points of $\cF$. The \emph{expected dimension} at a point $[p] \in \cM$ represented by $p \in \cF^{-1}(0)$ is defined to be the negative of the Euler characteristic,
\[
  \Euler(d_p) := \dim H_p^0 - \dim H_p^1 + \dim H_p^2
\]
of a \emph{Fredholm complex},
\[
  0 \xrightarrow T_p\cG \xrightarrow{d_p^0} T_p\cX \xrightarrow{d_p^1} \cY
\]
with cohomology groups $H_p^i = \Ker d_p^i / \Ran d_p^{i-1}$ for $i=0,1,2$. In the examples cited, one may compute the Euler characteristic of the elliptic complex from a suitable version of the Atiyah--Singer Index Theorem (see Gilkey \cite{Gilkey2} or H\"ormander \cite{Hormander_v3}). The expected dimension is equal to the dimension of $\cM$ as an analytic manifold at $p$ if the action of $\cG$ is free at $p$ and if $p$ is a regular point of $\cF$. For all $\xi \in T_p\cG$ and $a \in \cX$, the differentials are defined by
\[
  d_p^0\xi := \d R_p(\id_\cG)\xi
  \quad\text{and}\quad
  d_p^1 a := \d\cF(p)a,
\]
where $d_p^0$ is the differential of the analytic, right multiplication map $R_p:\cG \ni u \mapsto u\cdot p \in \cX$. The complex is Fredholm in the sense that the subspaces $\Ran d_p^{i-1} \subset \Ker d_p^i$ are closed, where $d_p^{-1}$ and $d_p^2$ are the zero operators, and the quotients $H_p^i$ are finite-dimensional vector spaces over $\KK$ for $i=0,1,2$.

An application of the Implicit Mapping Theorem for analytic maps (see Donaldson and Kronheimer \cite[Section 4.2.4]{DK} for a description in the case of general Fredholm maps and \cite[Section 4.2.5]{DK} for an application to the anti-self-dual Yang--Mills equation) shows that an open neighborhood of $[p]$ in $\cM$ given by
\[
  F^{-1}(0)/\Stab(p),
\]
defines a germ of an analytic space, where $F: H_p^1 \supset D \to H_p^2$ is an analytic map on an open neighborhood $D$ of the origin in $H_p^1$ and $\Stab(p) := \{u \in \cG: u\cdot p=p\}$ is the isotropy group of $p$ in $\cG$. If the isotropy group is trivial, $\Stab(p) = \{\id_\cG\}$, then the Zariski tangent space to $\cM$ at $[p]$ is given by
\[
  T_{[p]}\cM \cong \Ker\d\cF(p)/\Ran \d R_p(\id_\cG) = \Ker d_p^1/\Ran d_p^0 = H_p^1
\]
and $\dim\Xi_p = \dim H_p^2$. The preceding outline, is described in more detail (for $\KK=\CC$) by Friedman and Morgan \cite[Sections 4.1.3 and 4.1.4]{FrM}, based in turn on an exposition due to Douady \cite{Douady_1966_sem_bourbaki, Douady_1966aifg} of ideas due to Kuranishi \cite{Kuranishi}. 

\section{Outline}
\label{sec:Outline}
While our ultimate goal in this work is to prove new results for complex analytic spaces, we shall often provide expositions of foundational results and intermediate results from the perspectives of both algebraic geometry (specifically, schemes) and analytic spaces (initially for real or complex analytic spaces and later for complex analytic spaces only). We do this because many of the results that we require are often discussed far more clearly and comprehensively in the category of schemes than in the category of analytic spaces. In \cite{Hironaka_1964-I-II}, Hironaka also treats the categories of schemes and (real or complex) analytic spaces more or less separately and, despite the well-known theorems due to Serre \cite{Serre_1956} that address the close relationship between algebraic geometry over $\CC$ and complex analytic spaces, we believe that it is mathematically preferable to provide parallel but separate treatments of key concepts rather than attempt any (necessarily cumbersome) synthesis. 

We begin in Chapter \ref{chap:Analytic_spaces} with a brief introduction to analytic spaces as used in our work. Chapter \ref{chap:Smooth_point_regular_point_dimension_scheme} focuses on algebraic geometry and provides definitions and results for smooth points, regular points, and dimension theory for schemes. In Chapter \ref{chap:Smooth_point_regular_point_dimension_analytic_space}, we discuss the concepts of smooth points, regular points, and dimension theory in the category of analytic spaces, by analogy as far as possible with our discussion of those topics in the category of schemes. Chapter \ref{chap:Blowups_properties_schemes} provides an introduction to blowups and a review of their properties for schemes. In Chapter \ref{chap:Blowups_properties_analytic_spaces}, we describe blowups and their properties in the category of analytic spaces, again preserving the analogy with our description in the category of schemes whenever possible. The preceding chapters provide the background required in Chapter \ref{chap:Resolution_singularities} to state the main results on monomialization of ideal sheaves and resolution of singularities in the categories of schemes and analytic spaces. In Chapter \ref{chap:Stable_manifold_theorems_Morse_and_Morse-Bott_functions}, we discuss definitions of critical sets, stable manifolds, and unstable manifolds for Morse--Bott functions, the stable manifold theorem for Morse--Bott functions, and the Morse--Bott decomposition of a smooth manifold. Chapter \ref{chap:BB_decomposition_algebraic_variety_scheme_complex_manifold} provides an introduction to Bia{\l}ynicki--Birula decompositions for smooth algebraic varieties and complex manifolds. In Chapters \ref{chap:BB_decomposition_blowup_complex_manifold_algebraic} and \ref{chap:BB_decomposition_blowup_complex_manifold_analytic}, respectively, we develop our results on Bia{\l}ynicki--Birula decompositions for blowups of complex manifolds and complex analytic spaces using two different approaches --- algebraic and analytic, respectively. Appendix \ref{chap:Technical_results_definitions} contains some technical results and definitions that are occasionally used in the body of our work, but which are of secondary importance.









\chapter{Analytic spaces}
\label{chap:Analytic_spaces}
While the literature for algebraic varieties and schemes is well established, for analytic spaces the available references are less comprehensive, at least for the aspects of the theory that we shall need, and less so again for real analytic spaces. Moreover, terminology can vary among the references that do exist. For that reason, we provide in this chapter a concise summary of the basic definitions and results that we shall need in our work and a reader's guide to the literature. In Section \ref{sec:Analytic_varieties_and_spaces}, we review the key definitions underlying the concepts of analytic spaces, analytic varieties, and analytic sets. Section \ref{sec:Explicit_description_morphisms_analytic_varieties} provides an explicit description of a morphism of analytic model spaces, by analogy with the explicit description of a morphism of affine algebraic varieties. In Section \ref{sec:Complexification_real_analytic_spaces}, we briefly discuss complexification of real analytic germs and real analytic spaces. Section \ref{sec:Real_part_complex_analytic_spaces} outlines the concept of the real part of a complex analytic space. 

\section{Analytic spaces, varieties, and sets}
\label{sec:Analytic_varieties_and_spaces}
One can define analytic spaces over any field that is complete with respect to a valuation (see Onishchik \cite{Onishchik}). Write $\RR_+ = [0,\infty)$ and recall that a \emph{valuation} (or \emph{norm}) on a field $\KK$ is a function $|\cdot|:\KK \to \RR_+$ with the following properties (see Berkovich \cite[Section 1.1.1]{Berkovich_ictp_lectures}) for all $a,b\in\KK$:
\begin{enumerate}
\item (positivity) $|a|= 0$ if and only if $a=0$;
\item (multiplicativity) $|ab|=|a|\cdot|b|$;
\item (triangle axiom) $|a+b| \leq |a|+|b|$.
\end{enumerate}
Any valuation $|\cdot|:\KK \to \RR_+$ defines a metric on $\KK$ with respect to which the distance between two elements $a,b \in \KK$ equals $|a-b|$. The \emph{completion} of $\KK$ with respect to this metric is a field $\widehat\KK$ which contains $\KK$ and is provided with a valuation $|\cdot|:\widehat\KK \to \RR_+$ that extends that on $\KK$. A field $\KK$ with a valuation is \emph{complete} if it is complete as a metric space.

The \emph{Archimedean axiom} says that, for any non-zero $x \in \KK$, there is a positive integer $n$ such that $|nx| > 1$. An \emph{Archimedean field} is one in which this axiom holds, such as the real numbers, $\RR$, and the complex numbers, $\CC$ (see Payne \cite[Section 1.3]{Payne_2015}). The only complete Archimedean fields are $\RR$ or $\CC$, with their usual valuation $|\cdot|$ or a power $|\cdot|^\eps$, for $0<\eps<1$ (see \cite[Section 1.1.4]{Berkovich_ictp_lectures}, \cite[Section 1.3]{Payne_2015}).

A \emph{non-Archimedean field} is a complete normed field for which the Archimedean axiom does not hold (see \cite[Section 1.1.4]{Berkovich_ictp_lectures}, \cite[Section 1.4]{Payne_2015}). The theory of analytic spaces over non-Archimedean fields (also known as \emph{Berkovich spaces} --- see Ducros, Favre, and Nicaise \cite{Ducros_Favre_Nicaise_berkovich_spaces_and_applications}) has been developed by Berkovich  \cite{Berkovich_spectral_theory_analytic_geometry_non-archimedean_fields, Berkovich_1993, Berkovich_integration_one-forms_p-adic_analytic_spaces, Berkovich_ictp_lectures}. Bosch, G\"untzer, and Remmert \cite{Bosch_Guntzer_Remmert_non-archimedean_analysis} provide an extensive reference for analysis over non-Archimedean fields. Lecture notes by Jonsson \cite{Jonsson_topics_algebraic_geometry_berkovich_spaces, Jonsson_lecture_notes_on_berkovich_spaces} and articles by Conrad and Temkin \cite{Conrad_2008, Conrad_Temkin_2009}, Ducros \cite{Ducros_2015}, Payne \cite{Payne_2015}, and Temkin \cite{Temkin_2015} provide further background, development, applications, and motivations for Berkovich spaces. See also Grauert and Remmert \cite{Grauert_Remmert_analytic_local_algebras}. 

In our discussion of analytic spaces in this monograph, we shall always assume that $\KK$ is $\RR$ or $\CC$, with their standard valuations, and we refer the reader to the preceding references for insights into possible extensions of our results from the category of analytic spaces over $\KK=\RR$ or $\CC$ to one where $\KK$ is a non-Archimedean field. 

\begin{defn}[Analytic model space]
\label{defn:Analytic_model_space}
(See Acquistapace, Broglia, and Fernando \cite[Section 1.1, ``Way 2'']{Acquistapace_Broglia_Fernando_topics_global_real_analytic_geometry} for $\KK=\CC$ and \cite[Definition 2.1, second item]{Acquistapace_Broglia_Fernando_topics_global_real_analytic_geometry} for $\KK=\RR$, Fischer \cite[Section 0.14]{Fischer_complex_analytic_geometry} for $\KK=\CC$, Grauert and Remmert \cite[Section 1.1.2]{Grauert_Remmert_coherent_analytic_sheaves} for $\KK=\CC$, Guaraldo, Macr\`\i, and Tancredi \cite[Section 2.1, Definition 1.1, p. 11]{Guaraldo_Macri_Tancredi_topics_real_analytic_spaces} for $\KK=\RR$ or $\CC$, and Hironaka \cite[Chapter 0, Section 2, p. 119]{Hironaka_1964-I-II} for $\KK=\RR$ or $\CC$.)  
Let $D \subset \KK^n$ be a domain, $\sO_D$ be the sheaf of germs of analytic functions on $D$, and $\sI$ be an ideal sheaf in $\sO_D$ of \emph{finite type}, so for every point $x \in D$ there are an open neighborhood $U \subset D$ of $x$ and analytic functions $f_1,\ldots,f_k \in \sO(U)$, such that the sheaf $\sI$ is generated over $U$ by $f_1,\ldots,f_k$, that is,
\[
  \sI_U = \sO_U f_1 + \cdots + \sO_U f_k.
\]
The quotient sheaf $\sO_D/\sI$ is a sheaf of rings on $D$ and has topological support $Y = \supp(\sO_D/\sI)$, namely the set of all points $x\in D$ where $(\sO_D/\sI)_x \neq 0$, or equivalently where $\sI_x \neq \sO_{D,x}$. In a neighborhood $U$ of $x$, one has
\[
  Y \cap U = \{y\in U: f_1(y) = \cdots = f_k(y) = 0\}.
\] 
The restriction
\[
  \sO_Y := (\sO_D/\sI)\restriction Y
\]
of $\sO_D/\sI$ to $Y$ is a sheaf of rings on $Y$ and $(Y, \sO_Y)$ is a $\KK$-ringed space. One calls $(Y, \sO_Y)$ a \emph{analytic model space over $\KK$} (or \emph{$\KK$-analytic model space}) (in $D$) and $(Y, \sO_Y)$ is the \emph{analytic model space} (or \emph{$\KK$-analytic model space}) defined by an ideal $\sI \subset \sO_D$ of finite type.
\end{defn}

See B\u{a}nic\u{a} and St\u{a}n\u{a}\c{s}il\u{a} \cite{Banica_Stanasila_algebraic_methods_global_theory_complex_spaces} for an approach to complex analytic spaces that emphasizes links with algebraic geometry. In order to help readers understand differences between and properties in common for real and complex analytic sets or varieties, we recall the

\begin{defn}[Coherent sheaf]
\label{defn:Coherent_sheaf}
(See Acquistapace, Broglia, and Fernando \cite[Definition 1.1]{Acquistapace_Broglia_Fernando_topics_global_real_analytic_geometry}, Fischer \cite[Section 0.5]{Fischer_complex_analytic_geometry}, Grauert and Remmert \cite[Annex, Section 3.3, p. 235]{Grauert_Remmert_coherent_analytic_sheaves}, Guaraldo, Macr\`\i, and Tancredi \cite[Section 1.2, Definition 2.1, p. 4]{Guaraldo_Macri_Tancredi_topics_real_analytic_spaces}, Noguchi \cite[Definition 2.4.4]{Noguchi_analytic_function_theory_several_variables}, or the Stacks Project \cite[\href{https://stacks.math.columbia.edu/tag/01BV}{Definition 01BV}]{stacks-project}.)  
A sheaf $\sF$ on a ringed space $(X,\sO_X)$ is called $\sO_X$-\emph{coherent} if it is a sheaf of $\sO_X$-modules of \emph{finite presentation}, that is, it obeys the following two conditions:
\begin{enumerate}
\item\label{item:Finite_type} The sheaf is of \emph{finite type}, that is, for each $x \in X$ there exists an open neighborhood
$U$ of $x$ and finitely many sections $s_1,\ldots,s_k$ on $U$ generating the fiber $\sF_y$ for each $y \in U$.
\item\label{item:Relation_finite_type} The sheaf is of \emph{relation finite type}, that is, for each open set $U \subset X$ and for each finite number of sections $t_1,\ldots,t_l \in \sF(U)$, the sheaf of relations among them, that is, the kernel of the following homomorphism of sheaves is a sheaf of finite type:
\[
  \sigma:\sO_X^l \restriction U \ni (a_1,\ldots,a_l) \mapsto a_1t_1 + \cdots + a_lt_l \in \sF(U).  
\]    
\end{enumerate}
\end{defn}

Next, we recall the

\begin{thm}[Oka's Coherence Theorem]
\label{thm:Oka_coherence_theorem}
(See Abhyankar\footnote{In his monograph, Abhyankar allows $\KK$ to be any complete valued field, though he occasionally also requires $\KK$ to be algebraically closed \cite[pp. vii--viii]{Abhyankar_local_analytic_geometry}.} \cite[Theorem 15.2]{Abhyankar_local_analytic_geometry} for $\KK=\CC$ or $\RR$, Grauert and Remmert \cite[Section 2.5.2, Theorem, p. 58]{Grauert_Remmert_coherent_analytic_sheaves} for $\KK=\CC$, H\"ormander \cite[Theorem 6.4.1]{Hormander_introduction_complex_analysis_several_variables} for $\KK=\CC$, Narasimhan \cite[Section 4.2, Theorem 3, p. 73]{Narasimhan_introduction_theory_analytic_spaces} for $\KK=\CC$ or $\RR$, or Noguchi \cite[Theorem 2.5.1]{Noguchi_analytic_function_theory_several_variables} for $\KK=\CC$.)
Let $\KK=\RR$ or $\CC$ and $n$ be a positive integer. If $\Omega$ is an open set in $\KK^n$ and $\sO_\Omega$ denotes the sheaf of germs of $\KK$-analytic functions on $\Omega$, then $\sO_\Omega$ is a coherent sheaf of rings.  
\end{thm}

\begin{rmk}[Significance of Oka's Coherence Theorems]
\label{rmk:Significance_Oka_coherence_theorems}
For a detailed discussion of the history and significance of Oka's Coherence Theorems, we refer the reader to Noguchi \cite[Chapter 2, Historical Supplements, p. 61, Chapter 9, and Appendix Kiyoshi Oka]{Noguchi_analytic_function_theory_several_variables}.
\end{rmk}

\begin{rmk}[Coherence of analytic model spaces]
\label{rmk:Analytic_model_space_coherence}
For verifications that $\sO_Y$ in the Definition \ref{defn:Analytic_model_space} of an analytic model space $(Y,\sO_Y)$ is a coherent sheaf of $\sO_D$-modules, we refer to Acquistapace, Broglia, and Fernando \cite[Section 1.1, pp. 7--8]{Acquistapace_Broglia_Fernando_topics_global_real_analytic_geometry} for $\KK=\CC$ and \cite[Remark 2.2, pp. 23]{Acquistapace_Broglia_Fernando_topics_global_real_analytic_geometry} for $\KK=\RR$, Grauert and Remmert \cite[Section 2.5.3, Proposition, p. 60]{Grauert_Remmert_coherent_analytic_sheaves} for $\KK=\CC$, and Guaraldo, Macr\`\i, and Tancredi \cite[Section 2.1, Definition 1.1, p. 11]{Guaraldo_Macri_Tancredi_topics_real_analytic_spaces} for $\KK=\RR$ or $\CC$.

According to Theorem \ref{thm:Oka_coherence_theorem}, the sheaf $\sO_D$ is a coherent sheaf of rings and, in particular, $\sO_D$ is a coherent sheaf of $\sO_D$-modules (see Grauert and Remmert \cite[Annex, Section 3.2, Example 1, p. 235]{Grauert_Remmert_coherent_analytic_sheaves}). Since $\sI$ in Definition \ref{defn:Analytic_model_space} is a subsheaf (of ideals) of the sheaf $\sO_D$ of relation finite type, then $\sI$ is of relation finite type as well (see Grauert and Remmert \cite[Annex, Section 3.2, Example 2, p. 235]{Grauert_Remmert_coherent_analytic_sheaves}). Since $\sI$ is a finitely generated $\sO_D$-module, it is of finite type and thus a coherent sheaf of $\sO_D$-modules. Consequently, $\sO_Y = (\sO_D/\sI)\restriction Y$ in Definition \ref{defn:Analytic_model_space} is a coherent sheaf of $\sO_D$-modules as a consequence of Grauert and Remmert \cite[Annex, Section 4.1, Consequence 6, p. 238]{Grauert_Remmert_coherent_analytic_sheaves}, together with Fischer \cite[Section 0.11, Lemma, p. 8]{Fischer_complex_analytic_geometry} or Grauert and Remmert \cite[Annex, Section 4.3, Remark and Proposition, p. 239]{Grauert_Remmert_coherent_analytic_sheaves} (the extension principle for coherent sheaves). 
\end{rmk}

For the extension principle cited in Remark \ref{rmk:Analytic_model_space_coherence}, one is given a closed subset $Y$ of a topological space $X$ and $\iota:Y\to X$ denotes the injection. Given a sheaf $\sO_Y$ of rings on $Y$ and an $\sO_Y$-module $\sF$ on $Y$, the image sheaf $\iota_*\sF$ is an $\iota_*\sO_Y$-module on $X$. The extension principle states that $\sF$ is $\sO_Y$-coherent if and only if $\iota_*\sF$ is $\iota_*\sO_Y$-coherent. In Definition \ref{defn:Analytic_model_space}, we have $\sO_Y = (\sO_{\KK^n}/\sI)\restriction Y$ and thus $\iota_*\sO_Y = \sO_{\KK^n}/\sI$. Since $\sO_{\KK^n}/\sI$ is $(\sO_{\KK^n}/\sI)$-coherent (that is, a coherent sheaf of rings), then $\sO_Y$ is $\sO_Y$-coherent (that is, a coherent sheaf of rings). See also Grauert and Remmert \cite[Section 1.2.4, Extension Principle, p. 17]{Grauert_Remmert_coherent_analytic_sheaves}.

The model space $(Y,\sO_Y)$ in Definition \ref{defn:Analytic_model_space} is a \emph{ringed space}: see Fischer \cite[Section 0.1, p. 1]{Fischer_complex_analytic_geometry}, G\"ortz and Wedhorn \cite[Definition 2.29]{Gortz_Wedhorn_algebraic_geometry_v1}, Grauert and Remmert \cite[Section 1.1.0, p. 1]{Grauert_Remmert_coherent_analytic_sheaves}, Hartshorne \cite[Chapter II, Section 2, Definition, p. 72]{Hartshorne_algebraic_geometry}, Shafarevich \cite[Section 5.3.1, Definition 5.1, p. 24]{Shafarevich_v2}, the Stacks Project  \cite[\href{https://stacks.math.columbia.edu/tag/0091}{Definition 0091}]{stacks-project}, and Vakil \cite[Example 2.2.13]{Vakil_foundations_algebraic_geometry}.

More specifically, the model space $(Y,\sO_Y)$ in Definition \ref{defn:Analytic_model_space} is a \emph{locally ringed space} since, for each $x\in Y$, the stalk $\sO_{Y,x}$ is a local ring, with unique maximal ideal $\fm_{Y,x}$: see Fischer \cite[Section 0.1, p. 1, and Section 0.14, p. 9]{Fischer_complex_analytic_geometry}, G\"ortz and Wedhorn \cite[Definition 2.30]{Gortz_Wedhorn_algebraic_geometry_v1}, Grauert and Remmert \cite[Section 1.1.3, p. 5]{Grauert_Remmert_coherent_analytic_sheaves}, Hartshorne \cite[Chapter II, Section 2, Definition, p. 72]{Hartshorne_algebraic_geometry}, and Vakil \cite[Example 2.2.13]{Vakil_foundations_algebraic_geometry}.

More specifically still, $(Y,\sO_Y)$ is a $\KK$-\emph{ringed space} since $\sO_Y$ is a \emph{sheaf of local $\KK$-algebras}: see Fischer \cite[Section 0.1, p. 2, and Section 0.14, p. 9]{Fischer_complex_analytic_geometry} and Grauert and Remmert \cite[Section 1.1.3, p. 5]{Grauert_Remmert_coherent_analytic_sheaves} for $\KK=\CC$ and Guaraldo, Macr\`\i, and Tancredi \cite[Section 1.1, Definition 1.1, p. 1 and Section 2.1, Definition 1.4, p. 13]{Guaraldo_Macri_Tancredi_topics_real_analytic_spaces} for $\KK=\CC$ or $\RR$.


\begin{rmk}[Verification that $\KK$-analytic model spaces are $\KK$-ringed spaces]
\label{rmk:Analytic_model_space_K-ringed_space}  
To explain why the model space $(Y,\sO_Y)$ in Definition \ref{defn:Analytic_model_space} is a $\KK$-ringed space and not merely a ringed space, we review the relevant commutative algebra. A ring $A$ is \emph{local} if it has a unique maximal ideal; if $\fm$ is the maximal ideal of a local ring $A$, then $\kappa := A/\fm$ is the \emph{residue field} (see Atiyah and Macdonald \cite[Chapter 1, p. 4]{Atiyah_Macdonald_introduction_commutative_algebra} or Matsumura \cite[Section 1.1, p. 3]{Matsumura_commutative_ring_theory}). If $\KK$ is any field, then the ring $\KK[[x_1,\ldots,x_n]]$ of \emph{formal power series} (see Zariski and Samuel \cite[Section 7.1, p. 129]{Zariski_Samuel_communtative_algebra_II}) is a local ring with maximal ideal $\fm_0$ generated by the indeterminates $x_1,\ldots,x_n$ (see Zariski and Samuel \cite[Section 7.1, Corollary 1, p. 131]{Zariski_Samuel_communtative_algebra_II} or Nagata\footnote{Noting that any field is a local ring.} \cite[Chapter 2, Section 15, Corollary 15.4, p. 50]{Nagata_local_rings}). More generally, if $p = (p_1,\ldots,p_n) \in \KK^n$ is any point, then the kernel $\fm_{\KK^n,p}$ of the evaluation homomorphism of rings,
\[
  \ev_p:\KK[[x_1-p_1,\ldots,x_n-p_n]] \ni f \mapsto f(p_1,\ldots,p_n) \in \KK
\]
is an ideal and is maximal since the quotient is a field (see Atiyah and Macdonald \cite[Chapter 1, p. 3]{Atiyah_Macdonald_introduction_commutative_algebra}):
\[
  \KK[[x_1-p_1,\ldots,x_n-p_n]]/\fm_{\KK^n,p} \cong \KK.
\]
When $\KK=\RR$ or $\CC$, we may replace the role of the ring of formal power series in the preceding discussion by that of the ring $\KK\{x_1-p_1,\ldots,x_n-p_n\}$ of \emph{convergent power series} to show that it is also a local ring with maximal ideal $\fm_{\KK^n,p}$ (given by the kernel of the evaluation homomorphism). In particular, the stalk $\sO_{\KK^n,p}$ at a point $p \in \KK^n$ of the sheaf $\sO_{\KK^n}$ of $\KK$-analytic functions on $\KK^n$ is a local ring with maximal ideal $\fm_{\KK^n,p}$ and $\sO_{\KK^n,p}/\fm_{\KK^n,p} \cong \KK$.

For any ideal $\fb \subset \sO_{\KK^n,p}$, there is a bijection between the ideals of $\sO_{\KK^n,p}$ that contain $\fb$ and the ideals of the quotient ring $\sO_{\KK^n,p}/\fb$.
Hence, $\fm_{\KK^n,p}/\fb$ is the unique maximal ideal of $\sO_{\KK^n,p}/\fb$ since $\fm_{\KK^n,p}$ is the unique maximal ideal of $\sO_{\KK^n,p}$ and thus $\sO_{\KK^n,p}/\fb$ is a local ring. In particular, since $\sO_Y = \sO_D/\sI\restriction Y$ and $\sO_{D,p} = \sO_{\KK^n,p}$, we see that $\sO_{Y,p} = \sO_{D,p}/\sI_p$ is a local ring for any point $p \in Y \subset D$. Let $\fm_{Y,p} \subset \sO_{Y,p}$ denote the maximal ideal and note that we must have $\fm_{Y,p} = \fm_{\KK^n,p}/\sI_p$, so that the residue field of $\sO_{Y,p}$ is given by
\[
  \sO_{Y,p}/\fm_{Y,p}
  =
  \left(\sO_{\KK^n,p}/\sI_p\right)/\left(\fm_{\KK^n,p}/\sI_p\right)
  \cong
  \sO_{\KK^n,p}/\fm_{\KK^n,p} \cong \KK.
\]
Hence, $(Y,\sO_Y)$ is a $\KK$-ringed space, as claimed. \qed
\end{rmk}

\begin{defn}[Morphisms of ringed spaces]
\label{defn:Morphisms_ringed_spaces}
(For the definition of a \emph{morphism of ringed spaces}, see Fischer \cite[Section 0.3, p. 2]{Fischer_complex_analytic_geometry}, G\"ortz and Wedhorn \cite[Definition 2.29]{Gortz_Wedhorn_algebraic_geometry_v1}, Grauert and Remmert \cite[Section 1.1.4, p. 6]{Grauert_Remmert_coherent_analytic_sheaves}, Shafarevich \cite[Section 5.3.1, Definition 5.2, p. 24]{Shafarevich_v2}, and the Stacks Project  \cite[\href{https://stacks.math.columbia.edu/tag/0091}{Definition 0091}]{stacks-project}; for the definition of a \emph{morphism of $\CC$-ringed spaces}, see Fischer \cite[Section 0.3, p. 2]{Fischer_complex_analytic_geometry}, Grauert and Remmert \cite[Section 1.1.4, p. 6]{Grauert_Remmert_coherent_analytic_sheaves}, and Guaraldo, Macr\`\i, and Tancredi \cite[Section 1.1, Definition 1.5, p. 3]{Guaraldo_Macri_Tancredi_topics_real_analytic_spaces}.)
A \emph{morphism of ringed spaces}, $(X,\sO_X) \to (Y,\sO_Y)$, is a pair $(\varphi,\varphi^\sharp)$ comprising a continuous map $\varphi:X\to Y$ of topological spaces and a homomorphism $\varphi^\sharp:\sO_Y \to \varphi_*\sO_X$ of sheaves of rings on $Y$. For $p\in X$, the ring homomorphism
\begin{equation}
  \label{eq:Ring_homomorphism_stalks}
  \varphi_p^\sharp:\sO_{Y,\varphi(p)} \to \sO_{X,p}
\end{equation}
is defined to be the composition of the canonical homomorphisms.
\[
  \sO_{Y,\varphi(p)} \to (\varphi_*\sO_X)_{\varphi(p)} \to \sO_{X,p}.
\]  
A \emph{morphism of locally ringed spaces}, $(X,\sO_X) \to (Y,\sO_Y)$, is a morphism of ringed spaces such that the ring homomorphism \eqref{eq:Ring_homomorphism_stalks} is \emph{local} for every $p\in X$, that is,
\[
  \varphi_p^\sharp\left(\fm_{Y,\varphi(p)}\right) \subset \fm_{X,p}.
\]
A \emph{morphism of $\KK$-ringed spaces}, $(X,\sO_X) \to (Y,\sO_Y)$, is a morphism $(\varphi,\varphi^\sharp)$ of ringed spaces
where $\varphi^\sharp$ is furthermore a homomorphism of sheaves of $\KK$-algebras. 
\end{defn}

Morphisms of ringed spaces may be composed in a canonical way (see Fischer \cite[Section 0.3, p. 2]{Fischer_complex_analytic_geometry}), yielding categories of \emph{ringed spaces}, \emph{locally ringed spaces}, and \emph{$\KK$-ringed spaces} and corresponding definitions of \emph{epimorphisms}, \emph{monomorphisms}, and \emph{isomorphisms}. One has the

\begin{lem}[Monomorphisms, epimorphisms, and isomorphisms of $\KK$-ringed spaces]
\label{lem:Fischer_0-3}  
(See Fischer \cite[Section 0.4, p. 3]{Fischer_complex_analytic_geometry}.)  
Let $(\varphi,\varphi^\sharp):(X,\sO_X) \to (Y,\sO_Y)$ be a morphism of $\KK$-ringed spaces.
\begin{enumerate}
\item If $\varphi$ is injective and $\varphi_p^\sharp$ is surjective for every $p \in X$, then $(\varphi,\varphi^\sharp)$ is a monomorphism.
\item If $\varphi$ is surjective and $\varphi_p^\sharp$ is injective for every $p \in X$, then $(\varphi,\varphi^\sharp)$ is an epimorphism.  
\item $(\varphi,\varphi^\sharp)$ is an isomorphism of $\KK$-ringed spaces if and only if $\varphi$ is a homeomorphism and $\varphi_p^\sharp$ is an isomorphism for every $p \in X$.  
\end{enumerate}
\end{lem}    

\begin{defn}[Analytic space]
\label{defn:Analytic_space}
(See Abhyankar \cite[Chapter 7, Section 43.2]{Abhyankar_local_analytic_geometry}, Acquistapace, Broglia, and Fernando \cite[Definition 1.2]{Acquistapace_Broglia_Fernando_topics_global_real_analytic_geometry} for $\KK=\CC$ and \cite[Definition 2.1]{Acquistapace_Broglia_Fernando_topics_global_real_analytic_geometry} for $\KK=\RR$, Fischer \cite[Section 0.14]{Fischer_complex_analytic_geometry}, Grauert and Remmert \cite[Section 1.1.5]{Grauert_Remmert_coherent_analytic_sheaves} for $\KK=\CC$, Guaraldo, Macr\`\i, and Tancredi \cite[Section 2.1, Definition 1.4, p. 13]{Guaraldo_Macri_Tancredi_topics_real_analytic_spaces} for $\KK=\RR$ or $\CC$, Hironaka \cite[Chapter 0, Section 2, pp. 119--120]{Hironaka_1964-I-II} for $\KK=\RR$ or $\CC$, and Narasimhan \cite[Section 4.1, p. 64]{Narasimhan_introduction_theory_analytic_spaces} for $\KK=\CC$.)
Let $\KK=\RR$ or $\CC$. Let $(X, \sO_X)$ be a $\KK$-ringed space such that $X$ is a Hausdorff and paracompact\footnote{Paracompactness is required by Acquistapace, Broglia, and Fernando but not by Grauert and Remmert or Fischer.}. One calls $(X, \sO_X)$ an \emph{analytic space over $\KK$} (or \emph{$\KK$-analytic space} or \emph{complex} (respectively, \emph{real}) \emph{analytic space}) if every point of $X$ has an open neighborhood $U$ such that the open $\KK$-ringed subspace $(U, \sO_U)$ of $(X, \sO_X)$ is isomorphic (as a $\KK$-ringed space) to an analytic model space over $\KK$ as in Definition \ref{defn:Analytic_model_space}. (Some authors write $|X|$ or $\supp\sO_X$ to denote the underlying topological space.) 
\end{defn}

When $\KK=\RR$ in Definition \ref{defn:Analytic_space}, Acquistapace, Broglia, and Fernando \cite{Acquistapace_Broglia_Fernando_topics_global_real_analytic_geometry} call $(X, \sO_X)$ a \emph{real $C$-analytic space} while Guaraldo, Macr\`\i, and Tancredi \cite{Guaraldo_Macri_Tancredi_topics_real_analytic_spaces} call it a \emph{real analytic space}, as we do in this monograph.

\begin{rmk}[Coherence of analytic spaces]
\label{rmk:Analytic_space_coherence}  
The structure sheaf $\sO_X$ of every complex analytic space $X$ is coherent, as a consequence of that property for complex analytic model spaces noted in Remark \ref{rmk:Analytic_model_space_coherence} for $\KK=\CC$ or $\RR$, or directly by Fischer \cite[Section 0.14, Theorem, p. 11]{Fischer_complex_analytic_geometry} or Grauert and Remmert \cite[Section 2.5.3, Proposition, p. 60]{Grauert_Remmert_coherent_analytic_sheaves}. Narasimhan \cite[Section 4.3, Theorem 6, p. 80]{Narasimhan_introduction_theory_analytic_spaces} also proves coherence, for $\KK=\CC$, for what he refers to as a ``complex space'', but his definition \cite[Section 4.1, p. 64]{Narasimhan_introduction_theory_analytic_spaces} corresponds to that of a \emph{reduced complex analytic space} in the sense of Definition \ref{defn:Analytic_variety}.
\end{rmk}

\begin{defn}[Morphisms of analytic spaces]
\label{defn:Morphism_analytic_spaces}
(See Fischer \cite[Section 0.14, p. 9]{Fischer_complex_analytic_geometry} or Grauert and Remmert \cite[Section 1.1.4, p. 6]{Grauert_Remmert_coherent_analytic_sheaves} for the case $\KK=\CC$.)  
A \emph{morphism of $\KK$-analytic spaces}, $(\varphi,\varphi^\sharp):(X,\sO_X) \to (Y,\sO_Y)$, is a morphism of $\KK$-ringed spaces.
\end{defn}

One thus obtains a category of \emph{$\KK$-analytic spaces} (see Grauert and Remmert \cite[Section 1.1.4, p. 7]{Grauert_Remmert_coherent_analytic_sheaves} for the case $\KK=\CC$).

A fundamental difference between the cases $\KK=\RR$ and $\CC$ is explained by Acquistapace, Broglia, and Fernando \cite[Section 1.2]{Acquistapace_Broglia_Fernando_topics_global_real_analytic_geometry} and Guaraldo, Macr\`\i, and Tancredi \cite[Chapter II, Remark 1.3]{Guaraldo_Macri_Tancredi_topics_real_analytic_spaces}. To explain this difference, we recall the

\begin{defn}[Analytic set]
\label{defn:Analytic_set}
(See Abhyankar \cite[Chapter 1, Section 2, p. 8]{Abhyankar_local_analytic_geometry} for $\KK=\RR$ or $\CC$, Acquistapace, Broglia, and Fernando \cite[Section 1.1]{Acquistapace_Broglia_Fernando_topics_global_real_analytic_geometry} for $\KK=\CC$ and \cite[Definition 2.1]{Acquistapace_Broglia_Fernando_topics_global_real_analytic_geometry} for $\KK=\RR$, Demailly \cite[Definition 4.1, p. 91]{Demailly_complex_analytic_differential_geometry} for $\KK=\CC$, Fischer \cite[Section 0.14, p. 10]{Fischer_complex_analytic_geometry}, Grauert and Remmert \cite[Section 4.1.1, p. 76]{Grauert_Remmert_coherent_analytic_sheaves} for $\KK=\CC$, Guaraldo, Macr\`\i, and Tancredi \cite[Section 2.1, Definition 1.2, p. 11]{Guaraldo_Macri_Tancredi_topics_real_analytic_spaces} for $\KK=\RR$ or $\CC$, H\"ormander \cite[Definition 6.5.1, p. 167]{Hormander_introduction_complex_analysis_several_variables} for $\KK=\CC$, and Narasimhan \cite[Chapter 1, Definition 1, p. 5]{Narasimhan_introduction_theory_analytic_spaces} for $\KK=\RR$ or $\CC$.)
Let $\KK=\RR$ or $\CC$ and $(X,\sO_X)$ be an analytic space over $\KK$. A subset $S \subset X$ is called \emph{analytic at a point $p\in X$} if there exist an open neighborhood $U\subset X$ of $p$ and finitely many analytic functions $f_1, \ldots, f_k \in \sO_X(U)$ such that
\[
  S \cap U = \{y\in U: f_1(y) = \cdots = f_k(y) = 0\},
\]
where the number $k$ may depend on $U$. The set $S$ is called \emph{analytic in $X$} if $S$ is analytic at every point of $X$. Guaraldo, Macr\`\i, and Tancredi replace $X$ by an open subset $D \subset \KK^n$ and instead call $S$ a \emph{closed analytic subvariety} of $D$.
\end{defn}

Fischer \cite[Section 0.14, p. 10]{Fischer_complex_analytic_geometry} defines a subset $S\subset X$ of a complex analytic space to be an analytic set if there is a coherent ideal $\sI \subset \sO_X$ such that $S = \supp(\sO_X/\sI)$. By Theorem \ref{thm:Oka-Cartan}, the requirement that $\sI$ have finite type is equivalent to $\sI$ being coherent.

\begin{defn}[Coherent real analytic sets and germs]
\label{defn:Coherent_real_analytic_set}
(See Narasimhan \cite[Section 5.1, Definition 2, p. 93]{Narasimhan_introduction_theory_analytic_spaces} for $\KK=\RR$ or $\CC$.)
Let $S$ be a real analytic set in an open set $U \subset \RR^n$. Then $S$ is \emph{coherent} if the sheaf of germs of real analytic functions vanishing on $S$ is a coherent sheaf of $\sO_U$-modules. If $p\in S$ is a point, one says that the germ $S_p$ is \emph{coherent} if it is induced by a coherent real analytic set.
\end{defn}

We shall encounter the following concept in the contexts of analytic spaces and schemes. Recall first that a ring $R$ is \emph{reduced} if it has no non-zero nilpotent elements (see G\"ortz and Wedhorn \cite[Appendix B, p. 554]{Gortz_Wedhorn_algebraic_geometry_v1}). 

\begin{defn}[Reduced locally ringed space]
\label{defn:Reduced_locally_ringed_space}  
A locally ringed space $(X,\sO_X)$ is called \emph{reduced} if all its local rings, $\sO_{X,p}$ for $p\in X$, are reduced rings.
\end{defn}  

Let $(X,\sO_X)$ be an analytic space. Closely following Fischer \cite[Sections 0.15 and 0.16, pp. 11--12]{Fischer_complex_analytic_geometry} (see also  Grauert and Remmert \cite[Section 4.2.5, pp. 86--87, and Section 4.3, pp. 87--90]{Grauert_Remmert_coherent_analytic_sheaves}) but allowing $\KK=\RR$ as well as $\KK=\CC$, we recall that the \emph{nilradical}
\[
  \sN_X = \sqrt{0} \subset \sO_X
\]
is defined to be the sheaf of ideals associated to the presheaf
\[
  U \mapsto \{f \in \sO_X(U): f^k = 0, \text{ for some } k \in \NN\}.
\]
(Recall that if $\fa$ is an ideal in a ring $A$, then $\sqrt{a} = \{a \in A: a^k \in \fa\text{ for some positive integer } k\}$ by Atiyah and Macdonald \cite[p. 8]{Atiyah_Macdonald_introduction_commutative_algebra}.) Thus, for every $x \in X$ we have
\[
  \sN_{X,x} = \{f \in \sO_{X,x}: f^k = 0, \text{ for some } k \in \NN\}.
\]
Let $\sC_X$ denote the sheaf of continuous $\KK$-valued functions on the topological space $X$. If $U \subset X$ is an open subset, there is a canonical homomorphism
\[
  \sO_X(U) \ni f \mapsto \tilde f \in \sC_X(U),
\]
where $\tilde f(x) := f(x)$ for $x \in X$ and this defines a canonical homomorphism $\vartheta:\sO_X \to \sC_X$.

\begin{thm}[Coherence of the nilradical of a complex analytic space]
\label{thm:Nilradical_analytic_space}
(See Fischer \cite[Section 0.15, Theorem, p. 11]{Fischer_complex_analytic_geometry} and Grauert and Remmert \cite[Section 4.2.5, Lemma, p. 86 and Proposition, p. 87]{Grauert_Remmert_coherent_analytic_sheaves})  
If $(X,\sO_X)$ is an analytic space over $\KK=\CC$ or $\RR$, then $\sN_X = \Ker(\vartheta:\sO_X \to \sC_X)$ and, if $\KK=\CC$, then $\sN_X \subset \sO_X$ is a coherent ideal. 
\end{thm}

When $\KK=\CC$, the coherence of $\sN_X$ in Theorem \ref{thm:Nilradical_analytic_space} follows from Theorem \ref{thm:Oka-Cartan}, which is false for $\KK=\RR$ (see Remark \ref{rmk:Ideal_sheaf_real_analytic_set_need_not_be_coherent}).

Let $(X,\sO_X)$ be an analytic space over $\KK=\CC$ or $\RR$. Since $\sN_{X,p} \neq \sO_{X,x}$ for every
$x \in X$, then $\supp (\sO_X/\sN_X) = X$ and
\[
  \sO_X^{\red} := \sO_X/\sN_X.
\]
Observe that $\sO_X/\sN_X$ is isomorphic to the image of the homomorphism $\vartheta:\sO_X \to \sC_X$. Following Fischer \cite[Section 0.16, p. 12]{Fischer_complex_analytic_geometry} and Grauert and Remmert \cite[Section 4.3.2, p. 88]{Grauert_Remmert_coherent_analytic_sheaves}), we note that $(X,\sO_X^{\red})$ is a closed analytic subspace of $(X,\sO_X)$ (see the forthcoming Definition \ref{defn:Analytic_subspace}), called its \emph{reduction}, and $(X,\sO_X)$ is said to be \emph{reduced} if $\sO_X^{\red} = \sO_X$. (Grauert and Remmert \cite[Section 4.3.2, p. 88]{Grauert_Remmert_coherent_analytic_sheaves} abbreviate $(X,\sO_X^{\red})$ as $\red X$, whereas we prefer the abbreviation $X^{\red}$). Equivalently, $(X,\sO_X)$ is reduced if it satisfies Definition \ref{defn:Reduced_locally_ringed_space}. 

\begin{exmp}[Non-reduced schemes or analytic spaces]
\label{exmp:Fischer_0-14_page_10}  
See Fischer \cite[Section 0.14, p. 10, Example]{Fischer_complex_analytic_geometry}, in the category of complex analytic spaces, and Eisenbud and Harris \cite[Section II.3.1, Example II-9, p. 58]{Eisenbud_Harris_geometry_schemes}, in the category of schemes, for examples of $n$-\emph{fold points}. Consider the complex analytic space $(X,\sO_X)$, where $X = \supp(\sO_{\CC}/(z^n)) = \{0\} \subset \CC$ and $\sO_X = \sO_\CC/(z^n) \restriction \{0\} = \CC + \CC\eps + \cdots + \CC\eps^{n-1}$, for an integer $n \geq 1$, and $\eps^n = 0$.
Thus, $\sN_{X,0} = 0$ when $n=1$ but $\sN_{X,0} \neq 0$ when $n>1$, while $\sO_X^{\red} = \sO_\CC/I_1\restriction \{0\}$ and $X_{\red} = (\{0\},\sO_X^{\red})$ for all integers $n \geq 1$.
\qed
\end{exmp}  

See Eisenbud and Harris \cite[Sections II.3.1 and II.3.2]{Eisenbud_Harris_geometry_schemes} for more involved examples and exercises in the category of schemes than that of Example \ref{exmp:Fischer_0-14_page_10}. We have the

\begin{lem}[Characterization of reduced analytic spaces]
\label{lem:Characterization_reduced_analytic_space}
(See Fischer \cite[Section 0.16, Lemma, p. 12]{Fischer_complex_analytic_geometry} and Grauert and Remmert \cite[Section 4.3.3, Criterion of Reducedness, p. 89]{Grauert_Remmert_coherent_analytic_sheaves} for $\KK=\CC$ and Guaraldo, Macr\`\i, and Tancredi \cite[Section 1.1, Definition 1.2, p. 2]{Guaraldo_Macri_Tancredi_topics_real_analytic_spaces} for $\KK$-ringed spaces with $\KK=\CC$ or $\RR$.)
Let $(X,\sO_X)$ be an analytic space over $\KK=\CC$ or $\RR$. Then the following conditions are
equivalent:
\begin{enumerate}
\item $X$ is reduced;
\item $\sN_X = 0$;
\item The canonical homomorphism $\vartheta:\sO_X \to \sC_X$ is injective.
\end{enumerate}
\end{lem}

\begin{defn}[Analytic variety]
\label{defn:Analytic_variety}
(See Acquistapace, Broglia, and Fernando \cite[Definition 1.2, Way 1]{Acquistapace_Broglia_Fernando_topics_global_real_analytic_geometry} for $\KK=\CC$ and \cite[Definition 2.1]{Acquistapace_Broglia_Fernando_topics_global_real_analytic_geometry} for $\KK=\RR$, Guaraldo, Macr\`\i, and Tancredi \cite[Section 2.1, Definition 1.4, p. 13]{Guaraldo_Macri_Tancredi_topics_real_analytic_spaces} for $\KK=\RR$ or $\CC$.)
An \emph{analytic variety} is a $\KK$-ringed space $(X,\sO_X)$ that is locally isomorphic, as a $\KK$-ringed space, to a local model $(Y,\sO_Y)$, where $U\subset\RR^n$ is a domain and $Y \subset U \subset \RR^n$ is the zero set of finitely many real analytic functions in $\sO_{\RR^n}(U)$ and $\sO_Y$ is the quotient sheaf of the sheaf $\sO_U$ of germs of real analytic functions on $U$ by the ideal sheaf $\sI_Y$ of germs of analytic functions vanishing on $Y$. 
\end{defn}

\begin{rmk}[Terminology for analytic varieties and reduced analytic spaces]
\label{rmk:Varieties_and_reduced_analytic_spaces}
Instead of the term `analytic variety' in Definition \ref{defn:Analytic_variety}, Acquistapace, Broglia, and Fernando employ the term `reduced complex analytic space' \cite[Definition 1.2, Way 1]{Acquistapace_Broglia_Fernando_topics_global_real_analytic_geometry} when $\KK=\CC$ and `reduced real analytic space' \cite[Definition 2.1]{Acquistapace_Broglia_Fernando_topics_global_real_analytic_geometry} when $\KK=\RR$. We prefer the term `analytic variety' employed by Guaraldo, Macr\`\i, and Tancredi, since the qualifier `reduced' is applicable to any ringed space (see \cite[Section 1.1, Definition 1.2, p. 2]{Guaraldo_Macri_Tancredi_topics_real_analytic_spaces} and the introduction to their monograph, especially \cite[p. vii]{Guaraldo_Macri_Tancredi_topics_real_analytic_spaces}).
\end{rmk}

Guaraldo, Macr\`\i, and Tancredi \cite[Section 2.1, p. 13]{Guaraldo_Macri_Tancredi_topics_real_analytic_spaces} note that the reduced $\KK$-ringed space $(X,\sO_X^{\red})$ associated to a $\KK$-analytic space $(X,\sO_X)$ is an analytic variety, called the \emph{associated analytic variety}. 

\begin{defn}[Ideal sheaf defined by an analytic set in an analytic space]
\label{defn:Ideal_sheaf_for_analytic_subset}  
(See Grauert and Remmert in \cite[Section 4.1.2, p. 77]{Grauert_Remmert_coherent_analytic_sheaves} or  Narasimhan in \cite[Section 3.1, p. 31 and Section 4.1, p. 64]{Narasimhan_introduction_theory_analytic_spaces}.)  
Let $\KK=\RR$ or $\CC$ and $(X,\sO_X)$ be an analytic space over $\KK$. If $S \subset X$ is an analytic set, then $\iota(S) \subset \sO_X$ denotes the \emph{sheaf of ideals of the analytic set} defined by the pre-sheaf,
\[
  X \supset U \mapsto \left\{f \in \sO_X(U): S \cap U \subset f^{-1}(0) \right\} \subset \sO_X(U),
\]
where $U \subset X$ is an open subset. 
\end{defn}  

\begin{thm}[Oka--Cartan Coherence Theorem]
\label{thm:Oka-Cartan}
(See Grauert and Remmert \cite[Section 4.2, Fundamental Theorem, p. 84]{Grauert_Remmert_coherent_analytic_sheaves} and Narasimhan \cite[Section 4.3, Theorem 5, p. 77]{Narasimhan_introduction_theory_analytic_spaces}.)  
Let $(X,\sO_X)$ be a complex analytic space. If $S \subset X$ is a complex analytic set, then the sheaf of ideals, $\iota(S)$, is a coherent $\sO_X$-sheaf.
\end{thm}

\begin{rmk}[The Oka--Cartan Coherence Theorem does not hold for real analytic sets]
\label{rmk:Ideal_sheaf_real_analytic_set_need_not_be_coherent}  
When $\KK=\RR$, the ideal sheaf $\sI_Y$ in Definition \ref{defn:Analytic_variety} need not to be coherent since Theorem \ref{thm:Oka-Cartan} does not hold when $\CC$ is replaced by $\RR$: see Acquistapace, Broglia, and Fernando \cite[Remark 2.1 and Example 2.9]{Acquistapace_Broglia_Fernando_topics_global_real_analytic_geometry} and Guaraldo, Macr\`\i, and Tancredi \cite[Section 2.1, Remark 1.3, p. 12]{Guaraldo_Macri_Tancredi_topics_real_analytic_spaces} for a discussion of counterexamples when $\KK=\RR$, including
\begin{align}
  \label{eq:Whitney_umbrella}
  \{(x_1,x_2,x_3) \in \RR^3: x_1^2-x_3x_2^2 = 0\} &\subset \RR^3,
  \\
  \label{eq:Cartan_umbrella}
  \{(x,y,z) \in \RR^3: x_1^3-x_3(x_1^2+x_2^2) = 0\} &\subset \RR^3,
\end{align}
called the \emph{Whitney umbrella} \cite[Section 4, p. 163]{Whitney_1943}
and the \emph{Cartan umbrella} \cite[Section 8, Remarque, p. 93]{Cartan_1957}, respectively. Narasimhan \cite[Section 5.3]{Narasimhan_introduction_theory_analytic_spaces} discusses the Cartan umbrella and four more striking examples due to Cartan \cite[Section 11]{Cartan_1957} and Bruhat and Cartan \cite{Bruhat_Cartan_1957a, Bruhat_Cartan_1957b}.
\end{rmk}

\begin{defn}[Locally closed subspace of a topological space]
\label{defn:Locally_closed_subspace_topological_space} 
(See Bourbaki \cite[Chapter I, Section 3.3, Definition 2, p. 38]{Bourbaki_general_topology_chapters_1-4}.)  
A subset $Y$ of a topological space $X$ is said to be \emph{locally closed at a point} $p \in Y$ if there is an open neighborhood $U$ of $p$ in $X$ such that $Y \cap U$ is a closed subset of the subspace $U$. Moreover, $Y$ is \emph{locally closed} in $X$ if it is locally closed at each point $p \in Y$.
\end{defn}

\begin{defn}[Open analytic subspace, closed analytic subspace, and locally closed analytic subspace]
\label{defn:Analytic_subspace}
(See Fischer \cite[Section 0.14, p. 10]{Fischer_complex_analytic_geometry} and Grauert and Remmert \cite[Section 1.2.2, pp. 14--15]{Grauert_Remmert_coherent_analytic_sheaves} for $\KK=\CC$ and Bierstone and Milman \cite[Section 2.1, p. 804]{Bierstone_Milman_1989} and Guaraldo, Macr\`\i, and Tancredi \cite[Section 2.1, Definition 1.9, p. 15]{Guaraldo_Macri_Tancredi_topics_real_analytic_spaces} for $\KK=\CC$ or $\RR$; compare De Jong and Pfister \cite[Definition 6.1.26 (3), p. 235]{DeJong_Pfister_local_analytic_geometry}.)
Let $(X,\sO_X)$ be an analytic space over $\KK=\CC$ or $\RR$. An analytic space $(Y,\sO_Y)$ is called an \emph{open analytic subspace} of an analytic space $X$ if $Y$ is an open subset of $X$ and $\sO_Y = \sO_X\restriction Y$, and one writes $Y \subset X$.

An analytic space $(Y,\sO_Y)$ is a \emph{closed analytic subspace} of $(X,\sO_X)$ if there is a coherent ideal $\sI \subset \sO_X$ such that $Y = \supp(\sO_X/\sI)$ and $\sO_Y = \sO_X/\sI\restriction Y$. There is a canonical analytic morphism determined by the injection, denoted by $Y \hookrightarrow X$.

An analytic space $(Y,\sO_Y)$ is a \emph{locally closed analytic subspace} of $(X,\sO_X)$ if for each point $p \in Y$, there is an open neighborhood $U$ of $p$ in $X$ such that $(Y\cap U,\sO_Y\restriction Y\cap U)$ is a closed analytic subspace of $(U,\sO_X\restriction U)$. \qed
\end{defn}  

Grauert and Remmert only assume in \cite[Section 1.2.2, pp. 14--15]{Grauert_Remmert_coherent_analytic_sheaves} that $\sI$ is of finite type but Theorem \ref{thm:Oka-Cartan} then implies that $\sI$ is coherent for $\KK=\CC$. In \cite[Definition 6.1.26 (3), p. 235]{DeJong_Pfister_local_analytic_geometry}, De Jong and Pfister say that $(Y,\sO_Y)$ a \emph{closed analytic subspace} of $(X,\sO_X)$ if there exists a map of locally ringed spaces $(\iota,\iota^\sharp):(Y,\sO_Y) \to (X,\sO_X)$ such that $\iota:Y\to X$ is injective, and $\iota^\sharp:\sO_X\to\sO_Y$ is surjective.

\section{Explicit description of morphisms of affine algebraic varieties and analytic model spaces}
\label{sec:Explicit_description_morphisms_analytic_varieties}
In this section \ref{sec:Explicit_description_morphisms_analytic_varieties}, we give a useful explicit description of morphisms of analytic model spaces, by analogy with prior explicit descriptions of morphisms of affine algebraic varieties. According to Shafarevich \cite[Section 5.3.1, Example 5.19, p. 29]{Shafarevich_v2}, there is a one-to-one correspondence between \emph{regular maps}\footnote{That is, morphisms \cite[pp. 25 and 47]{Shafarevich_v1} in the case of affine or quasiprojective varieties} \emph{of quasiprojective varieties} over $\KK$ and \emph{morphisms of schemes}. Shafarevich \cite[Section 5.3.1, Definition 5.4, p. 28]{Shafarevich_v2} defines a morphism of schemes to be a morphism of the corresponding ringed spaces. 

By analogy with Shafarevich \cite[Section 1.2.2, p. 25, Definition]{Shafarevich_v1} for regular functions on affine algebraic varieties (see also Milne \cite[Section 2j, p. 48]{Milne_algebraic_geometry} and Mumford \cite[Section 3, Definition 1, p. 12]{Mumford_red_book_varieties_schemes}) and \cite[Section 1.4.2, p. 46]{Shafarevich_v1} for regular functions on quasiprojective algebraic varieties, we make the

\begin{defn}[Regular function on a $\KK$-analytic model space]
\label{defn:Regular_function_analytic_variety}
For $\KK=\RR$ or $\CC$, let $(X,\sO_X)$ be a $\KK$-analytic model space as in Definition \ref{defn:Analytic_model_space}, so we are given a domain a domain $D \subset \KK^n$, an ideal sheaf $\sI \subset \sO_D$ of finite type, and $X = \supp(\sO_D/\sI) \subset D$ with structure sheaf $\sO_X = (\sO_D/\sI)\restriction X$. A function $f:X \to \KK$ is \emph{regular} if there is a $\KK$-analytic function $\tilde f:D\to \KK$ such that $f(x) = \tilde f(x)$ for all $x\in X$.
\end{defn}

By analogy with Shafarevich \cite[Section 1.2.3, p. 27, Definition]{Shafarevich_v1} for regular maps of affine algebraic varieties (see also Milne \cite[Section 1i, p. 47]{Milne_algebraic_geometry} and Mumford \cite[Section 3, Definition 1, p. 12]{Mumford_red_book_varieties_schemes}) and \cite[Section 1.4.2, p. 47, Definition]{Shafarevich_v1} for regular maps of quasiprojective algebraic varieties, we make the

\begin{defn}[Regular map of $\KK$-analytic model spaces]
\label{defn:Regular_map_analytic_variety}
Continue the notation of Definition \ref{defn:Regular_function_analytic_variety} and let $(Y,\sO_Y)$ be a $\KK$-analytic model space defined by a domain $D' \subset \KK^m$, ideal sheaf $\sJ \subset \sO_{D'}$ of finite type, so $Y = \supp(\sO_{D'}/\sJ) \subset D'$ with structure sheaf $\sO_Y = (\sO_{D'}/\sJ)\restriction Y$. A continuous map $\varphi:X \to Y$ is a \emph{regular map of $\KK$-analytic model spaces} if each component function $\varphi_i = y_i\circ\varphi:X\to\KK$ of $\varphi = (\varphi_1,\ldots,\varphi_m):X \to Y\subset\KK^m$ is a regular function in the sense of Definition \ref{defn:Regular_function_analytic_variety}, where $y_i$ is the coordinate function
\[
  y_i:\KK^m \ni (b_1,\ldots,b_m) \mapsto b_i \in \KK,
\]
for $i=1,\ldots,m$. 
\end{defn}

By analogy with Grauert and Remmert \cite[Sections 1.1.4 and 1.1.5, pp. 6--7]{Grauert_Remmert_coherent_analytic_sheaves} for morphisms of complex analytic spaces, we have the

\begin{defn}[Morphism of $\KK$-analytic model spaces]
\label{defn:Morphism_analytic_space}
Continue the notation of Definition \ref{defn:Regular_map_analytic_variety}. A \emph{morphism of $\KK$-analytic model spaces} $(\varphi,\varphi^\sharp):(X,\sO_X)\to(Y,\sO_Y)$ is a morphism of $\KK$-ringed spaces and thus comprises a continuous map $\varphi:X\to Y$ and a $\KK$-algebra homomorphism $\varphi^\sharp:\sO_Y \to \varphi_*\sO_X$, where the image sheaf $\varphi_*\sO_X$ is defined by $\varphi_*\sO_X(V) := \sO_X(\varphi^{-1}(V))$ for open subsets $V \subset Y$.
\end{defn}

It is well-known that there is a one-to-one correspondence between regular maps of affine algebraic varieties $X \subset \KK^n$ and $Y \subset \KK^m$ and homomorphisms of coordinate rings of affine algebraic varieties (or algebraic sets) $\KK[X] := \KK[x_1,\ldots,x_n]/\sI$ and $\KK[Y] := \KK[y_1,\ldots,y_m]/\sJ$ (see Milne \cite[Section 2a, p. 35, and Section 2i, p. 47 and Section 3g]{Milne_algebraic_geometry}), where $\sI \subset \KK[x_1,\ldots,x_n]$ and $\sJ \subset \KK[y_1,\ldots,y_m]$ are finitely-generated ideals and, more generally, between regular maps of quasiprojective algebraic varieties $X$ and $Y$ over $\KK$ and their associated schemes over $\KK$. See Milne \cite[Proposition 3.26 and Section 3g]{Milne_algebraic_geometry}, Shafarevich \cite[Section 5.3.1, p. 29, Example 5.19]{Shafarevich_v2}, and Mumford \cite[Section 3, Propositions 1 and 2, pp. 13--14]{Mumford_red_book_varieties_schemes} for an explicit description of morphisms of affine algebraic varieties.

\begin{defn}[Morphism of $\KK$-ringed spaces]
\label{defn:Morphism_K-ringed_spaces}
(See Milne \cite[Section 3b, p. 58 and Section 3d, p. 62]{Milne_algebraic_geometry}.)  
Let $\KK$ be field, $X$ and $Y$ be topological spaces, and $\sO_X$ and $\sO_Y$ be sheaves of $\KK$-algebras over $X$ and $Y$, respectively. A \emph{morphism of $\KK$-ringed spaces} $(\varphi,\varphi^\sharp):(X,\sO_X) \to (Y,\sO_Y)$ comprises a continuous map $\varphi:X \to Y$ and a homomorphism of $\KK$-algebras $\varphi^\sharp:\sO_Y\to \varphi_*\sO_X$ defined so that for every open set $V \subset Y$, there is a homomorphism of $\KK$-algebras
\[
  \sO_Y(V) \ni g \mapsto g\circ\varphi \in \sO_X(\varphi^{-1}(V))
\]
and these homomorphisms are compatible with restriction to smaller open subsets. A morphism of ringed spaces is an \emph{isomorphism} if it is bijective and its inverse is also a morphism of ringed spaces (in particular, it is a homeomorphism).
\end{defn}  

\begin{prop}[Explicit description of morphisms of affine algebraic varieties]
\label{prop:Milne_3-26_affine}
(See Milne \cite[Section 3g, Proposition 3.26, p. 65]{Milne_algebraic_geometry}.)  
Let $\KK$ be a field and $X \subset \KK^m$ and $X \subset \KK^n$ be affine algebraic varieties. The following conditions on a continuous map $\varphi:X\to Y$ are equivalent:
\begin{enumerate}
\item
  \label{item:varphi_sharp_morphism_ringed_spaces_affine}
  $(\varphi,\varphi^\sharp)$ is a morphism of ringed spaces $(X,\sO_X)\to (Y,\sO_Y)$.
\item
  \label{item:varphi_components_regular_functions_affine}
  The components $\varphi_1,\ldots,\varphi_m$ of $\varphi$ are regular functions on $X$ for $i=1,\ldots,m$.
\item
  \label{item:pullback_functions_in_coordinate_ring_affine}
  $g \in \KK[Y] \implies g\circ\varphi \in \KK[X]$.
\end{enumerate}
\end{prop}

We have the following analogue in the category of $\KK$-analytic model spaces of Proposition \ref{prop:Milne_3-26_affine} for morphisms of affine algebraic varieties over $\KK$. In this analogy, the coordinate ring $\KK[X]$ for an affine algebraic variety $X = \VV(f_1,\ldots,f_r) \subset \KK^n$ corresponds to $\sO_D/\sI$ for a $\KK$-analytic model space $(X,\sO_X)$ in a domain $D \subset \KK^n$ and ideal $\sI \subset \sO_D$ with generators $f_1,\ldots,f_r \in \sO_D(D) = \sO_{\KK^n}(D)$, where $X = \cosupp\sI$ and $\sO_X = (\sO_D/\sI)\restriction X$. 

\begin{prop}[Explicit description of morphisms of analytic spaces]
\label{prop:Milne_3-26}
Let $\KK=\RR$ or $\CC$ and $(X,\sO_X)$ and $(Y,\sO_Y)$ be $\KK$-analytic model spaces as in Definition \ref{defn:Morphism_analytic_space}. Then the following conditions on a continuous map $\varphi:X\to Y$ are equivalent:
\begin{enumerate}
\item
  \label{item:varphi_sharp_morphism_ringed_spaces}
  $(\varphi,\varphi^\sharp)$ is a morphism of ringed spaces $(X,\sO_X)\to (Y,\sO_Y)$;
\item
  \label{item:varphi_components_regular_functions}
  The components $\varphi_i$ of $\varphi$ are regular functions on $X$ for $i=1,\ldots,m$;
\item
  \label{item:pullback_functions_in_coordinate_ring}
  $g \in \sO_{D'}/\sJ \implies g\circ\varphi \in \sO_D/\sI$.
\end{enumerate}
\end{prop}

\begin{proof}
We adapt the proof of Proposition \ref{prop:Milne_3-26_affine} by Milne for affine algebraic varieties. Given a ring $A$, Milne uses $\Spec(A)$ to denote the set of prime ideals in $A$ and $\spm(A)$ to denote the set of maximal ideals $A$ (see Milne \cite[Section 1a, p. 13]{Milne_algebraic_geometry} in his proof of Proposition \ref{prop:Milne_3-26_affine}.

Consider \eqref{item:varphi_sharp_morphism_ringed_spaces} $\implies$ \eqref{item:varphi_components_regular_functions}. If $y_i:\KK^m\to\KK$ are the coordinate functions as in Definition \ref{defn:Regular_map_analytic_variety} for $i=1,\ldots,m$, then $y_i \in \sO_{\KK^m} \subset \sO_{D'}$ and the restriction to $Y = \cosupp\sJ$ of its image in the quotient $y_i \in \sO_{D'}/\sJ$ defines a regular function $y_i \in \sO_Y = (\sO_{D'}/\sJ)\restriction Y$ as in Definition \ref{defn:Regular_function_analytic_variety}. By assumption for Item \eqref{item:varphi_sharp_morphism_ringed_spaces}, $\varphi^\sharp:\sO_Y \to \varphi_*\sO_X$ is a $\KK$-algebra homomorphism and by Definition \ref{defn:Morphism_K-ringed_spaces}, for each open set $V \subset Y$ and $U \subset \varphi^{-1}(V) \subset X$, this $\KK$-algebra homomorphism is induced by
\[
  \sO_Y(V) \ni g \mapsto g\circ\varphi \in \sO_X(U).
\]
By definition, the components of $\varphi:X\to Y\subset D' \subset \KK^m$ are given by $\varphi_i = y_i\circ\varphi$ and so $\varphi_i = y_i\circ\varphi \in \sO_X(U)$, for $i=1,\ldots,m$. Hence, $\varphi_i$ is a regular function on $X$ by Definition \ref{defn:Regular_function_analytic_variety}, for $i=1,\ldots,m$, and this verifies Item \eqref{item:varphi_components_regular_functions}.

Consider \eqref{item:varphi_components_regular_functions} $\implies$ \eqref{item:pullback_functions_in_coordinate_ring}. If well-defined on all of $\sO_{D'}/\sJ$, the map
\[
  \alpha:\sO_{D'}/\sJ \ni g \mapsto g\circ\varphi \in \sO_D/\sI
\]
would be a $\KK$-algebra homomorphism from the $\KK$-algebra of all regular ($\KK$-analytic) functions $Y\to \KK$ into the $\KK$-algebra of all regular ($\KK$-analytic) functions $X\to\KK$. Item \eqref{item:varphi_components_regular_functions} asserts that $\alpha$ maps the coordinate functions $y_i$ on $Y \subset D' \subset \KK^m$ to $\alpha(y_i) = y_i\circ\varphi = \varphi_i \in \sO_D/\sI$, for $i=1,\ldots,m$. But the set $\{y_1,\ldots,y_m\}$ of coordinate functions on $D' \subset \KK^m$ generates $\sO_{D'}$ and thus $\sO_{D'}/\sJ$ as $\KK$-algebras, so the homomorphism $\alpha$ maps all of $\sO_{D'}/\sJ$ into $\sO_D/\sI$.

Consider \eqref{item:pullback_functions_in_coordinate_ring} $\implies$ \eqref{item:varphi_sharp_morphism_ringed_spaces}.
The homomorphism $\alpha:\sO_{D'}/\sJ \to \sO_D/\sI$ defines a map $\spm(\sO_{D'}/\sJ)\to\spm(\sO_D/\sI)$, and it remains to show that this coincides with $\varphi$ when we identify $\spm(\sO_D/\sI)$ with $X$ and identify $\spm(\sO_{D'}/\sJ)$ with $Y$. Let $p \in X$ and $q = \varphi(p) \in Y$ and let $\fm_p$ and $\fm_q$ be the (maximal) ideals of elements of $\sO_D/\sI$ and $\sO_{D'}/\sJ$ that are zero at $p$ and $q$, respectively. Then, for $g \in \sO_{D'}/\sJ$ and noting that $\alpha(g) = g\circ\varphi$, we have
\[
  \alpha(g) \in \fm_p \iff g(\varphi(p)) = 0 \iff g(q) = 0 \iff g \in \fm_q.
\]
Therefore, $\alpha^{-1}(\fm_p) = \fm_q$, which is what we needed to show. This completes the proof of Proposition \ref{prop:Milne_3-26}.
\end{proof}

\section{Complexification of a real analytic variety or space}
\label{sec:Complexification_real_analytic_spaces}
We refer to Acquistapace, Broglia, and Fernando \cite[Section 2.A]{Acquistapace_Broglia_Fernando_topics_global_real_analytic_geometry}, Fensch \cite[Section 2.1]{Fensch_1966}, Guaraldo, Macr\`\i, and Tancredi \cite[Chapter 3]{Guaraldo_Macri_Tancredi_topics_real_analytic_spaces}, Hironaka \cite[Chapter 0, Section 2, pp. 119--120]{Hironaka_1964-I-II}, and Narasimhan \cite[Section 5.1]{Narasimhan_introduction_theory_analytic_spaces} for details concerning complexification of real analytic germs and real analytic spaces. 

Following Acquistapace, Broglia, and Fernando \cite[Section 2.A.1]{Acquistapace_Broglia_Fernando_topics_global_real_analytic_geometry}, we assume that $\RR^n$ is canonically embedded in $\CC^n$ as the fixed-point set of complex conjugation. We describe when a given real analytic space can be viewed as a subspace of a complex analytic space, first locally and then globally.

\begin{prop}[Complexification of the germ of a real analytic set at a point]
\label{prop:Acquistapace_Broglia_Fernando_2-3}  
(See Acquistapace, Broglia, and Fernando \cite[Proposition 2.3]{Acquistapace_Broglia_Fernando_topics_global_real_analytic_geometry} and Guaraldo, Macr\`\i, and Tancredi \cite[Section 3.1, Remark 1.3, p. 41 and Section 3.2, Proposition 2.1, p. 44]{Guaraldo_Macri_Tancredi_topics_real_analytic_spaces}.)
If $A_x$ is the germ of a real analytic set at a point $x \in \RR^n \subset \CC^n$, then there
exists in $\CC^n$ a unique germ of a complex analytic set $B_x$ such that
\begin{enumerate}
\item $A_x \subset B_x$ and $B_x \cap \RR^n = A_x$;
\item If the germ at $x$ of holomorphic function vanishes on $A_x$,  then it vanishes on $B_x$;
\item $B_x$ is minimal with respect to the preceding two properties;
\item If $\sI_x$ is the ideal of germs of functions vanishing on $A_x$ and $\sJ_x$ the one of
germs in $\sO_{\CC^n,x} = \sO_{\RR^n,x}\otimes_\RR\CC$ vanishing on $B_x$, then $\sJ_x = \sI_x\otimes_\RR\CC$.
\end{enumerate}
\end{prop}

\begin{defn}[Complexification of the germ of a real analytic set at a point]
\label{defn:Acquistapace_Broglia_Fernando_2-4}  
(See Acquistapace, Broglia, and Fernando \cite[Definition 2.4]{Acquistapace_Broglia_Fernando_topics_global_real_analytic_geometry}.)
Let $A_x$ be the germ of a real analytic set at a point $x \in \RR^n \subset \CC^n$. The \emph{complexification} of $A_x$ is the germ at $x\in\CC^n$ of a complexification of a representative of $A_x$.
\end{defn}

\begin{prop}[Necessary and sufficient condition for the germ of a real analytic set at a point to be coherent]
\label{prop:Acquistapace_Broglia_Fernando_2-6}  
(See Acquistapace, Broglia, and Fernando \cite[Proposition 2.6]{Acquistapace_Broglia_Fernando_topics_global_real_analytic_geometry} and Guaraldo, Macr\`\i, and Tancredi \cite[Section 3.2, Proposition 2.8, p. 47]{Guaraldo_Macri_Tancredi_topics_real_analytic_spaces}.)
Let $A_x$ be the germ of a real analytic at a point $x \in \RR^n \subset \CC^n$ and $B$ be a complex analytic set in an open neighborhood in $\CC^n$ of $x$ such that $B_x$ is the complexification of $A_x$. Then $A_x$ is coherent if and only if for all $y$ close to $x$, the germ $B_y$ is the complexification of $A_y$.
\end{prop}

Following Acquistapace, Broglia, and Fernando \cite[Section 2.A.3]{Acquistapace_Broglia_Fernando_topics_global_real_analytic_geometry}, a complex analytic space $X$ can be viewed as a real analytic space by identifying $\CC$ with $\RR^2$: one calls this structure the \emph{underlying real structure on $X$} and denotes this space by $X^\RR$. In this sense, a real analytic space $Y$ can be viewed as a (closed, real analytic) subspace of a complex analytic space $X$ if $Y$ is a (closed, real analytic) subspace of $X^\RR$.

\begin{defn}[Complexification of a real analytic space or variety]
\label{defn:Acquistapace_Broglia_Fernando_2-11}  
(See Acquistapace, Broglia, and Fernando \cite[Definition 2.11]{Acquistapace_Broglia_Fernando_topics_global_real_analytic_geometry} and Guaraldo, Macr\`\i, and Tancredi \cite[Section 3.3, Definition 3.1 and Remark 3.2, p. 51]{Guaraldo_Macri_Tancredi_topics_real_analytic_spaces}.)
Let $X$ be a real analytic space or variety\footnote{In the terminology of Acquistapace, Broglia, and Fernando, $X$ could be a reduced real analytic space or a $C$-analytic space, although in their \cite[Definition 2.11]{Acquistapace_Broglia_Fernando_topics_global_real_analytic_geometry} it is simply called a `real analytic space'.} (in the sense of Definitions \ref{defn:Analytic_space} or \ref{defn:Analytic_variety}, respectively, with $\KK=\RR$). A complex analytic space or variety $\widetilde X$ (in the sense of Definitions \ref{defn:Analytic_space} or \ref{defn:Analytic_variety}, respectively, with $\KK=\CC$) is called a \emph{complexification} of $X$ if
\begin{enumerate}
\item $X$ is a closed real analytic set in $\widetilde X$, considered as a real analytic variety; and
\item For any $x \in X$, the germ $\widetilde X_x$ is the complexification of $X_x$.
\end{enumerate}
\end{defn}

If a real analytic variety $(X,\sO_X)$ has a complexification, then Proposition \ref{prop:Acquistapace_Broglia_Fernando_2-6} implies that it is a coherent real analytic variety. The following theorem yields a converse of this fact.

\begin{thm}[Existence of a complexification for a real analytic space]
\label{thm:Acquistapace_Broglia_Fernando_2-12}  
(See Acquistapace, Broglia, and Fernando \cite[Theorem 2.12]{Acquistapace_Broglia_Fernando_topics_global_real_analytic_geometry} and Guaraldo, Macr\`\i, and Tancredi \cite[Section 3.3, Theorem 3.3, p. 52]{Guaraldo_Macri_Tancredi_topics_real_analytic_spaces}.)
If $(X,\sO_X)$ is a real analytic space (in the sense of Definition \ref{defn:Analytic_space} with $\KK=\RR$), then there exists a complex analytic space $(Z,\sO_Z)$ (in the sense of Definition \ref{defn:Analytic_space} with $\KK=\CC$) such that
\begin{enumerate}
\item $X$ embeds in $Z$ as a closed subset;
\item For all $x\in X$, one has $\sO_{Z,x} = \sO_{X,x}\otimes_\RR\CC$.
\end{enumerate}
Moreover if $Z_1$ and $Z_2$ are two such complex analytic spaces, the identity map $\id:(X,\sO_X) \to (X,\sO_X)$ extends to an isomorphism
\[
  \phi:(U_1,\sO_{Z_1}\restriction U_1) \to (U_2,\sO_{Z_2}\restriction U_2)
\]
between two open subsets of $Z_1$ and $Z_2$ that are open neighborhoods of the images of $X$.
\end{thm}

\begin{rmk}
\label{rmk:Acquistapace_Broglia_Fernando_2-13}  
The complex analytic space $(Z,\sO_Z)$ in Theorem \ref{thm:Acquistapace_Broglia_Fernando_2-12} depends on the structure sheaf $\sO_X$ of the real analytic space $X$ (see Acquistapace, Broglia, and Fernando \cite[Remark 2.13]{Acquistapace_Broglia_Fernando_topics_global_real_analytic_geometry}). If $(X,\sO_X)$ is a reduced real analytic space, then its complexification $(Z,\sO_Z)$ is a reduced complex analytic space (see Guaraldo, Macr\`\i, and Tancredi \cite[Section 3.3, Theorem 3.3, p. 52]{Guaraldo_Macri_Tancredi_topics_real_analytic_spaces}).
\end{rmk}

\begin{cor}[Existence of a complexification for a real analytic space]
\label{cor:Guaraldo_Macri_Tancredi_3-5}  
(See Guaraldo, Macr\`\i, and Tancredi \cite[Section 3.3, Corollary 3.5, p. 54]{Guaraldo_Macri_Tancredi_topics_real_analytic_spaces}.)
A real analytic variety (in the sense of Definition \ref{defn:Analytic_variety}) is coherent if and only if it admits a complexification.
\end{cor}

\section{Real part of a complex analytic space}
\label{sec:Real_part_complex_analytic_spaces}
We refer to Acquistapace, Broglia, and Fernando \cite[Section 2.B]{Acquistapace_Broglia_Fernando_topics_global_real_analytic_geometry} and Guaraldo, Macr\`\i, and Tancredi \cite[Chapter 4]{Guaraldo_Macri_Tancredi_topics_real_analytic_spaces} for details concerning taking the real part of a complex analytic space.

Following Acquistapace, Broglia, and Fernando \cite[Section 2.B]{Acquistapace_Broglia_Fernando_topics_global_real_analytic_geometry}, those authors observe that there are at least two ways to view a real analytic space as a closed subset of a complex analytic space. One approach is to consider local complex model spaces $(Z,\sO_Z) \subset (\Omega,\sO_\Omega) \subset (\CC^n,\sO_{\CC^n})$ such that $\Omega$ is stable under complex conjugation and take $X = Z\cap\RR^n$. Another approach is to equip the complex analytic space with an anti-holomorphic involution $\sigma:Z\to Z$ in such a way that $X$ is the fixed point set $Z^\sigma$ of $\sigma$.

\begin{defn}[Anti-holomorphic involution and real part of a complex analytic space]
\label{defn:Acquistapace_Broglia_Fernando_2-16}  
(See Acquistapace, Broglia, and Fernando \cite[Definition 2.16]{Acquistapace_Broglia_Fernando_topics_global_real_analytic_geometry} and Guaraldo, Macr\`\i, and Tancredi \cite[Section 4.1, Definition 1.1 and Remark 1.2, p. 60]{Guaraldo_Macri_Tancredi_topics_real_analytic_spaces}.)
We define the following concepts.
\begin{enumerate}
\item Let $Z, W$ be two complex analytic spaces. A continuous map $\varphi:W \to Z$ is
\emph{anti-holomorphic} if for any $z \in Z$ there exist open neighborhoods $U$ of $z$ and $V$ of $\varphi(z)$ and realizations of $U$, respectively $V$, as local models $U' \subset \CC^n$ and $V' \subset \CC^m$ in such a way that $\varphi':U' \to V'$ is the restriction of an antiholomorphic map
between open subsets of $\CC^n$ and $\CC^m$.

\item A map $\sigma:Z \to Z$ is an \emph{anti-involution} if it is anti-holomorphic and $\sigma\circ\sigma = \id$.

\item Let $Z$ be a complex analytic space and $X$ be a real analytic subspace of $Z^\RR$. We call $X$ the \emph{real part} of $Z$ if there exists an open covering $U_i$ of $X$ and realizations $\varphi_i: U_i  \to U_i' \subset \CC^n$ in such a way that $\varphi_i(U_i\cap X) = \varphi_i(U_i)\cap\RR^n$.
\end{enumerate}
\end{defn}

\begin{rmk}
\label{rmk:Acquistapace_Broglia_Fernando_2-17}  
(See Acquistapace, Broglia, and Fernando \cite[Remark 2.17]{Acquistapace_Broglia_Fernando_topics_global_real_analytic_geometry}.)
Given a real analytic space $(X,\sO_X)$, its $\sO_X$-complexification is constructed from its real local models (see Theorem \ref{thm:Acquistapace_Broglia_Fernando_2-12}, its proof in \cite{Acquistapace_Broglia_Fernando_topics_global_real_analytic_geometry}, and Remark \ref{rmk:Acquistapace_Broglia_Fernando_2-13}). Hence, $(X,\sO_X)$ is always the real part
of its $\sO_X$-complexification when the local models are invariant under complex conjugation.
\end{rmk}

We have the following 

\begin{thm}[Characterization of real analytic spaces]
\label{thm:Acquistapace_Broglia_Fernando_2-18}  
(See Acquistapace, Broglia, and Fernando \cite[Theorem 2.18]{Acquistapace_Broglia_Fernando_topics_global_real_analytic_geometry}.)
Let $X$ be a (connected) real analytic set. The following conditions are equivalent.
\begin{enumerate}
\item $X$ is a real analytic space.

\item $X$ is a real part of a complex analytic space.

\item $X$ is a fixed part of a complex analytic space under an anti-involution.
\end{enumerate}
\end{thm}

\section{Gluing analytic spaces}
\label{sec:Gluing_analytic_spaces}
It is worth recalling that, as in the case of schemes (see the forthcoming Section \ref{sec:Gluing_schemes}), one may construct $\KK$-analytic spaces by gluing. We adapt the discussions by Fischer \cite[Section 0.24, p. 20]{Fischer_complex_analytic_geometry} and Grauert and Remmert \cite[Section 1.1.7, p. 10]{Grauert_Remmert_coherent_analytic_sheaves} in the category of complex analytic spaces.

Let $X$ be a topological space and $\KK$ be $\RR$ or $\CC$. Any $\KK$-analytic space $(U, \sO_U)$, where $U$ is an open set in $X$, is called a $\KK$-\emph{analytic chart} on $X$. A pair comprising a family $\{(U_i, \sO_i)\}_{i\in I}$ of $\KK$-analytic charts on $X$ and a family $\{\theta_{ij}\}_{i,j\in I}$ of $\KK$-algebra isomorphisms
\[
  \theta_{ij}:\sO_i\restriction U_i\cap U_j \to \sO_i\restriction U_i\cap U_j, \quad\text{for all } i, j \in I,
\]
is called a $\KK$-\emph{analytic atlas} on $X$ if $\{(U_i, \sO_i)\}_{i\in I}$ is a covering of $X$ and the isomorphisms, called the \emph{gluing isomorphisms of the atlas}, obey the compatibility conditions:
\[
  \theta_{ij}\circ\theta_{jk} = \theta_{ik} \quad\text{on } U_i\cap U_j \cap U_k.
\]
One has the

\begin{lem}[Gluing lemma for analytic spaces]
\label{lem:Gluing_lemma_analytic_spaces}  
(See Fischer \cite[Section 0.24, p. 20]{Fischer_complex_analytic_geometry} or Grauert and Remmert \cite[Section 1.1.7, Gluing Lemma, p. 10]{Grauert_Remmert_coherent_analytic_sheaves} for $\KK=\CC$.)
Let $\KK$ be $\RR$ or $\CC$ and $\{(U_i, \sO_i)\}_{i\in I}$ and $\{\theta_{ij}\}_{i,j\in I}$ be a $\KK$-analytic atlas on a Hausdorff topological space $X$. Then there exists, up to isomorphism, a unique $\KK$-analytic space $(X,\sO_X)$ and $\KK$-algebra isomorphisms $\varphi_i:\sO_X\restriction U_i \to \sO_i$, for all $i\in I$, such that the following \emph{gluing equations} are obeyed for all $i,j \in I$:
\[
  \theta_{ij} = \varphi_i\circ\varphi_j^{-1} \quad\text{on } U_i\cap U_j.
\]
\end{lem}

\chapter{Smooth points, regular points, and dimensions for schemes}
\label{chap:Smooth_point_regular_point_dimension_scheme}
In this chapter, we provide concise guide to the basic definitions and results for schemes that we shall need in our work. Section \ref{sec:Schemes} contains a summary of the key definitions underlying the concept of a scheme. In Section \ref{sec:Analyticification_schemes_over_complex_number_field}, we describe the concept of analyticification of schemes over the complex number field, $\CC$. Section \ref{sec:Dimensions_schemes} contains a guide to dimension theory for schemes. In Section \ref{sec:Smooth_point_regular_point_schemes}, we review the concepts of smooth points and regular points for schemes.

\section{Schemes}
\label{sec:Schemes}
In this section, we review some of the terminology one encounters in the definition and application of schemes. Further details are provided by G\"ortz and Wedhorn \cite{Gortz_Wedhorn_algebraic_geometry_v1}, Hartshorne \cite{Hartshorne_algebraic_geometry}, and Vakil \cite{Vakil_foundations_algebraic_geometry}.

Given a ring $R$, one defines \cite[Chapter 2, Equation (2.1.1), p. 42]{Gortz_Wedhorn_algebraic_geometry_v1}
\[
  \Spec R := \{\fp \subset R: \fp \text{ is a prime ideal}\}.
\]
If $\fa \subset R$ is an ideal, then $\VV(\fa)$ is the set of prime ideals of $R$ containing $\fa$ \cite[Chapter 2, p. 42]{Gortz_Wedhorn_algebraic_geometry_v1}. The \emph{spectrum} of $R$ is the set $\Spec R$ endowed with the \emph{Zariski topology}, whose closed sets are given by $\VV(\fa)$, where $\fa \subset R$ is an ideal \cite[Chapter 2, Definition 2.2, p. 42]{Gortz_Wedhorn_algebraic_geometry_v1}. We refer to \cite[Section 2.10, p. 59]{Gortz_Wedhorn_algebraic_geometry_v1} for the definition of the structure sheaf $\sO_{\Spec R}$. If $R$ is a field or any ring with a single prime ideal, then $\Spec R$ consists of a single point.  

A locally ringed space $(X,\sO_X)$ is an \emph{affine scheme} if there exists a ring $R$ such that $(X,\sO_X)$ is isomorphic to $(\Spec R,\sO_{\Spec R})$; a morphism of affine schemes is a morphism of locally ringed spaces \cite[Chapter 2, Definition 2.34, p. 60]{Gortz_Wedhorn_algebraic_geometry_v1}.

A \emph{scheme} is a locally ringed space $(X,\sO_X)$ (usually denoted simply by $X$) which admits an open covering $X = \cup_i U_i$ such that all locally ringed spaces $(U_i,\sO_X\restriction U_i)$ are \emph{affine schemes}; a \emph{morphism of schemes} is a morphism of locally ringed space \cite[Chapter 3, Definition 3.1, p. 68]{Gortz_Wedhorn_algebraic_geometry_v1}.

A scheme $X$ is \emph{locally Noetherian} if it admits an affine open cover $X = \cup_iU_i$ such that all the affine coordinate rings $\sO_X(U_i)$ are Noetherian; if in addition $X$ is quasi-compact\footnote{In algebraic geometry, a topological space $X$ is called \emph{quasi-compact} if every open covering of $X$ has a finite subcover and \emph{compact} if in addition $X$ is Hausdorff (see G\"ortz and Wedhorn \cite[Chapter 1, Definition 1.22, p. 14]{Gortz_Wedhorn_algebraic_geometry_v1}).}, then $X$ is \emph{Noetherian} \cite[Chapter 3, Definition 3.18, p. 76]{Gortz_Wedhorn_algebraic_geometry_v1}.

For a scheme $S$, a \emph{scheme $X$ over $S$} means a morphism $X \to S$ of schemes; a \emph{scheme $X$ over a commutative ring $R$} means a morphism $X \to \Spec R$. 
If $k$ is a field, then a \emph{$k$-scheme} is a scheme $X$ together with a morphism $X \to \Spec k$.

Let $\KK$ be a field, and let $X \to \Spec \KK$ be a $\KK$-scheme. One calls $X$ a \emph{$\KK$-scheme locally of finite type} or says that $X$ is \emph{locally of finite type over $\KK$} if there is an affine open cover $X = \cup_i U_i$ such that for all $i$, one has that $U_i = \Spec A_i$ is the spectrum of a finitely generated $\KK$-algebra $A_i$; one says that $X$ is of \emph{finite type over $\KK$} if $X$ is locally of finite type and quasi-compact \cite[Chapter 3, Definition 3.30, p. 81]{Gortz_Wedhorn_algebraic_geometry_v1}.

\section{Analyticification of schemes over the complex number field}
\label{sec:Analyticification_schemes_over_complex_number_field}
Any scheme of finite type over a complete normed field $\KK$ naturally determines an analytic space over $\KK$; this correspondence between schemes and analytic spaces over $\KK$ for reduced complex spaces was first explored by Serre \cite{Serre_1956}. See Hartshorne \cite[Appendix B.1, p. 438]{Hartshorne_algebraic_geometry} and Neeman \cite{Neeman_algebraic_and_analytic_geometry} for expositions of some of those ideas.

If $X$ is a scheme of finite type over $\CC$, one may define the associated complex analytic space $X^{\mathrm{an}}$ following the outline provided by Hartshorne \cite[Appendix B.1, p. 439]{Hartshorne_algebraic_geometry}. (If $A$ is a ring and a scheme $X$ over $A$ can be covered by affine open sets $\Spec B_i$ where each $B_i$ is a finitely generated $A$-algebra, then $X$ is \emph{locally of finite type over $A$}; furthermore, if $X$ is quasicompact (see Vakil \cite[Section 3.6.5]{Vakil_foundations_algebraic_geometry}), then $X$ is of \emph{finite type over $A$}. See Shafarevich \cite[Section 5.3.5, Definition, p. 37]{Shafarevich_v2} and Vakil \cite[Section 5.3.6]{Vakil_foundations_algebraic_geometry}.)

Following Hartshorne, cover $X$ by open affine subsets $Y_i := \Spec A_i$. Each $A_i$ is an algebra of finite type over $\CC$, so $A_i \cong \CC[x_1,\ldots,x_n]/(f_1,\ldots,f_q)$, where $f_1,\ldots,f_q$ are polynomials in the indeterminates $x_1,\ldots,x_n$ and $(f_1,\ldots,f_q)$ denotes the ideal with those generators in the polynomial ring $\CC[x_1,\ldots,x_n]$. The $f_i$ are also holomorphic functions on $\CC^n$, so the zero set
\[
  \{x \in \CC^n: f_1(x) = \cdots = f_q(x) = 0\}
\]
defines a complex analytic subspace $Y_i^{\mathrm{an}}\subset \CC^n$. The scheme $X$ is obtained by gluing the open sets $Y_i \subset X$, so the same gluing data can be used to glue the analytic spaces  $Y_i^{\mathrm{an}}$ to produce an analytic space $X^{\mathrm{an}}$. (For the construction of complex analytic spaces by \emph{gluing}, see Grauert and Remmert \cite[Section 1.1.7, p. 10]{Grauert_Remmert_coherent_analytic_sheaves} or Fischer \cite[Section 0.24. p. 20]{Fischer_complex_analytic_geometry} and for the construction of schemes by gluing, see Hartshorne \cite[Chapter II, Section 2, Example 2.3.5, p. 75, and Exercise 2.12, p. 80]{Hartshorne_algebraic_geometry} or Shafarevich \cite[Section 5.3.2]{Shafarevich_v2}.)

Well-known results due to Chow \cite{Chow_WL_1949}, Serre \cite{Serre_1956}, and Grothendieck and Raynaud \cite{Grothendieck_Raynaud_SGA1} compare the categories of analytic spaces and schemes of finite type over $\CC$. For more recent published expositions of their work, we refer to Hartshorne \cite[Appendix B]{Hartshorne_algebraic_geometry} and Neeman \cite{Neeman_algebraic_and_analytic_geometry}; for useful concise, unpublished introductions to their work, we refer to the lecture notes by Halleck--Dub\'e \cite{Halleck-Dube_gaga}, Kedlaya \cite{Kedlaya_2009_gaga}, Warner \cite{Warner_chows_theorem}, and Zhao \cite{Zhao_2013}. The article by Huisman \cite{Huisman_2002} discusses the extent which these results carry over when $\CC$ is replaced by $\RR$.

\section{Dimensions for schemes}
\label{sec:Dimensions_schemes}
In this section, we discuss the meaning of dimension of a scheme.

\begin{defn}[Krull dimension]
\label{defn:Krull_dimension}
(See Abhyankar \cite[Chapter 2, Section 17, Definition, p. 142]{Abhyankar_local_analytic_geometry}, Atiyah and Macdonald \cite[Chapter 8, p. 89]{Atiyah_Macdonald_introduction_commutative_algebra}, Hartshorne \cite[Chapter I, Section 1, Definition, p. 6]{Hartshorne_algebraic_geometry}, Gopalakrishnan \cite[Section 8.2, p. 198 for rings and p. 199 for modules]{Gopalakrishnan_commutative_algebra}, G\"ortz and Wedhorn \cite[Section B.12, p. 569]{Gortz_Wedhorn_algebraic_geometry_v1}, Matsumura \cite[Section 2.5, p. 30]{Matsumura_commutative_ring_theory}, Raghavan, Balwant Singh, and Sridharan \cite[Section 3.2, p. 56]{Raghavan_BalwantSingh_Sridharan_homological_methods_commutative_algebra} (for rings and modules), the Stacks Project  \cite[\href{https://stacks.math.columbia.edu/tag/00KE}{Definition 00KE}]{stacks-project}, or Vakil \cite[Definition 11.1.1]{Vakil_foundations_algebraic_geometry}.)  
Let $R$ be a commutative ring. The \emph{Krull dimension} of $R$ is the supremum of the lengths of strictly increasing chains of prime ideals of $R$,
\begin{equation}
  \label{eq:Krull_dimension}
  \dim R := \sup\left\{l \in \NN: \fp_0 \subsetneq \fp_1 \subsetneq \cdots \subsetneq \fp_l\right\},
\end{equation}
where $\fp_j \subsetneq R$ is a prime ideal for $j=0,\ldots,l$.
\end{defn}

According to Atiyah and Macdonald \cite[Corollary 11.15, p. 121]{Atiyah_Macdonald_introduction_commutative_algebra} (see also Gopalakrishnan \cite[Section 8.2, Corollary 2, p. 201]{Gopalakrishnan_commutative_algebra} or Vakil \cite[Theorem 12.2.1]{Vakil_foundations_algebraic_geometry}), one has the inequality
\begin{equation}
  \label{eq:Krull_dimension_leq_dimension_Zariski_cotangent_space}
  \dim R \leq \dim_\kappa\fm/\fm^2,
\end{equation}
for any Noetherian local ring $(R,\fm,\kappa)$.

\begin{rmk}[Noetherian rings]
\label{rmk:Noetherian_rings}
Recall that a commutative ring $R$ is \emph{Noetherian} if it satisfies the \emph{ascending chain condition}, that is, for every sequence $\fa_1 \subseteq \fa_2 \subseteq \cdots \subseteq \fa_n \subseteq \cdots$ of ideals in $R$, there is a largest element (see Atiyah and Macdonald \cite[Chapter 7, p. 80]{Atiyah_Macdonald_introduction_commutative_algebra} for equivalent conditions for $R$ to be Noetherian). If $R$ is commutative Noetherian ring, then the rings $R[x_1,\ldots,x_n]$ of polynomials and $R[[x_1,\ldots,x_n]]$ of formal power series are also Noetherian (see Atiyah and Macdonald \cite[Corollary 7.6, p. 81, and Corollary 10.27, p. 113]{Atiyah_Macdonald_introduction_commutative_algebra}). In particular, if $\KK$ is any field, then $\KK[x_1,\ldots,x_n]$ and $\KK[[x_1,\ldots,x_n]]$ are Noetherian rings. Zariski and Samuel note that their proof that $\KK[[x_1,\ldots,x_n]]$ is a Noetherian ring (see \cite[Chapter VII, Section 1, Theorem $4'$, pp. 147--148]{Zariski_Samuel_communtative_algebra_II}) also shows that the ring $\KK\{x_1,\ldots,x_n\}$ of convergent power series (for $\KK=\RR$ or $\CC$) is Noetherian.

Here, $\KK\{x_1,\ldots,x_n\}$ is the ring of power series with positive radius of convergence and can be identified with the ring $\sO_{\KK^n,0}$ of germs of analytic functions at the origin in $\KK^n$. However, if $\sO_{\KK^n}(U)$ is the ring of analytic functions on a domain $U\subseteq\KK^n$, then $\sO_{\KK^n}(U)$ is not Noetherian since not every ideal in $\sO_{\KK^n}(U)$ is finitely generated.
\end{rmk}  

\begin{defn}[Regular local ring]
\label{defn:Regular_local_ring}
(See Atiyah and Macdonald \cite[Theorem 11.22]{Atiyah_Macdonald_introduction_commutative_algebra}, G\"ortz and Wedhorn \cite[Definition/Proposition B.76]{Gortz_Wedhorn_algebraic_geometry_v1}, Hartshorne \cite[Chapter I, Section 5, Definition, p. 32]{Hartshorne_algebraic_geometry}, and Vakil \cite[Definition 12.2.3]{Vakil_foundations_algebraic_geometry}.)
Let $R$ be a Noetherian local ring with maximal ideal $\fm$ and residue field $\kappa = R/\fm$. Then $R$ is a \emph{regular local ring} if
\begin{equation}
  \label{eq:Krull_dimension_equals_dimension_Zariski_cotangent_space_regular_local_ring}
  \dim R = \dim_\kappa\fm/\fm^2,
\end{equation}
where $\dim R$ is the Krull dimension \eqref{eq:Krull_dimension} of $R$.
\end{defn}

The geometric meaning of Krull dimension is revealed by the

\begin{thm}[Noether normalization theorem]
\label{thm:Noether_normalization}
(See Eisenbud \cite[Section 8.2.1, Theorem A1, p. 221]{Eisenbud_commutative_algebra} and G\"ortz and Wedhorn \cite[Theorem 5.15 and Corollary 5.17, pp. 127--128]{Gortz_Wedhorn_algebraic_geometry_v1}.)  
If $R$ is an affine ring over a field $\KK$ (that is, $R$ is a finitely generated $\KK$-algebra) and $\fp_0 \subsetneq \fp_1 \subsetneq \cdots \fp_{d-1}\subsetneq\fp_d$ is a a chain of prime ideals of $R$, maximal in the sense that no further prime ideals can be inserted into the chain, then there is a subring $S$ of $R$ with $S = \KK[x_1,\ldots,x_d]$ such that $R$ is a finitely generated $S$-module and and $\fp_i\cap S = (x_1,\ldots,x_i)$ for $i=0,\ldots,d$. 
\end{thm}

\begin{rmk}[Noether normalization and geometric meaning of Krull dimension]
\label{rmk:Noether_normalization_geometric_meaning_Krull_dimension}
(See G\"ortz and Wedhorn \cite[Remark 5.16, p. 127]{Gortz_Wedhorn_algebraic_geometry_v1}, and Mumford \cite[Corollary 2.29, p. 36]{Mumford_algebraic_geometry_I}, \cite[Chapter I, Section 7, p. 42]{Mumford_red_book_varieties_schemes}.)  
Theorem \ref{thm:Noether_normalization} has a well-known geometric interpretation. Suppose that $R = \KK[y_1,\ldots,y_n]/\fa$, where $\fa$ is a finitely generated ideal such that the coordinate ring $R$ has dimension $d$ and so $d$ is the dimension of the affine variety $X \subset \AAA_\KK^n$ over $\KK$ defined by $\fa$. The inclusion map $S \to R$ induces a surjective finite morphism of affine varieties $X \to \AAA_\KK^d$. Hence, any affine variety is a branched covering of affine space.   
\end{rmk}

\begin{thm}[Krull principal ideal theorem]
\label{thm:Krull_principal_ideal}
(See Eisenbud \cite[Section 8.2.1, Theorem B, p. 222]{Eisenbud_commutative_algebra}.)
If $R$ is a Noetherian local ring with maximal ideal $\fm$, then $\dim R$ is the minimal number $d$ such that there exist $d$ elements $f_1, \ldots, f_d \in \fm$ not all contained in any prime ideal other than $\fm$.
\end{thm}

One has the following alternative, well-known characterization of Krull dimension.

\begin{thm}[Krull dimension and degree of the Hilbert polynomial for a Noetherian local ring]
\label{thm:Hilbert_polynomial}
(See Eisenbud \cite[Section 8.2.1, Theorem C, p. 223]{Eisenbud_commutative_algebra}.)
Let $R$ be a Noetherian local ring with maximal ideal $\fm$, and denote the \emph{Hilbert function} by\footnote{The notation for the Hilbert function and polynomial are redefined from the statement in Eisenbud \cite[Section 8.2.1, Theorem C, p. 223]{Eisenbud_commutative_algebra} so the degree of $HP(n)$ matches that of $\chi_\fm(n)$ in Atiyah and Macdonald \cite[Chapter 11, pp. 118--119]{Atiyah_Macdonald_introduction_commutative_algebra}.}
\[
  HF(n) = \dim_{R/\fm}\fm^{n-1}/\fm^n, \quad\text{for all } n \in \NN.
\]  
For large $n$, $HF(n)$ agrees with a polynomial $HP(n)$, and
\[
  \dim R = \deg HP.
\]
\end{thm}

\begin{rmk}[Hilbert polynomial for finitely-generated module over a Noetherian local ring]
\label{rmk:Hilbert_polynomial_finitely-generated_module_over_Noetherian_local_ring}  
The definition of the Hilbert function and polynomial for a Noetherian local ring $R$ in the statement of Theorem \ref{thm:Hilbert_polynomial} can be generalized to the case of a finitely generated $R$-module $M$. See Atiyah and Macdonald \cite[Chapter 11, pp. 118--119]{Atiyah_Macdonald_introduction_commutative_algebra}, Gopalakrishnan \cite[Section 8.1, Definition, p. 194]{Gopalakrishnan_commutative_algebra}, Matsumura \cite[Chapter 5, Section 13, Theorem 13.2, p. 94 and Corollary, p. 95]{Matsumura_commutative_ring_theory}, and Raghavan, Balwant Singh, and Sridharan \cite[Section 3.2, p. 57]{Raghavan_BalwantSingh_Sridharan_homological_methods_commutative_algebra}. 
\end{rmk}

We recall the key

\begin{lem}[Dimension formula]
\label{lem:Atiyah_Macdonald_corollary_11-18}  
(See Atiyah and Macdonald \cite[Corollary 11.18, p. 122]{Atiyah_Macdonald_introduction_commutative_algebra}.)
Let $(R,\fm)$ be a Noetherian local ring. If $x \in \fm$ is an element that is not a zero-divisor, then
\[
  \dim R/(x) = \dim R - 1.
\]  
\end{lem}

If $\KK$ is a field, then it has exactly two ideals, $(0)$ and $\KK$, and only $(0)$ is prime, so $\dim\KK = 0$ according to the Definition \ref{defn:Krull_dimension} of Krull dimension. Of central importance in our applications is the following well-known

\begin{thm}[Krull dimension of a polynomial or formal power series ring over a Noetherian ring]
\label{thm:Krull_dimension_polynomial_power_series_rings}
(See Bruns and Herzog \cite[Theorem A.12 and Corollary A.13, p. 416]{Bruns_Herzog_cohen-macaulay-rings}, Gopalakrishnan \cite[Section 8.2, Corollary 9, p. 203]{Gopalakrishnan_commutative_algebra}, or Matsumura \cite[Theorem 15.4, p. 117]{Matsumura_commutative_ring_theory}.)  
If $R$ is a Noetherian ring, $x_1,\ldots,x_n$ are indeterminates over $R$, and $R[x_1,\ldots,x_n]$ and $R[[x_1,\ldots,x_n]]$ are the polynomial and formal power series rings with indeterminates $x_1,\ldots,x_n$ over $R$, then their Krull dimensions are given by
\begin{equation}
  \label{eq:Krull_dimension_polynomial_power_series_rings_over_Noetherian_ring}
  \dim R[x_1,\ldots,x_n] = \dim R[[x_1,\ldots,x_n]] = \dim R + n.
\end{equation}
In particular, if $R$ is replaced by a field $\KK$, then the Krull dimensions are given by
\begin{equation}
  \label{eq:Krull_dimension_polynomial_power_series_rings_over_field}
  \dim \KK[x_1,\ldots,x_n] = \dim \KK[[x_1,\ldots,x_n]] = n.
\end{equation}
\end{thm}

\begin{rmk}[Krull dimension and regularity of the rings of formal and convergent power series]
\label{rmk:Krull_dimension_polynomial_formal_convergent_power_series_ring}
For any field $\KK$ and integer $n\geq 1$, we noted in Remark \ref{rmk:Analytic_model_space_K-ringed_space} that the ring $\KK[[x_1,\ldots,x_n]]$ of formal power series is a local ring, with unique maximal ideal $\fm = (x_1,\ldots,x_n)$. Since
\begin{equation}
  \label{eq:Dimension_fm_mod_fm2_formal_convergent_power_series_ring_equals_n}
  \dim_\KK\fm/\fm^2 = n,
\end{equation}
then the equality \eqref{eq:Krull_dimension_polynomial_power_series_rings_over_field} in Theorem \ref{thm:Krull_dimension_polynomial_power_series_rings} implies that 
\[
  \dim \KK[[x_1,\ldots,x_n]] = \dim_\KK\fm/\fm^2,
\]
and so $\KK[[x_1,\ldots,x_n]]$ is a regular local ring by Definition \ref{defn:Regular_local_ring}.

If $\KK=\RR$ or $\CC$, then the Krull dimension of the ring $\KK\{x_1,\ldots,x_n\}$ of convergent power series is 
\begin{equation}
  \label{eq:Krull_dimension_polynomial_convergent_power_series_ring}
  \dim \KK\{x_1,\ldots,x_n\} = n.
\end{equation}
This is most easily seen from the proof by Gopalakrishnan of \cite[Section 8.2, Corollary 9, p. 203]{Gopalakrishnan_commutative_algebra}, replacing the role of formal power series rings with convergent power series rings. Because $\fm = (x_1,\ldots,x_n)$ is the maximal ideal of $\KK\{x_1,\ldots,x_n\}$, then $\dim \KK\{x_1,\ldots,x_n\} \leq n$ by \eqref{eq:Krull_dimension_leq_dimension_Zariski_cotangent_space} since the equality \eqref{eq:Dimension_fm_mod_fm2_formal_convergent_power_series_ring_equals_n} continues to hold. On the other hand, $\dim \KK\{x_1,\ldots,x_n\} \geq n$ by considering the chain of prime ideals
\[
  (x_1) \subsetneq (x_1,x_2) \subsetneq \cdots \subsetneq (x_1,\ldots,x_n).
\]
Combining these two inequalities proves the equality \eqref{eq:Krull_dimension_polynomial_convergent_power_series_ring}. We also noted in Remark \ref{rmk:Analytic_model_space_K-ringed_space} that $\KK\{x_1,\ldots,x_n\}$ is a local ring, with unique maximal ideal $\fm = (x_1,\ldots,x_n)$, and because
\[
  \dim \KK\{x_1,\ldots,x_n\} = \dim_\KK\fm/\fm^2,
\]
then $\KK\{x_1,\ldots,x_n\}$ is a regular local ring by Definition \ref{defn:Regular_local_ring}.
\end{rmk}  

\begin{thm}[Chevalley's Theorem]
\label{thm:Chevalley_theorem}
(See Atiyah and Macdonald\footnote{Atiyah and Macdonald assume that $R$ is Noetherian.} \cite[Chapter 11, Theorem 11.14, p. 121]{Atiyah_Macdonald_introduction_commutative_algebra}, Zariski and Samuel \cite[Chapter 8, Section 9, Theorem 20, p. 288]{Zariski_Samuel_communtative_algebra_II}, or Abhyankar \cite[Chapter 3, Section 17.2, p. 142]{Abhyankar_local_analytic_geometry}; compare the Stacks Project \cite[\href{https://stacks.math.columbia.edu/tag/00KQ}{Proposition 00KQ}]{stacks-project}.)  
If $(R,\fm)$ is a Noetherian local ring, then $\dim R$ is equal to the smallest integer $d$ for which there exist $d$ elements of $R$ that generate an ideal which is primary for $\fm$. 
\end{thm}

An ideal $\fa$ in a ring $R$ is called \emph{primary} if $\fa \neq R$ and if for any $x,y \in R$,
\[
  xy \in \fa \implies \text{either } x \in \fa \text{ or } y^n \in \fa \text{ for some } n > 0.
\]
(See Atiyah and Macdonald \cite[Chapter 4, pp. 50--51]{Atiyah_Macdonald_introduction_commutative_algebra} or Zariski and Samuel \cite[Chapter 3, Section 9, p. 152]{Zariski_Samuel_communtative_algebra_I}.) If $\fa \subset R$ is a primary ideal, then its radical $\fb = \sqrt{\fa}$ is called the \emph{associated prime ideal} of $\fa$ and one says that $\fa$ \emph{is primary for} $\fb$ or that $\fa$ is $\fb$-primary. The conclusion of Theorem \ref{thm:Chevalley_theorem} means that there there exist $x_1,\ldots,x_d \in R$ such that the ideal generated by these $d$ elements obeys $(x_1,\ldots,x_d) = \sqrt{\fm}$ (see Zariski and Samuel \cite[Chapter 8, Section 9, Remark (2), p. 291]{Zariski_Samuel_communtative_algebra_II}).

If $R$ is a local Noetherian ring with maximal ideal $\fm$, then an ideal $\fa\subset R$ such that $\fa = \sqrt{\fm}$ is also called an \emph{ideal of definition of $R$} (see the Stacks Project \cite[\href{https://stacks.math.columbia.edu/tag/07DU}{Definition 07DU}]{stacks-project}), as well as an $\fm$-primary ideal as above. If $x_1,\ldots,x_d \in \fm$ generate an $\fm$-primary ideal, where $d=\dim R$ (the Krull dimension \eqref{eq:Krull_dimension} of $R$), then $\{x_1,\ldots,x_d\}$ is called a \emph{system of parameters} of $R$ (see Atiyah and Macdonald \cite[Chapter 11, p. 122]{Atiyah_Macdonald_introduction_commutative_algebra} or Matsumura \cite[Chapter 5, Section 14, p. 104]{Matsumura_commutative_ring_theory}).

One has the following generalization of Lemma \ref{lem:Atiyah_Macdonald_corollary_11-18}.

\begin{lem}[Generalized dimension formula]
\label{lem:Stacks_Project_10-60-13}
(See Matsumura \cite[Theorem 14.1]{Matsumura_commutative_ring_theory} or the Stacks Project \cite[\href{https://stacks.math.columbia.edu/tag/02IE}{Lemma 02IE}]{stacks-project}.)
Let $(R,\fm)$ be a Noetherian local ring. Suppose $x_1,\ldots,x_d \in \fm$ generate an ideal of definition and $d=\dim R$. Then $\dim(R/(x_1,\ldots,x_i)) = d-i$ for all $i=1,\ldots,d$.
\end{lem}

\begin{rmk}[Artinian rings and dimension]
\label{rmk:Explanation_Grauert_Remmert_definition_dimension_complex_analytic_space}
Recall that a ring $R$ is called \emph{Artinian} (respectively, \emph{Noetherian}) if it satisfies the descending (respectively, ascending) chain condition on ideals (see Atiyah and Macdonald \cite[Chapter 6, p. 76]{Atiyah_Macdonald_introduction_commutative_algebra}). If $R$ is Noetherian, then so is a quotient $R/\fa$ for any ideal $\fa$ (see Atiyah and Macdonald \cite[Proposition 6.6, p. 76]{Atiyah_Macdonald_introduction_commutative_algebra}). According to the Stacks Project \cite[\href{https://stacks.math.columbia.edu/tag/00KJ}{Proposition 00KJ}]{stacks-project}, a ring is Artinian if and only if it is Noetherian and has dimension zero. Thus, if $R$ is a Noetherian ring, then the quotient $R/(x_1,\ldots,x_d)$ is Noetherian and so $\dim(R/(x_1,\ldots,x_d)) = 0$ if and only if $R/(x_1,\ldots,x_d)$ is Artinian.
\end{rmk}

\begin{rmk}[Intrinsic characterization of the minimal number of generators of an ideal in a local ring]
\label{rmk:Nakayama_lemma_minimal_number_generators_local_ring}  
As a consequence of Nakayama's Lemma (see Atiyah and Macdonald \cite[Proposition 2.6, p. 21]{Atiyah_Macdonald_introduction_commutative_algebra} or Matsumura \cite[Theorem 2.2, p. 8]{Matsumura_commutative_ring_theory}) and Matsumura \cite[Theorem 2.3, p. 8]{Matsumura_commutative_ring_theory}, one has the following well-known, useful characterization: The \emph{minimal number of generators} for a finitely generated ideal $I \subset R$ in a local ring $(R,\fm,\kappa)$ is equal to the dimension of $I/(\fm I)$ as a vector space over $\kappa$ (see De Jong and Pfister \cite[Corollary 1.3.6, p. 18]{DeJong_Pfister_local_analytic_geometry} or Pirnes \cite[Corollary 3.7, p. 23]{Pirnes_minimal_number_generators_ideals_commutative_rings}). 
\end{rmk}

As we discuss in the forthcoming Section \ref{sec:Dimension_analytic_space}, several other concepts of dimension are used in the theory of analytic spaces and but some may defined in more generality. For this purpose, we need the

\begin{defn}[Length of a module over a ring]
\label{defn:Length_module_over_ring}  
(See Atiyah and Macdonald \cite[Chapter 6, p. 76]{Atiyah_Macdonald_introduction_commutative_algebra} and the Stacks Project \cite[\href{https://stacks.math.columbia.edu/tag/02LY}{Definition 02LY}]{stacks-project}; compare Gopalakrishnan \cite[Section 3.4, Definition, p. 81]{Gopalakrishnan_commutative_algebra} and Raghavan, Balwant Singh, and Sridharan \cite[Section 1.5, p. 21]{Raghavan_BalwantSingh_Sridharan_homological_methods_commutative_algebra} for more restrictive definitions.)
Let $R$ be a ring. For any $R$-module $M$, one defines the \emph{length} of $M$ over $R$ by the formula
\[
  \length_R(M) = \sup\left\{n \in \ZZ: (0) = M_0 \subsetneq M_1 \subsetneq \cdots \subsetneq M_n = M,
    \ M_i \neq M_{i+1} \text{ for } i = 0,\ldots,n\right\}.
\]
\end{defn}



\begin{defn}[Chevalley dimension of a module over a Noetherian local ring]
\label{defn:Chevalley_dimension_module_over_Noetherian_local_ring}  
(See Ash \cite[Section 5.3.2]{Ash_course_commutative_algebra}, Gopalakrishnan \cite[Section 8.2, Definition, p. 200]{Gopalakrishnan_commutative_algebra} and Raghavan, Balwant Singh, and Sridharan \cite[Section 3.2, p. 57]{Raghavan_BalwantSingh_Sridharan_homological_methods_commutative_algebra}.) Let $(R,\fm)$ be a Noetherian local ring and $M$ be a finitely-generated $R$-module. If $M\neq (0)$, the \emph{Chevalley dimension} of $M$ is
\begin{equation}
  \label{eq:Chevalley_dimension_module_over_Noetherian_local_ring}
  s(M) := \inf\left\{ r \in \ZZ: \exists\ x_1,\ldots,x_r \in \fm \text{ such that } \length_R\left(M/M(x_1,\ldots,x_r)\right) < \infty \right\},
\end{equation}
where $(x_1,\ldots,x_r) \subset R$ is the ideal generated by $x_1,\ldots,x_r$ and if $M=(0)$, then $s(M) = -1$. IfIf $M=R$, then $s(R)$ is the Chevalley dimension of the ring $R$.
\end{defn}

\begin{thm}[Equivalence of Chevalley and Krull dimensions for finitely-generated modules over Noetherian local rings]
\label{thm:Gopalakrishnan_1_page_200}  
(See  Ash \cite[Section 5.3.3]{Ash_course_commutative_algebra}, Gopalakrishnan \cite[Section 8.2, Theorem 1, p. 200]{Gopalakrishnan_commutative_algebra} and Raghavan, Balwant Singh, and Sridharan \cite[Section 3.2, Theorem 3.10, p. 57]{Raghavan_BalwantSingh_Sridharan_homological_methods_commutative_algebra}.) If $(R,\fm)$ is a Noetherian local ring and $M$ be a finitely-generated $R$-module, then $\dim M = s(M)$ and, in particular, if $M=R$, then
\[
  \dim R = s(R).
\]  
\end{thm}

In the special case where the module is the Noetherian local ring itself, Theorem \ref{thm:Gopalakrishnan_1_page_200} can be obtained via De Jong and Pfister \cite[Remark 4.2.15, p. 145]{DeJong_Pfister_local_analytic_geometry}.

\begin{rmk}[Rings of formal or convergent power series are regular local rings]
\label{rmk:Example_formal_power_series_regular_local_ring}  
If $\KK$ is a field, then the ring $\KK[[x_1,\ldots,x_n]]$ of formal power series is a regular local ring of dimension $n$ by Remark \ref{rmk:Krull_dimension_polynomial_formal_convergent_power_series_ring} (or Zariski and Samuel \cite[Section 8.11, Example, p. 301]{Zariski_Samuel_communtative_algebra_II} or Abhyankar \cite[Chapter 2, Proposition 10.7, p. 80]{Abhyankar_local_analytic_geometry}). Similarly\footnote{Abhyankar allows $\KK$ to be any complete, non-discrete valued field \cite[Chapter 1, Section 2, p. 6]{Abhyankar_local_analytic_geometry}.}, for $\KK=\RR$ or $\CC$, the ring $\KK\{x_1,\ldots,x_n\}$ of convergent power series in the indeterminates $x_1,\ldots,x_n$ is a regular local ring of dimension $n$ by Remark \ref{rmk:Krull_dimension_polynomial_formal_convergent_power_series_ring} (or Abhyankar \cite[Chapter 1, Section 2, p. 7 and Chapter 2, Proposition 10.7, p. 80]{Abhyankar_local_analytic_geometry} or Ruiz \cite[Chapter II, Section 1, Lemma 1.9 and its proof, p. 20]{Ruiz_basic_theory_power_series}).
\end{rmk}  

\begin{rmk}[Maximal ideals in the ring of polynomials over a field]
\label{rmk:Example_polynomial_ring_not_regular_local_ring}
(See Conrad \cite[Theorem 3.1]{Conrad_maximal_ideal_polynomial_ring}.)
Let $\KK$ be a field and $\KK[x_1,\ldots,x_n]$ denote the ring of polynomials in $n$ indeterminates ($n\geq 1$). If $p = (p_1,\ldots,p_n) \in \KK^n$ is any point, then the kernel $\fm_p$ of the evaluation homomorphism of rings
\[
  \ev_p:\KK[x_1,\ldots,x_n] \ni f \mapsto f(p_1,\ldots,p_n) \in \KK
\]
is an ideal and is maximal since the quotient $\KK[x_1,\ldots,x_n]/\fm_p \cong \KK$ is a field. Thus, $\KK[x_1,\ldots,x_n]$ is not a local ring since it has more than one maximal ideal. Conversely, if $\KK$ is algebraically closed, then by the Hilbert Nullstellensatz (see the forthcoming Theorem \ref{thm:Hilbert_Nullstellensatz}) every maximal ideal in $\KK[x_1,\ldots,x_n]$ has the form $(x_1-p_1,\ldots,x_n-p_n)$ for some $p = (p_1,\ldots,p_n) \in \KK^n$. 
\end{rmk}

In order to use expected dimension to give a lower bound for dimension, we shall need the following estimate from commutative algebra of the Krull dimension (see Definition \ref{defn:Krull_dimension}) of a quotient ring. While Lemma \ref{lem:Krull_dimension_quotient_ring} is a corollary of Bruns and Herzog \cite[Proposition A.4, p. 414]{Bruns_Herzog_cohen-macaulay-rings} (by taking the finite $R$-module $M$ to be $R$ itself), but we include a more explicit proof given its importance.

\begin{lem}[Krull dimension of a quotient ring]
\label{lem:Krull_dimension_quotient_ring}
Suppose $R$ is a commutative Noetherian local ring with maximal ideal $\fm$. If $I \subset R$ is an ideal with generators $x_1,\ldots,x_l$, then the Krull dimension of $R/I$ obeys
\begin{equation}
  \label{eq:Krull_dimension_quotient_ring}
  \dim R/I  \geq \dim R - l.
\end{equation}
\end{lem}

\begin{proof}
The Krull dimension of a Noetherian local ring $R$ with maximum ideal $\fm$ is equal to the minimum number of generators of an $\fm$-primary ideal in $R$ by Atiyah and Macdonald \cite[Theorem 11.14]{Atiyah_Macdonald_introduction_commutative_algebra}. Thus, noting that $I \subset \fm$ since $\fm$ is maximal, suppose $\fb \subset R/I$ is an $(\fm/I)$-primary ideal with minimum number $k := \dim R/I$ of generators $\bar y_1,\ldots,\bar y_k$. (Recall that the property of $\fb$ being $(\fm/I)$-primary means that $\fm/I = \sqrt{\fb}$ by Atiyah and Macdonald \cite[pp. 50--51]{Atiyah_Macdonald_introduction_commutative_algebra}.) Let $y_1,\ldots,y_k \in R$ be lifts of the generators $\bar y_1,\ldots,\bar y_k \in R/I$. We claim that $\fa := (x_1,\ldots,x_l,y_1,\ldots,y_k)$ is an $\fm$-primary ideal in $R$. If $x \in \fa$, then we may write
\[
  x = a_1x_1 + \cdots + a_lx_l + b_1y_1 + \cdots + b_ky_k,
\]
where $a_i, b_j \in R$ for $i=1,\ldots,l$ and $j=1,\ldots,k$. By definition of $\sqrt{\fb}$ in Atiyah and Macdonald \cite[p. 8]{Atiyah_Macdonald_introduction_commutative_algebra}, there is a least positive integer $N$ such that $\bar y_j^N \in \fm/I$ for $j=1,\ldots,k$. Observe that $y_j^N \in \fm$ for $j=1,\ldots,k$. Since $I \subset \fm$, then $x_i \in \fm$ for $i=1,\ldots,l$ and thus $y_0 := a_1x_1 + \cdots + a_lx_l \in \fm$. Consequently, writing $b_0=1$ and applying the Multinomial Theorem (see Olver, Lozier, Boisvert, and Clark \cite[Section 26.4 (ii), Equation (26.4.9), p. 620]{Olver_Lozier_Boisvert_Clark})
\begin{align*}
  x^{(k+1)N}
  &= \left(b_0y_0 + b_1y_1 + \cdots + b_ky_k\right)^{(k+1)N}
  \\
  &= \sum_{n_0+n_1+\cdots+n_k = (k+1)N}\binom{(k+1)N}{n_0,n_1,\ldots,n_k}\prod_{j=0}^k b_j^{n_j}y_j^{n_j},
\end{align*}
where the sum is taken over all combinations of nonnegative integer indices $n_0$ through $n_k$ such that their sum is equal to $(k+1)N$. If $n_0=n_1=\cdots=n_k$, a choice that realizes the minimal possible power $n_j$ for some index $j$ subject to the constraint $n_0+\cdots+n_k = (k+1)N$, then $n_j = N$ for $j=0,1,\ldots,k$. Therefore, using $y_j^N \in \fm$ for $j=0,1,\ldots,k$, we obtain
\[
  x^{(k+1)N} \in \fm.
\]
Thus, $\sqrt{\fa} = \fm$ by definition of $\sqrt{\fa}$ in Atiyah and Macdonald \cite[p. 8]{Atiyah_Macdonald_introduction_commutative_algebra}. Hence, $\fa$ is an $\fm$-primary ideal in $R$ as claimed. The definition $k = \dim R/I$, the characterization of the Krull dimension $\dim R$ as the minimum number of generators of an $\fm$-primary ideal in $R$, and the fact that $\fa = (x_1,\ldots,x_l,y_1,\ldots,y_k)$ is an $\fm$-primary ideal in $R$ yield the following equality and inequality, 
\[
  l + \dim R/I = l + k \geq \dim R.
\]
This verifies inequality \eqref{eq:Krull_dimension_quotient_ring} and completes the proof of Lemma \ref{lem:Krull_dimension_quotient_ring}.
\end{proof}

\section{Smooth points and regular points for schemes}
\label{sec:Smooth_point_regular_point_schemes}
By definition of a \emph{scheme} (see Eisenbud and Harris \cite[Section 1.2, p. 21]{Eisenbud_Harris_geometry_schemes}, G\"ortz and Wedhorn \cite[Definition 3.1]{Gortz_Wedhorn_algebraic_geometry_v1}, Hartshorne \cite[Chapter II, Section 2, Definition, p. 74]{Hartshorne_algebraic_geometry}, Shafarevich \cite[Section 5.3.1, p. 28]{Shafarevich_v2}, and Vakil \cite[Definition 4.3.1]{Vakil_foundations_algebraic_geometry}), the stalk $\sO_{X,p}$ in Definition \ref{defn:Regular_point_scheme} is a local ring with unique maximal ideal $\fm_p$ and \emph{residue field} $\kappa(p) := \sO_{X,p}/\fm_p$ \cite[Section 2.9, p. 57 and Section B.1, pp. 554--555]{Gortz_Wedhorn_algebraic_geometry_v1}. If $\sO_{X,p}$ is a finitely generated $\KK$-algebra, then $\kappa(p)$ is a finite algebraic extension of $\KK$ by Atiyah and Macdonald \cite[Corollary 7.10]{Atiyah_Macdonald_introduction_commutative_algebra}; in particular, if $\KK$ is an algebraically closed field, then $\kappa(p) \cong \KK$. See also G\"ortz and Wedhorn \cite[Example 2.32]{Gortz_Wedhorn_algebraic_geometry_v1} and Shafarevich \cite[Section 5.1.2, p. 7]{Shafarevich_v2}. Recall that if $\KK$ is an arbitrary field and $X$ is a scheme over $\KK$ that is locally of finite type (see Section \ref{sec:Schemes}), then $p\in X$ is a $\KK$-\emph{rational point} if $\kappa(p)\cong \KK$ (see \cite[Section 5.1, p. 121]{Gortz_Wedhorn_algebraic_geometry_v1}). If a locally ringed space $(X,\sO_X)$ is further restricted to be a \emph{$\KK$-ringed space} then, by definition, $\kappa(p)$ is isomorphic to $\KK$ as a $\KK$-algebra (see Fischer \cite[Section 0.1, p.1]{Fischer_complex_analytic_geometry}, Grauert and Remmert \cite[Section 1.1.3, p. 5]{Grauert_Remmert_coherent_analytic_sheaves}, or Guaraldo, Macr\`\i, and Tancredi \cite[Section 1.1, Definition 1.1, p. 1]{Guaraldo_Macri_Tancredi_topics_real_analytic_spaces}).

\begin{defn}[Regular point of a scheme]
\label{defn:Regular_point_scheme}
(See G\"ortz and Wedhorn \cite[Definition 6.24, p. 161]{Gortz_Wedhorn_algebraic_geometry_v1}, Hartshorne \cite[Chapter II, Section 6, Remark 6.11.1A, p. 142]{Hartshorne_algebraic_geometry}, and Vakil \cite[Definition 12.2.3]{Vakil_foundations_algebraic_geometry}.)  
Let $(X,\sO_X)$ be a locally Noetherian scheme. A point $p \in X$ is \emph{regular} if the local ring $\sO_{X,p}$ is regular and the scheme $X$ is \emph{regular} if it is regular at every point. 
\end{defn}

We next recall the 

\begin{defn}[Smooth point of a scheme]
\label{defn:Smooth_point_scheme}
(See G\"ortz and Wedhorn \cite[Definition 6.14, p. 156]{Gortz_Wedhorn_algebraic_geometry_v1}, Hartshorne \cite[Chapter III, Section 10, Definition, p. 268]{Hartshorne_algebraic_geometry}, and Vakil \cite[Definition 12.2.6]{Vakil_foundations_algebraic_geometry}.)  
Let $X$ be a scheme over a field $\KK$ and $d\geq 0$ be an integer. Then $X$ is \emph{smooth (of relative dimension $d$ over $\KK$) at a point $p \in X$} if it has an open affine neighborhood of the form
\[
  U = \Spec\KK[x_1,\ldots,x_n]/(f_1,\ldots,f_{n-d})
\]
such that
\begin{equation}
  \label{eq:Jacobi_criterion_scheme}
  \Corank_{\kappa(p)} J_{f_1,\ldots,f_{n-d}}(p) = d,
\end{equation}
where $\kappa(p)$ is the residue field and the \emph{corank} of an $(n-d)\times n$ matrix $M$ is defined by $\Corank M := n - \Rank M$, and
\begin{equation}
  \label{eq:Jacobian_matrix_scheme}
  J_{f_1,\ldots,f_{n-d}}(p) := \left(\frac{\partial f_i}{\partial x_j}(p)\right) \in \Hom_\KK(\KK^n,\KK^{n-d})
\end{equation}
is the \emph{Jacobian matrix} defined by $f_1,\ldots,f_{n-d}$ at the point $p$. One says that $X$ is \emph{smooth (of relative dimension $d$ over $\KK$)} if $X$ is smooth (of relative dimension $d$ over $\KK$) at every point $p \in X$.
\end{defn}

Definition \ref{defn:Regular_point_scheme} is a specialization of the indicated definition in G\"ortz and Wedhorn obtained by taking $Y=\Spec\KK$ (see G\"ortz and Wedhorn \cite[Section 6.8, p. 157]{Gortz_Wedhorn_algebraic_geometry_v1}). The concepts of smooth point and regular point of a scheme are compared in the following

\begin{thm}[Comparison of smoothness and regularity for schemes]
\label{thm:Gortz_Wedhorn_6-28}
(See G\"ortz and Wedhorn \cite[Theorem 6.28]{Gortz_Wedhorn_algebraic_geometry_v1} and Vakil \cite[Theorem 12.2.10]{Vakil_foundations_algebraic_geometry}.)  
Let $X$ be a scheme locally of finite type over a field $\KK$, and $p \in X$ be a closed point\footnote{A point $p$ in a topological space is a \emph{closed point} if $\{p\}\subset X$ is a closed subset \cite[Definition 3.6.8]{Vakil_foundations_algebraic_geometry}.}, and $d \geq 0$ be an integer. Then the following are equivalent:
\begin{enumerate}
\item\label{item:X_smooth_scheme_at_point_dimension_d} If $X$ is smooth of relative dimension $d$ at $p$.
\item The following equalities hold:
  \begin{equation}
    \label{eq:Dimension_Zariski_tangent_space_equals_dimension_local_ring_scheme}
  \dim_{\kappa(p)}T_pX = \dim\sO_{X,p} = d,
\end{equation}
where $T_pX$ is the Zariski tangent space to $X$ at the point $p$.
\end{enumerate}
If either one (and thus both) of the preceding conditions are satisfied, then
\begin{enumerate}
\setcounter{enumi}{2}    
\item\label{item:X_regular_scheme_at_point_dimension_d} The local ring $\sO_{X,p}$ is regular and has dimension $d$.
\end{enumerate}  
Furthermore, if
\begin{inparaenum}[\itshape a\upshape)]
\item $\kappa(p)=\KK$, or
\item $\KK$ is perfect,
\end{inparaenum}  
then the final condition implies the other ones.
\end{thm}

For the definitions of the \emph{Zariski cotangent space} and \emph{Zariski tangent space}, respectively,
\begin{subequations}
  \label{eq:Zariski_cotangent_and_tangent_space}
  \begin{align}
    \label{eq:Zariski_cotangent_space}
    T_p^*X &:= \fm_p/\fm_p^2,
    \\
    \label{eq:Zariski_tangent_space}
    T_pX &:= \Hom_{\kappa(p)}(\fm_p/\fm_p^2,\kappa(p)),
  \end{align}
\end{subequations}
to a scheme $(X,\sO_X)$ at a point $p$ and equivalent definitions (especially for $\KK$-ringed spaces, using derivations as in \cite{Fischer_complex_analytic_geometry} and \cite{Guaraldo_Macri_Tancredi_topics_real_analytic_spaces}), we refer to Fischer \cite[Section 2.1, pp. 77--78]{Fischer_complex_analytic_geometry}, Hartshorne \cite[Chapter II, Section 2, Exercise 2.8, p. 80]{Hartshorne_algebraic_geometry}, G\"ortz and Wedhorn \cite[Definition 6.2]{Gortz_Wedhorn_algebraic_geometry_v1}, Guaraldo, Macr\`\i, and Tancredi \cite[Definition 2.1, pp. 18--19]{Guaraldo_Macri_Tancredi_topics_real_analytic_spaces}, Shafarevich \cite[Section 5.1.2, p. 9]{Shafarevich_v2}, the Stacks Project \cite[\href{https://stacks.math.columbia.edu/tag/0B2C}{Definition 0B2C}]{stacks-project}, and Vakil \cite[Definition 12.1.1]{Vakil_foundations_algebraic_geometry}.

\begin{rmk}[Distinction between the regular point and smooth point for a scheme over an imperfect field]
\label{rmk:Vakil_12-2-11}
G\"ortz and Wedhorn \cite[Corollary 6.32, Remark 6.33, and Example 6.34]{Gortz_Wedhorn_algebraic_geometry_v1} and Vakil \cite[Section 12.2.11]{Vakil_foundations_algebraic_geometry} underline the distinction between the concepts of `regularity' and `smoothness'. Vakil describes an example where the field $\KK$ is not perfect and regularity does not imply smoothness. Recall that a field of characteristic zero is perfect (see Lang \cite[Section 5.6, p. 252]{LangAlgebra}).
\end{rmk}

If $X$ is a topological space, then its \emph{topological dimension} $\dim X$ is the supremum of all lengths $l$ of chains $X_0 \supsetneq X_1 \subsetneq \cdots \supsetneq X_l$ of irreducible closed subsets of $X$ (see G\"ortz and Wedhorn \cite[Definition 5.5]{Gortz_Wedhorn_algebraic_geometry_v1}, the Stacks Project \cite[\href{https://stacks.math.columbia.edu/tag/0055}{Definition 0055}]{stacks-project},and Vakil \cite[Definition 11.1.1]{Vakil_foundations_algebraic_geometry}). (A nonempty topological space X is called \emph{irreducible} if $X$ cannot be expressed as the union of two proper closed subsets; a nonempty subset $Z$ of $X$ is called \emph{irreducible} if $Z$ is irreducible when we endow it with the induced topology \cite[Definition 1.14]{Gortz_Wedhorn_algebraic_geometry_v1}.) We recall the important

\begin{thm}[Equality of topological dimension of a scheme and dimension of the local ring at a closed point]
\label{thm:Gortz_Wedhorn_5-22}
(See G\"ortz and Wedhorn \cite[Theorem 5.22]{Gortz_Wedhorn_algebraic_geometry_v1}.)
Let $X$ be an irreducible scheme over a field $\KK$ that is locally of finite type. If $p \in X$ is a closed point, then $\dim\sO_{X,p} = \dim X$.
\end{thm}

Theorems \ref{thm:Gortz_Wedhorn_6-28} and \ref{thm:Gortz_Wedhorn_5-22} yield the

\begin{cor}[Smoothness and dimensionality]
\label{cor:Gortz_Wedhorn_6-29}
(See G\"ortz and Wedhorn \cite[Corollary 6.29]{Gortz_Wedhorn_algebraic_geometry_v1}.)
Let $X$ be an irreducible scheme locally of finite type over a field $\KK$ and let $p \in X$ be a $\KK$-rational point\footnote{By \cite[Proposition 3.33]{Gortz_Wedhorn_algebraic_geometry_v1}, a point $p \in X$ is closed if and only if the residue field $\kappa(p)$ is a finite extension of $\KK$.} of $X$). Then $X$ is smooth over $\KK$ at $p$ if and only if $\dim_\KK T_pX = \dim X$.
\end{cor}

\section{Gluing schemes}
\label{sec:Gluing_schemes}
As we noted for the category of analytic spaces (see Section \ref{sec:Gluing_analytic_spaces}), one can construct schemes by an analogous method of gluing. We follow the outline provided by Hartshorne \cite[Chapter II, Section 2, Exercise 2.12, p. 80]{Hartshorne_algebraic_geometry}.

Let $\{X_i\}_{i\in I}$ be a (possibly infinite) family of schemes. For each $i\neq j$, let $U_{ij} \subset X_i$ be an open subset and let it have the induced scheme structure (see Hartshorne \cite[Chapter II, Section 2, Exercise 2.2, p. 79]{Hartshorne_algebraic_geometry}). Suppose also that for each $i\neq j$, we are given an isomorphism of schemes
\[
  \theta_{ij}: U_{ij} \to U_{ji}
\]
such that the following conditions are obeyed:
\begin{enumerate}
\item $\theta_{ji} = \theta_{ij}^{-1}$, for all $i,j \in I$, and
\item $\theta_{ij}(U_{ij}\cap U_{ik}) = U_{ji}\cap U_{jk}$ and $\theta_{ik} = \theta_{jk}\circ\theta_{ij}$ on $U_{ij}\cap U_{ik}$, for all $i,j,k \in I$.
\end{enumerate}
Then there are a scheme $X$ and morphisms $\psi_i:X_i\to X$, for each $i\in I$, such that the following hold: 
\begin{enumerate}
\item $\psi_i$ is an isomorphism from $X_i$ onto an open subscheme of $X$, for each $i\in I$,
\item The images $\psi_i(X_i)$, for all $i\in I$, cover $X$
\item $\psi_i(U_{ij}) = \psi_i(X_i)\cap \psi_j(X_j)$, for all $i,j\in I$,
\item $\psi_i = \psi_j\circ\theta_{ij}$ on $U_{ij}$, for all $i,j\in I$.  
\end{enumerate}  
One says that $X$ is obtained by gluing the schemes along the isomorphisms $\theta_{ij}$.

\chapter{Smooth points, regular points, and dimensions for analytic spaces}
\label{chap:Smooth_point_regular_point_dimension_analytic_space}
Our discussion in this chapter of smooth points, regular points, and dimensions for analytic spaces partially mirrors our discussion of the corresponding topics for schemes in Chapter \ref{chap:Smooth_point_regular_point_dimension_scheme}. However, unlike in the category of schemes, many of the corresponding results in the category of analytic spaces are difficult to find in the literature and due to their central importance in our work, we include proofs as needed. Some simplifications are enabled by the fact that analytic spaces are $\KK$-ringed spaces, so the residue field is equal to $\KK$ at every point in an analytic space. Section \ref{sec:Dimension_analytic_space} contains a development of dimension theory for analytic spaces. In Section \ref{sec:Embedding_dimension_analytic_space}, we review the concept of embedding dimension for a point in an analytic space. In Section \ref{sec:Smooth_point_regular_point_analytic_space}, we discuss the concepts of smooth points and regular points for analytic spaces. Section \ref{sec:Analyticity_dimension_singular_set} contains an exposition of results on analyticity and dimension of the singular set of an analytic space. In Section \ref{sec:Algebraic_analytic_expected_dimensions}, we develop the relationship between the Krull dimension, analytic dimension, and expected dimension at a point in an analytic space. We conclude in Section \ref{sec:Real_analytic_sets_counterexamples} with a review of how real analytic sets differ from complex analytic sets in important ways and recall well-known counterexamples that highlight some those differences.

\section{Dimensions for analytic spaces}
\label{sec:Dimension_analytic_space}
The \emph{Krull dimension} $\dim\sO_{X,p}$ of the local ring at a point $p$ of an analytic space $(X,\sO_X)$ over $\KK=\RR$ or $\CC$ is as in Definition \ref{defn:Krull_dimension}. If $U \subset \KK^n$ is an open set and $X \subset U $ is an analytic set, Abhyankar \cite[Chapter 5, Section 29.4, pp. 232--233]{Abhyankar_local_analytic_geometry} takes $\dim\sO_{X,p}$ as his definition of the dimension of $X$ at a point $p\in X$. In order to discuss the geometric significance of the Krull dimension for analytic spaces, we shall require some definitions. 

\begin{defn}[Irreducible germs and analytic sets that are irreducible at a point]
\label{defn:Irreducible_germ_analytic_set}  
Let $\KK=\RR$ or $\CC$, and $D \subset \KK^n$ be a domain, $S \subset D$ be an analytic set in the sense of Definition \ref{defn:Analytic_set}, and $p \in S$ be a point. If $S_p$ is the germ of $S$ at $p$, then $\sI(S_p) \subset \sO_{D,p}$ denotes the ideal of (germs of) analytic functions in $\sO_{D,p} = \sO_{\KK^n,p}$ which vanish on $S_p$ (see Narasimhan \cite[Section 3.1, p. 31]{Narasimhan_introduction_theory_analytic_spaces}). An analytic germ $S_p$ is \emph{irreducible} if whenever there are two analytic germs $S_{1p}, S_{2p}$ with $S_p = S_{1p}\cup S_{2p}$, then one must have $S_{ip} = S_p$ for $i=1$ or $2$ (see Narasimhan \cite[Section 3.1, p. 31]{Narasimhan_introduction_theory_analytic_spaces}). An analytic set $S$ is \emph{irreducible at a point $p\in S$} if the stalk $\iota(S)_p$ is a prime ideal in $\sO_{\KK^n,p}$ (see Grauert and Remmert \cite[Section 4.1.3, p. 78]{Grauert_Remmert_coherent_analytic_sheaves}), noting that $\iota(S)_p = \sI(S_p)$, where $\iota(S) \subset \sO_D$ is the ideal sheaf of $S$ as in Definition \ref{defn:Ideal_sheaf_for_analytic_subset}.
\end{defn}  

According to Narasimhan \cite[Section 3.1, Lemma 1, p. 31]{Narasimhan_introduction_theory_analytic_spaces}, an analytic germ $S_p$ is irreducible if and only if $\sI(S_p)$ is a prime ideal and so the definitions due to Grauert and Remmert (for $\KK=\CC$) and due to Narasimhan in Definition \ref{defn:Irreducible_germ_analytic_set} are equivalent. Furthermore, according to Narasimhan \cite[Section 3.1, Proposition 1, p. 32]{Narasimhan_introduction_theory_analytic_spaces} (or Grauert and Remmert \cite[Section 4.1.3, Local Decomposition Lemma, p. 79]{Grauert_Remmert_coherent_analytic_sheaves}, when $\KK=\CC$), an analytic germ $S_p$ can be expressed as a finite union of irreducible analytic germs (the \emph{irreducible components} of $S_p$), uniquely up to ordering.

To complete our discussion of the geometric significance of the Krull dimension, we need to review a version of the Noether normalization theorem for convergent power series rings. (See Theorem \ref{thm:Noether_normalization} for the algebraic version for polynomial rings.) We begin with the

\begin{defn}[Analytic algebras and their morphisms]
\label{defn:Analytic_algebra}
(See Grauert and Remmert \cite[Section 2.0.1, p. 77]{Grauert_Remmert_analytic_local_algebras} for a field $\KK$ as below or De Jong and Pfister \cite[Definition 3.2.8 (1), p. 89]{DeJong_Pfister_local_analytic_geometry} and Ebelin \cite[Section 2.6, Definition, p. 81]{Ebeling_functions_several_complex_variables_singularities} for $\KK=\CC$.)  
Let $\KK$ be a complete, valued field with infinitely many elements and let $\KK\{x_1,\ldots,x_n\}$ denote the ring of convergent power series (see Grauert and Remmert \cite[Section 1.3.1, p. 27]{Grauert_Remmert_analytic_local_algebras}) in the indeterminates $x_1,\ldots,x_n$. If $\fa \subsetneq \KK\{x_1,\ldots,x_n\}$ is an ideal, then $A := \KK\{x_1,\ldots,x_n\}/\fa$ is called an \emph{analytic algebra} (or \emph{analytic $\KK$-algebra}). A $\KK$-algebra homomorphism $\varphi:A\to B$ between two analytic $\KK$-algebras is called an \emph{analytic homomorphism}. An analytic homomorphism $\varphi:A\to B$ of analytic $\KK$-algebras gives $B$ the structure of an $A$-module and $\varphi$ is called a \emph{finite analytic homomorphism}\footnote{See also Narasimhan \cite[Chapter 2, p. 10]{Narasimhan_introduction_theory_analytic_spaces} or the Stacks Project \cite[\href{https://stacks.math.columbia.edu/tag/01WG}{Section 01WG}]{stacks-project}.} if $B$ is a finitely generated $A$-module (see Grauert and Remmert \cite[Section 2.2.2, p. 89]{Grauert_Remmert_analytic_local_algebras}).
\end{defn}

Grauert and Remmert observe \cite[Section 2.0.1, p. 78]{Grauert_Remmert_analytic_local_algebras} that analytic algebras form a category, with morphisms given by analytic homomorphisms.

\begin{thm}[Analytic algebras are Noetherian local rings]
\label{thm:Analytic_algebras_are_Noetherian_local_rings}
(See Grauert and Remmert \cite[Section 2.0.1, Satz 1, p. 77]{Grauert_Remmert_analytic_local_algebras} for $\KK$ as in Definition \ref{defn:Analytic_algebra} and Ebelin \cite[Section 2.6, Proposition 2.29, p. 81]{Ebeling_functions_several_complex_variables_singularities} for $\KK=\CC$.)  
Each analytic algebra\footnote{In the sense of Definition \ref{defn:Analytic_algebra}.} $A$ is a Noetherian local ring with residue field $A/\fm \cong \KK$, where $\fm \subset A$ is the maximal ideal.
\end{thm}

We have the following analogue of Theorem \ref{thm:Noether_normalization}, with the polynomial ring being replaced by the ring of convergent power series. 

\begin{thm}[Noether normalization theorem for analytic algebras]
\label{thm:Noether_normalization_analytic_algebras} 
(See Grauert and Remmert \cite[Section 2.2.2, p. 90]{Grauert_Remmert_analytic_local_algebras} and De Jong and Pfister \cite[Theorem 2.2.9, p. 61]{DeJong_Pfister_local_analytic_geometry} for $\KK$ as in Definition \ref{defn:Analytic_algebra}, Narasimhan \cite[Section 3.1, Proposition 2, p. 32]{Narasimhan_introduction_theory_analytic_spaces} for $\KK=\RR$ or $\CC$, and Ebelin \cite[Section 2.10, Proposition 2.44, p. 101]{Ebeling_functions_several_complex_variables_singularities} for $\KK=\CC$.)
Let $A$ be an analytic algebra as in Definition \ref{defn:Analytic_algebra}. Then there are a non-negative integer $d$ and an analytic, finite, injective homomorphism of analytic $\KK$-algebras,
\[
  \varphi:\KK\{x_1,\ldots,x_d\} \to A.
\]
\end{thm}

\begin{rmk}[Related versions of Theorem \ref{thm:Noether_normalization_analytic_algebras}]
\label{rmk:Proof_Noether_normalization_analytic_algebras}
Atiyah and Macdonald state \cite[Chapter 11, Exercise 2, p. 125]{Atiyah_Macdonald_introduction_commutative_algebra} that if $A$ is a \emph{complete} Noetherian local ring with system of parameters $x_1,\ldots,x_d$, then the homomorphism $\KK[[t_1,\ldots,t_d]]$ given by $t_i\mapsto x_i$ for $1\leq i\leq d$ is injective and that $A$ is a finitely generated module over $\KK[[t_1,\ldots,t_d]]$. This statement is close, though not identical to Theorem \ref{thm:Noether_normalization_analytic_algebras}. Nagata proves a version of Theorem \ref{thm:Noether_normalization_analytic_algebras} as \cite[Chapter 7, Theorem 45.5, p. 194]{Nagata_local_rings} for the ring $\KK\{t_1,\ldots,t_d\}$ and a field $\KK$ with a multiplicative valuation \cite[Chapter 7, Section 45, p. 190]{Nagata_local_rings}, but does not prove or state that $d = \dim A$. For $\KK=\CC$, both Fischer \cite[Section 3.1, Theorem 1, p. 131]{Fischer_complex_analytic_geometry} and Huybrechts \cite[Theorem 1.1.30, p. 20]{Huybrechts_2005} state variants of Theorem \ref{thm:Noether_normalization_analytic_algebras}, though without proof. Remark \ref{rmk:Noether_normalization_geometric_meaning_Krull_dimension} on the geometric meaning of Theorem \ref{thm:Noether_normalization} applies here too for Theorem \ref{thm:Noether_normalization_analytic_algebras}.
\end{rmk}

\begin{defn}[Weierstrass dimension of an analytic algebra]
\label{defn:Weierstrass_dimension_analytic_algebra}
Let $A$ be an analytic algebra as in Definition \ref{defn:Analytic_algebra}. Then the \emph{Weierstrass dimension} of $A$ is the least integer $d$ such that there exists a Noether normalization homomorphism $\KK\{x_1,\ldots,x_d\} \to A$ in the sense of Theorem \ref{thm:Noether_normalization_analytic_algebras}.
\end{defn}

If $(X,p)$ is a germ of a complex analytic space $(X,\sO_X)$ with local ring $(\sO_{X,p},\fm_{X,p})$, then the \emph{Weierstrass dimension} of $(X,p)$ is the least integer $d$ such that there exists a Noether normalization homomorphism $\CC\{x_1,\ldots,x_d\} \to \sO_{X,p}$ in the sense of Theorem \ref{thm:Characterization_integer_Noether_normalization_analytic_algebras}, where $\CC\{x_1,\ldots,x_d\}$ is the ring of convergent power series in the indeterminates $x_1,\ldots,x_d$ with coefficients in $\CC$. See De Jong and Pfister \cite[Corollary 3.3.19, p. 101, and Definition 4.1.3, p. 128]{DeJong_Pfister_local_analytic_geometry} or Ebeling \cite[Section 2.10, Proposition 2.44, p. 101 and Definition, p. 102]{Ebeling_functions_several_complex_variables_singularities}. 

Grauert and Remmert remark \cite[Section 2.2.2, p. 90]{Grauert_Remmert_analytic_local_algebras} that the integer $d$ in Theorem \ref{thm:Noether_normalization_analytic_algebras} is equal to the dimension of the analytic $\KK$-algebra $A$, by which they mean its Chevalley dimension. We outline their justification of this assertion. We begin with the following specialization of Definition \ref{defn:Chevalley_dimension_module_over_Noetherian_local_ring}.

\begin{defn}[Chevalley dimension of a module over an analytic algebra]
\label{defn:Chevalley_dimension_module_over_analytic_algebra}  
(See Grauert and Remmert \cite[Section II.4.3, p. 110]{Grauert_Remmert_analytic_local_algebras} for $\KK$ as in Definition \ref{defn:Analytic_algebra} or De Jong and Pfister \cite[Definition 4.1.3, pp. 128--129]{DeJong_Pfister_local_analytic_geometry} and Ebelin \cite[Proposition 2.45, p. 102]{Ebeling_functions_several_complex_variables_singularities} for $\KK=\CC$.) Let $A$ be an analytic algebra as in Definition \ref{defn:Analytic_algebra} and $M$ be a non-zero $A$-module. The \emph{Chevalley dimension} of $M$ is
\begin{equation}
  \label{eq:Chevalley_dimension_module_over_analytic_algebra}
  s(M) := \inf\left\{l \in \ZZ: \exists\ f_1,\ldots,f_l \in A \text{ such that } \dim_\KK\left(M/M(f_1,\ldots,f_l)\right) < \infty \right\},
\end{equation}
where $(f_1,\ldots,f_l) \subset A$ is the ideal generated by $f_1,\ldots,f_l$ and $\dim_\KK(\cdot)$ computes the dimension of $M/M(f_1,\ldots,f_l)$ as a vector space over $\KK$. One calls $f_1,\ldots,f_d$ a \emph{system of parameters} of $M$. If $M=A$, then $s(A)$ is the Chevalley dimension of the ring $A$.
\end{defn}

In \cite[Section 5.1.2, p. 95]{Grauert_Remmert_coherent_analytic_sheaves}, Grauert and Remmert define the Chevalley dimension of $\sO_{X,p}$ of the stalk of a complex analytic space $(X,\sO_X)$ at a point $p \in X$ to be the minimal number $d$ of germs $f_{1,p},\ldots,f_{d,p} \in \fm_{X,p}$ such that $\sO_{X,p}/(f_{1,p},\ldots,f_{d,p})\sO_{X,p}$ is finite-dimensional as a complex vector space. The integer $d$ in Theorem \ref{thm:Noether_normalization_analytic_algebras} can be characterized by the following

\begin{thm}[Equivalence of the Chevalley and Weierstrass dimensions of an analytic algebra]
\label{thm:Characterization_integer_Noether_normalization_analytic_algebras} 
(See Grauert and Remmert \cite[Section 2.5.2, Satz 4, p. 120]{Grauert_Remmert_analytic_local_algebras} for $\KK$ as in Definition \ref{defn:Analytic_algebra} and Ebelin \cite[Proposition 2.45, p. 102]{Ebeling_functions_several_complex_variables_singularities} for $\KK=\CC$.)
Let $A$ be an analytic algebra as in Definition \ref{defn:Analytic_algebra}. If $\fm \subset A$ is the maximal ideal, then $f_1,\ldots,f_d \in \fm$ is a system of parameters for $A$ (and $d$ is the Chevalley dimension of $A$) if the map
\[
  \varphi:\KK\{x_1,\ldots,x_d\} \to A,
\]
given by $x_i \mapsto f_i$ for $i=1,\ldots,d$ defines an analytic homomorphism that is finite and injective.
\end{thm}

We conclude by recalling the relationship between the Weierstrass dimension of an analytic algebra in Definition \ref{defn:Weierstrass_dimension_analytic_algebra}, its Chevalley dimension in Definition \ref{defn:Chevalley_dimension_module_over_analytic_algebra}, and its Krull dimension in Definition \ref{defn:Krull_dimension}. We first have the

\begin{lem}[Chevalley dimension as upper bound for the length of chain of prime ideals]
\label{lem:Chevalley_dimension_upper_bound_length_chain_prime_ideals} 
(See Grauert and Remmert \cite[Section 2.6.1, Folgerung, p. 127]{Grauert_Remmert_analytic_local_algebras}.)
Let $A$ be an analytic algebra as in Definition \ref{defn:Analytic_algebra}. Every chain of prime ideals $\fp_0 \subsetneq \fp_1 \subsetneq \cdots \subsetneq \fp_l \subset A$ has length $l \leq s(A)$, where $s(A)$ is the Chevalley dimension \eqref{eq:Chevalley_dimension_module_over_analytic_algebra} of $A$. 
\end{lem}

\begin{lem}[Chevalley dimension as lower bound for the length of chain of prime ideals]
\label{lem:Chevalley_dimension_lower_bound_length_chain_prime_ideals} 
(See Grauert and Remmert \cite[Section 2.6.2, Korollar, p. 130]{Grauert_Remmert_analytic_local_algebras} or De Jong and Pfister \cite[Lemma 4.1.7, p. 131]{DeJong_Pfister_local_analytic_geometry}.)
Let $A$ be an analytic algebra as in Definition \ref{defn:Analytic_algebra}. Then there exists a chain of prime ideals $\fp_0 \subsetneq \fp_1 \subsetneq \cdots \subsetneq \fp_l \subset A$ of length $l \geq s(A)$, where $s(A)$ is the Chevalley dimension \eqref{eq:Chevalley_dimension_module_over_analytic_algebra} of $A$. 
\end{lem}

Consequently, Lemmas \ref{lem:Chevalley_dimension_upper_bound_length_chain_prime_ideals} and \ref{lem:Chevalley_dimension_lower_bound_length_chain_prime_ideals} yield the

\begin{cor}[Equivalence of the Chevalley, Krull, and Weierstrass dimensions of an analytic algebra]
\label{cor:Chevalley_equal_to_Krull_dimension_analytic_algebra} 
(See Grauert and Remmert \cite[Section 2.6.2, Bemerkung, p. 130]{Grauert_Remmert_analytic_local_algebras} for $\KK$ as in Definition \ref{defn:Analytic_algebra} or De Jong and Pfister \cite[Theorem 4.1.9, p. 131]{DeJong_Pfister_local_analytic_geometry} for $\KK=\CC$.)
If $A$ is an analytic algebra as in Definition \ref{defn:Analytic_algebra}, then its Chevalley dimension \eqref{eq:Chevalley_dimension_module_over_analytic_algebra}, Krull dimension \eqref{eq:Krull_dimension}, and Weierstrass dimension are equal:
\begin{equation}
  \label{eq:2}
  s(A) = d = \dim A,
\end{equation}
where $d$ is the integer in Definition \ref{defn:Weierstrass_dimension_analytic_algebra}.
\end{cor}

Grauert and Remmert use the following result to motivate their definition of the \emph{analytic dimension} of a complex analytic space at a point.

\begin{lem}
\label{lem:Equivalence_analytic_algebraic_dimensions_complex_analytic_space}  
(See Grauert and Remmert \cite[Section 5.1.2, Lemma, p. 95]{Grauert_Remmert_coherent_analytic_sheaves}.)
Let $(X,\sO_X)$ be a complex analytic space and $p\in X$ be a point. The following integers are equal:
\begin{enumerate}
\item\label{item:Analytic_dimension_point_complex_space} $k\geq 0$ is the smallest integer such that there is an open neighborhood $U\subset X$ of $p$ and $f_1,\ldots,f_k \in \sO_X(U)$ such that $f_1^{-1}(0)\cap\cdots\cap f_k^{-1}(0) = \{p\}$.
\item\label{item:Algebraic_dimension_point_complex_space} $l\geq 0$ is the smallest integer such that there are $g_1,\ldots,g_l \in \fm_p$ with $\sO_{X,p}/(g_1,\ldots,g_l)$ an Artinian ring (see Remark \ref{rmk:Explanation_Grauert_Remmert_definition_dimension_complex_analytic_space}).
\end{enumerate}
\end{lem}

\begin{rmk}[Equivalence of analytic and algebraic dimensions of a complex analytic space at a point]
\label{rmk:Equivalence_analytic_algebraic_dimensions_complex_analytic_space}
In Grauert and Remmert \cite[Section 5.1.1, first paragraph and Definition, p. 93 and Section 5.1.2, first paragraph, p. 95 and Dimension Formula, p. 96]{Grauert_Remmert_coherent_analytic_sheaves}, the integer $k$ in Lemma \ref{lem:Equivalence_analytic_algebraic_dimensions_complex_analytic_space} is defined by them to be the \emph{analytic dimension} of $X$ at a point $p$, while the integer $l$ is shown by them to be equal to the \emph{Chevalley dimension} of $\sO_{X,p}$ as in Definition \ref{defn:Chevalley_dimension_module_over_analytic_algebra} and which they call the \emph{algebraic dimension} of $X$ at $p$. Grauert and Remmert prove Lemma \ref{lem:Equivalence_analytic_algebraic_dimensions_complex_analytic_space} using the R\"uckert Nullstellensatz (Theorem \ref{thm:Ruckert_Nullstellensatz}) and that in turn requires the field $\KK$ to be algebraically closed.
\end{rmk}

\section{Embedding dimensions for analytic spaces}
\label{sec:Embedding_dimension_analytic_space}
We begin the followong observation.

\begin{rmk}[Comparison of the dimension of the Zariski tangent space and Krull dimension of the local ring at a point]
\label{rmk:Comparison_dimension_Zariski_tangent_space_and_Krul_dimension_local_ring_at_point}  
Let $(X,\sO_X)$ be an analytic space over $\KK=\RR$ or $\CC$. If $p \in X$ is a point, then the inequality \eqref{eq:Krull_dimension_leq_dimension_Zariski_cotangent_space}, definitions \eqref{eq:Zariski_cotangent_and_tangent_space} of the Zariski cotangent and tangent spaces, and fact that $\sO_{X,p}$ is a Noetherian local ring (by Theorem \ref{thm:Analytic_algebras_are_Noetherian_local_rings} or the argument below) yield the inequality
\begin{equation}
  \label{eq:Krull_dimension_leq_dimension_Zariski_tangent_space}
  \dim\sO_{X,p} \leq \dim_\KK T_pX,
\end{equation}
with equality if and only if $\sO_{X,p}$ is also regular.

The ring of convergent power series $\KK\{x_1,\ldots,x_n\}$ is Noetherian (see Abhyankar \cite[Chapter 2, Section 10, p. 81]{Abhyankar_local_analytic_geometry} or Zariski and Samuel \cite[Section 7.1, Theorem $4'$, p. 148]{Zariski_Samuel_communtative_algebra_II}). Consequently, for any ideal $\fa \subset \KK\{x_1,\ldots,x_n\}$, the quotient $\KK\{x_1,\ldots,x_n\}/\fa$ is also Noetherian by Atiyah and Macdonald \cite[Proposition 6.6]{Atiyah_Macdonald_introduction_commutative_algebra}. Moreover, writing $p=(p_1,\ldots,p_n)$, we see that $\sO_{D,p} = \sO_{\KK^n,p} = \KK\{x_1-p_1,\ldots,x_n-p_n\}$ is a Noetherian ring (see also Narasimhan \cite[Chapter 2, Theorem 4, p. 27]{Narasimhan_introduction_theory_analytic_spaces}). In particular, if $(X,\sO_X)$ is a local analytic model space with $\sO_X = \sO_D/\sI\restriction X$, then $\sO_{X,p} = \sO_{D,p}/\sI_p$ is a Noetherian ring and we recall from Remark \ref{rmk:Analytic_model_space_K-ringed_space} that it is also a local ring.
\end{rmk}

To define the concept of embedding dimension, we shall need the following elementary lemma and corresponding definition.

\begin{lem}[Closed analytic subspaces of analytic spaces]
\label{lem:Closed_analytic_subspace_of_analytic_space}  
(See Fischer \cite[Section 0.14, p. 10]{Fischer_complex_analytic_geometry} or Grauert and Remmert \cite[Section 1.2.2, p. 14]{Grauert_Remmert_coherent_analytic_sheaves} for $\KK=\CC$.)  
Let $\KK=\RR$ or $\CC$ and $(X,\sO_X)$ be a $\KK$-analytic space and $\sI \subset \sO_X$ be a coherent\footnote{Grauert and Remmert only require $\sI$ to be of finite type while Fischer and  and Guaraldo, Macr\`\i, and Tancredi further require that $\sI$ be coherent, as we do here.} sheaf of ideals. Then the induced $\KK$-ringed subspace $(Y,\sO_Y)$ of $(X,\sO_X)$, where $Y := \supp(\sO_X/\sI) \subset X$ and $\sO_Y := (\sO_X/\sI)\restriction Y$, is a $\KK$-analytic space and the injection morphism $\iota:(Y,\sO_Y)\to (X,\sO_X)$ is a $\KK$-analytic map.
\end{lem}

\begin{defn}[Closed analytic subspaces of analytic spaces]
\label{defn:Closed_analytic_subspace_of_analytic_space}
(See Fischer \cite[Section 0.14, p. 10]{Fischer_complex_analytic_geometry}, Grauert and Remmert \cite[Section 1.2.2, p. 14]{Grauert_Remmert_coherent_analytic_sheaves} for $\KK=\CC$ and Guaraldo, Macr\`\i, and Tancredi \cite[Section 2.9, Definition 1.9, p. 15]{Guaraldo_Macri_Tancredi_topics_real_analytic_spaces} for $\KK=\RR$ or $\CC$; see Guaraldo, Macr\`\i, and Tancredi \cite[Section 1.3, Definition 3.1, p. 7]{Guaraldo_Macri_Tancredi_topics_real_analytic_spaces} more generally for embeddings of ringed spaces.)    
In the setting of Lemma \ref{lem:Closed_analytic_subspace_of_analytic_space}, one calls $(Y,\sO_Y)$ a \emph{closed $\KK$-analytic subspace} of $(X,\sO_X)$. 
\end{defn}

We recall the important

\begin{defn}[Embedding dimension]
\label{defn:Embedding_dimension}
(See Fischer \cite[Section 2.3]{Fischer_complex_analytic_geometry} and Grauert and Remmert \cite[Section 6.1.1, p. 113]{Grauert_Remmert_coherent_analytic_sheaves} for $\KK=\CC$ or Guaraldo, Macr\`\i, and Tancredi \cite[Chapter 2, Remark  2.2, p. 19]{Guaraldo_Macri_Tancredi_topics_real_analytic_spaces} for $\KK=\RR$ or $\CC$.)
Let $(X, \sO_X)$ be a $\KK$-analytic space. The \emph{embedding dimension} $\embdim_x X$ of $X$ at a point $x$ is the smallest integer $e\geq 0$ such that an open neighborhood $U \subset X$ is $\KK$-bianalytic to a closed $\KK$-analytic subspace of a domain in $\KK^e$. 
\end{defn}

Locally, a $\KK$-analytic space $X$ is isomorphic to $\KK$-analytic model space, thus a closed $\KK$-analytic subspace of an open subset of $\KK^n$, so the embedding dimension is well-defined at each point of $X$. See Guaraldo, Macr\`\i, and Tancredi \cite[Section 1.3, Definition 3.8, p. 9]{Guaraldo_Macri_Tancredi_topics_real_analytic_spaces} for the definition of an embedding of $\KK$-ringed spaces. The following result gives an important geometric property of the Zariski tangent space.

\begin{prop}[Equality of embedding and Zariski tangent space dimensions] 
\label{prop:Equality_dimension_zariski_tangent_space_and_embedding_dimension}  
(See Fischer \cite[Section 2.3, Proposition, p. 79]{Fischer_complex_analytic_geometry}, Grauert and Remmert \cite[Section 6.1.2, Proposition, p. 115]{Grauert_Remmert_coherent_analytic_sheaves}, or Guaraldo, Macr\`\i, and Tancredi \cite[Remark  2.2.2]{Guaraldo_Macri_Tancredi_topics_real_analytic_spaces}.)  
For every point $p$ of a $\KK$-analytic space $(X, \sO_X)$, one has
\begin{equation}
  \label{eq:Embedding_equal_to_dimension_Zariski_tangent_space}
  \embdim_pX = \dim T_pX.
\end{equation}
\end{prop}

\begin{proof}
We shall give a proof of Proposition \ref{prop:Equality_dimension_zariski_tangent_space_and_embedding_dimension} in a style that is closer to that in differential topology (for example, Lee \cite{Lee_john_smooth_manifolds}). By Definition \ref{defn:Analytic_space}, each point in $X$ has an open neighborhood such that the restriction of $(X,\sO_X)$ to that neighborhood is $\KK$-bianalytic to a $\KK$-analytic model space in a domain in $\KK^n$ for some $n$. Hence, we may assume without loss of generality that $(X,\sO_X)$ is a $\KK$-analytic model space in a domain $D$ in $\KK^n$, with defining ideal sheaf $\sI\subset \sO_D$ and $X = \supp(\sO_D/\sI)$ and $\sO_X = \sO_D/\sI\restriction X$. Assume that $\sI$ is generated by $f_1,\ldots,f_r\in\sO_D(D)$ and $F = (f_1,\ldots,f_r):D\to\KK^r$ is the corresponding $\KK$-analytic map and $T_pX \cong \Ker\d F(p) \subset \KK^n$ by Lemma \ref{lem:Identification_extrinsic_cotangent_space_with_Zariski_cotangent_space} and the definition \eqref{eq:Zariski_cotangent_and_tangent_space} of Zariski cotangent and tangent spaces. Finally, by Definition \ref{defn:Embedding_dimension} of embedding dimension, we may assume without loss of generality that $n = \embdim_pX$.

Let $\Xi\subset \KK^r$ be a linear subspace such that $\KK^r = \Ran\d F(p)\oplus \Xi$ as a direct sum of $\KK$-vector spaces and observe, after possibly shrinking $D$, that $S := F^{-1}(\Xi)$ is an embedded open $\KK$-analytic submanifold of $\KK^n$ since $F$ is transverse to $\Xi \subset \KK^r$ at the point $F(p)=0\in\KK^r$ by construction. Moreover,
\begin{align*}
  T_pS  &= (\d F(p))^{-1}(T_{F(p)}\Xi) = (\d F(p))^{-1}(\Xi)
  \\
                    &= (\d F(p))^{-1}(0) = \Ker\d F(p) \cong T_pX.
\end{align*}  
Because $\d F(p):\KK^n/\Ker\d F(p) \to \Ran \d F(p)$ is an isomorphism of $\KK$-vector spaces and $\Ker\d F(p) \cong T_pX$, then
\[
   n - \dim_\KK \Ker\d F(p) = \dim_\KK\left(\KK^n/\Ker\d F(p)\right) = \dim_\KK\Ran \d F(p).
\]
Therefore, if $s := \dim_\KK\Ran \d F(p)$, we obtain
\begin{equation}
  \label{eq:Equality_dimensions_smooth_manifold_Zariski_tangent_space_nullity_Jacobian}
  \dim_\KK T_pS = \dim_\KK T_pX = \dim_\KK \Ker\d F(p) = n - s.
\end{equation}
Clearly, $(X, \sO_X)$ is a closed $\KK$-analytic subspace (in the sense of Definition \ref{defn:Closed_analytic_subspace_of_analytic_space}) of the $\KK$-analytic manifold $(S,\sO_S)$ and which has dimension $\dim S = n-s$ by \eqref{eq:Equality_dimensions_smooth_manifold_Zariski_tangent_space_nullity_Jacobian}. By Definition \ref{defn:Embedding_dimension} of embedding dimension, we thus have
\[
  \embdim_pX \leq \dim S = n-s = \dim_\KK T_pX.
\]
But $\embdim_pX = n$ by assumption and so we must have $s = 0$ and thus $\dim_\KK T_pX = n$. This completes the proof of Proposition \ref{prop:Equality_dimension_zariski_tangent_space_and_embedding_dimension}.
\end{proof}

\section{Smooth points and regular points for analytic spaces}
\label{sec:Smooth_point_regular_point_analytic_space}
By analogy with Definition \ref{defn:Regular_point_scheme} for schemes, we have the\footnote{This definition and terminology are not used by Fischer \cite{Fischer_complex_analytic_geometry}, Grauert and Remmert \cite{Grauert_Remmert_coherent_analytic_sheaves}, Guaraldo, Macr\`\i, and Tancredi \cite{Guaraldo_Macri_Tancredi_topics_real_analytic_spaces}, or Narasimhan \cite{Narasimhan_introduction_theory_analytic_spaces}.}

\begin{defn}[Regular point of an analytic space]
\label{defn:Regular_point_analytic_space}
(See Abhyankar\footnote{Abhyankar and Hironaka use the term `simple point'.} \cite[Chapter 5, Section 30.16, Definition (4), p. 240, and Section 44.6, Definition (4), p. 404]{Abhyankar_local_analytic_geometry} and Hironaka \cite[Chapter 0, Section 2, p. 121]{Hironaka_1964-I-II} for $\KK=\RR$ or $\CC$ and see Acquistapace, Broglia, and Fernando\footnote{However, these authors give a different and separate definition of `regular point', involving the Jacobian of generators of $\sO_{X,x}$, for a reduced complex analytic space \cite[Section 1.A.1, Definition 1.5, p. 10]{Acquistapace_Broglia_Fernando_topics_global_real_analytic_geometry}.} \cite[Section 1.A.2, p. 15]{Acquistapace_Broglia_Fernando_topics_global_real_analytic_geometry} for $\KK=\CC$ and \cite[Section 2.E, Definition 2.31, p. 40]{Acquistapace_Broglia_Fernando_topics_global_real_analytic_geometry} for $\KK=\RR$.)
Let $\KK=\RR$ or $\CC$ and $(X,\sO_X)$ be an analytic space over $\KK$ in the sense of Definition \ref{defn:Analytic_space}. A point $x \in X$ is \emph{regular} if the local ring $\sO_{X,x}$ is regular and, otherwise, $x$ is called \emph{non-regular}. The analytic space $X$ is \emph{regular} if it is regular at every point. 
\end{defn}

By analogy with Definition \ref{defn:Regular_point_scheme} for schemes, we have the

\begin{defn}[Smooth point of an analytic space]
\label{defn:Smooth_point_analytic_space}
(See Fischer\footnote{Fischer uses the term `non-singular'.} \cite[Section 0.14, p. 9]{Fischer_complex_analytic_geometry} for $\KK=\CC$, Grauert and Remmert\footnote{Grauert and Remmert use the term `smooth', but also use the terms `regular' and `simple' as synonyms.} \cite[Section 1.1.5, p. 8 or Section 6.2.1, p. 116]{Grauert_Remmert_coherent_analytic_sheaves} for $\KK=\CC$, or, for an analytic variety,  Guaraldo, Macr\`\i, and Tancredi \cite[Definition 2.7, p. 22]{Guaraldo_Macri_Tancredi_topics_real_analytic_spaces} for $\KK=\CC$ or $\RR$.)
Let $\KK=\RR$ or $\CC$ and $(X,\sO_X)$ be an analytic space over $\KK$ in the sense of Definition \ref{defn:Analytic_space}. A point $p \in X$ is called \emph{smooth} (or \emph{non-singular}) \emph{(of dimension $d$)} if there is an open neighborhood $U \subset X$ of $p$ such that $(U,\sO_X|_U)$ is isomorphic to an analytic model space $(W,\sO_W)$ in the sense of Definition \ref{defn:Analytic_model_space}, where $W \subset \CC^d$ is an open subset and, otherwise, $p$ is called \emph{non-smooth} (or \emph{singular}). If $X$ is smooth (of dimension $d$) at every point, then $X$ is an \emph{analytic manifold (of dimension $d$)} over $\KK$.
\end{defn}

It will be useful in our applications to specialize Definition \ref{defn:Smooth_point_analytic_space} to the case of an analytic set.

\begin{defn}[Smooth point an analytic set]
\label{defn:Smooth_point_analytic_set}
(See Acquistapace, Broglia, and Fernando\footnote{Acquistapace, Broglia, and Fernando define `smooth' for a `$C$-analytic set', that is, a real analytic set that admits a structure, which may not be unique, as a real analytic space in the sense of Definition \ref{defn:Analytic_space}.} \cite[Section 2.E, Definition 2.31, p. 40]{Acquistapace_Broglia_Fernando_topics_global_real_analytic_geometry} for $\KK=\RR$, Demailly \cite[Definition 4.23, p. 98]{Demailly_complex_analytic_differential_geometry} for $\KK=\CC$, Gunning and Rossi \cite[Chapter III, Section C, Definition 2, pp. 110--111]{Gunning_Rossi_analytic_functions_several_complex_variables} for $\KK=\CC$, or Narasimhan \cite[Chapter III, Section 1, Definition 2, p. 36]{Narasimhan_introduction_theory_analytic_spaces} for $\KK=\RR$ or $\CC$.)
Let $\KK=\RR$ or $\CC$ and $S\subset \Omega$ be an analytic set in the sense of Definition \ref{defn:Analytic_set}, where $\Omega \subset \KK^n$ is an open set. A point $p \in S$ is called \emph{smooth} (or \emph{non-singular}) \emph{of dimension $d$} if there is an open neighborhood $U \subset \Omega$ of $p$ such that $S\cap U$ is an embedded analytic submanifold of dimension $d$ in $U$ and, otherwise, $p$ is called \emph{non-smooth} (or \emph{singular}).
\end{defn}

Definition \ref{defn:Smooth_point_analytic_space} of a smooth point of an analytic space is apparently different from Definition \ref{defn:Smooth_point_scheme} of a smooth point of a scheme. However, these definitions can be seen to be equivalent in the category of analytic spaces via the

\begin{lem}[Jacobian matrix rank condition for a smooth point of an analytic space]
\label{lem:Jacobi_criterion_smooth_point}  
If $(X,\sO_X)$ is an analytic space over a field $\KK=\RR$ or $\CC$, then a point $p \in X$ is smooth of dimension $d$ in the sense of Definition \ref{defn:Smooth_point_analytic_space} if and only if there exists an open neighborhood $U\subset X$ of $p$ such that $(U,\sO_X|_U)$ is isomorphic to an analytic model space $(Y,\sO_Y)$, where $\sO_Y = \sO_D/(f_1,\ldots,f_{n-d})$ for a domain $D\subset\KK^n$ and (after identifying $p\in X$ with its image in $Y$) the Jacobian matrix
\begin{equation}
  \label{eq:Jacobian_matrix_analytic_space}
  J_{f_1,\ldots,f_{n-d}}(p) = \left(\frac{\partial f_i}{\partial x_j}(p)\right)_{(n-d)\times n}
  \in \Hom_\KK(\KK^n,\KK^{n-d}),
\end{equation}
obeys the corank condition,
\begin{equation}
  \label{eq:Jacobi_criterion_analytic_space}
  \Corank_\KK J_{f_1,\ldots,f_{n-d}}(p) = d,
\end{equation}
or, equivalently, the rank condition,
\[
  \Rank_\KK J_{f_1,\ldots,f_{n-d}}(p) = n-d.
\]
\end{lem}

\begin{proof}
If $p\in X$ is smooth in the sense of Definition \ref{defn:Smooth_point_analytic_space}, then $(U,\sO_X|_U) \cong (D,\sO_D)$, so $(f_1,\ldots,f_{n-d}) = (0)$ and thus $J_{f_1,\ldots,f_{n-d}}(p) = 0 \in \Hom_\KK(\KK^n,0)$ has rank zero and corank $d$, where $d=n$, and the Jacobian matrix rank condition is trivially obeyed.

Conversely, if $(U,\sO_X|_U) \cong (Y,\sO_Y)$, where $\sO_Y = \sO_D/(f_1,\ldots,f_{n-d})$ and $J_{f_1,\ldots,f_{n-d}}(p)$ has rank $n-d$ over $\KK$, then, after possibly shrinking $D$, the $\KK$-ringed space $(Y,\sO_Y)$ is a $\KK$-analytic manifold of dimension $d$ as a consequence of the Implicit Function Theorem for $\KK$-analytic maps: see, for example, Griffiths and Harris \cite[Section 0.1, pp. 19--20]{GriffithsHarris}, Huybrechts \cite[Proposition 1.1.11 and Corollary 1.1.12]{Huybrechts_2005}, or Noguchi \cite[Theorem 1.2.41, Definition 1.2.44, and Theorem 1.2.45]{Noguchi_analytic_function_theory_several_variables} for $\KK=\CC$ and Abhyankar \cite[Chapter 2, Section 10.8, p. 84]{Abhyankar_local_analytic_geometry} or Krantz and Park \cite[Theorem 6.1.2]{Krantz_Parks_implicit_function_theorem} for $\KK=\RR$ or $\CC$.

Indeed, if we define a $\KK$-analytic map by $F = (f_1,\ldots,f_{n-d}):\KK^n \supset D \to \KK^{n-d}$, then $Y = F^{-1}(0)$ and there are an open neighborhood $D'$ of the origin in $\KK^d$ and a $\KK$-analytic map $g = (g_{d+1},\ldots,g_n):\KK^d \supset D' \to \KK^{n-d}$ such that $(D',g(D')) \subset D \subset \KK^d\times \KK^{n-d}$ and
\[
  F(x_1,\ldots,x_d,g_{d+1}(x_1,\ldots,x_d),\ldots,g_n(x_1,\ldots,x_d)) = 0,
  \quad\text{for all } (x_1,\ldots,x_d) \in D'.
\]
Therefore, $(Y,\sO_Y)$ is a $\KK$-analytic manifold of dimension $d$ and so $(U,\sO_X|_U)$ is a $\KK$-analytic manifold of dimension $d$ and $p$ is a smooth point of $X$ by Definition  \ref{defn:Smooth_point_analytic_space}.
\end{proof}  

Let $(X,\sO_X)$ be an analytic space over $\KK=\RR$ or $\CC$. Recall from Definition \ref{defn:Regular_local_ring} that $\sO_{X,p}$ is a regular local ring if
\[
  \dim\sO_{X,p} = \dim_\KK \fm_p/\fm_p^2,
\]
and recall from \eqref{eq:Zariski_cotangent_and_tangent_space} that $T_p^*X = \fm_p/\fm_p^2$ and $T_pX = \Hom_\KK(\fm_p/\fm_p^2,\KK)$ are the Zariski cotangent and tangent spaces, respectively, to $X$ at the point $p$. We note the

\begin{lem}[Dimensions of a complex analytic space and its reduction are equal]
(See Grauert and Remmert \cite[Section 5.1.3, Observation, p. 96]{Grauert_Remmert_coherent_analytic_sheaves} for the statement below or De Jong and Pfister \cite[Lemma 4.1.6, p. 130]{DeJong_Pfister_local_analytic_geometry} for the Krull dimension of any Noetherian local ring.)  
Let $(X,\sO_X)$ be a complex analytic space. Then the dimensions of $X$ and $X^{\red}$ coincide
at all points, so $\dim_p X = \dim_p X^{\red}$ for all $p \in X$.
\end{lem}

We have the following analogue of Theorem \ref{thm:Gortz_Wedhorn_6-28} for analytic spaces.

\begin{thm}[Comparison of smoothness and regularity for analytic spaces]
\label{thm:Comparison_smoothness_regularity_analytic_spaces}
(See Abhyankar\footnote{In general, Abhyankar allows $\KK$ to be any complete (non-discrete) valued field as he notes in \cite[Preface, p. viii]{Abhyankar_local_analytic_geometry}, but in \cite[Chapter 5, p. 230]{Abhyankar_local_analytic_geometry} for Sections 30--39 and \cite[Chapter 7, p. 394]{Abhyankar_local_analytic_geometry} for Sections 44--46, he further assumes that $\KK$ is algebraically closed, though some of the results remain true without that assumption.} \cite[Chapter 5, Section 33.18, p. 279 and Chapter 7, Section 44.19, p. 408]{Abhyankar_local_analytic_geometry}, Acquistapace, Broglia, and Fernando \cite[Definition 1.6, Proposition 1.7, and Theorem 1.12]{Acquistapace_Broglia_Fernando_topics_global_real_analytic_geometry} for $\KK=\CC$ and reduced complex analytic spaces, Fischer \cite[Section 0.22, Corollary 1, p. 18]{Fischer_complex_analytic_geometry}, De Jong and Pfister \cite[Theorem 4.3.10, p. 153]{DeJong_Pfister_local_analytic_geometry}, Hironaka \cite[Remark 5.5, Definition 5.6, and page 5.1]{Hironaka_intro_real-analytic_sets_maps}, and Grauert and Remmert \cite[Section 6.2.1, Criterion of Smoothness, p. 116]{Grauert_Remmert_coherent_analytic_sheaves} for $\KK=\CC$.)
Let $\KK=\RR$ or $\CC$ and $(X,\sO_X)$ be an analytic space over $\KK$ and $p\in X$ be a point and $d\geq 0$ be an integer. Then the following are equivalent:
\begin{enumerate}
\item\label{item:X_smooth_analytic_space_at_point_dimension_d}
  The analytic space $X$ is smooth of dimension $d$ at $p$ in the sense of Definition  \ref{defn:Smooth_point_analytic_space}.
\item\label{item:X_regular_analytic_space_at_point_dimension_d}
  The local ring $\sO_{X,p}$ is regular and has Krull dimension $d$, so the following equalities hold:
\begin{equation}
  \label{eq:Dimension_Zariski_tangent_space_equals_dimension_local_ring_analytic_space}
    \dim_\KK T_pX = \dim\sO_{X,p} = d,
\end{equation}
where $T_pX$ is the Zariski tangent space to $X$ at the point $p$.
\end{enumerate}
\end{thm}

Demailly provides the following closely related result.

\begin{prop}[Characterization of a smooth point of a complex analytic space in terms of algebraic property 
of the local ring]
\label{prop:Demailly_4-32}  
(See Demailly \cite[Proposition 4.32, p. 101]{Demailly_complex_analytic_differential_geometry}.)  
Let $(S, x)$ be a germ of a complex analytic set $S$ of dimension $d = \dim\sO_{S,x}$ and let $\fm_{S,x} \subset \sO_{S,x}$ be the maximal ideal of holomorphic functions that vanish at $x$. Then $\fm_{S,x}$ cannot have less than $d$ generators and $\fm_{S,x}$ has $d$ generators if and only if $x$ is a smooth point as in Definition \ref{defn:Smooth_point_analytic_space}.
\end{prop}  

\begin{proof}[Proof of Theorem \ref{thm:Comparison_smoothness_regularity_analytic_spaces} for complex analytic spaces]
While implicit in the cited references, the variations among definitions, terminology, and emphases in their expositions obscure the conclusion and thus we summarize the argument here from the indicated references when $\KK=\CC$ and translate their terminology as needed. While Grauert and Remmert assume that $\KK=\CC$, Abhyankar allows in \cite[Sections 33 and 44]{Abhyankar_local_analytic_geometry} that $\KK$ be an algebraically closed, complete, non-discrete valued field, although he notes that algebraic closure is not required for every result stated by him in those sections. Hironaka allows $\KK=\RR$ or $\CC$, but only briefly outlines the argument in \cite[Chapter 5]{Hironaka_intro_real-analytic_sets_maps}, a limited circulation monograph that is difficult to obtain. 

Grauert and Remmert \cite[Section 6.2.1, Criterion of Smoothness, p. 116]{Grauert_Remmert_coherent_analytic_sheaves} prove the equivalence of \eqref{item:X_smooth_analytic_space_at_point_dimension_d} and the fact that the following identity holds,
\[
  \dim_\KK T_pX = \dim_pX = d,
\]
where ``$\dim_pX$'' is the \emph{analytic dimension} of $X$ at the point $p$, as defined by Grauert and Remmert in Item \eqref{item:Analytic_dimension_point_complex_space} of Lemma \ref{lem:Equivalence_analytic_algebraic_dimensions_complex_analytic_space}. But Corollary \ref{cor:Chevalley_equal_to_Krull_dimension_analytic_algebra}, Lemma \ref{lem:Equivalence_analytic_algebraic_dimensions_complex_analytic_space}, and Remark \ref{rmk:Equivalence_analytic_algebraic_dimensions_complex_analytic_space} also imply that this analytic dimension is equal to the Krull dimension $\dim\sO_{X,p}$ and this yields the equivalence of Item \eqref{item:X_smooth_analytic_space_at_point_dimension_d} and the fact that the identity \eqref{eq:Dimension_Zariski_tangent_space_equals_dimension_local_ring_analytic_space} holds. Because $T_pX = (\fm_p/\fm_p^2)^*$ by definition \eqref{eq:Zariski_cotangent_and_tangent_space}, then 
\[
  \dim\sO_{X,p} = \dim_\CC(\fm_p/\fm_p^2)
\]
and so $\sO_{X,p}$ is a regular local ring by Definition \ref{defn:Regular_local_ring}, since its residue field $\kappa(p) = \sO_{X,p}/\fm_p$ is isomorphic to $\CC$ by virtue of the fact that $(X,\sO_X)$ is a $\CC$-ringed space.

Fischer \cite[Section 0.22, Corollary 1, p. 18]{Fischer_complex_analytic_geometry} outlines a proof for $\KK=\CC$ of equivalence between the condition that $X$ is smooth at $p$ in the sense of Definition  \ref{defn:Smooth_point_analytic_space} and that the local ring $\sO_{X,p}$ is regular, though he does not prove the indicated equalities of dimensions.

Abyhyankar obtains the equivalence of Items \eqref{item:X_regular_analytic_space_at_point_dimension_d} and \eqref{item:X_smooth_analytic_space_at_point_dimension_d} through his equivalence of Items (i) and (v) in \cite[Chapter 5, Section 33.18, p. 279]{Abhyankar_local_analytic_geometry} or through \cite[Chapter 7, Section 44.19, p. 408]{Abhyankar_local_analytic_geometry}.

When $(X,\sO_X)$ is a reduced complex analytic space, Acquistapace, Broglia, and Fernando obtain the equivalence of Items \eqref{item:X_regular_analytic_space_at_point_dimension_d} and \eqref{item:X_smooth_analytic_space_at_point_dimension_d} through their \cite[Theorem 1.12]{Acquistapace_Broglia_Fernando_topics_global_real_analytic_geometry} and the equivalence of the property of $p$ being a ``regular point'' in the sense of their \cite[Definition 1.6]{Acquistapace_Broglia_Fernando_topics_global_real_analytic_geometry} and the property of $X$ being smooth at $p$ by their \cite[Proposition 1.7]{Acquistapace_Broglia_Fernando_topics_global_real_analytic_geometry}.
\end{proof}

Before giving a self-contained proof of Theorem \ref{thm:Comparison_smoothness_regularity_analytic_spaces} that applies simultaneously to $\KK=\RR$ or $\CC$, we shall need a definition and a preparatory lemma. With the exception of the Stacks Project \cite[\href{https://stacks.math.columbia.edu/tag/0094}{Section 0094}]{stacks-project}, standard references for sheaf theory tend to only define a morphism between a pair of sheaves $\sF$ and $\sG$ over a given topological space, $X$, and not a morphism from a sheaf $\sF$ over $X$ to a sheaf $\sG$ over another topological space, $Y$. For this reason, we include the

\begin{defn}[Morphism of sheaves of modules over ringed spaces]
\label{defn:Morphism_sheaves_modules_over_ringed_spaces}
(See the paragraph following
\cite[\href{https://stacks.math.columbia.edu/tag/0097}{Lemma 0097}]{stacks-project} in the Stacks Project.)  
Let $(\varphi,\varphi^\sharp):(X,\sO_X) \to (Y,\sO_Y)$ be a morphism of ringed spaces, $\sF$ be a sheaf of $\sO_X$-modules over $X$, and $\sG$ be a sheaf of $\sO_Y$-modules over $Y$. One calls $\alpha:\sG\to\sF$ a $\varphi$-\emph{morphism of sheaves of modules} if for every open subset $V \subset Y$, then
\begin{equation}
  \label{eq:varphi-morphism_presheaves_modules}
  \alpha_V:\sG(V) \to \sF(\varphi^{-1}(V))
\end{equation}
is a morphism of $\sO_Y(V)$-modules, where $\sF(\varphi^{-1}(V))$ inherits an $\sO_Y(V)$-module structure via the morphism of rings\footnote{See the Stacks Project \cite[\href{https://stacks.math.columbia.edu/tag/0091}{Definition 0091}]{stacks-project} for the definition of a morphism of ringed spaces.}
\begin{equation}
  \label{eq:varphi-morphism_presheaves_ringed_spaces}
  \varphi_V^\sharp:\sO_Y(V) \to \sO_X(\varphi^{-1}(V)).
\end{equation}
When $(X,\sO_X) = (Y,\sO_Y)$ and $(\varphi,\varphi^\sharp)$ is the identity morphism, one calls $\alpha:\sG\to\sF$ a \emph{morphism of sheaves of modules} (see the Stacks Project \cite[\href{https://stacks.math.columbia.edu/tag/0077}{Definition 0077}]{stacks-project} or Grauert and Remmert \cite[Annex, Section 0.1, p. 223 and Section 2.2, p. 229]{Grauert_Remmert_coherent_analytic_sheaves}).
\end{defn}

\begin{rmk}[Morphism of sheaves of modules over ringed spaces and pushforward sheaves]
\label{rmk:Morphism_sheaves_modules_over_ringed_spaces_via_pushforward}  
In the paragraph following
\cite[\href{https://stacks.math.columbia.edu/tag/0097}{Lemma 0097}]{stacks-project} in the Stacks Project, it is noted that a $\varphi$-morphism of sheaves of modules as in Definition \ref{defn:Morphism_sheaves_modules_over_ringed_spaces},
\[
  \alpha:\sG \to \sF,
\]  
can be canonically identified with a morphism of sheaves of $\sO_Y$-modules over $Y$,
\[
  \alpha_Y:\sG \to \varphi_*\sF,
\]
where $\varphi_*\sF$ is the pushforward of $\sF$ to $Y$ and is a sheaf of $\sO_Y$-modules (see the Stacks Project \cite[\href{https://stacks.math.columbia.edu/tag/0095}{Definition 0095}]{stacks-project}). Indeed,
\[
  (\varphi_*\sF)(V) := \sF(\varphi^{-1}(V)), \quad\text{for any open subset } V \subset Y,
\]
and this defines the pushforward sheaf $\varphi_*\sF$ as a sheaf of sets over a topological space $Y$, following the Stacks Project \cite[\href{https://stacks.math.columbia.edu/tag/008C}{Section 008C} and \href{https://stacks.math.columbia.edu/tag/008D}{Lemma 008D}]{stacks-project}. Moreover, by the Stacks Project \cite[\href{https://stacks.math.columbia.edu/tag/008S}{Lemma 008S}]{stacks-project}, we observe that $\varphi_*\sF$ is a (pre-)sheaf of $\varphi_*\sO_X$-modules via the multiplication map
\[
  (\varphi_*\sO_X)(V) \times (\varphi_*\sF)(V)
  =
  \sO_X(\varphi^{-1}(V)) \times \sF(\varphi^{-1}(V)) \to \sF(\varphi^{-1}(V)) = (\varphi_*\sF)(V).
\]
Thus, as noted in Definition \ref{defn:Morphism_sheaves_modules_over_ringed_spaces}, the morphism \eqref{eq:varphi-morphism_presheaves_ringed_spaces} of rings yields 
the multiplication map
\[
  \sO_Y(V) \times (\varphi_*\sF)(V)
  =
  \sO_Y(V) \times \sF(\varphi^{-1}(V)) \to \sF(\varphi^{-1}(V)) = (\varphi_*\sF)(V),
\]
and so $\varphi_*\sF$ may be viewed as a (pre-)sheaf of $\sO_Y$-modules. 
\end{rmk}

The forthcoming lemma generalizes Fischer \cite[Section 0.6, Corollary 2, p. 5]{Fischer_complex_analytic_geometry}, who restricts to the case where $(X,\sO_X) = (Y,\sO_Y)$ and omits the proof. The forthcoming lemma also generalizes Grauert and Remmert \cite[Section A.3.1, Statement 2, p. 234]{Grauert_Remmert_coherent_analytic_sheaves}, who restrict to the case where $(X,\sO_X) = (Y,\sO_Y)$ and $\alpha_x$ is an epimorphism but only assume that $\sF$ has finite type. Lastly, the forthcoming lemma generalizes G\"ortz and Wedhorn \cite[Lemma 6.17]{Gortz_Wedhorn_algebraic_geometry_v1}, who restrict to the case where $X$ and $Y$ are schemes locally of finite type over a field $\KK$ and $\sF=\sO_X$ and $\sG=\sO_Y$.

\begin{lem}[Morphism of coherent sheaves of modules that induces an epimorphism (monomorphism, isomorphism) of stalks at a point]
\label{lem:Morphism_coherent_sheaves_modules_implies_isomorphism_stalks_at_point}  
Let $(\varphi,\varphi^\sharp):(X,\sO_X) \to (Y,\sO_Y)$ be a morphism of ringed spaces, $\sF$ be a coherent sheaf of $\sO_X$-modules over $X$, and $\sG$ be a coherent sheaf of $\sO_Y$-modules over $Y$, and $\alpha:\sG\to\sF$ be a $\varphi$-morphism of sheaves of modules as in Definition \ref{defn:Morphism_sheaves_modules_over_ringed_spaces}. If $\alpha_x:\sG_{\varphi(x)}\to\sF_x$ is an epimorphism (respectively, monomorphism or isomorphism) of modules, then $\alpha$ has the same property as a morphism of sheaves of modules from some open neighborhood $V\subset Y$ of $\varphi(x)$ to the corresponding open neighborhood $\varphi^{-1}(V)\subset X$ of $x$.     
\end{lem}

\begin{proof}
As noted in Remark \ref{rmk:Morphism_sheaves_modules_over_ringed_spaces_via_pushforward}, the $\varphi$-morphism $\alpha:\sG\to\sF$ of sheaves of modules may be canonically identified with a morphism of sheaves of $\sO_Y$-modules over $Y$,
\[
  \alpha_Y:\sG \to \varphi_*\sF.
\]
Suppose first that $\alpha_x:\sG_{\varphi(x)}\to\sF_x$ is an epimorphism of $\sO_{Y,\varphi(x)}$-modules, where $\sF_x$ is an $\sO_{Y,\varphi(x)}$-module via the morphism of rings (see Definition \ref{defn:Morphism_sheaves_modules_over_ringed_spaces}) 
\[
  \varphi_x^\sharp:\sO_{Y,\varphi(x)} \to \sO_{X,x}.
\]
Consequently, $\alpha_{Y,\varphi(x)}:\sG_{\varphi(x)} \to (\varphi_*\sF)_{\varphi(x)}$ is an epimorphism of $\sO_{Y,\varphi(x)}$-modules. Since both $\sG$ and $\varphi_*\sF$ are sheaves of $\sO_Y$-modules over $Y$, we can apply the following result proved by Grauert and Remmert as \cite[Section A.4.2, Consequence 3, p. 237]{Grauert_Remmert_coherent_analytic_sheaves}: If a sequence of coherent sheaves of $\sO_X$-modules over a ringed space $(X,\sO_X)$,
\[
  \sR \xrightarrow{\alpha} \sS \xrightarrow{\beta} \sT
\]
is exact at a point $x\in X$, then there is an open neighborhood $U\subset X$ of $x$ such that the following sequence of sheaves of $\sO_X|_U$-modules over $U$ is exact:
\[
  \sR_U \xrightarrow{\alpha_U} \sS_U \xrightarrow{\beta_U} \sT_U
\]
Hence, there is an open neighborhood $V\subset Y$ of $\varphi(x)$ such that the following morphism of sheaves of $\sO_Y|_V$-modules over $V$ is an epimorphism: 
\[
  \alpha_Y:\sG_V \to (\varphi_*\sF)_V.
\]
Therefore, by the canonical identification, the following $\varphi$-morphism of sheaves of modules is an epimorphism:
\[
  \alpha:\sG_V \to \sF_{\varphi^{-1}(V)}.
\]
By applying the same argument with epimorphism replaced by monomorphism, we obtain the analogous conclusion for monomorphisms and by combining the conclusions from both steps, we obtain the analogous conclusion for isomorphisms. 
\end{proof}

\begin{cor}[Morphism of coherent ringed spaces that induces an epimorphism (monomorphism, isomorphism) of stalks at a point]
\label{cor:Morphism_coherent_ringed_spaces_implies_isomorphism_stalks_at_point}  
Let $(\varphi,\varphi^\sharp):(X,\sO_X) \to (Y,\sO_Y)$ be a morphism of ringed spaces, where $\sO_X$ and $\sO_Y$ are coherent sheaves (see Definition \ref{defn:Coherent_sheaf}). If $\varphi_x:\sO_x\to\sO_{\varphi(x)}$ is an epimorphism (respectively, monomorphism or isomorphism) of rings, then $(\varphi,\varphi^\sharp)$ has the same property as a morphism of ringed spaces on some open neighborhood $V\subset Y$ of $\varphi(x)$ and corresponding open neighborhood $\varphi^{-1}(V)\subset X$ of $x$.     
\end{cor}

\begin{proof}
Apply Lemma \ref{lem:Morphism_coherent_sheaves_modules_implies_isomorphism_stalks_at_point} with $\sF = \sO_X$ and $\sG = \sO_Y$.
\end{proof}

\begin{rmk}[Application to morphisms of analytic manifolds]
\label{rmk:Application_morphism_analytic_manifolds}  
Suppose that $(\varphi,\varphi^\sharp):X \to Y$ is a morphism of analytic manifolds over $\KK=\RR$ or $\CC$ and that $\varphi_x^\sharp:\sO_{Y,\varphi(x)}\to\sO_{X,x}$ is an isomorphism of rings. By Oka's Theorem \ref{thm:Oka_coherence_theorem}, the sheaves $\sO_X$ and $\sO_Y$ are coherent in the sense of Definition \ref{defn:Coherent_sheaf}. We can thus apply Corollary \ref{cor:Morphism_coherent_ringed_spaces_implies_isomorphism_stalks_at_point} to conclude that there is an open neighborhood $V\subset Y$ of $\varphi(x)$ such that  $\varphi_V^\sharp:\sO_Y|_V\to\sO_X|_{\varphi^{-1}(V)}$ is an isomorphism of sheaves of rings and hence that $(\varphi,\varphi^\sharp):(\varphi^{-1}(V), \sO_X|_{\varphi^{-1}(V)}) \to (V,\sO_Y|_V)$ is an isomorphism of analytic manifolds (see Definition \ref{defn:Analytic_space}). 
\end{rmk}

Recall that if $D\subset \KK^n$ is a domain and $p\in D$ a point and $f \in \sO_D$, then its differential at $p$ is
\begin{equation}
  \label{eq:Differential_function_at_point}
  \d f(p) := \left(\frac{\partial f}{\partial x_1}(p),\ldots,\frac{\partial f}{\partial x_n}(p)\right) \in \Hom_\KK(\KK^n,\KK),
\end{equation}
and we thus obtain a linear map
\begin{equation}
  \label{eq:Differential_map_at_point}
  \mathrm{d}_p:\sO_D \ni f \mapsto \d f(p) \in (\KK^n)^*.
\end{equation}
We have the well-known

\begin{lem}[Identification of the extrinsic cotangent space at point with the maximal ideal modulo its square]
\label{lem:Identification_extrinsic_cotangent_space_with_Zariski_cotangent_space}
(Compare Shafarevich \cite[Chapter 2, Section 1.3, Theorem 2.1, p. 87]{Shafarevich_v1} or Fischer \cite[Section 2.2, p. 78]{Fischer_complex_analytic_geometry}.)  
Let $(X, \sO_X)$ be a $\KK$-analytic model space in a domain $D$ of $\KK^n$ in the sense of Definition \ref{defn:Analytic_model_space} and $p\in X$ be a point. If $\sO_X = \sO_D/(f_1,\ldots,f_r) \restriction X$, then the differential $\mathrm{d}_p$ defines an isomorphism of $\KK$-vector spaces, 
\begin{equation}
  \label{eq:Zariski_cotangent_space_isomorphic_to_dual_kernel_differential_defining_map}
  \fm_{X,p}/\fm_{X,p}^2 \cong (\Ker J_{f_1,\ldots,f_r}(p))^*,
\end{equation}
where $\fm_{X,p} \subset \sO_{X,p}$ is the maximal ideal and
\[
  J_{f_1,\ldots,f_r}(p)
  =
  (\d f_1(p),\ldots,\d f_r(p))
  =
  \begin{pmatrix}
    \frac{\partial f_1}{\partial x_1}(p) & \cdots & \frac{\partial f_1}{\partial x_n}(p)
    \\
    \vdots & & \vdots
    \\
    \frac{\partial f_r}{\partial x_1}(p) & \cdots & \frac{\partial f_r}{\partial x_n}(p)
  \end{pmatrix}
  \in \Hom_\KK(\KK^n,\KK^r).
\]
\end{lem}

\begin{rmk}[Expected dimension]
\label{rmk:Expected_dimension}  
Continue the notation of Lemma \ref{lem:Identification_extrinsic_cotangent_space_with_Zariski_cotangent_space}. We call
\begin{equation}
  \label{eq:Expected_dimension_analytic_model_space}
  \expdim_p X := n-r
\end{equation}
the \emph{expected dimension} of $(X,\sO_D/(f_1,\ldots,f_r)\restriction X)$ at the point $p$, noting that $D\subset\KK^n$ is a domain. If $p \in X$ is a smooth point of dimension $d$ in the sense of Definition \ref{defn:Smooth_point_analytic_space}, then Lemma  \ref{lem:Jacobi_criterion_smooth_point} implies that we can take $r=n-d$ and the isomorphism \eqref{eq:Zariski_cotangent_space_isomorphic_to_dual_kernel_differential_defining_map} and the definition \eqref{eq:Zariski_cotangent_and_tangent_space} of Zariski tangent and cotangent spaces yield
\[
  \expdim_pX = \dim_\KK T_pX.
\]
Unlike the dimension over $\KK$ of the Zariski tangent space $T_pX$ or the Krull dimension of the local ring $\sO_{X,p}$, expected dimension is not intrinsically defined.
\end{rmk}

\begin{proof}[Proof of Lemma \ref{lem:Identification_extrinsic_cotangent_space_with_Zariski_cotangent_space}]
It is convenient to define the analytic map
\[
  F = (f_1,\ldots,f_r):\KK^n\supset D\to \KK^r,
\]
and thus write $\d F(p) = J_{f_1,\ldots,f_r}(p)$. To show that $\mathrm{d}_p$ is surjective, let $\alpha \in (\Ker \d F(p))^*$. Because $\Ker \d F(p) \subset \KK^n$ is a linear subspace, then any basis for $\Ker \d F(p)$ can be extended to a basis for $\KK^n$ and thus gives
\[
  \KK^n = \Ker \d F(p) \oplus (\Ker \d F(p))^\perp
\]
as a direct sum of vector spaces, where $(\Ker \d F(p))^\perp$ is defined to be the linear span of the additional basis vectors. Define $f \in \Hom(\KK^n,\KK)$ by setting $f=\alpha-\alpha(p)$ on $\Ker \d F(p)$ and $f=0$ on $(\Ker \d F(p))^\perp$. Observe that $\d f(p) = \alpha$ and $f(p)=0$ and thus $[f]_p\in \fm_{X,p}$, for any open neighborhood $U\subset X$ of $p$ and pair $(U,f)$ representing the germ $[f]_p$ of $f$ at $p$.

To see that $\Ker\mathrm{d}_p = \fm_{X,p}^2$, suppose $[g]_p\in\fm_{X,p}$ obeys $\d g(p)=0$ on $\Ker \d F(p)$. Because
\[
  \Ker \d F(p) = \Ker \d f_1(p) \cap \cdots \cap\Ker \d f_r(p),
\]
there are constants $\lambda_1,\ldots,\lambda_r\in\KK$ such that
\[
  \d g(p) = \lambda_1\d f_1(p) + \cdots + \lambda_r\d f_r(p) \quad\text{on } \KK^n.
\]
Define $h := g - \lambda_1f_1 + \cdots + \lambda_rf_r$ and observe that $\d h(p) = 0$, so $h \in \fm_{X,p}$, and ${\mathrm d}^2 h(p) = 0$, and thus $h \in \fm_{X,p}^2$, with the pair $(U,h)$ representing the germ $[h]_p$ of $h$ at $p$. Here,
\[
  {\mathrm d}^2 h(p) := \left(\frac{\partial^2 h}{\partial x_i\partial x_j}(p)\right) \in \Hom(\KK^n\otimes\KK^n,\KK),
\]
again equivalently viewed as an element of $\KK^n\otimes\KK^n$ or $(\KK^n)^*\otimes(\KK^n)^*$ using the canonical bases. 
\end{proof}

We are finally ready to give the

\begin{proof}[Proof of Theorem \ref{thm:Comparison_smoothness_regularity_analytic_spaces}]
For $\KK=\RR$, one approach might be to use \emph{base change}, transforming the question of equivalence of Items \eqref{item:X_smooth_analytic_space_at_point_dimension_d} and \eqref{item:X_regular_analytic_space_at_point_dimension_d} for $\KK=\RR$ to the corresponding question for $\KK=\CC$, which we have already treated. However, the references cited in Section \ref{sec:Complexification_real_analytic_spaces} --- which discuss complexification of a real analytic space --- and Section \ref{sec:Real_part_complex_analytic_spaces} --- which discuss taking the real part of a complex analytic space --- do not contain all the results we need to provide a concise argument. Instead, we shall adapt the proof of Theorem \ref{thm:Gortz_Wedhorn_6-28} for schemes over $\KK$, which G\"ortz and Wedhorn allow to be any field. Abhyankar's proof of Theorem \ref{thm:Comparison_smoothness_regularity_analytic_spaces} requires that $\KK$ be algebraically closed through his appeal in \cite[Chapter 5, Section 33.17, p. 279]{Abhyankar_local_analytic_geometry} to the R\"uckert Nullstellensatz (see the forthcoming Theorem \ref{thm:Ruckert_Nullstellensatz}).

In the argument that we present here, we do not require $\KK$ to be algebraically closed and thus allow $\KK=\RR$ or $\CC$. We may assume without loss of generality that $(X,\sO_X)$ is a local model space, so $\sO_X = \sO_D/\sI\restriction X$ for a domain $D\subset\KK^n$ and $X = \supp(\sO_D/\sI)$, with $\sI = (f_1,\ldots,f_r) \subset \sO_D$. 

  Assume Item \eqref{item:X_smooth_analytic_space_at_point_dimension_d}, namely that $X$ is smooth of dimension $d$ at $p$ in the sense of Definition \ref{defn:Smooth_point_analytic_space}.  We shall adapt the proof by  G\"ortz and Wedhorn of \cite[Lemma 6.26]{Gortz_Wedhorn_algebraic_geometry_v1} from the category of schemes to show that $\sO_{X,p}$ is a regular local ring of dimension $d$ and thus Item \eqref{item:X_regular_analytic_space_at_point_dimension_d} holds. (See Abhyankar \cite[Chapter 2, Section 10.11, p. 89]{Abhyankar_local_analytic_geometry} for a proof of the analogous result in the category of analytic spaces.) By Lemma \ref{lem:Jacobi_criterion_smooth_point}, we may assume that $m=n-d$ for this direction. Consider the analytic map $F = (f_1,\ldots,f_{n-d}):\KK^n \supset D \to \KK^{n-d}$ and its differential (Jacobian matrix) at $p\in D$,
\[
  \d F(p) = J_{f_1,\ldots,f_{n-d}}(p) = \left(\frac{\partial f_i}{\partial x_j}(p)\right)_{(n-d)\times n} \in \Hom_\KK(\KK^n,\KK^{n-d}).
\]
Because $X$ is smooth at $p$ of dimension $d$, Lemma \ref{lem:Jacobi_criterion_smooth_point} implies that 
\begin{equation}
  \label{eq:Corank_Jacobian_smooth_point}
  \Corank_\KK \d F(p) = d,
\end{equation}
where $\Corank_\KK \d F(p) = n - \Rank_\KK \d F(p)$. Observe that $X = F^{-1}(0)$ and that the natural morphisms of analytic spaces,
\[
  X \hookrightarrow D \xrightarrow{F} \KK^{n-d},
\]
yields an exact sequence of rings (noting that $F_p^\sharp[g]_p = [g\circ F]_p = [g\circ (f_1,\ldots f_{n-d})]_p$, for all $[g]_p \in \sO_{\KK^{n-d},F(p)}$),
\[
  \sO_{\KK^{n-d},F(p)} \xrightarrow{F_p^\sharp} \sO_{D,p} \to \sO_{D,p}/(f_1,\ldots f_{n-d})_p = \sO_{X,p} \to 0, 
\]
where $(f_1,\ldots f_{n-d})_p = ([f_1]_p,\ldots [f_{n-d}]_p)$ is the ideal in $\sO_{D,p}$ generated by $[f_1]_p,\ldots [f_{n-d}]_p$. By taking the maximal ideals of these local rings modulo their squares, we obtain an exact sequence of vector spaces,
\[
  \fm_{\KK^{n-d},F(p)}/ \fm_{\KK^{n-d},F(p)}^2 \to \fm_{\KK^n,p}/\fm_{\KK^n,p}^2 \to \fm_{X,p}/\fm_{X,p}^2 \to 0 
\]
where $\fm_{\KK^{n-d},F(p)} \subset \sO_{\KK^{n-d},F(p)}$ and $\fm_{\KK^n,p} \subset \sO_{\KK^n,p} = \sO_{D,p}$ and $\fm_{X,p} \subset \sO_{X,p}$ are the maximal ideals, that is, by definition of the Zariski cotangent spaces \eqref{eq:Zariski_cotangent_and_tangent_space},
\[
  T_{F(p)}^*\KK^{n-d} \to T_p^*\KK^n \to T_p^*X \to 0 
\]
By dualizing, we obtain an exact sequence of vector spaces,
\[
  0 \to T_pX \to T_p\KK^n \xrightarrow{\d F(p)} T_{F(p)}\KK^{n-d},
\]
where $T_pX = (\fm_{X,p}/\fm_{X,p}^2)^* \cong \Ker \d F(p)$ by Lemma \ref{lem:Identification_extrinsic_cotangent_space_with_Zariski_cotangent_space}.

Now the matrix $\d F(p)$ has rank $n-d$, so the images of $f_1,\ldots,f_{n-d} \in \sO_D$ in $T_p^*\KK^n = \fm_{\KK^n,p}/\fm_{\KK^n,p}^2$ are linearly independent over $\KK$. But $\sO_{\KK^n,p}$ is a regular local ring and by G\"ortz and Wedhorn \cite[Proposition B.77 (3), p. 572]{Gortz_Wedhorn_algebraic_geometry_v1} or Matsumura \cite[Chapter 5, Theorem 14.2, p. 105]{Matsumura_commutative_ring_theory}, it follows that $\sO_{X,p}$ is a regular local ring of dimension $d$ and thus Item \eqref{item:X_regular_analytic_space_at_point_dimension_d} holds.

Conversely, assume that Item \eqref{item:X_regular_analytic_space_at_point_dimension_d} holds. We shall adapt the proof by  G\"ortz and Wedhorn of \cite[Lemma 6.27]{Gortz_Wedhorn_algebraic_geometry_v1} to show that Item \eqref{item:X_smooth_analytic_space_at_point_dimension_d} holds. Since $\dim T_pX = d$ by \eqref{eq:Dimension_Zariski_tangent_space_equals_dimension_local_ring_analytic_space} and Lemma \ref{lem:Identification_extrinsic_cotangent_space_with_Zariski_cotangent_space} provides that 
\[
  T_pX = \Ker J_{f_1,\ldots,f_r}(p),
\]
then $\dim \Ker J_{f_1,\ldots,f_r}(p) = d$, where the Jacobian matrix is given by
\[
  J_{f_1,\ldots,f_r}(p) = \left(\frac{\partial f_i}{\partial x_j}(p)\right)_{r\times n} \in \Hom_\KK(\KK^n,\KK^r).
\]
Since $\KK^n/\Ker J_{f_1,\ldots,f_r}(p) \cong \Ran J_{f_1,\ldots,f_r}(p) \subseteq \KK^r$, we have
\[
  \Rank_\KK J_{f_1,\ldots,f_r}(p) = n-d \leq r.
\]
Therefore we may assume, if necessary after renumbering the $f_i$ in the presentation $\sO_X = \sO_D/(f_1,\ldots,f_r)$, that the first $n-d$ rows of the Jacobian matrix $J_{f_1,\ldots,f_r}(p)$ are linearly independent. Let
\[
  Y := \supp(\sO_D/(f_1,\ldots,f_{n-d}))
\]
and note that $p\in X \subseteq Y \subset D \subset \KK^n$. For convenience, we define an analytic map by $F = (f_1,\ldots,f_{n-d}):\KK^n \supset D \to \KK^{n-d}$ and observe that $Y = F^{-1}(0)$. The corresponding Jacobian matrix,
\[
  \d F(p) = J_{f_1,\ldots,f_{n-d}}(p) = \left(\frac{\partial f_i}{\partial x_j}(p)\right) \in \Hom_\KK(\KK^n,\KK^{n-d}),
\]
has rank $n-d$ over $\KK$ and thus $Y$ is smooth at $p$ (in the sense of Definition \ref{defn:Smooth_point_analytic_space}) by Lemma \ref{lem:Jacobi_criterion_smooth_point}. By the fact that Item \eqref{item:X_smooth_analytic_space_at_point_dimension_d} $\implies$ Item \eqref{item:X_regular_analytic_space_at_point_dimension_d}, we know that the local ring $\sO_{Y,p}$ is regular of dimension $d$, where $\sO_Y := \sO_D/(f_1,\ldots,f_{n-d})$. The homomorphism of rings
\begin{equation}
  \label{eq:Homomorphism_OYx_to_OXx}
  \iota_p^\sharp:\sO_{Y,p} \to \sO_{X,p}
\end{equation}
induced by the inclusion $\iota:X \hookrightarrow Y$ is then a surjection of $d$-dimensional rings. We claim that the epimorphism \eqref{eq:Homomorphism_OYx_to_OXx} is actually an isomorphism of rings

Indeed\footnote{Following the argument by A. Stonestrom in \url{https://math.stackexchange.com/questions/4208721} for an epimorphism of rings $A\to B$ with the same finite Krull dimension, where $A$ is a domain.}, the epimorphism of rings \eqref{eq:Homomorphism_OYx_to_OXx} induces an isomorphism of rings,
\[
  \sO_{Y,p}/\Ker\iota_p^\sharp \cong \sO_{X,p}.
\]
Therefore, the prime ideals of $\sO_{X,p}$ are in one-to-one correspondence with the prime ideals of $\sO_{Y,p}$ that contain $\Ker\iota_p^\sharp$. Now choose any strictly ascending chain of prime ideals in $\sO_{Y,p}$ that contain $\Ker\iota_p^\sharp$ and is of maximal possible length $d$, the Krull dimension of $\sO_{Y,p}$ as in Definition \ref{defn:Krull_dimension}:
\[
  \Ker\iota_p^\sharp \subseteq \fp_0 < \fp_1 < \cdots < \fp_d = \fm_{Y,p}.
\]
Since $\dim\sO_{X,p} = d = \dim\sO_{Y,p}$, this chain must also be a strictly ascending chain of prime ideals of $\sO_{X,p}$ that has maximal possible length. But $\sO_{Y,p}$ is a regular local ring and is thus an integral domain (see Atiyah and Macdonald \cite[Chapter 11, Lemma 11.23, p. 123]{Atiyah_Macdonald_introduction_commutative_algebra} or Matsumura \cite[Chapter 5, Theorem 14.3, p. 106]{Matsumura_commutative_ring_theory}) and therefore $(0)$ is a prime ideal (see Atiyah and Macdonald \cite[Chapter 1, p. 3]{Atiyah_Macdonald_introduction_commutative_algebra} or Matsumura \cite[Section 1.1, p. 2]{Matsumura_commutative_ring_theory}), so if $\fp_0$ were non-zero, then we could add a link $(0) = \fp_{-1} < \fp_0$ to the chain, contradicting maximality. So $\fp_0$ must be zero, which forces $\Ker\iota_p^\sharp$ to be zero, and hence $\iota_p^\sharp$ is a ring isomorphism, as claimed.

Thus, $(\iota,\iota^\sharp):(X,\sO_X) \to (Y,\sO_Y)$ is a morphism of $\KK$-ringed spaces such that $\iota_p^\sharp:\sO_{Y,p}\to\sO_{X,p}$ is a ring isomorphism. The sheaves $\sO_X$ and $\sO_Y$ are coherent by Remark \ref{rmk:Analytic_model_space_coherence}, so Corollary \ref{cor:Morphism_coherent_ringed_spaces_implies_isomorphism_stalks_at_point} implies that there is an open neighborhood $V\subset Y$ of $p$ such that $(\iota,\iota^\sharp):(V\cap X,\sO_{V\cap X}) \to (V,\sO_Y|_V)$ is an isomorphism of $\KK$-ringed spaces, noting that $\iota^{-1}(V) = V\cap X$. In particular, $(V\cap X,\sO_{V\cap X})$ is isomorphic to $(V,\sO_Y|_V)$ as an analytic space (see Definition \ref{defn:Analytic_space}) and therefore $X$ is smooth at $p$ (in the sense of Definition \ref{defn:Smooth_point_analytic_space}) and so Item \eqref{item:X_smooth_analytic_space_at_point_dimension_d} holds. This completes the proof of Theorem \ref{thm:Comparison_smoothness_regularity_analytic_spaces}.
\end{proof}

\section{Analyticity and dimension of the singular set of an analytic space}
\label{sec:Analyticity_dimension_singular_set}
We begin with the

\begin{defn}[Singular locus of an analytic space]
\label{defn:Singular_locus_analytic_space}
(See Fischer \cite[Section 2.15, p. 96]{Fischer_complex_analytic_geometry} for $\KK=\CC$, Grauert and Remmert \cite[Section 6.2.1, p. 116]{Grauert_Remmert_coherent_analytic_sheaves} for $\KK=\CC$, and Narasimhan \cite[Section 3.1, Definition 2, p. 36]{Narasimhan_introduction_theory_analytic_spaces} for $\KK=\RR$ or $\CC$.)  
Let $(X,\sO_X)$ be an analytic space over $\KK=\RR$ or $\CC$. Then the \emph{singular locus} $X_\sing$ of $X$ is the set of all points $p \in X$ such that $p$ is not a smooth point of $X$ in the sense of Definition \ref{defn:Smooth_point_analytic_space}. (By Theorem \ref{thm:Comparison_smoothness_regularity_analytic_spaces}, this is equivalent to the condition that $p$ is not a regular point of $X$ in the sense of Definition \ref{defn:Regular_point_analytic_space}.) Let $X_\sm := X \less X_\sing$ denote the subset of smooth points.
\end{defn}

For properties of the singular locus of an analytic space, we shall primarily consider the case $\KK=\CC$, which is better behaved.

\begin{thm}[Analyticity of the singular locus of a complex analytic space]
\label{thm:Analyticity_singular_locus_complex_analytic_space}  
(See Abhyankar \cite[Section 33.24, p. 283 and Section 44.20, p. 409]{Abhyankar_local_analytic_geometry},
Demailly \cite[Theorem 4.31, p. 100]{Demailly_complex_analytic_differential_geometry},
Fischer \cite[Section 2.15, Corollary, p. 96]{Fischer_complex_analytic_geometry}, Grauert and Remmert \cite[Section 6.2.2, p. 117]{Grauert_Remmert_coherent_analytic_sheaves}, Guaraldo, Macr\`\i, and Tancredi \cite[Section 2.3, Theorem 3.3, p. 25]{Guaraldo_Macri_Tancredi_topics_real_analytic_spaces}, Narasimhan \cite[Section 3.2, Theorem 6, p. 56 and Corollary 1, p. 58 and Section 4.1, Proposition 2, p. 66]{Narasimhan_introduction_theory_analytic_spaces}.) 
Let $(X,\sO_X)$ be a complex analytic space as in Definition \ref{defn:Analytic_space}. Then the following hold:
\begin{enumerate}
\item\label{item:Singular_locus_analytic_set} The singular locus $X_\sing$ is an analytic set in $X$ in the sense of Definition \ref{defn:Analytic_set}.
\item\label{item:Singular_locus_lower_dimension_X_reduced} If $X$ is reduced, then $X_\sing$ is nowhere dense in $X$ and $\dim_pX_\sing < \dim_pX$ for all points $p \in X$, where $\dim_pX$ is either one of the integers in Lemma \ref{lem:Equivalence_analytic_algebraic_dimensions_complex_analytic_space}.
\end{enumerate}  
\end{thm}

Regarding Item \eqref{item:Singular_locus_lower_dimension_X_reduced} in Theorem \ref{thm:Analyticity_singular_locus_complex_analytic_space}, it is useful to keep the following in mind:

\begin{lem}[Dimension at a point of a complex analytic space and its reduction]
\label{lem:Dimension_point_complex_analytic_space_reduction}
(See Grauert and Remmert \cite[Section 5.1.3, Observation, p. 96]{Grauert_Remmert_coherent_analytic_sheaves}.)
Let $(X,\sO_X)$ be a complex analytic space as in Definition \ref{defn:Analytic_space} and $X_\red = (X,\sO_X^\red)$ be its reduction as in Grauert and Remmert \cite[Section 4.3.2, p. 88]{Grauert_Remmert_coherent_analytic_sheaves}. Then $\dim_p X = \dim_p X_\red$ for all $p \in X$, where $\dim_pX$ is either one of the integers in Lemma \ref{lem:Equivalence_analytic_algebraic_dimensions_complex_analytic_space}.
\end{lem}

The proof of Theorem \ref{thm:Analyticity_singular_locus_complex_analytic_space} due to Grauert and Remmert relies on their \cite[Section 3.3.4, Open Projection Lemma, p. 71 and Section 5.4.4, Proposition, p. 108]{Grauert_Remmert_coherent_analytic_sheaves} and the \emph{Active Lemma} (Lemma \ref{lem:Active}) and its consequences, including Lemma \ref{lem:Ritt}, which we now describe. Grauert and Remmert use \cite[Section 3.3.4, Open Projection Lemma, p. 71]{Grauert_Remmert_coherent_analytic_sheaves} to prove the

\begin{thm}[Existence of finite open holomorphic maps]
\label{thm:Existence_finite_open_holomorphic_maps}  
(See Grauert and Remmert \cite[Section 3.3.4, Existence Theorem, p. 71]{Grauert_Remmert_coherent_analytic_sheaves}.)
Let $(X,\sO_X)$ be a complex analytic space as in Definition \ref{defn:Analytic_space} and that is irreducible at a point $p \in X$. Then there exists an open neighborhood $U \subset X$ of $p$ and a finite and
open holomorphic map $f: U \to V$ onto a connected domain $V \subset \CC^d$.
\end{thm}

\begin{defn}[Active section]
\label{defn:Active_section}
(See Grauert and Remmert \cite[Section 5.2.1, p. 97]{Grauert_Remmert_coherent_analytic_sheaves} or Ebelin \cite[Section 2.10, Definition, p. 104]{Ebeling_functions_several_complex_variables_singularities}.)
Let $(X,\sO_X)$ be a complex analytic space as in Definition \ref{defn:Analytic_space}. A \emph{section} $g \in \sO_X(X)$ is \emph{active} if and only if every \emph{germ} $g_p \in \sO_{X,p}$ is \emph{active} in the sense that for every $f_p \in \sO_{X,p}$ with $f_pg_p \in \sN_{X,p}$ we have $f_p \in \sN_{X,p}$, where $\sN_X \subset \sO_X$ is the nilradical sheaf as in Theorem \ref{thm:Nilradical_analytic_space}.
\end{defn}

\begin{lem}[Criterion of activity]
\label{lem:Criterion_activity}
(See Grauert and Remmert \cite[Section 5.2.2, p. 98]{Grauert_Remmert_coherent_analytic_sheaves}.)  
Given a holomorphic function $g \in \sO_X(X)$ and a point $p$ in its zero set $N(g)$, the following statements are equivalent:
\begin{enumerate}
\item $g$ is active at $p$.
\item There is an open neighborhood $V \subset X$ of $p$ such that $N(g)\cap V$ is nowhere
dense in $V$.
\item There is a decomposition $U = A_1 \cup \cdots \cup A_s$ of an open neighborhood $U\subset X$ into local prime components at $p$ such that $N(g) \cap A_j$ is not a neighborhood of $p$ in $A_j$ for $1 \leq j \leq s$.
\end{enumerate}  
\end{lem}

\begin{lem}[Active Lemma]
\label{lem:Active}
(See Grauert and Remmert \cite[Section 5.2.4, p. 100]{Grauert_Remmert_coherent_analytic_sheaves}, De Jong and Pfister \cite[Corollary 4.1.10, p. 132]{DeJong_Pfister_local_analytic_geometry}, or Ebelin \cite[Proposition 2.46, p. 104]{Ebeling_functions_several_complex_variables_singularities}.)
Let $(X,\sO_X)$ be a complex analytic space as in Definition \ref{defn:Analytic_space}. If $g \in \sO_X(X)$ is active at $p \in N(g)$ in the sense of Definition \ref{defn:Active_section}, then
\[
  \dim_p N(g) = \dim_p X - 1.
\]
\end{lem}

Lemma \ref{lem:Active} is an analytic analogue of Lemma \ref{lem:Atiyah_Macdonald_corollary_11-18} in commutative ring theory and an analytic analogue of Vakil \cite[Exercise 11.2.G, p. 311]{Vakil_foundations_algebraic_geometry} for a scheme of finite type over a field $k$ and a closed irreducible subset. Grauert and Remmert apply Lemma \ref{lem:Active} to prove the

\begin{lem}[Ritt's Lemma]
\label{lem:Ritt}
(See Grauert and Remmert \cite[Section 5.3.1, p. 102]{Grauert_Remmert_coherent_analytic_sheaves}.)
Let $(X,\sO_X)$ be a complex analytic space as in Definition \ref{defn:Analytic_space} and $A$ be an analytic set in $X$. Then the following statements are equivalent:
\begin{enumerate}
\item $\dim_p A < \dim_p X$ for all $p\in A$,
\item $A$ is nowhere dense in $X$.
\end{enumerate}  
\end{lem}

The following result may be derived from Theorem \ref{thm:Analyticity_singular_locus_complex_analytic_space}, but we shall include the separate statement and proof because of its importance in our applications and the fact that its statement does not require $X$ to be reduced as in Item \eqref{item:Singular_locus_lower_dimension_X_reduced} in Theorem \ref{thm:Analyticity_singular_locus_complex_analytic_space}.

\begin{thm}[Krull dimension and nearby smooth points for a complex analytic set]
\label{thm:Narasimhan_section_3-1_theorem_1_p_41}    
(See Demailly \cite[Theorem 4.31 and proof, p. 100]{Demailly_complex_analytic_differential_geometry} or Narasimhan \cite[Section 3.1, Proposition 6, p. 40 and Theorem 1, p. 41]{Narasimhan_introduction_theory_analytic_spaces}.)  
Let $S$ be a complex analytic subset (as in Definition \ref{defn:Analytic_set}) of an open set $\Omega \subset \CC^n$. If $p \in S$ is a point with $\dim\sO_{S,p} = d$ and\footnote{We have added the hypothesis that $d>0$ since the statement becomes vacuous when $d=0$.} $d > 0$, then every open neighborhood of $p$ in $S$ contains points at which $S$ is smooth of dimension $d$ as in Definition \ref{defn:Smooth_point_analytic_set}. In particular, the set of smooth points of $S$ is dense in $S$.
\end{thm}

\begin{proof}
If $p \in S_\sm$, then by Definition \ref{defn:Smooth_point_analytic_set} there is an open neighborhood $U \subset S$ of $p$ such that $U$ is a complex manifold, say of dimension $k$, so $p$ is a smooth point of dimension $k$. Since $p$ is a smooth point, it is also a regular point and $k = d$ by Theorem \ref{thm:Comparison_smoothness_regularity_analytic_spaces}, so the conclusions immediately follow in this case.
  
If $p \in S_\sing$, then the conclusion follows by Demailly's proof \cite[Theorem 4.31, p. 100]{Demailly_complex_analytic_differential_geometry}. By the argument in the proof \cite[Theorem 4.31, p. 100]{Demailly_complex_analytic_differential_geometry}, we assume without loss of generality that $p=0\in\CC^n$ and that $(S,0)$ is an irreducible germ. After possibly shrinking $\Omega$, we may suppose that the sheaf $\sI_S$ has generators $\{g_1,\ldots,g_N\}$. For each point $x \in S_\sm\cap \Omega$, there are an open neighborhood $U\subset\Omega$ of $x$ and $f_j \in \sO_{\CC^n}(U)$ for $j = 1,\ldots,n-d$ such that 
\[
  S\cap U = \{z\in U: f_1(z) = \cdots = f_{n-d}(z) = 0\}
\]
and the differentials $\d f_j$, for $j = 1,\ldots,n-d$, are linearly independent. Since
\[
  f_j \in g_1\sO_U + \cdots g_N\sO_N, \quad\text{for } j = 1,\ldots,n-d,
\]
where $\sO_U := \sO_{\CC^n} \restriction U$, then one can extract a subset of generators $\{g_{j_1},\ldots,g_{j_{n-d}}\}$ that has a Jacobian matrix
\[
  J_{g_{j_1},\ldots,g_{j_{n-d}}}(x)
  =
  \left(\frac{\partial g_{j_a}}{\partial z_k}(x)\right)_{(n-d)\times n}
  \in \Hom_\CC(\CC^n,\CC^{n-d}),
\]
as in \eqref{eq:Jacobian_matrix_analytic_space} of rank $n-d$ at $x$. Therefore, $S_\sing$ is defined by the holomorphic equations on $\Omega$, 
\begin{multline}
  \label{eq:Defining_equations_analytic_space_singular_set}
  g_j = 0, \quad\text{for } j=1,\ldots,N,
  \quad\text{and}
  \\
  \det J_{g_{j_1},\ldots,g_{j_{n-d}}} = 0, \quad\text{for all }
  J \subset \{1,\ldots,N\} \text{ with } |J| = n-d.
\end{multline}
Thus $S_\sm = S \less S_\sing$ and $S_\sing$ is nowhere dense in $S$. (Grauert and Remmert apply Lemma \ref{lem:Criterion_activity} or \ref{lem:Ritt} to reach this conclusion.)
\end{proof}  

Since an analytic space is locally isomorphic to an analytic set, Theorem \ref{thm:Narasimhan_section_3-1_theorem_1_p_41} yields the

\begin{cor}[Algebraic dimension and nearby smooth points for a complex analytic space]
\label{cor:Generic_smoothness_complex_analytic_space}
Let $(X,\sO_X)$ be a complex analytic space as in Definition \ref{defn:Analytic_space}. If point $p \in X$ has $\dim\sO_{X,p} = d$, then every open neighborhood $U\subset X$ of $p$ contains points at which $X$ is smooth of dimension $d$ as in Definition \ref{defn:Smooth_point_analytic_set} and the set of smooth points of $X\cap U$ is dense in $X\cap U$.
\end{cor}

\begin{rmk}[Interpretation of $d$ in Theorem \ref{thm:Narasimhan_section_3-1_theorem_1_p_41}]
\label{rmk::Narasimhan_section_3-1_theorem_1_p_41_meaning_of_d}
Given a non-zero prime ideal $\fp \subsetneq \sO_{\CC^n,0}$, Narasimhan in 
\cite[Section 3.1, Definition 3, p. 41]{Narasimhan_introduction_theory_analytic_spaces} defines $d$ in \cite[Section 3.1, Theorem 1, p. 41]{Narasimhan_introduction_theory_analytic_spaces} via \cite[Section 3.1, Proposition 2, p. 32 and Proposition 3, p. 35]{Narasimhan_introduction_theory_analytic_spaces} to be the largest integer in the range $0 \leq d < n$ such that the ring homomorphism
\[
  \varphi:\sO_{\CC^d,0} \to \sO_{\CC^n,0}
\]
is injective and ensures that $\sO_{\CC^n,0}/\fp$ a finite $\sO_{\CC^d,0}$-module. In particular, $d$ is the Krull dimension \eqref{eq:Krull_dimension} of $\sO_{\CC^n,0}/\fp$ and that explains our interpretation of the dimension of the germ $(S,p)$ in \cite[Section 3.1, Theorem 1, p. 41]{Narasimhan_introduction_theory_analytic_spaces} as $\dim\sO_{S,p}$. This is also the integer $d$ of Theorem \ref{thm:Noether_normalization_analytic_algebras}.
The dimension of an arbitrary complex analytic germ $(S,p)$ is the maximum dimension of the irreducible components $(S_\nu,p)$ of $(S,p)$.
\end{rmk}  

\begin{rmk}[Restriction to $\KK=\CC$ in Theorem \ref{thm:Narasimhan_section_3-1_theorem_1_p_41}]
\label{rmk::Narasimhan_section_3-1_theorem_1_p_41_restriction_to_C}
The proof of Theorem \ref{thm:Narasimhan_section_3-1_theorem_1_p_41} relies in part on the Weierstrass Preparation Theorem (through the proof of Theorem \ref{thm:Noether_normalization_analytic_algebras}), which is valid for $\KK=\RR$ or $\CC$: see Abhyankar \cite[Chapter 2, Section 10.3, p. 74]{Abhyankar_local_analytic_geometry} or Narasimhan \cite[Chapter 2, Theorem 2, p. 15]{Narasimhan_introduction_theory_analytic_spaces} for $\KK=\RR$ or $\CC$ and Chirka \cite[Section 1.1]{Chirka_1989}, Griffiths and Harris \cite[Section 0.1, p. 8]{GriffithsHarris}, Noguchi \cite[Theorem 2.1.3]{Noguchi_analytic_function_theory_several_variables}, Gunning and Rossi \cite[Section 2.B, Theorem 2, p. 68]{Gunning_Rossi_analytic_functions_several_complex_variables}, or H\"ormander \cite[Theorem 6.1.1]{Hormander_introduction_complex_analysis_several_variables} for $\KK=\CC$.

However, while Narasimhan explicitly allows $\KK=\RR$ or $\CC$ in his statement of Theorem \ref{thm:Narasimhan_section_3-1_theorem_1_p_41} (see \cite[Section 3.1, p. 31]{Narasimhan_introduction_theory_analytic_spaces}), Example \ref{exmp:Hypersurface} illustrates that the result is not literally true when $\KK=\RR$, at least within the category of real analytic spaces rather than algebraic geometry. See also the forthcoming Remark \ref{rmk:Analyticity_singular_locus_real_analytic_space}.
\end{rmk}

\begin{rmk}[Definition of dimension of an analytic set in terms of nearby smooth points]
\label{rmk:Equivalent_definition_dimension_analytic_set}
If $X$ is an analytic variety over a field $\KK$ that is irreducible at a point $p$, Guaraldo, Macr\`\i, and Tancredi \cite[Section 1.2, Definition 2.7, p. 22]{Guaraldo_Macri_Tancredi_topics_real_analytic_spaces} define the \emph{dimension of $X$ at $p$} to be the largest $d$ such that every open neighborhood of $p$ contains a point that is smooth of dimension $d$ in the sense of Definition \ref{defn:Smooth_point_analytic_space}. Guaraldo, Macr\`\i, and Tancredi remark \cite[Section 1.2, Remark 2.9, p. 22]{Guaraldo_Macri_Tancredi_topics_real_analytic_spaces} that the integer $d$ is equal to the Krull dimension of the local ring, $\dim \sO_{X,p}$, but they do not provide a proof or a precise reference. When $\KK=\CC$, Theorem \ref{thm:Narasimhan_section_3-1_theorem_1_p_41} confirms that this integer $d$ is equal to the integer $d$ in Theorem \ref{thm:Noether_normalization_analytic_algebras} and thus equal to $\dim \sO_{X,p}$.
\end{rmk}  

We include the following comments and useful basic results for later ease of reference.

\begin{rmk}[Consequences of finite type property and coherence of the structure sheaf $\sO_Y$]
\label{rmk:Consequences_finite_type_property_coherence_structure_sheaf}
Suppose that $(Y,\sO_Y)$ is a local model space as in Definition \ref{defn:Analytic_model_space} over $\KK=\RR$ or $\CC$, so $Y \subset D$ is the support of $\sO_D/\sI$ with a domain $D \subset \KK^n$ and ideal $\sI \subset \sO_D$ with generators $f_1,\ldots,f_r$ and $\sO_Y = (\sO_D/\sI)\restriction Y$. According to Remark \ref{rmk:Nakayama_lemma_minimal_number_generators_local_ring}, for a point $p\in Y$, the \emph{minimal number} $r$ of sections $f_j$ of $\sO_D$ such that their germs $f_{j,p} \in \sO_{D,p}$ generate the ideal $\sI_p$ may be characterized as the dimension of the vector space $\sI_p/(\fm_p\sI_p)$ over the residue field $\kappa_p = \sO_{D,p}/\fm_p \cong \KK$:
\[
  r = \dim_\KK \sI_p/(\fm_p\sI_p).
\]
According to Remark \ref{rmk:Analytic_model_space_coherence}, the subsheaf $\sI \subset \sO_D$ (of ideals of the sheaf of rings of germs of $\KK$-analytic functions on $D$) is a coherent sheaf of $\sO_D$-modules and, in particular, is of finite type in the sense of Definition \ref{defn:Coherent_sheaf} (Item \eqref{item:Finite_type}). Hence, assuming only that the germs $f_{1,p},\ldots,f_{r,p} \in \sO_{D,p}$ generate $\sI_p$, by Grauert and Remmert \cite[Annex 3.1, Lemma, p. 233]{Grauert_Remmert_coherent_analytic_sheaves} there exists an open neighborhood $U\subset D$ of $p$ such the sections $f_1,\ldots,f_r$ of $\sO_D\restriction U$ generate the sheaf of ideals $\sI\restriction U$.

Furthermore, because the following sequence of coherent sheaves is exact,
\[
  (f_{1,p},\ldots,f_{r,p}) \xrightarrow{} \sO_{D,p} \xrightarrow{} \sO_{X,p} = \sO_{D,p}/(f_{1,p},\ldots,f_{r,p})
\]
there is an open neighborhood $U\subset D$ of $p$ such that the following sequence is also exact, 
\[
  (f_1,\ldots,f_r) \xrightarrow{} \sO_U \xrightarrow{} \sO_{Y\cap U} = \sO_U/(f_1,\ldots,f_r),
\]
noting that $\sO_{Y\cap U,p} = \sO_{Y,p}$ (see Grauert and Remmert \cite[Annex 4.2, Consequence 3, p. 237]{Grauert_Remmert_coherent_analytic_sheaves}).
\end{rmk}

\begin{lem}[Monotonicity of dimension for complex analytic sets]
\label{lem:Grauert_Remmert_note_5-1-3}  
(See Grauert and Remmert \cite[Note 5.1.3, p. 96]{Grauert_Remmert_coherent_analytic_sheaves}.)
Let $(X,\sO_X)$ be a complex analytic space as in Definition \ref{defn:Analytic_space}. If $A, B$ are analytic subsets of $X$ with $A \subset B$, then for all $p \in X$,
\[
  \dim_pA \leq \dim_pB.
\]  
\end{lem}

\begin{lem}[Dimension of the product of two complex analytic spaces]
\label{lem:Grauert_Remmert_5-3-1_dimension_product}
(See Grauert and Remmert \cite[Product Formula 5.3.1, p. 101]{Grauert_Remmert_coherent_analytic_sheaves}.)
Let $(X,\sO_X)$ and $(Y,\sO_Y)$ be complex analytic spaces as in Definition \ref{defn:Analytic_space}. If $(p, q) \in X \times Y$, then
\[
  \dim_{(p,q)} X \times Y = \dim_pX + \dim_q Y.
\]
\end{lem}

Lemma \ref{lem:Grauert_Remmert_5-3-1_dimension_product} is the analytic analogue of Vakil \cite[Exercise 11.2.E, p. 311]{Vakil_foundations_algebraic_geometry} for products of schemes of finite type over a field $k$.

\begin{lem}[Dimension of the intersection of two complex analytic sets]
\label{lem:Grauert_Remmert_5-3-1_dimension_intersection}
(See Grauert and Remmert \cite[Intersection Inequality 5.3.1, p. 102]{Grauert_Remmert_coherent_analytic_sheaves}.)
Let $(X,\sO_X)$ be a complex analytic space as in Definition \ref{defn:Analytic_space}. If $A, B$ are analytic subsets of $X$ and $p\in A\cap B$, which is a smooth point of $X$ in the sense of Definition \ref{defn:Smooth_point_analytic_set}, then
\[
  \dim_p (A \cap B) \geq \dim_p A + \dim_p B - \dim_p X.
\]
\end{lem}

Lemma \ref{lem:Grauert_Remmert_5-3-1_dimension_intersection} is the analytic analogue of Shafarevich \cite[Chapter I, Section 6.1, Theorem 1.24, p. 75]{Shafarevich_v1} for the codimension of the intersection (if non-empty) of two irreducible quasiprojective varieties in projective space $\PP^n$ over a field $k$ or Vakil \cite[Exercise 11.3.E, p. 317]{Vakil_foundations_algebraic_geometry} for the case of algebraic varieties in affine space $k^n$ over a field $k$. 

\begin{thm}[Analyticity of the singular locus of a real analytic space]
\label{thm:Analyticity_singular_locus_real_analytic_space}  
(See Guaraldo, Macr\`\i, and Tancredi \cite[Section 4.1, Theorem 1.5, p. 61]{Guaraldo_Macri_Tancredi_topics_real_analytic_spaces} and Narasimhan \cite[Section 5.2, Proposition 18, p. 105]{Narasimhan_introduction_theory_analytic_spaces}.)
Let $(X,\sO_X)$ be a real analytic variety. Then the following hold:
\begin{enumerate}
\item There exists a closed analytic subvariety of $X$, of codimension at least one, which contains the singular locus $X_\sing$.
\item If $X$ is coherent, then $X_\sing$ is a closed analytic subvariety of $X$, of codimension at least one.
\end{enumerate}  
\end{thm}

\begin{rmk}[Counterexamples to Theorem \ref{thm:Analyticity_singular_locus_real_analytic_space} for arbitrary real analytic sets]
\label{rmk:Analyticity_singular_locus_real_analytic_space}  
See Guaraldo, Macr\`\i, and Tancredi \cite[Section 2.3, Remark 3.5, (2) and (3), p. 26 and Section 4.1, Remark 1.6, p. 62]{Guaraldo_Macri_Tancredi_topics_real_analytic_spaces} and Narasimhan \cite[Section 5.2, p. 105 and Section 5.3, pp. 106--109]{Narasimhan_introduction_theory_analytic_spaces} for counterexamples to Theorem \ref{thm:Analyticity_singular_locus_real_analytic_space} for arbitrary real analytic sets.
\end{rmk}




\section{Algebraic, analytic, and expected dimensions for analytic spaces}
\label{sec:Algebraic_analytic_expected_dimensions}
We use Lemma \ref{lem:Krull_dimension_quotient_ring} to show that expected dimension gives a lower bound for the algebraic dimension of a point in an analytic space.

\begin{lem}[Algebraic dimension greater than or equal to expected dimension at a point in an analytic space]
\label{lem:Algebraic_dimension_pX_geq_expdim_pX}  
Let $(X,\sO_X)$ be an analytic space over $\KK=\RR$ or $\CC$ as in Definition \ref{defn:Analytic_space}. If $p \in X$ is a point, then
\begin{equation}
  \label{eq:Algebraic_dimension_pX_geq_expdim_pX}
  \dim\sO_{X,p} \geq \expdim_pX,
\end{equation}
where $\dim\sO_{X,p}$ is the algebraic or Krull dimension \eqref{eq:Krull_dimension} of $X$ at $p$ and $\expdim_pX$ is the expected dimension of $X$ at $p$ as in Remark \ref{rmk:Expected_dimension}.
\end{lem}

\begin{proof}
We may assume without loss of generality that $(X,\sO_X)$ is an analytic model space as in Definition \ref{defn:Analytic_model_space}, where $X \subset D$ and $D \subset \KK^n$ is a domain, $\sI \subset \sO_D$ is an ideal, $X = \cosupp\sI$, and $\sO_X = (\sO_D/\sI)\restriction X$. We write $\sI = (f_1,\ldots,f_r)$, where $f_j:D\to\KK$ are analytic functions for $j=1,\ldots,r$. Observe that $\expdim_pX = n-r$, since $X$ is given by the zero locus of the analytic map $F = (f_1,\ldots,f_r):\KK^n\supset D \to \KK^r$. We have
\[
  \sO_{X,p} = (\sO_D/\sI)_p = \left(\sO_D/(f_1,\ldots,f_r)\right)_p = \sO_{D,p}/(f_{1,p},\ldots,f_{r,p}).
\]
Moreover, $\sO_{D,p}$ is a Noetherian local ring by Theorem \ref{thm:Analytic_algebras_are_Noetherian_local_rings} or the argument described in Remark \ref{rmk:Comparison_dimension_Zariski_tangent_space_and_Krul_dimension_local_ring_at_point}. Therefore, Lemma \ref{lem:Krull_dimension_quotient_ring} gives
\[
  \dim\sO_{X,p} = \dim\sO_{D,p}/(f_{1,p},\ldots,f_{r,p}) \geq \dim\sO_{D,p} - r. 
\]
But $\dim\sO_{D,p} = n$ by Remark \ref{rmk:Example_formal_power_series_regular_local_ring} since $\sO_{D,p}= \KK\{x_1-p_1,\ldots,x_n-p_n\}$, for $p = (p_1,\ldots,p_n)$,
and $\expdim_pX = n-r$, so these observations and the preceding inequality yield the inequality \eqref{eq:Algebraic_dimension_pX_geq_expdim_pX}. 
\end{proof}

Next we apply Lemma \ref{lem:Algebraic_dimension_pX_geq_expdim_pX} to show that expected dimension gives a lower bound for the dimension of an analytic space at any smooth point.

\begin{lem}[Analytic dimension greater than or equal to the expected dimension of a complex analytic space]
\label{lem:Dim_q_Xsmooth_geq_expdim_pX}
Let $(X,\sO_X)$ be a complex analytic space as in Definition \ref{defn:Analytic_space}. If $p \in X$ is a point with $\expdim_pX > 0$, then there is an open neighborhood $U\subset X$ of $p$ such that $X_\sm\cap U$ is non-empty of dimension $\dim (X_\sm\cap U) = \dim\sO_{X,p}$ and
\begin{equation}
  \label{eq:Dim_q_Xsmooth_geq_expdim_pX}
  \dim (X_\sm\cap U) \geq \expdim_pX, 
\end{equation}
where $X_\sm$ is as in Definition \ref{defn:Smooth_point_analytic_space}.
\end{lem}

\begin{proof}
Equation \eqref{eq:Algebraic_dimension_pX_geq_expdim_pX} in Lemma \ref{lem:Algebraic_dimension_pX_geq_expdim_pX} implies that
\[
  \dim\sO_{X,p} \geq \expdim_pX > 0.
\]
Hence, Theorem \ref{thm:Narasimhan_section_3-1_theorem_1_p_41} provides an open neighborhood $U\subset X$ of $p$ such that $X_\sm\cap U$ is non-empty and
\[
  \dim(X_\sm\cap U) = \dim\sO_{X,p}.
\]
The equality \eqref{eq:Dim_q_Xsmooth_geq_expdim_pX} now follows by combining the preceding inequality and equality.
\end{proof}

Lemma \ref{lem:Dim_q_Xsmooth_geq_expdim_pX} may become vacuous in the category of real analytic spaces, as Examples \ref{exmp:Hypersurface} and \ref{exmp:Intersection_hypersurfaces} illustrate. 

\begin{exmp}[Hypersurface]
\label{exmp:Hypersurface}  
Consider the ideal $\sI = (f) \subset \sO_{\KK^2}$, where $f(x,y) = x^2+y^2$ and $\KK=\RR$ or $\CC$. Consider $p = (0,0) \in \KK^2$ and observe that $\fm_p = (x,y) \subset \sO_{\KK^2,p} = \KK\{x,y\}$ is the maximal ideal. According to Theorem \ref{thm:Krull_dimension_polynomial_power_series_rings}, we have $\dim \KK\{x,y\} = 2$. Clearly, $f$ is not a zero divisor in $\fm_p$ (in the sense of rings), so Atiyah and Macdonald \cite[Corollary 11.18, p. 122]{Atiyah_Macdonald_introduction_commutative_algebra} implies that
\[
  \dim \KK\{x,y\}/(x^2+y^2) = \dim \KK\{x,y\} - 1 = 1.
\]
Since $\sO_{X,p} = \sO_{\KK^2,p}/\sI_p$, we obtain $\dim\sO_{X,p} = 1$. Denote $\VV_\KK(f) = \{(x,y) \in \KK^2: f(x,y) = 0\}$ and observe that
\[
  \expdim_p\VV_\KK(f) = 1, \quad\text{for all } p \in \VV_\KK(f),
\]
for $\KK=\RR$ or $\CC$.

If $\KK=\RR$, then $X = \VV_\RR(f) = \{(0,0)\}$ and $X_\sm = \varnothing$, since $\d f(x,y) = (2x,2y)$ and $\d f(0,0) = (0,0)$, thus $(0,0)$ is a singular point of $\VV_\RR(f) \subset \RR^2$. Theorem \ref{thm:Narasimhan_section_3-1_theorem_1_p_41} and Lemma \ref{lem:Dim_q_Xsmooth_geq_expdim_pX} do not hold since $U = \{p\}$ and there are no smooth points $q \in U$.

If $\KK=\CC$, then $X = \VV_\CC(f) = \{(x,\pm ix): x \in \CC\}$ and $X_\sm = \VV_\CC(f)\less\{(0,0\}$, since $(0,0)$ is the only singular point of $\VV_\CC(f) \subset \CC^2$, so $X_\sm$ is a complex manifold of dimension $1$. Theorem \ref{thm:Narasimhan_section_3-1_theorem_1_p_41} and Lemma \ref{lem:Dim_q_Xsmooth_geq_expdim_pX} hold in this case too, since $\dim_q X = 1$ for all $q \in X_\sm$.
\qed
\end{exmp}

Example \ref{exmp:Hypersurface} leads to the following simple generalization.

\begin{exmp}[Complete intersection]
\label{exmp:Intersection_hypersurfaces}
Consider the ideal $\sI = (f_1,f_2) \subset \sO_{\KK^3}$, where $f_1(x,y,z) = x^2-y^2-z^2$ and $f_2(x,y,z) = x$ and $\KK=\RR$ or $\CC$. Consider $p = (0,0,0) \in \KK^3$ and observe that $\fm_p = (x,y,z) \subset \sO_{\KK^3,p} = \KK\{x,y,z\}$ is the maximal ideal. According to Theorem \ref{thm:Krull_dimension_polynomial_power_series_rings}, we have $\dim \KK\{x,y,z\} = 3$. We claim that
\[
  \dim \KK\{x,y,z\}/(x^2-y^2-z^2,x) = 1.
\]
Indeed, the ring epimorphism $\KK\{x,y,z\} \ni f(x,y,z) \mapsto f(0,y,z) \in \KK\{y,z\}$ has kernel $\sJ = x\KK\{x,y,z\} \subset \KK\{x,y,z\}$ and so we obtain a ring isomorphism,
\[
  \KK\{x,y,z\}/(x^2-y^2-z^2,x) \ni [f(x,y,z)] \mapsto [f(0,y,z)] \in \KK\{y,z\}/(y^2+z^2).
\]
Since $\dim \KK\{y,z\}/(y^2+z^2) = 1$ by Example \ref{exmp:Hypersurface}, the claim follows. Since $\sO_{X,p} = \sO_{\KK^3,p}/\sI_p$, we obtain $\dim\sO_{X,p} = 1$. Denote $\VV_\KK(\sI) = \{(x,y,z) \in \KK^3: f_1(x,y,z) = 0 = f_2(x,y,z)\}$ and observe that
\[
  \expdim_p\VV_\KK(\sI) = 1, \quad\text{for all } p \in \VV_\KK(\sI),
\]
for $\KK=\RR$ or $\CC$.

If $\KK=\RR$, then $X = \VV_\RR(\sI) = \{(0,0,0)\}$ and $X_\sm = \varnothing$, since $\d f_1(x,y,z) = (2x,2y,0)$ and $\d f(0,0,0) = (1,0,0)$, thus $(0,0,0)$ is a singular point of $\VV_\RR(\sI) \subset \RR^3$. Theorem \ref{thm:Narasimhan_section_3-1_theorem_1_p_41} and Lemma \ref{lem:Dim_q_Xsmooth_geq_expdim_pX} do not hold in this case since $U = \{p\}$ and there are no smooth points $q \in U$.

If $\KK=\CC$, then $X = \VV_\CC(\sI) = \{(0,y,\pm iy): y \in \CC\}$ and $X_\sm = \VV_\CC(\sI)\less\{(0,0,0\}$, since $(0,0,0)$ is the only singular point of $\VV_\CC(\sI) \subset \CC^3$, so $X_\sm$ is a complex manifold of dimension $1$. Theorem \ref{thm:Narasimhan_section_3-1_theorem_1_p_41} and Lemma \ref{lem:Dim_q_Xsmooth_geq_expdim_pX} hold in this case too, since $\dim_q X = 1$ for all $q \in X_\sm$.
\qed
\end{exmp}

Both Examples \ref{exmp:Hypersurface} and \ref{exmp:Intersection_hypersurfaces} illustrate a basic distinction between real and complex algebraic or analytic sets, due to the failure of the Hilbert Nullstellensatz (see the forthcoming Theorem \ref{thm:Hilbert_Nullstellensatz}) for algebraic varieties over fields that are not algebraically closed the or the R\"uckert Nullstellensatz (see the forthcoming Theorem \ref{thm:Ruckert_Nullstellensatz}) for analytic sets over $\KK$ for $\KK=\RR$ rather than $\KK=\CC$. See Cox, Little, and O'Shea \cite[Section 9.3, p. 491]{Cox_Little_OShea_ideals_varieties_algorithms} for a similar discussion in the context of dimensions of real and complex algebraic varieties.

\section{Real analytic sets and counterexamples}
\label{sec:Real_analytic_sets_counterexamples}
Real analytic sets differ from complex analytic sets in important ways and in this section, we review counterexamples that highlight some important differences. To set these examples in context, we first recall the

\begin{defn}[$C$-analytic set]
\label{defn:C-analytic set}  
(See Narasimhan \cite[Section 5.2, Definition 6, p. 104]{Narasimhan_introduction_theory_analytic_spaces}.)  
Let $S$ be a subset of an open set $U \subset \RR^n$ and view $\RR^n$ as a subset of $\CC^n$. Then $S$ is called a \emph{$C$-analytic set} if there exists an open set $\tilde U \subset \CC^n$ such that $\tilde U\cap\RR^n = U$ and a complex analytic set $\tilde S \subset \tilde U$ such that $\tilde S \cap \CC^n = S$.
\end{defn}

\begin{prop}[Conditions for a subset of an open subset of Euclidean space to be $C$-analytic]
\label{prop:_Narasimhan_section_5-2_proposition_15}  
(See Narasimhan \cite[Section 5.2, Proposition 15, p. 104]{Narasimhan_introduction_theory_analytic_spaces} or Cartan \cite[Section 10, Proposition 15, p. 96]{Cartan_1957}.)  
Let $S$ be a subset of an open set $U \subset \RR^n$. Then $S$ is $C$-analytic if and only if one of the following conditions is obeyed:
\begin{enumerate}
\item There are finitely many real analytic functions $f_1,\ldots,f_m \in \sO_U(U)$ such that $S = \{x \in U: f_1(x) = \cdots = f_m(x) = 0\}$.
\item There is a coherent sheaf of ideals $\sI \subset \sO_U$ such that $S = \supp(\sO_U/\sI)$, that is, $S$ is the set of zeros of a coherent sheaf of ideals.
\end{enumerate}
\end{prop}

Note that a coherent analytic set is $C$-analytic, but the converse is not true in general (see Narasimhan \cite[Section 5.2, p. 104]{Narasimhan_introduction_theory_analytic_spaces}). We now review the examples.

\begin{figure}
	\centering
	\includegraphics[width=0.7\linewidth]{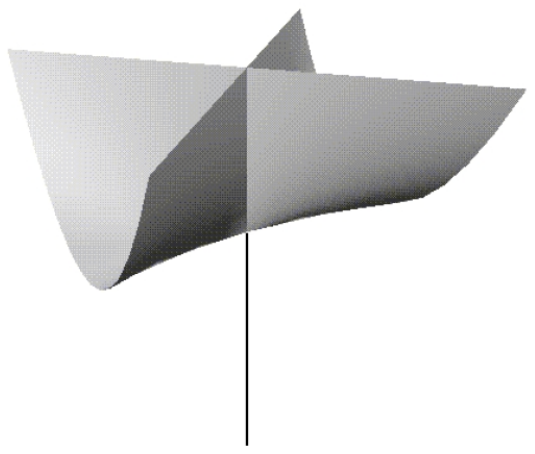}
	\caption[Whitney umbrella]{Whitney umbrella as the zero locus of all points $(x,y,z)$ in $\RR^3$ such that $y^2-zx^2=0$ (from Massey and Tr\'ang \cite[Example 5.21]{Massey_Le_2007})}
	\label{fig:Whitney_umbrella}
      \end{figure}
      
\begin{exmp}[Whitney Umbrella]
\label{exmp:Whitney_umbrella} 
(See Acquistapace, Broglia, and Fernando \cite[Section 2.A.1, Example 2.9 (1), p. 25]{Acquistapace_Broglia_Fernando_topics_global_real_analytic_geometry}, Guaraldo, Macr\`\i, and Tancredi \cite[Section 2.3, Remark 3.5 (2), p. 26]{Guaraldo_Macri_Tancredi_topics_real_analytic_spaces}, and Massey and Tr\'ang \cite[Example 5.21]{Massey_Le_2007}.)
For the subvariety $S$ defined in \eqref{eq:Whitney_umbrella}, the singular locus of $S$ is the set $\{(x_1,x_2,x_3) \in \RR^3: x_1=0, x_2=0, x_3 \leq 0\}$, which is not a subvariety of $S$. Moreover, the set $S$ is not coherent at the point $(0,0,0)$ in the sense of Definition \ref{defn:Coherent_real_analytic_set} (see Guaraldo, Macr\`\i, and Tancredi \cite[Section 2.3, Remark 3.5 (2), p. 26]{Guaraldo_Macri_Tancredi_topics_real_analytic_spaces}). Indeed, Acquistapace, Broglia, and Fernando observe that $S$ is irreducible but not equidimensional\footnote{Acquistapace, Broglia, and Fernando use the term equidimensional whereas Grauert and Remmert \cite[Section 5.4.2, p. 106]{Grauert_Remmert_coherent_analytic_sheaves} use the term \emph{pure dimensional} for a complex analytic space $X$ to mean that $\dim_pX = \dim X$ for all $p \in X$; they say that a complex analytic space $X$ is \emph{pure dimensional at a point $p\in X$} if there exists an open neighborhood $U\subset X$ of $p$ such that $U$ is pure dimensional.}, so $S$ cannot be a coherent real analytic set by Narasimhan \cite[Section 5.1, Proposition 7, p. 95]{Narasimhan_introduction_theory_analytic_spaces}. See Figure \ref{fig:Whitney_umbrella}.

Following Massey and Tr\'ang \cite[Example 5.21]{Massey_Le_2007}, we note the two-dimensional portion above the $x_1x_2$-plane, and that when $x_3< 0$, the only points of $S$ are on the
$x_3$-axis. While $S$ is the union of the one-dimensional $x_3$-axis (the ``handle of the umbrella'' or ``tail'') and the two-dimensional ``umbrella'' portion, the set $S \cap \{(x_1,x_2,x_3): x_3 \geq 0\}$ is not an analytic set and $S$ is irreducible; see Massey and Tr\'ang \cite[Corollary 5.17]{Massey_Le_2007}.

However, $S$ is $C$-analytic by Definition \ref{defn:C-analytic set}. In contrast, a complex analytic space $X$ is pure dimensional at all points where $X$ is irreducible (see Grauert and Remmert \cite[Section 5.4.2, p. 106]{Grauert_Remmert_coherent_analytic_sheaves}).
\qed
\end{exmp}

\begin{rmk}[Resolution of singularities for the Whitney Umbrella]
\label{rmk:Whitney_umbrella_resolution}  
Resolution of singularities for the Whitney Umbrella is described by Annala \cite[Example 3.6]{Annala_resolution_of_singularities} over a field $\KK$ of characteristic zero, in detail by Faber and Hauser \cite[pp. 382--385]{Faber_Hauser_2010} over $\CC$, briefly by Koll\'ar \cite[Example 3.6.1, p. 123]{Kollar_lectures_resolution_singularities} over $\KK$, and in the Mathematics Stack Exchange answer \cite{MathStackExchange_Blowing-up-the-whitney-umbrella-over-the-z-axis} over $\KK$.
\qed
\end{rmk}
 
\begin{figure}
	\centering
	\includegraphics[width=0.7\linewidth]{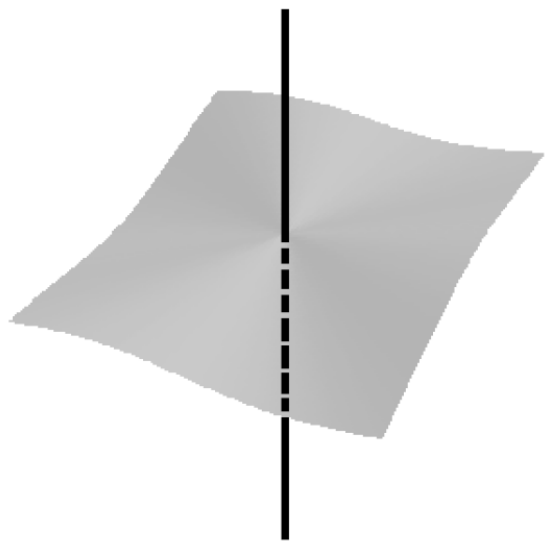}
	\caption[Cartan umbrella]{Cartan umbrella as the zero locus of all points $(x,y,z)$ in $\RR^3$ such that $z(x^2+y^2)-x^3=0$ (from Massey and Tr\'ang \cite[Example 5.26]{Massey_Le_2007})}
	\label{fig:Cartan_umbrella}
      \end{figure}
      
\begin{exmp}[Cartan's Umbrella]
\label{exmp:Cartan_umbrella} 
(See Acquistapace, Broglia, and Fernando \cite[Section 2.A.1, Example 2.9 (2), p. 25]{Acquistapace_Broglia_Fernando_topics_global_real_analytic_geometry}, Bochnak, Coste, and Roy \cite[Section 3.1, Example 3.1.2 (d), p. 60]{Bochnak_Coste_Roy_real_algebraic_geometry}, Guaraldo, Macr\`\i, and Tancredi \cite[Section 2.1, Remark 1.3, p. 12]{Guaraldo_Macri_Tancredi_topics_real_analytic_spaces}, Massey and Tr\'ang \cite[Example 5.26]{Massey_Le_2007}, and Narasimhan \cite[Section 5.3, Example 1, p. 106]{Narasimhan_introduction_theory_analytic_spaces}.)
For the real analytic set $S\subset \RR^3$ defined in \eqref{eq:Cartan_umbrella}, the sheaf of ideals $\iota(S) \subset \sO_{\RR^3}$ is generated by at the origin by the function $g(x_1,x_2,x_3) := x_3(x_1^2+x_2^2)-x_1^3$. Guaraldo, Macr\`\i, and Tancredi observe that this may be proved by using complexification (see Section \ref{sec:Complexification_real_analytic_spaces}) and applying the R\"uckert Nullstellensatz (see the forthcoming Theorem \ref{thm:Ruckert_Nullstellensatz}) for complex analytic spaces to the ideal of $\sO_{\CC^3}$ generated by the holomorphic function $\tilde g(z_1,z_2,z_3) := z_3(z_1^2+z_2^2)-z_1^3$.

In an open neighborhood of any point $(0,0,t) \in \RR^3$ with $t\neq 0$, the set $S$ reduces to the line $x_1=x_2=0$ (the one-dimensional ``handle'' or ``tail'') and $\iota(S)$ is generated at such a point by the functions $x_1,x_2$. See Figure \ref{fig:Cartan_umbrella}. Therefore, in an open neighborhood of the origin, $\iota(S)$ cannot be generated by $g$ and thus $\iota(S)$ is not of finite type by Guaraldo, Macr\`\i, and Tancredi \cite[Section 1.2, Proposition 2.3, p. 5]{Guaraldo_Macri_Tancredi_topics_real_analytic_spaces} or Grauert and Remmert \cite[Annex 3.1, Lemma, p. 233]{Grauert_Remmert_coherent_analytic_sheaves} (a result due to Serre \cite{Serre_1955}) and thus not coherent by Definition \ref{defn:Coherent_sheaf}.

The real analytic set $S\subset \RR^3$ is irreducible at $(0,0,0)$ and it has dimension two at this point; nevertheless, in every open neighborhood of $(0,0,0) \in \RR^3$ there exist points at which $S$ has dimension one (see \cite[Section 2.3, Remark 3.5 (1), p. 25]{Guaraldo_Macri_Tancredi_topics_real_analytic_spaces}) and \cite[Section 5.3, Example 1, p. 106]{Narasimhan_introduction_theory_analytic_spaces}). Narasimhan observes \cite[Section 5.3, Example 1, p. 106]{Narasimhan_introduction_theory_analytic_spaces} that because $S$ is not pure dimensional at $(0,0,0)$, then $S$ cannot be a coherent real analytic set by \cite[Section 5.1, Proposition 7, p. 95]{Narasimhan_introduction_theory_analytic_spaces}. Moreover, $S$ is irreducible but $S \less S_\sing$ is not connected (see Guaraldo, Macr\`\i, and Tancredi \cite[Section 2.3, Remark 3.5 (1), p. 25]{Guaraldo_Macri_Tancredi_topics_real_analytic_spaces}). However, $S$ is $C$-analytic by Definition \ref{defn:C-analytic set}.
\qed
\end{exmp}

\begin{exmp}
\label{exmp:Quartic_not_coherent} 
(See Acquistapace, Broglia, and Fernando \cite[Section 2.A.1, Example 2.9 (3), p. 25]{Acquistapace_Broglia_Fernando_topics_global_real_analytic_geometry}.)
Consider the set
\begin{equation}
  \label{eq:Quartic}
  S := \left\{(x_1,x_2,x_3)\in\RR^3: x_3(x_1+x_2)(x_1^2+x_2^2) - x_1^4 = 0\right\}.
\end{equation}
Then $S$ has dimension two everywhere, but one can show via Proposition \ref{prop:Acquistapace_Broglia_Fernando_2-6} that the ideal sheaf $\iota(S)$ is not coherent at  $(0,0,0)$.
\qed
\end{exmp}

\begin{exmp}[A closed analytic set in $\RR^3$ which does not admit any coherent structure]
\label{exmp:Closed_analytic_subset_3-space_not_admitting_a_coherent_structure} 
(See Acquistapace, Broglia, and Fernando \cite[Section 2.D, Example 2.26 (1), p. 37]{Acquistapace_Broglia_Fernando_topics_global_real_analytic_geometry} and Cartan \cite[Section 11, pp. 97--98]{Cartan_1957}.)  
Define
\[
  a(t)
  :=
  \begin{cases}
    \exp(1/(t^2-1)), &\text{for } -1 < t < 1,
    \\
    0, &\text{for } t \leq -1 \text{ or } t \geq 1.
  \end{cases}
\]
Let $S \subset \RR^3$ be the zero locus of the smooth function
\[
  g(x_1,x_2,x_3) := a(x_3)x_1^3 - x_3(x_1^2 + x_2^2).
\]
One can verify that $S$ is an analytic set in $\RR^3$ and that if $f \in \sO_{\RR^3}(\RR^3)$, then $f \equiv 0$ on $\RR^3$. Consequently, $S$ does not admit any coherent structure.
\qed
\end{exmp}

\begin{exmp}[A closed analytic set in $\RR^3$ which does not admit any coherent structure and which has a compact two-dimensional part]
\label{exmp:Closed_analytic_subset_3-space_not_admitting_a_coherent_structure_compact_2D} 
(See Acquistapace, Broglia, and Fernando \cite[Section 2.D, Example 2.26 (2), p. 37]{Acquistapace_Broglia_Fernando_topics_global_real_analytic_geometry} and Cartan \cite[Section 11, p. 99]{Cartan_1957}.)
If in Example \ref{exmp:Closed_analytic_subset_3-space_not_admitting_a_coherent_structure}, one takes the function
\[
  g(x_1,x_2,x_3) := x_3^2(1-2x_3^2)(x_1^2+x_2^2)-(x_1^4+x_2^4)a(x_3)
\]
instead of $g(x_1,x_2,x_3) = a(x_3)x_1^3 - x_3(x_1^2 + x_2^2)$, then the zero locus $S$ is an analytic set whose two-dimensional part is compact. As in Example \ref{exmp:Closed_analytic_subset_3-space_not_admitting_a_coherent_structure}, any analytic function vanishing on $S$ necessarily vanishes on $\RR^3$.
\qed
\end{exmp}

\begin{exmp}[A closed analytic set in $\RR^3$ which does not admit any coherent structure and which is compact]
\label{exmp:Closed_analytic_subset_3-space_not_admitting_a_coherent_structure_compact} 
(See Acquistapace, Broglia, and Fernando \cite[Section 2.D, Example 2.26 (2), p. 37]{Acquistapace_Broglia_Fernando_topics_global_real_analytic_geometry}.)
Furthermore, if instead of $g(x) = a(x_3)x_1^3 - x_3(x_1^2 + x_2^2)$ in Example \ref{exmp:Closed_analytic_subset_3-space_not_admitting_a_coherent_structure}, one takes the following function
\[
  g(x_1,x_2,x_3) := \left(1-4(x_1^2 + x_2^2 + x_3^2)\right)
  \left((x_1^2+x_3^2-1)^2 + x_2^2\right) - \left((x_1^2+x_3^2-1)^4 + x_2^4\right)a(x_3),
\]
one obtains a compact analytic set $S$ as the zero locus of $g$ and which has the same properties as those in Examples \ref{exmp:Closed_analytic_subset_3-space_not_admitting_a_coherent_structure} or \ref{exmp:Closed_analytic_subset_3-space_not_admitting_a_coherent_structure_compact_2D}.
\qed
\end{exmp}

\chapter{Blowups and their properties for schemes}
\label{chap:Blowups_properties_schemes}
In this chapter on blowups we provide the necessary background for our exposition of resolution of singularities for schemes and algebraic varieties and for our review of results on the Birula--Bia{\l}ynicki decomposition for algebraic varieties. In Section \ref{sec:Subschemes_and_algebraic_subvarieties}, we review definitions of subschemes and algebraic subvarieties. Section \ref{sec:Fiber_products_base_change} provides a brief introduction to the concepts of fiber product and base change. Section \ref{sec:Blowing_up_algebraic_varieties_schemes} contains a discussion of Cartier divisors, definitions and constructions of blowups of algebraic varieties and schemes, properties of blowups, and strict and total transforms. We conclude in Section \ref{sec:Functorial_properties_blowups_schemes} with a review of the concept of flatness in the category of schemes and discuss functorial property of blowups of schemes.

\section{Subschemes and algebraic subvarieties}
\label{sec:Subschemes_and_algebraic_subvarieties}
We begin by stating a few definitions that we shall need prior to our discussion of blowups.

\begin{defn}[Open subscheme and open immersion (or embedding)]
\label{defn:Open_subscheme_and_open_immersion_or_embedding}  
(See G\"ortz and Wedhorn \cite[Proposition--Definition 3.2 and Definition 3.40]{Gortz_Wedhorn_algebraic_geometry_v1}; see also Eisenbud and Harris \cite[Section I.2.1, p. 23]{Eisenbud_Harris_geometry_schemes}, Hartshorne \cite[Chapter II, Section 3, Definition, p. 85]{Hartshorne_algebraic_geometry}, Hauser \cite[Remark 2.19]{Hauser_2014}, the Stacks Project \cite[\href{https://stacks.math.columbia.edu/tag/01IO}{Definition 01IO}]{stacks-project}, and Vakil \cite[Section 7.1.1, p. 202]{Vakil_foundations_algebraic_geometry}.)
Let $X$ be a scheme. If $U \subseteq X$ is an open subset, then the locally ringed space $(U,\sO_{X|U})$ is a scheme and $U$ is called an \emph{open subscheme} of $X$. If $U$ is an  affine scheme, then $U$ is called an affine open subscheme. A morphism $\iota:Y \to X$ of schemes is called an \emph{open immersion}\footnote{Vakil uses the term \emph{open embedding} and \emph{open immersion} interchangeably \cite[Section 7.1.1, p. 202]{Vakil_foundations_algebraic_geometry}.}  if the underlying continuous map is a homeomorphism of $Y$ onto an open subset $U$ of $X$ and the sheaf homomorphism $\sO_X \to \iota_*\sO_Y$ induces an isomorphism $\sO_{X|U} \cong \iota_*\sO_Y$ of sheaves on $U$.
\end{defn}

\begin{defn}[Closed immersion (or embedding) of schemes and closed subscheme]
\label{defn:Closed_subscheme_and_closed_immersion_or_embedding}  
(See G\"ortz and Wedhorn \cite[Definition 3.41]{Gortz_Wedhorn_algebraic_geometry_v1}; see also Eisenbud and Harris \cite[Definition I.27]{Eisenbud_Harris_geometry_schemes}, Hartshorne \cite[Chapter II, Section 3, Definition, p. 85]{Hartshorne_algebraic_geometry}, Hauser \cite[Definition 2.18]{Hauser_2014}, Shafarevich \cite[Section 5.3.3, Definition 5.5, p. 32]{Shafarevich_v2}, the Stacks Project \cite[\href{https://stacks.math.columbia.edu/tag/01IO}{Definition 01IO}]{stacks-project}, and Vakil \cite[Definition 8.1.1]{Vakil_foundations_algebraic_geometry}.)
Let $X$ be a scheme.
\begin{enumerate}
\item A \emph{closed subscheme} of $X$ is given by a closed subset $Z \subseteq X$ and an ideal sheaf $\sJ \subseteq \sO_X$ such that $Z = \{x \in X : (\sO_X/\sJ)_x \neq 0\}$ and $(Z, (\sO/\sJ)\restriction Z)$ is a scheme.
\item A morphism $\iota : Z \to X$ of schemes is called a \emph{closed immersion}\footnote{Vakil prefers the term \emph{closed embedding}, noting that differential geometric notion of closed immersion is closer to what algebraic geometers call an \emph{unramified morphism}, which Vakil defines in \cite[Section 21.6]{Vakil_foundations_algebraic_geometry}.} if the underlying continuous map is a homeomorphism between $Z$ and a closed subset of $X$, and the sheaf homomorphism $\iota^\flat:\sO_X \to \iota_*\sO_Z$ is surjective, where $\iota^\flat$ denotes the canonical projection $\sO_X \to \sO_X/\sJ = \iota_*\sO_Z$.
\end{enumerate}
\end{defn}

G\"ortz and Wedhorn remark \cite[Section 3.16, p. 87]{Gortz_Wedhorn_algebraic_geometry_v1} that open and closed subschemes are special cases of the notion of (locally closed) subscheme, though Eisenbud and Harris, Hartshorne, Shafarevich, the Stacks Project, and Vakil only define subschemes that are open or closed subschemes in the sense of Definitions \ref{defn:Open_subscheme_and_open_immersion_or_embedding} or \ref{defn:Closed_subscheme_and_closed_immersion_or_embedding}.

\begin{defn}[(Locally closed) subscheme]
\label{defn:Locally_closed_subscheme}  
(See G\"ortz and Wedhorn \cite[Definition 3.43]{Gortz_Wedhorn_algebraic_geometry_v1}.)
\begin{enumerate}
\item Let $X$ be a scheme. A \emph{subscheme} of $X$ is a scheme $(Y, \sO_Y)$ such that $Y \subset X$ is a locally closed subset and such that $Y$ is a closed subscheme of the open subscheme $U \subseteq X$, where $U$ is the largest open subset of $X$ which contains $Y$ and in which $Y$ is
closed (that is, $U$ is the complement of $\bar Y \less Y$). There is a natural morphism of
schemes $Y \to X$.
\item An \emph{immersion} $\iota: Y \to X$ is a morphism of schemes whose underlying continuous map is a homeomorphism of $Y$ onto a locally closed subset of $X$ and such that for all $y \in Y$ the ring homomorphism $\iota_y^\sharp:\sO_{X,\iota(y)} \to \sO_{Y,y}$ between the local rings is surjective.
\end{enumerate}
\end{defn}

\begin{defn}[Algebraic variety]
\label{defn:Algebraic_variety}  
(See the Stacks Project \cite[\href{https://stacks.math.columbia.edu/tag/020D}{Definition 020D}]{stacks-project} and Vakil \cite[Definition 10.1.7]{Vakil_foundations_algebraic_geometry}; compare Hartshorne \cite[Chapter II, Section 4, Definition, p. 105]{Hartshorne_algebraic_geometry}.)  
An \emph{algebraic variety} over a field $\KK$ is a reduced, separated scheme of finite type over $\KK$.
\end{defn}

A scheme $X$ is \emph{reduced} if it satisfies Definition \ref{defn:Reduced_locally_ringed_space} (see G\"ortz and Wedhorn \cite[Chapter 3, Definition 3.26(a), p. 79]{Gortz_Wedhorn_algebraic_geometry_v1}). We refer the reader to G\"ortz and Wedhorn \cite[Chapter 9, Definition 9.7, p. 233]{Gortz_Wedhorn_algebraic_geometry_v1} or Hartshorne \cite[Chapter II, Section 4, p. 96]{Hartshorne_algebraic_geometry} for the definition of a \emph{separated scheme}.

G\"ortz and Wedhorn \cite[Example 9.17, p. 236]{Gortz_Wedhorn_algebraic_geometry_v1} and Hartshorne \cite[Chapter II, Section 4, Definition, p. 105]{Hartshorne_algebraic_geometry} define an \emph{(abstract) variety} to be an integral, separated scheme of finite type over an algebraically closed field $\KK$. Varieties (in the sense of Definition \ref{defn:Algebraic_variety}) form a category: morphisms of varieties are just morphisms as schemes (see Vakil \cite[Section 10.1.7]{Vakil_foundations_algebraic_geometry}). We recall that a scheme $X$ is \emph{integral} if it is reduced and irreducible (see G\"ortz and Wedhorn \cite[Chapter 3, Definition 3.26(b), p. 79]{Gortz_Wedhorn_algebraic_geometry_v1}). 

\begin{defn}[(Algebraic) subvariety]
\label{defn:Algebraic_subvariety}    
(See the Stacks Project \cite[\href{https://stacks.math.columbia.edu/tag/0AZ9}{Section 0AZ9}]{stacks-project} and Vakil \cite{Vakil_foundations_algebraic_geometry}.)
A \emph{subvariety} of a variety $X$ is a reduced, locally closed subscheme of $X$ (which
one can show is itself a variety). An \emph{open subvariety} of $X$ is an open subscheme of $X$. (One can show that reducedness is automatic in this case.) A \emph{closed subvariety} of $X$ is a reduced, closed subscheme of $X$.  
\end{defn}

Classical definitions of algebraic varieties that are \emph{not} based on the definition of a variety as a scheme are provided by Hartshorne (see \cite[Section 1.1, Definition, p. 3]{Hartshorne_algebraic_geometry} for an affine algebraic variety and \cite[Section 1.2, Definition, p. 10]{Hartshorne_algebraic_geometry} for a projective algebraic variety and quasiprojective variety) and by Shafarevich (see \cite[Section 1.2.1, p. 23]{Shafarevich_v1} for a closed subset of affine space, $\KK^n$, and \cite[Section 1.4.1, p. 41]{Shafarevich_v1} for a closed subset of projective space, $\PP^n$, and \cite[Section 1.4.1, p. 46]{Shafarevich_v1} for a quasiprojective variety).

According to Hartshorne \cite[Chapter II, Section 2, Proposition 2.6, p. 78]{Hartshorne_algebraic_geometry} and Shafarevich \cite[Section 5.3.1, Example 5.19, p. 29]{Shafarevich_v2}, to every affine, projective, or quasiprojective variety, one may associate a scheme; see also Hartshorne \cite[Chapter II, Section 3, Example 3.2.1, p. 84]{Hartshorne_algebraic_geometry}.

Definitions of (closed) subvarieties that are not based on the definition of a variety as a scheme are provided by Hartshorne \cite[Chapter I, Section 3, Exercise 3.10, p. 21]{Hartshorne_algebraic_geometry}, Milne \cite[Section 3h, p. 68]{Milne_algebraic_geometry} (closed subvariety of an affine algebraic variety) and Shafarevich \cite[Section 1.4.1, p. 46]{Shafarevich_v1} (subvariety of a quasiprojective variety).

\section{Fiber products and base change}
\label{sec:Fiber_products_base_change}
For further explanation of the concepts of \emph{fiber products} and \emph{base change} in algebraic geometry, we refer to G\"ortz and Wedhorn \cite[Chapter 4, p. 107]{Gortz_Wedhorn_algebraic_geometry_v1}, the Stacks Project \cite[\href{https://stacks.math.columbia.edu/tag/01JW}{Sections 01JW} and \href{https://stacks.math.columbia.edu/tag/02WE}{02WE}]{stacks-project}, Vakil \cite[Sections 9.3 and 9.4]{Vakil_foundations_algebraic_geometry}, and especially Milne \cite[Section 5.i]{Milne_algebraic_geometry}, whose exposition in the category of algebraic varieties we closely follow here (see G\"ortz and Wedhorn \cite{Gortz_Wedhorn_algebraic_geometry_v1} or Vakil \cite{Vakil_foundations_algebraic_geometry} in the category of schemes). For the category of complex analytic spaces, we refer to Fischer \cite{Fischer_complex_analytic_geometry} and for techniques to adapt constructions in the category of complex analytic spaces to the category of real analytic spaces, we refer Guaraldo, Macr\`\i, and Tancredi \cite{Guaraldo_Macri_Tancredi_topics_real_analytic_spaces}. We refer to Milne \cite[Section 5.33]{Milne_algebraic_geometry} for comments regarding differences in fiber products that depend on whether is working in the category of algebraic varieties or algebraic schemes.

Following Milne \cite[Section 5.i]{Milne_algebraic_geometry}, let $\varphi:V\to S$ and $\psi:W\to S$ be regular maps of algebraic varieties over a field $\KK$. The set
\[
  V\times_S W := \left\{ (v,w) \in V\times W: \varphi(v) = \psi(w) \right\}
\]
is closed in $V \times W$, having a canonical structure of an algebraic variety (see Milne \cite[Section 5.e]{Milne_algebraic_geometry}). The algebraic variety $V\times_S W$ is called the \emph{fiber product} of $V$ and $W$ over $S$. Note that if $S$ consists of a single point, then $V\times_S W = V\times W$.

We write $\varphi'$ for the map $V\times_S W \ni (v,w) \mapsto w \in W$ and $\psi'$ for the map $V\times_S W \ni (v,w) \mapsto v \in V$. One then has a commutative diagram:
\[
  \begin{tikzcd}
  V \times_S W \arrow[r, "\varphi' "] \arrow[d, "\psi' "] & W \arrow[d, "\psi"] \\
  V \arrow[r, "\varphi"] &S
\end{tikzcd}
\]
The system $(V \times_S W,\varphi',\psi')$ has the following universal property: for any regular maps $\alpha:T\to V$ and $\beta:T\to W$ such that $\varphi\alpha = \psi\beta$, there is a unique regular map $(\alpha,\beta):T\to V \times_S W$ such that the following diagram commutes:
\[
  \begin{tikzcd}
  T
  \arrow[drr, bend left, "\beta"]
  \arrow[ddr, bend right, "\alpha"]
  \arrow[dr, dotted, "{(\alpha,\beta)}" description] & & \\
    & V \times_S W \arrow[r, "\varphi' "] \arrow[d, "\psi' "]
      & W \arrow[d, "\psi"] \\
& V \arrow[r, "\varphi"] &S
\end{tikzcd}
\]
In other words,
\[
  \Hom(T,V\times_SW) \cong \Hom(T,V)\times_{\Hom(T,S)}\Hom(T,W).
\]
Indeed, there is a unique such map of sets, namely, $T \ni t \mapsto (\alpha(t),\beta(t))$, which is regular because it is regular as a map into $V\times W$. 

The map $\varphi'$ in the preceding diagrams is called the \emph{base change} of $\varphi$ with respect to $\psi$. For any point $s \in S$, the base change of $\varphi:V\to S$ with respect to $s \hookrightarrow S$ is the map $\varphi^{-1}(s) \to s$ induced by $\varphi$, namely the fibre of $V$ over $s$. If $f:V\to S$ is a regular map and $U \subset S$ is a subvariety, then $V\times_SU$ is the inverse image $f^{-1}(U)$ of $U$ in $V$ (see Milne \cite[Example 5.31]{Milne_algebraic_geometry}).

For fiber products in the category of complex analytic spaces, see Section \ref{sec:Fiber_products_analytic_spaces}.

\section{Blowups of algebraic varieties and schemes}
\label{sec:Blowing_up_algebraic_varieties_schemes}
In Section \ref{subsec:Existence_uniqueness_blowups}, we review results for the existence and uniqueness of the blowup of a scheme along a closed subscheme. Section \ref{subsec:Strict_transforms_subschemes_blowups} provides definitions of the strict transform of a subscheme in the blowup of a scheme. We conclude in Section \ref{subsec:Smoothness_blowups_strict_transforms} with a summary of results concerning the smoothness of blowups and strict transforms.

\subsection{Existence and uniqueness of blowups}
\label{subsec:Existence_uniqueness_blowups}
We first review the definition and properties of Cartier divisors in the category of algebraic varieties before proceeding to their definition in the category of schemes.

\begin{defn}[Cartier divisor in an algebraic variety]
\label{defn:Hauser_4-3}
(See Hauser \cite[Definition 4.1, p. 17]{Hauser_2014}; see also Shafarevich \cite[Section 3.1.2, Definition, p. 153]{Shafarevich_v1}.)
A closed subvariety $Z$ of an irreducible variety $X$ is a \emph{Cartier divisor in $X$ at a point $p \in Z$} if $Z$ can be defined locally at $p$ by a single equation $h = 0$ for some non-zero element $h \in \sO_{X,p}$. If $X$ is not assumed to be irreducible, $h$ is required to be a non-zero divisor of $\sO_{X,p}$ (this excludes the possibility that $Z$ is a component or a union of components of $X$). The subvariety $Z$ is a \emph{Cartier divisor in $X$} if it is a Cartier divisor at each of its points. The empty subvariety is considered as a Cartier divisor. 
\end{defn}

\begin{rmk}[Properties of Cartier divisors in algebraic varieties]
\label{rmk:Properties_Cartier_divisors_algebraic_varieties}
(See Hauser \cite[Lecture IV, p. 17]{Hauser_2014}.)  
A (nonempty) Cartier divisor $Z\subset X$ has codimension one, that is, $\dim_pZ = \dim_pX - 1$ for all $p\in Z$, but the converse is false (see Fischer \cite[Section 0.45, Example, p. 42]{Fischer_complex_analytic_geometry}). When $X$ is non-singular and irreducible, a closed subvariety $Y\subset X$ with codimension one is locally defined at any point $p\in Y$ by a single non-trivial equation, that is, an equation given by a non-zero and non-invertible element $h \in \sO_{X,p}$. This need not hold for singular varieties (see Hauser \cite[Example 4.21, p. 22]{Hauser_2014}).

The complement $X \less Z$ of a Cartier divisor $Z\subset X$ is dense in $X$. Cartier divisors are, essentially, the largest closed and properly contained subvarieties of $X$. If $Z$ is a Cartier divisor in $X$, then the ideal $\sI$ defining $Z$ in $X$ is called \emph{locally principal} (see Hartshorne \cite[Chapter II, Section 6, Proposition 6.13, p. 144]{Hartshorne_algebraic_geometry} or Eisenbud and Harris \cite[Section 3.2.5, p. 117]{Eisenbud_Harris_geometry_schemes}). 
\end{rmk}

We now turn to the definition of Cartier divisors in the category of schemes.

\begin{defn}[Cartier subscheme]
\label{defn:Eisenbud_Harris_4-15}
(See Eisenbud and Harris \cite[Definition IV.15, p. 165]{Eisenbud_Harris_geometry_schemes}; see also G\"ortz and Wedhorn \cite[Definition 11.26, p. 304]{Gortz_Wedhorn_algebraic_geometry_v1}, Hartshorne \cite[Chapter II, Section 6, Definition, p. 141]{Hartshorne_algebraic_geometry}, the Stacks Project \cite[\href{https://stacks.math.columbia.edu/tag/01WR}{Definition 01WR}]{stacks-project}, and Vakil \cite[Section 8.4.1]{Vakil_foundations_algebraic_geometry}.)  
Let $X$ be a scheme and $Y \subset X$ a subscheme. Then $Y$ is a \emph{Cartier subscheme} (or \emph{divisor}) in $X$ if it is locally the zero locus of a single nonzero divisor; that is, if for all $p \in X$ there is an affine neighborhood $U = \Spec A$ of $p$ in $X$ such that $Y \cap U = \VV(f) \subset U$ for some nonzero divisor $f \in A$. 
\end{defn}

\begin{rmk}[On Definition \ref{defn:Eisenbud_Harris_4-15} of a Cartier subscheme]
\label{rmk:On_definition_Cartier_subscheme}
G\"ortz and Wedhorn define both general and \emph{effective} Cartier divisors \cite[Definition 11.26, (1) and (4)]{Gortz_Wedhorn_algebraic_geometry_v1}. Hartshorne includes a definition of effective Cartier divisor in \cite[Section 2.6, Definition, p. 145]{Hartshorne_algebraic_geometry}. Vakil does not discuss general Cartier divisors, only effective Cartier divisors \cite[Section 0.2, p. 17, and Section 8.4]{Vakil_foundations_algebraic_geometry}. Similarly, the Stacks Project focuses on effective Cartier divisors \cite[\href{https://stacks.math.columbia.edu/tag/01WQ}{Section 01WQ}]{stacks-project}.
\end{rmk}

There are two approaches to defining the blowup of a scheme along a subscheme:
\begin{inparaenum}
\item \emph{non-constructive}, in which the blowup is (uniquely) characterized by a universal mapping property and one subsequently establishes its existence, and
\item \emph{constructive}, in which the blowup is defined explicitly and one subsequently establishes the universal mapping property and its uniqueness. 
\end{inparaenum}
We begin with the first definition of the blowup of a scheme.

\begin{defn}[Blowup of a scheme along a subscheme via characterization by a universal property]
\label{defn:Eisenbud_Harris_4-16}
(See Eisenbud and Harris \cite[Definition 4.16, p. 165]{Eisenbud_Harris_geometry_schemes}; see also G\"ortz and Wedhorn \cite[Definition 13.90, p. 413]{Gortz_Wedhorn_algebraic_geometry_v1}, Hironaka \cite[Chapter 0, Section 2, p. 129]{Hironaka_1964-I-II}, and Vakil \cite[Section 22.2]{Vakil_foundations_algebraic_geometry}. See Hauser \cite[Definition 4.1, p. 17]{Hauser_2014} for the definition of the blowup of an algebraic variety.)  
Let $X$ be a scheme and $Z \subset X$ be a closed subscheme as in Definition \ref{defn:Closed_subscheme_and_closed_immersion_or_embedding}. The \emph{blowup of $X$ along $Z$}, denoted $\pi : \Bl_Z(X) \to X$, is the morphism to $X$ characterized by the following properties:
\begin{enumerate}
\item\label{item:Cartier_subscheme} The inverse image\footnote{See Remark \ref{rmk:Inverse_images_divisors} for a reference to a definition.} $\pi^{-1}(Z)$ of $Z$ is a Cartier subscheme in $\Bl_Z(X)$.
\item\label{item:Universal_mapping_property_blowup_scheme} $\pi : \Bl_Z(X) \to X$ is universal with respect to this property; that is, if $\tau : W \to X$ is any morphism such that $\tau^{-1}(Z)$ is a Cartier subscheme in $Z$, then there is a unique morphism $\varphi: W \to \Bl_Z(X)$ such that $\tau = \pi\circ \varphi$:
\[
  \begin{tikzcd}
    W \arrow[r, "\varphi", dashed] \arrow[rd, "\tau"'] &\Bl_Z(X) \arrow[d, "\pi"]
    \\
    &X
  \end{tikzcd}
\]  
\end{enumerate}
The morphism $\pi$ is the \emph{blowup map}, the inverse image $E_ZX := \pi^{-1}(Z)$ of $Z$ in $\Bl_Z(X)$ is the \emph{exceptional divisor} of the blowup, and $Z$ is the \emph{center} of the blowup. We write $\widetilde X$ (or $X'$) and $E$ rather than $\Bl_Z(X)$ and $E_ZX$ when their meanings are clear from the context.
\end{defn}

\begin{rmk}[On Definition \ref{defn:Eisenbud_Harris_4-16} of the blowup]
\label{rmk:On_definition_blowup}
G\"ortz and Wedhorn and also Vakil require that $E$ be an effective Cartier divisor in their version of Definition \ref{defn:Eisenbud_Harris_4-16}. 
\end{rmk}

\begin{rmk}[Uniqueness and existence of the blowup in Definition \ref{defn:Eisenbud_Harris_4-16}]
\label{rmk:Existence_blowup}
The universal properties in Definition \ref{defn:Eisenbud_Harris_4-16} \emph{uniquely} characterize the blowup $\pi : \Bl_Z(X) \to X$ of a scheme along a subscheme as noted by Eisenbud and Harris \cite[Chapter IV, p. 165]{Eisenbud_Harris_geometry_schemes}, G\"ortz and Wedhorn \cite[p. 413]{Gortz_Wedhorn_algebraic_geometry_v1}, and Vakil \cite[p. 599]{Vakil_foundations_algebraic_geometry}.

\emph{Existence} of the blowup in Definition \ref{defn:Eisenbud_Harris_4-16} is provided by Eisenbud and Harris \cite[Theorem IV.23, p. 170]{Eisenbud_Harris_geometry_schemes} or, alternatively, using existence of the blowup of an affine scheme along a closed subscheme provided by Eisenbud and Harris \cite[Proposition IV.22, p. 169]{Eisenbud_Harris_geometry_schemes} and deducing the existence of blowups in general by gluing \cite[Chapter IV, p. 170]{Eisenbud_Harris_geometry_schemes}. See also G\"ortz and Wedhorn \cite[Proposition 13.92, p. 415]{Gortz_Wedhorn_algebraic_geometry_v1}, Hauser \cite[Definition 4.7, p. 18, and Theorems 4.18 and 4.19, pp. 20--21]{Hauser_2014}, and Vakil \cite[Theorem 22.3.2]{Vakil_foundations_algebraic_geometry} for existence. 
\end{rmk}

We now turn to the second definition of the blowup of a scheme, by explicit construction.

\begin{defn}[Blowup of a scheme along a subscheme via explicit construction]
\label{defn:Hartshorne_page_163}
(See Fulton \cite[Section B.6, p. 435]{Fulton2}, Grothendieck \cite[Section 8, p. 152]{Grothendieck_EGAII}, Hartshorne \cite[Chapter II, Section 7, Definition, p. 163, and Proposition 2.7.14]{Hartshorne_algebraic_geometry}, and the Stacks Project \cite[\href{https://stacks.math.columbia.edu/tag/01OG}{Definition 01OG}]{stacks-project}.)
Let $X$ be a Noetherian scheme and $\sI$ be a coherent sheaf of ideals on $X$. Consider the sheaf of graded algebras $\sS := \oplus_{d=0}\sI^d$, where $\sI^d$ is the $d$-th power of $\sI$ and $\sI^0 := \sO_X$. Define the \emph{blowup of $X$ with respect to the coherent sheaf of ideals $\sI$} to be
\[
  \Bl_Z(X) \text{ or } \widetilde X := \bProj\sS,
\]
where the scheme $\bProj\sS$ and morphism $\pi:\bProj\sS \to X$ are defined by Hartshorne in \cite[Section 2.7, Construction, p. 160]{Hartshorne_algebraic_geometry}, and $E := \pi^{-1}(Z)$ is the \emph{exceptional divisor}. If $Z$ is the closed subscheme of $X$ corresponding to $\sI$, then $\widetilde X$ is the \emph{blowup of $X$ along $Z$} or \emph{with center} $Z$.
\end{defn}

\begin{lem}[Properties and uniqueness of the blowup in Definition \ref{defn:Hartshorne_page_163}]
\label{lem:Hartshorne_propositions_2-7-13_and_2-7-14}  
(See Hartshorne \cite[Chapter II, Section 7, Propositions 7.13 and 7.14, p. 164]{Hartshorne_algebraic_geometry} and the Stacks Project \cite[\href{https://stacks.math.columbia.edu/tag/02OS}{Lemma 02OS}]{stacks-project}; see also G\"ortz and Wedhorn \cite[Proposition 13.91 (3), p. 414]{Gortz_Wedhorn_algebraic_geometry_v1} for Item \eqref{item:pi_isomorphism_over_X_less_Z}.)  
If $X$ is a scheme and $Z\subset X$ is a closed subscheme, then the blowup $\pi:\widetilde X \to X$ along $Z$ in Definition \ref{defn:Hartshorne_page_163} has the following properties:
\begin{enumerate}
\item\label{item:pi_isomorphism_over_X_less_Z} $\pi:\pi^{-1}(X\less Z) \to X\less Z$ is an isomorphism.
\item The exceptional divisor $E$ is an effective Cartier divisor on $\widetilde X$.
\item The blowup $\pi:\widetilde X \to X$ obeys the universal mapping property in Item \eqref{item:Universal_mapping_property_blowup_scheme} in Definition \ref{defn:Eisenbud_Harris_4-16}.  
\end{enumerate}
\end{lem}

\begin{rmk}[Explicit constructions of blowups for algebraic varieties]
\label{rmk:Definition_existence_uniqueness_blowup_algebraic_varieties}
Explicit constructions of the blowup $\widetilde X$ are provided by Hartshorne \cite[Chapter I, Section 4, pp. 28--30]{Hartshorne_algebraic_geometry}, Hauser \cite[Section 4]{Hauser_2014}, and Shafarevich \cite[Sections 2.4.1--2.4.4]{Shafarevich_v1}.
\end{rmk}

\begin{prop}[Properties of the blowup of an algebraic variety with respect to an ideal]
\label{prop:Hartshorne_2-7-16}
(See Hartshorne \cite[Chapter II, Section 7, Proposition 7.16, p. 166]{Hartshorne_algebraic_geometry}.)
Let $X$ be an algebraic variety over an algebraically closed field $\KK$, let $\sI \subset \sO_X$ be a nonzero coherent sheaf of ideals on $X$, and $\pi:\widetilde X\to X$ be the blowup with respect to $\sI$. Then the following hold:
\begin{enumerate}
\item $\widetilde X$ is an algebraic variety.
\item $\pi$ is a birational, proper, surjective morphism.
\item If $X$ is quasi-projective (respectively, projective) over $\KK$, then the same is true for $\widetilde X$, and $\pi$ is a projective morphism.  
\end{enumerate}
\end{prop}  

\begin{rmk}[Blowups in resolutions of singularities]
\label{rmk:Properties_blowups}  
Koll\'ar \cite[Section 3.3, p. 136]{Kollar_lectures_resolution_singularities} provides a summary of certain properties of blowups that are relevant to monomialization and resolution of singularities. One says that a blowup is \emph{trivial} if $Z$ is a Cartier divisor in $X$ and in this case, $\pi : \widetilde X \to X$ is an isomorphism by the universal property of blowups (see also  Hauser \cite[Remark 4.5, p. 18]{Hauser_2014}, the Stacks Project \cite[\href{https://stacks.math.columbia.edu/tag/0807}{Lemma 0807}]{stacks-project}, and Vakil \cite[Observation 22.2.1]{Vakil_foundations_algebraic_geometry}). In resolution of singularities, one includes the possibility that $Z = \varnothing$, which is known as the \emph{empty} blowup. Empty blowups arise when one restricts a blowup sequence to an open subset $U \subset X$ and the center of the blowup is disjoint from $U$. However, in statements of theorems on resolution of singularities, one does not include empty blowups in the final blowup sequences.
\end{rmk}

\subsection{Strict transforms of subschemes in blowups}
\label{subsec:Strict_transforms_subschemes_blowups}
In this subsection, we review definitions of the strict transform of a subscheme in the blowup of a scheme.

\begin{defn}[Total transform in the blowup of a scheme along a closed subscheme]
\label{defn:Total_transform}  
(See Hauser \cite[Definition 6.1, p. 30]{Hauser_2014} for algebraic varieties and Eisenbud and Harris \cite[Section 4.2, p. 168]{Eisenbud_Harris_geometry_schemes}, Hauser \cite[p. 24]{Hauser_2000_survey}, Hironaka \cite[Chapter 0, Section 5, p. 142]{Hironaka_1964-I-II}, or Vakil \cite[Section 22.2]{Vakil_foundations_algebraic_geometry} for schemes.)
Let $\pi : \widetilde X \to X$ denote the blowup of an a scheme $X$ along a closed subscheme $Z \subset X$, with exceptional divisor $E := \pi^{-1}(Z)$ defined by the principal ideal $\sI_E$ of $\sO_{\widetilde X}$. Let $\pi^* : \sO_X \to \sO_{\widetilde X}$ be the dual homomorphism of $\pi$. Let $Y \subset X$ be a closed subscheme and $\sI \subset \sO_X$ be an ideal. The inverse image\footnote{See Remark \ref{rmk:Inverse_images_divisors} for a reference to a definition.} $Y^* := \pi^{-1}(Y)$ of $Y$ and the extension $\sI^* = \pi^*\sI = \sI\cdot\sO_{\widetilde X}$ of $\sI$ are called the \emph{total transform} of $Y$ and $\sI$ under $\pi$.
\end{defn}

\begin{defn}[Strict transform in the blowup of a scheme along a closed subscheme]
\label{defn:Strict_transform}
(See Hauser \cite[Definition 6.2, p. 30]{Hauser_2014} or Koll\'ar \cite[Definition 2.1, p. 68]{Kollar_lectures_resolution_singularities} for algebraic varieties and see Eisenbud and Harris \cite[Section 4.2, p. 168]{Eisenbud_Harris_geometry_schemes}, G\"ortz and Wedhorn \cite[Chapter 13, p. 416]{Gortz_Wedhorn_algebraic_geometry_v1}, Hartshorne \cite[Chapter II, Section 7, Definition, p. 165]{Hartshorne_algebraic_geometry}, Hauser \cite[p. 24]{Hauser_2000_survey}, or Vakil \cite[Section 22.2, p. 600]{Vakil_foundations_algebraic_geometry} for schemes; see also Hironaka \cite[Chapter 0, Section 2, p. 130]{Hironaka_1964-I-II}.)
Continue the notation of Definition \ref{defn:Total_transform}. The scheme theoretic closure\footnote{Zariski closure in the case of algebraic varieties.} $\overline{\pi^{-1}(Y \less Z)}$ of $\pi^{-1}(Y \less Z)$ in $\widetilde X$ is the \emph{strict} (or \emph{proper} or \emph{birational}) \emph{transform} of $Y$ under $\pi$ and denoted by\footnote{We generally write $Y'$ in the context of resolution of singularities and write $\widetilde Y$ in the context of a single blowup, as opposed to a composition of blowups. Hauser occasionally writes $Y^s$, but this notation appears uncommon.} $\widetilde Y$ or $Y'$, depending on the context.
\end{defn}

In the forthcoming Definition \ref{defn:Strict_transform_and_total_transform_algebraic_varieties}, we define the \emph{total transform} and \emph{strict transform} in the context of resolution of singularities in the category of algebraic varieties, where the resolution morphism is a composition of blowup morphisms. Given the importance of the concept in Definition \ref{defn:Strict_transform}, it is worth recalling the

\begin{defn}[Scheme theoretic closure and density]
\label{defn:Scheme_theoretic_closure_density}  
(See Grothendieck \cite[Definition 11.10.2, p. 171]{Grothendieck_EGAIV-3}, Hartshorne \cite[Chapter II, Section 3, Exercise 3.11(d), p. 92]{Hartshorne_algebraic_geometry}, and the Stacks Project \cite[\href{https://stacks.math.columbia.edu/tag/01RB}{Definition 01RB}]{stacks-project}; see also G\"ortz and Wedhorn \cite[Remark 10.31, p. 255]{Gortz_Wedhorn_algebraic_geometry_v1} and Vakil \cite[Section 8.3.8, p. 239]{Vakil_foundations_algebraic_geometry}.)  
Let $X$ be a scheme and $U\subset X$ be an open subscheme. The scheme theoretic image of the morphism $U \to X$ is the \emph{scheme theoretic closure} of $U$ in $X$. One says that $U$ is \emph{scheme theoretically dense} in $X$ if for every open subset $V \subset X$, the scheme theoretic closure of $U\cap V$ in $V$ is equal to $V$.
\end{defn}

It is useful to note the following relation between topological and scheme theoretic density.

\begin{lem}[Scheme theoretic and topological density]
\label{lem:Scheme_theoretic_toplogical_density}  
(See the Stacks Project \cite[\href{https://stacks.math.columbia.edu/tag/056D}{Lemma 056D}]{stacks-project}; compare G\"ortz and Wedhorn \cite[Remark 10.32, p. 255]{Gortz_Wedhorn_algebraic_geometry_v1} and Vakil \cite[Section 8.3.9, p. 239]{Vakil_foundations_algebraic_geometry}.)
Let $X$ be a reduced scheme and $U \subset X$ be an open subscheme. Then the following are equivalent:
\begin{enumerate}
\item $U$ is topologically dense in $X$.
\item The scheme theoretic closure of $U$ in $X$ is $X$.
\item $U$ is scheme theoretically dense in $X$.
\end{enumerate}
\end{lem}

Vakil provides the following intuition underlying Definition \ref{defn:Scheme_theoretic_closure_density} and Lemma \ref{lem:Scheme_theoretic_toplogical_density}.

\begin{rmk}[Induced reduced subscheme structure on a closed subset of a scheme]
\label{rmk:Vakil_8-3-9}
(See Vakil \cite[Section 8.3.9, p. 239]{Vakil_foundations_algebraic_geometry}; compare Eisenbud and Harris \cite[Section V.1.1, pp. 209--213]{Eisenbud_Harris_geometry_schemes}.)  
Suppose $X^\set$ is a topologically closed subset of a scheme $Y$. Then one can define a canonical
scheme structure $X$ on $X^\set$ that is reduced. Informally, $X$ is cut out
by those functions whose values are zero at all the points of $X^\set$. On the affine
open set $\Spec A$ of $Y$, if the set $X^\set$ corresponds to the radical ideal $\sI = \sI(X^\set)$ (see Vakil \cite[Section 3.7 p. 127]{Vakil_foundations_algebraic_geometry}, the scheme $X$ corresponds to $\Spec(A/I)$. Alternatively, one could define $X$ to be the smallest closed subscheme whose underlying set contains $X^\set$. According to Vakil \cite[Section 3.7 p. 127]{Vakil_foundations_algebraic_geometry}, these definitions yield the same scheme $X$. The construction is called the \emph{induced reduced subscheme structure} on the closed subset $X^\set$.  
\end{rmk}

In \cite[\href{https://stacks.math.columbia.edu/tag/080C}{Section 080C}]{stacks-project}, the Stacks Project authors provide a different approach to the definition of strict transform than that of Definition \ref{defn:Strict_transform} and we summarize their method now (in our notation). Let $X$ be a scheme, $Z\subset X$ be a closed subscheme, and $\pi:\widetilde X \to X$ be blowup of $X$ along $Z$, and $E = \pi^{-1}(Z)\subset \widetilde X$ be the exceptional divisor. We consider a scheme $Y$ over $X$, with morphism $f: Y \to X$, and form the Cartesian diagram,
\[
   \begin{tikzcd}
    \pr_{\widetilde X}^{-1}(E) \arrow[d] \arrow[r] &Y\times_X\widetilde X \arrow[d, "\pr_{\widetilde X}"]   \arrow[r, "\pr_Y"] &Y  \arrow[d, "f"]
    \\
    E \arrow[r] &\widetilde X \arrow[r, "\pi"] & X
   \end{tikzcd}
\]
where $\pr_Y$ and $\pr_{\widetilde X}$ are the projections onto the first and second factors, respectively, in the fiber product, $Y\times_X\widetilde X$. Since $E$ is an effective Cartier divisor (see the Stacks Project \cite[\href{https://stacks.math.columbia.edu/tag/02OS}{Lemma 02OS}]{stacks-project}), one sees that $\pr_{\widetilde X}^{-1}(E)$ is locally principal (see the Stacks Project \cite[\href{https://stacks.math.columbia.edu/tag/053P}{Lemma 053P}]{stacks-project}).

\begin{defn}[Strict transform in the blowup of a scheme along a closed subscheme via fiber products]
\label{defn:Strict_transform_stacks_project}
(See the Stacks Project \cite[\href{https://stacks.math.columbia.edu/tag/080D}{Definition 080D}]{stacks-project}.)  
With $Z\subset X$ and $f:Y\to X$ as above, then the \emph{strict transform} of $Y$ is the closed subscheme $\widetilde Y \subset Y\times_X\widetilde X$ cut out by the quasi-coherent ideal of sections of $\sO_{Y\times_X\widetilde X}$ supported on $\pr_{\widetilde X}^{-1}(E)$.
\end{defn}

The Stacks Project authors identify the strict transform more explicitly in the

\begin{lem}[Strict transforms as blowups or base changes]
\label{lem:Strict_transforms_blowups_or_base_changes}
(See the Stacks Project \cite[\href{https://stacks.math.columbia.edu/tag/080E}{Lemma 080E}]{stacks-project} and \cite[\href{https://stacks.math.columbia.edu/tag/080F}{Lemma 080F}]{stacks-project}.)
Continue the notation of Definition \ref{defn:Strict_transform_stacks_project}.
\begin{enumerate}
\item\label{item:Strict_transform_as_blowup} The strict transform $\widetilde Y$ of $Y$ is the blowup of $Y$ along the closed subscheme $f^{-1}(Z)$ of $Y$.
\item\label{item:Strict_transform_as_base_change} If $Y$ is flat over $X$ at all points lying over $Z$ (see the forthcoming Definition \ref{defn:Flat_modules_and_flat_morphisms_ringed_spaces}, then the strict transform of $Y$ is equal to the \emph{base change} $Y\times_X\widetilde X$.
\end{enumerate}  
\end{lem}

Hauser proves an analogue of Item \eqref{item:Strict_transform_as_blowup} in Lemma \ref{lem:Strict_transforms_blowups_or_base_changes} in the category of algebraic varieties as \cite[Proposition 5.1, p. 27]{Hauser_2014} and, in particular, establishes that the blowup $\pi_Y:\widetilde Y \to Y$ is equal to the restriction of the projection $\pr_Y:Y\times_X\widetilde X$ to the closed subvariety $\widetilde Y \subset Y\times_X\widetilde X$. Hauser also proves a special case of Item \eqref{item:Strict_transform_as_blowup} in \cite[Corollary 5.2 (a), p. 28]{Hauser_2014} when $Y \subset X$ is a closed subvariety, so $f:Y\to X$ is inclusion and $f^{-1}(Z) = Z\cap Y$.

For the definition of inverse image in the category of complex analytic spaces, see Fischer \cite[Section 0.27]{Fischer_complex_analytic_geometry} and in the category of schemes, see Remark \ref{rmk:Inverse_images_divisors} for references to definitions.

\begin{rmk}[Strict transform in the blowup of a scheme along a closed subscheme via colon ideals]
\label{rmk:Strict_transform_colon_ideals}
In our applications, it will be useful to have a more explicit construction of the strict transform in the blowup of a scheme along a closed subscheme than that provided by Definitions \ref{defn:Strict_transform} or \ref{defn:Strict_transform_stacks_project}. For this purpose, we recall the a construction provided by Hauser \cite[p. 335]{Hauser_2003}, \cite[Section 6, pp. 30--31]{Hauser_2014} and for further discussion, we refer the reader to Appendix \ref{subsec:Sheaves_ideals_strict_transforms_resolution_blowups}. Continue the notation of Definition \ref{defn:Strict_transform} and let $\sI_U$ denote the restriction of an ideal $\sI \subset \sO_X$ to the open set $U := X \less Z$. Define $\widetilde U := \pi^{-1}(U) \subset \widetilde X$ and let $\tau: \widetilde U \to U$ be the restriction of $\pi:\widetilde X \to X$ to $\widetilde U$. One defines the strict transform of the ideal $\sI$ by
\[
  \tilde\sI := \tau^*(\sI_U)\cap\sO_{\widetilde X}.
\]
If the ideal $\sI$ defines the closed subscheme $Y \subset X$, then $\tilde\sI$ defines the closed subscheme $\widetilde Y \subset \widetilde X$. One can show that $\tilde\sI$ is given by a union of \emph{colon ideals},
\begin{equation}
  \label{eq:Strict_transform_ideal_as_union_colon_ideals}
  \tilde\sI = \bigcup_{i\geq 0}(\sI^*:\sI_E^i),
\end{equation}
where $\sI^* \subset \sO_{\widetilde X}$ is the total transform of $\sI$ (see Definition \ref{defn:Total_transform}) and $\sI_E \subset \sO_{\widetilde X}$ is the ideal defining the exceptional divisor $E = \pi^{-1}(Z) \subset \widetilde X$. The expression for $\tilde\sI$ in terms of colon ideals can be understood with the aid of Hauser \cite[Lemma 6.3, p. 31]{Hauser_2014}.

Let $\tilde p \in \widetilde X$ be a point and $h \in \sO_{\widetilde X,\tilde p}$ be a germ that locally defines $E \subset \widetilde X$. Near $\tilde p$, one has
\[
  \tilde\sI = (\tilde f: f \in \sI),
\]
where the strict transform $\tilde f$ of $f$ is defined at $\tilde p$, up to multiplication by invertible elements in $\sO_{\widetilde X}$, by
\[
  f^* = h^k\cdot\tilde f,
\]
with maximal exponent $k := \ord_Ef^*$ given by the \emph{order} of $f^*$ along $E$ (see the forthcoming Definition \ref{defn:Order}), and $f^* := \pi^*f \in \sO_{\widetilde X}$. By abuse of notation, the preceding identity can be written as
\[
  \tilde f = h^{-k}\cdot f^*.
\]
This expression tells us that the strict transform $\widetilde Y \subset \widetilde X$ of $Y$ is cut out by an ideal $\tilde \sI \subset \sO_{\widetilde X}$ of sections that are supported on $E \subset \widetilde X$. Thus, $h$ is not a zero divisor in $\sO_{\widetilde Y,\tilde p}$, where $\sO_{\widetilde Y} = (\sO_{\widetilde X}/\tilde\sI)\restriction \widetilde Y$.
\end{rmk}

Suppose that $\sI$ in Remark \ref{rmk:Strict_transform_colon_ideals} is generated locally by elements $f_1, \ldots, f_k$ of $\sO_X$. Hauser notes that $\widetilde \sI$ \emph{contains} contains the ideal $(\tilde f_1,\ldots,\tilde f_k)$ generated by the strict transforms $\tilde f_1,\ldots,\tilde f_k$ of $f_1, \ldots, f_k$ but the inclusion can be strict (see \cite[Remark 6.4, p. 31 and Example 6.22, p. 32]{Hauser_2014}).

\begin{defn}[Order]
\label{defn:Order}  
(See Hauser \cite[Definition 8.4, p. 36]{Hauser_2014}.)
Let $Y$ be a subvariety of a not necessarily regular ambient variety $X$ defined by an ideal $\sI$ and let $Z$ be an irreducible subvariety of $X$ defined by a prime ideal $\sJ$. The \emph{order} of $Y$ or $\sI$ in $X$ along $Z$ or with respect to $\sJ$ is the maximal integer $k = \ord_Z(X) = \ord_Z(\sI)$ such that $\sI_Z \subset \sJ_Z^k$, where $\sI_Z = \sI ·\cdot \sO_{X,Z}$ and $\sJ_Z = \sJ\cdot\sO_{X,Z}$ denote the ideals generated by $\sI$ and $\sJ$ in the localization $\sO_{X,Z}$ of $X$ along $Z$. If $Z = \{p\}$ is a point in $X$, then the order of $Y$ and $\sI$ at $p$ is denoted by $\ord_p(Y) = \ord_p(\sI)$ or $\ord_{\fm_p}(X) = \ord_{\fm_p}(\sI)$.
\end{defn}

\subsection{Smoothness of blowups and strict transforms}
\label{subsec:Smoothness_blowups_strict_transforms}
In this subsection, we consider the question of smoothness of blowups and strict transforms.

\begin{thm}[Smoothness of the blowup of a smooth variety along a closed, smooth subvariety]
\label{thm:Hartshorne_2-8-24a}
(See Hartshorne \cite[Chapter II, Section 8, Theorem 8.24 (a), p. 186]{Hartshorne_algebraic_geometry}, Hauser \cite[Proposition 5.3, p. 28]{Hauser_2014} by taking $Y=X$, or Vakil \cite[Theorem 22.3.10, p. 607]{Vakil_foundations_algebraic_geometry}.)
Let $X$ be a smooth variety over a field $\KK$. If $Z \subset X$ is a closed, smooth subvariety, then the blowup $\widetilde X$ of $X$ along $Z$ is regular.
\end{thm}  

\begin{cor}[Smoothness of the strict transform of a closed, smooth subvariety]
\label{cor:Smoothness_strict_transform_subvariety}
(See also Hauser \cite[Proposition 5.3, p. 28]{Hauser_2014}.)
Continue the notation of Theorem \ref{thm:Hartshorne_2-8-24a} and assume further that $Y \subset X$ is a closed, smooth subvariety. If $Z\cap Y$  is a closed, smooth subvariety of $Y$, then the strict transform $\widetilde Y$ as in Definitions \ref{defn:Strict_transform} or \ref{defn:Strict_transform_stacks_project} is a smooth variety.
\end{cor}

\begin{proof}
According to Item \eqref{item:Strict_transform_as_blowup} in Lemma \ref{lem:Strict_transforms_blowups_or_base_changes} (or Hauser \cite[Corollary 5.2 (a), p. 28]{Hauser_2014}), the strict transform $\widetilde Y$ is equal to the blowup of $Y$ along $Z\cap Y$. Since $Y$ is regular and the blowup center $Z\cap Y$ is regular, then $\widetilde Y$ is regular by Theorem \ref{thm:Hartshorne_2-8-24a}.
\end{proof}

Analogues of Theorem \ref{thm:Hartshorne_2-8-24a} and Corollary \ref{cor:Smoothness_strict_transform_subvariety} also hold in the category of schemes.

\begin{lem}[Smoothness of the blowup of a smooth scheme along a closed, smooth subscheme]
\label{lem:Hartshorne_2-8-24a_schemes}  
(See the Stacks Project \cite[\href{https://stacks.math.columbia.edu/tag/0FUT}{Lemma 0FUT}]{stacks-project}.)  
Let $S$ be a scheme and $Z \to X$ be a closed immersion of schemes smooth over $S$. If $\pi:\widetilde X \to X$ is the blowup along $Z$ with exceptional divisor $E \subset \widetilde X$, then $\widetilde X$ and $E$ are smooth over $S$. The morphism $p:E\to Z$ is canonically isomorphic to the projective space bundle
\[
  \PP(\sI/\sI^2) \to Z,
\]  
where $\sI \subset X$ is the ideal sheaf of $Z$. 
\end{lem}

According to the Stacks Project \cite[\href{https://stacks.math.columbia.edu/tag/01IO}{Definition 01IO}]{stacks-project}, a morphism of schemes is a \emph{closed immersion (of schemes)} if it is a closed immersion of locally ringed spaces. Moreover, according to Stacks Project \cite[\href{https://stacks.math.columbia.edu/tag/01HK}{Definition 01HK}]{stacks-project}, a morphism $\iota:Z\to X$ of locally ringed spaces is a \emph{closed immersion of locally ringed spaces} if
\begin{inparaenum}
\item the map $\iota$ is a homeomorphism of $Z$ onto a closed subset of $X$,
\item the map $\sO_X\to\iota_*\sO_Z$ is surjective and, if $\sI$ denotes the kernel,
\item the $\sO_X$-module $\sI$ is locally generated by sections.
\end{inparaenum}

\begin{cor}[Smoothness of the strict transform of a smooth, closed subscheme]
\label{cor:Smoothness_strict_transform_subscheme}
Let $X$ be a smooth scheme and $Z \subset X$ be a closed, smooth subscheme. If $Y \subset X$ is a smooth, closed subscheme and $Z\cap Y$ is a closed, smooth subscheme of $Y$, then the strict transform $\widetilde Y$ as in Definitions \ref{defn:Strict_transform} or \ref{defn:Strict_transform_stacks_project} is a smooth scheme.
\end{cor}

\begin{proof}
According to Item \eqref{item:Strict_transform_as_blowup} in Lemma \ref{lem:Strict_transforms_blowups_or_base_changes}, the strict transform $\widetilde Y$ is equal to the blowup of $Y$ along $Z\cap Y$. Since $Y$ is smooth and the blowup center $Z\cap Y$ is smooth, then $\widetilde Y$ is smooth by Lemma \ref{lem:Hartshorne_2-8-24a_schemes}.
\end{proof}

\section{Flatness and functorial property of blowups of schemes}
\label{sec:Functorial_properties_blowups_schemes}
Before discussing the functorial property of blowups, we recall the following important definition. Let $f : X \to Y$ be a morphism of ringed spaces and let $\sF$ be an $\sO_X$-module. For each point $p \in X$, the $\sO_{X,p}$-module $\sF_p$ is endowed via the homomorphism $f_p^\sharp : \sO_{Y,f(p)} \to \sO_{X,p}$ with the structure of an $\sO_{Y,f(p)}$-module. For the notion of \emph{flat module over a ring}, see G\"ortz and Wedhorn \cite[Section (B.4), Definition/Proposition B.16, p. 559]{Gortz_Wedhorn_algebraic_geometry_v1}. For a more thorough discussion of flat morphisms, see G\"ortz and Wedhorn \cite[Chapter 14, p. 428]{Gortz_Wedhorn_algebraic_geometry_v1} and the Stacks Project \cite[\href{https://stacks.math.columbia.edu/tag/01U2}{Section 01U2}]{stacks-project}. G\"ortz and Wedhorn remind us that although it is not possible to fully describe the algebraic concept of flatness geometrically, the following intuition is useful: If $f : X \to S$ is a flat morphism of schemes, then the fibers $f^{-1}(s)$ form a continuously varying family of algebraic varieties, as $s$ varies in $S$.

\begin{defn}[Flat modules and flat morphisms]
\label{defn:Flat_modules_and_flat_morphisms_ringed_spaces}  
(See G\"ortz and Wedhorn \cite[Definition 7.38, p. 197 or Definition 14.1, p. 428]{Gortz_Wedhorn_algebraic_geometry_v1}
Let $f : X \to Y$ be a morphism of ringed spaces and let $\sF$ be an $\sO_X$-module.
\begin{enumerate}
\item\label{item:Flat_module_over_morphism_ringed_spaces} The $\sO_X$-module $\sF$ is called \emph{flat over $Y$ at $p$} if $\sF_p$ is a flat $\sO_{Y,f(p)}$-module; $\sF$ is called \emph{flat over $Y$} if $\sF$ is flat over $Y$ at all points $p \in X$.
\item\label{defn:Flat_morphisms_ringed_spaces} One says that $f$ is a \emph{flat morphism}, or that \emph{$X$ is flat over $Y$}, if $\sO_X$ is flat over $Y$.
\end{enumerate}
\end{defn}

In Definition \ref{defn:Flat_modules_and_flat_morphisms_ringed_spaces}, one considers $\sF_p$ to be $\sO_{Y,f(p)}$-module via the canonical homomorphism
\[
  f_p^\sharp:\sO_{Y,f(p)} \to \sO_{X,p}.
\]
For every ringed space $(X,\sO_X)$, locally free $\sO_X$-modules are flat (see G\"ortz and Wedhorn \cite[Section 7.18, p. 198]{Gortz_Wedhorn_algebraic_geometry_v1}). There is also the following converse.

\begin{prop}[Equivalent conditions for flat modules]
\label{prop:Equivalant_conditions_flat_modules}  
(See G\"ortz and Wedhorn \cite[Proposition 7.41, p. 198]{Gortz_Wedhorn_algebraic_geometry_v1}.)  
Let $(X,\sO_X)$ be a locally ringed space and let $\sF$ be an $\sO_X$-module. Then the following assertions are equivalent.
\begin{enumerate}
\item $\sF$ is locally free of finite type.
\item $\sF$ is of finite presentation and $\sF_p$ is a free $\sO_{X,p}$-module for all $p \in X$.
\item $\sF$ is flat and of finite presentation.
\end{enumerate}
\end{prop}

Flat morphisms between schemes have many desirable properties, of which the following are the most relevant to applications in our work.

\begin{prop}[Inverse image of an effective Cartier divisor under a flat morphism]
\label{prop:Inverse_image_effective_Cartier_divisor_flat_morphism}  
(See G\"ortz and Wedhorn \cite[Corollary 11.51, p. 316]{Gortz_Wedhorn_algebraic_geometry_v1}.) 
Let $f : X \to Y$ be a flat morphism of schemes and $Z \subset Y$ be an effective Cartier divisor as in Definition \ref{defn:Eisenbud_Harris_4-15} and Remark \ref{rmk:On_definition_Cartier_subscheme}. Then the inverse image $f^{-1}(Z)$ (as a subscheme) is an effective Cartier divisor in $X$ and this divisor is the inverse image $f^*(Z)$ (as a divisor).
\end{prop}

\begin{prop}[Properties of flat modules and morphisms of schemes]
\label{prop:Equivalant_conditions_flat_modules}  
(See G\"ortz and Wedhorn \cite[Proposition 14.3, p. 429]{Gortz_Wedhorn_algebraic_geometry_v1}.) 
Flatness is stable under base change (as in Section \ref{sec:Fiber_products_base_change}) and under composition of morphisms of schemes.
\end{prop}


We now recall the functorial property of blowups; we refer the reader to G\"ortz and Wedhorn \cite[Diagram (4.4.7), p. 101]{Gortz_Wedhorn_algebraic_geometry_v1} for the concept of a \emph{Cartesian diagram} appearing in Item \ref{item:Blowup_commutes_flat_base_change} of Proposition \ref{prop:Functorial_property_blowup_scheme} below.

\begin{prop}[Blowing up commutes with flat base change]
\label{prop:Functorial_property_blowup_scheme}
(See G\"ortz and Wedhorn \cite[Proposition 13.91, p. 414]{Gortz_Wedhorn_algebraic_geometry_v1}; see also Eisenbud and Harris \cite[Proposition IV.21, p. 167]{Eisenbud_Harris_geometry_schemes}, Hartshorne \cite[Chapter II, Section 7, Corollary 7.15, p. 165]{Hartshorne_algebraic_geometry}, and the Stacks Project \cite[\href{https://stacks.math.columbia.edu/tag/0805}{Lemma 0805}]{stacks-project}. See Hauser \cite[Proposition 5.1, p. 27]{Hauser_2014} for algebraic varieties.)  
Let $X$ be a scheme, let $Z$ be a closed subscheme of $X$, and let $\pi_X: \Bl_Z(X) \to X$ be the blowup of $X$ along $Z$. If $f : Y \to X$ is a morphism of schemes, then the following hold:
\begin{enumerate}
\item\label{item:Functorial_property_blowup_scheme} There exists a unique morphism\footnote{See Remark \ref{rmk:Inverse_images_divisors} for a reference to a definition of $f^{-1}(Z)$.} $\Bl_Z(f) : \Bl_{f^{-1}(Z)}(Y) \to \Bl_Z(X)$ such that the following diagram commutes:
    \begin{equation}
      \label{eq:Goertz_Wedhorn_13-19-1}
      \begin{tikzcd}
    \Bl_{f^{-1}(Z)}(Y) \arrow[r, "\Bl_Z(f)"] \arrow[d, "\pi_Y"] &\Bl_Z(X) \arrow[d, "\pi_X"]
    \\
    Y \arrow[r, "f"] &X
  \end{tikzcd}
\end{equation}
where $\pi_Y: \Bl_{f^{-1}(Z)}(Y) \to Y$ is the blowup of $Y$ along $f^{-1}(Z)$.
\item\label{item:Blowup_commutes_flat_base_change} If $f$ is flat in the sense of Definition \ref{defn:Flat_modules_and_flat_morphisms_ringed_spaces}, then the diagram \eqref{eq:Goertz_Wedhorn_13-19-1} is Cartesian, that is, the induced morphism to the fiber product (see Section \ref{sec:Fiber_products_base_change})
\[
  r:\Bl_{f^{-1}(Z)}(Y) \to \Bl_Z(X) \times_X Y
\]
is an isomorphism.
\end{enumerate}
\end{prop}

\begin{rmk}[Blowups do not commute with arbitrary base change]
Vakil notes that when the flatness condition is dropped, then blowing up need not commute with base change (see Vakil \cite[Exercise 24.2.P, p. 652]{Vakil_foundations_algebraic_geometry}). The often-cited phrase \emph{commutes with base change} may be understood as follows: If $f:Y \to X$ is a flat morphism of schemes and $Z \subset Y$ is a closed embedding of schemes, then there is a canonical isomorphism of schemes \cite[Exercise 24.2.P (a), p. 652]{Vakil_foundations_algebraic_geometry}:
\[
  \Bl_{Z\times_X Y}(Y) \cong \Bl_Z(X) \times_X Y,
\]
noting that $f^{-1}(Z) = Z\times_X Y$. However, the preceding relation does not hold for arbitrary morphisms $f:Y \to X$ of schemes \cite[Exercise 24.2.P (b), p. 652]{Vakil_foundations_algebraic_geometry}.
\end{rmk}

\chapter{Blowups and their properties for analytic spaces}
\label{chap:Blowups_properties_analytic_spaces}
In this chapter on blowups we provide the necessary background for our exposition of resolution of singularities for analytic spaces and for our development of results on the Birula--Bia{\l}ynicki decomposition for complex analytic spaces. In Section \ref{sec:Blowups_vector_spaces_along_linear_subspaces}, we review the model case of the blowup of a vector space along a linear subspace. Section \ref{sec:Blowups_analytic_manifolds_along_embedded_analytic_submanifolds} contains a discussion of blowups of analytic manifolds along embedded analytic submanifolds. In Section \ref{sec:Functorial_quivariance_properties_manifolds}, we describe functorial properties of blowups of manifolds. Section \ref{sec:Blowing_up_analytic_spaces} discusses blowups of analytic spaces. We conclude in Section \ref{sec:Exceptional_divisor_and_strict_transform_complex_analytic_model_space} with a descrioption of the exceptional divisor and strict transform for the blowup of an analytic model space along a linear subspace.

\section{Blowups of vector spaces along linear subspaces}
\label{sec:Blowups_vector_spaces_along_linear_subspaces}
We primarily follow the discussion due to Huybrechts \cite[Section 2.5, p. 98]{Huybrechts_2005}, except that we allow $\KK=\RR$ or $\CC$ and also include details from Griffiths and Harris \cite[Chapter 4, Section 6, pp. 602--605]{GriffithsHarris} and Voisin \cite[Section 3.3.3, p. 98]{Voisin_hodge_theory_complex_algebraic_geometry_I}. See also Hauser \cite[Appendix C, p. 394]{Hauser_2003}, \cite[Definition 4.12, p. 19]{Hauser_2014}, and \cite[Proposition 5.4, 28]{Hauser_2014} with proof in \cite{Hauser_2010}.

Adapting Huybrechts \cite[Section 2.5, Example 2.5.2, p. 99]{Huybrechts_2005}, we first recall the construction of the blowup $\Bl_{\KK^m}(\KK^n)$, where $0\leq m < n$, of $\KK^n$ along the coordinate subspace
\begin{equation}
  \label{eq:Coordinate_subspace_blowup_center}
  \KK^m := \left\{(z_1,\ldots,z_n) \in \KK^n:z_{m+1} = \cdots = z_n = 0\right\}
\end{equation}
defined by the standard coordinates $z = (z_1,\ldots,z_n)$ on $\KK^n$. We denote $\KK^0 = (0)$. We let $x = (x_{m+1}:\cdots:x_n)$ denote the homogeneous coordinates on the projective space $\PP(\KK^{n-m})$ and define
\begin{equation}
  \label{eq:Blowup_linear_subspace_equations}
  \Bl_{\KK^m}(\KK^n)
  := \left\{(x,z): z_ix_j = z_jx_i, \text{ for } i, j = m + 1,\ldots,n\right\}
    \subset
   \PP(\KK^{n-m})\times\KK^n.
\end{equation}    
Hence, the \emph{blowup of $\KK^n$ along the linear subspace $\KK^m\subset\KK^n$} is the incidence variety
\begin{equation}
  \label{eq:Blowup_linear_subspace_incidence_variety}
  \Bl_{\KK^m}(\KK^n) = \left\{(\ell,z) \in \PP(\KK^{n-m})\times\KK^n: z \in \langle\KK^m,\ell\rangle\right\},
\end{equation}
where $\ell \in \PP(\KK^{n-m})$ is represented by a line in the complementary linear subspace
\begin{equation}
  \label{eq:Coordinate_subspace_orthogonal_complement_blowup_center}
  \KK^{n-m} := \left\{(z_1,\ldots,z_n) \in \KK^n:z_1 = \cdots = z_m = 0\right\}
\end{equation}
and $\langle\KK^m,\ell\rangle$ is the linear span in $\KK^n$ of the subspace $\KK^m$ and line $\ell$.

The analytic projection map onto the first factor in \eqref{eq:Blowup_linear_subspace_equations} defined through the description \eqref{eq:Blowup_linear_subspace_incidence_variety} of $\Bl_{\KK^m}(\KK^n)$ as an incidence variety
\begin{equation}
  \label{eq:Blowup_linear_subspace_projection_KPn-m-1_incidence_variety}
  \pr_{\PP(\KK^{n-m})}:\Bl_{\KK^m}(\KK^n) \ni (\ell,z) \mapsto \ell \in \PP(\KK^{n-m})
\end{equation}
can also be constructed from the description \eqref{eq:Blowup_linear_subspace_equations} of $\Bl_{\KK^m}(\KK^n)$ in terms of its defining equations
\begin{equation}
  \label{eq:Blowup_linear_subspace_projection_KPn-m-1_equations}
  \pr_{\PP(\KK^{n-m})}:\Bl_{\KK^m}(\KK^n) \ni (x,z) \mapsto x \in \PP(\KK^{n-m}).
\end{equation}
The projection map \eqref{eq:Blowup_linear_subspace_projection_KPn-m-1_incidence_variety} (equivalently \eqref{eq:Blowup_linear_subspace_projection_KPn-m-1_equations}) exhibits $\Bl_{\KK^m}(\KK^n)$ as an analytic $\KK^{m+1}$-bundle over $\PP(\KK^{n-m})$, where the fiber over $\ell \in \PP(\KK^{n-m})$ is just
\[
  \pr_{\PP(\KK^{n-m})}^{-1}(\ell) = \langle\KK^m,\ell\rangle.
\]
Thus, $\Bl_{\KK^m}(\KK^n)$ is an analytic manifold. Moreover, the projection map onto the second factor in \eqref{eq:Blowup_linear_subspace_equations},
\begin{equation}
  \label{eq:Blowup_linear_subspace_projection_Kn}
  \pi = \pr_{\KK^n}:\Bl_{\KK^m}(\KK^n) \ni (x,z) \mapsto z \in \KK^n,
\end{equation}
is an isomorphism over $\KK^n\less\KK^m$ and $\pi^{-1}(\KK^m)$ is canonically isomorphic to $\PP(N_{\KK^n/\KK^m})$, where the normal bundle $N_{\KK^n/\KK^m}$ of the embedding $\KK^m\hookrightarrow\KK^n$ is canonically isomorphic to the product bundle $\KK^m\times\KK^{n-m}$ over $\KK^m$ and, in particular, is an embedded analytic submanifold of $\Bl_{\KK^m}(\KK^n)$ of dimension $n-1$. The linear subspace $\KK^m \subset \KK^n$ is the \emph{blowup center} and the hypersurface $\PP(N_{\KK^n/\KK^m})$ in $\Bl_{\KK^m}(\KK^n)$ is the \emph{exceptional divisor}.

The blowup $\Bl_{\KK^m}(\KK^n) \subset\PP(\KK^{n-m})\times\KK^n$ may be covered by $n-m$ coordinate domains
\begin{equation}
  \label{eq:Blowup_linear_subspace_coordinate_patch}
  U_j := \{((x_{m+1}:\cdots:x_n),z) \in \Bl_{\KK^m}(\KK^n): x_j \neq 0\}, \quad\text{for } j = m+1,\ldots,n,
\end{equation}
with $n$ analytic coordinate functions
\begin{equation}
  \label{eq:Blowup_linear_subspace_coordinates}
  w_k(j)
  :=
  \begin{cases}
    z_k, &\text{for } k = 1, \ldots,m,
    \\
    x_k/x_j = z_k/z_j, &\text{for } k = m+1, \ldots,\hat\jmath,\ldots,n,
    \\
    z_j, &\text{for } k = j,
  \end{cases}
\end{equation}
where $\hat\jmath$ indicates that the index $j$ is omitted from the index set $\{m+1,\ldots,n\}$,
and resulting analytic coordinate charts
\begin{equation}
  \label{eq:Blowup_linear_subspace_chart}
  \varphi_j:\Bl_{\KK^m}(\KK^n)\supset U_j \to \KK^n, \quad\text{for } j = m+1,\ldots,n.
\end{equation}
The coordinates $w_k(j)$ for $k \in \{m+1, \ldots,\hat\jmath,\ldots,n\}$ are Euclidean coordinates on each fiber $\pi^{-1}(z) \cong \PP(\KK^{n-m})$ of the exceptional divisor.

We may use the local coordinates \eqref{eq:Blowup_linear_subspace_coordinates} to describe the restriction $\pi:\Bl_{\KK^m}(\KK^n)|_{U_j} \to \KK^n$ of the blowup map $\pi$ in \eqref{eq:Blowup_linear_subspace_projection_Kn} to the domain $U_j \subset \Bl_{\KK^m}(\KK^n)$ of the coordinate chart:
\begin{multline}
  \label{eq:Blowup_linear_subspace_projection_Kn_local_chart}  
  \pi\circ\varphi_j^{-1}:\KK^n \supset \varphi_j(U_j) \ni (w_1,\ldots,w_n)
  \mapsto (z_1,\ldots,z_n) \in \KK^n,
  \\
  \text{where }
  z_k
  =
  \begin{cases}
    w_k, &\text{for } k = 1, \ldots,m,
    \\
    w_jw_k, &\text{for } k = m+1, \ldots,\hat\jmath,\ldots,n,
    \\
    w_j, &\text{for } k = j,
  \end{cases}
\end{multline}
for $j=m+1,\ldots,n$, where $\hat\jmath$ again indicates that the index $j$ is omitted from the index set $\{m+1,\ldots,n\}$. The map $\pi\circ\varphi_j^{-1}:\varphi_j(U_j)\less\{w_j=0\} \to \KK^n$ is an analytic embedding onto its image and $E_j := \varphi_j(U_j)\cap\{(w_1,\ldots,w_n)\in\KK^n:w_j=0\}$ is the exceptional divisor of the blowup map $\pi$ in \eqref{eq:Blowup_linear_subspace_projection_Kn} with respect to the local coordinate chart $\varphi_j$.

\section{Blowups of analytic manifolds along embedded analytic submanifolds}
\label{sec:Blowups_analytic_manifolds_along_embedded_analytic_submanifolds}
The construction of $\Bl_{\KK^m}(\KK^n)$ in Section \ref{sec:Blowups_vector_spaces_along_linear_subspaces} extends to define the blowup $\widetilde X = \Bl_Z(X)$ of an arbitrary analytic manifold $X$ of dimension $n$ along an arbitrary embedded analytic submanifold $Z \hookrightarrow X$ of dimension $m$. We choose an atlas $\{(V_\alpha,\varphi_\alpha)\}_{\alpha\in\sA}$ for $X$ with analytic coordinate charts
\[
  \varphi_\alpha: X \supset V_\alpha \to \varphi_\alpha(V_\alpha) \subset \KK^n
\]
that are compatible with the embedded submanifold $Z$ in the sense that
\[
  \varphi_\alpha(V_\alpha\cap Z) = \varphi_\alpha(V_\alpha) \cap\KK^m.
\]
We let $\pi : \Bl_{\KK^m}(\KK^n) \to \KK^n$ be the blowup of $\KK^n$ along $\KK^m$ as constructed above and denote by
\[
  \pi_\alpha:\widetilde X_\alpha \to \varphi_\alpha(V_\alpha)
\]
its restriction to the open subset $\varphi_\alpha(V_\alpha) \subset \KK^n$, so
\[
  \widetilde X_\alpha := \pi^{-1}(\varphi_\alpha(V_\alpha)), \quad\text{for all } \alpha\in\sA,
\]
and defined $\pi_\alpha := \pi\restriction \widetilde X_\alpha$, for all $\alpha\in\sA$. Huybrechts and Griffiths and Harris show that all the blowups over the open subsets $\varphi_\alpha(V_\alpha)$ naturally glue together to give the global blowup, $\pi:\Bl_Z(X)\to X$, thus proving the

\begin{prop}[Blowup of an analytic manifold along a submanifold]
\label{prop:Huybrechts_2-5-3}  
(See Huybrechts \cite[Proposition 2.5.3, p. 100]{Huybrechts_2005} or Voisin \cite[Lemma 3.33 and Definition 3.23, p. 80]{Voisin_hodge_theory_complex_algebraic_geometry_I} for $\KK=\CC$.)
Let $\KK=\RR$ or $\CC$ and $Z$ be an analytic embedded submanifold of an analytic manifold $X$. Then there exist an analytic manifold $\widetilde X = \Bl_Z(X)$, the \emph{blowup of $X$ along $Z$}, and an analytic map $\pi : \widetilde X \to X$ such that, if $E := \pi^{-1}(Z)$, then $\pi : \widetilde X \less E \to X\less Z$ is an analytic diffeomorphism and $\pi : E \to Z$ is isomorphic to the canonical projection $\PP(N_{Z/X}) \to Z$.
\end{prop}

\begin{defn}[Blowup center and exceptional divisor]
\label{defn:Huybrechts_2-5-4}    
(See Huybrechts \cite[Definition 2.5.4, p. 100]{Huybrechts_2005}.)
The hypersurface $E = \pi^{-1}(Z) = \PP(N_{Z/X}) \subset \Bl_Z(X)$ in Proposition \ref{prop:Huybrechts_2-5-3}  is the \emph{exceptional divisor} of the blowup morphism $\pi:\Bl_Z(X)\to X$, while $Z$ is its \emph{blowup center}.
\end{defn}

We shall need the

\begin{prop}[K\"ahler metric on the blowup of a compact complex K\"ahler manifold along a complex submanifold]
\label{prop:Voisin_3-24}  
(See Voisin \cite[Proposition 3.24, p. 80]{Voisin_hodge_theory_complex_algebraic_geometry_I}.)
Let $Z$ be a compact, complex embedded submanifold of a complex K\"ahler manifold $X$. Then the blowup $\Bl_Z(X)$ of $X$ along $Z$ is a K\"ahler manifold and is compact if $X$ is compact. 
\end{prop}

We also recall the useful

\begin{prop}[K\"ahler metric on complex submanifold of a complex K\"ahler manifold]
\label{prop:Hubrechts_3-1-10}  
(See Ballmann \cite[Example 4.10 (3), p. 43]{Ballmann_lectures_Kaehler_manifolds} or Huybrechts \cite[Proposition 3.1.10]{Huybrechts_2005}.)
Let $Z$ be a complex embedded submanifold of a complex manifold $X$. Then the restriction of a K\"ahler metric on $X$ to $Z$ yields a K\"ahler metric on $Z$.
\end{prop}

We may apply Proposition \ref{prop:Hubrechts_3-1-10} to give an elementary proof of existence of a K\"ahler metric that is quite different from that provided by Voisin for \cite[Proposition 3.24, p. 80]{Voisin_hodge_theory_complex_algebraic_geometry_I} in the simple case where $X=\CC^n$ and $Z=\CC^m$ is a coordinate subspace, so non-compact. Phong and Sturm \cite[Lemma 4, p. 427]{Phong_Sturm_2014} prove a version of Proposition \ref{prop:Voisin_3-24} that allows the embedded complex manifold $Z$ and complex K\"ahler manifold $X$ to be non-compact, but only asserting the existence of a K\"ahler metric on a compact subset of $\Bl_Z(X)$.

\begin{prop}[K\"ahler metric on the blowup of a complex vector space along a complex linear subspace]
\label{prop:Voisin_3-24_vector_space}  
Let $Z$ be a complex linear subspace of a complex vector space $X$ with a Hermitian inner product. Then the blowup $\Bl_Z(X)$ of $X$ along $Z$ is a K\"ahler manifold.
\end{prop}

\begin{proof}
We extend Huybrechts' exposition in \cite[Example 2.5.2, p. 99]{Huybrechts_2005} of the construction of $\Bl_{\CC^m}(\CC^n)$ to give a K\"ahler metric on $\Bl_Z(X)$ determined by the Hermitian inner product on $X$. Let $Z^\perp\subset X$ denote the orthogonal complement of $Z$ in $X$ with respect to the Hermitian inner product on $X$. The complex projective space $\PP(Z^\perp)$ has the Fubini--Study metric described, for example, by Huybrechts \cite[Example 3.1.9 (i), p. 117]{Huybrechts_2005} and, of course, the given Hermitian inner product is a K\"ahler metric on $X$. The product $\PP(Z^\perp)\times X$ is a product of two complex K\"ahler manifolds and thus itself a complex K\"ahler manifold. By its definition in \eqref{eq:Blowup_linear_subspace_incidence_variety} with $\KK=\CC$ and the development in \cite[Example 2.5.2, p. 99]{Huybrechts_2005}, the subset
\[
   \Bl_Z(X) \subset \PP(Z^\perp)\times X
\]
is an embedded complex submanifold and thus inherits a K\"ahler metric by Proposition \ref{prop:Hubrechts_3-1-10}.
\end{proof}  

\section{Functorial properties of blowups of manifolds}
\label{sec:Functorial_quivariance_properties_manifolds}
The references for blowups of complex analytic spaces in Section \ref{sec:Blowing_up_analytic_spaces} or blowups of complex manifolds along embedded complex submanifolds in Section \ref{sec:Blowups_analytic_manifolds_along_embedded_analytic_submanifolds} do not discuss analogues of the functorial properties of blowups of schemes along closed subschemes described in Section \ref{sec:Functorial_properties_blowups_schemes}. However, we have the following partial analogue. Let $M$ be a smooth manifold and $A$ a closed submanifold of $M$. If $x \in A$ is a point and the tangent spaces of $M$ and $A$ at $x$ are denoted by $T_xM$ and $T_xA$, respectively, then the \emph{normal space} of $A$ at $x$ is defined to be $N_xA := T_xM/T_xA$.

\begin{prop}[Functorial property of blowups for smooth or real analytic manifolds]
\label{prop:Functorial_property_blowup_real-analytic_smooth_manifold}
(See Arone and Kankaanrinta \cite[Theorem 4.1, p. 829]{Arone_Kankaanrinta_2010} and Akbulut and King \cite[Lemma 2.1, p. 51]{Akbulut_King_1985}.)
Let $M$ and $N$ be smooth (respectively, real analytic) manifolds with neat closed smooth (respectively, real analytic) submanifolds $A$ and $B$, respectively. Let $f : M \to N$ be a smooth (respectively, real analytic) map such that $f^{-1}(B) = A$. If the map $N_xA \to N_{f(x)}B$ induced by $f$ is injective for all $x \in A$, then $f$ induces a unique smooth (respectively, real analytic) map $\Bl_Z(f)$ between the blowups such that the following diagram commutes:
\begin{equation}
  \label{eq:Functorial_property_blowup_real-analytic_smooth_manifold}
  \begin{tikzcd}
    \Bl_A(M) \arrow[r, "\Bl_Z(f)"] \arrow[d, "\pi_M"] &\Bl_B(N) \arrow[d, "\pi_N"]
    \\
    M \arrow[r, "f"] &N
  \end{tikzcd}
\end{equation}
\end{prop}

Following Arone and Kankaanrinta \cite[Section 5, p. 831]{Arone_Kankaanrinta_2010}, one can apply Proposition \ref{prop:Functorial_property_blowup_real-analytic_smooth_manifold} to deduce equivariance properties of blowups of smooth or real analytic manifolds. Let $G$ be a Lie group and let $M$ be a smooth (respectively, real analytic manifold) on which $G$ acts smoothly (respectively, real analytically). The action of $G$ on $M$ is called \emph{proper} if the following map is proper:
\[
  G \times M \ni (g,x) \to (gx,x) \in M \times M
\]
It is well-known that the action of $G$ on $M$ is proper if and only if for every two
points $x$ and $y$ in $M$ there are neighborhoods $U$ and $V$ of $x$ and $y$, respectively,
such that the closure of the set
\[
  G(U,V) = \{g \in G | g(U) \cap V \neq \varnothing\}
\]
is compact. We recall the

\begin{thm}[Equivariance property of blowups for smooth or real analytic $G$-manifolds]
\label{thm:Equivariance_property_blowups_smooth_real-analytic_G-manifolds}  
(See Arone and Kankaanrinta \cite[Theorem 5.1, p. 831]{Arone_Kankaanrinta_2010}.)
Let $G$ be a Lie group and let $M$ be a smooth (respectively, real analytic) $G$-manifold. If $A$ is a closed, neat, smooth (respectively, real analytic) $G$-invariant submanifold of $M$, then the
blowup $\Bl_A(M)$ is a smooth (respectively, real analytic) $G$-manifold and the canonical projection
$\pi: \Bl_A(M) \to M$ is a smooth (respectively, real analytic) $G$-equivariant map. If the action of $G$ on $M$ is proper, then also the action of $G$ on $\Bl_A(M)$ is proper.
\end{thm}

We are not aware of prior statements of versions of Proposition \ref{prop:Functorial_property_blowup_real-analytic_smooth_manifold} or Theorem \ref{thm:Equivariance_property_blowups_smooth_real-analytic_G-manifolds} in the category of complex manifolds. However, we shall see that Proposition \ref{prop:Functorial_property_blowup_complex_manifold} below is a special case of the forthcoming, more general Proposition \ref{prop:Functorial_property_blowup_complex_analytic_space}. 

\begin{prop}[Functorial property of blowups for complex manifolds]
\label{prop:Functorial_property_blowup_complex_manifold}
Let $X$ and $Y$ be complex manifolds without boundary with closed, complex embedded submanifolds $Z$ and $W$, respectively.
If $f : X \to Y$ is a holomorphic map such that $f^{-1}(W) = Z$, then $f$ induces a unique holomorphic map $\Bl_Z(f)$ between the blowups such that the following diagram commutes:
\begin{equation}
  \label{eq:Functorial_property_blowup_complex_manifold}
  \begin{tikzcd}
    \Bl_Z(X) \arrow[r, "\Bl_Z(f)"] \arrow[d, "\pi_X"] &\Bl_W(Y) \arrow[d, "\pi_Y"]
    \\
    X \arrow[r, "f"] &Y
  \end{tikzcd}
\end{equation}
\end{prop}

While the proof of Theorem \ref{thm:Equivariance_property_blowups_complex_manifolds} below could be obtained by adapting the proof of Theorem \ref{thm:Equivariance_property_blowups_smooth_real-analytic_G-manifolds} for real analytic or smooth manifolds, we shall see in Section \ref{sec:Blowing_up_analytic_spaces} that it can also be proved as a corollary of Proposition \ref{prop:Functorial_property_blowup_complex_manifold}.

\begin{thm}[Equivariance property of blowups for complex $G$-manifolds]
\label{thm:Equivariance_property_blowups_complex_manifolds}  
Let $G$ be a Lie group and let $X$ be a complex $G$-manifold. If $Z$ is a closed, complex, embedded, $G$-invariant submanifold of $X$, then the blowup $\Bl_Z(X)$ is a complex $G$-manifold and the canonical projection $\pi: \Bl_Z(X) \to X$ is a holomorphic $G$-equivariant map. If the action of $G$ on $X$ is proper, then also the action of $G$ on $\Bl_Z(X)$ is proper.
\end{thm}

In the forthcoming Lemma \ref{lem:C*_equivariance_blowup_vector_space_along_linear_subspace}, we prove a simple special case of Theorem \ref{thm:Equivariance_property_blowups_complex_manifolds} by direct calculation for linear $G$ actions on a complex vector space when $G = S^1$ or $\CC^*$ and the blowup of the vector space along a $G$-invariant complex linear subspace. Before proceeding to the precise statement and proof, we recall some well-known facts from representation theory for Lie groups.

\begin{lem}[Classification of complex representations of the circle group]
\label{lem:Classification_complex_representations_circle_group}
(See, for example, Feehan and Leness \cite[Section 3.2]{Feehan_Leness_introduction_virtual_morse_theory_so3_monopoles}.)  
Let $V$ be a complex vector space of dimension $n\geq 1$ and $\rho:S^1 \to \GL(V)$ be a linear representation of the circle group. Then there is an $S^1$-equivariant isomorphism $\varphi:V\to\CC^n$ of complex vector spaces and a representation $\varrho: S^1 \to \U(n)$ such that $\rho = \varphi^{-1}\circ\varrho\circ\varphi$ and
\begin{equation}
  \label{eq:Circle_matrix_representation}
  \varrho(e^{i\theta})
  =
  \begin{pmatrix}
    e^{il_1\theta} &0            &\cdots &0 \\
    0            &e^{il_2\theta} &       & \\
    \vdots       &             &\ddots       &\vdots \\
    &            &      & \\
    0 &          &\cdots            & e^{il_n\theta}
  \end{pmatrix},
  \quad\text{for all } e^{i\theta} \in S^1,
\end{equation}
where the integers $l_k\in\ZZ$, for $k=1,\ldots,n$, are uniquely determined up to ordering by the representation $(\rho,V)$. The $S^1$-equivariant isomorphism $\varphi:V\to\CC^n$ of complex vector spaces determines an $S^1$-invariant, Hermitian inner product on $V$, with respect to which $\rho:S^1 \to \U(V)$ is a unitary representation.
\end{lem}

Lemma \ref{lem:Classification_complex_representations_circle_group} leads to the observations in the following

\begin{rmk}[Classification of complex representations of $\CC^*$]
\label{rmk:Classification_complex_representations_C*}
Let $V$ be a complex vector space of dimension $n\geq 1$ and $\rho_\CC:\CC^* \to \GL(V)$ be a linear representation of the group $\CC^*$. The standard embedding $S^1 \subset \CC^*$ and the representation $\rho_\CC:\CC^* \to \GL(V)$ uniquely determine a representation $\rho:S^1 \to \GL(V)$ and from Lemma \ref{lem:Classification_complex_representations_circle_group}, this is equivalent to a unitary, matrix representation $\varrho:S^1\to\U(n)$ of the form \eqref{eq:Circle_matrix_representation} via an $S^1$-equivariant isomorphism $\varphi:V\to\CC^n$ of complex vector spaces, so $\rho = \varphi^{-1}\circ\varrho\circ\varphi$. The unitary representation $\varrho$ admits a complexification with the same integer weights,
\begin{equation}
  \label{eq:C*_matrix_representation}
  \varrho_\CC(\lambda)
  =
  \begin{pmatrix}
    \lambda^{l_1} &0            &\cdots &0 \\
    0            &\lambda^{l_2} &       & \\
    \vdots       &             &\ddots       &\vdots \\
    &            &      & \\
    0 &          &\cdots            & \lambda^{l_n}
  \end{pmatrix},
  \quad\text{for all } \lambda \in \CC^*,
\end{equation}
and from the forthcoming Example \ref{exmp:Universal_complexification_circle_group}, this complexification is unique up to isomorphism of linear representations of $\CC^*$. In particular, after possibly replacing $\varphi$ by its composition with such an isomorphism, we see that $\rho_\CC$ and $\varrho_\CC$ are related by $\rho_\CC = \varphi^{-1}\circ\varrho_\CC\circ\varphi$, so $\varphi:V\to\CC^n$ becomes a $\CC^*$-equivariant isomorphism of complex vector spaces.
\end{rmk}

Before proceeding to the statement and proof of the forthcoming Lemma \ref{lem:C*_equivariance_blowup_vector_space_along_linear_subspace}, we recall the elementary

\begin{lem}[Group invariance of the orthogonal complement of an invariant subspace]
\label{lem:HilbertSpaceClosedUnderAction}
(See Br\"ocker and tom Dieck \cite[Chapter II, Proof of Proposition 1.9, p. 68]{BrockertomDieck} or Knapp \cite[Chapter IV, Equation (4.5), p. 239]{Knapp_1986}.)  
Let $\rho:G\to \U(H)$ be a unitary representation of a group $G$ on a complex Hilbert space $H$. If $V \subset H$ is a complex linear subspace that is $\rho(G)$-invariant, then the orthogonal complement $V^\perp$ of $V$ in $H$ is also $\rho(G)$-invariant.
\end{lem}

If one replaced the unitary group $\U(H)$ by its complexification $\GL(H)$ in Lemma \ref{lem:HilbertSpaceClosedUnderAction} then, naturally, one can no longer expect that if $V$ is $\rho(G)$-invariant, then $V^\perp$ will be $\rho(G)$-invariant. However, the following observation will suffice for our application.

\begin{lem}[$\CC^*$ invariance of an orthogonal complement of a $\CC^*$ invariant subspace]
\label{lem:HilbertSpaceClosedUnderC*Action}
Let $X$ be a finite-dimensional complex vector space and $\rho_\CC:\CC^* \to \GL(X)$ be a representation. If $Z\subset X$ is a $\CC^*$-invariant linear subspace, then there is a Hermitian inner product $\langle\cdot,\cdot\rangle$ on $X$ such that the orthogonal complement $Z^\perp$ of $Z$ in $X$ is also $\CC^*$-invariant.
\end{lem}

\begin{proof}
We claim that for some choice of Hermitian inner product $\langle\cdot,\cdot\rangle$ on $X$, we have
\begin{equation}
  \label{eq:rho_C_commutes_with_complex_conjugation}
  \langle w,\rho_\CC(\lambda)v\rangle = \langle \rho_\CC(\bar\lambda) w, v\rangle,
  \quad\text{for all } v, w \in X \text{ and } \lambda \in \CC^*.
\end{equation}
Let $\rho:S^1 \to \GL(X)$ be the linear representation determined by resttriction of $\rho_\CC$ to $S^1 \subset \CC^*$. By Lemma \ref{lem:Classification_complex_representations_circle_group}, the representation $\rho:S^1 \to \GL(X)$ has the form \eqref{eq:Circle_matrix_representation} and $X$ may be given an $S^1$-invariant Hermitian inner product $\langle\cdot,\cdot\rangle$ via the $S^1$-equivariant isomorphism $X \cong \CC^n$ and the standard Hermitian inner product on $\CC^n$. Let $\{e_1,\ldots,e_n\}$ be the orthogonal basis for $X$ induced by this isomorphism and the standard basis for $\CC^n$. By Remark \ref{rmk:Classification_complex_representations_C*}, the complexification $\rho_\CC$ of $\rho$ has the form \eqref{eq:C*_matrix_representation}. Thus, for any $v = \sum_{j=1}^na_je_j \in X$ and $\lambda \in \CC^*$, we see that
\[
  \rho_\CC(\lambda)v = \sum_{j=1}^na_j\lambda^{l_j} e_j
\]
and so, for any $w = \sum_{i=1}^nb_ie_i \in X$, we obtain
\begin{align*}
  \langle \rho_\CC(\lambda)^\dagger w, v\rangle
  &= \langle w,\rho_\CC(\lambda)v\rangle
  = \sum_{i,j=1}^n \langle b_ie_i, a_j\lambda^{l_j} e_j\rangle
  \\
  &= \sum_{i,j=1}^n b_i\bar\lambda^{l_j}\bar a_j \langle e_i,e_j\rangle
  = \sum_{i,j=1}^n b_i\bar\lambda^{l_j}\bar a_j \delta_{ij}
  \\
  &= \sum_{i=1}^n \bar\lambda^{l_i} b_i\bar a_i
  = \sum_{i,j=1}^n \langle b_i\bar\lambda^{l_i}e_i, a_je_j\rangle
  \\
  &= \langle \rho_\CC(\bar\lambda)w, v\rangle,
\end{align*}
and this proves \eqref{eq:rho_C_commutes_with_complex_conjugation}

Hence, if $v \in Z^\perp$ and $w \in Z$, then $\rho_\CC(\bar\lambda) w \in Z$ because $Z$ is $\CC^*$-invariant and therefore
\[
  \langle w,\rho_\CC(\lambda)v\rangle
  = \langle \rho_\CC(\lambda)^\dagger w, v\rangle
  = \langle \rho_\CC(\bar\lambda)w, v\rangle = 0,
\]
where the second equality follows from \eqref{eq:rho_C_commutes_with_complex_conjugation}. Since $w\in Z$ and $\lambda \in \CC^*$ are arbitrary, the linear subspace $Z^\perp \subset X$ is also $\CC^*$-invariant and we obtain a representation $\rho_\CC:\CC^* \to \Gl(Z^\perp)$. 
\end{proof}

Lemma \ref{lem:HilbertSpaceClosedUnderC*Action} is a special case of the following more general result in linear algebra.

\begin{lem}[Existence of an invariant complement for linear subspace that is invariant for a diagonalizable linear operator]
\label{lem:Existence_invariant_complement_invariant_linear_subspace}
Let $X$ be a finite-dimensional complex vector space and $T \in \End(X)$ be a linear operator. If $Z\subset X$ is a $T$-invariant linear subspace, that is, $T(Z) \subset Z$, then there exists a $T$-invariant linear subspace $Z' \subset X$ such that $X = Z\oplus Z'$ as a direct sum of complex vector spaces.
\end{lem}

\begin{lem}[$\CC^*$-equivariance property of the blowup of a complex vector space with a linear $\CC^*$ action and a $\CC^*$-invariant, linear subspace]
\label{lem:C*_equivariance_blowup_vector_space_along_linear_subspace}  
Let $X$ be a finite-dimensional, complex vector space, $Z \subsetneq X$ be a complex linear subspace, and $\rho_\CC:\CC^* \to \GL(X)$ be a homomorphism such that $Z$ is invariant under the induced action of $\CC^*$ on $X$. If $\pi:\Bl_Z(X)\to X$ is the blowup of $X$ along $Z$ constructed in Section \ref{sec:Blowups_analytic_manifolds_along_embedded_analytic_submanifolds}, then there is a unique homomorphism $\tilde\rho_\CC:\CC^* \to \Aut(\Bl_Z(X))$ such that, for each $\lambda\in\CC^*$ and corresponding $\rho_\CC(\lambda) \in \GL(X)$ and $\tilde\rho_\CC(\lambda) = \Bl_Z(\rho_\CC(\lambda)) \in \Aut(\Bl_Z(X))$, the following diagram commutes,
\begin{equation}
\label{eq:Commutative_diagram_pi_circ_rho_equals_tilderho_circ_pi}
\begin{tikzcd}
\Bl_Z(X) \arrow[d, "\pi"'] \arrow[r, "\tilde\rho_\CC(\lambda)"] &\Bl_Z(X) \arrow[d, "\pi"] 
\\
X \arrow[r, "\rho_\CC(\lambda)"'] &X
\end{tikzcd}
\end{equation}
and the blowup morphism $\pi$ is $\CC^*$-equivariant. 
\end{lem}  

\begin{proof}
By hypothesis, the subspace $Z \subset X$ is $\CC^*$-invariant and thus also the orthogonal complement $Z^\perp \subset X$ is $\CC^*$-invariant by Lemma \ref{lem:HilbertSpaceClosedUnderC*Action}. Hence, there are uniquely determined $\CC^*$ actions on $\PP(Z^\perp)$ and $\PP(Z^\perp)\times X$. The incidence relations \eqref{eq:Blowup_linear_subspace_incidence_variety} are invariant under the linear action of $\CC^*$ on $X$ and the subspaces $Z$ and $Z^\perp$ and hence the $\CC^*$ action on $\PP(Z^\perp)\times X$ restricts to a $\CC^*$ action on the embedded complex submanifold $\Bl_Z(X) \subset \PP(Z^\perp)\times X$. 
\end{proof}

\begin{rmk}[Equivariant tubular neighborhoods for invariant submanifolds of $G$-manifolds]
\label{rmk:Equivariant_tubular_neighborhoods_invariant_submanifolds_G-manifolds}  
Given a unitary linear representation $\rho:S^1 \to \U(n)$ of the circle or its complexification $\rho_\CC:\CC^* \to \GL(n,\CC)$ whose actions on $\CC^n$ preserve a coordinate subspace $\CC^m \hookrightarrow \CC^n$ and hence preserve its normal space $\CC^n/\CC^m \cong \CC^{n-m}$, the actions of $S^1$ and $\CC^*$ naturally lift to the blowup $\Bl_{\CC^m}(\CC^n)$. Furthermore, given a linear representation $\rho:G \to \GL(n,\CC)$ of a Lie group, a similar (if less transparent) argument would show that the action of $G$ on $\CC^n$ naturally lifts to the blowup $\Bl_{\CC^m}(\CC^n)$ of $\CC^n$ along a $G$-invariant coordinate subspace $\CC^m$ and that the normal space $\CC^n/\CC^m$ is $G$-invariant.
	
	More generally, given an orthogonal, smooth action of a compact Lie group $G$ on a smooth Riemannian manifold $M$ and a $G$-invariant embedded submanifold $S \hookrightarrow M$, we recall from Bredon \cite[Section VI.2, Theorem 2.2, p. 306]{Bredon} that $S$ has an open, $G$-invariant tubular neighborhood constructed via an open, $G$-invariant neighborhood of the zero section $S$ of the normal bundle $N_{M/S} \cong (TM\restriction S)/TS$ of $S$ in $M$. The normal bundle $N_{M/S}$ is a smooth $G$-vector bundle in the sense of Bredon \cite[Section VI.2, p. 303]{Bredon}, that is, the group $G$ acts smoothly on the total space $N_{M/S}$, linearly on the fibers of $N_{M/S}$, and the bundle projection $\pi:N_{M/S} \to S$ is $G$-equivariant. See Kankaanrinta \cite{Kankaanrinta_2007} for extensions of these results to proper actions of Lie groups on smooth manifolds.
\end{rmk}

\section{Fiber products of analytic spaces}
\label{sec:Fiber_products_analytic_spaces}
We closely follow Fischer \cite[Section 0.25, p. 21]{Fischer_complex_analytic_geometry} (but see also Grauert and Remmert \cite[Section 1.2.6, pp. 18--20]{Grauert_Remmert_coherent_analytic_sheaves}), who assumes that $\KK=\CC$ whereas we also allow $\KK=\RR$. Let $(X_j,\sO_{X_j})$ for $j=1,2$ and $(Y,\sO_Y)$ be $\KK$-analytic spaces and let $\varphi_j:X_j\to Y$ for $j=1,2$ be morphisms of $\KK$-analytic spaces. A $\KK$-analytic space $X_1\times_Y X_2$ together with morphisms of $\KK$-analytic spaces for $j=1,2$,
\[
  \pr_j:X_1\times_Y X_2 \to X_j
\]
such that
\[
  \varphi_1\circ\pr_1 = \varphi_2\circ\pr_2
\]
is called a \emph{fiber product} of $X_1$ and $X_2$ over $Y$ (or over $\varphi_1$ and $\varphi_2$) if it has the following universal property: Given any $\KK$-analytic space $X$ together with morphisms of $\KK$-analytic spaces $\psi_j:X \to X_j$ for $j=1,2$ such that
\[
  \varphi_1\circ\psi_1 = \varphi_2\circ\psi_2,
\]
there is a unique morphism $\psi:X\to X_1\times_Y X_2$ of $\KK$-analytic spaces such that the following diagram commutes:
\begin{equation}
  \label{eq:Fiber_products_analytic_spaces}
  \begin{tikzcd}
  X
  \arrow[drr, bend left, "\psi_1"]
  \arrow[ddr, bend right, "\psi_2"']
  \arrow[dr, dotted, "\psi" description] & & \\
    & X_1\times_Y X_2 \arrow[r, "\pr_1"'] \arrow[d, "\pr_2"]
      & X_1 \arrow[d, "\varphi_1"] \\
& X_2 \arrow[r, "\varphi_2"'] &Y
\end{tikzcd}
\end{equation}
A commutative diagram morphisms of $\KK$-analytic spaces
\begin{equation}
  \label{eq:Cartesian_square_analytic_spaces}
  \begin{tikzcd}
   X \arrow[r, " "] \arrow[d, " "] & X_1 \arrow[d, " "] \\
  X_2 \arrow[r, " "] &Y
\end{tikzcd}
\end{equation}
is called a \emph{Cartesian square} if $X$ is a fiber product of $X_1$ and $X_2$ over $Y$.

When $Y$ is a simple point, a fiber product of $X_1$ and $X_2$ is called a \emph{direct product} and denoted by $X_1 \times X_2$. Its universal property then is the following: Given any $\KK$-analytic space $X$ together with morphisms $\psi_j:X \to X_j$ of $\KK$-analytic spaces for $j=1,2$, there is a unique morphism $\psi=(\psi_1,\psi_2):X \to X_1 \times X_2$ of $\KK$-analytic spaces such that the following diagram commutes:
\[
  \begin{tikzcd}
  &X_1 \\
  X \arrow[ur, "\psi_1"] \arrow[dr, "\psi_2"'] \arrow[r, dotted, "\psi" description]
  & X_1 \times X_2 \arrow[u, "\pr_1"'] \arrow[d, "\pr_2"] \\
  & X_2
  \end{tikzcd}
\]  
An exposition of fiber products, for $\KK=\CC$, is given by Fischer in \cite[Sections 0.26--0.32, pp. 22--30]{Fischer_complex_analytic_geometry}. The universal property implies the uniqueness of a fiber product up to isomorphisms. Fischer's development does not require $\KK$ to be algebraically closed and is valid for $\KK=\RR$ as well.

\begin{prop}[Existence of fiber products]
\label{prop:Fiber_product_closed_subspace_product}  
(See Fischer \cite[Section 0.32, Corollary, p. 29]{Fischer_complex_analytic_geometry}.)
Let $(X_j,\sO_{X_j})$, for $j=1,2$, and $(Y,\sO_Y)$ be $\KK$-analytic spaces and $\varphi_j:X_j \to Y$, for $j=1,2$, be $\KK$-analytic morphisms. Then the following hold:
\begin{enumerate}
\item The fiber product $X_1\times_Y X_2$ exists and is a closed subspace of $X_1\times X_2$.
\item For any point $p \in X_2$, the fiber $(X_1 \times_Y X_2)_p$ of $X_1\times_Y X_2$ over $p$ is isomorphic to the fiber $(X_1)_{\varphi_2(p)}$.
\end{enumerate}
\end{prop}

In the category of sets, we have
\[
  X_1\times_Y X_2 = \left\{ (x_1,x_2) \in X_1\times X_2: \varphi_1(x_1) = \varphi_2(x_2) \right\},
\]
just as in Section \ref{sec:Fiber_products_base_change} in the category of algebraic varieties.

\section{Flatness of holomorphic maps of complex analytic spaces}
\label{sec:Flatness_holomorphic_maps_complex_analytic_spaces}
For further details concerning the concept of flatness in the category of complex analytic spaces, we refer the reader to Douady \cite{Douady_1968enseignement_math} and Fischer \cite[Sections 3.10--3.22, pp. 146--161]{Fischer_complex_analytic_geometry}. For the definitions of a flat morphism $f:X\to Y$ of ringed spaces and a flat $\sO_X$-module $\sF$, see Definition \ref{defn:Flat_modules_and_flat_morphisms_ringed_spaces}. In particular, these definitions apply in the category of holomorphic maps $f$ of complex analytic spaces $(X,\sO_X$ and $(Y,\sO_Y)$ and to $\sO_X$-modules $\sF$ (see Fischer \cite[Section 3.11, p. 146]{Fischer_complex_analytic_geometry}).

The following algebraic lemma will be useful in subsequent discussions of the inverse images of hypersurfaces by flat holomorphic maps. Recall that an element $a$ of a ring $R$ is a \emph{zero divisor} if the map $R \ni x \mapsto ax \in R$ is not injective, that is, there exists a non-zero element $y \in R$ such that $ya = 0$. Recall furthermore that a \emph{local homomorphism of local rings} is a homomorphism $\varphi:R\to S$ of rings such that $R$ and $S$ are local rings and $\varphi(\fm_R) \subset \fm_S$, where $\fm_R$ and $\fm_S$ are the maximal ideals in $R$ and $S$, respectively (see the Stacks Project \cite[\href{https://stacks.math.columbia.edu/tag/07BI}{Definition 07BI}]{stacks-project}).

\begin{lem}[Homomorphisms of rings and flatness]
\label{lem:Local_homomorphism_rings}  
(See Fischer \cite[Section 3.12, Lemma, p. 148]{Fischer_complex_analytic_geometry}.)
If $\varphi:R \to S$ is a homomorphism of rings such that $S$ is $R$-flat\footnote{That is, $S$ is a flat $R$-module, where the homomorphism $\varphi:R\to S$ is used to give $S$ the structure of an $R$-module.} and $a \in R$, then the following hold:
\begin{enumerate}
\item\label{item:varphi(a)_zero_divisor_implies_a_zero_divisor}
  If $\varphi(a)$ is a zero divisor, then $a$ is a zero divisor.
\item\label{item:a_zero_divisor_implies_varphi(a)_zero_divisor}
  If in addition $\varphi:R \to S$ is a local homomorphism of local rings, then $\varphi$ is injective
and thus if $a$ is a zero divisor, then $\varphi(a)$ is a zero divisor.
\end{enumerate}
\end{lem}

\begin{rmk}[Application of Lemma \ref{lem:Local_homomorphism_rings}]
\label{rmk:Local_homomorphism_rings_application}
In our work, we shall appeal to Lemma \ref{lem:Local_homomorphism_rings} in the following setting. Assume that $(f,f^\sharp):(X,\sO_X)\to (Y,\sO_Y)$ is a flat morphism of $\KK$-analytic spaces. For a point $p \in X$, let $f_p^\sharp:\sO_{Y,f(p)} \to \sO_{X,p}$ be the canonical homomorphism of local rings. Because $f$ is flat, $\sO_{X,p}$ is a flat $\sO_{Y,f(p)}$-module by Definition \ref{defn:Flat_modules_and_flat_morphisms_ringed_spaces}, where the homomorphism $f_p^\sharp:\sO_{Y,f(p)} \to \sO_{X,p}$ provides $\sO_{X,p}$ with the structure of an $\sO_{Y,f(p)}$-module. Item \eqref{item:varphi(a)_zero_divisor_implies_a_zero_divisor} in Lemma \ref{lem:Local_homomorphism_rings} implies that if $h \in \sO_{Y,f(p)}$ is a \emph{not} a zero divisor, then $f_p^\sharp(h) \in \sO_{X,p}$ is \emph{not} a zero divisor. (Although we shall not need this additional fact in our work, we note that $f_p^\sharp:\sO_{Y,f(p)} \to \sO_{X,p}$ is a local homomorphism of local rings, and so if $g \in \sO_{Y,f(p)}$ and $f_p^\sharp(g) \in \sO_{X,p}$ is \emph{not} a zero divisor, then $g$ is a \emph{not} a zero divisor.) 
\end{rmk}  

As in Section \ref{sec:Functorial_properties_blowups_schemes}, the property of flatness of a holomorphic map is important because it ensures strong continuity conditions upon the fibration induced by that map. In this context, one has the

\begin{prop}[Flatness of the canonical projection from a product of complex analytic spaces]
\label{prop:Canonical_projection_product_complex_analytic_spaces_flat}
(See Fischer \cite[Section 3.17, p. 155]{Fischer_complex_analytic_geometry}.)
If $(X,\sO_X)$ and $(Y,\sO_Y)$ are complex analytic spaces, then the canonical projection $\pr_Y:X\times Y \to Y$ is flat.
\end{prop}

As in the category of schemes, the property that flatness is preserved by base change will be used in our discussion of blowups of complex analytic spaces.

\begin{prop}[Properties of flat holomorphic maps of complex analytic spaces]
\label{prop:Properties_flat_holomorphic_maps_complex_analytic_spaces}
(See Fischer \cite[Section 3.15, Proposition, p. 152]{Fischer_complex_analytic_geometry}.)
Assume that we are given a Cartesian diagram of holomorphic maps of complex analytic spaces as in \eqref{eq:Cartesian_square_analytic_spaces}
\[
  \begin{tikzcd}
   X' \arrow[r, "\chi"] \arrow[d, "{f'}"] & X \arrow[d, "f"] \\
   Y' \arrow[r, "\psi"] &Y
\end{tikzcd}
\]
and thus $X'$ is the fiber product $X\times_Y Y'$ and $f'$ is the base change of $f$ with respect to $\psi$ (as in Section \ref{sec:Fiber_products_base_change}).
Let $p' \in X'$ and $p := \chi(p') \in X$. If $f$ is \emph{flat in $p \in X$} (that is, $\sO_{X,p}$ is a flat $\sO_{Y,f(p)}$ module), then $f'$ is flat in $p' \in X'$ (that is, $\sO_{X',p'}$ is a flat $\sO_{Y',f'(p')}$ module). In particular, if $f:X\to Y$ is flat, then $f':X'\to Y'$ is flat. 
\end{prop}

There are sufficient conditions for a holomorphic map of complex analytic spaces to be flat.

\begin{prop}[Flatness of finite holomorphic maps of complex analytic spaces]
\label{prop:Flatness_finite_holomorphic_maps_complex_analytic_spaces}
(See Fischer \cite[Section 3.13, Proposition, p. 149]{Fischer_complex_analytic_geometry}.)
Let $(X,\sO_X)$ and $(Y,\sO_Y)$ be complex analytic spaces and $f:X\to Y$ be a finite holomorphic map. Then $f$ is flat if and only if $f_*\sO_X$ is locally free.
\end{prop}

\begin{prop}[Flatness of holomorphic maps of complex manifolds]
\label{prop:Flatness_holomorphic_maps_complex_manifolds}
(See Fischer \cite[Section 3.20, Corollary, p. 158]{Fischer_complex_analytic_geometry}.)
Let $X$ and $Y$ be connected, complex manifolds and $f:X\to Y$ be a holomorphic map. Then the following are equivalent:
\begin{enumerate}
\item $f$ is flat.
\item $f$ is open.
\item Every fiber of $f$ is of pure dimension $\dim X - \dim Y$.
\end{enumerate}
\end{prop}

We refer the reader to Fischer \cite[Chapter 3]{Fischer_complex_analytic_geometry} for additional sufficient conditions for a holomorphic map of complex analytic spaces to be flat.

\section{Blowups of analytic spaces}
\label{sec:Blowing_up_analytic_spaces}
Following Griffiths and Harris \cite[Chapter 0, Section 1, pp. 12--14, Section 2, pp. 20--22, and Chapter 1, Section 1, pp. 129--130]{GriffithsHarris} (who consider \emph{complex} analytic subvarieties of smooth \emph{complex} manifolds), let $M$ be a (real or complex) analytic (not necessarily compact) manifold of dimension $d \geq 1$ and $V \subset M$ be a closed analytic subspace. A closed analytic subspace $V \subset M$ is called \emph{irreducible} if $V$ cannot be written as the union of two closed analytic subspace, $V_1, V_2 \subset M$, with $V_i \neq V$ for $i=1,2$. One calls $V \subset M$ a \emph{closed analytic subspace of dimension $d-1$} if $V$ is a \emph{analytic hypersurface}, that is, for any point $p \in V$, then $U\cap V = f^{-1}(0)$, for some open neighborhood, $U \subset M$ of $p$, and some analytic function, $f$, on $U$ \cite[Chapter 1, Section 2, p. 20]{GriffithsHarris}. We then recall the

\begin{defn}[Divisor on an analytic manifold]
(See Griffiths and Harris \cite[Chapter 1, Section 1, p. 130]{GriffithsHarris}.)
A \emph{divisor} $D$ on an analytic manifold $M$ is a locally finite, formal linear combination,
\[
D = \sum_{i} a_iV_i,
\]
of irreducible, analytic hypersurfaces of $M$, where $a_i \in \ZZ$.
\end{defn}

For blowing up in the category of complex analytic spaces, we refer to Fischer \cite{Fischer_complex_analytic_geometry} and for techniques to adapt constructions in the category of complex analytic spaces to the category of real analytic spaces, we refer to Guaraldo, Macr\`\i, and Tancredi \cite{Guaraldo_Macri_Tancredi_topics_real_analytic_spaces}. For blowing up in the category of complex  complex manifolds, we refer to Griffiths and Harris \cite[Chapter 4, Section 6, pp. 602--605]{GriffithsHarris}, Huybrechts \cite[Section 2.5]{Huybrechts_2005}, and Voisin \cite[Section 3.3.3]{Voisin_hodge_theory_complex_algebraic_geometry_I}.

\begin{defn}[Hypersurface in a complex analytic space]
\label{defn:Hypersurface_complex_analytic_space}  
(See Fischer \cite[Section 0.45, p. 42]{Fischer_complex_analytic_geometry}.)  
Let $X$ be a complex analytic space with a closed subspace $Y \subset X$, defined by a coherent ideal $\sI \subset \sO_X$. Then $Y$ is a \emph{hypersurface} in $X$ if for every $p \in Y$ there is a non-zero-divisor $h \in \sO_{X,p}$ such that $\sI_p = h\cdot \sO_{X,p}$.
\end{defn}

Recall that an $\sO_X$-module $\sF$ is called \emph{invertible} if for every $p \in X$ there are an open neighborhood $U\subset X$ of $p$ and an isomorphism $\sO_X\restriction U \to \sF\restriction U$ \cite[Section 0.45, p. 42]{Fischer_complex_analytic_geometry}. Then $Y$ is a hypersurface if and only if $\sI$ is invertible (see Fischer \cite[Section 0.45, p. 42]{Fischer_complex_analytic_geometry}). The Active Lemma \ref{lem:Active} implies that
\[
  \dim_pY = \dim_pX - 1
\]
for every $p \in Y$, when $Y$ is a hypersurface in $X$, but the converse is false (see Fischer \cite[Section 0.45, Example, p. 42]{Fischer_complex_analytic_geometry}). We recall the following generalization of Theorem \ref{prop:Huybrechts_2-5-3} from the category of complex manifolds and complex embedded submanifolds to complex analytic spaces and closed complex analytic subspaces and analogue of Definition \ref{defn:Eisenbud_Harris_4-16} in the categories of algebraic varieties and schemes.

\begin{thm}[Blowup of a complex analytic space]
\label{thm:Fischer_4-1}  
(See Fischer \cite[Theorem 4.1, p. 162]{Fischer_complex_analytic_geometry}, following Douady \cite{Douady_1968seminaire_bourbaki} and Hironaka and Rossi \cite{Hironaka_Rossi_1964}; see also Hironaka \cite[Chapter 0, Section 2, p. 129]{Hironaka_1964-I-II}.)  
Let $(X,\sO_X)$ be a complex analytic space with a closed complex analytic subspace\footnote{As usual, called the \emph{center}.} $Z \hookrightarrow X$. Then there exist a complex analytic space $\Bl_Z(X)$ and a holomorphic map $\pi:\Bl_Z(X) \to X$ with the following properties:
\begin{enumerate}
\item\label{item:Fischer_theorem_4-1_a} $\pi$ is proper.
\item\label{item:Fischer_theorem_4-1_b} $E := \pi^{-1}(Z) \hookrightarrow \Bl_Z(X)$ is a hypersurface as in Definition \ref{defn:Hypersurface_complex_analytic_space}.
\item\label{item:Fischer_theorem_4-1_c} $\pi$ is universal with respect to Item \eqref{item:Fischer_theorem_4-1_b}, that is, if $\tau: W \to X$ is a holomorphic map of complex analytic spaces such that $\tau^{-1}(Z) \hookrightarrow W$ is a hypersurface, then there is a unique holomorphic map $\varphi: W \to \Bl_Z(X)$ such that the following diagram commutes:
  \[
    \begin{tikzcd}
    W \arrow[r, "\varphi", dashed] \arrow[rd, "\tau"'] &\Bl_Z(X) \arrow[d, "\pi"]
    \\
    &X
  \end{tikzcd}
  \]
\item\label{item:Fischer_theorem_4-1_d} The restriction of $\pi$ to $\Bl_Z(X)\less E \to X\less Z$ is biholomorphic.
\item\label{item:Fischer_theorem_4-1_e} If $X$ is a manifold and $Z$ is a submanifold, then $\Bl_Z(X)$ is a manifold.
\end{enumerate}
\end{thm}

Fischer notes that blowing up is the principal tool in the resolution of singularities of complex spaces \cite[p. 168]{Fischer_complex_analytic_geometry}. If the hypothesis that $Z$ in Item \ref{item:Fischer_theorem_4-1_e} of Theorem \ref{thm:Fischer_4-1} be a submanifold is relaxed to allow $Z$ to be an arbitrary complex analytic subspace, then $\widetilde X$ may have singularities (see Fischer \cite[p. 168]{Fischer_complex_analytic_geometry}).

By analogy with the uniqueness assertion in Remark \ref{rmk:Existence_blowup}, we have the

\begin{rmk}[Uniqueness of the blowup in Theorem \ref{thm:Fischer_4-1}]
\label{rmk:Uniqueness_blowup_analytic_spaces}
The universal properties in Theorem \ref{thm:Fischer_4-1} \emph{uniquely} characterize the blowup $\pi: \Bl_Z(X) \to X$ of complex analytic space along a closed subpace. To see this, suppose that $\varpi:\widetilde X \to X$ is a holomorphic map of analytic spaces that also has the properties in Items \eqref{item:Fischer_theorem_4-1_b} and \eqref{item:Fischer_theorem_4-1_c}. Hence, there are unique holomorphic maps $\varphi: \Bl_Z(X) \to \widetilde X$ and $\psi: \widetilde X \to \Bl_Z(X)$ such that the following diagram commutes:
\[
  \begin{tikzcd}[every arrow/.append style={shift left}]
   \widetilde X \arrow[d, "\varpi"'] \arrow[r, "\psi", dashed] &\Bl_Z(X) \arrow[l,"\varphi",dashed] \arrow[d, "\pi"]
   \\
   X  \arrow[r, "\id"] &X
  \end{tikzcd}
\]
Therefore, by choosing the map $\tau:W \to X$ in Item \eqref{item:Fischer_theorem_4-1_c} to be $\pi: \Bl_Z(X) \to X$, its uniqueness assertion implies that $\psi\circ\varphi = \id_{\Bl_Z(X)}$ and by choosing the map $\tau:W \to X$ in Item \eqref{item:Fischer_theorem_4-1_c} to be $\varpi:\widetilde X \to X$, its uniqueness assertion implies that $\varphi\circ\psi = \id_{\widetilde X}$. Thus, we obtain a commutative diagram
\[
  \begin{tikzcd}
   \widetilde X \arrow[rd, "\varpi"'] \arrow[r, "\psi"] &\Bl_Z(X) \arrow[d, "\pi"]
   \\
   &X
  \end{tikzcd}
\]
where $\psi:\widetilde X \cong \Bl_Z(X)$ is an isomorphism of complex analytic spaces.
\end{rmk}

\begin{rmk}[On the definition of blowups in complex analytic spaces]
\label{rmk:Fischer_4-1}
Fischer proves Theorem \ref{thm:Fischer_4-1}, in particular the existence of the blowup $\pi:\widetilde X \to X$, using calculations in local coordinates. By comparing the combination of Items \eqref{item:Fischer_theorem_4-1_b} and \eqref{item:Fischer_theorem_4-1_c} in Theorem \ref{thm:Fischer_4-1} with the Definition \ref{defn:Eisenbud_Harris_4-16} for the blowup in the category of schemes, we see that this combination could serve as the analogous definition of blowup in the category of complex analytic spaces, with uniqueness of the blowup following immediately just as in the case of schemes while existence is proved by Fisher in \cite[pp. 162--168]{Fischer_complex_analytic_geometry}.
\end{rmk}

We recall the

\begin{prop}[Existence of the inverse image of a $\KK$-analytic subspace]
\label{prop:Fischer_0-27}  
(See Fischer \cite[Proposition 0.27, p. 23]{Fischer_complex_analytic_geometry} and Grauert and Remmert \cite[Section 1.2.6, p. 19]{Grauert_Remmert_coherent_analytic_sheaves} for $\KK=\CC$.)
Let $\KK=\RR$ or $\CC$ and $(X,\sO_X)$ and $(Y,\sO_Y)$ be $\KK$-analytic spaces, $\varphi:Y\to X$ be a morphism of $\KK$-analytic spaces, and $Z$ be an open (respectively, closed) $\KK$-analytic subspace of $X$. Then there is a Cartesian square (as in Section \ref{sec:Fiber_products_analytic_spaces})
\[
  \begin{tikzcd}
    \varphi^{-1}(Z) \arrow[r, "\jmath"] \arrow[d, "\psi"'] &Y \arrow[d, "\varphi"]
    \\
    Z \arrow[r, "\iota"] &X
  \end{tikzcd}
\]
where $\varphi^{-1}(Z)$ is an open (respectively, closed) $\KK$-analytic subspace of $Y$ and $\jmath, \iota$
denote the canonical injections. In particular, there is an morphism of $\KK$-analytic spaces:
\[
  \sO_{\varphi^{-1}(Z)} \cong \psi^*\sO_Z.
\]
\end{prop}

Recall that if $(f,f^\sharp):(Y,\sO_Y) \to (X,\sO_X)$ is a morphism of $\KK$-analytic spaces and $\sF$ is an $\sO_X$-module, then the image sheaf $f_*\sF$ is well-defined and an $f_*\sO_X$-sheaf; via the map $f^\sharp:\sO_X\to f_*\sO_Y$, the sheaf becomes an $\sO_X$-module, called the \emph{analytic image sheaf} of $\sF$ (with respect to the map $(f,f^\sharp)$ (for further details, see Fischer \cite[Section 0.2, p. 6]{Fischer_complex_analytic_geometry} and Grauert and Remmert \cite[Section 1.2.5, pp. 18--19]{Grauert_Remmert_coherent_analytic_sheaves} for $\KK=\CC$).

Furthermore, we recall that there is  natural construction leading to an adjoint of the functor $f_*$. Given an $\sO_Y$-module $\sG$, one may construct an $\sO_X$-module $f^*\sG$, called the \emph{analytic inverse image sheaf} of $\sG$ (for further details, see Fischer \cite[Section 0.10, p. 6]{Fischer_complex_analytic_geometry} and Grauert and Remmert \cite[Section 1.2.6, pp. 18--19]{Grauert_Remmert_coherent_analytic_sheaves} for $\KK=\CC$). One can show that \cite[Section 0.10, p. 6]{Fischer_complex_analytic_geometry}
\[
  \Hom_X(f^*\sG,\sF) = \Hom_Y(\sG,f_*\sF).
\]
We refer the reader to the indicated references for further properties.

\begin{defn}[Inverse image of a $\KK$-analytic subspace]
\label{defn:Fischer_0-27} 
In Proposition \ref{prop:Fischer_0-27}, the $\KK$-analytic space $\varphi^{-1}(Z)$ is called the \emph{inverse image} of $Z$. If $Z = (\{p\},\CC)$ is a simple point, then $\varphi^{-1}(Z)$ is called the \emph{fiber} of $\varphi$ over $p$ and denoted by $Y_p$.
\end{defn}  

Theorem \ref{thm:Fischer_4-1} and Proposition \ref{prop:Fischer_0-27} yield the following analogue in the category of complex analytic spaces of Proposition \ref{prop:Functorial_property_blowup_scheme} (for schemes), Proposition \ref{prop:Functorial_property_blowup_real-analytic_smooth_manifold} (for real analytic manifolds), and Proposition \ref{prop:Functorial_property_blowup_complex_manifold} (for complex manifolds).

\begin{prop}[Functorial property of blowups for complex analytic spaces]
\label{prop:Functorial_property_blowup_complex_analytic_space}
Let $(Y,\sO_Y)$ and $(X,\sO_X)$ be complex analytic spaces and let $(Z,\sO_Z)$ be a closed, complex analytic subspace of $(X,\sO_X)$. If $f : Y \to X$ is a morphism of complex analytic spaces, then the following hold:
\begin{enumerate}
\item\label{item:Functorial_property_blowup_complex_analytic_space} There exists a unique morphism\footnote{See Definition \ref{defn:Fischer_0-27} for the construction of $f^{-1}(Z)$.} $\Bl_Z(f) :  
\Bl_{f^{-1}(Z)}(Y)\to \Bl_Z(X)$ such that the following diagram commutes:
\begin{equation}
  \label{eq:Functorial_property_blowup_complex_analytic_space}
  \begin{tikzcd}
    \Bl_{f^{-1}(Z)}(Y) \arrow[r, "\Bl_Z(f)"] \arrow[d, "\pi_Y"] &\Bl_Z(X) \arrow[d, "\pi_X"]
    \\
    Y \arrow[r, "f"] &X
  \end{tikzcd}
\end{equation}
\item\label{item:Blowup_commutes_flat_base_change_complex_analytic_space} If $f$ is flat in the sense of Definition \ref{defn:Flat_modules_and_flat_morphisms_ringed_spaces}, then the diagram \eqref{eq:Functorial_property_blowup_complex_analytic_space} is Cartesian as in \eqref{eq:Cartesian_square_analytic_spaces}, that is, the induced morphism to the fiber product (see Section \ref{sec:Fiber_products_analytic_spaces})
\[
  r:\Bl_{f^{-1}(Z)}(Y) \to \Bl_Z(X) \times_X Y
\]
is an isomorphism of complex analytic spaces such that $\pi_Y = \pr_2\circ r$, where $\pr_2:\Bl_Z(X) \times_X Y \to Y$ is projection onto the second factor.
\end{enumerate}
\end{prop}

\begin{proof}
We shall adapt the proofs of Proposition \ref{prop:Functorial_property_blowup_scheme} due to G\"ortz and Wedhorn and to Eisenbud and Harris in the category of schemes and due to Hauser in the category of algebraic varieties. We provide details in addition to those given by the preceding authors.

Proposition \ref{prop:Fischer_0-27} yields a closed, complex analytic subspace $f^{-1}(Z)$ of $(Y,\sO_Y)$, so Theorem \ref{thm:Fischer_4-1} provides a complex analytic space $\Bl_{f^{-1}(Z)}(Y)$, a holomorphic map $\pi_Y:\Bl_{f^{-1}(Z)}(Y) \to Y$ comprising the blowup of $Y$ along the center $f^{-1}(Z)$, and a hypersurface as in Definition \ref{defn:Hypersurface_complex_analytic_space}:
\[
  \pi_Y^{-1}(f^{-1}(Z)) = (f\circ\pi_Y)^{-1}(Z) \hookrightarrow \Bl_{f^{-1}(Z)}(Y).
\]
Hence, Item \eqref{item:Fischer_theorem_4-1_c} in Theorem \ref{thm:Fischer_4-1} gives us a unique holomorphic map $\varphi: \Bl_{f^{-1}(Z)}(Y) \to \Bl_Z(X)$ such that the following diagram commutes:
\[
  \begin{tikzcd}
   \Bl_{f^{-1}(Z)}(Y) \arrow[r, "\varphi", dashed] \arrow[rd, "f\circ\pi_Y"'] &\Bl_Z(X) \arrow[d, "\pi_X"]
   \\
   &Y
  \end{tikzcd}
\]
We choose $\Bl_Z(f) := \varphi$ and thus obtain the commutative diagram \eqref{eq:Functorial_property_blowup_complex_analytic_space}.

We will show that the diagram \eqref{eq:Functorial_property_blowup_complex_analytic_space} is Cartesian by proving that the following induced morphism is an isomorphism:
\begin{equation}
  \label{eq:r_morphism_BlfinvZ(Y)_to_BlZ(X)_timesX_Y}
  r: \Bl_{f^{-1}(Z)}(Y) \to \Bl_Z(X)\times_X Y.
\end{equation}
The morphism $r$ is induced by applying the universal property of fiber products of analytic spaces (see Section \ref{sec:Fiber_products_analytic_spaces}) summarized in the diagram \eqref{eq:Fiber_products_analytic_spaces} to give the following commutative diagram:
\begin{equation}
  \label{eq:Fiber_products_and_blowups_analytic_spaces}
  \begin{tikzcd}
  \Bl_{f^{-1}(Z)}(Y)
  \arrow[drr, bend left, "\Bl_Z(f)"]
  \arrow[ddr, bend right, "\pi_Y"']
  \arrow[dr, dotted, "r" description] & & \\
    & \Bl_Z(X)\times_X Y \arrow[r, "\pr_1"'] \arrow[d, "\pr_2"]
      & \Bl_Z(X) \arrow[d, "\pi_X"] \\
& Y \arrow[r, "f"'] &X
\end{tikzcd}
\end{equation}
Here, $\pr_1:\Bl_Z(X)\times_X Y\to \Bl_Z(X)$ denotes the projection onto the first factor and $\pr_2:\Bl_Z(X)\times_X Y\to Y$ denotes the projection onto the second factor in the fiber product.

According to \eqref{eq:Fiber_products_analytic_spaces}, the following diagram commutes:
\begin{equation}
  \label{eq:Fiber_product_blowupX_timesX_Y}
  \begin{tikzcd}
  \Bl_Z(X)\times_X Y \arrow[r, "\pr_1"] \arrow[d, "\pr_2"] &\Bl_Z(X) \arrow[d, "\pi_X"] \\
  Y \arrow[r, "f"] &X
  \end{tikzcd}
\end{equation}
Since $f$ is flat by hypothesis, Proposition \ref{prop:Properties_flat_holomorphic_maps_complex_analytic_spaces} applies to the diagram \eqref{eq:Fiber_product_blowupX_timesX_Y} to show that $\pr_1$ is flat.

Let $E_X = \pi_X^{-1}(Z)$ denote the exceptional divisor on $X$. Because $E_X \subset \Bl_Z(X)$ is a hypersurface in the sense of Definition \ref{defn:Hypersurface_complex_analytic_space} by Item \eqref{item:Fischer_theorem_4-1_b} in Theorem \ref{thm:Fischer_4-1}, then $E_X$ is defined by a locally principal ideal $\sI_{E_X}$ in $\sO_{\Bl_Z(X)}$, and for any point $p \in E_X \subset \Bl_Z(X)$, there are an open neighborhood $U \subset \Bl_Z(X)$ of $p$ and an element $h \in \sO_{\Bl_Z(X)}\restriction U$ such that $\sI_{E_X}\restriction U = h\sO_{\Bl_Z(X)}\restriction U$, where (the germ at $p$ of) $h$ is not a zero divisor in $\sO_{\Bl_Z(X),p}$.

We introduce the abbreviation $\widehat Y := \Bl_Z(X)\times_X Y$ and define $E_Y := \pr_1^{-1}(E_X) \subset \widehat Y$. Since $E_X \subset \Bl_Z(X)$ is defined by a locally principal ideal and the morphism $\pr_1:\widehat Y \to \Bl_Z(X)$ is defined by restriction of the projection $\pr_1:\Bl_Z(X)\times Y \to \Bl_Z(X)$ to the closed subspace $\widehat Y \subset \Bl_Z(X)\times Y$ (see Proposition \ref{prop:Fiber_product_closed_subspace_product}), then $E_Y \subset \widehat Y$ is also defined by an ideal $\sI_{E_Y} \subset \sO_{\widehat Y}$ that is locally principal. To see this, observe that $E_X\cap U$ is the zero locus in $U$ of $h \in \sO_{\Bl_Z(X)}\restriction U$ and so
\[
  E_Y\cap \pr_1^{-1}(U) = E_Y\cap(U\times_X Y) = (E_X\cap U)\times_X Y
\]
is the zero locus in $\pr_1^{-1}(U) = U\times_X Y$ of $\pr_1^*h \in \sO_{\widehat Y}\restriction \pr_1^{-1}(U)$. Thus, we have
\[
  \sI_{E_Y} \restriction \pr_1^{-1}(U) = (\pr_1^*h)\sO_{\widehat Y}\restriction \pr_1^{-1}(U).
\]  
For any point $\hat p \in E_Y$ with $\pr_1(\hat p) = p \in E_X$, the morphism $\pr_1:\widehat Y \to \Bl_Z(X)$ defines a homomorphism of (local) rings, $\pr_{1,p}^\sharp:\sO_{\Bl_Z(X),p} \to \sO_{\widehat Y,\hat p}$. Since $h \in \sO_{\Bl_Z(X),p}$ is not a zero divisor and $\sO_{\widehat Y,\hat p}$ is a flat $\sO_{\Bl_Z(X),p}$-module via the homomorphism $\pr_{1,p}^\sharp$ since $\pr_1$ is flat, Lemma \ref{lem:Local_homomorphism_rings} (see Remark \ref{rmk:Local_homomorphism_rings_application}) implies that $\hat h := \pr_{1,p}^\sharp(h)$ is not a zero divisor in $\sO_{\widehat Y,\hat p}$, where we note that $\pr_{1,p}^\sharp(h)$ is the germ at $\hat p$ of $\pr_1^*h \in \sO_{\widehat Y}\restriction \pr_1^{-1}(U)$. Thus, $\sI_{E_Y,\hat p} = \hat h\sO_{\widehat Y,\hat p}$ and so $E_Y$ is a hypersurface in $\widehat Y$ in the sense of Definition \ref{defn:Hypersurface_complex_analytic_space}.

We now observe that
\begin{align*}
  E_Y &= \pr_1^{-1}(E_X) \quad\text{(by definition of $E_Y$)}
  \\
      &= \pr_1^{-1}(\pi_X^{-1}(Z)) \quad\text{(by definition of $E_X$)}
  \\
      &= (\pi_X\circ\pr_1)^{-1}(Z)
  \\
      &= (f\circ\,\pr_2)^{-1}(Z)
        \quad\text{(because \eqref{eq:Fiber_product_blowupX_timesX_Y} is commutative)}
  \\
      &= \pr_2^{-1}(f^{-1}(Z)).
\end{align*}  
Hence, we have a holomorphic map $\pr_2:\Bl_Z(X)\times_X Y \to Y$ such that $E_Y = \pr_2^{-1}(f^{-1}(Z))$ is a hypersurface in $\Bl_Z(X)\times_X Y$ in the sense of Definition \ref{defn:Hypersurface_complex_analytic_space}. Therefore, by the universal property of the blowup (see Item \eqref{item:Fischer_theorem_4-1_c} in Theorem \ref{thm:Fischer_4-1}), there exists a unique morphism
\begin{equation}
  \label{eq:s_morphism_BlZ(X)_timesX_Y_to_BlfinvZ(Y)}
  s: \Bl_Z(X)\times_X Y \to \Bl_{f^{-1}(Z)}(Y)
\end{equation}
such that the following diagram commutes:
\begin{equation}
  \label{eq:Commutative_diagram_defining_s_morphism}
  \begin{tikzcd}
    \Bl_Z(X)\times_X Y \arrow[r, "s", dashed] \arrow[rd, "\pr_2"'] &\Bl_{f^{-1}(Z)}(Y) \arrow[d, "\pi_Y"]
    \\
    &Y
  \end{tikzcd}
\end{equation}
By combining the commutative diagrams \eqref{eq:Fiber_products_and_blowups_analytic_spaces} and \eqref{eq:Commutative_diagram_defining_s_morphism}, we obtain the augmented commutative diagram:
\begin{equation}
  \label{eq:Commutative_diagram_defining_r_and_s_morphisms}
  \begin{tikzcd}
  \Bl_{f^{-1}(Z)}(Y)
  \arrow[drr, bend left, "\Bl_Z(f)"]
  \arrow[ddr, bend right, "\pi_Y"']
  \arrow[dr, "r" description, bend right] & & \\
    & \Bl_Z(X)\times_X Y \arrow[ul, "s" description, bend right] \arrow[r, "\pr_1"'] \arrow[d, "\pr_2"]
      & \Bl_Z(X) \arrow[d, "\pi_X"] \\
& Y \arrow[r, "f"'] &X
\end{tikzcd}
\end{equation}
The universal property of the blowup $\pi_Y:\Bl_{f^{-1}(Z)}(Y) \to Y$ given by Item \eqref{item:Fischer_theorem_4-1_c} in Theorem \ref{thm:Fischer_4-1} implies that $s\circ r =  \id$ on $\Bl_{f^{-1}(Z)}(Y)$ (see Remark \ref{rmk:Uniqueness_blowup_analytic_spaces}). The universal property of the fiber product given by Section \ref{sec:Fiber_products_analytic_spaces} implies that $r\circ s = \id$ on $\Bl_Z(X)\times_X Y$. Thus, $r$ in \eqref{eq:r_morphism_BlfinvZ(Y)_to_BlZ(X)_timesX_Y} is an isomorphism of complex analytic spaces such that $\pi_Y = \pr_2\circ r$ as in \eqref{eq:Fiber_products_and_blowups_analytic_spaces} and this completes the proof of Proposition \ref{prop:Functorial_property_blowup_complex_analytic_space}.
\end{proof}

Proposition \ref{prop:Functorial_property_blowup_complex_analytic_space} immediately yields as corollaries the following proofs of results stated earlier.

\begin{proof}[Proof of Proposition \ref{prop:Functorial_property_blowup_complex_manifold}]
It is enough to observe that Item \eqref{eq:Functorial_property_blowup_complex_analytic_space} in Proposition \ref{prop:Functorial_property_blowup_complex_analytic_space} yields Proposition \ref{prop:Functorial_property_blowup_complex_manifold}.  
\end{proof}

\begin{proof}[Proof of Theorem \ref{thm:Equivariance_property_blowups_complex_manifolds}]
We may apply Proposition \ref{prop:Functorial_property_blowup_complex_manifold} by choosing $Y = X$ and $f = g$, for any $g \in G$. The fact that $\Bl_Z(X)$ is a complex manifold is given by Proposition \ref{prop:Huybrechts_2-5-3}. Proposition \ref{prop:Functorial_property_blowup_complex_manifold} ensures that each element $g \in G \subset \Aut(X)$ (the group of biholomorphic maps of $X$ to itself) lifts to a unique element $\Bl_Z(g) \in \Aut(\Bl_Z(X))$ such that the following diagram commutes:
\[
  \begin{tikzcd}
    \Bl_Z(X) \arrow[r, "\Bl_Z(g)"] \arrow[d, "\pi_X"] &\Bl_Z(X) \arrow[d, "\pi_X"]
    \\
    X \arrow[r, "g"] &X
  \end{tikzcd}
\]
In other words, the blowup $\Bl_Z(X)$ is a complex $G$-manifold and the canonical projection $\pi_X: \Bl_Z(X) \to X$ is a holomorphic, $G$-equivariant map. If the action of $G$ on $X$ is proper, then so is the action of $G$ on $\Bl_Z(X)$ because the blowup map $\pi_X$ is proper and the preceding diagram commutes.
\end{proof}

We have the following analogue in the category of complex analytic spaces of Lemma \ref{lem:Strict_transforms_blowups_or_base_changes} and of Hauser \cite[Proposition 5.1, p. 27]{Hauser_2014} in the category of algebraic varieties.

\begin{prop}[Properties of blowups of complex analytic spaces under base change]
\label{prop:Properties_base_changes_and_strict_transforms_blowup_complex_analytic_space}  
Let $(Y,\sO_Y)$ and $(X,\sO_X)$ be complex analytic spaces and let $(Z,\sO_Z)$ be a closed, complex analytic subspace of $(X,\sO_X)$. If $f : Y \to X$ is a morphism\footnote{This morphism is called the \emph{base change} in algebraic geometry.} of complex analytic spaces and $\widetilde Y$ is the Zariski closure of $\pr_2^{-1}(Y\less f^{-1}(Z))$ in the fiber product $\Bl_Z(X)\times_X Y$, where $\pr_2:\Bl_Z(X)\times_X Y \to Y$ is projection onto the second factor, then the restriction $\pr_2:\widetilde Y\to Y$ of the projection $\pr_2:\Bl_Z(X)\times_X Y \to Y$ to the closed subspace $\widetilde Y \subset \Bl_Z(X)\times_X Y$ is equal to the blowup $\pi_Y:\Bl_{f^{-1}(Z)}(Y)\to Y$.
\end{prop}

\begin{proof}
We adapt the proof due to Hauser of \cite[Proposition 5.1, p. 27]{Hauser_2014}.

First, we claim that $\pr_2^{-1}(f^{-1}(Z))$ is a hypersurface in $\widetilde Y$. To prove this claim, we argue partly as in the proof of Proposition \ref{prop:Functorial_property_blowup_complex_analytic_space}. Let $E_X = \pi_X^{-1}(Z)$ denote the exceptional divisor on $X$ and observe that
\begin{align*}
  \pr_2^{-1}(f^{-1}(Z))
  &= (f\circ \pr_2)^{-1}(Z)
  \\
  &= (\pi_X\circ\pr_1)^{-1}(Z) \quad\text{(because \eqref{eq:Fiber_product_blowupX_timesX_Y} is commutative)}
  \\
  &= \pr_1^{-1}(\pi_X^{-1}(Z)) 
  \\
  &= \pr_1^{-1}(E_X) \quad\text{(by definition of $E_X$)},
\end{align*}
where $\pr_1:\widetilde Y \to \Bl_Z(X)$ is the restriction of the projection $\pr_1:\Bl_Z(X)\times_X Y \to \Bl_Z(X)$ to the closed subspace $\widetilde Y \subset \Bl_Z(X)\times_X Y$. Because $E_X \subset \Bl_Z(X)$ is a hypersurface in the sense of Definition \ref{defn:Hypersurface_complex_analytic_space}, then $E_X$ is defined by a locally principal ideal $\sI_{E_X}$ in $\sO_{\Bl_Z(X)}$, and for any point $p \in E_X \subset \Bl_Z(X)$, there are an open neighborhood $U \subset \Bl_Z(X)$ of $p$ and an element $h \in \sO_{\Bl_Z(X)}\restriction U$ such that $\sI_{E_X}\restriction U = h\sO_{\Bl_Z(X)}\restriction U$, where $h$ is not a zero divisor in $\sO_{\Bl_Z(X),p}$.

Since $E_X \subset \Bl_Z(X)$ is defined by a locally principal ideal and the morphism $\pr_1:\widetilde Y \to \Bl_Z(X)$ is defined by restriction of the projection $\pr_1:\Bl_Z(X)\times_X Y \to \Bl_Z(X)$ to the closed subspace $\widetilde Y \subset \Bl_Z(X)\times Y$, then $\pr_1^{-1}(E_X) \subset \widetilde Y$ is also defined by an ideal $\sI_{\pr_1^{-1}(E_X)} \subset \sO_{\widetilde Y}$ that is locally principal. To see this, observe that $E_X\cap U$ is the zero locus in $U$ of $h \in \sO_{\Bl_Z(X)}\restriction U$ and so $\pr_1^{-1}(E_X)\cap \pr_1^{-1}(U)$ is the zero locus in $\pr_1^{-1}(U)$ of $\pr_1^*h \in \sO_{\widetilde Y}\restriction \pr_1^{-1}(U)$. Thus, we have
\[
  \sI_{\pr_1^{-1}(E_X)} \restriction \pr_1^{-1}(U)
  = (\pr_1^*h)\sO_{\widetilde Y}\restriction \pr_1^{-1}(U).
\]
For any point $\tilde p \in \pr_1^{-1}(E_X)$ with $\pr_1(\tilde p) = p \in E_X$, the morphism $\pr_1:\widetilde Y \to \Bl_Z(X)$ defines a homomorphism of (local) rings, $\pr_{1,p}^\sharp:\sO_{\Bl_Z(X),p} \to \sO_{\widetilde Y,\tilde p}$. According to Remark \ref{rmk:Strict_transform_colon_ideals}, the germ $\tilde h := (\pr_1^*h)_p = \pr_{1,p}^\sharp(h)$ is not a zero divisor in $\sO_{\widetilde Y,\tilde p}$. Thus, $\sI_{\pr_1^{-1}(E_X),\tilde p} = \tilde h\sO_{\widetilde Y,\tilde p}$ and so $\pr_2^{-1}(f^{-1}(Z)) = \pr_1^{-1}(E_X)$ is a hypersurface in $\widetilde Y$ in the sense of Definition \ref{defn:Hypersurface_complex_analytic_space}.

Second, we claim that $\pr_2:\widetilde Y\to Y$ obeys the universal property in Theorem \ref{thm:Fischer_4-1}. To prove this claim, let $\tau:W \to Y$ be a morphism of complex analytic spaces such that $\tau^{-1}(f^{-1}(Z)) = (f\circ\tau)^{-1}(Z)$ is a hypersurface in $W$. As in the proof of Proposition \ref{prop:Functorial_property_blowup_complex_analytic_space},
the universal property of blowups provided by Items \eqref{item:Fischer_theorem_4-1_b} and \eqref{item:Fischer_theorem_4-1_c} in Theorem \ref{thm:Fischer_4-1} yields a unique morphism $\rho:W\to \Bl_Z(X)$ of complex analytic spaces such that $\pi_X\circ\rho = f\circ\tau:W \to X$. Similarly, the universal property of fiber products yields a unique morphism $\sigma:W\to \Bl_Z(X)\times_X Y$ of complex analytic spaces (see Section \ref{sec:Fiber_products_analytic_spaces}) such that the following diagram commutes:
\begin{equation}
  \label{eq:Fiber_products_and_blowups_analytic_spaces_base_change_strict_transform}
  \begin{tikzcd}
  W
  \arrow[drr, bend left, dashed, "\rho"]
  \arrow[ddr, bend right, "\tau"']
  \arrow[dr, dotted, "\sigma" description] & & \\
    & \Bl_Z(X)\times_X Y \arrow[r, "\pr_1"'] \arrow[d, "\pr_2"]
      & \Bl_Z(X) \arrow[d, "\pi_X"] \\
& Y \arrow[r, "f"'] &X
\end{tikzcd}
\end{equation}
It remains to prove that $\sigma(W) \subset \widetilde Y$. Since $\tau^{-1}(f^{-1}(Z))$ is a hypersurface
divisor in $W$ by assumption, its complement $W\less \tau^{-1}(f^{-1}(Z))$ is dense in $W$. But
\[
  W\less \tau^{-1}\left((f^{-1}(Z)\right)
  = \tau^{-1}\left(Y\less f^{-1}(Z)\right)
  = (\pr_2\circ\sigma)^{-1}\left(Y\less f^{-1}(Z)\right)
  = \sigma^{-1}\left(\pr_2^{-1}\left(Y\less f^{-1}(Z)\right)\right),
\]
and thus
\[
  \sigma\left(W\less \tau^{-1}\left(f^{-1}(Z)\right)\right) \subset \pr_2^{-1}\left(Y\less f^{-1}(Z)\right).
\]
Because $\widetilde Y$ is the Zariski closure of $\pr_2^{-1}(Y\less f^{-1}(Z))$ in $\Bl_Z(X)\times_X Y$, we see that $\sigma(W) \subset \widetilde Y$ and this completes the proof of the second claim. Hence, by Remark \ref{rmk:Uniqueness_blowup_analytic_spaces},
there is a unique isomorphism $\psi:\widetilde Y \cong \Bl_{f^{-1}(Z)}(Y)$ of complex analytic spaces such that $\pr_2 = \pi_Y\circ\psi$ and this completes the proof of Proposition \ref{prop:Properties_base_changes_and_strict_transforms_blowup_complex_analytic_space}.
\end{proof}  

In his definition of strict transform \cite[Chapter 0, Section 2, p. 130]{Hironaka_1964-I-II}, Hironaka allows for analytic spaces over $\KK=\RR$ or $\CC$, but his approach is rather different from that of Definition \ref{defn:Strict_transform}. The most natural analogue of Definition \ref{defn:Strict_transform} is given by the

\begin{defn}[Strict transform in the blowup of an analytic space along a closed analytic subspace]
\label{defn:Strict_transform_analytic_space}
(See Bierstone and Milman \cite[Section 2.5, p. 807]{Bierstone_Milman_1989}.)
Let $\KK=\RR$ or $\CC$ and $(X,\sO_X)$ be an analytic space over $\KK$ as in Definition \ref{defn:Analytic_space} and $(Z,\sO_Z)$ be a closed analytic subspace of $(X,\sO_X)$ as in Definition \ref{defn:Analytic_subspace} and $\pi:\widetilde X\to X$ be the blowup of $X$ along $Z$ as in Theorem \ref{thm:Fischer_4-1}, denoting $\widetilde X = \Bl_Z(X)$ for convenience. If $(Y,\sO_Y)$ is a closed analytic subspace of $(X,\sO_X)$, then the \emph{strict} (or \emph{proper}) \emph{transform} $(\widetilde Y,\sO_{\widetilde Y})$ of $(Y,\sO_Y)$ under the blowup morphism $\pi:\widetilde X \to X$ is the smallest closed analytic subspace of $(\widetilde X,\sO_{\widetilde X})$ such that $\pi$ induces an isomorphism $\widetilde Y\less E \cong Y\less Z$ of analytic spaces, where $E = \pi^{-1}(Z)$ is the exceptional divisor in $\widetilde X$.
\end{defn}

\begin{rmk}[Explicit construction of the strict transform in the blowup of an analytic space along a closed analytic subspace]
\label{rmk:Strict_transform_analytic_space_construction}
When $X$ is an analytic manifold and $Z$ is an embedded analytic submanifold of $X$ in Definition \ref{defn:Strict_transform_analytic_space}, Bierstone and Milman provide an explicit construction of the strict transform $\widetilde Y$ in \cite[Section 2.5, pp. 806--807]{Bierstone_Milman_1989} and state, but do not prove, that under these hypotheses, their construction gives the smallest closed analytic subspace of $\widetilde X$ such that $\pi$ induces an isomorphism $\widetilde Y\less E \cong Y\less Z$.
\end{rmk}  

We have the following analogue in the category of analytic spaces of Hauser \cite[Corollary 5.2 (a), p. 28]{Hauser_2014} in the category of algebraic varieties and corollary of Proposition
\ref{prop:Properties_base_changes_and_strict_transforms_blowup_complex_analytic_space}.

\begin{cor}[Strict transform of a closed, complex analytic subspace under blowup of a complex analytic space along a closed, complex analytic subspace]
\label{cor:Strict_transform_closed_subspace_under_blowup_analytic_space_along_subspace}
Let $(X,\sO_X)$ be a complex analytic space, $Z \subsetneq X$ be a closed, complex analytic subspace, and $\pi_X:\Bl_Z(X)\to X$ be the blowup of $X$ along $Z$ as in Theorem \ref{thm:Fischer_4-1}. If $Y \subseteq X$ is a closed, complex analytic subspace and $Y\cap Z = \iota^{-1}(Z)$, where $\iota: Y \hookrightarrow X$ denotes inclusion, then the strict transform of $Y$ in $\Bl_Z(X)$ as in Definition \ref{defn:Strict_transform_analytic_space} is given by
\begin{equation}
  \label{eq:Strict_transform_analytic_space_Y_is_blowup_along_YcapZ}
  \widetilde Y = \Bl_{Y\cap Z}(Y),
\end{equation}
where $\pi_Y:\Bl_{Y\cap Z}(Y)\to Y$ is the blowup of $Y$ along $Y\cap Z$ as in Theorem \ref{thm:Fischer_4-1}. Moreover, there is an embedding
\begin{equation}
  \label{eq:Blowup_embedding}
  \Bl_Z(\iota):\Bl_{Y\cap Z}(Y) \hookrightarrow \Bl_Z(X)
\end{equation}
such that the following diagram commutes:
\begin{equation}
  \label{eq:Blowup_complex_analytic_subspace_and_embedding_commutative_diagram}
  \begin{tikzcd}
    \Bl_{Y\cap Z}(Y) \arrow[r, "\Bl_Z(\iota) "] \arrow[d, "\pi_Y"]
    &\Bl_Z(X)\arrow[d, "\pi_X"]
    \\
    Y \arrow[r, "\iota"] & X
  \end{tikzcd}
\end{equation}
\end{cor}

\begin{proof}
By applying Proposition \ref{prop:Properties_base_changes_and_strict_transforms_blowup_complex_analytic_space} with $f=\iota$, we obtain a morphism \eqref{eq:Blowup_embedding} and a commutative diagram \eqref{eq:Blowup_complex_analytic_subspace_and_embedding_commutative_diagram}. Because $f=\iota$, the fiber product $\Bl_Z(X)\times_X Y$ is given as a set by (see Section \ref{sec:Fiber_products_analytic_spaces})
\[
  \Bl_Z(X)\times_X Y = \left\{(\tilde x, y) \in \Bl_Z(X)\times Y: \pi_X(\tilde x) = y\right\}
\]
and so $f^{-1}(Z) = \iota^{-1}(Z) = Y\cap Z$ and $Y \less f^{-1}(Z) =  Y\less Z$. In particular, the Zariski closure of $\pr_2^{-1}(Y\less Z)$ in $\Bl_Z(X)\times_X Y$ is equal to the Zariski closure of $\pi_X^{-1}(Y\less Z)$ in $\Bl_Z(X)$, namely the strict transform, $\widetilde Y$. Thus, Proposition \ref{prop:Properties_base_changes_and_strict_transforms_blowup_complex_analytic_space} yields the equality \eqref{eq:Strict_transform_analytic_space_Y_is_blowup_along_YcapZ} of complex analytic spaces. The morphism $\Bl_Z(\iota)$ is an embedding since $\Bl_{Y\cap Z}(Y) = \widetilde Y$, a closed, complex analytic subspace of $\Bl_Z(X)$.
\end{proof}

We obtain the following analogue of Corollaries \ref{cor:Smoothness_strict_transform_subvariety} and \ref{cor:Smoothness_strict_transform_subscheme} as an almost immediate consequence of Proposition \ref{prop:Huybrechts_2-5-3} and Corollary \ref{cor:Strict_transform_closed_subspace_under_blowup_analytic_space_along_subspace}.

\begin{cor}[Smoothness of the strict transform of an embedded complex submanifold]
\label{cor:Smoothness_strict_transform_submanifold}
Let $X$ be a complex manifold and $Z \subset X$ be an embedded complex submanifold. If $Y$ is an embedded complex submanifold of $X$ and the intersection $Z\cap Y$ is an embedded complex submanifold of $Y$, then the strict transform $\widetilde Y$ as in Definition \ref{defn:Strict_transform_analytic_space} is an embedded complex submanifold of the blowup $\widetilde X$ of $X$ along $Z$.
\end{cor}

\begin{proof}
According to Corollary \ref{cor:Strict_transform_closed_subspace_under_blowup_analytic_space_along_subspace}, the strict transform $\widetilde Y$ is equal to the blowup of $Y$ along $Z\cap Y$. Since $Y$ is a complex manifold and the blowup center $Z\cap Y$ is an embedded complex submanifold of $Y$, then the strict transform $\widetilde Y$ is complex manifold by Proposition \ref{prop:Huybrechts_2-5-3}.
\end{proof}

\section[Exceptional divisor and strict transform for blowup of analytic model space]{Exceptional divisor and strict transform for the blowup of an analytic model space along a linear subspace}
\label{sec:Exceptional_divisor_and_strict_transform_complex_analytic_model_space}
The proof of the forthcoming Theorem \ref{thm:Strict_transform_complex_analytic_model_space} requires the Hilbert Nullstellensatz in the category of complex algebraic varieties and the R\"uckert Nullstellensatz in the category of complex analytic spaces. We recall their statements here.

\begin{thm}[Hilbert Nullstellensatz]
\label{thm:Hilbert_Nullstellensatz}
(See Atiyah and Macdonald \cite[Chapter 7, Exercise 14, p. 85]{Atiyah_Macdonald_introduction_commutative_algebra}, Eisenbud \cite[Section 1.6, Theorem 1.6, p. 33]{Eisenbud_commutative_algebra}, Hartshorne \cite[Chapter I, Section 1, Theorem 1.3A, p. 4]{Hartshorne_algebraic_geometry}, and Zariski and Samuel \cite[Chapter VII, Section 3, Theorem 14, p. 164]{Zariski_Samuel_communtative_algebra_II}.)
Let $k$ be a field with algebraically closed extension $K$ and $n\geq 1$ be an integer, $\fa \subset k[x_1,\ldots,x_n]$ be an ideal, $\VV(\fa) \subset K^n$ denote the affine algebraic variety defined\footnote{As in Zariski and Samuel \cite[Chapter VII, Section 3, p. 160]{Zariski_Samuel_communtative_algebra_II}.} by $\fa$, and $I(V) \subset k[x_1,\ldots,x_n]$ denote the ideal defined\footnote{As in Zariski and Samuel \cite[Chapter VII, Section 3, p. 160]{Zariski_Samuel_communtative_algebra_II}.} by an affine algebraic variety $V \subset K^n$. Then
\[
  I(\VV(\fa)) = \sqrt{\fa},
\]
where $\sqrt{\fa}$ is the radical ideal\footnote{As in Zariski and Samuel \cite[Chapter III, Section 7, Definition 2, p. 147]{Zariski_Samuel_communtative_algebra_II}.} of $\fa$, or equivalently: If $f\in k[x_1,\ldots,x_n]$ and $\fa$ is generated by $f_1,\ldots,f_r\in k[x_1,\ldots,x_n]$
and $f$ vanishes at every common zero in $K^n$ of $f_1,\ldots,f_r$, then there exist a positive integer $\rho$ and $a_1\ldots,a_r \in k[x_1,\ldots,x_n]$ such that
\[
  f^\rho = a_1f_1 + \cdots + a_rf_r.
\]
\end{thm}

\begin{rmk}[Interpretation and proof of the Hilbert Nullstellensatz]
\label{rmk:Interpretation_proof_Hilbert_Nullstellensatz}
Tao \cite{Tao_2007_hilbert_nullstellensatz} provides a very useful discussion and proof of Theorem \ref{thm:Hilbert_Nullstellensatz}.
\end{rmk}  

The analogue of Theorem \ref{thm:Hilbert_Nullstellensatz} in the category of analytic spaces is given by

\begin{thm}[Analytic or R\"uckert Nullstellensatz]
\label{thm:Ruckert_Nullstellensatz}
(See Abhyankar \cite[Chapter V, Theorem 30.12, p. 238]{Abhyankar_local_analytic_geometry} for the case of a complete nondiscrete valued field $K$ that is algebraically closed and for $\CC$, see Grauert and Remmert \cite[Section 3.2.2, Theorem, p. 67, and Section 4.1.5, Theorem, p. 82]{Grauert_Remmert_coherent_analytic_sheaves}, Guaraldo, Macr\`\i, and Tancredi \cite[Section 2.2, Theorem 2.10 p. 12]{Guaraldo_Macri_Tancredi_topics_real_analytic_spaces}, or Huybrechts \cite[Section 1.1, Proposition 1.1.29, p. 19]{Huybrechts_2005}.)
Let $n\geq 1$ be an integer, $\fa \subset \sO_{\CC^n,0}$ be an ideal, $\VV(\fa) \subset \CC^n$ denote the analytic germ\footnote{As in Huybrechts \cite[Section 1.1, Definitions 1.1.21 and 1.1.22, p. 18]{Huybrechts_2005}.} defined by $\fa$, and $I(V) \subset \sO_{\CC^n,0}$ denote the ideal defined by an analytic germ\footnote{As in Huybrechts \cite[Section 1.1, Definition 1.1.24, p. 18]{Huybrechts_2005}.} $V \subset \CC^n$ at the origin. Then
\[
  I(\VV(\fa)) = \sqrt{\fa},
\]
where $\sqrt{\fa}$ is the radical ideal\footnote{As in Huybrechts \cite[Section 1.1, p. 19]{Huybrechts_2005}.} of $\fa$.
\end{thm}

While Theorem \ref{thm:Ruckert_Nullstellensatz} and Theorem \ref{thm:Strict_transform_complex_analytic_model_space} below are stated for the category of analytic spaces over $\CC$, one might ask whether they could be reformulated along the lines of Theorem \ref{thm:Hilbert_Nullstellensatz} for analytic spaces over any complete nondiscrete valued field $k$ with an algebraically closed extension $K$.

\begin{thm}[Exceptional divisor and strict transform for the blowup of a complex analytic model space along a linear subspace]
\label{thm:Strict_transform_complex_analytic_model_space}
Let $(Y,\sO_Y)$ be a complex analytic model space as in Definition \ref{defn:Analytic_model_space}, so there is a domain $D$ around the origin in $X = \CC^n$, a finitely-generated ideal $\sI \subset \sO_D$, a topological space $Y = \cosupp\sI$, and a structure sheaf $\sO_Y = (\sO_D/\sI)\restriction Y$. Let
\begin{equation}
  \label{eq:Complex_coordinate_subspace_blowup_center}
  Z := \{(z_1,\ldots,z_n) \in \CC^n: z_j=0, \text{ for } j=m+1,\ldots,n\} \cong \CC^m,
\end{equation}
denote the coordinate subspace in \eqref{eq:Coordinate_subspace_blowup_center} with $\KK=\CC$,
where $0\leq m < n$, and let $\pi_X:\Bl_Z(X)\to X$ denote the blowup of $X$ along the center $Z$ as in Section \ref{sec:Blowups_analytic_manifolds_along_embedded_analytic_submanifolds}. Let 
\begin{equation}
  \label{eq:Complex_coordinate_subspace_blowup_center_orthogonal_complement}
  Z^\perp := \left\{(z_1,\ldots,z_n) \in \CC^n: z_j=0, \text{ for } j=1,\ldots,m\right\} \cong \CC^{n-m}
\end{equation}
denote the coordinate subspace given by the orthogonal complement in $\CC^n$ of $Z$ in \eqref{eq:Complex_coordinate_subspace_blowup_center}. Then in each coordinate domain $U_j \subset \Bl_Z(X)$ as in \eqref{eq:Blowup_linear_subspace_coordinate_patch}, the topological supports of the exceptional divisor $E \subset \Bl_Z(X)$, the strict transform\footnote{In the sense of Definition \ref{defn:Strict_transform_analytic_space}.} $\Bl_{Z\cap Y}(Y)$ of the blowup, and the intersection $E\cap\Bl_{Z\cap Y}(Y)$ are defined by the forthcoming analytic varieties \eqref{eq:Exceptional_divisor_support_coordinate_domain}, \eqref{eq:Strict_transform_support_coordinate_domain}, and \eqref{eq:Intersection_exceptional_divisor_and_support_coordinate_domain}, respectively, that are explicitly determined by local generators of the ideal $\sI$.
\end{thm}

\begin{proof}
By replacing $D$ with a possibly smaller open neighborhood $U$ of the origin and relabeling, we may assume without loss of generality that $\sI = (f_1,\ldots,f_r)$, where $f_k \in \sO_X(D)$ for $k = 1,\ldots,r$. 
Because $Z \subset Y$ by hypothesis and $Z = \VV(z_{m+1},\ldots,z_n) \subset \CC^n$ by \eqref{eq:Complex_coordinate_subspace_blowup_center} and $Y = \VV(f_1,\ldots,f_r)$ by hypothesis, Theorem  \ref{thm:Ruckert_Nullstellensatz} yields positive integers $\rho_k$ and holomorphic functions $a_{k,m+1},\ldots,a_{k,n} \in \sO_X(D)$ such that $f_k^{\rho_k} \in I(Z) = z_{m+1}\sO_X+\cdots+z_n\sO_X$ and thus
\begin{equation}
  \label{eq:Ruckert_Nullstellensatz_Y_Z}
  f_k^{\rho_k} = z_{m+1}a_{k,m+1} + \cdots + z_na_{k,n}, \quad\text{for } k = 1,\ldots,r.
\end{equation}
We may use the local coordinate form \eqref{eq:Blowup_linear_subspace_projection_Kn_local_chart} of the restriction $\pi_X:\Bl_Z(X)|_{U_j} \to X$ of the blowup morphism $\pi_X$ in \eqref{eq:Blowup_linear_subspace_projection_Kn} to the domain $U_j \subset \Bl_Z(X)$ of the coordinate chart $\varphi_j:\Bl_Z(X)|_{U_j} \cong \varphi_j(U_j) \subset X$, for each $j=m+1,\ldots,n$, to describe the total transform $\pi_X^*\sI$ of the ideal $\sI$,
\[
  \pi_X^*\sI = (\pi_X^*f_1,\ldots, \pi_X^*f_r) \subset \sO_{\widetilde X},
\]
where $\widetilde X = \Bl_Z(X)$ and $\pi_X^*f_k \in \sO_{\widetilde X}(\pi_X^{-1}(D))$ for $k = 1,\ldots,r$. With respect to the chart $(U_j,\varphi_j)$ on $\Bl_Z(X)$ for $j=m+1,\ldots,n$ and $k = 1,\ldots,r$, we have
\begin{align*}
  (\varphi_j)^{-1}(\pi_X^*f_k)(w_1,\ldots,w_n)
  &=
  \left(\pi_X\circ\varphi_j^{-1}\right)^*f_k(w_1,\ldots,w_n)
  \\
  &=
  f_k\left(\pi_X\left(\varphi_j^{-1}(w_1)\right),\ldots,\pi_X\left(\varphi_j^{-1}(w_n)\right)\right)
  \\
  &= f_k\left(w_1,\ldots,w_m,w_j(w_{m+1},\ldots,w_{j-1},1,w_{j+1},\ldots,w_n)\right),
  \\
  &\qquad\text{for all } (w_1,\ldots,w_n) \in \varphi_j(U_j) \subset \CC^n.
\end{align*}
By applying the identity \eqref{eq:Ruckert_Nullstellensatz_Y_Z} for $j=m+1,\ldots,n$ and $k = 1,\ldots,r$, we see that
\begin{multline*}
  f_k^{\rho_k}\left(w_1,\ldots,w_m,w_j(w_{m+1},\ldots,w_{j-1},1,w_{j+1},\ldots,w_n)\right)
  \\
  =
  w_j\sum_{i=m+1, i\neq j}^n w_ia_{k,i}\left(w_1,\ldots,w_m,w_{m+1},\ldots,w_{j-1},1,w_{j+1},\ldots,w_n\right)
  \\
  + w_ja_{k,j}\left(w_1,\ldots,w_m,w_{m+1},\ldots,w_{j-1},1,w_{j+1},\ldots,w_n\right).
\end{multline*}
Since $f_k\in\CC\{z_1,\ldots,z_n\}$ defines a generator of $\sI \subset \sO_X$ for each $k = 1,\ldots,r$, it is not identically zero. Hence, for each $j=m+1,\ldots,n$ and $k = 1,\ldots,r$, there exist a unique maximum positive integer $k_j$ and a non-zero holomorphic function $g_{k,j} \in \CC\{w_1\ldots,w_{j-1},w_{j+1},\ldots,w_n\}$ such that
\begin{multline}
  \label{eq:Blowup_generator}
  f_k^{\rho_k}\left(w_1,\ldots,w_m,w_j(w_{m+1},\ldots,w_{j-1},1_j,w_{j+1},\ldots,w_n)\right)
  \\
  =
  w_j^{k_j}g_{k,j}(w_1,\ldots,w_{j-1},1,w_{j+1},\ldots,w_n).
\end{multline}
The topological support of the restriction $\widetilde Y\cap U_j$ of the strict transform $\widetilde Y = \Bl_{Z\cap Y}(Y) \subset \Bl_Z(X)$ (see Corollary \ref{cor:Strict_transform_closed_subspace_under_blowup_analytic_space_along_subspace}) to the coordinate neighborhood $U_j \subset \Bl_Z(X)$ is thus defined for each $j=m+1,\ldots,n$ by the variety $\VV(g_{1,j}\circ\varphi_j,\ldots,g_{r,j}\circ\varphi_j) \subset U_j$ corresponding to the ideal
\[
  (g_{1,j},\ldots,g_{r,j}) \subset \CC\{w_1\ldots,w_{j-1},w_{j+1},\ldots,w_n\}
  \subset \CC\{w_1\ldots,w_n\}. 
\]
Indeed, to see this we observe that if $p\in Y$, then $f_k(p)=0$ for $k = 1,\ldots,r$ and hence $f_k^{\rho_k}(p) = 0$ for $k = 1,\ldots,r$. If $\tilde p \in \pi_Y^{-1}(p) \in \Bl_{Z\cap Y}(Y)$, then we must have $\tilde p \in \Bl_{Z\cap Y}(Y)\cap U_j$ for some $j=j(\tilde p)$, because the blowup morphism $\pi_Y = \pi_X\restriction Y:\Bl_{Z\cap Y}(Y)\to Y$ is surjective, the neighborhoods $U_j$ cover $\Bl_Z(X)$, and the neighborhoods $\Bl_{Z\cap Y}(Y)\cap U_j$ cover $\Bl_{Z\cap Y}(Y)$. Therefore, writing $\varphi_j(\tilde p) = (w_1(\tilde p),\ldots,\tilde w_n(p)) \in \CC^n$, the identities \eqref{eq:Blowup_generator} yield
\[
  w_j^{k_j}(\tilde p)
  g_{k,j}(w_1(\tilde p),\ldots,w_{j-1}(\tilde p),1,w_{j+1}(\tilde p),\ldots,w_n(\tilde p)) = 0,
\]
and so, for $k = 1,\ldots,r$,
\[
  w_j^{k_j}(\tilde p) = 0 \quad\text{or}\quad
  g_{k,j}(w_1(\tilde p),\ldots,w_{j-1}(\tilde p),1,w_{j+1}(\tilde p),\ldots,w_n(\tilde p)) = 0.
\]
The equation
\begin{equation}
  \label{eq:Exceptional_divisor_support_equation_in_local_coordinates}
  w_j = 0
\end{equation}
defines points in an open neighborhood
\begin{equation}
  \label{eq:Exceptional_divisor_support_coordinate_domain}
  E\cap U_j = \VV(w_j\circ\varphi_j)
\end{equation}
in the topological support of the exceptional divisor $E \subset \Bl_Z(X)$ in $U_j$ for the blowup morphism $\pi_X:\Bl_Z(X)\to X$, while the equations
\begin{equation}
  \label{eq:Strict_transform_support_equations_in_local_coordinates}
  g_{k,j}(w_1,\ldots,w_{j-1},1,w_{j+1},\ldots,w_n) = 0, \quad\text{for } k = 1,\ldots,r,
\end{equation} 
define points in an open neighborhood
\begin{equation}
  \label{eq:Strict_transform_support_coordinate_domain}
  \Bl_{Z\cap Y}(Y)\cap U_j = \VV(g_{1,j}\circ\varphi_j,\ldots,g_{r,j}\circ\varphi_j)
\end{equation}
of the topological support of the strict transform $\Bl_{Z\cap Y}(Y)$ in $U_j \subset \Bl_Z(X)$, for $j=m+1,\ldots,n$. The variety
\begin{equation}
  \label{eq:Intersection_exceptional_divisor_and_support_coordinate_domain}
  E\cap\Bl_{Z\cap Y}(Y)\cap U_j = \VV(w_j\circ\varphi_j,g_{1,j}\circ\varphi_j,\ldots,g_{r,j}\circ\varphi_j)
\end{equation}
is the topological support of the intersection of the exceptional divisor $E$ and the strict transform $\Bl_{Z\cap Y}(Y)$ in $U_j \subset \Bl_Z(X)$, for $j=m+1,\ldots,n$. This completes the proof of Theorem \ref{thm:Strict_transform_complex_analytic_model_space}.
\end{proof}

\chapter{Monomialization of ideal sheaves and resolution of singularities}
\label{chap:Resolution_singularities}
In this chapter, we describe some of the main results due to Hironaka on resolution of singularities. We begin in Section \ref{sec:Normal_crossing_divisors} by reviewing the key concepts of divisors with normal and simple normal crossings, continuing in Section \ref{sec:Ideals_simple_normal_crossing_divisors} with an analogous discussion for ideals with normal and simple normal crossings. Section \ref{sec:Main_theorems_monomialization_resolution_singularities_schemes} contains the statements of Hironaka's main theorems on monomialization and resolution of singularities for schemes while Section \ref{sec:Principalization_and_resolution_singularities_analytic_spaces} contains the statements of his corresponding results for analytic spaces. In Section \ref{sec:Equivariant_resolution_of_singularities}, we discuss equivariant resolution of singularities.

\section{Divisors with normal and simple normal crossings}
\label{sec:Normal_crossing_divisors}
We begin in Section \ref{subsec:Divisors_normal_simple_normal_crossings_algebraic_varieties} by reviewing the concepts of divisors with normal and simple normal crossings on algebraic varieties. We continue in Section \ref{subsec:Divisors_normal_crossings_schemes} by reviewing the more general concept of divisors with normal crossings on schemes. We conclude in Section \ref{subsec:Divisors_normal_crossings_analytic_spaces} by briefly commenting on the definition of divisors with normal crossings on analytic spaces.

\subsection{Divisors with normal and simple normal crossings on algebraic varieties}
\label{subsec:Divisors_normal_simple_normal_crossings_algebraic_varieties}
We begin by recalling the

\begin{defn}[Simple normal crossing divisor on an algebraic variety]
\label{defn:Subvariety_having simple_normal_crossings_with_divisor}
(See Koll{\'a}r \cite[Definition 3.24, p. 137]{Kollar_lectures_resolution_singularities}.)
Let $X$ be a smooth algebraic variety of dimension $d\geq 1$. One says that $E = \sum_i E_i$ is a \emph{simple normal crossing divisor} on $X$ if each $E_i$ is smooth and for each point $p \in X$ one can choose local coordinates $x_1,\ldots,x_d$ in the maximal ideal $\fm_p$ of the local ring $\sO_p$ such that for each $i$ the following hold:
\begin{enumerate}
\item Either $p \notin E_i$ or $E_i\cap U = \{q \in U: x_{j_i}(q) = 0\}$ in an open neighborhood $U \subset X$ of $p$ for some $j_i$, and

\item $j_i \neq j_{i'}$ if $i \neq i'$.
\end{enumerate}
An algebraic subvariety $Z \subset X$ has \emph{simple normal crossings} with $E$ if one can choose $x_1,\ldots,x_d$ as above such that in addition
\begin{enumerate}
\setcounter{enumi}{2}
\item $Z = \{q \in U: x_{j_1}(q) = \cdots = x_{j_s}(q) = 0\}$ for some $j_1,\ldots,j_s$ and open neighborhood $U$ of $p$.
\end{enumerate}
\end{defn}

In particular, $Z$ is smooth in Definition \ref{defn:Subvariety_having simple_normal_crossings_with_divisor} and some of the $E_i$ are allowed to contain $Z$. Koll{\'a}r also gives the following, more elementary definition that serves, in part, to help compare the concepts of simple normal crossing divisor (as used by \cite{Kollar_1999, Wlodarczyk_2008}) and normal crossing divisor (as used by \cite{Bierstone_Milman_1997}), in the context of resolution of singularities. 

\begin{defn}[Simple normal crossing divisor on an algebraic variety]
\label{defn:Simple_normal_crossing_divisor}
(See Koll{\'a}r \cite[Definition 1.44, p. 30]{Kollar_lectures_resolution_singularities}.) Let $X$ be a smooth algebraic variety of dimension $d\geq 1$ and $E \subset X$ a divisor. One calls $E$ a \emph{simple normal crossing divisor} if every irreducible component of $E$ is smooth and all intersections are transverse. That is, for every point $p \in E$ we can choose local coordinates $x_1,\ldots,x_d$ on an open neighborhood $U \subset X$ of $p$ and $m_i \in \ZZ\cap [0,\infty)$ for $i=1,\ldots,d$ such that $U\cap E = \{q \in U: \prod_{i=1}^d x^{m_i}(q) = 0 \}$.
\end{defn}

\begin{rmk}[Normal crossing divisor on an algebraic variety]
\label{rmk:Normal_crossing_divisor}
(See Koll{\'a}r \cite[Remark 1.45., p. 30]{Kollar_lectures_resolution_singularities}.)
Continuing the notation of Definition \ref{defn:Simple_normal_crossing_divisor},
one calls $E$ a \emph{normal crossing divisor} if for every $p \in E$ there are local analytic or formal coordinates, $x_1,\ldots,x_d$, and natural numbers $m_1,\ldots,m_d$ such that $U\cap E = \{q \in U: \prod_{i=1}^d x^{m_i}(q) = 0 \}$.

For example, the nodal curve $y^2 = x^3 + x^2$ is a normal crossing divisor in $\CC^2$, but not a simple normal crossing divisor because we can write
\[
  y^2 = x^3 + x^2 = \left(y - x\sqrt{1+x}\right)\left(y + x\sqrt{1+x}\right)
\]
as a power series, but $y^2 - x^3 - x^2$ is irreducible as a polynomial.
\end{rmk}

Definitions of (simple) normal crossing divisors on algebraic varieties are also provided by Cutkosky \cite[Exercise 3.13 (2), p. 29]{Cutkosky_resolution_singularities}, Hartshorne \cite[Chapter V, Section 3, Remark 3.8.1, pp. 390--391]{Hartshorne_algebraic_geometry}, and Lazarsfeld \cite[Definition 4.1.1, p. 238]{Lazarsfeld_positivity_algebraic_geometry_v1}. See Hauser \cite[Definition 3.15, p. 14]{Hauser_2014} for definitions of normal crossings points and simple normal crossings points in affine algebraic varieties.

\subsection{Divisors with normal crossings on schemes}
\label{subsec:Divisors_normal_crossings_schemes}
One has analogues of the preceding definitions in the category of schemes.

\begin{defn}[Strict normal crossings divisor on a scheme]
\label{defn:Strict_normal_crossing_divisor_scheme}  
(See the Stacks Project \cite[\href{https://stacks.math.columbia.edu/tag/0BI9}{Definition 0BI9}]{stacks-project}.)  
Let $X$ be a locally Noetherian scheme. A \emph{strict normal crossings divisor} on $X$ is an effective Cartier divisor $D\subset X$ such that for every point $p\in D$, the local ring $\sO_{X,p}$ is regular and there exists a regular system of parameters $x_1,\ldots,x_d \in \fm_p$ and an integer $1\leq r\leq d$ such that $D$ is cut out by $x_1,\ldots,x_r$ in $\sO_{X,p}$.
\end{defn}

The following lemma provides an equivalent characterization of a strict normal crossings divisor on a scheme.

\begin{lem}[Equivalent characterization of a strict normal crossings divisor on a scheme]
\label{lem:Equivalent_characterization_strict_normal_crossings_divisor_scheme}
(See the Stacks Project \cite[\href{https://stacks.math.columbia.edu/tag/0BIA}{Lemma 0BIA}]{stacks-project}.)
Let $X$ be a locally Noetherian scheme, $D\subset X$ be an effective Cartier divisor, and $D_i \subset D$, for $i\in I$, be its irreducible components viewed as reduced closed subschemes of $X$. Then the following are equivalent:
\begin{enumerate}
\item $D$ is a strict normal crossings divisor.
\item $D$ is reduced, each $D_i$ is an effective Cartier divisor, and for $J\subset I$ finite the scheme theoretic intersection $D_J = \cap_{j\in J}D_j$ is a regular scheme each of whose irreducible components has codimension $|J|$ in X.
\end{enumerate}
\end{lem}

\begin{defn}[Normal crossings divisor on a scheme]
\label{defn:Normal_crossing_divisor_scheme}
(See the Stacks Project \cite[\href{https://stacks.math.columbia.edu/tag/0BSF}{Definition 0BSF}]{stacks-project}; see also Hironaka \cite[Chapter 0, Section 5, Definition 2, p. 141]{Hironaka_1964-I-II}.)  
Let $X$ be a locally Noetherian scheme. A \emph{normal crossings divisor} on $X$ is an effective Cartier divisor $D\subset X$ such that for every point $p \in D$ there exists an \'etale morphism $U \to X$ with $p$ in the image and $D\times_XU$ is a strict normal crossings divisor on $U$.
\end{defn}

For example, $D=\VV(x^2+y^2)$ is a normal crossings divisor but not a strict one on $\Spec(\RR[x,y])$ because after pulling back to the \'etale cover $\Spec(\CC[x,y])$, one obtains $(x-iy)(x+iy)=0$ (see the Stacks Project \cite[\href{https://stacks.math.columbia.edu/tag/0CBN}{Section 0CBN}]{stacks-project}). See the Stacks Project \cite[\href{https://stacks.math.columbia.edu/tag/02GI}{Definition 02GI}]{stacks-project} for a definition of an \emph{\'etale} morphism of schemes and \cite[\href{https://stacks.math.columbia.edu/tag/0215}{Definition 0215}]{stacks-project} for the definition of an \emph{\'etale cover}. The following lemma provides an equivalent characterization of a normal crossings divisor on a scheme.

\begin{lem}[Equivalent characterization of a normal crossings divisor on a scheme]
\label{lem:Equivalent_characterization_normal_crossings_divisor_scheme}
(See the Stacks Project \cite[\href{https://stacks.math.columbia.edu/tag/0CBR}{Lemma 0CBR}]{stacks-project}.)
Let $X$ be a locally Noetherian scheme and $D\subset X$ be a closed subscheme. Then the following are equivalent:
\begin{enumerate}
\item $D$ is a normal crossings divisor on $X$,
\item $D$ is reduced, the normalization $\nu:D^\nu \to D$ is unramified, and for any integer $n \geq 1$ the scheme
\[
   Z_n = D^\nu\times_D\cdots \times_D D^\nu \less\left\{(p_1,\ldots,p_n):p_i\neq p_j\text{ for some } i\neq j\right\}
\]
is regular and the morphism $Z_n\to X$ is a local complete intersection morphism whose conormal sheaf is locally free of rank $n$.
\end{enumerate}
\end{lem}

\subsection{Divisors with normal crossings on analytic spaces}
\label{subsec:Divisors_normal_crossings_analytic_spaces}
Definitions \ref{defn:Subvariety_having simple_normal_crossings_with_divisor} and \ref{defn:Simple_normal_crossing_divisor} for algebraic varieties extend to the categories of analytic spaces, where $\sO_p$ is then the local ring of analytic functions; see, for example, Bierstone and Milman \cite[Section 2.1, p. 804]{Bierstone_Milman_1989} and Hironaka \cite[Chapter 0, Section 5, Definition 2, p. 141]{Hironaka_1964-I-II} or Koll{\'a}r \cite[Section 3.44]{Kollar_lectures_resolution_singularities}. In the category of analytic spaces, Remark \ref{rmk:Normal_crossing_divisor} implies that the concepts of simple normal crossing divisor and normal crossing divisor coincide.

\section{Ideals with normal and simple normal crossings}
\label{sec:Ideals_simple_normal_crossing_divisors}
Before stating the main versions of resolution of singularities, we recall some definitions from Cutkosky \cite[pp. 40--41]{Cutkosky_resolution_singularities} and
Koll{\'a}r \cite[Note on Terminology 3.16, p. 143]{Kollar_lectures_resolution_singularities}.

\begin{defn}[Principalization or monomialization of an ideal sheaf]
\label{defn:Principalization_monomialization_ideal_sheaf}
(See Cutkosky \cite[Section 4.2, p. 40]{Cutkosky_resolution_singularities} or Koll{\'a}r \cite[Notation 3.15, p. 134 and p. 135]{Kollar_lectures_resolution_singularities}.)    
Let $X$ be a smooth scheme and $\sI \subset \sO_X$ be an ideal sheaf. A \emph{principalization} (or \emph{monomialization}) of $\sI$ is a proper birational morphism $\pi : \widetilde{X} \to X$ such that $\widetilde{X}$ is smooth and the inverse image ideal sheaf
\[
\pi^*\sI \subset \sO_{\widetilde{X}}
\]
is a locally principal ideal sheaf.
\end{defn}

\begin{rmk}[Inverse images of sheaves]
\label{rmk:Inverse_images_sheaves}
Let $f:Y\to X$ be a continuous map of topological spaces and $\sG$ be a sheaf over $X$. We refer to G\"ortz and Wedhorn \cite[Chapter 2, pp. 54--55]{Gortz_Wedhorn_algebraic_geometry_v1} for the definition of the inverse image sheaf $f^{-1}\sG$ over $X$. If $f$ is the inclusion of $Y$ as a subspace of $X$, they also write $\sG\restriction Y$ for $f^{-1}\sG$.  
\end{rmk}

\begin{rmk}[Inverse images of divisors and subschemes]
\label{rmk:Inverse_images_divisors}
Let $f:Y\to X$ be a morphism of schemes and $D$ be a Cartier divisor on $X$. We refer to G\"ortz and Wedhorn \cite[Chapter 11, Definition 11.49, p. 315]{Gortz_Wedhorn_algebraic_geometry_v1} for the definition of the \emph{pullback} or \emph{inverse image} $f^*D$ over $X$ and refer to \cite[Chapter 11, Proposition 11.50, p. 316]{Gortz_Wedhorn_algebraic_geometry_v1} for two criteria that indicate when it is possible to form the inverse image of a divisor. Let $\iota:Z\to X$ be an immersion of schemes. We refer to \cite[Chapter 4, Section 4.11, p. 112]{Gortz_Wedhorn_algebraic_geometry_v1} for the definition of the \emph{inverse image} $f^{-1}(Z)$ of $Z$ under $f$, a subscheme of $Y$. The authors note that $f^{-1}(Z)$ is a closed subscheme of $X$ if $Z$ is a closed subscheme of $Y$ and by \cite[Chapter 4, Proposition 4.20, p. 103]{Gortz_Wedhorn_algebraic_geometry_v1}, one knows that if $Z$ is an open subscheme of $Y$, then $f^{-1}(Z)$ is an open subscheme of $X$.
\end{rmk}

\begin{rmk}[Inverse image ideal sheaf]
\label{rmk:Inverse_image_ideal_sheaf}  
(See Koll\'ar \cite[Notation 3.15, p. 134]{Kollar_lectures_resolution_singularities}.)  
Let $f : Y \to X$ be a morphism of schemes and $\sI \subset \sO_X$ be an ideal sheaf. in Definition \ref{defn:Principalization_monomialization_ideal_sheaf}, Koll\'ar uses $f^*\sI$ to denote the \emph{inverse image ideal sheaf of $\sI$}, that is, the ideal sheaf generated by the pullbacks of local sections of $\sI$. (It is denoted by $f^{-1}\sI\cdot\sO_Y$ or $\sI\cdot\sO_Y$ by Hartshorne \cite[Chapter II, Section 7, p. 163]{Hartshorne_algebraic_geometry}, as opposed to the usual sheaf-theoretic pullback, also denoted by $f^*\sI$, and which may be different --- see Hartshorne \cite[Chapter II, Section 7, Caution 7.12.2, p. 163]{Hartshorne_algebraic_geometry}.
\end{rmk}

\begin{defn}[Simple normal crossings or monomial ideal at a point]
\label{defn:Simple_normal_crossings_monomial_ideal_at_point}  
(See Cutkosky \cite[Section 4.2, p. 40]{Cutkosky_resolution_singularities} or Koll{\'a}r \cite[Note on Terminology 3.16, p. 134]{Kollar_lectures_resolution_singularities}.)    
Let $X$ be a smooth scheme and $\sI \subset \sO_X$ be a locally principal ideal. One says that $\sI$ has \emph{simple normal crossings} (or is \emph{monomial}) \emph{at a point} $p \in X$ if there exist regular parameters $\{x_1,\ldots,x_d\} \subset \sO_{X,p}$ such that
\[
\sI_p = x_1^{m_1}\cdots x_n^{m_d}\sO_{X,p},
\]
for some $m_i\in\NN$, with $i=1,\ldots,d$.
\end{defn}

\begin{defn}[Locally monomial ideal sheaf]
\label{defn:Locally_monomial_ideal_sheaf}
(See Koll{\'a}r \cite[Note on Terminology 3.16, p. 143]{Kollar_lectures_resolution_singularities}.)  
Suppose that $X$ is a smooth algebraic variety and $\sI \subset \sO_X$ is an ideal sheaf. One says that $\sI$ is \emph{locally monomial} if one of the following equivalent conditions hold:
\begin{enumerate}
\item $\sI$ is monomial at every point $p \in X$ as in Definition \ref{defn:Simple_normal_crossings_monomial_ideal_at_point}.
\item $\sI$ is the ideal sheaf of a \emph{simple normal crossing divisor} as in Definition \ref{defn:Subvariety_having simple_normal_crossings_with_divisor}.
\end{enumerate}
\end{defn}

\section[Theorems on monomialization and resolution of singularities for schemes]{Main theorems on monomialization and resolution of singularities for schemes}
\label{sec:Main_theorems_monomialization_resolution_singularities_schemes}
See Hironaka \cite{Hironaka_1964-I-II} for the first proof of monomialization and resolution of singularities for schemes of arbitrary dimension over a field of characteristic zero. See also Hironaka \cite{Hironaka_infinitely_near_singular_points} and Aroca, Hironaka, Vicente \cite{Aroca_Hironaka_Vicente_theory_maximal_contact, Aroca_Hironaka_Vicente_desingularization_theorems} for related results specific to complex analytic spaces and real analytic spaces. 

For simpler proofs and strengthenings of Hironaka's theorems, see Abramovich \cite{Abramovich_2011}, Bierstone and Milman \cite{Bierstone_Milman_1997}, Bravo and Villamayor \cite{Bravo_Villamayor_2003}, Encinas and Hauser \cite{Encinas_Hauser_2002}, Encinas and Villamayor \cite{Encinas_Villamayor_1998, Encinas_Villamayor_2000}, Villamayor \cite{Villamayor_1989, Villamayor_1992}, and W{\l}odarczyk \cite{Wlodarczyk_2005, Wlodarczyk_2008}.

For expository introductions to resolution of singularities and its proof, see Abramovich \cite{Abramovich_2018icm}, Bierstone and Milman \cite{Bierstone_Milman_1999}, Cutkosky \cite{Cutkosky_resolution_singularities}, Faber and Hauser \cite{Faber_Hauser_2010}, Goward \cite{Goward_2005}, Hauser \cite{Hauser_2003}, Hironaka \cite{Hironaka_1963}, and Koll{\'a}r \cite{Kollar_2005arxiv, Kollar_lectures_resolution_singularities}, and Lipman \cite{Lipman_1975}.

\begin{thm}[Monomialization of an ideal sheaf over a scheme]
\label{thm:Monomialization_ideal_sheaf_algebraic_scheme}
(See Koll{\'a}r \cite[Theorem 3.21, p. 136, Theorem 3.26, p. 138, and Theorem  3.35, p. 135]{Kollar_lectures_resolution_singularities}.)
Let $X$ be a smooth scheme of finite type over a field of characteristic zero and $\sI \subset \sO_X$ be an ideal sheaf that is not zero on any irreducible component of $X$. If $E$ is a simple normal crossing divisor on $X$, then there is a sequence of smooth blowup morphisms
\begin{equation}
  \label{eq:Composition_blow-up_morphisms}
  \Pi : X' = X_r \xrightarrow{\pi_{r-1}} X_{r-1} \xrightarrow{\pi_{r-2}} \cdots \xrightarrow{\pi_1} X_1 \xrightarrow{\pi_0} X_0 = X
\end{equation}
such that the following  hold:
\begin{enumerate}
\item The centers of the blowups are smooth and have simple normal crossings with $E$;  
\item The pullback $\Pi^*\sI \subset \sO_{X'}$ is the ideal sheaf of a simple normal crossing divisor;
\item $\Pi : X' \to X$ is an isomorphism over $X \less \cosupp\sI$, where $\cosupp\sI$ (or $\supp(\sO_X/\sI)$) is the cosupport of $\sI$;
\item The assignment to $X$ of a blowup sequence is functorial in that it commutes with smooth morphisms in the sense of the forthcoming Definition \ref{defn:Kollar_3-34-1} and, when $E=\varnothing$, commutes with closed embeddings in the sense of the forthcoming Definition \ref{defn:Kollar_3-34-3}.  
\end{enumerate}
\end{thm}

\begin{rmk}[Specialization of monomialization to an ideal sheaf over an algebraic variety]
If $X$ is a smooth algebraic variety, then so also is $X'$ and $\Pi$ is birational and projective --- see Koll{\'a}r \cite[Theorem 3.21, p. 136 and Theorem 3.26, p. 138]{Kollar_lectures_resolution_singularities}. In the relatively simple \cite[Theorem 3.17, p. 135]{Kollar_lectures_resolution_singularities}, Koll{\'a}r establishes principalization but not monomialization.
\end{rmk}  

\begin{rmk}[Other versions of monomialization of an ideal sheaf over a scheme]
\label{rmk:Monomialization_ideal_sheaf_algebraic_scheme_other_versions}
Versions of Theorem \ref{thm:Monomialization_ideal_sheaf_algebraic_scheme} are provided by Bierstone and Milman \cite[Theorem 1.10, p. 216, and p. 221]{Bierstone_Milman_1997}, Hironaka \cite[Chapter 0, Section 5, Main Theorem II, p. 142]{Hironaka_1964-I-II}, and W{\l}odarczyk \cite[Theorem 1.0.1, p. 781]{Wlodarczyk_2005}. 
\end{rmk}

As a consequence of Theorem \ref{thm:Monomialization_ideal_sheaf_algebraic_scheme} one obtains strong resolution of singularities for schemes of finite type over a field of characteristic zero.

\begin{thm}[Strong resolution of singularities for schemes of finite type over a field of characteristic zero]
\label{thm:Resolution_singularities_algebraic_scheme}
(See Koll{\'a}r \cite[Theorem 3.27, p. 139 and Theorem 3.36, p. 146]{Kollar_lectures_resolution_singularities}.)
If $X$ is a scheme of finite type over a field of characteristic zero, then there is a sequence of smooth blowup morphisms as in \eqref{eq:Composition_blow-up_morphisms} such the following hold:
\begin{enumerate}
\item $X'$ is smooth;
\item $\Pi:X' \to X$ is an isomorphism over the smooth locus $X\less X_\sing$; 
\item $\Pi^{-1}(X_\sing)$ is a divisor with simple normal crossings;
\item The assignment to $X$ of a blowup sequence is functorial in that it commutes with smooth morphisms in the sense of the forthcoming Definition \ref{defn:Kollar_3-34-1}.
\end{enumerate}
\end{thm}

\begin{rmk}[Other versions of strong resolution of singularities for schemes of finite type over a field of characteristic zero]
\label{rmk:Resolution_singularities_algebraic_scheme_other_versions}
Versions of Theorem \ref{thm:Resolution_singularities_algebraic_scheme} are provided by Bierstone and Milman \cite[Theorem 1.6, p. 215, and p. 221]{Bierstone_Milman_1997}, Hironaka \cite[Chapter 0, Section 3, Main Theorem I, p. 132, and Chapter 0, Section 4, Main Theorem I${}^*$, p. 138]{Hironaka_1964-I-II}, and W{\l}odarczyk \cite[Theorem 1.0.3, p. 782]{Wlodarczyk_2005}. In \cite[Theorem 1.0.3, p. 782]{Wlodarczyk_2005}, W{\l}odarczyk states that $\Pi:X'\to X$ is a proper birational morphism; $\Pi$ is functorial with respect to smooth morphisms in the sense that for any smooth morphism $\phi:Y\to X$, there is a natural lifting $\phi':Y' \to X'$ which is a smooth morphism; and $\Pi:X'\to X$ is an isomorphism over $X\less X_\sing$.
\end{rmk}  

\begin{thm}[Embedded resolution of singularities for a scheme of finite type over a field of characteristic zero]
\label{thm:Embedded_resolution_of_singularities_algebraic_scheme}
(See Hauser \cite[p. 329 and pp. 334--335]{Hauser_2003} when $X$ is a smooth scheme and $Y$ is a subscheme and W{\l}odarczyk \cite[Theorem 1.0.2, p. 781]{Wlodarczyk_2005} when $X$ is a smooth variety and $Y$ is a subvariety.) 
Let $Y$ be a reduced singular scheme of finite type over a field of characteristic zero and $Y\hookrightarrow X$ be a closed embedding of $Y$ into a smooth scheme $X$. Then there there is a sequence of smooth blowup morphisms as in \eqref{eq:Composition_blow-up_morphisms} such the following hold:
\begin{enumerate}
\item\label{item:Embedded_resolution_of_singularities_algebraic_scheme_smooth_centers}
The blowup morphisms $\pi_i:X_i\to X_{i-1}$ have smooth closed centers $Z_{i-1}\subset X_{i-1}$ for $i=1,\ldots,r$, where $X_0 := X$ and $X' := X_r$;
\item\label{item:Embedded_resolution_of_singularities_algebraic_scheme_exceptional_divisors_snc}  
The exceptional divisor $E_i$ of the induced morphism $\Pi_i = \pi_1\circ\cdots\circ\pi_i:X_i\to X_0$ has only simple normal crossings and $Z_i$ has simple normal crossings with $E_i$, for $i=1,\ldots,r$; here, each $E_i$ is the inverse image in $X_i$ of the first $i$ blowup centers $Z_0,\ldots,Z_{i-1}$ under the preceding blowup morphisms for $i=1,\ldots,r$.
\item\label{item:Embedded_resolution_of_singularities_algebraic_scheme_centers_disjoint_Ysm}
If $Y_i \subset X_i$ denotes the strict transform of $Y$, then the center $Z_i\subset X_i$ is disjoint from the preimage $\Pi_i^{-1}(Y_\sm) \subset Y_i \subset X_i$ of the open subset $Y_\sm \subset Y$ of points where $Y$ is smooth for $i=0,\ldots,r-1$ and we write $\Pi_0 = \pi_0 = \id_{X_0}$;
\item\label{item:Embedded_resolution_of_singularities_algebraic_scheme_strict_transform_smooth}
The strict transform $Y' := Y_r$ of $Y$ is smooth and has only simple normal crossings with the exceptional divisor $E' := E_r$;
\item\label{item:Embedded_resolution_of_singularities_algebraic_scheme_commutes_morphisms_embeddings}
The morphism $\Pi:(X',Y')\to (X,Y)$ commutes with smooth morphisms (see the forthcoming Definition \ref{defn:Kollar_3-34-1}) and smooth embeddings of the ambient scheme $X$ (see the forthcoming Definition \ref{defn:Kollar_3-34-3});
\item\label{item:Embedded_resolution_of_singularities_algebraic_scheme_isomorphism_off_exceptional_divisor}
The morphism $\Pi:X'\to X$ is proper and induces an isomorphism $X'\less E' \to X\less Z$ outside the final exceptional divisor $E' \subset X'$, where $Z \subset X$ denotes the image of $E'$, that is, the image of all the intermediate blowup centers $Z_i$ for $i=0,\ldots,r-1$.
\end{enumerate}
\end{thm}

\begin{rmk}[Other versions of embedded resolution of singularities for schemes of finite type over a field of characteristic zero]
\label{rmk:Embedded_resolution_singularities_algebraic_scheme_other_versions}
Versions of Theorem \ref{thm:Embedded_resolution_of_singularities_algebraic_scheme} are provided by Bierstone and Milman \cite[Theorem 1.6, p. 215, and p. 221]{Bierstone_Milman_1997} and Hironaka \cite[Chapter 0, Section 5, Corollary 3, p. 146]{Hironaka_1964-I-II}. Hauser \cite[p. 329]{Hauser_2003} provides a useful guide to the literature on resolution of singularities in terms of key properties:
\begin{inparaenum}
\item\label{item:Hauser_E1} \emph{Explicitness.} $\Pi$ is a composition of blowups of $X$ in regular closed centers $Z_i$ transversal to the exceptional loci.
\item\label{item:Hauser_E2} \emph{Embeddedness.} The strict transform $Y'$ of $X$ is regular and transversal to the exceptional locus in $X'$.
\item\label{item:Hauser_E3} \emph{Excision.} The morphism $\Pi:Y' \to Y$ is independent of the embedding of $Y$ in $X$.
\item\label{item:Hauser_E4} \emph{Equivariance.} $\Pi$ commutes with smooth morphisms $X_1 \to X$, embeddings $X \to X_2$, and field extentions.
\item\label{item:Hauser_E5} \emph{Effectiveness}. The centers of blowups are given as the top locus of a local upper semicontinuous invariant of $Y$.
\end{inparaenum}
Existence with Properties \eqref{item:Hauser_E1} and \eqref{item:Hauser_E2} was proved by Hironaka \cite{Hironaka_1964-I-II}. Constructive proofs of resolution of singularities with Properties \eqref{item:Hauser_E1}, \eqref{item:Hauser_E2}, \eqref{item:Hauser_E3}, \eqref{item:Hauser_E4}, and \eqref{item:Hauser_E5} were provided by Villamayor \cite{Villamayor_1989, Villamayor_1989}, Bierstone and Milman \cite{Bierstone_Milman_1997}, Encinas and Villamayor \cite{Encinas_Villamayor_1998, Encinas_Villamayor_2000}, Encinas and Hauser \cite{Encinas_Hauser_2002}, and Bravo and Villamayor \cite{Bravo_Villamayor_2003}. For algorithmic implementations of and computer programs for resolution of singularities for algebraic varieties, we refer the reader to Bodn{\'a}r and Schicho \cite{Bodnar_Schicho_2000jsc, Bodnar_Schicho_2000bb, Bodnar_Schicho_2001}. See Hauser \cite{Hauser_2003, Hauser_2014} for very extensive bibliographic details and many additional references.
\end{rmk}

\begin{rmk}[Location of the blowup centers]
\label{rmk:Location_blow-up_centers}
For Item \eqref{item:Embedded_resolution_of_singularities_algebraic_scheme_centers_disjoint_Ysm} in Theorem \ref{thm:Embedded_resolution_of_singularities_algebraic_scheme}, W{\l}odarczyk emphasizes in \cite[Theorem 1.0.2 (b), p. 781]{Wlodarczyk_2005} that the blowup centers $Z_i$ need not be contained in $Y_i$ or disjoint from $(Y_i)_\sm$. Moreover, in \cite[Theorem 1.0.2 (b), p. 781]{Wlodarczyk_2005}, W{\l}odarczyk simply writes $Y_\sm \subset Y_i$ rather than $\Pi_i^{-1}(Y_\sm) \subset Y_i$ as we do: he uses the isomorphism $\Pi_i:X_i\less E_i \to X_0\less E$ to identify $\Pi_i^{-1}(Y_\sm) \cong Y_\sm$ for $i=0,\ldots,r$, where $E \subset X_0$ here denotes the image of $E_i \subset X_i$, that is, the image of all the intermediate blowup centers $Z_j$ for $j=0,\ldots,i-1$ when $i\geq 1$ and $E=Z_0$ when $i=0$.

When $i=0$, Item \eqref{item:Embedded_resolution_of_singularities_algebraic_scheme_centers_disjoint_Ysm} in Theorem \ref{thm:Embedded_resolution_of_singularities_algebraic_scheme} implies that $Z_0\subset X_0$ is disjoint from $\pi_0^{-1}(Y_\sm) = Y_\sm \subset Y_0 \subset X_0$ and similarly, if $i=1$, it implies that $Z_1\subset X_1$ is disjoint from $\pi_1^{-1}(Y_\sm) \subset Y_1 \subset X_1$. In particular, $\Pi:\Pi^{-1}(Y_\sm) \to Y_\sm$ is an isomorphism, $\Pi^{-1}(Y_\sm) \subset Y' \subset X'$ is disjoint from the preimages in $X$ of the blowup centers $Z_i$ for $i=0,\ldots,r-1$, and $Y_\sm \subset Y \subset X$ is disjoint from the preimages in $X$ of the blowup centers $Z_i$ for $i=0,\ldots,r-1$.

We may compare the preceding remarks with those of Koll{\'a}r \cite[Section 3.3, p. 120, and Section 3.5, p. 122]{Kollar_lectures_resolution_singularities} who implies that, in the setting of Theorem \ref{thm:Resolution_singularities_algebraic_scheme} for strong resolution, one has $\Pi_i(Z_i) \subset X_\sing$ for $i=0,\ldots,r-1$, where $X = X_\sm \sqcup X_\sing$.
\end{rmk} 

\begin{defn}[Strict transform and total transform in resolution of singularities in the category of algebraic varieties]
\label{defn:Strict_transform_and_total_transform_algebraic_varieties}
In the specialization of Theorem \ref{thm:Embedded_resolution_of_singularities_algebraic_scheme} to the category of algebraic varieties, the \emph{strict transform} $Y'$ is the Zariski closure of $\Pi^{-1}(Y \less Z)$ in $X'$ (see Bierstone and Milman \cite[Section 3, The strict transform, p. 237]{Bierstone_Milman_1997} and \cite[Section 1.6, p. 50]{Bierstone_Milman_1999}, and Hauser \cite[Chapter 0, Section $-1$, p. 335]{Hauser_2003}).  The preimage $Y^* := \Pi^{-1}(Y)$ is the \emph{total transform} of $Y$ (see Hartshorne \cite[Chapter V, Section 5, p. 410]{Hartshorne_algebraic_geometry}, Hauser \cite[Chapter 0, Section $-1$, p. 334]{Hauser_2003}, and Koll\'ar \cite[Section 3.3, p. 138 and Definition 3.65, p. 163]{Kollar_lectures_resolution_singularities}).
\end{defn}

To illustrate some of the preceding ideas, we recall the following simple example. For additional examples, we refer the reader to Hauser \cite{Hauser_2014} and Hauser and Faber \cite{Faber_Hauser_2010}.

\begin{exmp}[Resolution of the quadratic cone]
\label{exmp:Quadratic_cone}  
Following Smith et al. \cite[Section 7.1, Example, pp. 103--106]{Smith_Kahanpaa_Kekalainen_Traves_invitation_algebraic_geometry}, consider the quadratic cone $Y = \VV_\CC(x^2+y^2-z^2) \subset X = \CC^3$ and its blowup $\Bl_0(Y) \subset \Bl_0(X)$ at the origin $p_0=0\in X$. For the coordinate domain $U_z \subset \Bl_0(X)$ with local coordinates $(u,v,z)$ corresponding to the chart $\varphi_z:U_z\to\CC^3$, the blowup morphism $\pi_X:\Bl_0(X)\to X$ is given by $(u,v,z) \mapsto (x,y,z) = (zu,zv,z)$ and thus $\Bl_0(Y)\cap U_z$ is isomorphic to $\VV_\CC(u^2+v^2-1)$ and $E\cap U_z$ is isomorphic to $\VV_\CC(z)$. In particular, $\pi_X^{-1}(p_0)\cap \Bl_0(Y) \cap U_z$ is isomorphic to $\VV_\CC(u^2+v^2-1)\cap \VV_\CC(z) \subset \CC^3$.

Suppose that $\rho:S^1\to\Gl(3,\CC)$ is a unitary representation as in \eqref{eq:Circle_matrix_representation} and that $\rho$ restricts to an $S^1$ action on $Y$. Hence, $l_1=l_2=l_3=l$ for some integer $l\in\ZZ$ and the induced action of $e^{i\theta}\in S^1$ on $U_z$ is given by $(u,v,z) \mapsto (u,v,e^{il\theta}z)$, so each point in the variety $\VV_\CC(u^2+v^2-1,z) \subset \CC^3$ is a fixed point of the induced $S^1$ action on $\Bl_0(Y)$.

Denoting $F(x,y,z) := x^2+y^2-z^2$, we see that $DF(x,y,z) = (2x,2y,-2z)$ for all $(x,y,z) \in \CC^3$ and so $\Ran DF(x,y,z)=\CC$ unless $(x,y,z)=(0,0,0)$, in which case $T_0Y = \CC^3$, while $\dim_\CC T_pY = 2$ for all $p\in Y\less\{0\}$ and all such points are smooth.

Note that intuition provided by illustrations in the case $\KK=\RR$ can be misleading when we consider the corresponding varieties over $\KK=\CC$. For example, $\VV_\RR(x^2+y^2-z^2,z) = (0,0,0) \in \RR^3$ whereas $\VV_\CC(x^2+y^2-z^2,z) = \{(x,\pm i x,0): x \in \CC\} \subset \CC^3$, which is smooth of complex dimension one away from the origin.
\qed
\end{exmp}

\section{Main theorems on resolution of singularities for analytic spaces}
\label{sec:Principalization_and_resolution_singularities_analytic_spaces}
For versions of resolution of singularities for analytic spaces, we refer to Bierstone and Milman \cite{Bierstone_Milman_1997}, Hironaka \cite{Hironaka_1964-I-II}, Koll{\'a}r \cite{Kollar_lectures_resolution_singularities} and W{\l}odarczyk \cite{Wlodarczyk_2008}. Koll{\'a}r \cite[p. 135 and Section 3.44]{Kollar_lectures_resolution_singularities} observes that a good resolution method should also work for complex, real or $p$-adic analytic spaces. Indeed, Hironaka \cite{Hironaka_1964-I-II} explicitly allows for analytic spaces, as do Bierstone and Milman \cite{Bierstone_Milman_1997} and W{\l}odarczyk \cite{Wlodarczyk_2008}. Following Hironaka \cite[p. 111]{Hironaka_1964-I-II} and W{\l}odarczyk \cite[p. 34]{Wlodarczyk_2008}, we shall restrict our attention to the case $\KK=\RR$ or $\CC$ and refer the reader to the literature for analytic spaces over other fields. We begin with the following analogue of Theorem \ref{thm:Resolution_singularities_algebraic_scheme}.

\begin{thm}[Strong resolution of singularities for analytic spaces]
\label{thm:Resolution_singularities_analytic_space}
(See Koll{\'a}r \cite[Theorem 3.45]{Kollar_lectures_resolution_singularities}.)
Let $\KK$ be a locally compact field of characteristic zero. If $X$ is a $\KK$-analytic space, then the conclusions of Theorem \ref{thm:Resolution_singularities_algebraic_scheme} continue to hold.
\end{thm}

\begin{rmk}[Other versions of strong resolution of singularities for analytic spaces]
\label{rmk:Resolution_singularities_analytic_space_other_versions}
Versions of Theorem \ref{thm:Resolution_singularities_analytic_space} are provided by Bierstone and Milman \cite[Theorem 1.6, p. 215, and p. 221]{Bierstone_Milman_1997} for a field $\KK$ with a complete valuation, Hironaka \cite[Chapter 0, Section 6, Main Theorem I${}'$, p. 151 and Main Theorem I${}''(n)$, p. 155]{Hironaka_1964-I-II} for $\KK=\CC$ and \cite[Chapter 0, Section 6, Main Theorem I${}''(n)$, p. 158]{Hironaka_1964-I-II} for $\KK=\RR$, and W{\l}odarczyk \cite[Theorem 2.0.1, p. 34]{Wlodarczyk_2008} for $\KK=\RR$ or $\CC$. In \cite[Theorem 2.0.1, p. 34]{Wlodarczyk_2008}, W{\l}odarczyk states that $\Pi:X'\to X$ is a bimeromorphic proper morphism and that $\Pi$ is functorial with respect to local analytic isomorphisms in the following sense: for any local analytic isomorphism $\phi:Y\to X$, there is a natural lifting $\phi':Y' \to X'$ which is a local analytic isomorphism. In this context, \emph{proper} means that the preimage by $\Pi$ of any compact subset of $X$ is compact in $X'$ (see Abramovich \cite[Section 1.2]{Abramovich_2018icm}, Hauser \cite[Chapter 0, Section $-1$, p. 334]{Hauser_2003}, and Hironaka \cite[Introduction, p. 111]{Hironaka_1964-I-II}); for other definitions of properness, see Cutkosky \cite[Section 2.3, Definition 2.12, p. 9]{Cutkosky_resolution_singularities} and Hartshorne \cite[Chapter II, Section 4, p. 100, Definition]{Hartshorne_algebraic_geometry}.
\end{rmk}

Koll{\'a}r \cite[Section 3.44, p. 149]{Kollar_lectures_resolution_singularities} asserts that the proofs in \cite[Chapter 3]{Kollar_lectures_resolution_singularities} extend (locally) to analytic spaces over locally compact fields. Thus, we have the following analogue of Theorem \ref{thm:Monomialization_ideal_sheaf_algebraic_scheme}.

\begin{thm}[Monomialization of an ideal sheaf over a smooth analytic space]
\label{thm:Monomialization_ideal_sheaf_analytic_space}
(See Koll{\'a}r \cite[Section 3.44, p. 149]{Kollar_lectures_resolution_singularities}.)
Let $\KK$ be a locally compact field of characteristic zero. If $X$ is a smooth $\KK$-analytic space, then the conclusions of Theorem \ref{thm:Monomialization_ideal_sheaf_algebraic_scheme} continue to hold.
\end{thm}

\begin{rmk}[Other versions of monomialization of an ideal sheaf over a smooth analytic space]
\label{rmk:Monomialization_ideal_sheaf_analytic_space_other_versions}
Versions of Theorem \ref{thm:Monomialization_ideal_sheaf_analytic_space} are provided by Bierstone and Milman \cite[Theorem 1.10, p. 216]{Bierstone_Milman_1997} for a field $\KK$ with a complete valuation, Hironaka \cite[Chapter 0, Section 7, Main Theorem II${}'$(N), p. 156]{Hironaka_1964-I-II} for $\KK=\CC$ and \cite[Chapter 0, Section 7, Main Theorem II${}''(N)$, p. 158]{Hironaka_1964-I-II} for $\KK=\RR$, and W{\l}odarczyk \cite[Theorem 2.0.3, p. 35]{Wlodarczyk_2008} for $\KK=\RR$ or $\CC$. In \cite[Theorem 2.0.3, p. 35]{Wlodarczyk_2008}, W{\l}odarczyk states that $\Pi:X'\to X$ is a proper morphism and that $\Pi$ commutes with local analytic isomorphisms and embeddings of ambient varieties.
\end{rmk}

Lastly, we have the following analogue of Theorem \ref{thm:Embedded_resolution_of_singularities_algebraic_scheme}.

\begin{thm}[Embedded resolution of singularities for an analytic space]
\label{thm:Embedded_resolution_of_singularities_analytic_space}
(See W{\l}odarczyk \cite[Theorem 2.0.2, p. 34]{Wlodarczyk_2008}.) 
Let $\KK=\RR$ or $\CC$ and $Y$ be a $\KK$-analytic subspace of a $\KK$-analytic manifold $X$. Then there exist a $\KK$-analytic manifold $X'$, a simple normal crossing locally finite divisor $E$ on $X'$, bimeromorphic proper morphism $\Pi:X'\to X$ such that the strict transform $Y'\subset X'$ is smooth and has simple normal crossings with the divisor $E$; the support of the divisor $E$ is the the exceptional locus of $\Pi$; and the morphism $\Pi$ locally factors into a sequence of blowups at smooth centers. In particular, all the conclusions of Theorem \ref{thm:Embedded_resolution_of_singularities_algebraic_scheme} continue to hold.
\end{thm}

\begin{rmk}[Other versions of embedded resolution of singularities for an analytic space]
\label{rmk:Embedded_resolution_singularities_analytic_space_other_versions}
Versions of Theorem \ref{thm:Resolution_singularities_analytic_space} are provided by Bierstone and Milman \cite[Theorem 1.6, p. 215 and p. 221]{Bierstone_Milman_1997} for a field $\KK$ with a complete valuation, Hauser \cite[p. 329 and Chapter 0, Section $-1$, pp. 334--335]{Hauser_2003} and \cite[Definition 7.6, p.34]{Hauser_2014}, Hironaka \cite[Chapter 0, Section 7, Main Theorem I${}'$, p. 151 and Main Theorem I${}''(n)$, p. 155]{Hironaka_1964-I-II} for $\KK=\CC$ and Hironaka \cite[Chapter 0, Section 7, Main Theorem I${}''(n)$, p. 158]{Hironaka_1964-I-II} for $\KK=\RR$, and Koll{\'a}r \cite[Section 3.1]{Kollar_lectures_resolution_singularities} for a locally compact field  $\KK$ of characteristic zero.
\end{rmk}

Lastly, we have the following analogue of Definition \ref{defn:Strict_transform_and_total_transform_algebraic_varieties} and generalization of Definition \ref{defn:Strict_transform_analytic_space} from the case of one blowup morphism to a composition of blowup morphisms.

\begin{defn}[Strict transform and total transform in resolution of singularities in the category of analytic spaces]
\label{defn:Strict_transform_and_total_transform_analytic_spaces}
(See Bierstone and Milman \cite[Section 2.5, p. 807]{Bierstone_Milman_1989}.)  
In the setting of Theorem \ref{thm:Embedded_resolution_of_singularities_algebraic_scheme}, the \emph{strict transform} $Y'$ of $Y \subset X$ is the smallest closed analytic subspace of $X'$ such that the resolution morphism $\Pi:X'\to X$ induces an isomorphism $Y'\less E \cong Y\less Z$ of analytic spaces, while the preimage $Y^* := \Pi^{-1}(Y)$ is the \emph{total transform} of $Y$.
\end{defn}

\section{Equivariant blowups and resolution of singularities}
\label{sec:Equivariant_resolution_of_singularities}
As noted by Koll{\'a}r \cite[Section 3.4.1, p. 121, and Proposition 3.9.1, p. 125]{Kollar_lectures_resolution_singularities}, where $X$ is assumed to be a scheme, the action of an algebraic group $G$ on $X$ lifts to an action on its functorial resolution $X'$. See Abramovich and Wang \cite{Abramovich_Wang_1997}, Bierstone, Grigoriev, Milman and W{\l}odarczyk \cite[Theorems 1.0.1, 1.02, and 1.03, p. 195]{Bierstone_Grigoriev_Milman_Wlodarczyk_2011}, Encinas and Villamayor \cite{Encinas_Villamayor_2000}, Hauser \cite[p. 329 and Chapter 0, Section $-1$, p. 336]{Hauser_2003}, Reichstein and Youssin \cite{Reichstein_Youssin_2000, Reichstein_Youssin_2002pams, Reichstein_Youssin_2002pjm}, and Villamayor \cite[Corollary 7.6.3, p. 669]{Villamayor_1992}, \cite{Villamayor_2014} for expositions and proofs of results related to equivariant resolution of singularities. In particular, Encinas and Villamayor \cite{Encinas_Villamayor_2000} and Reichstein and Youssin \cite{Reichstein_Youssin_2000} discuss equivariant resolution of singularities in detail. We shall ultimately apply these equivariance properties (for the groups $G=S^1$ or $\CC^*$) acting on a complex $G$-manifold when resolving the singularities of a $G$-invariant, closed, complex analytic subspace.

In order to precisely describe the equivariance properties of resolution of singularities, we begin by recalling the following

\begin{defn}[Blowup sequence functors for schemes]
\label{defn:Kollar_3-31}  
(See Koll\'ar \cite[Definition 3.31, p. 142]{Kollar_lectures_resolution_singularities}.)
A \emph{blowup sequence functor} is a functor $\sB$ whose
\begin{enumerate}
\item Inputs are triples $(X,\sI,E)$, where $X$ is a scheme, $\sI \subset \sO_X$ an ideal sheaf that is nonzero on every irreducible component, and $E$ a divisor on $X$ with ordered index set, and
\item Outputs are blowup sequences
  \begin{equation}
    \label{eq:Kollar_3-31_blowup_sequence}
    \begin{tikzcd}
      \Pi: X_r \arrow[r, "\pi_{r-1}"] &X_{r-1} \arrow[r, "\pi_{r-2}"]
      &\cdots \arrow[r, "\pi_1"] &X_1 \arrow[r, "\pi_0"] &X_0
      \\
      & Z_{r-1} \arrow[hookrightarrow, u] &\cdots &Z_1 \arrow[hookrightarrow, u]
      &Z_0 \arrow[hookrightarrow, u]
  \end{tikzcd}
  \end{equation}
with specified blowup centers $Z_i$ for $i=0,\ldots,r-1$ and $X_0 := X$ while $X' := X_r$. Here the length of the sequence $r$, the schemes $X_i$, and the centers $Z_i$ all depend on $(X,\sI,E)$. 
\end{enumerate}
\end{defn}

If each $Z_i$ in Definition \ref{defn:Kollar_3-31} is smooth, then a nontrivial blowup $\pi_i:X_{i+1} \to X_i$ uniquely determines $Z_i$, so one can omit $Z_i$ from the notation; however, different centers may lead to the same birational map \cite[p. 142]{Kollar_lectures_resolution_singularities}.  

\begin{defn}[Blowup sequence functors for schemes commute with smooth morphisms]
\label{defn:Kollar_3-34-1}  
(See Koll\'ar \cite[Paragraph 3.34.1, p. 144]{Kollar_lectures_resolution_singularities}.)  
A blowup sequence functor $\sB$ \emph{commutes with smooth morphisms} if $\sB$ commutes with every smooth surjective morphism $h:Y\to X$ of schemes,
\[
  \sB(Y,h^*\sI,h^{-1}(E)) = h^*\sB(X,\sI,E),
\]
and for every smooth morphism $h$, the blowup sequence $\sB(Y,h^*\sI,h^{-1}(E))$ is obtained from
the pullback $h^*\sB(X,\sI,E)$ of the sequence $\sB(X,\sI,E)$ by deleting every blowup $h^*\pi_i$ whose
center is empty and reindexing the resulting blowup sequence.
\end{defn}

\begin{defn}[Blowup sequence functors for schemes commute with closed embeddings]
\label{defn:Kollar_3-34-3}  
(See Koll\'ar \cite[Paragraph 3.34.3, p. 145]{Kollar_lectures_resolution_singularities}.)  
A blowup sequence functor $\sB$ \emph{commutes with closed embeddings} if
\[
  \sB(X,\sI_X,E) = j_*\sB(Y,\sI_Y,E|_Y),
\]  
whenever
\begin{itemize}
\item $j:Y\hookrightarrow X$ is a closed embedding of smooth schemes,
\item $0 \neq \sI_Y \subset \sO_Y$ and $0 \neq \sI_X \subset \sO_X$ are ideal sheaves such that $\sO_X/\sI_X = j_*(\sO_Y/\sI_Y)$, and
\item $E$ is a simple normal crossing divisor on $X$ such that $E|_Y$ is also a simple normal crossing divisor on $Y$.
\end{itemize}
\end{defn}

As a consequence of the fact that blowup sequences commute with smooth morphisms in Theorems \ref{thm:Monomialization_ideal_sheaf_algebraic_scheme}, \ref{thm:Resolution_singularities_algebraic_scheme}, and \ref{thm:Embedded_resolution_of_singularities_algebraic_scheme} in the sense of Definition \ref{defn:Kollar_3-34-1}, one immediately obtains the following corollaries. However, we emphasize that (as a consequence of Definition \ref{defn:Kollar_3-34-1}) each blowup center $Z_{i-1} \subset X_{i-1}$ is $G$-invariant and each blowup morphism $\pi_i:X_i\to X_{i-1}$ is $G$-equivariant, for $i=1,\ldots,r$, as in Encinas and Villamayor \cite[Paragraph 2.11, p. 155]{Encinas_Villamayor_2000}. See Reichstein and Youssin \cite[Theorem 1.1, p. 1019]{Reichstein_Youssin_2000} for one example of such a statement when $G$ is an algebraic group acting algebraic variety $X$ and Encinas and Villamayor \cite{Encinas_Villamayor_2000} for many more examples.

\begin{cor}[Equivariance of monomialization of an ideal sheaf over a scheme]
\label{cor:Monomialization_ideal_sheaf_algebraic_scheme_equivariant}
(See Bierstone, Grigoriev, Milman, and W{\l}odarczyk \cite[Theorem 1.0.1, p. 195]{Bierstone_Grigoriev_Milman_Wlodarczyk_2011}.)  
Continue the hypotheses of Theorem \ref{thm:Monomialization_ideal_sheaf_algebraic_scheme}. Then the action of an algebraic group $G$ on the scheme $X$ lifts to an action of $G$ on its functorial resolution $X'$ and the morphism $\Pi:(X',\sI') \to (X,\sI)$ is $G$-equivariant, where $\sI' := \Pi^*\sI$. In particular, each center $Z_{i-1} \subset X_{i-1}$ is $G$-invariant, $G$ acts on each blowup $X_i$, and each blowup morphism $\pi_i:X_i\to X_{i-1}$ is $G$-equivariant, for $i=1,\ldots,r$.
\end{cor}

\begin{proof}
Applying Koll\'ar \cite[Proposition 3.9.1, p. 125]{Kollar_lectures_resolution_singularities} to Theorem \ref{thm:Monomialization_ideal_sheaf_algebraic_scheme} yields the conclusion.
\end{proof}

In exactly the same manner, one obtains

\begin{cor}[Equivariance of strong resolution of singularities for schemes]
\label{cor:Resolution_singularities_algebraic_scheme_equivariant}
(See Bierstone, Grigoriev, Milman, and W{\l}odarczyk \cite[Theorem 1.0.3, p. 196]{Bierstone_Grigoriev_Milman_Wlodarczyk_2011}.)  
Continue the hypotheses of Theorem \ref{thm:Resolution_singularities_algebraic_scheme}. Then the action of an algebraic group $G$ on the scheme $X$ lifts to an action of $G$ on its functorial resolution $X'$ and the morphism $\Pi:X' \to X$ is $G$-equivariant. In particular, each blowup center $Z_{i-1} \subset X_{i-1}$ is $G$-invariant, $G$ acts on each blowup $X_i$, and each blowup morphism $\pi_i:X_i\to X_{i-1}$ is $G$-equivariant, for $i=1,\ldots,r$.
\end{cor}

\begin{cor}[Equivariance of embedded resolution of singularities for schemes]
\label{cor:Embedded_resolution_of_singularities_algebraic_scheme_equivariant}
(See Bierstone, Grigoriev, Milman, and W{\l}odarczyk \cite[Theorem 1.0.2, p. 195]{Bierstone_Grigoriev_Milman_Wlodarczyk_2011} or Encinas and Villamayor \cite[Theorem 3.1, Item 2.d (iii), p. 156]{Encinas_Villamayor_2000}.)  
Continue the hypotheses of Theorem \ref{thm:Embedded_resolution_of_singularities_algebraic_scheme}. Then the action of an algebraic group $G$ on the scheme $X$ lifts to an action of $G$ on its functorial resolution $X'$ and the morphism $\Pi:(X',Y') \to (X,Y)$ is $G$-equivariant. In particular, each blowup center $Z_{i-1} \subset X_{i-1}$ is $G$-invariant, $G$ acts on each blowup $X_i$, and each blowup morphism $\pi_i:X_i\to X_{i-1}$ is $G$-equivariant, for $i=1,\ldots,r$.
\end{cor}

\chapter{Critical sets, stable manifolds, and unstable manifolds for Morse--Bott functions}
\label{chap:Stable_manifold_theorems_Morse_and_Morse-Bott_functions}
In Section \ref{sec:Critical_points_well-defined}, we give a preliminary definition of a critical point of a regular function on an analytic space. In Section \ref{sec:Critical_sets_stable_unstable_manifolds}, we specialize to the case of an analytic space that is a closed subspace of an analytic manifold and give a more refined definition of critical point and critical set in an analytic space. We also recall definitions of stable manifolds and unstable manifolds of a Morse--Bott function on a smooth manifold. Section \ref{sec:Stable_manifold_theorems_Morse_and_Morse-Bott_functions} contains statements of the stable manifold theorems for Morse and Morse--Bott functions. We conclude in Section \ref{sec:Morse-Bott_decomposition_smooth_manifold} with a review of the Morse--Bott decomposition of a smooth manifold.

\section{Definition of a critical point of a regular function}
\label{sec:Critical_points_well-defined}
In this section, we give a definition (see Lemma \ref{lem:Definition_critical_point_regular_function_well-defined}) of a critical point of a regular function on an analytic space and show that it is, in fact, well-defined. However, because the Hessian of a regular function on an analytic space at a critical point need not be well-defined, we shall subsequently describe a less general, but simpler concept (see Definition \ref{defn:Critical_point}) that serves as our standard definition of critical point throughout our work.

Let $\KK=\RR$ or $\CC$ and $(X,\sO_X)$ be a $\KK$-analytic model space as in Definition \ref{defn:Analytic_model_space}, so we are given a domain $D \subset \KK^n$, an ideal sheaf $\sI \subset \sO_D$ of finite type, and $X = \supp(\sO_D/\sI) \subset D$ with structure sheaf $\sO_X = (\sO_D/\sI)\restriction X$. Recall that a function $f:V\to \KK$ on an algebraic subset $V\subset\KK^n$ over a field $\KK$ is \emph{regular} if there is a polynomial function $\tilde f:\KK^n\to\KK$ such that $f=\tilde f$ on $V$ (see Shafarevich \cite[Section 1.2.2]{Shafarevich_v1} or Smith \cite[Definition 3.1]{Smith_notes_algebraic_geometry}). By analogy, we say that a function $f:U\to\KK$ defined on an open neighborhood $U \subset X$ of a point $p \in X$ is \emph{analytic} if there is an analytic function $\tilde f:\tilde U \to \KK$ on an open neighborhood $\tilde U \subset \KK^n$ of $p$ such that $U = X\cap \tilde U$ and $f = \tilde f$ on $U$. We may thus interpret the analytic function $f:X\to\KK$ as the restriction to $X \subset D$ of an analytic function $\tilde f:D\to\KK$.

We refer to
\begin{inparaenum}[\itshape a\upshape)]
\item Hartshorne \cite[Chapter I, Section 3, Definition, p. 14]{Hartshorne_algebraic_geometry} for the definition of a function that is regular at a point and regular on a quasi-affine variety $Y \subset \KK^n$
\item Hartshorne \cite[Chapter I, Section 3, Definition, p. 14]{Hartshorne_algebraic_geometry} for the definition of a function that is regular at a point and regular on a quasi-projective variety $Y \subset \PP(\KK)^n$;
\item Hartshorne \cite[Chapter I, Section 3, Definition, p. 15]{Hartshorne_algebraic_geometry} for the definition of the ring $\sO_Y$ of regular functions on a variety $Y$ over $\KK$;
\item Hartshorne \cite[Chapter II, Section 1, Example 1.0.1, p. 62]{Hartshorne_algebraic_geometry} for the definition of the sheaf $\sO_X$ of regular functions on a variety $X$ over $\KK$; and 
\item Eisenbud and Harris \cite[Section 1.2, p. 22]{Eisenbud_Harris_geometry_schemes} for the definition of a regular function on an open set $U\subset X$ and global regular function on a scheme $X$ over a field $\KK$.
\end{inparaenum}

We note that a global regular function, $f\in \sO_X(X)$, is the same as a morphism, $f:X\to\KK$, of algebraic varieties or schemes or analytic spaces over $\KK$ (see Greul and Pfister \cite[Section A.6]{Greuel_Pfister_singular_introduction_commutative_algebra}).
The forthcoming Definition \ref{defn:Critical_point} provides a simpler definition of critical point that applies when $X$ is a closed, $\KK$-analytic subspace of a $\KK$-analytic manifold, but the definition in the following lemma is more general and of some independent interest.

\begin{lem}[Critical points of regular functions are well-defined]
\label{lem:Definition_critical_point_regular_function_well-defined}  
Let $(X,\sO_X)$ be an analytic space over a field $\KK=\RR$ or $\CC$ as in Definition \ref{defn:Analytic_space}. If $f:X\to\KK$ is a regular function and $p \in X$ is a \emph{critical point} for $f$ in the sense that
\[
  T_pX \subseteq \Ker\left(\d{\tilde f}(p):T_pD \to \KK\right),
\]
for some open neighborhood $U\subset X$ of $p$ such that $(\varphi,\varphi^\sharp):(U,\sO_X\restriction U) \to (Y,\sO_Y)$ is an isomorphism of $\KK$-analytic spaces, where $D \subset \KK^n$ is a domain, $p \in X$ is identified with $\varphi(p) \in D$, and $\sI \subset \sO_D$ is a finitely generated ideal sheaf, $Y = \cosupp\sI$, and $\sO_Y = (\sO_D/\sI)\restriction Y$, and $\tilde f \in \sO_D(D)$ is such that $\tilde f = (\varphi^{-1})^\sharp(f) = f\circ\varphi^{-1} \in (\varphi^{-1})_*\sO_Y$, then the property of $p$ being a critical point of $f$ is independent of the choice of local lift $\tilde f$ of $(\varphi^{-1})^\sharp(f) \in (\varphi^{-1})_*\sO_Y$.
\end{lem}  

\begin{proof}
Let $\tilde g \in \sO_D(D)$ be another local lift of $(\varphi^{-1})^\sharp(f)$ such that $\tilde g = f\circ\varphi^{-1} \in (\varphi^{-1})_*\sO_Y$, where $\varphi^\sharp:\sO_Y \to \varphi_*(\sO_X\restriction U)$ and $(\varphi^{-1})^\sharp:\sO_X\restriction U \to (\varphi^{-1})_*\sO_Y$. Thus, $\tilde f$ and $\tilde g$ are two $\KK$-analytic functions on $D$ such that
\[
  \tilde g = \tilde f + \sum_{j=1}^l a_j f_j,
\]
where $\sI = f_1\sO_D+\cdots+f_l\sO_D$ and $a_j, f_j \in \sO_D(D)$ for $j=1,\ldots,l$. Because $p$ is a critical point of $f$ then, by the definition in the lemma,
\[
  \Ker \d{\tilde f}(p) \supseteq T_pX = \Ker \d f_1(p)\cap\cdots\cap \Ker \d f_l(p).
\]
We compute that
\[
  \d{\tilde g}(x) = \d{\tilde f}(x) + \sum_{j=1}^l a_j(x)\d f_j(x) + f_j(x)\d a_j(x),
  \quad\text{for all } x \in D.
\]
Because $f_j(x)=0$ for all $x\in Y$ and $j=1,\ldots,l$ and $p\in Y$ in particular, we see that
\[
  \d{\tilde g}(p) = \d{\tilde f}(p) + \sum_{j=1}^l a_j(p)\d f_j(p).
\]
Therefore,
\[
  \d{\tilde g}(p) = \d{\tilde f}(p) \quad\text{on } T_pX = \bigcap_{j=1}^l\Ker \d f_j(p)
\]
and so $\Ker \d{\tilde g}(p) \supseteq T_pX$ if and only if $\Ker \d{\tilde f}(p) \supseteq T_pX$. Hence, the definition of a critical point of $f$ is independent of the choice of representative analytic function $\tilde f$, for a given isomorphism $(\varphi,\varphi^\sharp):(U,\sO_X\restriction U) \to (Y,\sO_Y)$ of $\KK$-analytic spaces.
\end{proof}

\section{Critical sets, stable manifolds, and unstable manifolds}
\label{sec:Critical_sets_stable_unstable_manifolds}
Our development of Morse theory for analytic spaces uses the restriction of an analytic function $f$ on a analytic manifold $M$ to a possibly singular analytic subspace $X\subset M$, so we shall adopt a definition of `critical point' of the restriction of $f$ to $X$ that is more specific to this situation than that of Lemma \ref{lem:Definition_critical_point_regular_function_well-defined}. One approach to a definition is given by Goresky and MacPherson \cite[Section 2.1]{GorMacPh}. They assume given a Whitney stratification of a real analytic space $X$ and call $p\in S$ a critical point of $f_X := f|_X$ if $T_pS \subset \Ker(df(p):T_pM\to\RR)$, where $S\subset X$ is the smooth stratum of $X$ that contains $p$. Another definition is provided by Wilkin \cite[Section 1]{Wilkin_2019arxiv} and is described further below, while Massey \cite{Massey_2000} provides a more topological definition. We find it convenient to instead adopt the following definition, as this fits most naturally with our application.

\begin{defn}[Critical point of the restriction of an analytic function on an analytic manifold to an analytic subspace]
\label{defn:Critical_point}  
Let $M$ be an analytic manifold over $\KK=\RR$ or $\CC$, and $f:M\to\KK$ be an analytic function, and $X \subset M$ be a closed analytic subspace as in Definition \ref{defn:Analytic_subspace}, and let $f_X:X\to\KK$ be given by $f_X = f \restriction X$. Then $p \in X$ is a \emph{critical point} for $f_X:X\to\KK$ if $T_pX \subseteq \Ker(df(p):T_pM\to\KK)$, where $T_pX \subseteq T_pM$ is the \emph{Zariski tangent space} for $X$ at $p$.
\end{defn}  

It is useful to recall some definitions and terminology and results from dynamical systems that is relevant to our approach to Morse theory on analytic spaces. We rely on the discussions in Wilkin \cite[Section 1]{Wilkin_2019arxiv} and Palis and de Melo \cite[Section 2.6]{Palis_deMelo_geometric_theory_dynamical_systems}.

Let $(M,g)$ be a real analytic Riemannian manifold with the metric topology and $X \subset M$ be a closed real analytic subspace. Given an analytic function $f: M \to \RR$, the gradient vector field $\grad f$ is well-defined on $M$. Suppose that for each point $x \in X$ there exists an open interval $(-\eps,\eps)$ (depending on $x$) such that the flow $\varphi(t,x)$ of $\grad_g f$ with initial condition $x$ exists,
\begin{equation}
  \label{eq:Gradient_flow_equation}
  \frac{d\varphi}{dt}(t,x) = -\grad_g f(\varphi(t,x)), \quad t \in (-\eps,\eps),
\end{equation}
depends continuously on the initial condition $x$, and $\varphi(t,x) \in X$ for all $t\in(-\eps,\eps)$. (These conditions are satisfied when the flow is generated by a group action which preserves the subspace $X$ and moment map flows form an important class of such examples.)

Wilkin and Palis and de Melo define a critical point of $f_X = f\restriction X$ to be a fixed point of this flow on $X$. Let $\Crit f_X \subset X$ denote the subset of all critical points of $f_X$. Given a critical value $c \in \RR$ and a connected component $X_C \subseteq \Crit f_X \cap f_X^{-1}(c)$ of the critical set, define $X_C^+$, $X_C^-$ to be the corresponding \emph{stable} (or \emph{ascending}) and \emph{unstable}  (or \emph{descending}) \emph{sets} with respect to the flow by
\begin{subequations}
\label{eq:Morse-Bott_function_stable_unstable_set}  
\begin{align}
  \label{eq:Morse-Bott_function_stable_set}  
  X_C^+ &:= \left \{ x \in X: \lim_{t\to\infty} \varphi(t,x) \in X_C \right\},
  \\
  \label{eq:Morse-Bott_function_unstable_set}  
  X_C^- &:= \left \{ x \in X: \lim_{t\to-\infty} \varphi(t,x) \in X_C \right\}.
\end{align}
\end{subequations}
We write $X_p^\pm$ when the set $X_C$ is replaced by a point $p\in X_C$ in these definitions. See Wilkin \cite[Theorem 1.1 and Corollary 1.2]{Wilkin_2019arxiv} for general criteria when the desirable \cite[Conditions 1--4]{Wilkin_2019arxiv} hold, together with Kirwan \cite{Kirwan_cohomology_quotients_symplectic_algebraic_geometry} and Wilkin \cite{Wilkin_2019} for more specific examples. Specifically, Wilkin \cite[Corollary 1.2]{Wilkin_2019arxiv} requires only that $f:X\to\RR$ be \emph{proper} (that is, the preimage under $f$ of compact sets are compact), in addition to being real analytic.

\begin{figure}
	\centering
	\includegraphics[width=0.7\linewidth]{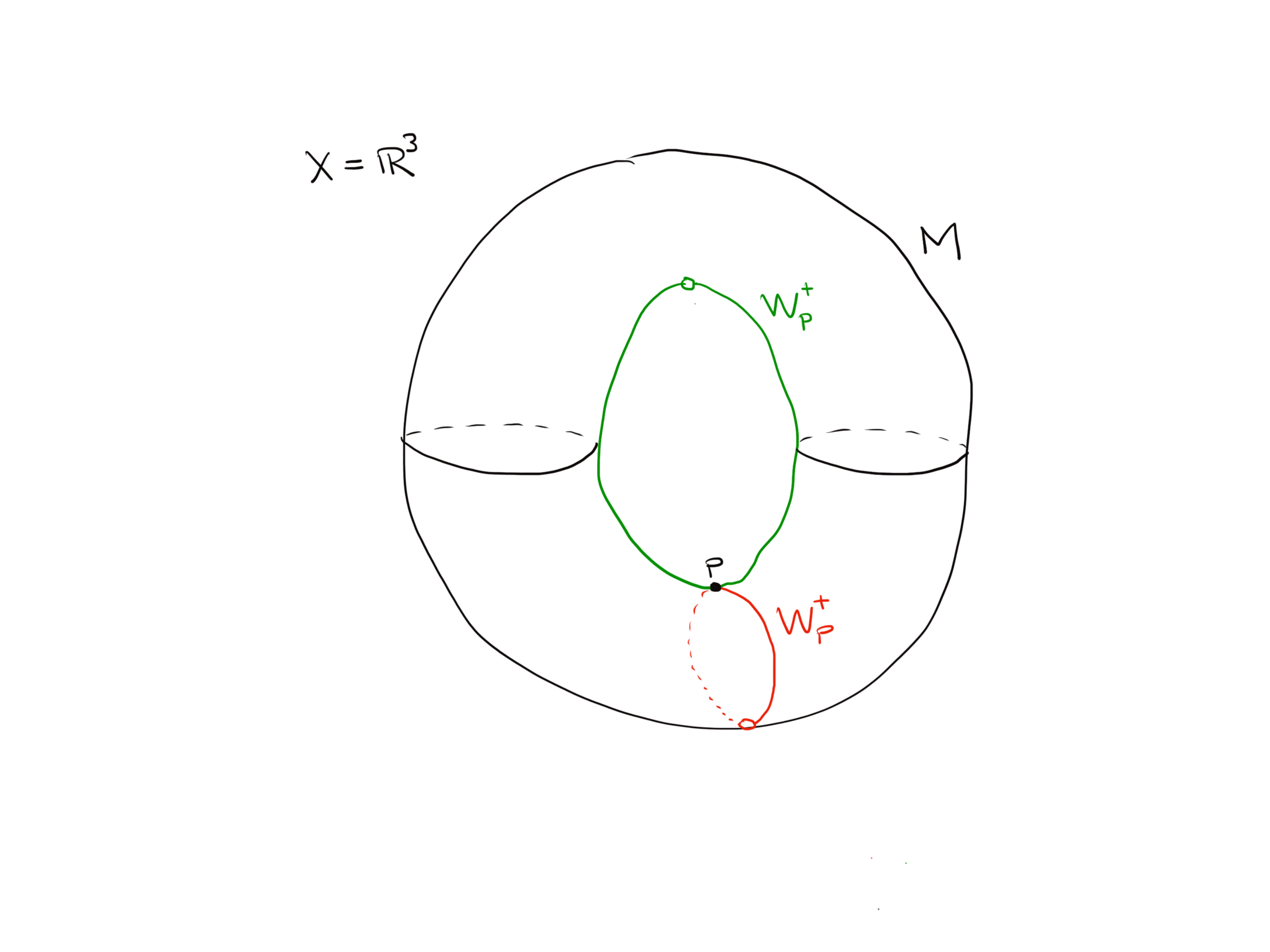}
	\caption[Stable and unstable manifolds for a saddle point in the torus]{Stable and unstable manifolds $X_p^\pm$ for the height function around a saddle point $p$ in the torus $X$ in $M=\RR^3$}
	\label{fig:stable_unstable_manifolds_torus_saddle_point}
\end{figure}

Figure \ref{fig:stable_unstable_manifolds_torus_saddle_point} illustrates the stable and unstable manifolds $X_p^\pm$ for the height function $f_X$ around a saddle point $p$ in the torus $X$ in $M=\RR^3$, where $X_p^\pm$ are both given by copies of $S^1\less\{p\}$, that is, one-dimensional open disks.

\begin{figure}
	\centering
	\includegraphics[width=0.7\linewidth]{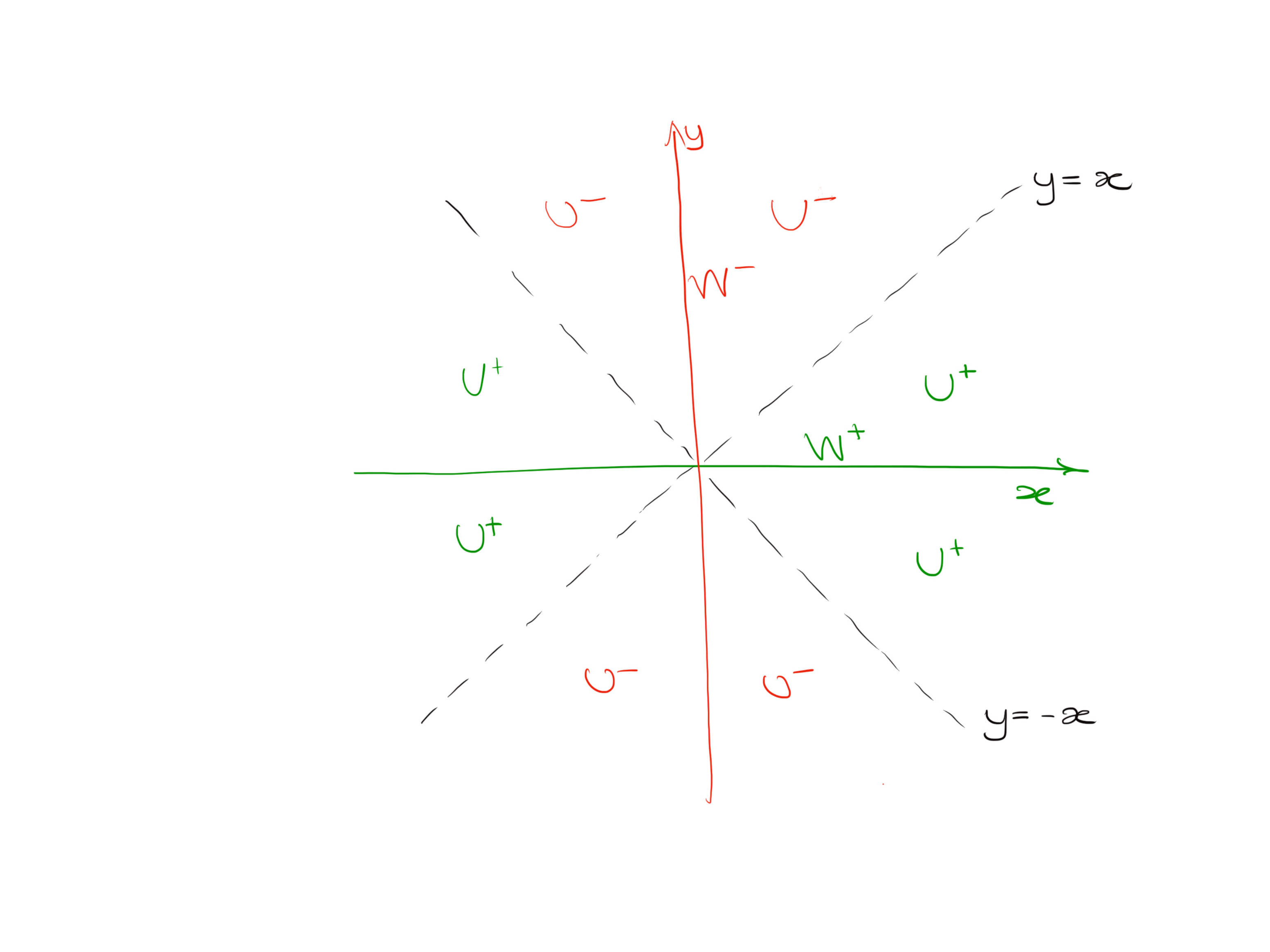}
	\caption[Open neighborhoods of the stable and unstable manifolds in the plane]{Open neighborhoods $U^\pm$ in $\RR^2$ of the stable and unstable manifolds $X_p^\pm$ for the gradient flow of $f(x,y)=\frac{1}{2}(x^2-y^2)$ near the origin}
	\label{fig:open_neighborhoods_unstable_manifolds}
\end{figure}

\begin{exmp}[Flow lines near a saddle point for the height function in $\RR^3$]
\label{exmp:Gradient_flow_near_saddle_point_height_function}  
Suppose $M=\RR^3$ and that $X\subset M$ is the graph of $z=f(x,y)=\frac{1}{2}(x^2-y^2)$ for $(x,y)\in\RR^2$ and consider the critical point $p=(0,0)$. We have $f(x,y)<0$ for all points $(x,y)$ with $|x|<|y|$ and $f(x,y)>0$ for all points $(x,y)$ with $|x|>|y|$. Noting that
\[
  \Hess f(0,0) = \begin{pmatrix} 1 & 0 \\ 0 & -1 \end{pmatrix},
\]
we see that the $x$-axis (the stable manifold $X_p^+$ for the (negative) gradient flow of $f$ near the origin,
\[
  (\dot x(t),\dot y(t)) = -f'(x(t),y(t)), \quad\text{for } t \in \RR,
\]
and the line through eigenvector $(1,0)$ with eigenvalue $+1$) has an open neighborhood
\[
  U_+ := \{(x,y)\in\RR^2:|y|<|x|\} = \{(r,\theta)\in\RR^2: r\in [0,\infty) \text{ and } |\theta|<\pi/4\}
\]
such that flow lines from $(0,0)$ into $U_+$ \emph{increase} in height, $f$. Indeed, for any direction $u=(\cos\theta,\sin\theta)$ with $|\theta|<\pi/4$ and thus $|\cos\theta|<1/\sqrt{2}$, then the point $(x(t),y(t)) = (t\cos\theta,t\sin\theta)$ for $t \in \RR$ has height
\[
  f(x(t),y(t)) = \frac{t^2}{2}(\cos^2\theta-\sin^2\theta) = t^2\left(\cos^2\theta-\frac{1}{2}\right) > 0
\]
that \emph{increases from zero} for $t$ increasing or decreasing from zero. Similarly, the $y$-axis (the unstable manifold $X_p^-$ for the gradient flow of $f$ near the origin and the line through eigenvector $(0,1)$ with eigenvalue $-1$) has an open neighborhood
\[
  U_- := \{(x,y)\in\RR^2:|y|>|x|\} = \{(r,\theta)\in\RR^2: r\in [0,\infty) \text{ and } \pi/4<|\theta|<3\pi/4\}
\]
such that flow lines from $(0,0)$ into $U_-$ \emph{decrease} in height, $f$. Indeed, for any direction $u=(\cos\theta,\sin\theta)$ with $\pi/4<|\theta|<3\pi/4$ and thus $|\cos\theta|>1/\sqrt{2}$, then the point $(x(t),y(t)) = (t\cos\theta,t\sin\theta)$ for $t \in \RR$ has height
\[
  f(x(t),y(t)) = t^2\left(\cos^2\theta-\frac{1}{2}\right) < 0
\]
that \emph{decreases from zero} for $t$ increasing or decreasing from zero. Flow lines from $(0,0)$ along the lines $y=\pm x$ (equivalently, $\theta=\pm\pi/4$ or $\pm 3\pi/4$) remain at height zero. See Figure \ref{fig:open_neighborhoods_unstable_manifolds} for illustrations of the neighborhoods $U^\pm$ in $\RR^2$.
\qed
\end{exmp}

\begin{figure}
	\centering
	\includegraphics[width=0.7\linewidth]{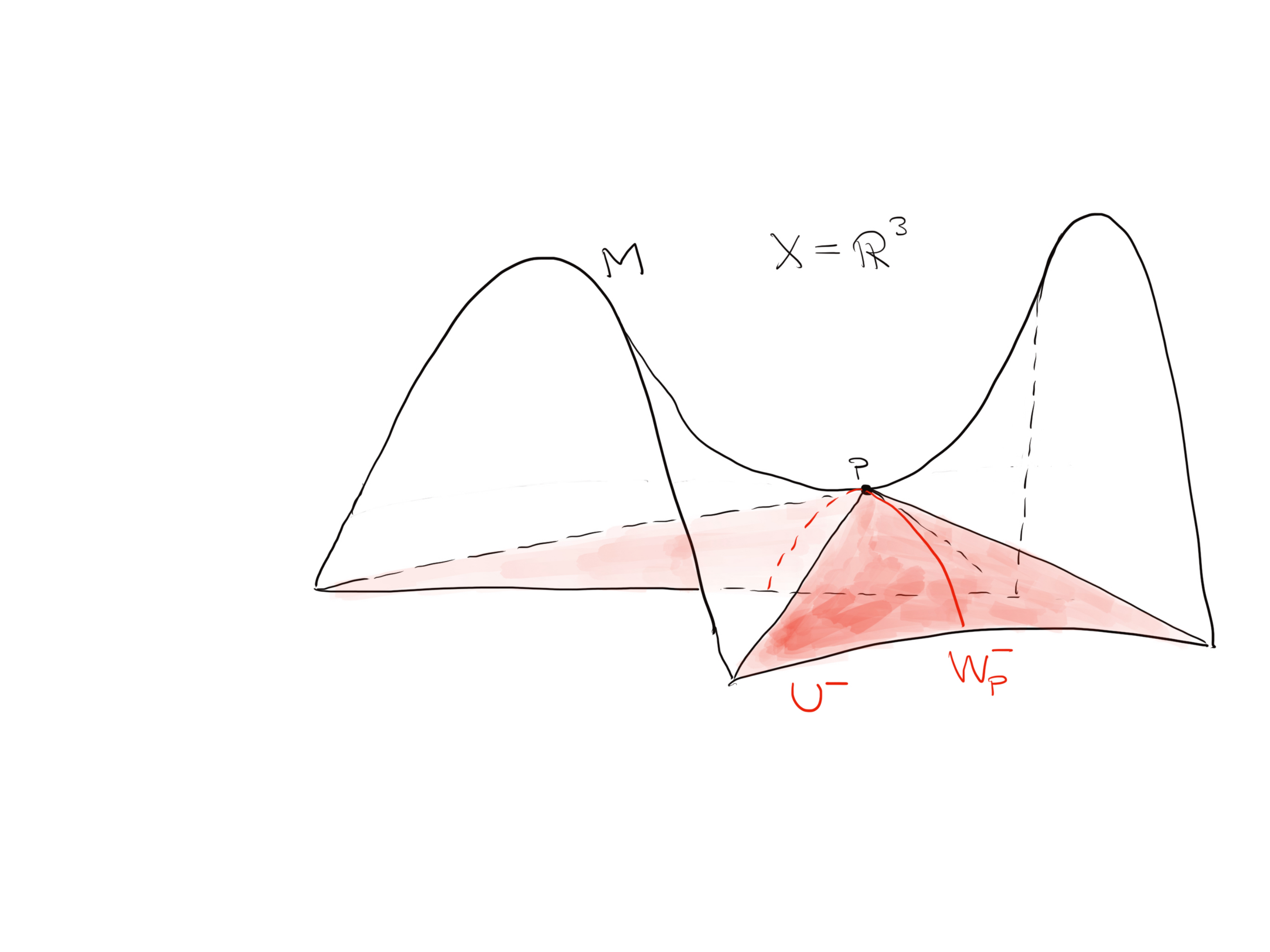}
	\caption[Open neighborhoods of the unstable manifold in a torus]{Open neighborhood $U^-$ of the unstable manifold $X_p^-$ for the gradient flow of the height function near a saddle point $p$ in the torus $X$ in $M=\RR^3$}
	\label{fig:open_neighborhood_unstable_manifold_saddle_point}
\end{figure}

See Figure \ref{fig:open_neighborhood_unstable_manifold_saddle_point} for an illustration of the corresponding neighborhood $U^-$ of the unstable manifold $X_p^-$ for the gradient flow of the height function near a saddle point $p$ in the torus $X$ in $M=\RR^3$. Example \ref{exmp:Gradient_flow_near_saddle_point_height_function} generalizes to give a description of the open neighborhoods of the stable and unstable manifolds near a critical point. 

\section{Stable manifold theorems for Morse and Morse--Bott functions}
\label{sec:Stable_manifold_theorems_Morse_and_Morse-Bott_functions}
We begin with the key

\begin{thm}[Stable and unstable manifolds for the gradient flow of a Morse function near a critical point]
\label{thm:Stable_manifold_Morse_function}
(See Audin and Damian \cite[Proposition 2.1.5, p. 28]{Audin_Damian_morse_theory_floer_homology}, Banyaga and Hurtubise \cite[Theorem 4.2]{Banyaga_Hurtubise_lectures_morse_homology}, Cohen, Iga, and Norbury \cite[Theorem 2 in Section 6.3]{Cohen_Iga_Norbury_morse_theory}, and Nicolaescu \cite[Proposition 2.24]{Nicolaescu_morse_theory}.)
Let $(M,g)$ be a closed, finite-dimensional, smooth, Riemannian manifold of dimension $m \geq 2$ and $f:M\to\RR$ be a smooth Morse function with critical point $p\in M$ with Morse index $\lambda_p^-$ and Morse co-index $\lambda_p^+ := m-\lambda_p^-$. Then the tangent space $T_pM$ to $M$ at $p$ splits as an orthogonal direct sum,
\[
  T_pM = T_p^+M \oplus T_p^-M,
\]
where the Hessian $\Hess f(p) \in \End(T_pM)$ is positive definite on $T_p^+M$ and negative definite on $T_p^-M$. Moreover, the stable and unstable manifolds $M_p^\pm$ are surjective images of smooth embeddings,
\[
  \iota^\pm: T_p^\pm M \to M_p^\pm \subset M,
\]
and the manifolds $M_p^\pm$ are smoothly embedded open balls of dimension $\lambda_p^\pm$.
\end{thm}

Extensions of Theorem \ref{thm:Stable_manifold_Morse_function} from Morse functions to Morse--Bott functions are described by Cohen, Iga, and Norbury \cite[Theorem 15.3]{Cohen_Iga_Norbury_morse_theory}, Zhou \cite[Proposition 8.4 and Lemma 8.6]{Zhou_zhengyi_2019arxiv_v2},
and Austin and Braam \cite[Proposition 3.2 and Theorem A.9]{Austin_Braam_1995} who provide a proof based on Hirsch, Pugh, and Shub \cite[Theorems 1.1 and 4.1]{Hirsch_Pugh_Shub_invariant_manifolds}, one of the most general statements of the stable manifold theorem for dynamical systems. (Palis and de Melo \cite[Theorem 2.6.2]{Palis_deMelo_geometric_theory_dynamical_systems} provide a statement and proof of a simpler version of the stable manifold theorem for dynamical systems.) We refer the reader to the articles by Naito \cite{Naito_1988}, \cite{Naito_Kozono_Maeda_1990}, \cite{Naito_Kozono_Maeda_1989} and his collaborators for analogues of the stable manifold theorem for gradient flow on Banach manifolds, especially Yang--Mills gradient flow.

\begin{thm}[Morse--Bott Lemma for functions on Euclidean space]
\label{thm:Morse-Bott_Lemma}
(See Feehan \cite[Theorem 2.10 and Remark 2.12]{Feehan_lojasiewicz_inequality_all_dimensions_morse-bott}.)
Let $U \subset \RR^m$ be an open neighborhood of the origin and $f:U\to\RR$ be a $C^{p+2}$ function ($p\geq 1$) such that $f(0) = 0$. If $f$ is Morse--Bott at the origin as in Definition \ref{maindefn:Morse-Bott_function}, then, after possibly shrinking $U$, there are an open neighborhood of the origin $V \subset \RR^m$ and a $C^p$ diffeomorphism, $V \ni y \mapsto x = \Phi(y) \in U$ with $\Phi(0)=0$ and $D\Phi(0)=\id_{\RR^m}$, such that
\begin{equation}
\label{eq:Morse-Bott_function}
f(\Phi(y)) = \frac{1}{2}\langle y, Ay\rangle_{\RR^m}, \quad\text{for all } y \in V,
\end{equation}
where
\[
  A := f''(0) = (f\circ\Phi)''(0) \in \End(\RR^m)
\]
is symmetric. If $f$ is analytic, then $\Phi$ is analytic. There is an element $P \in \SO(m)$ such that $P^\intercal A P = \Lambda$, where $\Lambda$ is diagonal with $\lambda^+$ positive and $\lambda^-$ negative eigenvalues and the eigenvalue zero with multiplicity $c=m-\lambda^+-\lambda^-$ and, after replacing the coordinates $y$ by $y=P^{-1}\tilde y$ on $\RR^m$ and relabeling,
\[
  f(\Phi(y)) = \frac{1}{2}\sum_{i=1}^{\lambda^+} \nu_i y_i^2 - \frac{1}{2}\sum_{i=\lambda^+ + 1}^{\lambda^+ + \lambda^-} \nu_i y_i^2,
  \quad\text{for all } y \in V,
\]
where $\nu_i>0$ for $i=1,\ldots,m-c$.
\end{thm}

More generally, we have the

\begin{thm}[Stable and unstable manifolds for the gradient flow of a Morse--Bott function near a critical submanifold]
\label{thm:Stable_manifold_Morse-Bott_function}
(See Cohen, Iga, and Norbury \cite[Theorem 15.3]{Cohen_Iga_Norbury_morse_theory} and Austin and Braam \cite[Proposition 3.2 and Theorem A.9]{Austin_Braam_1995}.)
Let $(M,g)$ be a closed, finite-dimensional, smooth, Riemannian manifold of dimension $m\geq 2$ and $f:M\to\RR$ be a smooth Morse--Bott function with connected, smooth critical submanifold $C\subset M$ with Morse--Bott index $\lambda_C^-$ and co-index $\lambda_C^+$. Then the maps
\[
  \pi^\pm: M_C^\pm \to C
\]
are smooth fiber bundles with fibers diffeomorphic to open balls $B^{\lambda^\pm} \subset \RR^{\lambda^\pm}$. Moreover as bundles they are isomorphic to the negative and positive normal bundles of $C$ in $M$,
\[
  p^\pm: N_{C/M}^\pm \to C.
\]
\end{thm}

Theorem \ref{thm:Morse-Bott_Lemma} implies that, after a change of coordinates, the stable and unstable manifolds of gradient flow for $f$ and its critical submanifold on a small enough open neighborhood of the origin $p\in \RR^m$ take the simple forms indicated in Example \ref{exmp:Gradient_flow_near_saddle_point_height_function} as coordinate subspaces of $\RR^m=\RR^{\lambda^+}\times\RR^{\lambda^-}\times\RR^c$:
\begin{equation}
\label{eq:(Un)stable_and_critical_submanifolds_Morse-Bott_function}
\begin{aligned}
  M_p^+ &= \{y\in U: y_i=0, \text{for } i=1,\ldots,\lambda^+\} = \RR^{\lambda^+}\times\{0\}\times\{0\},
  \\
  M_p^- &= \{y\in U: y_i=0, \text{for } i=\lambda^++1,\ldots,\lambda^+ + \lambda^-\} = \{0\}\times\RR^{\lambda^-}\times\{0\},
  \\
  C = M^0 &= \{y\in U: y_i=0, \text{for } i=m-c+1,\ldots,m\} = \{0\}\times\{0\}\times\RR^c.
\end{aligned}
\end{equation}
We thus have the

\begin{lem}[Open cone neighborhoods of stable and unstable submanifolds]
\label{lem:Cone_neighborhoods_(un)stable_submanifolds}
Continue the hypotheses of Theorem \ref{thm:Morse-Bott_Lemma} and notation in \eqref{eq:(Un)stable_and_critical_submanifolds_Morse-Bott_function}, let $\{e_i\}_{i=1}^m$ denote the standard basis of unit vectors for $\RR^m$, and relabel $f\circ\Phi$ as $f$. If
\begin{align*}
  U^+ &:= \left\{y\in\RR^m: |\angle\, y, M_p^+| < \pi/4\right\},
  \\
  U^- &:= \left\{y\in\RR^m: |\angle\, y, M_p^-| < \pi/4\right\},
\end{align*}
then $f(y)<f(0)$ for all $y\in U^-\less\{0\}$ and $f(y)>f(0)$ for all $y\in U^+\less\{0\}$ and \emph{second-order gradient flow}, $\ddot y(t) = f''(y(t)f'(y(t))$, for $f$ emanating from the origin in any direction in $U^-$ or $U^+$ is strictly decreasing or increasing, respectively.
\end{lem}

In Lemma \ref{lem:Cone_neighborhoods_(un)stable_submanifolds}, we define the \emph{second-order gradient flow} for $f$ starting at the origin in a direction $v\in U^\pm$ by taking the time derivative of the \emph{first-order (negative) gradient flow} (or \emph{method of steepest descent}),
\begin{equation}
\label{eq:Gradient_flow_Euclidean_space}  
\dot y(t) = -f'(y(t)), \quad\text{for } t \in [0,\eps), \quad y(0)=y_0
\end{equation}
starting at a point $y_0\in U$ to give $\ddot y(t) = -f''(y(t))\dot y(t)  = f''(y(t))f'(y(t)$, that is, a solution to the second-order system,
\begin{equation}
\label{eq:Time_derivative_gradient_flow_Euclidean_space} 
  \ddot y(t) = f''(y(t))f'(y(t)) \quad\text{for } t \in [0,\eps),
  \quad y(0)=y_0, \quad \dot y(0)=v,
\end{equation}
for some initial point $y_0\in U$ and initial direction $v$. The angle $\theta = \angle\, v, W$ in Lemma \ref{lem:Cone_neighborhoods_(un)stable_submanifolds} between a unit vector $v\in\RR^m$ and a subspace $W\subset\RR^m$ is defined by $\cos\theta := \langle v, \pi_W v \rangle_{\RR^m}$, where $\pi_W v$ is the orthogonal projection of $v$ onto $W$. Note that if we define
\[
  F(y) := \frac{1}{2}|f'(y)|^2, \quad y \in U,
\]
then $F'(y) = f''(y)f'(y)$ and thus $y(t)$ in \eqref{eq:Time_derivative_gradient_flow_Euclidean_space} may be viewed as a solution to the \emph{second-order gradient system},
\begin{equation}
\label{eq:Second-order_gradient_flow_Euclidean_space}
  \ddot y(t) = F'(y(t)) \quad\text{for } t \in [0,\eps),
  \quad y(0)=0, \quad \dot y(0)=v.
\end{equation}
See, for example, Boulmezaoud, Cieutat, and Daniilidis \cite{Boulmezaoud_Cieutat_Daniilidis_2018} for a comparison of solutions to \eqref{eq:Gradient_flow_Euclidean_space} and \eqref{eq:Second-order_gradient_flow_Euclidean_space}.

We illustrate an application of Lemma \ref{lem:Cone_neighborhoods_(un)stable_submanifolds} in the following

\begin{exmp}[Flow lines near a saddle point for the height function in $\RR^3$ (continued)]
\label{exmp:Second-order_gradient_flow_near_saddle_point_height_function} 
For instance, in Example \ref{exmp:Gradient_flow_near_saddle_point_height_function} we have $\grad f(x,y) = (x,-y)$ and we consider
\begin{equation}
\label{eq:Second_order_gradient_flow_saddle}  
  \ddot x(t) = -\dot x(t), \quad \dot x(0) = \cos\theta, x(0) = 0
  \quad\text{and}\quad
  \ddot y(t) = \dot y(t), \quad \dot y(0) = \sin\theta, y(0) = 0.
\end{equation}
This integrates to give
\[
  (\dot x(t), \dot y(t)) = (e^{-t}\cos\theta, e^t\sin\theta), \quad\text{for } t \in \RR
\]
with $(x(0),y(0)) = (0,0)$, and thus
\begin{align*}
  (x(t),y(t))
  &= (-e^{-t}\cos\theta, e^t\sin\theta) - (-\cos\theta, \sin\theta)
  \\
  &= (-(e^{-t}-1)\cos\theta, (e^t-1)\sin\theta), \quad\text{for } t \in \RR.
\end{align*}
Hence,
\[
  f(x(t),y(t)) = \frac{1}{2}\left((e^{-t}-1)^2\cos^2\theta - (e^t-1)^2\sin^2\theta\right),
  \quad\text{for } t \in \RR.
\]
Using $e^t = 1 + t + \frac{1}{2}t^2 + \cdots$, we see that for $t \in (-\eps,\eps)$ and $\eps\in (0,1]$ small, we have
\[
  (x(t),y(t)) \approx t(\cos\theta, \sin\theta), \quad\text{for } t \in (-\eps,\eps)
\]
and
\[
  f(x(t),y(t)) \approx \frac{t^2}{2}\left(\cos^2\theta - \sin^2\theta\right),
  \quad\text{for } t \in (-\eps,\eps),
\]
and so the flow lines for \eqref{eq:Second_order_gradient_flow_saddle} behave as claimed near the origin in Example \ref{exmp:Gradient_flow_near_saddle_point_height_function} when \emph{started at the origin}. Note that our solution $(x(t),y(t))$ to \eqref{eq:Second_order_gradient_flow_saddle} is not literally a solution to negative gradient flow starting at the origin, $(\dot x(t), \dot y(t)) = -(x(t),-y(t))$ with $(x(0),y(0)) = (0,0)$, but rather obeys
\[
  (\dot x(t), \dot y(t)) = (-x(t),y(t)) + (\cos\theta, \sin\theta), \quad\text{for } t \in \RR,
  \quad (x(0),y(0)) = (0,0),
\]
as one can see by combining the preceding equations.
\qed
\end{exmp}

\section{Morse--Bott decomposition of a smooth manifold}
\label{sec:Morse-Bott_decomposition_smooth_manifold}
In this section, we review the decomposition of a smooth manifold into \emph{stable} or \emph{unstable} subsets defined by a Morse--Bott function, broadly following Bott \cite[p. 341]{Bott_1988}, \cite[p. 104]{Bott_1988} and Smale \cite{Smale_1960}, \cite[Theorem 2.3, p. 753]{Smale_1967}, both of whom pay tribute to ideas of Thom \cite{Thom_1949}.

Let $(M,g)$ be a finite-dimensional, real analytic Riemannian manifold and $f:M\to\RR$ be a real analytic Morse--Bott function in the sense of Definition \ref{maindefn:Morse-Bott_function}. Unlike the cited references, we allow the critical set $M^0 := \Crit f \subset M$ to comprise a union of critical submanifolds that may be positive-dimensional. We let $\{M_\alpha^0\}$ denote the connected components of $M^0$ and note that each $M_\alpha^0$ is an embedded, real analytic submanifold of $M$ by Definition \ref{maindefn:Morse-Bott_function} and our assumption that $f$ is real analytic and Morse--Bott. An assumption of real analyticity of $f$ does not appear in the cited papers by Bott or Smale but, we include it here in order to ensure convergence of its gradient flow (in the presence of additional assumptions of compactness, such as an assumption that $M$ is closed) by applying results originating in the work of {\L}ojasiewicz \cite{Lojasiewicz_1965, Lojasiewicz_1984} and Simon \cite{Simon_1983}. Following \eqref{eq:Morse-Bott_function_stable_unstable_set}, we denote
\begin{align*}  
  M_\alpha^+ &:= \left \{ x \in M: \lim_{t\to\infty} \varphi(t,x) \in M_\alpha^0 \right\},
  \\
  M_\alpha^- &:= \left \{ x \in M: \lim_{t\to-\infty} \varphi(t,x) \in M_\alpha^0 \right\},
\end{align*}
where $\varphi:\RR\times M \to M$ is the (negative) gradient flow \eqref{eq:Gradient_flow_equation} defined by $f$. If $M_\alpha^0$ consists of an isolated point $p_\alpha\in M$, then Theorem \ref{thm:Stable_manifold_Morse_function} ensures that the subsets $M_\alpha^\pm\cap U$ are (smoothly) embedded submanifolds of an open neighborhood of $U$ of $p_\alpha$. More generally, if $M_\alpha^0$ is a positive-dimensional submanifold, then Theorem \ref{thm:Stable_manifold_Morse-Bott_function} ensures that the subsets $M_\alpha^\pm \subset M$ are (smoothly) embedded submanifolds and that the associated fiber bundles $\pi_\alpha^\pm: M_\alpha^\pm \to M_\alpha^0$ are isomorphic to the positive and negative normal bundles $N_{M_\alpha^0/M}^\pm \to M_\alpha^0$ of $M_\alpha^0$ in $M$, so that
\begin{equation}
  \label{eq:Morse-Bott_splitting_tangent_bundlle_along_component_critical_submanifold}
  TM\restriction M_\alpha^0 = TM_\alpha^0 \oplus N_{M_\alpha^0/M}^+ \oplus N_{M_\alpha^0/M}^-
\end{equation}
and $N_{M_\alpha^0/M}^\pm \cong N_{M_\alpha^0/M_\alpha^\pm}$. In the presence of a suitable assumption of compactness (for example, that $M$ is closed), then the limits $\lim_{t\to\infty}\varphi(t,x)$ and $\lim_{t\to -\infty}\varphi(t,x)$ exist and belong to the critical subset $M^0\subset M$ for every point $x \in M$. Consequently, when $M^0$ comprises a set of isolated points, one obtains the two stratifications
\begin{equation}
  \label{eq:Morse-Bott_decomposition_smooth_manifold}
  M = \bigsqcup_\alpha M_\alpha^+ = \bigsqcup_\alpha M_\alpha^-
\end{equation}
indexed by the critical points of $f$ (see Bott \cite[p. 341]{Bott_1988}, \cite[p. 104]{Bott_1988} and Smale \cite[Theorem 2.3, p. 753]{Smale_1967}). More generally, these stratifications are indexed by the connected components of $M^0$, some or all of which may be positive-dimensional.

\chapter{Group actions on schemes and analytic spaces}
\label{chap:Group_actions_analytic_spaces_algebraic_schemes}
In this chapter, we review results for group actions on schemes and analytic spaces and establish some extensions that will become important in our review of Bia{\l}ynicki--Birula decompositions for algebraic varieties and complex analytic spaces and our development of certain extensions of those results. We begin in Section \ref{sec:Algebraic_group_actions_algebraic_schemes} by reviewing the concepts of algebraic groups and their actions on algebraic schemes. Section \ref{sec:Analytic_group_actions_analytic_spaces} contains the analogous discussion for analytic groups and their actions on analytic spaces. In Section \ref{sec:Linearization_theorems_compact_group_actions_manifolds_complex_analytic_spaces}, we review results on linearizations for actions of compact groups on manifolds and complex analytic spaces and establish extensions that we shall apply when discussing induced group actions on blowups. Section \ref{sec:Linearization_theorems_noncompact_group_actions_manifolds} provides an exposition of Weyl's unitary trick and the corresponding linearizations for actions of noncompact groups on complex manifolds. In Section \ref{sec:Holomorphic_maps_vector_spaces_with_circle_actions} we discuss holomorphic maps of complex vector spaces with circle actions and the circle actions that they induce on the codomains of the given holomorphic maps. We conclude in Section \ref{sec:Real_analytic_maps_vector_spaces_with_circle_actions} by discussing the considerably weaker results available for real analytic  maps of complex  vector spaces with circle actions.

\section{Algebraic groups and their actions on algebraic schemes}
\label{sec:Algebraic_group_actions_algebraic_schemes}
We follow Milne \cite[Sections 1a and 1f]{Milne_algebraic_groups_cup} for the definition of an algebraic group and its action on an algebraic scheme. For more informal introductions to group actions on analytic varieties or schemes, see Brion \cite{Brion_introduction_actions_algebraic_groups}, Dr\'ezet \cite{Drezet_2000}, Prasad \cite{Prasad_lectures_algebraic_groups}, and Serganova \cite{Serganova_lie_groups}.

Before recalling Milne's definitions, we review his terminology. An \emph{affine algebraic scheme} over $\KK$ is a $\KK$-ringed space isomorphic to $\Spm(A)$ for some $\KK$-algebra $A$; a \emph{morphism} (or \emph{regular map}) of affine algebraic schemes over $\KK$ is a morphism of $\KK$-ringed spaces --- it is automatically a morphism of locally ringed spaces (see Milne \cite[Appendix A.7, p. 568]{Milne_algebraic_groups_cup}). For $A = \KK[x_1,\ldots,x_n]$, we recall that Milne denotes the topological space $V$ of maximal ideals in $A$ by $\spm(A)$ and denotes the $\KK$-ringed space $(V,\sO_V)$ by $\Spm(A)$ (see Milne \cite[Section 3.e, p. 63]{Milne_algebraic_geometry}).

An \emph{algebraic scheme} over a field $\KK$ (or \emph{algebraic $\KK$-scheme}) is a scheme of finite type over $\KK$; an algebraic scheme is an \emph{algebraic variety} if it is geometrically reduced and separated (see Milne \cite[p. 3]{Milne_algebraic_groups_cup}). Let $(X,\sO_X)$ be a $\KK$-ringed space. An open subset $U$ of $X$ is said to be \emph{affine} if $(U,\sO_X)$ is an affine algebraic scheme over $\KK$; an algebraic scheme over $\KK$ is a $\KK$-ringed space $(X,\sO_X)$ that admits a finite covering by open affines; a \emph{morphism of algebraic schemes} (also called a \emph{regular map}) over $\KK$ is a morphism of $\KK$-ringed spaces (see Milne \cite[Appendix A.11, p. 569]{Milne_algebraic_groups_cup}).

\begin{defn}[Algebraic group]
\label{defn:Algebraic_group}  
(See Milne \cite[Section 1a, Definition 1.1, p. 6]{Milne_algebraic_groups_cup} and
Conrad \cite[Definitions 1.1.1 and 1.1.5]{Conrad_linear_algebraic_groups}.)  
Let $\KK$ be field and denote $\star = \Spm(\KK)$. Let $G$ be an algebraic scheme over $\KK$ and let $m:G\times G\to G$ be a regular map. The pair $(G,m)$ is an \emph{algebraic group over $\KK$} if there exist regular maps
\[
  e:\star \to G \quad\text{and}\quad \inv: G \to G,
\]
such that the following diagrams commute:
\begin{equation}
  \label{eq:Algebraic_group_associativity}
  \begin{tikzcd}
  G\times G\times G \arrow[r, "\id\times m"] \arrow[d, "m\times\id"] &G\times G \arrow[d, "m"]
  \\
  G\times G \arrow[r, "m"] &G
  \end{tikzcd}
\end{equation}
and
\begin{equation}
  \label{eq:Algebraic_group_identity}
  \begin{tikzcd}
    \star\times G \arrow[r, "e\times \id"] \arrow[rd, "\cong"'] &G\times G \arrow[d, "m"]
    &G\times \star \arrow[l, "\id\times e"'] \arrow[ld, "\cong"]
    \\
    &G
  \end{tikzcd}
\end{equation}
and
\begin{equation}
  \label{eq:Algebraic_group_inverse}
  \begin{tikzcd}
  G \arrow[d] \arrow[r, "(\inv{,}\id)"] &G\times G \arrow[d, "m"] &G \arrow[l, "(\id{,}\inv)"'] \arrow[d]
  \\
  \star \arrow[r, "e"] &G  &\star \arrow[l, "e"']
  \end{tikzcd}
\end{equation}
When $G$ is an algebraic variety, one calls $(G,m)$ a \emph{group variety}, and when $G$ is an affine scheme, one calls $(G,m)$ an \emph{affine algebraic group} (or \emph{linear algebraic group}).
\end{defn}

The \emph{multiplicative group} $\GG_m(\KK) = \GL(1,\KK)$ or \emph{torus} of dimension one over $\KK$ is represented by $\sO(\GG_m) = \KK[T,T^{-1}] \subset \KK[T]$ (see Milne \cite[Examples 2.2 and 2.8, pp. 40--41]{Milne_algebraic_groups_cup} and Conrad \cite[Examples 1.4.1, 1.4.4, and 1.4.5]{Conrad_linear_algebraic_groups}). The special orthogonal group
\[
  \SO(2,\RR) = \left\{\begin{pmatrix} a & b \\ -b & a \end{pmatrix}: a, b \in \RR \text { and } a^2+b^2 = 1 \right\}
\]
is an $\RR$-torus under the isomorphism
\[
  \SO(2,\RR) \ni \begin{pmatrix} a & b \\ -b & a \end{pmatrix} \mapsto a+ib \in \CC^*.
\]
See Conrad \cite[Definition 4.1.1 and Example 4.1.2]{Conrad_linear_algebraic_groups}.

\begin{defn}[Action of an algebraic group on a scheme]
\label{defn:Action_algebraic_group}  
(See Milne \cite[Section 1f, p. 26]{Milne_algebraic_groups_cup}.)
An \emph{action} of an algebraic group $G$ on an algebraic scheme $X$ is a regular map $\mu:G\times X \to X$ such that the following diagrams commute:
\begin{equation}
  \label{eq:Algebraic_group_action_associativity}
  \begin{tikzcd}
  G\times G\times X \arrow[r, "\id\times \mu"] \arrow[d, "m\times\id"] &G\times X \arrow[d, "\mu"]
  \\
  G\times X \arrow[r, "\mu"] &G
  \end{tikzcd}
\end{equation}
and
\begin{equation}
  \label{eq:Algebraic_group_action_identity}
  \begin{tikzcd}
    \star\times X \arrow[r, "e\times \id"] \arrow[rd, "\cong"'] &G\times X \arrow[d, "\mu"]
    \\
    &X
  \end{tikzcd}
\end{equation}
\end{defn}

The action $\mu$ of a group $G$ on a scheme $X$, as in Definition \ref{defn:Action_algebraic_group}, is called \emph{trivial} if it factors through the projection $G \times X \to X$ (see Milne \cite[Section 7a, p. 138]{Milne_algebraic_groups_cup}.).

\begin{thm}[Existence of a fixed subscheme]
\label{thm:Milne_7-1}
(See Milne \cite[Theorem 7.1, p. 138]{Milne_algebraic_groups_cup}.)
Let $\mu:G\times X \to X$ be an action of an algebraic group $G$ on a separated algebraic scheme $X$ over a field $\KK$. Then there exists a largest closed subscheme $X^G$ of $X$ on which $G$ acts trivially
\end{thm}

We refer to Milne \cite[Section 7b, p. 138]{Milne_algebraic_groups_cup} for a more precise (that is, functorial) statement of Theorem \ref{thm:Milne_7-1}. Regarding properties of fixed-point subschemes, we recall the

\begin{thm}[Smoothness of the fixed subscheme]
\label{thm:Milne_13-1}  
(See Milne \cite[Theorem 13.1, p. 254]{Milne_algebraic_groups_cup}.)  
Let $G$ be a linearly reductive group variety acting on a smooth variety $X$ over a field $\KK$. Then the fixed-point scheme $X^G$ is smooth.
\end{thm}

In particular, tori are linearly reductive (see Milne \cite[Section 13a, p. 254]{Milne_algebraic_groups_cup}).

\section{Analytic groups and their actions on analytic spaces}
\label{sec:Analytic_group_actions_analytic_spaces}
For definitions of transformation groups of $\KK$-analytic spaces for $\KK=\RR$ or $\CC$, we adapt Akhiezer \cite[Sections 1.1 and 1.2]{Akhiezer_lie_group_actions_complex_analysis} for the case $\KK=\CC$.
Guaraldo, Macr\`\i, and Tancredi provide a brief discussion of the action of real Lie groups on real analytic spaces in \cite[Section 8.1]{Guaraldo_Macri_Tancredi_topics_real_analytic_spaces}.

Let $(X,\sO_X)$ be a $\KK$-analytic space as in Definition \ref{defn:Analytic_space}. An \emph{$\KK$-analytic automorphism} of $X$ is a $\KK$-analytic isomorphism of  $(X,\sO_X)$ onto $(X,\sO_X)$ and the set of all such automorphisms is denoted by $\Aut(X)$. Composition of $\KK$-analytic mappings gives a group structure to $\Aut(X)$ (see Akhiezer \cite[Section 1.1, p. 6]{Akhiezer_lie_group_actions_complex_analysis}).

Let $G$ be an abstract group and assume that $X$ has a countable topology. Given a homomorphism $\Phi: G \to \Aut(X)$ of groups, one says that $\Phi$ defines an action of $G$ (or $G$-action) on $X$ and the pair $(G,\Phi)$ is called a \emph{transformation group of $X$}. The automorphism $\Phi(g)$ of $(X,\sO_X)$ is denoted by $g$ and, given an open subset $V\subset X$, one writes $f\circ g$ instead of $f\circ\Phi(g)$ for $f\in\sO_X(V)$ when this causes no ambiguity. A $G$-action on a $\KK$-analytic space $(X, \sO_X)$ induces an action on the underlying topological space $X$ (see Akhiezer \cite[Section 1.2, p. 6]{Akhiezer_lie_group_actions_complex_analysis}).

If $G$ is a topological group, then the $G$-action defined by $\Phi$ is said to be \emph{continuous} if for each open relatively compact subset $U \Subset X$ and for each open subset $V \subset X$, the following hold:
\begin{enumerate}
\item $W := \{g \in G: g\cdot C \subset V\}$, where $C = \bar U$, is open in $G$;
\item for any $f \in \sO_X(U)$, the map
  \begin{equation}
    \label{eq:Action_topological_group_structure_sheaf}
    W \ni g \mapsto f\circ g|_U \in \sO_X(U)
  \end{equation}
is continuous when $\sO_X(U)$ carries the canonical Fr\'echet topology.
\end{enumerate}
If these conditions are fulfilled, $(G,\Phi)$ called a \emph{topological transformation group} of the $\KK$-analytic space $X$. Note that the first condition means that the mapping of topological
spaces $G \times X \ni (g,x) \mapsto gx \in X$ is continuous (see Akhiezer \cite[Section 1.2, p. 7]{Akhiezer_lie_group_actions_complex_analysis}).

Let $G$ be a real (respectively, complex) Lie group and assume now that $(X,\sO_X)$ is a complex analytic space. The $G$-action defined by $\Phi$ is said to be \emph{real analytic} (respectively, \emph{complex analytic} or \emph{holomorphic}) if it is continuous and, in addition, the map \eqref{eq:Action_topological_group_structure_sheaf} is a \emph{real analytic} (respectively, \emph{complex analytic} or \emph{holomorphic}) map from $W$ to the Fr\'echet space $\sO_X(V)$. In this situation, $(G,\Phi)$ is called a  \emph{real} (respectively, \emph{complex}) \emph{(Lie) transformation group} of the complex analytic space $(X,\sO_X)$ (see Akhiezer \cite[Section 1.2, p. 7]{Akhiezer_lie_group_actions_complex_analysis}). 

When $(X,\sO_X)$ is a real analytic space, we restrict $G$ to be a real Lie group and, as above, define the $G$-action to be real analytic and $(G,\Phi)$ to be a real (Lie) transformation group of the real analytic space $(X,\sO_X)$.

Suppose we are given two $G$-actions of the same type (that is, both continuous, real or complex analytic or holomorphic) on complex analytic spaces $(X,\sO_X)$ and $(Y,\sO_Y)$. Let $\Phi: G \to \Aut(X)$ and $\Psi:G \to \Aut(Y)$ be the corresponding group homomorphisms. A complex analytic mapping $\varphi: X \to Y$ is a \emph{morphism of $G$-actions} or a \emph{$G$-equivariant map} if for every $g \in G$ the following diagram commutes:
\[
  \begin{tikzcd}
    X \arrow[r, "\Phi(g)"] \arrow[d, "\varphi"] &X \arrow[d, "\varphi"]
    \\
    Y \arrow[r, "\mu"] &Y
  \end{tikzcd}
\]
The group actions described above are sometimes called \emph{global}. It is also useful to consider \emph{local actions} and those are described in Section \ref{sec:Local_transformation_group_analytic_space}, to which we refer for details.

Lastly, we recall some relevant results on automorphism groups for manifolds and analytic spaces.

\begin{thm}[The automorphism group of a compact, complex analytic space]
\label{thm:Automorphism_group_compact_complex_analytic_space}  
(See Akhiezer \cite[Section 2.3, Theorem, p. 40]{Akhiezer_lie_group_actions_complex_analysis} or Kaup \cite{Kaup_1965}.)  
Let $(X,\sO_X)$ be a compact, complex analytic space. Then the automorphism group $\Aut(X)$ can be endowed with a structure of a complex Lie group so that the action of $\Aut(X)$ on $X$ is holomorphic. The Lie algebra of $\Aut(X)$ is isomorphic to the Lie algebra of all holomorphic vector fields on $X$.
\end{thm}

\begin{rmk}[Other versions of Theorem \ref{thm:Automorphism_group_compact_complex_analytic_space}]
\label{rmk:Automorphism_group_other_versions}  
When $X$ is a complex manifold, Theorem \ref{thm:Automorphism_group_compact_complex_analytic_space} is due to Bochner and Montgomery \cite{Bochner_Montgomery_1947}, while the general case stated is due to Kaup. When $X$ is a $C^k$ (for $k\geq 1$) or a real analytic manifold, the corresponding analogue of Theorem \ref{thm:Automorphism_group_compact_complex_analytic_space} is again due to Bochner and Montgomery \cite{Bochner_Montgomery_1945, Bochner_Montgomery_1946}. For noncompact manifolds, the conclusions of Theorem \ref{thm:Automorphism_group_compact_complex_analytic_space} do not hold: there is no complex connected group acting in a nontrivial manner with complex parameters on any bounded, complex Euclidean domain (see Bochner and Montgomery \cite[Section 3, Theorem 2, p. 663]{Bochner_Montgomery_1947}).
\end{rmk}

\begin{rmk}[Holomorphic Lie group actions on compact, complex projective K\"ahler manifolds]
\label{rmk:Holomorphic_Lie_group_actions_complex_projective_Kaehler_manifolds}  
According to Sommese \cite[Section II, Definition, p. 110]{Sommese_1975}, a complex connected Lie group $G$ \emph{acts projectively} on a compact K\"ahler manifold $X$ if $G$ acts holomorphically and the Lie algebra of holomorphic vector-fields that $G$ generates on $X$ is annihilated by every holomorphic one form. Furthermore, according to Sommese \cite[Remark II-B, p. 110]{Sommese_1975}, the preceding definition is justified by a theorem due to Blanchard \cite{Blanchard_1956} which states that if $X$ is projective and $G$ is as above, then $X$ can be embedded into some $\PP^n$ equivariantly with respect to a faithful representation of $G$ in $\mathrm{PGL}(n,\CC)$. The theorem is a consequence of a fixed-point theorem dexcribed in Sommese \cite{Sommese_1974} or Lieberman \cite[Section III]{Lieberman_1982}.
\end{rmk}

We can now state and prove the following generalization of Theorem \ref{thm:Equivariance_property_blowups_complex_manifolds} from the category of complex $G$-manifolds to complex analytic $G$-spaces.

\begin{thm}[Equivariance property of blowups for complex analytic $G$-spaces]
\label{thm:Equivariance_property_blowups_complex_analytic_spaces}  
Let $G$ be a Lie group and let $(X,\sO_X)$ be a complex analytic $G$-space. If $(Z,\sO_Z)$ is a closed, complex analytic, $G$-invariant subspace of $(X,\sO_X)$, then the blowup $\Bl_Z(X)$ is a complex analytic $G$-space and the canonical projection $\pi: \Bl_Z(X) \to X$ is a holomorphic $G$-equivariant morphism. If the action of $G$ on $X$ is proper, then also the action of $G$ on $\Bl_Z(X)$ is proper.
\end{thm}

\begin{proof}
We may apply Proposition \ref{prop:Functorial_property_blowup_complex_analytic_space} by choosing $Y = X$ and $f = g$, for any $g \in G$. The fact that $\Bl_Z(X)$ is a complex analytic space is given by Theorem \ref{thm:Fischer_4-1}. Proposition \ref{prop:Functorial_property_blowup_complex_analytic_space} ensures that each element $g \in G \subset \Aut(X)$ (the group of biholomorphic maps of $X$ to itself) lifts to a unique element $\Bl_Z(g) \in \Aut(\Bl_Z(X))$ such that the following diagram commutes:
\[
  \begin{tikzcd}
    \Bl_Z(X) \arrow[r, "\Bl_Z(g)"] \arrow[d, "\pi_X"] &\Bl_Z(X) \arrow[d, "\pi_X"]
    \\
    X \arrow[r, "g"] &X
  \end{tikzcd}
\]
In other words, the blowup $\Bl_Z(X)$ is a complex analytic $G$-space and the canonical projection $\pi_X: \Bl_Z(X) \to X$ is a holomorphic, $G$-equivariant mmorphism. If the action of $G$ on $X$ is proper, then so is the action of $G$ on $\Bl_Z(X)$ because the blowup morphism $\pi_X$ is proper and the preceding diagram commutes.
\end{proof}

\section{Linearization for actions of compact groups on manifolds and complex analytic spaces}
\label{sec:Linearization_theorems_compact_group_actions_manifolds_complex_analytic_spaces}
We first recall the

\begin{thm}[Bochner linearization theorem for actions of compact groups on $C^k$ manifolds]
\label{thm:Bochner_linearization_theorem}  
(See Bochner \cite{Bochner_1945}, Duistermaat and Kolk \cite[Theorem 2.2.1, p. 96]{DuistermaatLieGroups}, or Kydonakis \cite[Theorem 3.1]{Kydonakis_2014_note}.)
Let $M$ be a finite-dimensional, real-analytic manifold and $A$ be a continuous homomorphism from a compact topological group $G$ to $\Diff^k(M)$, where $k$ is a positive integer, $\infty$, or $\omega$ (real analytic). If $p \in M$ is a point such that $A(g)(p) = p$ for all $g \in G$, then there  exists a $G$-invariant open neighborhood $U$ of $p$ in $M$ and a $C^k$ diffeomorphism $\chi$ from $U$ onto an open neighborhood $V$ of the origin $0 \in T_pM$ such that
\[
  \chi(p) = 0 \quad\text{and}\quad D\chi(p) = \id \in \End(T_pM)
\]
and, noting that $A(g) \in \Diff^k(M)$ and $D(A(g))(p) \in \GL(T_pM)$ for all $g \in G$,
\[
  \chi(A(g)(x)) = D(A(g))(p)\chi(x), \quad\text{for all } g \in G \text{ and } x \in U.
\]
\end{thm}

The Bochner Linearization Theorem \ref{thm:Bochner_linearization_theorem} in the category of real $C^k$ real manifolds was motivated by the following theorem due to Cartan in the category of complex manifolds.

\begin{thm}[Cartan linearization theorem for actions of compact groups on complex manifolds]
\label{thm:Cartan_linearization_theorem}  
(See Bochner \cite[Item (i), p. 372]{Bochner_1945}, Cartan \cite{Cartan_1935} or \cite[Article 32, pp. 474--523]{Cartan_collected_works_I}, and Martin \cite[Theorems 3 and 7]{Martin_1944}.)
Let $X$ be an $n$-dimensional, complex manifold and $G$ be a compact group of automorphisms of $X$. If $p \in X$ is a fixed point of $G$, then there  exists a $G$-invariant open neighborhood $U$ of $p$ in $X$ and a biholomorphic map $\varphi$ from $U$ onto an open neighborhood $V$ of the origin $0 \in \CC^n$ such that $\varphi(p) = 0$ and 
\[
  \varphi(g(z)) \in \GL(n,\CC), \quad \text{for all } g \in G \text{ and } z \in U.
\]
\end{thm}

By automorphism of $X$ in Theorem \ref{thm:Cartan_linearization_theorem} one means a biholomorphic transformation. Following Akhiezer \cite[Section 2.2, p. 35]{Akhiezer_lie_group_actions_complex_analysis}, suppose now that $X$ is an arbitrary complex space, $G$ is a group acting on $X$, and $p \in X$ is a fixed point of $G$. Then $G$ acts naturally on the local algebra $\sO_{X,p}$ by algebra automorphisms. Namely, for $g \in G$, we have an automorphism $\sO_{X,p} \ni f_p \mapsto (f \circ g^{-1})_p \in \sO_{X,p}$, where $f$ is a holomorphic function defined on an open neighborhood $U$ of $p$. Given $g$, the germ $(f \circ g^{-1})_p$ depends only on the germ $f_p$ and not on its representative $f \in \sO_X(U)$. This representation of $G$ on $\sO_{X,p}$ induces the \emph{isotropy representation},
\begin{equation}
  \label{eq:Isotropy_representation}
  \tau_p:G \to \GL(T_pX),
\end{equation}
where $T_pX = (\fm_p/\fm_p^2)^*$ is the Zariski tangent space to $X$ at the point $p$ and $\fm_p$ is the maximal ideal in the local ring $\sO_{X,p}$. Kaup \cite{Kaup_1965} provides the following generalization of Theorem \ref{thm:Cartan_linearization_theorem} from the category of complex manifolds to complex analytic spaces.

\begin{thm}[Kaup linearization theorem for actions of compact groups on complex analytic spaces]
\label{thm:Kaup_linearization_theorem}
(See Akhiezer \cite[Section 2.2, p. 36, Linearization Theorem]{Akhiezer_lie_group_actions_complex_analysis} or Kaup \cite[Satz 5, p. 84]{Kaup_1965}.)
Let $K$ be a compact topological group acting continuously on a complex analytic space $X$. If $p \in X$ is a fixed point of $K$, then there exist a $K$-invariant neighborhood $U$ of $p$ in $X$, a $\tau_p(K)$-invariant neighborhood $V$ of the origin $0 \in T_pX$, and a closed $K$-equivariant embedding $\varphi: U \to V$.
\end{thm}

For our application, a variant (see the forthcoming Corollary \ref{cor:Cartan_linearization_theorem_equivariant}) of the Cartan Linearization Theorem \ref{thm:Cartan_linearization_theorem} will suffice for our application, but we include the statements of Theorems \ref{thm:Bochner_linearization_theorem} and \ref{thm:Kaup_linearization_theorem} for the sake of completeness, references to their proofs, and relevance of the isotropy representation \eqref{eq:Isotropy_representation}. The following corollary of the proof of Theorem \ref{thm:Kaup_linearization_theorem} will be especially useful and the same argument yields the corresponding refinement of Theorem \ref{thm:Bochner_linearization_theorem}. Of course, our proof of this refinement follows \mutatis for real $C^k$ manifolds and submanifolds as well.

\begin{cor}[Linearization for actions of compact groups on analytic manifolds and submanifolds]
\label{cor:Cartan_linearization_theorem_equivariant}
Let $\KK=\RR$ or $\CC$, and $X$ be a $\KK$-analytic manifold, $Z \subset X$ be an embedded, $\KK$-analytic submanifold, and $\Phi:G \to \Aut(X)$ be a continuous homomorphism from a compact topological group $G$ into the group $\Aut(X)$ of $\KK$-analytic automorphisms of $X$. If $p \in X$ is a fixed point under the action $G$ on $X$ and the submanifold $Z$ is $G$-invariant, then there exist
\begin{enumerate}
\item An open neighborhood $U$ of $p$ in $X$ that is invariant under the action of $\Phi(g) \in \Aut(X)$, for all $g\in G$;
\item An open neighborhood $V$ of the origin $0 \in T_pX$ that is invariant under the action of $\d(\Phi(g))(p) \in \GL(T_pX)$, for all $g\in G$; and
\item A $\KK$-bianalytic, $G$-equivariant map $\varphi:U\to V$
\end{enumerate}
such that the following hold:
\begin{itemize}
\item $\varphi(p) = 0$ and $\d\varphi(p) = \id \in \End(T_pX)$; and
\item $\varphi(U\cap Z) = V\cap T_pZ$.
\end{itemize}
\end{cor}

\begin{proof}
We shall modify the proof of Theorem \ref{thm:Bochner_linearization_theorem} given by Duistermaat and Kolk \cite[Section 2.2, pp. 96--97]{DuistermaatLieGroups} and the proof of 
Theorem \ref{thm:Kaup_linearization_theorem} given by Akhiezer \cite[Section 2.2, pp. 36--37]{Akhiezer_lie_group_actions_complex_analysis}. When $S$ is an embedded, real, smooth submanifold of dimension $k$ of a real, smooth manifold $M$ of dimension $m$ and $x_0\in S$ is a point, then there are an open neighborhood $U \subset M$ of $x_0$ and a smooth coordinate chart $\psi:M \subset U \to \RR^m$ such that $\psi(U\cap S) = \psi(U)\cap\RR^k$, where $\RR^k \hookrightarrow \RR^m$ has its standard embedding as a coordinate subspace (see Lee \cite[Theorem 5.8, p. 101]{Lee_john_smooth_manifolds}). By replacing the role in the proofs of \cite[Theorem 4.12, p. 81, and Theorem 5.8, p. 101]{Lee_john_smooth_manifolds} given by Lee of the Implicit Mapping Theorem for smooth maps of smooth manifolds by the role of the Implicit Mapping Theorem for $\KK$-analytic maps of $\KK$-analytic manifolds, one obtains the corresponding result for $\KK$-analytic manifolds. Hence, if $X$ and $Z$ have dimensions $n$ and $l$ over $\KK$, respectively, there are an open neighborhood $U \subset X$ of $p$ and a $\KK$-analytic coordinate chart $\phi:X \supset U \to \KK^n$ such that $\phi(U\cap Z) = \phi(U)\cap\KK^l$, where $\KK^l \hookrightarrow \KK^n$ has its standard embedding as a coordinate subspace.

The continuous homomorphism $\Phi:G \to \Aut(U)$ and the $\KK$-analytic diffeomorphism $\phi:X \supset U \to\phi(U) \subset \KK^n$ induce a $G$-action on the open subset $\phi(U) \subset \KK^n$, which we denote by $\phi:G \ni g \mapsto g\circ\phi^{-1} \in \Aut(\phi(U))$, such that $\phi(U)\cap\KK^l$ is $G$-invariant. The tangent space $T_pX$ and subspace $T_pZ \subset T_pX$ may be identified with $\KK^n$ and the coordinate subspace $\KK^l \subset \KK^n$, respectively, via the $\KK$-linear isomorphism $d\phi(p):T_pX \to \KK^n$, so that $\phi(p) = 0$ and $\d\phi(p) = \id \in \End(T_pX)$.  

For each $g \in G$, we let $\phi^g$ denote the composition
\[
  \begin{tikzcd}
    U \xrightarrow{\ \Phi(g)\ } U \xrightarrow{\ \phi\ } T_pX  \xrightarrow{\d(\Phi(g))(p)^{-1}} T_pX.
  \end{tikzcd}
\]
Then $g \mapsto \phi^g$ is a continuous function on $G$ with values in the Fr\'echet space $\An(U, T_pX)$ of $\KK$-analytic maps from $U$ into $T_pX$. By averaging $\phi^g$ over $g \in G$, we obtain a $\KK$-analytic map
\[
  \varphi := \int_G \phi^g\cdot d\mu(g) \in \An(U, T_pX),
\]
where $\mu$ is the normalized Haar measure on $G$. The invariance of $\mu$ implies that
\[
  \varphi\circ\Phi(g) = \d(\Phi(g))(p)\circ\varphi, \quad\text{for all } g \in G,
\]
that is, $\varphi:U \to T_pX$ is a $G$-equivariant map.

By hypothesis, $U\cap Z$ is $G$-invariant and because $\phi:U \to T_pX$ is $G$-equivariant, the image
\[
  \phi(U\cap Z) = \phi(U)\cap \KK^l = \phi(U)\cap T_pZ
\]
is also $G$-invariant. Since the $\KK$-analytic automorphism $\Phi(g)\in \Aut(U)$ restricts to a $\KK$-analytic automorphism $\Phi(g) \in \Aut(U\cap Z)$, then the isotropy representation, $\d(\Phi(g))(p) \in \End(T_pX)$, restricts to $\d(\Phi(g))(p) \in \End(T_pZ)$. Hence, each map $\phi^g:U\to T_pX$ restricts to a $\KK$-analytic map
\[
  U\cap Z \xrightarrow{\ \Phi(g)\ } U\cap Z \xrightarrow{\ \phi\ }
  T_pZ  \xrightarrow{\d(\Phi(g))(p)^{-1}} T_pZ.
\]
and therefore the average, $\varphi$, restricts to a $G$-equivariant, $\KK$-analytic map, $\varphi:U\cap Z \to T_pZ$.

We now compute the differential of $\varphi$ at the point $p$ using the expression
\[
  \varphi = \int_G \Phi(g) \circ \phi \circ \d(\Phi(g))(p)^{-1}\cdot d\mu(g)
\]
and the fact that $\d\phi(p) = \id \in \End(T_pX)$ to give
\[
  \d\varphi(p) = \int_G \d(\Phi(g))(p) \circ \d\phi(p) \circ \d(\Phi(g))(p)^{-1}\cdot d\mu(g)
   = \id \in \End(T_pX).
\]
The Inverse Mapping Theorem for $\KK$-analytic maps now implies that, after possibly shrinking $U$, the map $\varphi:X \supset U \to T_pX$ is a $\KK$-analytic embedding from the open neighborhood $U$ of $p$ in $X$ onto an open neighborhood $V = \varphi(U)$ of the origin in $T_pX$ and this completes the proof of Corollary \ref{cor:Cartan_linearization_theorem_equivariant}.
\end{proof}

While its proof is an immediate consequence of the definitions, the following elementary observation has useful applications and we include the statement here for completeness.

\begin{lem}[Fixed points of group actions on sets and invariant subsets]
\label{lem:Fixed_points_group_actions_on_subsets_and_invariant_subsets} 
Let $G$ be a group, $X$ be a set, and $\Phi:G\to \Aut(X)$ be a homomorphism from $G$ into the group $\Aut(X)$ of automorphisms of $X$. If $Y \subset X$ be a subset that is invariant under the action
\[
  G\times X \ni (g,x) \mapsto g\cdot x := \Phi(g)(x) \in X,
\]
then
\[
  Y^G = Y\cap X^G,
\]
where $Y^G := \{y\in Y: g\cdot y = y, \text{ for all } g \in G\}$ and $X^G := \{x\in X: g\cdot x = x, \text{ for all } g \in G\}$.
\end{lem}

\section[Linearization for actions of noncompact groups]{Weyl's unitary trick and linearization for actions of noncompact groups on complex manifolds}
\label{sec:Linearization_theorems_noncompact_group_actions_manifolds}
It is known that the Bochner Linearization Theorem \ref{thm:Bochner_linearization_theorem} can be extended from the case of compact Lie groups to noncompact Lie groups with compact real form by applying Weyl's `unitary trick' (see Hilgert and Neeb \cite[Theorem 5.2.10, p. 576]{Hilgert_Neeb_structure_geometry_lie_groups} or Varadarajan \cite[Lemma 4.11.13, p. 349]{Varadarajan}). For example, see Cairns and Ghys \cite[Theorem 2.6, p. 139]{Cairns_Ghys_1997}, Guillemin and Sternberg \cite{Guillemin_Sternberg_1968}, and Ku\v{s}nirenko \cite{Kusnirenko_1967}. We shall subsequently apply such an extension for the Cartan Linearization Theorem \ref{thm:Cartan_linearization_theorem} and in particular its Corollary \ref{cor:Cartan_linearization_theorem_equivariant} for the case of $G = \U(1)$ with $G_\CC = \CC^*$. However, because $G_\CC$ is noncompact one needs to broaden the concept of linearization.

\begin{defn}[Linearizable group action]
\label{defn:Cairns_Ghys_1-3}  
(See Cairns and Ghys \cite[Definition 1.3 p. 135]{Cairns_Ghys_1997}.)  
Let $G$ be a group, $M$ be a $C^1$ manifold, $p\in M$ be a point, and $\Phi:G\to\Aut(M)$ be a homomorphism from $G$ into the group $\Aut(M)$ of $C^1$ diffeomorphisms of $M$. One says that $\Phi$ is \emph{linearizable at $p$} if there are an open neighborhood $U$ of $p$ in $M$ and an open neighborhood $V$ of the origin in $T_pM$ and a homeomorphism $\varphi: U \to V$ such that $\varphi(p) = 0$ and for each $g \in G$ the maps
\[
  U \ni x \mapsto \varphi(\Phi(g)(x)) \in T_pM
  \quad\text{and}\quad U
  \ni x \mapsto \d(\Phi(g))(p)\varphi(x) \in T_pM,
\]
have the same germ at $p$, that is, there exists an open neighborhood $U_g \subset U$ of $p$ such that
\[
  \varphi\circ\Phi(g) = \d(\Phi(g))(p)\circ\varphi \quad\text{on } U_g.
\]  
If $M$ and $\varphi$ and $\Phi(g)$ are $C^k$ for all $g\in G$, where $k$ is a positive integer, $\infty$, or $\omega$, or $M$ is a complex manifold and the maps $\varphi$ and $\Phi(g)$ are biholomorphic, for all $g\in G$, then one says that $\Phi$ is $C^k$ or \emph{holomorphically linearizable}.  
\end{defn}

\begin{rmk}[Sumihiro's Theorem and proof of the decomposition Bia{\l}ynicki--Birula decomposition for smooth algebraic varieties]
\label{rmk:Sumihiro_theorem_and_proof_BB_decomposition}  
The cited results due to Cairns and Ghys \cite{Cairns_Ghys_1997}, Guillemin and Sternberg \cite{Guillemin_Sternberg_1968}, and Ku\v{s}nirenko \cite{Kusnirenko_1967} may be viewed, when specialized to the complexification $G_\CC = \CC^*\times\cdots\times\CC^*$ of $G=\U(1)\times\U(1)$, as analogues of a well-known result \cite[Corollary 2, p. 8]{Sumihiro_1974} due to Sumihiro in the context of algebraic groups acting on algebraic varieties: Let $X$ be a normal algebraic variety over an algebraically closed field $k$ and let $T$ be a torus group that acts regularly on $X$. Then, for any point $p \in X$, there is a $T$-stable affine open neighborhood of $p$ in $X$. This result was used by Konarski \cite{Konarski_1996} to give a new and elementary proof of the Bia{\l}ynicki--Birula decomposition for smooth algebraic varieties (Theorem \ref{thm:Milne_13-47}). 
\end{rmk}

We prove the forthcoming Theorem \ref{thm:Cartan_linearization_theorem_equivariant_noncompact} using Weyl's Unitary trick to pass from a compact, real Lie group to a noncompact, complex Lie group. Before proceeding to this application, we recall that a \emph{complex Lie group} is a real Lie group $G$ whose Lie algebra $\bg = L(G)$ is a complex Lie algebra and for which $\Ad(G) \subseteq \Aut_\CC(\bg)$, the group of complex linear automorphisms of g (see Hilgert and Neeb \cite[Definition 15.1.1 (a)]{Hilgert_Neeb_structure_geometry_lie_groups}). A homomorphism $\alpha: G_1 \to G_2$ of complex Lie groups is called \emph{holomorphic} if the induced homomorphism of Lie algebras, $L(\alpha):L(G_1)\to L(G_2)$, is complex linear; if $G_2 = \GL(V)$ for a complex vector space, then a holomorphic homomorphism $\alpha: G_1 \to G_2$ is also called a \emph{holomorphic representation} of $G_1$ on $V$ (see \cite[Definition 15.1.1 (b)]{Hilgert_Neeb_structure_geometry_lie_groups}). A subgroup $H$ of a complex Lie group $G$ is called a \emph{complex Lie subgroup} if $H$ is closed and its Lie algebra $L(H)$ is a complex subspace of $L(G)$ (see \cite[Definition 15.1.1 (c)]{Hilgert_Neeb_structure_geometry_lie_groups}).

\begin{defn}[Universal complexification of a Lie group]
\label{defn:Hilgert_Neeb_5-1-2}  
(See Hilgert and Neeb \cite[Definition 15.1.2, p. 566]{Hilgert_Neeb_structure_geometry_lie_groups})  
Let $G$ be a real Lie group. A pair $(\eta_G,G_\CC)$, comprising a complex Lie group $G_\CC$ and a morphism $\eta_G: G \to G_\CC$ of real Lie groups, is called a \emph{universal complexification} of $G$ if for every homomorphism $\alpha: G \to H$ into a complex Lie group $H$, there exists a unique holomorphic homomorphism $\alpha_\CC:G_\CC \to H$ with $\alpha_\CC \circ \eta_G = \alpha$. This is also called the \emph{universal property} of $G_\CC$.
\end{defn}

According to Hilgert and Neeb \cite[Remark 15.1.3, p. 566]{Hilgert_Neeb_structure_geometry_lie_groups}, $G_\CC$ is unique due to the universal property and by Hilgert and Neeb \cite[Theorem 15.1.4, p. 566]{Hilgert_Neeb_structure_geometry_lie_groups}, the group $G_\CC$ exists and has the five properties listed there. We shall often appeal to the following important example.

\begin{exmp}[Universal complexification of the circle group]
\label{exmp:Universal_complexification_circle_group}
(See Hilgert and Neeb \cite[Example 15.1.6 (a), p. 569]{Hilgert_Neeb_structure_geometry_lie_groups}.)  
If $G = S^1 = \{z\in \CC:|z|=1\}$ and $\eta_{S^1}: S^1 \hookrightarrow \CC^*$ is the standard inclusion, then $(\eta_{S^1},\CC^*)$ is a universal complexification of $S^1$.
\end{exmp}

\begin{thm}[Linearization for noncompact group actions on complex manifolds and submanifolds]
\label{thm:Cartan_linearization_theorem_equivariant_noncompact}
Let $X$ be a complex manifold, $Z \subset X$ be an embedded, complex submanifold, $G$ be a compact real Lie group with complexification $G_\CC$, and $\Phi:G_\CC \to \Aut(X)$ be a holomorphic homomorphism from $G_\CC$ onto a subgroup of the group $\Aut(X)$ of biholomorphic automorphisms of $X$. If $p \in X$ is a fixed point under the action $G_\CC$ on $X$ and the submanifold $Z$ is $G_\CC$-invariant, then there exist
\begin{enumerate}
\item An open neighborhood $U$ of $p$ in $X$ that is invariant under the action of $\Phi(g) \in \Aut(X)$, for all $g\in G$;
\item An open neighborhood $V$ of the origin $0 \in T_pX$ that is invariant under the action of $\d(\Phi(g))(p) \in \GL(T_pX)$, for all $g\in G$; and
\item A biholomorphic, $G$-equivariant map $\varphi:U\to V$
\end{enumerate}
such that the following hold:
\begin{itemize}
\item $\varphi(p) = 0$ and $\d\varphi(p) = \id \in \End(T_pX)$;
\item $\varphi(U\cap Z) = V\cap T_pZ$; and
\item $G_\CC$ is holomorphically linearizable at $p$ in the sense of Definition \ref{defn:Cairns_Ghys_1-3}.
\end{itemize}
\end{thm}

\begin{proof}
We continue the notation of the proof of Corollary \ref{cor:Cartan_linearization_theorem_equivariant}, noting that all assertions but the final one regarding holomorphic linearizability of $G_\CC$ at $p$ are provided by Corollary \ref{cor:Cartan_linearization_theorem_equivariant}. By analogy with the proof of Cairns and Ghys \cite[Theorem 2.6, p. 139]{Cairns_Ghys_1997}, we define
\[
  S := \left\{g \in G_\CC: \varphi\circ\Phi(g) = \d(\Phi(g))(p)\circ\varphi \text{ on } U_g\right\},
\]
where each $U_g \subseteq U$ is an open neighborhood of $p$ that may depend on $g$. By construction, $S$ is a complex Lie subgroup of $G_\CC$ that contains $G$. By Hilgert and Neeb \cite[Theorem 5.1.4 and Remark 15.1.3, p. 566]{Hilgert_Neeb_structure_geometry_lie_groups}, the universal complexification $G_\CC$ of $G$ exists (even without our hypothesis that $G$ is compact) and is unique up to isomorphism of complex Lie groups, so we must have $S = G_\CC$. Thus, $G_\CC$ is holomorphically linearizable at $p$ in the sense of Definition \ref{defn:Cairns_Ghys_1-3}.
\end{proof}

\section{Holomorphic  maps of complex vector spaces with circle actions}
\label{sec:Holomorphic_maps_vector_spaces_with_circle_actions}
To motivate the results of this section, we first recall the following generalization of Blanchard's Lemma \cite[Section 1.1]{Blanchard_1956} (for holomorphic transformation groups on complex analytic spaces) to the category of algebraic groups acting on schemes.

\begin{thm}[Blanchard's Lemma]
\label{thm:Blanchard_schemes}  
(See Brion \cite[Theorem 7.2.1]{Brion_2017} or Brion, Samuel, and Uma \cite[Proposition
4.2.1]{Brion_Samuel_Uma_lectures_structure_algebraic_groups_geometric_applications}.)
Let $G$ be a connected algebraic group, $X$ a $G$-scheme of finite type, $Y$ a scheme of finite type, and $\varphi : X \to Y$ a proper morphism such that $\varphi^\sharp: \sO_Y \to \varphi_*(\sO_X)$ is an isomorphism. Then there exists a unique action of $G$ on $Y$ such that $\varphi$ is $G$-equivariant.
\end{thm}

The following lemma is closer to the original version of Blanchard's Lemma in the setting of holomorphic transformation groups acting on complex analytic spaces; see also Bierstone \cite{Bierstone_1977} for related results.

\begin{lem}[Blanchard's lemma for Lie group actions on complex analytic spaces]
\label{lem:Blanchard_complex_analytic_spaces}
(See Akhiezer \cite[Section 2.4, p. 44, Lemmas 1 and 2, and Remark, p. 45]{Akhiezer_lie_group_actions_complex_analysis} and Blanchard \cite[Section 1.1]{Blanchard_1956}.)  
Let $\varphi:X \to Y$ be a proper, holomorphic map between reduced complex analytic spaces $X$ and $Y$ such that $\sO_Y = \varphi_*(\sO_X)$. If there is a real analytic action of a connected, real Lie group $G$ on $X$, then the following hold:
\begin{enumerate}
\item Each $g \in G$ induces a biholomorphic transformation of $Y$, the induced action of $G$ on $Y$ is real analytic, and the map $\varphi$ is $G$-equivariant;
\item If in addition $G$ is a complex Lie group and the action of $G$ on $X$ is holomorphic, then the induced action of $G$ on $Y$ is holomorphic.
\end{enumerate}  
\end{lem}

\begin{rmk}[Zariski's Connectedness Theorem and Blanchard's Lemma]
\label{rmk:Blanchard_lemma}  
The author is very grateful to Michel Brion for providing the author with the following explanation of the relationship between the hypotheses of Theorem \ref{thm:Blanchard_schemes} and Lemma \ref{lem:Blanchard_complex_analytic_spaces}. If $\varphi: X \to Y$ is a proper morphism of Noetherian schemes such that $\sO_Y = \varphi_*(\sO_X)$, then the fibers of $\varphi$ are connected by \emph{Zariski's Connectedness Theorem}. The latter result is proved for projective morphisms by Hartshorne as \cite[Chapter 3, Section 11, p. 279, Corollary 11.3]{Hartshorne_algebraic_geometry} and in the general case of proper morphisms by Grothendieck as \cite[Corollaire 4.3.2, p. 131]{Grothendieck_EGAIII-1}.

There is a partial converse to Zariski's Connectedness Theorem. If $\varphi : X \to Y$ is a proper morphism of algebraic varieties over a field $\KK$ that has characteristic zero and the fibers of $\varphi$ are connected and $Y$ is normal, then $\sO_Y = \varphi_*(\sO_X)$ (this follows from the existence of the Stein factorization).
\end{rmk}

In his proof of Lemma \ref{lem:Blanchard_complex_analytic_spaces}, Akhiezer relies on the following version of Zariski's Connectedness Theorem for complex analytic spaces; we include a proof due to Murayama since we are unable to find any other reference.

\begin{lem}[Zariski connectedness theorem for complex analytic spaces]
\label{lem:Zariski_connectedness_theorem_complex_analytic_spaces}
(See Murayama \cite{Murayama_2016mathoverflow_geometric_intuition_Stein_factorization_theorem}.)  
If $\varphi:X\to Y$ is a proper holomorphic map of complex analytic spaces such that $\sO_Y = \varphi_*(\sO_X)$, then the fibers $\varphi^{-1}(y)$ are connected for all $y\in Y$.
\end{lem}

\begin{proof}
Murayama adapts Hartshorne's proof of \cite[Chapter 3, Section 11, p. 279, Corollary 11.3]{Hartshorne_algebraic_geometry}. The hypotheses imply that $\varphi:X\to Y$ is surjective and so its fibers $\varphi^{-1}(y)$ are non-empty for all $y\in Y$. Suppose that $\varphi^{-1}(y)$ is disconnected. Then there exists an open neighborhood $U$ of $\varphi^{-1}(y)$ which is disconnected. Because $\varphi$ is a continuous, proper map of topological spaces, it is necessarily a closed map.
By shrinking the neighborhood if necessary, we can assume that $U = U_1\cup U_2$ has the form $\varphi^{-1}(V)$ for an open neighborhood of $y$ in $Y$ (because $\varphi$ is closed and by properties of closed maps of topological spaces, see for example, Grauert and Remmert \cite[Section 2.3.1, p. 46, Lemma]{Grauert_Remmert_coherent_analytic_sheaves}) and that $U_1\cap U_2 = \varnothing$. The section of $\sO_X$ which is equal to $1$ on $U_1$ and equal to $0$ on $U_2$ gives a section $f$ of $\sO_Y$ over $V$, by the identity $\sO_Y = \varphi_*(\sO_X)$, such that $f(y) = 1$ and $f(y) = 0$, a contradiction.
\end{proof}

A holomorphic map $(X,\sO_X)\to(Y,\sO_Y)$ of complex analytic spaces comprises a continuous map $\varphi:X\to Y$ topological spaces and a homomorphism $\varphi^\sharp:\sO_Y \to \varphi_*(\sO_X)$. If $V \subset Y$ is an open subset and $h \in \sO_Y(V)$, then $h\circ\varphi \in \sO_X(\varphi^{-1}(V))$ and this defines $\varphi_x^\sharp:\sO_{Y,\varphi(x)} \to \varphi_*(\sO_{X,x})$, for all $x \in X$.

In our application, we shall not easily be able to verify that $\varphi^\sharp:\sO_Y \cong \varphi_*(\sO_X)$ is actually an isomorphism, so we shall explore an alternative approach --- see the forthcoming Proposition \ref{prop:Equivariance_complex_analytic_maps_complex_Banach_spaces_circle_actions} and Lemma \ref{lem:Blanchard_holomorphic_map_from_domain_around_fixed_point_into_vector_space} --- to obtaining the conclusion of Blanchard's Lemma \ref{lem:Blanchard_complex_analytic_spaces} in the special case where $(X,\sO_X) = (D,\sO_D)$, for a domain $D \subset \CC^n$ that is invariant under the (possibly nonlinear) action of the circle group $G = S^1$. In the forthcoming Proposition \ref{prop:Equivariance_complex_analytic_maps_complex_Banach_spaces_circle_actions}, we allow the domain and codomain to be complex Hilbert spaces but restrict our attention to the case $G = S^1$ with a linear representation on the domain. (For the theory of Banach space representations of Lie groups, we refer to Lyubich \cite{Lyubich_introduction_theory_Banach_representations_groups}.)

\begin{prop}[Equivariance of complex analytic maps of complex Hilbert spaces with orthogonal circle actions]
\label{prop:Equivariance_complex_analytic_maps_complex_Banach_spaces_circle_actions}  
Let $W$ and $V$ be complex Hilbert spaces and $S^1\times W \to W$ be a real analytic circle action defined by a unitary representation $S^1\to\U(W)$ given by an orthogonal direct sum decomposition
\begin{equation}
  \label{eq:domain_decomposition}
  W = \bigoplus_{k=0}^mW_k
\end{equation}
into closed complex linear subspaces $W^k \subset W$, for $k=0,1,\ldots,m$, and integer weights $l_0=0$ and\footnote{The integers $l_k$ may be positive or negative when $k\neq 0$.} $l_1<l_2<\cdots<l_m$ such that
\begin{equation}
  \label{eq:domain_circle_action}
  e^{i\theta}\cdot x = \left(x_0,e^{il_1\theta}x_1,\ldots,e^{il_m\theta}x_m\right),
   \quad\text{for all } x = (x_0,x_1,\ldots,x_m) \in W \text{ and } \theta \in \RR,
\end{equation}
where $e^{il_k\theta}$ acts on $x_k$ by complex scalar multiplication, for $k=0,\ldots,m$. If $U \subset W$ is an $S^1$-invariant open neighborhood of the origin and $F:W\supset U\to V$ is a complex analytic map with Taylor series\footnote{See Whittlesey \cite{Whittlesey_1965} for an introduction to analytic maps of Banach spaces.} having radius of convergence $r>0$,
\begin{equation}
  \label{eq:Taylor_series}
  F(x) = \sum_{|\alpha|\geq 0}\frac{1}{\alpha!}D^\alpha F(0)x^\alpha, \quad\text{for all } x \in B_r(0) \subseteq U,
\end{equation}
where $\alpha\in\ZZ_{\geq 0}^{m+1}$ and $|\alpha| = \alpha_0+\alpha_1+\cdots+\alpha_m$ and $\alpha! := \alpha_0!\alpha_1!\cdots\alpha_m!$ and
\[
  D^\alpha F(0) := D_0^{\alpha_0}D_1^{\alpha_1}\cdots D_m^{\alpha_m} F(0) \in \Hom(\otimes^\alpha W,V)
\]
are bounded, multilinear operators, and $x^\alpha := x_0^{\alpha_0}x_1^{\alpha_1}\cdots x_m^{\alpha_m} \in \otimes^\alpha W$, and
\[
  \otimes^\alpha W := (\otimes^{\alpha_0} W_0) \otimes \cdots \otimes(\otimes^{\alpha_m} W_m),
  \quad\text{where } \otimes^{\alpha_k} W_k := \underbrace{W_k\otimes\cdots\otimes W_k}_{\alpha_k\text{ factors}},
\]
and for each $(w_0,\ldots,w_m) \in W$ and $k\in\{0,\ldots,m\}$, the partial derivative $D_kF(w_0,\ldots,w_m) \in \Hom(W_k,V)$ is the bounded linear operator that is uniquely defined by the condition
\begin{multline}
  \label{eq:Partial_derivative_map_Banach_spaces}
  \left\|F(w_0,\ldots,w_{k-1},w_k+x_k,w_{k+1},\ldots,w_m) \right.
    \\
    \left. - F(w_0,\ldots,w_m) - D_kF(w_0,\ldots,w_m)x_k\right\|_V = o\left(\|x_k\|_W\right),
    \quad\text{for all } x_k \in W_k,
\end{multline}
then the following hold.
\begin{enumerate}
\item\label{item:Circle_actions_on_domain_and_codomain_constrain_map}
  Assume that $S^1\times V \to V$ is a real analytic circle action defined by a unitary representation $S^1\to\U(V)$ given by an orthogonal direct sum decomposition
\begin{equation}
  \label{eq:codomain_decomposition}
  V = \bigoplus_{j=0}^n V_j
\end{equation}
into closed complex linear subspaces $V_j \subset V$, for $j=0,1,\ldots,n$, and integer weights $m_0=0$ and\footnote{The integers $m_j$ may be positive or negative when $j\neq 0$.} $m_1< m_2 < \cdots < m_n$ such that
\begin{equation}
  \label{eq:codomain_circle_action}
  e^{i\theta}\cdot y = \left(y_0,e^{im_1\theta}y_1,\ldots,e^{im_n\theta}y_n\right),
   \quad\text{for all } y = (y_0,y_1,\ldots,y_n) \in V \text{ and } \theta \in \RR,
\end{equation}
where $e^{im_j\theta}$ acts on $y_j$ by complex scalar multiplication for $j=0,\ldots,n$. If $F$ is $S^1$-equivariant with respect to the real analytic $S^1$ actions on $W$ and $V$ defined by the unitary representations \eqref{eq:domain_circle_action} and \eqref{eq:codomain_circle_action}, respectively, and $j \in \{0,\ldots,n\}$, then the subset $A_j \subset \ZZ_{\geq 0}^{m+1}$ of multiindices $\alpha = (\alpha_0,\ldots,\alpha_m)$ corresponding to nonzero, bounded operators
\[
  D^\alpha F_j(0) \in \Hom(\otimes^\alpha W,V_j)
\]
in the induced Taylor series for $F_j = \pi_j\circ F$, where $\pi_j:V\to V_j$ is the canonical projection, obeys the following relation:
\begin{equation}
  \label{eq:codomain_weights_equal_multiindex_dot_domain_weights}
  m_j = \alpha_1l_1+\cdots+\alpha_ml_m.
\end{equation}
Moreover, the subset $A_j$ is finite when $j\geq 1$.
\item\label{item:Circle_action_on_domain_and_map_induce_circle_action_on_codomain}
  Conversely, if $V$ is a Hilbert space and the operators $D^\alpha F(0) \in \Hom(\otimes^\alpha W,V)$  have closed range and are non-zero for at most finitely many multiindices $\alpha = (\alpha_0,\alpha_1,\ldots,\alpha_m) \in \ZZ_{\geq 0}^{m+1}$ with $\alpha_k>0$ for some $k\geq 1$, then the map $F$ and the unitary representation \eqref{eq:domain_circle_action} uniquely determine a decomposition \eqref{eq:codomain_decomposition} of $V$ as an orthogonal direct sum of a maximal closed complex linear subspace $V_0\subset V$ with subset $A_0 \subset \ZZ_{\geq 0}^{m+1}$ of multiindices $\alpha = (\alpha_0,0\ldots,0)$ and integer $m_0=0$, closed complex linear subspaces $V_j\subset V$ and finite subsets $A_j \subset \ZZ_{\geq 0}^{m+1}$ for $j=1,\ldots,n$ and integers $m_1 < \cdots <m_n$ such that the following holds: The multiindices $\alpha \in A_j$ for the induced Taylor series of $F_j = \pi_j\circ F$ obey the relations \eqref{eq:codomain_weights_equal_multiindex_dot_domain_weights} for $j=0,\ldots,n$ and the map $F$ is $S^1$-equivariant with respect to the real analytic $S^1$ actions on $W$ and $V$ defined by the unitary representations \eqref{eq:domain_circle_action} and \eqref{eq:codomain_circle_action}, respectively.
\end{enumerate}
\end{prop}

\begin{exmp}
To illustrate the scope of Proposition \ref{prop:Equivariance_complex_analytic_maps_complex_Banach_spaces_circle_actions}, consider the following elementary example. Suppose $W=\CC^3=V$ with
\[
  F(x_1,x_2,x_3) = (f_0(x_1,x_2,x_3),x_1^2 + x_1x_2, x_2x_3^3),
  \quad\text{for all } (x_1,x_2,x_3) \in \CC^3,
\]
where $f_0$ is an arbitrary holomorphic function. If
\[
  e^{i\theta}\cdot(x_1,x_2,x_3) := (e^{-i\theta}x_1,e^{-i\theta}x_2,e^{i\theta}x_3),
\]
so $m=2$ with $l_1=-1$ and $l_2=1$, and
\[
  e^{i\theta}\cdot(y_1,y_2,y_3) := (y_1,e^{-i2\theta}y_2,e^{i2\theta}y_3),
\]
so $n=3$ with $m_0=0$, $m_1=-2$, and $m_2=2$, then $F:V\to W$ is $S^1$-equivariant with respect to the indicated $S^1$ actions.
\qed
\end{exmp}  

\begin{rmk}[On the hypotheses of Proposition \ref{prop:Equivariance_complex_analytic_maps_complex_Banach_spaces_circle_actions}]
\label{rmk:Equivariance_complex_analytic_maps_complex_Banach_spaces_circle_actions}  
Our assumptions that the integers $m$ and $n$ in Proposition \ref{prop:Equivariance_complex_analytic_maps_complex_Banach_spaces_circle_actions} are finite implies that, when $j\neq 0$, the operators $D^\alpha F_j(0)$ are non-zero for at most finitely many multiindices $\alpha$. These assumptions are made for convenience and consistency with our applications. Indeed, 
in our example in Feehan and Leness \cite[Section 13.6]{Feehan_Leness_introduction_virtual_morse_theory_so3_monopoles}, $W$ and $V$ are complex Banach spaces, $m=n=3$, and $F$ is a quadratic polynomial map. However, one could generalize the hypothesis \eqref{eq:domain_decomposition} to allow a countable direct sum
\begin{equation}
  \label{eq:domain_decomposition_countable}
  W = \bigoplus_{k\in\ZZ}W_k
\end{equation}
and generalize \eqref{eq:domain_circle_action} by allowing a sequence of integer weights $\cdots < l_{-2}<l_{-1}<l_0=0<l_1<l_2<\cdots$ such that
\begin{multline}
  \label{eq:domain_circle_action_countable}
  e^{i\theta}\cdot x = \left(\ldots,e^{il_{-2}\theta}x_2,e^{il_{-1}\theta}x_1,x_0,e^{il_1\theta}x_1,e^{il_2\theta}x_2,\ldots\right),
  \\
  \text{for all } x = (\ldots,x_{-2},x_{-1},x_0,x_1,x_2,\ldots) \in W \text{ and } \theta \in \RR.
\end{multline}
Similarly, in Item \eqref{item:Circle_actions_on_domain_and_codomain_constrain_map}, we may generalize the hypotheses \eqref{eq:codomain_decomposition} and \eqref{eq:codomain_circle_action} to allow
\begin{equation}
  \label{eq:codomain_decomposition_countable}
  V = \bigoplus_{j\in\ZZ} V_j
\end{equation}
with a sequence of integer weights $\cdots<m_{-2}< m_{-1} < m_0=0<m_1<m_2<\cdots$ such that
\begin{multline}
  \label{eq:codomain_circle_action_countable}
  e^{i\theta}\cdot y =
  \left(\ldots,e^{im_{-2}\theta}y_2,e^{im_{-1}\theta}y_1,y_0,e^{im_1\theta}y_1,e^{im_2\theta}y_2,\ldots\right),
  \\
  \text{for all } y = (\ldots,y_{-2},y_{-1},y_0,y_1,y_2,\ldots) \in V \text{ and } \theta \in \RR.
\end{multline}
Finally, the relation \eqref{eq:codomain_weights_equal_multiindex_dot_domain_weights} is generalized to
\begin{equation}
  \label{eq:codomain_weights_equal_multiindex_dot_domain_weights_countable}
  m_j = \sum_{k\in\ZZ}\alpha_kl_k,
\end{equation}
with the understanding that, after excluding the case $j=0$ and $m_0=0$, the multiindex sequences $(\alpha_k)_{k\in\ZZ} \in A_j$ have finite support and so the sums in \eqref{eq:codomain_weights_equal_multiindex_dot_domain_weights_countable} are finite.
\end{rmk}  

\begin{proof}[Proof of Proposition \ref{prop:Equivariance_complex_analytic_maps_complex_Banach_spaces_circle_actions}]
Consider Item \eqref{item:Circle_actions_on_domain_and_codomain_constrain_map}. The map $F_j = \pi_j\circ F:U \to V$ is complex analytic with induced Taylor series
\begin{equation}
  \label{eq:Taylor_series_component}
  F_j(x) = \sum_{|\alpha|\geq 0}\frac{1}{\alpha!}D^\alpha F_j(0)x^\alpha,
  \quad\text{for all } x \in B_r(0).
\end{equation}
Observe that by our definition \eqref{eq:codomain_circle_action} of the unitary representation $S^1\to \U(V)$ we have
\[
  e^{i\theta}\cdot F_j(x) = e^{im_j\theta}F_j(x), \quad\text{for all } x \in U \text{ and } j = 0,\ldots,n.
\]
Each operator $D^\alpha F_j(0) \in \Hom(\otimes^\alpha W, V_j)$ is complex multilinear and so, because 
\[
  e^{i\theta}\cdot x = \left(x_0,e^{il_1\theta}x_1,\ldots,e^{il_m\theta}x_m\right),
   \quad\text{for all } x = (x_0,x_1,\ldots,x_m) \in W,
\]
by definition \eqref{eq:domain_circle_action} of the unitary representation $S^1\to \U(W)$ and writing $x_k^{\alpha_k} = (x_k,\ldots,x_k)$, repeated $\alpha_k$ times for $k=0,\ldots,m$, we obtain 
\begin{align*}
  D^\alpha F_j(0)\left(e^{i\theta}\cdot x\right)^\alpha
  &=
  D^\alpha F_j(0)\left(x_0,e^{il_1\theta}x_1,\ldots,e^{il_m\theta}x_m\right)^\alpha
  \\
  &=
  D^\alpha F_j(0)\left(x_0^{\alpha_0},e^{il_1\theta}x_1^{\alpha_1},\ldots,e^{il_m\theta}x_m^{\alpha_m}\right)
  \\
  &= e^{i(\alpha_1l_1+\cdots+\alpha_ml_m)\theta}
    D^\alpha F_j(0)\left(x_0^{\alpha_0},x_1^{\alpha_1},\ldots,x_m^{\alpha_m}\right)
  \\
  &\qquad\text{(by complex multilinearity)}
  \\
  &=
    e^{im_j\theta}D^\alpha F_j(0)\left(x_0^{\alpha_0},x_1^{\alpha_1},\ldots,x_m^{\alpha_m}\right)
    \quad\text{(by \eqref{eq:codomain_weights_equal_multiindex_dot_domain_weights})},
  \\
  &\qquad \text{for all } x \in W, \alpha = (\alpha_0,\ldots,\alpha_m) \in A_j, \text{ and } j = 0,\ldots,n,
\end{align*}
and therefore we obtain
\begin{equation}
  \label{eq:Circle_equivariance_derivative_operators}
  D^\alpha F_j(0)\left(e^{i\theta}\cdot x\right)^\alpha
  =
  e^{im_j\theta}D^\alpha F_j(0)x^\alpha,
  \quad \text{for all } x \in W, \alpha \in A_j, \text{ and } j = 0,\ldots,n.
\end{equation}
Consequently, noting that the ball $B_r(0)\subset W$ is $S^1$ invariant because the representation $S^1\to\U(W)$ given by \eqref{eq:domain_circle_action} is unitary, we obtain the desired equivariance of $F:B_r(0)\to V$ by definition \eqref{eq:codomain_circle_action} of the unitary representation $S^1\to \U(V)$, since
\[
   D^\alpha F(0)\left(e^{i\theta}\cdot x\right)^\alpha
  =
  e^{i\theta}\cdot D^\alpha F(0)x^\alpha,
  \quad \text{for all } x \in W \text{ and } \alpha \in \ZZ_{\geq 0}^{m+1}, 
\]
and therefore the the Taylor series \eqref{eq:Taylor_series} for $F$ yields
\[
   F\left(e^{i\theta}\cdot x\right)
  =
  e^{i\theta}\cdot F(x),
  \quad \text{for all } x \in W\cap B_r(0).
\]  
The constraints \eqref{eq:codomain_weights_equal_multiindex_dot_domain_weights} for $j = 0,\ldots,n$ and the fact that $n<\infty$ ensure that the operators $D^\alpha F_j(0)$ can be non-zero for at most finitely many multiindices $\alpha = (\alpha_0,\alpha_1,\ldots,\alpha_m) \in \ZZ_{\geq 0}^{m+1}$ with $\alpha_k>0$ for some $k\geq 1$ and hence the subset $A_j$ is finite when $j\geq 1$. This proves Item \eqref{item:Circle_actions_on_domain_and_codomain_constrain_map}.

Consider Item \eqref{item:Circle_action_on_domain_and_map_induce_circle_action_on_codomain}. We define linear subspaces of $V$ by 
\[
  V_j := \sum_{\alpha\in A_j} D^\alpha F(0)(\otimes^\alpha W),
  \quad \text{for } j = 0,\ldots,n,
\]
where $n$ is the largest integer such that $V_j$ is a non-zero subspace. Because the operators $D^\alpha F(0)$ have closed range for all $\alpha$ by hypothesis and because the preceding sum is also finite by hypothesis, then each subspace $V_j\subset V$ is closed and uniquely determined. We label the integers $m_j$ so that $m_0=0$ while $m_1<\cdots<m_n$. If $y \in V_j\cap V_p$ for $j \neq p$, and thus $m_j \neq m_p$, then there are vectors $x, w \in W$ and multiindices $\alpha \in A_j$, $\beta \in A_p$ such that
\[
  D^\alpha F(0)x^\alpha = y = D^\beta F(0)w^\beta.
\]
But for any $e^{i\theta}\in S^1$, circle equivariance of the map $F:B_r(0)\to V$ with respect to the representations \eqref{eq:domain_circle_action} of $S^1$ on $W$ and \eqref{eq:codomain_circle_action} of $S^1$ on $V$ yields the identities
\begin{align*}
  D^\alpha F(0)\left(e^{i\theta}\cdot x\right)^\alpha &= e^{im_j\theta}D^\alpha F(0)x^\alpha = e^{im_j\theta}y,
  \\
  D^\beta F(0)\left(e^{i\theta}\cdot w\right)^\beta &= e^{im_p\theta}D^\beta F(0)w^\beta = e^{im_p\theta}y,
\end{align*}
and because $m_j\neq m_p$, we must have $y=0$. Therefore, $V_j\cap V_p = (0)$ for all $j\neq p$ and so the linear sum of the subspaces $V_0,\ldots,V_n$ is a direct sum and a closed subspace of $V$,
\[
  \sum_{j=0}^n V_j = \bigoplus_{j=0}^n V_j \subseteq V.
\]
If the preceding direct sum is equal to $V$, we are done. If not equal to $V$, we proceed as follows. Suppose that $P, Q \subset V$ are closed subspaces with orthogonal complements $P^\perp, Q^\perp \subset V$, so $V = P\oplus P^\perp = Q\oplus Q^\perp$ as direct sums of Hilbert spaces. Then, $P+Q\subset V$ is a closed subspace with orthogonal complement $(P+Q)^\perp = P^\perp\cap Q^\perp$ and thus
\[
  V = (P+Q) \oplus\left(P^\perp\cap Q^\perp\right).
\]
By an extension of this principle to finite sums of closed linear subspaces and the fact that each subspace $V_j$ is closed in $V$ for $j=0,\ldots,n$ with orthogonal complement $V_j^\perp$, we may write
\[
  V = R\oplus\bigoplus_{j=0}^n V_j,
\]
where $R := V_0^\perp\cap \cdots\cap V_n^\perp$ is the orthogonal complement of the closed subspace $\oplus_{j=0}^n V_j$ in $V$. By assumption for this case, $R\neq(0)$ and we define the $S^1$ action on $R$ to be trivial, that is, $e^{i\theta}\cdot y := y$ for all $y\in R$, and relabel $V_0\oplus R$ as $V_0$, so that we obtain the orthogonal direct sum decomposition
\[
  V = \bigoplus_{j=0}^n V_j,
\]
as desired for \eqref{eq:codomain_decomposition}. This completes the proof of Item \eqref{item:Circle_action_on_domain_and_map_induce_circle_action_on_codomain} and thus Proposition \ref{prop:Equivariance_complex_analytic_maps_complex_Banach_spaces_circle_actions}.
\end{proof}

Proposition \ref{prop:Equivariance_complex_analytic_maps_complex_Banach_spaces_circle_actions} leads to the

\begin{cor}[Circle-equivariant complex analytic maps of complex Banach spaces preserve weight-sign decompositions]
\label{cor:Circle-equivariant_complex_analytic_maps_preserve_weight-sign_decompositions}  
Let $W$ and $V$ be complex Banach spaces endowed with real analytic circle actions having weight decompositions defined by the representations \eqref{eq:domain_decomposition}, \eqref{eq:domain_circle_action} and \eqref{eq:codomain_decomposition}, \eqref{eq:codomain_circle_action}, respectively, $U \subset W$ be an $S^1$-invariant open neighborhood of the origin, and $F:W\supset U\to V$ be a circle-equivariant, complex analytic map. If
\begin{subequations}
\label{eq:Weight-sign_decomposition_W_and_V}  
\begin{align}
  \label{eq:Weight-sign_decomposition_W}
  W &=W_0 \oplus W^+ \oplus W^-, \quad\text{where}\quad W^+ := \bigoplus_{\{k:l_k>0\}} W_k
  \quad\text{and}\quad W^- := \bigoplus_{\{k:l_k<0\}} W_k,
  \\
  \label{eq:Weight-sign_decomposition_V}
  V &= V_0 \oplus V^+ \oplus V^-, \quad\text{where}\quad V^+ := \bigoplus_{\{j:m_j>0\}} V_j \quad\text{and}\quad V^- := \bigoplus_{\{j:m_j<0\}} V_j,
\end{align}
\end{subequations}
are the weight-sign decompositions of $W$ and $V$, respectively, then
  \begin{equation}
  \label{eq:FW0pm_subset_V0pm_complex}
  F(W_0\cap B_r(0)) \subset V_0 \quad\text{and}\quad F(W^\pm\cap B_r(0)) \subset V^\pm,
\end{equation}
where $r$ is the radius of convergence of the Taylor series of $F$ around the origin.
\end{cor}

\begin{proof}
By the calculations leading to the identity \eqref{eq:Circle_equivariance_derivative_operators}, we see that if $x \in W_0\cap B_r(0)$, so $x = (x_0,0,\ldots,0)$ for $x_0\in W_0$, then
\begin{multline*}
  D^\alpha F_j(0)x^\alpha
  =
  D^\alpha F_j(0)\left(e^{i\theta}\cdot x\right)^\alpha
  =
  e^{im_j\theta}D^\alpha F_j(0)x^\alpha,
  \\
  \text{for all } x \in W_0\cap B_r(0), \theta\in\RR, \text{ and } j=0,\ldots,n.
\end{multline*}
The preceding identity is a tautology when $j=0$ and thus $m_0=0$, but when $j=1,\ldots,n$, it forces $D^\alpha F_j(0)x^\alpha = 0$ for all $\alpha \in A_j$ and $x\in W_0\cap B_r(0)$. Therefore, $F_j(x) = 0$ for all $x\in W_0\cap B_r(0)$ and $j=1,\ldots,n$ by the Taylor series \eqref{eq:Taylor_series_component} for $F_j(x)$. This yields the first inclusion in \eqref{eq:FW0pm_subset_V0pm_complex}.

For the remaining two cases, suppose first that the non-zero weights $l_k$ are either all positive or all negative for $k=1,\ldots,m$. Then the relations \eqref{eq:codomain_weights_equal_multiindex_dot_domain_weights} imply that the non-zero weights $m_j$ are either all positive or all negative as well for $j=1,\ldots,n$. Consequently, $F(x) \in V^\pm$ for all $x \in W^\pm\cap B_r(0)$ and this yields the second pair of inclusions in \eqref{eq:FW0pm_subset_V0pm_complex}, in the special case where $W^+=(0)$ or $W^-=(0)$.

In general, the non-zero weights $l_k$ do not all have the same sign if $m\geq 2$ and so there is an integer $q\in\{2,\ldots,m-1\}$ such that $l_1<\cdots<l_q<0<l_{q+1}<\cdots<l_m$. If $x \in W^-\cap B_r(0)$, so $x = (0,x_1,\ldots,x_q,0,\ldots,0)$ with $x_k\in W_k$ for $k=1,\ldots,q$, then the calculations leading to \eqref{eq:Circle_equivariance_derivative_operators} imply that
\begin{multline*}
  D^\alpha F_j(0)\left(e^{i\theta}\cdot x\right)^\alpha
  =
  D^\alpha F_j(0)\left(0,e^{il_1\theta}x_1,\ldots,e^{il_q\theta}x_q,0,\ldots,0\right)^\alpha
  \\
  =
  e^{i(\alpha_1l_1+\cdots+\alpha_ql_q)\theta}D^\alpha F_j(0)x^\alpha
  =
  e^{im_j\theta}D^\alpha F_j(0)x^\alpha,
  \\
  \text{for all } x \in W^-\cap B_r(0), \theta\in\RR, \text{ and } j=0,\ldots,n.
\end{multline*}
Because $l_k<0$ for $k=1,\ldots,q$, then $\alpha_1l_1+\cdots+\alpha_ql_q < 0$ for all $\alpha \in A_j$ and $j=0,\ldots,n$. Consequently, either $F_j(x)=0$ for $j=1,\ldots,n$ and all $x \in W^-\cap B_r(0)$ or there is a least integer $s\in\{1,\ldots,n-1\}$ such that
\[
  D^\alpha F_j(0)x^\alpha = 0, \quad\text{for all } x \in W^-\cap B_r(0), 
  \quad j=s+1,\ldots,n,
\]
while, for each $j=1,\ldots,s$, we may have $D^\alpha F_j(0)x^\alpha \neq 0$ for some $x \in W^-\cap B_r(0)$ and $\alpha\in A_j$ and
\[
  m_j = \alpha_1l_1+\cdots+\alpha_ql_q < 0.
\]
Consequently,
\[
  V^- = \bigoplus_{\{j:m_j<0\}} V_j = \bigoplus_{j=1}^s V_j,
\]
and $F_j(x) \in V^-$ for all $x\in W^-\cap B_r(0)$ and $j=1,\ldots,s$, while $F_j(x) = 0$ for all $x\in W^-\cap B_r(0)$ and $j=s+1,\ldots,n$ by the Taylor series \eqref{eq:Taylor_series_component} for $F_j(x)$. This implies that $F(x) \in V^-$ for all $x\in W^-\cap B_r(0)$ and hence yields the inclusion \eqref{eq:FW0pm_subset_V0pm_complex} into $V^-$. An essentially identical argument yields the inclusion \eqref{eq:FW0pm_subset_V0pm_complex} into $V^+$.
\end{proof}

We now return to state and prove the special case of Blanchard's Lemma \ref{lem:Blanchard_complex_analytic_spaces} promised at the beginning of this section.

\begin{lem}[Linear circle action induced by a holomorphic map on neighborhood of fixed point]
\label{lem:Blanchard_holomorphic_map_from_domain_around_fixed_point_into_vector_space}
Let $n,r$ be positive integers, $D \subset \CC^n$ be a domain, and $F:\CC^n\supset D \to \CC^r$ be a holomorphic map. If $\Phi:S^1 \to \Aut(D)$ is a homomorphism from $S^1$ into the group $\Aut(D)$ of biholomorphic transformations of $D$ that defines a real analytic action $S^1\times D\to D$ and the origin $0 \in D$ is a fixed point of this $S^1$ action then, after possibly shrinking $D$, the following hold:
\begin{enumerate}
\item
\label{item:Bochner_Cartan_Kaup_circle}
There are an open neighborhood $V$ of the origin in $\CC^n$ and a biholomorphic map $\chi:D \to V$ such that $\chi(0) = 0$ and $D\chi(0) = \id \in \U(n)$ and 
\[
  \chi(e^{i\theta}\cdot z) = D(\Phi(e^{i\theta}))(0)\chi(z),
  \quad\text{for all } e^{i\theta} \in S^1 \text{ and } z \in D,
\]
where $\rho:S^1 \ni e^{i\theta} \mapsto D(\Phi(e^{i\theta}))(0) \in \U(n)$ is the isotropy representation of $S^1$. In particular, $V$ is invariant under the linear $S^1$ action on $\CC^n$ defined by $\rho:S^1 \to \U(n)$ and the map $\chi$ is $S^1$-equivariant with respect to the possibly nonlinear $S^1$ action on $D$ and linear $S^1$ action on $V$.
\item
\label{item:Blanchard_linear_circle_action}
The induced holomorphic map $F\circ \chi^{-1}:V \to \CC^r$ uniquely determines a linear $S^1$ action on $\CC^r$ defined by a unitary representation $\varrho:S^1 \to \U(r)$ such that $F\circ \chi^{-1}$ is $S^1$ equivariant.
\item
\label{item:Blanchard_original_nonlinear_circle_action}
The holomorphic map $F:D \to \CC^r$ uniquely determines an $S^1$ action on $\CC^r$ defined by a unitary representation $\varrho:S^1 \to \U(r)$ such that $F$ is $S^1$-equivariant with respect to the given nonlinear $S^1$ action on $D$ and induced linear $S^1$ action on $\CC^r$.
\end{enumerate}  
\end{lem}

\begin{proof}
Item \eqref{item:Bochner_Cartan_Kaup_circle} follows from the Linearization Theorem \ref{thm:Cartan_linearization_theorem} or \ref{thm:Kaup_linearization_theorem} (due to Cartan or Kaup, respectively) and the implication that the complex representation $\rho:S^1\to\GL(n,\CC)$ provided by Theorem \ref{thm:Cartan_linearization_theorem} or \ref{thm:Kaup_linearization_theorem} is actually unitary, so $\rho:S^1\to\U(n)$, follows from the classification of complex, linear representations of $S^1$ provided by Br\"ocker and tom Dieck \cite[Chapter II, Proposition 1.9, p. 68, and Proposition 8.1, p. 107]{BrockertomDieck}. Item \eqref{item:Blanchard_linear_circle_action} of the lemma follows from Item \eqref{item:Circle_action_on_domain_and_map_induce_circle_action_on_codomain} of Proposition \ref{prop:Equivariance_complex_analytic_maps_complex_Banach_spaces_circle_actions}. Finally, Item \eqref{item:Blanchard_original_nonlinear_circle_action} is an immediate corollary of Items \eqref{item:Bochner_Cartan_Kaup_circle} and \eqref{item:Blanchard_linear_circle_action}.
\end{proof}

In Section \ref{sec:Analytic_group_actions_analytic_spaces}, we recalled the definition of a (global) group action on a $\KK$-analytic space but for applications, it is useful to consider a stronger concept that is motivated by many examples in gauge theory and by Lemma \ref{lem:Blanchard_holomorphic_map_from_domain_around_fixed_point_into_vector_space} when $G = S^1$.

\section{Real analytic  maps of complex  vector spaces with circle actions}
\label{sec:Real_analytic_maps_vector_spaces_with_circle_actions}
We now describe what can be said in the more delicate case of a real analytic map. For simplicity, we restrict our attention to the case of orthogonal representations of the circle group on finite-dimensional, real inner product spaces --- see Br\"ocker and tom Dieck \cite[Chapter 2, Proposition 8.5, p. 109]{BrockertomDieck} for such representations --- and based on that reference, we begin with the important

\begin{rmk}[Weight-sign ambiguity for orthogonal representations of $S^1$]
\label{rmk:Weight_sign_ambiguity_real_representations_circle}
Given an orthogonal representation $\rho:S^1\to\Or(W)$ with weight $+1$ on a real, two-dimensional, inner product space $W=\RR^2\cong\CC$, we may define an orthogonal representation $\bar\rho:S^1\to\Or(W)$ with weight $-1$ by setting
\[
  \bar\rho(e^{i\theta})\bar z := e^{-i\theta}\bar z = \overline{e^{i\theta}z}
  = \overline{\rho(e^{i\theta})z}, \quad\text{for all } z \in W.
\]
Complex conjugation $C:W\to W$ is thus an isometric isomorphism of real inner product spaces that obeys
\[
  \left(C\circ\rho(e^{i\theta})\right)(z) =  \left(\bar\rho(e^{i\theta})\circ C\right)(z),
  \quad\text{for all } z \in W.
\]
Hence, the orthogonal representations $\rho:S^1\to\Or(W)$ with weight $+1$ and $\bar\rho:S^1\to\Or(W)$ with weight $-1$ are isomorphic as \emph{real representations of $S^1$} by Br\"ocker and tom Dieck \cite[Chapter 2, Definition 1.4, p. 67]{BrockertomDieck}, since they are intertwined by the isomorphism $C:W\to W$ of real vector spaces.
\end{rmk}

Remark \ref{rmk:Weight_sign_ambiguity_real_representations_circle} shows that the weights of an orthogonal representation of $S^1$ on a real vector space are only determined up to sign by the representation. Moreover, this remark yields an example of a real analytic map $F:\RR^2\to\RR^2$ that is circle-equivariant with respect to a circle action $\rho_+:S^1\to\Or(2)$ on the domain $\RR^2$ with weight $+1$ and a circle action $\rho_-:S^1\to\Or(2)$ on the codomain $\RR^2$ with weight $-1$,
\[
  \left(F\circ\rho_+(e^{i\theta})\right)(x,y) = \left(\rho_-(e^{i\theta})\circ F\right)(x,y),
  \quad\text{for all } (x,y) \in \RR^2.
\]
and so this map $F$ does \emph{not} preserve signs of weight of the circle actions. As in Remark \ref{rmk:Weight_sign_ambiguity_real_representations_circle}, we may consider $\RR^2 = \CC$ as a real vector space and define $F(z) := \bar z$ (complex conjugation) for all $z\in\CC$ and define $\rho_\pm(e^{i\theta})z := e^{\pm i\theta}z$ (complex scalar multiplication) for all $\theta\in\RR$ and $z\in\CC$. Thus, unlike holomorphicity in Section \ref{sec:Holomorphic_maps_vector_spaces_with_circle_actions}, real analyticity of a map $F$ is insufficient to ensure that signs of weights of circle actions are preserved. We recall the

\begin{defn}[Pseudoholomorphic map]
\label{defn:Pseudoholomorphic_map}
(See Boothby, Kobayashi, and Wang \cite[Section 2, p. 329]{Boothby_Kobayashi_Wang_1963}, Cirici and Wilson \cite[Definition 2.2, p. 11]{Cirici_Wilson_2021}, Di Scala, Kasuya, and Zuddas \cite[Section 2, p. 3, Definition 7]{DiScala_Kasuya_Zuddas_2016}, or McDuff and Salamon \cite[Section 4.5, p. 180]{McDuffSalamonSympTop3}.)  
Let $(M, J)$ and $(N, J')$ be smooth, almost complex manifolds. A smooth map $F : M \to N$
is \emph{pseudoholomorphic} (or \emph{almost analytic} or a \emph{morphism of almost complex manifolds}) if $\d F\circ J = J'\circ\d F:TM \to TN$. Equivalently, $F$ is pseudo-holomorphic (or almost analytic or morphism of almost complex manifolds) if and only if the differential $\d F:TM \to TN$ is complex linear at each point. An \emph{isomorphism of almost complex manifolds} is a morphism of almost complex manifolds that is also a diffeomorphsim of smooth manifolds.
\end{defn}

\begin{rmk}[Pseudoholomorphic maps that are not holomorphic]
\label{rmk:Pseudoholomorphic_map_not_holomorphic}
Fern\'andez \cite{Fernandez_2015} provides many explicit examples of pseudoholomorphic maps $F:S^2 \to S^6$ that are not holomorphic when $S^6$ is equipped with its standard round metric and an orthogonal almost complex structure constructed with the aid of the octonions \cite[Section 1, p. 2437]{Fernandez_2015}. According to Blanchard \cite{Blanchard_1953} and LeBrun \cite{LeBrun_1987}, there are no orthogonal complex structures on $S^6$ and so the orthogonal almost complex structure employed by Fern\'andez cannot be integrable. Therefore, the examples constructed by Fern\'andez are strictly pseudoholomorphic in the sense that the codomain is strictly almost complex, with an almost complex structure that is not integrable. Fern\'andez employs techniques from twistor theory and its relation to harmonic maps, for which a general reference is Burstall and Rawnsley \cite{Burstall_Rawnsley_twistor_theory_riemannian_symmetric_spaces}.
\end{rmk}

Suppose now that we choose in addition an almost complex structure,
\[
  J = \begin{pmatrix} 0 & 1 \\ -1 & 0\end{pmatrix} \in \End(\RR^2),
\]
on the domain and codomain $\RR^2$, so $(\RR^2,J) \cong \CC$ with $e_2 = Je_1$ and $xe_1 + ye_2 = (x + yJ)e_1$, for all $x,y\in\RR$, and the indicated real linear isomorphism is defined by $\RR^2 \ni (x + yJ)e_1 \mapsto x+iy \in \CC$. If $F:\RR^2\to\RR^2$ is a pseudoholomorphic map in the sense of Definition \ref{defn:Pseudoholomorphic_map}, then
\[
  \d F(p)\circ J_p = J_p\circ \d F(p), \quad\text{for all } p \in \RR^2,
\]
where $T_p\RR^2 = \RR^2$ and $J_p = J$ for all $p\in\RR^2$. With respect to the identification $(\RR^2,J) \cong \CC$, the differential $\d F(p):\CC\to\CC$ is complex linear (since $F:\RR^2\to\RR^2$ is a pseudoholomorphic map) and thus
\[
  \d F(p)(e^{i\theta}v) = e^{i\theta}\d F(p)v,
  \quad\text{for all } p \in \RR^2, v \in T_p\RR^2 = \RR^2, \text{ and } e^{i\theta} \in S^1 \subset \CC^*,
\]
so $\d F(p)$ preserves the sign of the weight (equal to $+1$ in this example) of the circle action on $\RR^2$.

Given finite-dimensional, real Hilbert spaces $W$ and $V$ and orthogonal representations $\rho_W:S^1\to\Or(W)$ and $\rho_V:S^1\to\Or(V)$, the preceding simple examples suggest that we seek partial analogues for circle-equivariant, pseudoholomorphic, smooth maps $F:W\supset U \to V$ of the results in Section \ref{sec:Holomorphic_maps_vector_spaces_with_circle_actions} for real Hilbert spaces with circle-invariant almost complex structures.

\begin{prop}[Equivariance of pseudoholomorphic maps of almost complex vector spaces with circle actions]
\label{prop:Equivariance_real_analytic_maps_real_vector_spaces_circle_actions}
Let $W$ and $V$ be finite-dimensional, real Hilbert spaces with orthogonal almost complex structures and $S^1\times W \to W$ be a $C^1$ circle action defined by an orthogonal representation $S^1\to\Or(W)$ given by an orthogonal direct sum decomposition \eqref{eq:domain_decomposition} and integer weights $l_0=0$ and\footnote{Not assumed to be strictly increasing as in Proposition \ref{prop:Equivariance_complex_analytic_maps_complex_Banach_spaces_circle_actions}.} $l_1\leq l_2\leq \cdots\leq l_m$ such that the almost complex structure $J_W$ on $W$ is circle invariant. Each $W_k$ is assumed to be a real, two-dimensional, linear subspace for $k=1,\ldots,m$ and we assume that $J_W(W_k) \subset W_k$, so $W_k$ has an almost complex structure $J_{W_k} = J_W\restriction W_k$, for $k=1,\ldots,m$. We assume that $e^{i\theta}\in S^1$ acts by real matrix multiplication,
\begin{multline}
  \label{eq:domain_circle_action_real}
  e^{i\theta}\cdot x
  =
  \left(x_0,
    \begin{pmatrix}\cos l_1\theta & -\sin l_1\theta \\ \sin l_1\theta &\cos l_1\theta\end{pmatrix}
    \begin{pmatrix} x_1^1\\x_1^2 \end{pmatrix},
    \ldots,
    \begin{pmatrix}\cos l_m\theta & -\sin l_m\theta \\ \sin l_m\theta &\cos l_m\theta\end{pmatrix}
    \begin{pmatrix} x_m^1\\x_m^2 \end{pmatrix}
  \right),
  \\
  \text{for all } x = (x_0,x_1,\ldots,x_m) \in W \text{ and } \theta \in \RR,
 \end{multline}
where $x_k = x_k^1e_1^k + x_k^2e_2^k \in W_k$ with respect to a basis $\{e_1^k,e_2^k\}$ for $W_k$ with $e_2^k = J_We_1^k$, for $k=1,\ldots,m$. If $U \subset W$ is an $S^1$-invariant open neighborhood of the origin and $F:W\supset U\to V$ is a $C^1$ pseudoholomorphic map in the sense of Definition \ref{defn:Pseudoholomorphic_map}, then the following hold.
\begin{enumerate}
\item\label{item:Circle_actions_on_domain_and_codomain_constrain_map_real}
If $J_V$ is the almost complex structure on $V$, then we assume that $J_V$ is circle invariant and $J_V(V_j) \subset V_j$, so $V_j$ has an almost complex structure $J_{V_j} = J_V\restriction V_j$, for $j=0,\ldots,n$. Assume that $S^1\times V \to V$ is a $C^1$ circle action defined by an orthogonal representation $S^1\to\Or(V)$ given by an orthogonal direct sum decomposition \eqref{eq:codomain_decomposition} and integer weights $m_0=0$ and\footnote{Not assumed to be strictly increasing as in Proposition \ref{prop:Equivariance_complex_analytic_maps_complex_Banach_spaces_circle_actions}.} $m_1 \leq m_2 \leq \cdots \leq m_n$. Each $V_j$ is now a real, two-dimensional, linear subspace for $j=1,\ldots,n$ and $e^{i\theta}\in S^1$ acts by real matrix multiplication,
\begin{multline}
  \label{eq:codomain_circle_action_real}
  e^{i\theta}\cdot y
  =
  \left(y_0,
    \begin{pmatrix}\cos m_1\theta & -\sin m_1\theta \\ \sin m_1\theta &\cos m_1\theta\end{pmatrix}
    \begin{pmatrix} y_1^1\\y_1^2 \end{pmatrix},
    \ldots,
    \begin{pmatrix}\cos m_n\theta & -\sin m_n\theta \\ \sin m_n\theta &\cos m_n\theta\end{pmatrix}
    \begin{pmatrix} y_n^1\\y_n^2 \end{pmatrix}
  \right),
  \\
  \text{for all } y = (y_0,y_1,\ldots,y_n) \in V \text{ and } \theta \in \RR,
 \end{multline}
where $y_j = y_j^1f_1^j + y_j^2f_2^j \in V_j$ with respect to a basis $\{f_1^j,f_2^j\}$ for $V_j$ with $f_2^j = J_Vf_1^j$, for $j=1,\ldots,n$. If $F$ is circle-equivariant with respect to the circle actions on $W$ and $V$ defined by the orthogonal representations \eqref{eq:domain_circle_action_real} and \eqref{eq:codomain_circle_action_real}, respectively, and the partial derivative $D_kF_j(0) \in\Hom(W_k,V_j)$ is nonzero for some $k\in\{1,\ldots,m\}$ and $j \in \{1,\ldots,n\}$, then
\begin{equation}
   \label{eq:codomain_weight_equals_domain_weight_real}
   m_j = l_k,
\end{equation}
where $F_j := \pi_j\circ F$ and $\pi_j:V\to V_j$ is the canonical projection.
 
\item\label{item:Circle_action_on_domain_and_map_induce_circle_action_on_codomain_real}
Conversely, if the operator $DF(0) \in \Hom(W,V)$ is surjective then the orthogonal representation \eqref{eq:domain_circle_action_real} of $S^1$ on $W$, orthogonal almost complex structures $J_W$ and $J_V$, and map $F$ uniquely determine an orthogonal representation $\rho_V$ of $S^1$ on $V$ such that the relations \eqref{eq:codomain_weight_equals_domain_weight_real} hold for all $k\in\{1,\ldots,m\}$ and $j \in \{1,\ldots,n\}$ with $D_kF_j(0) \in\Hom(W_k,V_j)$ nonzero. 
\end{enumerate}
\end{prop}

\begin{proof}
Because $F:W\supset U\to V$ is equivariant with respect to linear $S^1$ actions on $W$ and $V$, then
\[
  DF(0)\left(e^{i\theta}\cdot x\right) = e^{i\theta}\cdot DF(0)x,
  \quad\text{for all } e^{i\theta} \in S^1 \text{ and } x \in W,
\]
that is, the operator $DF(0) \in \Hom(W,V)$ is also is equivariant with respect to the $S^1$ actions on $W$ and $V$. For $k=1,\ldots,m$, there are real linear isomorphisms
\[
  \iota_k:W_k \ni x_k = x_k^1e_1^k +  x_k^2Je_1^k \mapsto x_k^1 + ix_k^2 \in \CC.
\]
Each $e^{i\theta}\in S^1$ acts by real matrix multiplication on $W_k$ with respect to its basis $\{e_1^k, Je_1^k\}$ via \eqref{eq:domain_circle_action_real}, so that
\begin{multline*}
  e^{i\theta}\cdot x_k
  := \left(x_k^1\cos l_k\theta - x_k^2\sin l_k\theta\right)e_1^k
  + \left(x_k^1\sin l_k\theta + x_k^2\cos l_k\theta\right)e_2^k,
  \\
  \text{for all }  x_k = x_k^1e_1^k +  x_k^2Je_1^k \in W_k.
\end{multline*}
If we write $z_k = x_k^1 + ix_k^2 \in \iota_k(W_k) = \CC$, then
\begin{multline*}
  e^{il_k\theta}z_k
  = (\cos l_k\theta + i\sin l_k\theta)(x_k^1 + ix_k^2)
  \\
  = (x_k^1\cos l_k\theta - x_k^2\sin l_k\theta) + i(x_k^1\sin l_k\theta + x_k^2\cos l_k\theta)
 \in \CC.
\end{multline*}
In other words,
\[
  \iota_k\left(e^{i\theta}\cdot x\right) = e^{il_k\theta}\iota_k(x_k),
  \quad\text{for all }
  e^{i\theta} \in S^1 \text{ and } x_k\in W_k \text{ and } k=1,\ldots,m,
\]  
and so the following diagram commutes, with the top horizontal arrow indicating multiplication by real, two-by-two matrices and bottom horizontal arrow indicating complex scalar multiplication
\[
\begin{tikzcd}
    W_k \arrow[r, "e^{i\theta}\cdot"] \arrow[d, "\iota_k"'] &W_k \arrow[d, "\iota_k"]
    \\
    \CC \arrow[r, "e^{il_k\theta}\times"] &\CC
  \end{tikzcd}
\]
for $k=1,\ldots,m$. Similarly, the following diagram commutes, with the top horizontal arrow indicating multiplication by real, two-by-two matrices and bottom horizontal arrow indicating complex scalar multiplication: 
\[
\begin{tikzcd}
    V_j \arrow[r, "e^{i\theta}\cdot"] \arrow[d, "\iota_j"'] &V_j \arrow[d, "\iota_j"]
    \\
    \CC \arrow[r, "e^{im_j\theta}\times"] &\CC
  \end{tikzcd}
\]  
for $j=1,\ldots,n$. Suppose $D_kF_j(0) \in\Hom(W_k,V_j)$ is nonzero for some $k\in\{1,\ldots,m\}$ and $j \in \{1,\ldots,n\}$. The operator $DF(0) \in\Hom(W,V)$ is complex linear with respect to the almost complex structures $J_W$ on $W$ and $J_V$ on $V$, so the operator $D_kF_j(0) \in\Hom(W_k,V_j)$ is complex linear with respect to the almost complex structures $J_{W_k}$ on $W_k$ and $J_{V_j}$ on $V_j$. Thus
\[
  \iota_j\circ D_kF_j(0)\circ\iota_k^{-1} \in \End(\CC)
\]
is complex linear in the usual sense and hence $S^1$ equivariance,
\[
  \left(\iota_j\circ D_kF_j(0)\circ\iota_k^{-1}\right)\left(e^{il_k\theta}z_k\right)
  =
   e^{im_j\theta}\left(\iota_j\circ D_kF_j(0)\circ\iota_k^{-1}\right)z_k,
\]
yields the identity $m_j = l_k$ in \eqref{eq:codomain_weight_equals_domain_weight_real}. This yields Item \eqref{item:Circle_actions_on_domain_and_codomain_constrain_map_real}, while Item \eqref{item:Circle_action_on_domain_and_map_induce_circle_action_on_codomain_real} is immediate. This completes the proof of Proposition \ref{prop:Equivariance_real_analytic_maps_real_vector_spaces_circle_actions}.
\end{proof}

Item \eqref{item:Circle_actions_on_domain_and_codomain_constrain_map_real} in Proposition \ref{prop:Equivariance_real_analytic_maps_real_vector_spaces_circle_actions} immediately leads to the following partial analogue of Corollary \ref{cor:Circle-equivariant_complex_analytic_maps_preserve_weight-sign_decompositions}.

\begin{cor}[Differentials of circle-equivariant, pseudoholomorphic maps of almost complex vector spaces preserve weight-sign decompositions]
\label{cor:Circle-equivariant_real_analytic_maps_preserve_weight-sign_decompositions}
Let $W$ and $V$ be finite-dimensional, real Hilbert spaces endowed with circle actions defined by the orthogonal representations \eqref{eq:domain_decomposition}, \eqref{eq:domain_circle_action_real} on $W$ and \eqref{eq:codomain_decomposition}, \eqref{eq:codomain_circle_action_real} on $V$, respectively,
with integer weights $l_k$ associated to subspaces $W_k\subset W$, for $k=0,\ldots,m$, and integer weights $m_j$ associated to subspaces $V_j\subset V$, for $j=0,\ldots,n$. Assume further that $W$ and $V$ are endowed with circle-invariant almost complex structures that induce almost complex structures on subspaces $W_k\subset W$, for $k=1,\ldots,m$, and $V_j\subset V$, for $j=1,\ldots,n$, that are assumed to be of real dimensional two. If $U \subset W$ is an circle-invariant, open neighborhood of the origin and $F:W\supset U\to V$ is a $C^1$ circle-equivariant map that is pseudoholomorphic in the sense of Definition \ref{defn:Pseudoholomorphic_map} and the weight-sign decompositions of $W$ and $V$ are as in \eqref{eq:Weight-sign_decomposition_W_and_V}, then
\begin{equation}
\label{eq:FW0pm_subset_V0pm_real}
DF(p)W_0 \subset V_0 \quad\text{and}\quad DF(p)W^\pm \subset V^\pm.
\quad\text{for all } p \in U,
\end{equation}
\end{cor}

\chapter[Bia{\l}ynicki--Birula decompositions for algebraic varieties and complex manifolds]{Bia{\l}ynicki--Birula decompositions for smooth algebraic varieties and complex manifolds}
\label{chap:BB_decomposition_algebraic_variety_scheme_complex_manifold}
The Bia{\l}ynicki--Birula decomposition for a smooth algebraic variety or complex manifold provides the key geometric method in our work. We begin in Section \ref{sec:BB_decomposition_algebraic_variety_scheme} by reviewing the classical Bia{\l}ynicki--Birula decomposition for a $\GG_m$ action on a smooth algebraic variety. Section \ref{sec:BB_decomposition_complex_projective_space} provides a development of the Bia{\l}ynicki--Birula decompositions for $\CC^*$ actions on complex projective space. In Section \ref{sec:BB_decomposition_complex_manifold} we review the Bia{\l}ynicki--Birula decomposition for a $\CC^*$ action on a compact, complex K\"ahler manifold. Section \ref{sec:BB_decomposition_noncompact_complex_manifold} contains our results on the Bia{\l}ynicki--Birula decomposition for a $\CC^*$ action on a noncompact, complex manifold. In Section \ref{sec:Functorial_properties_Bialynicki-Birula_decomposition}, we discuss functorial properties of the Bia{\l}ynicki--Birula decomposition. We conclude in Section \ref{sec:Hamiltonian_functions_circle_actions} by discussing two elementary examples of a Hamiltonian function, one for a circle action on a complex vector space and another for a circle action on a complex projective space.

\section{Bia{\l}ynicki--Birula decompositions for $\GG_m$ actions on smooth algebraic varieties}
\label{sec:BB_decomposition_algebraic_variety_scheme}
Before stating the Bia{\l}ynicki--Birula decomposition, we review some definitions. Recall that an algebraic group $G$ over $k$ is a \emph{torus} if it becomes isomorphic to a product of copies of $\GG_m$ over some field containing $k$ and is \emph{split} if it is isomorphic to a product of copies of $\GG_m$ over $k$ (see Milne \cite[Section 12e, Definition 12.14, p. 236]{Milne_algebraic_groups_cup}).

\begin{defn}
\label{defn:Milne_13-5}  
(See Milne \cite[Definition 13.45, p. 273]{Milne_algebraic_groups_cup}.)   
An action of a torus $T$ on a scheme $X$ over $k$ is \emph{locally affine} if $X$ admits a covering by $T$-invariant open affine subschemes.
\end{defn}

\begin{defn}
\label{defn:Milne_7-26}  
(See Milne \cite[Definition 7.26, p. 145]{Milne_algebraic_groups_cup}.)   
An action of $G$ on an algebraic scheme $X$ over $k$ is said to be \emph{linear} if there exists a representation $r:G \to \GL(V)$ of $G$ on a finite-dimensional vector space $V$ and an equivariant immersion $X \hookrightarrow \PP(V)$.
\end{defn}

Let $X$ be a scheme equipped with an action of $\GG_m$ and let $p \in X(k)$, where $X(k)$ denotes the set of points $p \in |X|$ such that $\kappa(p)=k$, where $\kappa(p) = \sO_{X,p}/\fm_p$ is the residue field (see Milne \cite[Section 1a, p. 9]{Milne_algebraic_groups_cup}). If $p$ is fixed by $\GG_m$, then $\GG_m$ acts on the tangent space $T_pX$, which therefore decomposes into a direct sum
\[
  T_pX = \bigoplus_{i\in\ZZ} (T_pX)_i
\]
of eigenspaces, so $t \in T(k)$ acts on $(T_pX)_i$ as multiplication by $t^i$. Let
\begin{align*}
  T_p^+X &:= \bigoplus_{i>0}(T_pX)_i,
  \\
  T_p^-X &:= \bigoplus_{i<0}(T_pX)_i,
  \\
  T_p^0X &:= (T_pX)_0,
\end{align*}  
so that
\[
  T_pX = T_p^+X \oplus T_p^-X \oplus T_p^0X.
\]
Given a field $k$, let $k^a$ denote its algebraic closure and for an algebraic scheme $(X,\sO_X)$ over $k$, one lets $X$ denote the scheme and $|X|$ the underlying topological space of closed
points (see Milne \cite[Introduction, p. 3]{Milne_algebraic_groups_cup}). We then have the

\begin{thm}[Bia{\l}ynicki--Birula decomposition]
\label{thm:Milne_13-47}  
(See Milne \cite[Theorem 13.47, p. 272]{Milne_algebraic_groups_cup} or Bia{\l}ynicki--Birula \cite[Theorem 4.1, p. 492]{Bialynicki_1973} for the original statement.)
Let $X$ be a smooth algebraic variety over $k$ equipped with a locally affine action of $\GG_m$.
\begin{enumerate}
\item Let $Z$ be a connected component of $X^{\GG_m}$. There exist a unique smooth subvariety $X(Z)$ of $X$ such that
\[
  X(Z)(k^a) = \left\{p \in X(k^a):\lim_{t\to 0}tp \text{ exists and lies in } Z(k^a)\right\}
\]
and a unique regular map $\pi_Z:X(Z) \to Z$ sending $p \in X(Z)(k^a)$ to the limit $\lim_{t\to 0}tp \in Z(k^a)$.

\item The map $\pi_Z$ realizes $X(Z)$ as a fiber bundle over $Z$. More precisely, every point $p \in Z(k)$ has an open neighborhood $U$ such that the restriction of $\pi_Z$ to $\pi_Z^{-1}(U)$ is isomorphic over $U$ to the projection $U \times (T_p^+X)_\fa \to U$.

\item The topological space $|X|$ is a disjoint union of the locally closed subsets $|X(Z)|$ as $Z$ runs over the connected components of $X^{\GG_m}$.
\end{enumerate}  
\end{thm}

Recall that $V$ is a finite-dimensional vector space over a field $k$, then Milne denotes the corresponding algebraic group by $V_\fa$ (see Milne \cite[Section 2.6, p. 40, Remark 4.35, p. 97, or Section 10.9, p. 189]{Milne_algebraic_groups_cup}).

While one may consider the action of a torus $\TT$ on an algebraic variety, one can apply the following observation to reduce to the case of a one-parameter subgroup.

\begin{lem}
\label{lem:Abdellaoui_3-7}  
(See Abdellaoui \cite[Lemma 3.7, p. 8]{Abdellaoui_2013arxiv_v2}.)  
Let $X$ be a normal algebraic variety with a linear action of a torus $\TT$. If the subset $X^\TT\subset X$ of fixed points of $\TT$ is non-empty, then then there exists a one-parameter subgroup $\Lambda \subset \TT$ such that $X^\Lambda = X^\TT$.
\end{lem}

\section{Bia{\l}ynicki--Birula decompositions for $\CC^*$ actions on complex projective space}
\label{sec:BB_decomposition_complex_projective_space}
The following elementary lemma is surely well-known, but we include its statement and proof due to its relevance to our discussion and the absence of a reference known to us. A general result for \emph{integral cohomology projective spaces} is given by Bredon \cite[Chapter VII, Section 5, p. 393, Theorem 5.1]{Bredon}.

\begin{lem}[Fixed-point sets of circle actions on complex projective space]
\label{lem:Fixed-point_sets_circle_actions_complex_projective_space}  
Let $n\geq 2$ be an integer and
\begin{equation}
  \label{eq:Complex_linear_circle_representation_complex_projective_space}
  \rho_{\PP^n}:S^1\times\PP(\CC^n)\to\PP(\CC^n)
\end{equation}
be the circle action on complex projective space $\PP(\CC^n)$ induced by a unitary representation $\rho:S^1\to \U(n)$ as in \eqref{eq:Circle_matrix_representation}. The set\footnote{In the notation of Bredon \cite[Chapter I, Section 5, p. 44]{Bredon}.} $\PP(\CC^n)^{S^1}$ of fixed points of the circle action \eqref{eq:Complex_linear_circle_representation_complex_projective_space} is a finite union of complex projective subvarieties $\PP(L) \subset \PP(\CC^n)$, where
\begin{equation}
  \label{eq:Coordinate_subspace_defining_fixed-point_set}
  L = \{(z_1,\ldots, z_n) \in \CC^n: z_k = 0, \text{ for all } k \notin J\}
\end{equation}
is a coordinate subspace defined by a nonempty subset $J \subset \{1,\ldots,n\}$ such that $l_k=l$ for some $l\in\ZZ$ and all $k\in J$.
\end{lem}

\begin{rmk}[Fixed-point sets of $\CC^*$ actions on complex projective space]
\label{rmk:Fixed-point_sets_C*_actions_complex_projective_space}  
The statement and proof of Lemma \ref{lem:Fixed-point_sets_circle_actions_complex_projective_space} extend \mutatis to the case of a $\CC^*$ action
\begin{equation}
  \label{eq:Complex_linear_C*_representation_complex_projective_space}
  \rho_{\PP^n}^\CC:\CC^*\times\PP(\CC^n)\to\PP(\CC^n)
\end{equation}
on complex projective space $\PP(\CC^n)$ induced by the complexification $\rho_\CC:\CC^*\to \GL(n,\CC)$ in \eqref{eq:C*_matrix_representation} of the unitary representation $\rho:S^1\to \U(n)$ in \eqref{eq:Circle_matrix_representation}. Moreover, by the forthcoming Lemma \ref{lem:S1_and_C*_fixed_points_coincide_on_complex_analytic_space}, the fixed-point sets $\PP(\CC^n)^{S^1}$ and $\PP(\CC^n)^{\CC^*}$ necessarily coincide.
\end{rmk}

\begin{proof}[Proof of Lemma \ref{lem:Fixed-point_sets_circle_actions_complex_projective_space}]
If $J = \{1,\ldots,n\}$, then $L = \CC^n$ and $\PP(L) = \PP(\CC^n)$ and (by definition of $J$) the action of $S^1$ on $\PP(\CC^n)$ is trivial (since $l_1=\cdots=l_n=l$), so $\PP(\CC^n)^{S^1} = \PP(\CC^n)$.

If $J \subsetneq \{1,\ldots,n\}$ is a nonempty subset such that $l_k=l$ for some $l\in\ZZ$ and all $k\in J$, then $\rho(e^{i\theta})$ acts on $L$ by scalar multiplication by $e^{il\theta}$ and so the coordinate subspace $L$ in \eqref{eq:Coordinate_subspace_defining_fixed-point_set} defines a complex projective subspace $\PP(L) \subset \PP(\CC^n)^{S^1}$.

Conversely, suppose $q\in \PP(\CC^n)^{S^1}$. Thus, $q = (q_1:\cdots:q_n)$ is a fixed point of the circle action \eqref{eq:Complex_linear_circle_representation_complex_projective_space} and because that action is induced by a unitary representation as in \eqref{eq:Circle_matrix_representation}, we must have
\[
  \rho(e^{i\theta})(q_1,\ldots,q_n) = (e^{l_1\theta}q_1,\ldots,e^{l_n\theta}q_n)
  = e^{l\theta}(q_1,\ldots,q_n),
  \quad\text{for all } \theta \in \RR,
\]
for some $l\in\ZZ$. Let $J \subseteq \{1,\ldots,n\}$ be the subset of indices $k$ such that $q_k\neq 0$ and observe that the preceding identity implies that $l_k = l$, for all $k \in J$. Hence, $q \in \PP(L)$, where $L$ in \eqref{eq:Coordinate_subspace_defining_fixed-point_set} is the coordinate subspace defined by $J$. Since $\{1,\ldots,n\}$ has $2^n$ subsets, including $\varnothing$ and $\{1,\ldots,n\}$, then either $\PP(\CC^n)^{S^1} = \PP(\CC^n)$ (when the $S^1$ action is trivial) or $\PP(\CC^n)^{S^1}$ is a union over up to $2^n-2$ non-empty, distinct complex projective subspaces.
\end{proof}

\begin{lem}[Bia{\l}ynicki--Birula decompositions for $\CC^*$ actions on complex projective space]
\label{lem:BB_decomposition_complex_projective_space}  
Let $X$ be a complex vector space of dimension greater than or equal to two and $\PP(X)$ be the corresponding complex projective space. Let $\rho_\CC:\CC^* \to \GL(X)$ be a representation and let $\CC^*\times\PP(X)\to \PP(X)$ be the induced $\CC^*$ action. Then $\PP(X)$ admits a Bia{\l}ynicki--Birula decomposition in the sense of Definition \ref{maindefn:BB_decomposition_complex_manifold}.
\end{lem}

\begin{proof}
We may assume without loss of generality that $X = \CC^n$ for $n\geq 2$. From Remark \ref{rmk:Classification_complex_representations_C*}, we know that $\rho_\CC:\CC^* \to \GL(n,\CC)$ has the form \eqref{eq:C*_matrix_representation}. The connected components $\PP(\CC^n)_\alpha^0$ of the subset $\PP(\CC^n)^0$ of fixed points of the $\CC^*$ action on $\PP(\CC^n)$ are given explicitly by Lemma \ref{lem:Fixed-point_sets_circle_actions_complex_projective_space} and Remark \ref{rmk:Fixed-point_sets_C*_actions_complex_projective_space} and, in particular, are identified with projective subvarieties $\PP(L) \subset \PP(\CC^n)$, where $L \subset \CC^n$ is a coordinate subspace as in \eqref{eq:Coordinate_subspace_defining_fixed-point_set}.

If $z = (z_1,\ldots,z_n) \in X\less\{0\}$, then $\rho_\CC(\lambda)z = (\lambda^{l_1}z_1,\lambda^{l_2}z_2,\ldots,\lambda^{l_n}z_n)$ and we claim that $\lim_{z\to 0}\rho_\CC(\lambda)z \in \PP(\CC^n)_\alpha^0$, for some $\alpha$, that is, $[z] \in \PP(\CC^n)_\alpha^+$, for some $\alpha$. If $l_1=l_2=\cdots=l_n$, then the induced action on $\PP(\CC^n)$ is trivial and so we may assume without loss of generality that at there at least two distinct integer weights. By relabeling the coordinates on $\CC^n$, we may further assume without loss of generality that $l_1 \leq l_2\leq \cdots \leq l_n$. We first consider the case $l_1<l_2<\cdots<l_n$ and observe that
\begin{multline*}
  \left(\lambda^{l_1}z_1,\lambda^{l_2}z_2,\ldots,\lambda^{l_n}z_n\right)
  =
  \lambda^{l_1}\left(z_1,\lambda^{l_2-l_1}z_2,\ldots,\lambda^{l_n-l_1}z_n\right),
  \\
  \text{for all } \lambda \in \CC^*
  \text{ and } (z_1,\ldots,z_n) \in \CC^n\less\{0\},
\end{multline*}
so that
\begin{multline*}
  \left[\lambda^{l_1}z_1,\lambda^{l_2}z_2,\ldots,\lambda^{l_n}z_n\right]
  =
  \left[z_1,\lambda^{l_2-l_1}z_2,\ldots,\lambda^{l_n-l_1}z_n\right],
  \\
  \text{for all } \lambda \in \CC^*
  \text{ and } [z_1,\ldots,z_n] \in \PP(\CC^n).
\end{multline*}
From Lemma \ref{lem:Fixed-point_sets_circle_actions_complex_projective_space} and Remark \ref{rmk:Fixed-point_sets_C*_actions_complex_projective_space}, the $n$ distinct fixed points of this $\CC^*$ action on $\PP(\CC^n)$ are given by
\[
  \PP(\CC^n)_j^0 = [0,\ldots,0,1_j,0,\ldots,0], \quad\text{for } j = 1,\ldots,n,
\]
where $1_j$ denotes $1$ in the $j$-th entry. Denote the affine coordinate domains on $\PP(\CC^n)$ by
\begin{multline*}
  U_j := \{[z_1,\ldots,z_n] \in \PP(\CC^n): z_j \neq 0\}
  \\
  = \{[z_1,\ldots,z_{j-1},1,z_{j+1},\cdots,z_n] \in \PP(\CC^n)\}, \quad\text{for } j=1,\ldots,n.
\end{multline*}
(As an aside, we see that
\[
  \PP(\CC^n)\less U_j = \{[z_1,\ldots,z_{j-1},0,z_{j+1},\cdots,z_n] \in \PP(\CC^n)\} \cong \PP(\CC^{n-1}),
  \quad\text{for } j=1,\ldots,n,
\]
and therefore we can alternatively define the affine coordinate domains on $\PP(\CC^n)$ by
\[
  U_j = \PP(\CC^n)\less \{[z_1,\ldots,z_{j-1},0,z_{j+1},\cdots,z_n] \in \PP(\CC^n)\},
  \quad\text{for } j=1,\ldots,n.
\]
We use a similar definition in our discussion below of the general case where the $n$ weights need not be distinct.) We then observe that
\[
  \lim_{\lambda\to 0}\left[\lambda^{l_1}z_1,\lambda^{l_2}z_2,\ldots,\lambda^{l_n}z_n\right]
  =
  [1,0,\ldots,0],
  \quad\text{for all } [z_1,\ldots,z_n] \in U_1,
\]
since $l_j-l_1 > 0$ for $j=2,\ldots,n$ and $[z_1,0,\ldots,0] = [1,0,\ldots,0] \in U_1$. Thus,
\[
  \PP(\CC^n)_1^+ = U_1.
\]
Suppose $[z_1,\ldots,z_n] \in \PP(\CC^n)\less U_1$, so $z_1=0$ but $z_j\neq 0$ for some $j \in \{2,\ldots,n\}$. If $j=2$, then
\[
  \left[0,\lambda^{l_2}z_2,\lambda^{l_3}z_3,\ldots,\lambda^{l_n}z_n\right]
  =
  \left[0,z_2,\lambda^{l_3-l_2}z_3,\ldots,\lambda^{l_n-l_2}z_n\right]
\]
and therefore
\[
  \lim_{\lambda\to 0}\left[0,\lambda^{l_2}z_2,\lambda^{l_3}z_3,\ldots,\lambda^{l_n}z_n\right]
  =
  [0,1,0,\ldots,0],
  \quad\text{for all } [0,z_2,\ldots,z_n] \in U_2\cap \left(\PP(\CC^n)\less U_1\right).
\]
Consequently,
\[
  \PP(\CC^n)_2^+ = U_2\cap \left(\PP(\CC^n)\less U_1\right).
\]
Continuing in this way, we see that
\[
  \PP(\CC^n)_j^+ = U_j\cap \left(\PP(\CC^n)\less (U_1\cup\cdots\cup U_{j-1})\right),
  \quad\text{for } j = 2,\ldots,n,
\]
and, in particular, we obtain the plus decomposition,
\[
  \PP(\CC^n) = \bigsqcup_{j=1}^n \PP(\CC^n)_j^+.
\]
It is straightforward to check that the natural projections $\pi_j^+:\PP(\CC^n)_j^+ \to \PP(\CC^n)_j^0$ and sets $\PP(\CC^n)_j^+$ have the properties for a Bia{\l}ynicki--Birula plus decomposition enumerated in Definition \ref{maindefn:BB_decomposition_complex_manifold}.

To obtain the minus decomposition, we observe that
\begin{multline*}
  \left[\lambda^{l_1}z_1,\ldots,\lambda^{l_{n-1}}z_{n-1},\lambda^{l_n}z_n\right]
  =
  \left[z_1^{l_1-l_n},\ldots,\lambda^{l_{n-1}-l_n}z_{n-1},\ldots,z_n\right],
  \\
  \text{for all } \lambda \in \CC^*
  \text{ and } [z_1,\ldots,z_n] \in \PP(\CC^n).
\end{multline*}
Hence,
\[
  \lim_{\lambda\to \infty}\left[\lambda^{l_1}z_1,\lambda^{l_2}z_2,\ldots,\lambda^{l_n}z_n\right]
  =
  [0,\ldots,0,1],
  \quad\text{for all } [z_1,\ldots,z_n] \in U_n,
\]
since $l_j-l_n < 0$ for $j=1,\ldots,n-1$ and $[0,\ldots,0,z_n] = [0,\ldots,0,1] \in U_n$. Thus,
\[
  \PP(\CC^n)_n^- = U_n.
\]
Continuing in this way, we see that
\[
  \PP(\CC^n)_j^- = U_j\cap \left(\PP(\CC^n)\less (U_{j+1}\cup\cdots\cup U_n)\right),
  \quad\text{for } j = 1,\ldots,n-1,
\]
and, in particular, we obtain the minus decomposition,
\[
  \PP(\CC^n) = \bigsqcup_{j=1}^n \PP(\CC^n)_j^-.
\]
As in the plus case, one finds that the natural projections $\pi_j^-:\PP(\CC^n)_j^- \to \PP(\CC^n)_j^0$ and sets $\PP(\CC^n)_j^-$ have the properties for a Bia{\l}ynicki--Birula minus decomposition enumerated in Definition \ref{maindefn:BB_decomposition_complex_manifold}.

For the general case, where $l_1\leq \cdots \leq l_n$ but at least two weights are distinct, we write $\CC^n = V_1\times \cdots \times V_p$, where $p \geq 2$ and we relabel the distinct weights as $k_1 < \cdots < k_p$, where each $k_j$ is equal to one or more of the integers $l_i$. We now proceed as in the case of $n$ distinct weights, $l_1< \cdots < l_n$, but write $z_j \in V_j$ for $j = 1, \ldots,p$ and observe that the connected components of the set of fixed points of the $\CC^*$ action on $\PP(\CC^n)$ are given by
\[
  \PP(\CC^n)_j^0 = \PP(0\times\cdots 0\times V_j\times 0\times\cdots\times 0),
  \quad\text{for } j = 1,\ldots,p.
\]
We now define
\[
  U_j := \PP(\CC^n)\less\PP(V_1\times\cdots\times V_{j-1}\times 0\times V_{j+1}\times\cdots\times V_p),
  \quad\text{for } j = 1,\ldots,p,
\]
and thus obtain
\begin{align*}
  \PP(\CC^n)_1^+ &= U_1,
  \\
  \PP(\CC^n)_j^+ &= U_j \cap \PP(\CC^n)\less\left(U_1\cup \cdots\cup U_{j-1}\right),
  \quad\text{for } j = 2 = 1,\ldots,p.
\end{align*}
Hence, we arrive at the plus decomposition for $\PP(\CC^n)$ and a similar argument yields the minus decomposition. Just as before, these subsets and their natural projections have the properties for Bia{\l}ynicki--Birula plus and minus decompositions enumerated in Definition \ref{maindefn:BB_decomposition_complex_manifold}.
\end{proof}

\section[Bia{\l}ynicki--Birula decompositions for $\CC^*$ actions on compact K\"ahler manifolds]{Bia{\l}ynicki--Birula decompositions for $\CC^*$ actions on compact, complex K\"ahler manifolds}
\label{sec:BB_decomposition_complex_manifold}
We summarize the discussions of $\CC^*$ actions on complex K\"ahler manifolds due to Bia{\l}ynicki--Birula and Sommese \cite[Section 0, p. 777]{Bialynicki-Birula_Sommese_1983}, Carrell and Sommese \cite[Section III, p. 55]{Carrell_Sommese_1978ms}, \cite[Section 1c, p. 569]{Carrell_Sommese_1979cmh}, Fujiki \cite[Section 1, p. 798, and Section 2, p. 803]{Fujiki_1979}, and Yang \cite{Yang_2008}. In Section \ref{sec:BB_decomposition_noncompact_complex_manifold} we shall relax the hypothesis in Theorem \ref{thm:Bialynicki-Birula_decomposition_compact_complex_Kaehler_manifold} that the complex manifold be compact and remove the hypothesis that it be K\"ahler.

\begin{lem}[Fixed-point set of a complex Lie group acting on a complex manifold]
\label{lem:Fixed-point_set_complex_Lie_group_on_complex_manifold}
(See Fujiki \cite[Lemma 1.1 and Corollary 1.2, p. 799]{Fujiki_1979}.)
Let $(X,\sO_X)$ be a complex analytic space and $G$ a connected, complex Lie group acting biholomorphically on $X$. If $X^G := \{x\in X: g\cdot x = x, \text{ for all } g \in G\}$ denotes the fixed-point set of $G$ on $X$, then $X^G$ is an analytic subset of $X$. If $X$ is a complex manifold and $G$ is reductive, then $X^G$ is a complex submanifold of $X$.
\end{lem}

\begin{rmk}[Dimensions of fixed-point submanifolds]
\label{rmk:Dimensions_fixed-point_submanifolds}
Suppose in the setting of Lemma \ref{lem:Fixed-point_set_complex_Lie_group_on_complex_manifold} that $X$ is a complex manifold and $G$ is reductive, so $X^G$ is a complex submanifold of $X$. However, the connected components of $X^G$ may have different dimensions. In the forthcoming Example \ref{exmp:Bialynicki-Birula_decomposition_CP^1}, we consider $S^1$ actions on $\PP^1$ induced by $\CC^*$ actions. Extending this example to $\PP^2$, one sees that\footnote{See the discussion of dimensions of fixed-point submanifolds by Thomas Rot at \url{https://mathoverflow.net/questions/268069/fixed-point-set-of-smooth-circle-action}.} the $\CC^*$ action on $\PP^2$ induced by $[z_0,z_1,z_2] \mapsto [\lambda z_0,z_1,z_2]$, for all $\lambda \in \CC^*$ and $[z_0,z_1,z_2] \in \PP^2$, has fixed-point sets with different dimensions, namely the isolated point $[1,0,0]$ and the line $\PP^1 \subset \PP^2$ defined by the subset of all points $[0,z_1,z_2] \in \PP^2$. By extending this example to $\PP^n$ with $n\geq 2$, one can produce\footnote{See the discussion of dimensions of fixed-point submanifolds by Jason DeVito at \url{https://math.stackexchange.com/questions/2632738/when-is-the-dimension-of-the-fixed-point-set-of-an-isometry-defined}.} pairs of fixed-point submanifolds $\PP^a, \PP^b \subset \PP^n$ with $a+b = n-1$. More generally again, one can show\footnote{See the discussion of dimensions of fixed-point submanifolds by Moishe Kohan at \url{https://math.stackexchange.com/questions/2632738/when-is-the-dimension-of-the-fixed-point-set-of-an-isometry-defined}.}  that for every integer $n\geq 1$ there exists an $n$-dimensional connected compact Riemannian manifold $M$ and an isometry $\sigma:M\to M$ whose fixed-point set contains components of all possible dimensions, $0,1,2,\ldots,n-1$. For more examples and discussion, see Lemma \ref{lem:Fixed-point_sets_circle_actions_complex_projective_space} and the articles by Delzant \cite{Delzant_1988} and Yang \cite{Yang_1995}.
\end{rmk}

In particular, $\CC^*$ is a connected, complex, reductive Lie group and thus Lemma \ref{lem:Fixed-point_set_complex_Lie_group_on_complex_manifold} implies that if $\CC^*\times X \to X$ is a holomorphic $\CC^*$-action, then its fixed-point subset $X^0:=X^G$ is a complex submanifold of $X$. Assume that $X^0$ is nonempty with connected components $X_\alpha^0$, for $\alpha$ varying over an index set $\sA$. Each $X_\alpha^0$ is a complex submanifold of $X$ by Lemma \ref{lem:Fixed-point_set_complex_Lie_group_on_complex_manifold} for all $\alpha \in \sA$.

Assume further that $X$ is compact and K\"ahler. It is a basic fact (see Sommese \cite[Section II, Lemma II-A, p. 109]{Sommese_1975}) that for any point $x \in X$, the holomorphic map
\[
  \CC^* \ni \lambda \mapsto \lambda\cdot x \in X
\]
extends to a holomorphic map
\[
  \PP^1 \ni \lambda \mapsto \lambda\cdot x \in X.
\]
Thus, one observes (see Carrell and Sommese \cite[Section 1c, p. 569]{Carrell_Sommese_1979cmh}) that the following limits exist,
\[
  \lim_{\lambda\to 0} \lambda\cdot x \quad\text{and}\quad \lim_{\lambda\to \infty}\lambda\cdot x,
\]
and by the group action properties, both must lie in the fixed-point subset $X^0$. This yields two $\CC^*$-invariant decompositions of $X$, called the plus and minus decompositions as in Definition \ref{maindefn:BB_decomposition_complex_manifold}, with 
\[
  X = \bigsqcup_\alpha X_\alpha^+ \quad\text{and}\quad X = \bigsqcup_\alpha X_\alpha^-,
\]
where, as in \eqref{eq:Complex_manifold_Xalpha_plus_minus_submanifolds},
\[
  X_\alpha^+ := \left\{x:\lim_{\lambda\to 0}\lambda\cdot x \in X_\alpha^0 \right\}
  \quad\text{and}\quad
  X_\alpha^- := \left\{x:\lim_{\lambda\to \infty}\lambda\cdot x \in X_\alpha^0 \right\}.
\]
These decompositions were first obtained for smooth algebraic varieties by Bia{\l}ynicki--Birula (see Section \ref{sec:BB_decomposition_algebraic_variety_scheme}) and subsequently shown to exist for compact, complex K\"ahler manifolds by Carrell and Sommese \cite{Carrell_Sommese_1978ms} and Fujiki \cite{Fujiki_1979}. The properties of the decomposition are given by the following analogue of Theorem \ref{thm:Milne_13-47} when $X$ is a smooth algebraic variety.

\begin{thm}[Bia{\l}ynicki--Birula decomposition for a holomorphic $\CC^*$ action on a compact, complex K\"ahler manifold]
\label{thm:Bialynicki-Birula_decomposition_compact_complex_Kaehler_manifold}
(See Carrell and Sommese \cite[Proposition II, p. 55]{Carrell_Sommese_1978ms} and Fujiki \cite[Lemma 2.1, p. 803, and Theorem 2.2, p. 805]{Fujiki_1979}.)
Let $X$ be a compact, complex K\"ahler manifold. If $\CC^*\times X \to X$ is a holomorphic $\CC^*$ action with at least one fixed point, then $X$ has plus and minus Bia{\l}ynicki--Birula decompositions in the sense of Definition \ref{maindefn:BB_decomposition_complex_manifold} with finitely many components $X_\alpha^\pm$ as in \eqref{eq:Complex_manifold_Xalpha_plus_minus_submanifolds} for $\alpha\in\{1,\ldots,r\}$. Moreover, there is exactly one component $X_1^0$ (respectively, $X_r^0$) called the \emph{source} (respectively, the \emph{sink}) such that $TX\restriction X_1^0 = TX_1^0\oplus N_{X_1^0/X}^+$ (respectively, $TX\restriction X_r^0 = TX_r^0\oplus N_{X_r^0/X}^-$).
\end{thm}

Yang \cite[Theorem 4.12, p. 92]{Yang_2008} proves a version of Theorem \ref{thm:Bialynicki-Birula_decomposition_compact_complex_Kaehler_manifold} for holomorphic $\CC^*$ actions on complex K\"ahler manifolds that are not necessarily compact but where the circle action $S^1\times X\to X$ induced by the restriction of the $\CC^*$ action to $S^1\subset \CC^*$ has a Hamiltonian function that is proper and bounded below and the number of connected components of $X^{\CC^*}$ is finite.

As asserted by Theorem \ref{thm:Bialynicki-Birula_decomposition_compact_complex_Kaehler_manifold}, the normal bundle of $X_\alpha^0$ in $X_\alpha$ is a specific subbundle of the normal bundle of $X_\alpha^0$ in $X$ whose description we now provide in more detail. For each point $x \in X$, elements $\lambda \in \CC^*$ act linearly on the holomorphic tangent bundle $TX$ of $X$ via the differentials $d\rho(\lambda,x)$, where we denote
\[
  \rho:\CC^*\times X \ni (\lambda,x) \mapsto \rho(\lambda,x) = \lambda\cdot x \in X.
\]
When $x\in X^0$, it is well known that the resulting complex representation of $\CC^*$ on $T_xX$ is determined by the existence of a basis $v_1, \ldots, v_n$ of $T_xX$ and integers $m_1,\ldots,m_n$, such that $\lambda\cdot v_i = d\rho(\lambda,x)v_i = \lambda^{m_i}v_i$, for all $\lambda\in\CC^*$. The integers $m_i$ are the weights of the action of $\CC^*$ on $T_xX$. Therefore, for each $\alpha$, one has a canonical holomorphic direct-sum decomposition
\[
  TX \restriction X_\alpha^0 = TX_\alpha^0 \oplus N^+X_\alpha^0 \oplus N^-X_\alpha^0,
\]
where $N^+X_\alpha^0$ (respectively, $N^-X_\alpha^0$) is the holomorphic vector bundle over $X_\alpha^0$ whose fiber at $x$ is generated by the $v_i$ corresponding to positive (respectively, negative) weights $m_i$, and $T_xX_\alpha^0$ is generated by the $v_i$ corresponding to $m_i=0$. Then $N^+X_\alpha^0$ (respectively, $N^-X_\alpha^0$) is the normal bundle of $X_\alpha^0$ in $X_\alpha^+$ (respectively, $X_\alpha^-$) whose existence is asserted by Theorem \ref{thm:Bialynicki-Birula_decomposition_compact_complex_Kaehler_manifold}.

According to Carrell and Sommese \cite[Theorem 1]{Carrell_Sommese_1979cmh}, the \emph{Morse--Bott co-index} of $X_\alpha^0 \subset X$ for the Morse--Bott function $f$ in Frankel's Theorem  \ref{mainthm:Frankel_almost_Hermitian} (see also Carrell and Sommese \cite[Theorem, p. 567]{Carrell_Sommese_1979cmh} or Frankel \cite[Section 3, pp. 2--5]{Frankel_1959}) is given by
\[
  \lambda_\alpha^+ = \lambda_x^+(f) = \dim_\RR N_x^+X_\alpha^0
\]
for any $x \in X_\alpha^0$, where $f$ is the Hamiltonian function as in \eqref{eq:MomentMap}, so $df = \iota_\Theta\omega$ and $\omega$ is the circle-invariant K\"ahler form on $X$ given by \eqref{eq:Fundamental_two-form} and $\Theta$ is the smooth vector field on $X$ generated by the circle action (see the forthcoming \eqref{eq:Vector_field_generator_circle_action}). Consequently, the \emph{Morse--Bott index} and \emph{nullity} for $f$ at any point $x \in X_\alpha^0$ are given by
\begin{align*}
  \lambda_\alpha^- &= \lambda_x^-(f) = \dim_\RR N_x^-X_\alpha^0,
  \\
  \lambda_\alpha^0 &= \lambda_x^0(f) = \dim_\RR T_xX_\alpha^0.
\end{align*}
Associated to a $\CC^*$ action on $X$ are the circle and radial actions arising from the circle subgroup $S^1 \subset \CC^*$ and the radial subgroup $\RR^+ \subset \CC^*$ consisting of all positive real numbers. In \cite[Section IIa]{Carrell_Sommese_1979cmh}, Carrell and Sommese show that $f$ is a \emph{Lyupanov function}, that is, $f$ is strictly increasing along the radial orbits in $X$. They also prove some geometric consequences of the existence of this function. The two group actions give rise to a pair of vector fields $\Theta$ (as above) and $R$ on $X$ such that
\begin{equation}
  \label{eq:Almost_complex_structure_on_circular_vector_field_equals_radial_vector_field}
  J\Theta = R \quad\text{on X},
\end{equation}
where $J$ is the almost complex structure on $TX$. Indeed, one can define
\begin{subequations}
\label{eq:Vector_field_generator_circle_and_radial_actions}  
\begin{align}
  \label{eq:Vector_field_generator_circle_action} 
  \Theta_x &:= \left.\frac{d}{d\theta}(e^{i\theta}\cdot x)\right|_{\theta=0} \quad (\theta\in\RR),
  \\
  \label{eq:Vector_field_generator_radial_action} 
  R_x &:= \left.\frac{d}{dr}(r\cdot x)\right|_{r=1} \quad (r\in\RR^+),
          \quad\text{for all } x \in X.
\end{align}
\end{subequations}
Note that $\Zero (\Theta) = \Zero (R) = X^0$. The trajectories of $R$, that is, the rays $r\cdot x$, for $r\in(0,\infty)$, have limits in $X^0$ as $r\to 0$ or $r\to\infty$. The next lemma shows that these limits lie in different components of the fixed-point set $X^0\subset X$.

\begin{lem}
\label{lem:Carrell_Sommese_1979_lemma_1_and_corollaries_1_and_2}  
(See Carrell and Sommese \cite[Section IIa, Lemma 1 and Corollaries 1 and 2]{Carrell_Sommese_1979cmh})
Let $X$ be a compact, complex K\"ahler manifold with a holomorphic $\CC^*$ action $\CC^*\times X\to X$, K\"ahler form $\omega$, and smooth Hamiltonian function $f:X\to\RR$ such that $df = \iota_\Theta\omega$ as in \eqref{eq:MomentMap}, where $\Theta$ is as in \eqref{eq:Vector_field_generator_circle_action}. Then the following hold:
\begin{enumerate}
\item The function $f$ is decreasing along the trajectories of $R$, where $R$ is as in \eqref{eq:Vector_field_generator_radial_action}.
\item The Morse--Bott index of $f$ at $x\in X_\alpha^0$ is equal to $\dim_\RR N_{X_\alpha^0/X}^-|_x$ and, consequently, is even.
\item The source $X_1^0$ of $X$ is $\{x \in X: \text{$f$ assumes its absolute minimum at $x$}\}$.
\item The sink $X_r^0$ of $X$ is $\{x \in X: \text{$f$ assumes its absolute minimum at $x$}\}$.
\end{enumerate}
\end{lem}

By Lemma \ref{lem:Carrell_Sommese_1979_lemma_1_and_corollaries_1_and_2}, the fixed-point components $X_1^0, X_2^0,\ldots,X_r^0$ can be indexed in such a way that
\[
  f(X_1^0) < f(X_2^0) \leq \cdots \leq f(X_{r-1}^0) < f(X_r^0),
\]
as expected from Frankel \cite[Section 3, pp. 2--5]{Frankel_1959}. Observe that from the proof of \cite[Lemma 1]{Carrell_Sommese_1979cmh} by Carrell and Sommese, we have
\[
  (df)\xi = \iota_\Theta\omega(\xi) = \omega(\Theta,\xi) = -\omega(JR,\xi) = -g(R,\xi),
\]
for any $\xi \in C^\infty(TX)$ and using the convention \eqref{eq:Fundamental_two-form}, 
\[
  \omega(\eta,\xi) = g(\eta,J\xi), \quad\text{for all } \xi, \eta \in C^\infty(TX),
\]
and assumption that $(\omega,g,J)$ is a compatible triple. The definition \eqref{eq:DefineGradient}
\[
  (df)\xi = g(\grad_g f,\xi), \quad\text{for all } \xi \in C^\infty(TX)
\]
of the gradient of $f$ with respect to the Riemannian metric $g$ thus yields
\[
  g(\grad_g f,\xi) = -g(R,\xi), \quad\text{for all } \xi \in C^\infty(TX),
\]
and therefore
\begin{equation}
  \label{eq:Gradient_Hamiltonian_and_radial_rotation_vector_fields}
  \grad_g f = -R = -J\Theta.
\end{equation}
In particular, for $x\notin X^0$, we obtain $df(x)R = -g(R,R)<0$. Note that this leads to the opposite of the conclusion by Carrell and Sommese in \cite[Lemma 1]{Carrell_Sommese_1979cmh} regarding the behavior of $f$ along trajectories of $R$ because we used the convention $\omega = g(\cdot, J\cdot)$ and not $g(J\cdot, \cdot)$, so $f$ decreases rather than increases along trajectories of $R$.

\section[Bia{\l}ynicki--Birula decomposition for $\CC^*$ actions on noncompact K\"ahler manifolds]{Bia{\l}ynicki--Birula decompositions for $\CC^*$ actions on noncompact, complex manifolds}
\label{sec:BB_decomposition_noncompact_complex_manifold}
In this section, we shall relax the hypothesis in Theorem \ref{thm:Bialynicki-Birula_decomposition_compact_complex_Kaehler_manifold} that the complex manifold be compact, remove the hypothesis that it be K\"ahler, and conclude this section with the proof of Theorem \ref{mainthm:BB_decomposition_complex_manifold_C*_action}, our generalization of Theorem \ref{thm:Bialynicki-Birula_decomposition_compact_complex_Kaehler_manifold}. The following simple example highlights the differences between the compact and noncompact cases when constructing Bia{\l}ynicki--Birula decompositions.

\begin{exmp}[Bia{\l}ynicki--Birula decomposition for $\CC^2$]
\label{exmp:Bialynicki-Birula_decomposition_C^2}  
Let $\rho:\CC^*\to\GL(2,\CC)$ be a linear representation defined by $\CC^2 \ni (z,w) \mapsto (t^mz,t^nw) \in \CC^2$ for given integers $m, n \in \ZZ$ and all $t\in\CC^*$. We first consider $X = \CC^2$ with the $\CC^*$ action induced by the representation $\rho$ and let $p = (0,0) \in X$ denote the origin. If $m,n$ are both nonzero and $e_1 := (1,0)$ and $e_2 := (0,1) \in \CC^2$, then the fixed-point set is $F = \{p\}$ and, in the notation of Theorem \ref{thm:Bialynicki-Birula_decomposition_compact_complex_Kaehler_manifold}, one has the following splittings of the tangent space to $X$ at the origin,
\[
  T_pX
  =
  \begin{cases}
    N_p^+X = \CC^2, &\text{if } m, n > 0,
    \\
    N_p^-X = \CC^2, &\text{if } m, n < 0,
    \\
    N_p^+X\oplus N_p^-X, \text{ for } N_p^+X = \CC e_1,
    \text{ and } N_p^-X = \CC e_2, &\text{if } m > 0, n < 0.
  \end{cases}
\]
However, the global Bia{\l}ynicki--Birula decompositions of $X = \CC^2$ are given by 
\[
  X
  =
  \begin{cases}
    X^+ = \CC^2, &\text{if } m, n > 0,
    \\
    X^-\cup\{0\}, \text{ for } X^- = \CC^2\less\{0\}, &\text{if } m, n < 0,
    \\
    X^+\times X^-, \text{ for } X^+ = \CC e_1 \text{ and } X^- = \CC e_2, &\text{if } m > 0, n < 0,
  \end{cases}
\]
and thus do not follow the pattern predicted by Theorem \ref{thm:Bialynicki-Birula_decomposition_compact_complex_Kaehler_manifold}, where $X$ is assumed to be compact.
\qed
\end{exmp}

\begin{exmp}[Bia{\l}ynicki--Birula decomposition for $\PP(\CC^2)$]
\label{exmp:Bialynicki-Birula_decomposition_CP^1}
We let $\rho:\CC^*\to\GL(2,\CC)$ be as in Example \ref{exmp:Bialynicki-Birula_decomposition_C^2}, but we now consider $X = \PP^1 = \PP(\CC^2)$ with the $\CC^*$ action induced by the representation $\rho$ and coordinate charts
\begin{align*}
  \varphi_1:U_1 &= \{[z,w]:z\neq 0, w\in\CC\} \ni [z,w] \mapsto w/z \in \CC, 
                  \\
  \varphi_2:U_2 &= \{[z,w]:w\neq 0, z\in\CC\} \ni [z,w] \mapsto z/w \in \CC.
\end{align*}
If $m=n$, then the induced $\CC^*$ action on $X$ is trivial, so we shall assume that $m\neq n$, in which case we obtain the following fixed-point set:
\[
  X^0  =  \{p_1,p_2\} = \{[1,0], [0,1]\}.
\]
In the notation of Theorem \ref{thm:Bialynicki-Birula_decomposition_compact_complex_Kaehler_manifold}, one has the following splittings of the tangent space to $X$ at the point $p_1 = [1,0]$ (the analysis for $p_2 = [0,1]$ is identical):
\[
  T_{p_1}X
  =
  \begin{cases}
    N_{p_1}^+X \cong \CC &\text{if } m, n > 0 \text{ or }  m, n < 0,
    \\
    N_{p_1}^-X \cong \CC &\text{if } m > 0, n < 0 \text{ or }  m <0, n > 0.
  \end{cases}
\]
The global Bia{\l}ynicki--Birula decompositions of $X = \PP^1$ do follow the pattern predicted by Theorem \ref{thm:Bialynicki-Birula_decomposition_compact_complex_Kaehler_manifold}: 
\[
  X = X^+ = X^- = \PP^1,
\]
for all combinations of $m\neq n$.
\qed
\end{exmp}

Examples \ref{exmp:Bialynicki-Birula_decomposition_C^2} and \ref{exmp:Bialynicki-Birula_decomposition_CP^1} and their analogues in higher dimensions suggest that Theorem \ref{thm:Bialynicki-Birula_decomposition_compact_complex_Kaehler_manifold} describes the local structure of a complex manifold $X$ near a fixed point of a holomorphic $\CC^*$ action, but not necessarily its global structure when the hypothesis that $X$ is compact is relaxed.

The proofs of the forthcoming Proposition \ref{prop:Extension_C*_actions_real_analytic_almost_Hermitian_manifolds} and Theorem \ref{thm:Extension_C*_actions_Hermitian_complex_manifolds} have the virtue of clarifying the relationship between convergence of the limits in the Bia{\l}ynicki-Birula decomposition for manifolds with a $\CC^*$ action and solutions to gradient flow for a Hamiltonian function defined by the induced circle action and a circle-invariant non-degenerate two-form. The forthcoming Lemma \ref{lem:Barbalat} is well-known in the theory of Lyapunov stability for dynamical systems (for example, see Khalil \cite[Lemma 8.2, p. 323]{Khalil_1992})
and we shall use the variant in Lemma \ref{lem:Barbalat_variant} in our proof of Proposition \ref{prop:Extension_C*_actions_real_analytic_almost_Hermitian_manifolds}.

\begin{lem}[Barb\u{a}lat's lemma]
\label{lem:Barbalat}
(See Barb\u{a}lat \cite[p. 269]{Barbalat_1959}, Farkas and Wegner \cite[Theorem 1, p. 825, and Theorem 4, p. 826]{Farkas_Wegner_2016}, or Niculescu and Popovici \cite[Corollary 2, p. 159]{Niculescu_Popovici_2012}.)
If $E$ is a Banach space, $h : [0,\infty) \to E$ is uniformly continuous, and $\lim_{t\to \infty}\int_0^t h(s)\,ds$ exists, then $h(t) \to 0$ as $t \to \infty$.
\end{lem}

\begin{lem}[Variant of Barb\u{a}lat's lemma]
\label{lem:Barbalat_variant}
(See Farkas and Wegner \cite[Theorem 5, p. 827]{Farkas_Wegner_2016}.)  
Let $p \in [1,\infty)$ and $q \in (1,\infty]$. If $h : [0,\infty) \to \RR$ obeys $h \in L^p(0,\infty)$ and $\dot h \in L^q(0,\infty)$, then $h(t) \to 0$ as $t \to \infty$.
\end{lem}

\begin{prop}[Subsequential limits of solutions to gradient flow are critical points]
\label{prop:Subsequential_limits_solutions_gradient_flow_critical_points}  
Let $(X,g)$ be a $C^2$ Riemannian manifold \footnote{In the sense of Abraham, Marsden, and Ratiu \cite[Definition 5.2.12, p. 352]{AMR}.}, $f:X\to\RR$ be a $C^2$ function, $x \in X$ be a point, and $u \in C^0([0,\infty)\times X, X)\cap C^{1;0}((0,\infty)\times X, X)$ be a solution\footnote{The notation means that the solution is continuously differentiable with respect to $t\in(0,\infty)$ and continuous with respect to $(t,x)\in [0,\infty)\times X$.} to the gradient flow equation
\begin{equation}
  \label{eq:Gradient_flow_equation_proposition}
  \frac{du(t,x)}{dt} = \grad_g f(u(t,x)), \quad u(0,x) = x,
\end{equation}
where the gradient vector field $\grad_g f$ is defined by the differential $f':TX\to\RR$ and the metric $g$ as in \eqref{eq:DefineGradient}. If $\{t_k\}_{k=1}^\infty \subset [0,\infty)$ is sequence such that $t_k\to\infty$ as $k\to\infty$, the limit $u_\infty(x) := \lim_{k\to\infty} u(t_k,x) \in X$ exists, and $f'':T^2M \to \RR$ is bounded on an open neighborhood of $u_\infty(x)$, then $u_\infty(x)$ is a critical point of $f$, that is, $df(u_\infty(x)) = 0$.
\end{prop}  

\begin{proof}
By choosing a local coordinate chart $(U,\varphi)$ centered at $u_\infty(x)$ in $X$, so $\varphi:U\to T_{u_\infty(x)}$ is a $C^2$ embedding onto an open neighborhood of the origin, we may reduce to the case where $X$ is replaced by the Hilbert space $H = T_{u_\infty(x)}$ with inner product $\langle\cdot,\cdot\rangle_H = g_{u_\infty(x)}(\cdot,\cdot)$. By the Fundamental Theorem of Calculus, the gradient flow equation \eqref{eq:Gradient_flow_equation_proposition} is equivalent to the integral equation
\begin{equation}
  \label{eq:Gradient_flow_integral_equation}
  u(t'',x) - u(t',x) = \int_{t'}^{t''}\grad_g f(u(t,x))\,dt,
  \quad\text{for all } 0 \leq t' \leq t'' <\infty,
\end{equation}
where the right-hand side is a Bochner integral (see Sell and You \cite[Appendix C.1, p. 613]{Sell_You_2002}). In particular, for all integers $k,l\geq 1$, the identity \eqref{eq:Gradient_flow_integral_equation} yields
\[
  u(t_k,x) - u(t_l,x) = \int_{t_l}^{t_k}\grad_g f(u(t,x))\,dt.
\]
Allowing $l=0$ and setting $t_0 = 0$, taking the limit as $k\to\infty$, and applying our hypothesis that $u(t_k,x) \to u_\infty(x)$ in $H$ as $k\to\infty$ gives
\[
  u_\infty(x) - u(0,x) = \int_0^\infty\grad_g f(u(t,x))\,dt.
\]
Hence, setting $h(t) := \grad_g f(u(t,x))$ for the given point $x$ and all $t\in[0,\infty)$, we see that $h \in L^1(0,\infty;H)$.

By applying the Chain Rule, substituting the gradient flow equation \eqref{eq:Gradient_flow_equation_proposition}, and applying the Riesz isomorphism to $f'(u(t,x)) \in H^* \cong H$, noting that $f'(u) = g(\cdot,\grad_g f(u))$, we obtain
\begin{multline*}
  \frac{df(u(t,x))}{dt} = f'(u(t,x))\frac{du(t,x)}{dt} = f'(u(t,x))\grad_g f(u(t,x))
  \\
  = \|f'(u(t,x))\|_{H^*}^2 = \|\grad_g f(u(t,x))\|_H^2, \quad\text{for all } t > 0.
\end{multline*}
Hence, the Fundamental Theorem of Calculus yields
\begin{equation}
  \label{eq:Gradient_flow_energy_identity}
  f(u(t'',x)) - f(u(t',x)) = \int_{t'}^{t''}\|\grad_g f(u(t,x))\|_H^2\,dt,
  \quad\text{for all } 0 \leq t' \leq t'' <\infty.
\end{equation}
In particular, for all integers $k,l\geq 0$, the identity \eqref{eq:Gradient_flow_energy_identity} yields
\[
  f(u(t_k,x)) - f(u(t_l,x)) = \int_{t_l}^{t_k}\|\grad_g f(u(t,x))\|_H^2\,dt.
\]
Setting $l=0$, applying our hypotheses that $u(t_k,x) \to u_\infty(x)$ in $H$ as $k\to\infty$ and that $f$ is $C^0$ (in fact, $C^2$), and taking the limit as $k\to\infty$ of the preceding identity yields
\[
  f(u_\infty(x)) - f(u(0,x)) = \int_0^\infty\|\grad_g f(u(t,x))\|_H^2\,dt.
\]
The definition $h(t) = \grad_g f(u(t,x))$ and the gradient flow equation \eqref{eq:Gradient_flow_equation_proposition} (see Abraham, Marsden, and Ratiu \cite[Exercise 2.4H, p. 113]{AMR} or Zeidler \cite[Section 4.5, p. 141]{Zeidler_nfaa_v1} for expressions for second tangent maps) give
\[
  \dot h(t) = f''(u(t,x))\dot u(t,x) = f''(u(t,x))\grad_g f(u(t,x)), \quad\text{for all } t > 0.
\]
Hence, by defining the self-adjoint operator $\Hess_g f(u) \in \End(H)$ via the symmetric bilinear form $f''(u) = g(\cdot,\Hess_g f(u)\cdot) \in \Hom(H\otimes H,\RR) = \Hom(H,H^*)$, we obtain
\begin{align*}
  \|\dot h(t)\|_H
  &\leq \|f''(u(t,x))\|_{\Hom(H,H^*)}\|\grad_g f(u(t,x))\|_H
  \\
  &= \|\Hess_g f(u(t,x))\|_{\End(H)}\|\grad_g f(u(t,x))\|_H.
\end{align*}
Therefore, recalling that $f'':T^2X\to \RR$ is bounded over $U\subset X$ by hypothesis,
\[
  \int_0^\infty\|\dot h(t)\|_H^2\,dt
  \leq
  \sup_{t>0}\|\Hess_g f(u(t,x))\|_{\End(H)}^2\int_0^\infty\|\grad_g f(u(t,x))\|_H^2\,dt.
\]
We conclude that $\dot h \in L^2(0,\infty;H)$ and because $h \in L^1(0,\infty;H)$, Lemma \ref{lem:Barbalat_variant} implies that $h(t) \to 0$ in $H$ as $t \to \infty$. But our hypothesis on convergence in $H$ of the subsequence $\{u(t_k,x)\}_{k=1}^\infty$ to $u_\infty(x)$ and the fact that $f$ is $C^1$ (in fact, $C^2$) by hypothesis yields
\[
  h(t_k) = \grad_g f(u(t_k,x)) \to \grad_g f(u_\infty(x)), \quad\text{as } k \to \infty.
\]
Hence, $\grad_g f(u_\infty(x)) = 0$, as claimed.
\end{proof}

We shall apply Proposition \ref{prop:Subsequential_limits_solutions_gradient_flow_critical_points} to prove

\begin{prop}[Extension of real analytic $\CC^*$ actions on real analytic, almost Hermitian manifolds]
\label{prop:Extension_C*_actions_real_analytic_almost_Hermitian_manifolds}
Let $(X,g,J)$ be a real analytic, almost Hermitian manifold with fundamental two-form $\omega = g(\cdot,J\cdot)$ as in \eqref{eq:Fundamental_two-form} and let $\CC^* \times X \to X$ be a real analytic action of $\CC^*$ by real analytic diffeomorphisms of $X$. Assume that the real analytic circle action $S^1\times X \to X$ obtained by restricting the $\CC^*$ action to $S^1\subset \CC^*$ is Hamiltonian in the sense of \eqref{eq:MomentMap}, so $\iota_\Theta\omega = df$ on $X$ for a real analytic function $f:X\to\RR$, where the vector field $\Theta$ on $X$ given by \eqref{eq:Vector_field_generator_circle_action} is the generator of the circle action. Suppose $z \in X$ and that 
\[
  A_z:\CC^* \ni \lambda \mapsto \lambda\cdot z \in X
\]
is the corresponding orbit as in \eqref{eq:Orbit_point_under_C*_action_complex_manifold}. Let $D \subset \CC$ denote an open disk centered at the origin and let $D^* := D\less\{0\} \subset \CC$ denote the punctured disk.
\begin{enumerate}
\item\label{item:AD*_limit_lambda_to_zero_exists_fixed_point}
If $A_z(D^*)$ is contained in a compact subset of $X$, then $A_z$ extends to a $\CC^*$-equivariant, continuous map, $A_z:\CC \to X$, and the limit $z_0 := \lim_{\lambda\to 0}\lambda\cdot z$ exists and is a fixed point of the $\CC^*$ action on $X$;
\item\label{item:AC-D_limit_lambda_to_infinity_exists_fixed_point}
If $A_z(\CC^*\less D)$ is contained in a compact subset of $X$, then $A_z$ extends to a $\CC^*$-equivariant, continuous map, $A_z:\CC^*\cup\{\infty\}\to X$, and the limit $z_\infty := \lim_{\lambda\to \infty}\lambda\cdot z$ exists and is a fixed point of the $\CC^*$ action on $X$.  
\end{enumerate}
In particular, if $A_z(\CC^*)$ is contained in a compact subset of $X$, then $A_z$ extends to a $\CC^*$-equivariant, continuous map, $A_z:\PP^1 \to X$, and the images $z_0 = A_z([1,0])$ and $z_\infty = A_z([0,1])$ of the points $[1,0], [0,1] \in \PP^1$ are fixed points of the $\CC^*$ action on $X$.  
\end{prop}

\begin{proof}
In our proof, we focus on Case \eqref{item:AD*_limit_lambda_to_zero_exists_fixed_point}, since the proof of Case \eqref{item:AC-D_limit_lambda_to_infinity_exists_fixed_point} is virtually identical and the final conclusion follows immediately from Cases \eqref{item:AD*_limit_lambda_to_zero_exists_fixed_point} and \eqref{item:AC-D_limit_lambda_to_infinity_exists_fixed_point}. We first make the

\begin{claim}
\label{claim:BB_limit_exists_iff_gradient_flow_limit_exists}  
The Bia{\l}ynicki--Birula limit $\lim_{\lambda\to 0}\lambda\cdot z$ (respectively, $\lim_{\lambda\to \infty}\lambda\cdot z$) exists if and only if the limit $\lim_{t\to\infty}u(t,z)$ (respectively, $\lim_{t\to -\infty}u(t,z)$) exists for a solution $u(t,z)$ to the gradient flow equation for $f$ with initial (respectively, final) condition $u(z,0) = z$.
\end{claim}

\begin{proof}[Proof of Claim \ref{claim:BB_limit_exists_iff_gradient_flow_limit_exists}]
Observe that if $r=e^{-t}$, then $dr/dt = de^{-t}/dt = -e^{-t} = -r$ and so the Chain Rule gives
\begin{equation}
  \label{eq:Vector_field_generator_temporal_to_radial_conversion_chain_rule}
  \frac{d}{dt}(e^{-t}\cdot x) = -r\frac{d}{dr}(r\cdot x), \quad\text{for all } x \in X.
\end{equation}
Recall from \eqref{eq:Gradient_Hamiltonian_and_radial_rotation_vector_fields} that
\[
  \grad_g f(x) = -R_x, \quad\text{for all } x \in X,
\]
where the (real analytic) radial vector field $R_x = d(r\cdot x)/dr|_{r=1}$, for $x\in X$, is as in \eqref{eq:Vector_field_generator_radial_action} and $g$ is the (real analytic) Riemannian metric on $X$ associated to the given (real analytic) Hermitian metric. Define a (real analytic) vector field $T$ on $X$ by
\begin{equation}
  \label{eq:Vector_field_generator_translation_action}
  T_x := \left.\frac{d}{dt}(e^{-t}\cdot x)\right|_{t=0}, \quad\text{for all } x \in X,
\end{equation}
and observe that the identity \eqref{eq:Vector_field_generator_temporal_to_radial_conversion_chain_rule} yields
\begin{equation}
  \label{eq:Vector_field_generator_temporal_to_radial_conversion}
  T_x = -R_x, \quad\text{for all } x\in X.
\end{equation}
Define $u(t,z) := e^{-t}\cdot z$, so the gradient flow equation for $u:\RR\times X\to X$ implied by \eqref{eq:Gradient_Hamiltonian_and_radial_rotation_vector_fields} is given by
\begin{equation}
  \label{eq:Gradient_flow_from_C*_action}
  \frac{du}{dt}(\cdot,z) = \grad_g f(u(t,z)), \quad\text{for all } t \in \RR, \quad u(0,z) = z.
\end{equation}
Hence, if the limit $\lim_{\lambda\to 0}\lambda\cdot z$ exists, then the limit $\lim_{t\to \infty}u(t,z)$ exists since $u(t,z) = e^{-t}\cdot z$, for all $t\in\RR$. Conversely, if the limit $\lim_{t\to \infty}u(t,z)$ exists, we may write that limit as $u_\infty(z) \in X$ and observe that it is a critical point of the Hamiltonian function $f$ by Proposition \ref{prop:Subsequential_limits_solutions_gradient_flow_critical_points}. Since
\[
  df(u_\infty(z)) = 0 \in T_{u_\infty(z)}^*X,
\]
the relation \eqref{eq:MomentMap} implies that
\[
  \omega_{u_\infty(z)}(\Theta,\cdot) = 0 \in T_{u_\infty(z)}^*X.
\]
Because $\omega$ is non-degenerate by hypothesis, we obtain $\Theta_{u_\infty(z)} = 0$ and therefore $u_\infty(z)$ is a fixed point of the $S^1$ action on $X$ (see, for example, Feehan and Leness \cite[Lemma 3.3.5 (1)]{Feehan_Leness_introduction_virtual_morse_theory_so3_monopoles}). Hence, writing $\lambda = e^{-t+i\theta}$, let $\{t_k\}_{k=1}^\infty \subset [0,\infty)$ be any sequence such that $t_k\to\infty$ as $k\to\infty$ and let $\{\theta_l\}_{l=1}^\infty \subset [0,2\pi)$ be any sequence. Then, writing $\lambda_{k,l} = e^{-t_k+i\theta_l}$, we obtain
\[
  \lim_{\lambda_{k,l}\to 0}\lambda_{k,l}\cdot z
  = \lim_{k,l\to\infty}\left(e^{i\theta_l}e^{-t_k}\right)\cdot z
  = \lim_{l\to\infty}e^{i\theta_l}\cdot\lim_{k\to\infty}u(t_k,z)
  = \lim_{l\to\infty}e^{i\theta_l}\cdot u_\infty(z)
  = u_\infty(z),
\]
where the final equality follows from the fact that $u_\infty(z)$ is a fixed point of the $S^1$ action on $X$ and thus $e^{i\theta_l}\cdot u_\infty(z) = u_\infty(z)$, for all integers $l\geq 1$. Hence, the limit $\lim_{\lambda\to 0}\lambda\cdot z$ exists, since the preceding equalities hold for every sequence $\{\lambda_{k,l}\}_{k,l=1}^\infty \subset \CC^*$ such that $\lambda_{k,l} = e^{-t_k+i\theta_l} \to 0$ as $k\to\infty$ and $l\to\infty$.

To see that the limit $\lim_{\lambda\to \infty}\lambda\cdot z$ exists if and only if the limit $\lim_{t\to-\infty}u(t,z)$ exists, we may write $\mu:=1/\lambda$ for $\lambda\in\CC^*$ and $v(s,z) := u(-s,z)$ for $s \in [0,\infty)$ and observe that this second case follows from the first. This completes the proof of Claim \ref{claim:BB_limit_exists_iff_gradient_flow_limit_exists}.
\end{proof}

Convergence of the gradient flow $u(\cdot,z):[0,\infty)\to X$ to the limit $u_\infty(z) \in X$ as $t\to\infty$ thus follows from the two facts below:
\begin{enumerate}
\item\label{item:Hamiltonian_Lojasiewcicz_gradient_inequality}
  The Hamiltonian function $f$ obeys the {\L}ojasiewicz gradient inequality,
  \begin{equation}
    \label{eq:Lojasiewicz_gradient_inequality}
    \|\grad_g f(z)\|_g \geq C|f(z)-f(u_\infty(z))|^\theta, \quad\text{for all } z \in B_\sigma(u_\infty(z)),
  \end{equation}
on an open ball $B_\sigma(u_\infty(z))$ of radius $\sigma \in (0,1]$ centered at the critical point $u_\infty(z)$, for some constants $C \in [1,\infty)$ and $\theta \in [1/2,1)$.
  
\item\label{item:Gradient_flow_solution_subsequence_convergence}
  There are a critical point $u_\infty(z) \in X$ of $f$ and a strictly increasing subsequence $\{t_k\}_{k=1}^\infty \subset [0,\infty)$ such that $u(t_k,z) \to u_\infty(z)$ as $k\to\infty$.  
\end{enumerate}
To verify Item \eqref{item:Hamiltonian_Lojasiewcicz_gradient_inequality}, note that $f$ is real analytic by hypothesis (see also Remark \ref{rmk:Real_analyticity_Hamiltonian_function}) and see Feehan \cite[Theorem 1, p. 3277]{Feehan_lojasiewicz_inequality_all_dimensions} for the statement and a modern proof of the {\L}ojasiewicz gradient inequality using Resolution of Singularities or see {\L}ojasiewicz\footnote{The first page number refers to the version of {\L}ojasiewicz's original manuscript mimeographed by IHES while the page number in parentheses refers to the cited LATEX version of his manuscript prepared by M. Coste.} \cite[Proposition 1, page 92 (67)]{Lojasiewicz_1965} for his statement and original proof using the theory of semianalytic sets.

To verify Item \eqref{item:Gradient_flow_solution_subsequence_convergence}, observe that orbit $u(z,[0,\infty)) \subset X$ is a subset of the orbit $A_z(D^*) \subset X$, which is contained in a compact subset of $X$ by hypothesis. Hence, there is a strictly increasing sequence $\{t_k\}_{k=1}^\infty \subset [0,\infty)$ such that the sequence $\{u(t_k,z)\}_{k=1}^\infty \subset X$ converges to a limit $u_\infty(z)\in X$ as $k\to\infty$.

Convergence of a global solution $u(\cdot,z):[0,\infty)\to X$ to the gradient flow equation \eqref{eq:Gradient_flow_from_C*_action} to the subsequential limit $u_\infty(z)$ now follows from Feehan \cite[Theorem 1 = Theorem 24.14]{Feehan_yang_mills_gradient_flow_v4} and Items \eqref{item:Hamiltonian_Lojasiewcicz_gradient_inequality} and \eqref{item:Gradient_flow_solution_subsequence_convergence}. The technical hypotheses of \cite[Theorem 1 = Theorem 24.14]{Feehan_yang_mills_gradient_flow_v4} are satisfied exactly as in our proof of Feehan \cite[Theorem 5.1, p. 3299]{Feehan_lojasiewicz_inequality_all_dimensions}.

We can apply Claim \ref{claim:BB_limit_exists_iff_gradient_flow_limit_exists} to establish that the Bia{\l}ynicki--Birula limit $\lim_{\lambda\to 0}\lambda\cdot z$ exists and thus $A_z$ extends to a continuous map, $A_z:\CC\to X$, and by construction, the map is clearly $\CC^*$-equivariant.

Lastly, the limits $z_0$ and $z_\infty$ are necessarily fixed points of the $\CC^*$ action on $X$. Indeed, if $\mu \in \CC^*$, then
\[ 
  \mu\cdot z_0
  = \mu\cdot\left(\lim_{\lambda\to 0} \lambda\cdot z\right)
  = \lim_{\lambda\to 0} \mu\cdot(\lambda\cdot z)
  = \lim_{\lambda\to 0} (\mu\lambda)\cdot z
  = \lim_{\lambda'\to 0} \lambda'\cdot z
  = z_0,
\]
where the penultimate equality follows by writing $\lambda' := \mu\lambda \in \CC^*$. Hence, $z_0$ is a fixed point of the $\CC^*$ action and applying this argument \mutatis shows that $z_\infty$ is a fixed point of the $\CC^*$ action. This completes the proof of Proposition \ref{prop:Extension_C*_actions_real_analytic_almost_Hermitian_manifolds}.
\end{proof}

\begin{rmk}[Real analyticity of the Hamiltonian function]
\label{rmk:Real_analyticity_Hamiltonian_function}
By definition of the Hamiltonian function $f$ in Proposition \ref{prop:Extension_C*_actions_real_analytic_almost_Hermitian_manifolds}, we have
\[
  \Delta_gf = d^{*_g}\iota_\Theta\omega \quad\text{on } X,
\]
where $\Delta_g = d^{*_g}d$ is the Laplace operator, the function on the right-hand side is real analytic since the Riemannian metric $g$, the vector field $\Theta$, the almost complex structure $J$, and hence the fundamental two-form $\omega$ are all real analytic. Therefore, $f$ is a solution to a linear, second-order, elliptic partial differential equation with real analytic coefficients and is thus a real analytic function by Morrey and Nirenberg \cite{Morrey_Nirenberg_1957}.

We recall that the existence of a smooth Hamiltonian function $f$ for a smooth $S^1$ action $S^1\times X \to X$ is provided by Remark \ref{rmk:Application_Cartan_magic_formula_existence_Hamiltonian_functions} when $(X,\omega)$ is a smooth, closed manifold with $S^1$-invariant symplectic form $\omega$ and $H^1(X;\RR) = 0$.
\qed
\end{rmk}

Proposition \ref{prop:Extension_C*_actions_real_analytic_almost_Hermitian_manifolds} leads to the following corollary and strengthening of the fundamental observation \cite[Lemma II-A, p. 109]{Sommese_1975} due to Sommese, where his requirements that $X$ be compact or K\"ahler and that the $\CC^*$ action have at least one fixed point are all relaxed. (As an aside, we note that a well-known removable-singularities result for finite-energy harmonic maps from the punctured disk into a closed Riemannian manifold is due to Sacks and Uhlenbeck \cite[Theorem 3.6, p.  13]{Sacks_Uhlenbeck_1981}.)

\begin{thm}[Extension of $\CC^*$ actions on Hermitian complex manifolds]
\label{thm:Extension_C*_actions_Hermitian_complex_manifolds}
Continue the hypotheses of Proposition \ref{prop:Extension_C*_actions_real_analytic_almost_Hermitian_manifolds}. If $J$ is integrable and so $(X,J)$ is a complex manifold, then the extension $A_z:\CC \to X$ (respectively, $A_z:\CC^*\cup\{\infty\}\to X$) is a $\CC^*$-equivariant, holomorphic map. Moreover, if the orbit $A_z(\CC^*)$ of the map
\[
  A_z:\CC^*\ni \lambda \mapsto \lambda\cdot z \in X
\]
as in \eqref{eq:Orbit_point_under_C*_action_complex_manifold} is contained in a compact subset of $X$, then $A_z$ extends to a $\CC^*$-equivariant, holomorphic map,
\[
  A_z:\PP^1 \to X,
\]
and the images $z_0 = A_z([1,0])$ and $z_\infty = A_z([0,1])$ of the points $[1,0], [0,1] \in \PP^1$ are fixed points of the $\CC^*$ action on $X$.  
\end{thm}

\begin{proof}
Because the orbit $A_z(\CC^*)$ is contained in a compact subset of $X$, then Proposition \ref{prop:Extension_C*_actions_real_analytic_almost_Hermitian_manifolds} implies that the points $z_0$ and $z_\infty$ exist and are fixed points of the $\CC^*$ action on $X$. 

Assume that $X$ has complex dimension $n$ and choose a holomorphic coordinate chart $(U_0,\varphi)$ centered at $z_0$, so that $\varphi = (\varphi_1,\ldots,\varphi_n):X \supset U_0 \to \CC^n$ obeys $\varphi(z_0)=0$, and define $f_k = \varphi_k\circ A_z:\CC \supset D_0 \to \CC$, for $k = 1,\ldots, n$, for an open disk $D_0$ centered at the origin and chosen small enough that $A_z(\bar D_0) \subset U_0$. Since $A_z(0) = z_0$ and $\varphi_k(z_0) = 0$, we obtain $f_k(0) = 0$, for $k = 1,\ldots, n$. Each function $f_k$ is continuous on $\bar D_0$ and holomorphic on $D_0^*$, so Rad\'o's Theorem (see Rudin \cite[Theorem 12.14, p. 263]{RudinRealComplex} for functions on domains in $\CC$ or \cite[Theorem 5.1.7, p. 302]{Rudin_function_theory_ball_Cn} for domains in $\CC^n$) implies that each $f_k$ is holomorphic on $D_0$. Hence, the map $\varphi\circ A_z:D_0 \to \CC^n$ is holomorphic and so $A_z:D_0 \to X$ is holomorphic.

By writing $\mu = 1/\lambda$ for $\lambda \in \CC^*$ and using a holomorphic coordinate chart $(U_\infty,\psi)$ centered at $z_\infty$, so that $\psi:X \supset U_\infty \to \CC^n$ obeys $\psi(z_\infty)=0$, the same argument applies to prove that $A_z:D_\infty \to X$ is holomorphic, for an open disk $D_\infty$ centered at the point $\infty \in \PP^1 = \CC\cup\{\infty\}$ and chosen small enough that $A_z(D_\infty) \subset U_\infty$. The $\CC^*$-equivariance of the map $A_z:\PP^1\to X$ is given by Proposition \ref{prop:Extension_C*_actions_real_analytic_almost_Hermitian_manifolds}.
\end{proof}

We now give the

\begin{proof}[Proof of Theorem \ref{mainthm:BB_decomposition_complex_manifold_C*_action}]
Consider Item \eqref{item:BB_complex_manifold_limits_z0_zinfty_exist_fixed_points}. Proposition \ref{prop:Extension_C*_actions_real_analytic_almost_Hermitian_manifolds} implies that for every point $z \in X$, at least one of the limits $\lim_{\lambda\to 0}\lambda\cdot z$ or $\lim_{\lambda\to \infty}\lambda\cdot z$ exists, the corresponding extended map $A_z:\CC\to X$ or $A_z:\CC\cup\{\infty\}\to X$ is $\CC^*$-equivariant and continuous, and the subset $X^0 \subset X$ of fixed points is non-empty. This verifies Item \eqref{item:BB_complex_manifold_limits_z0_zinfty_exist_fixed_points}.

Consider Item \eqref{item:BB_complex_manifold_Az_holomorphic_BB_decomposition_exists}. Since $J$ is integrable by hypothesis, Theorem \ref{thm:Extension_C*_actions_Hermitian_complex_manifolds} implies that the map $A_z:\CC\to X$ or $A_z:\CC\cup\{\infty\}\to X$ is also holomorphic for each $z \in X$. Lemma \ref{lem:Fixed-point_set_complex_Lie_group_on_complex_manifold} implies that $X^0$ is a complex (embedded) submanifold of $X$
and so has at most countably connected components (of possibly different dimensions), $X_\alpha^0$ for $\alpha\in\sA$. (Indeed, every topological manifold of dimension $d$ is a disjoint union of at most countably many connected manifolds of dimension $d$ by Lee \cite[Problem 4.9, p. 123]{Lee_john_topological_manifolds} and we may apply this property to each submanifold of $X^0$ with complex dimension $d$ between zero and $\dim X$.)
Consequently, the subsets $X_\alpha^\pm$ in \eqref{eq:Complex_manifold_Xalpha_plus_minus_submanifolds} and the maps $\pi_\alpha^\pm:X_\alpha^\pm\to X_\alpha^0$ in \eqref{eq:pi_Xalpha_pm_to_Xalpha_0_bundle_maps} are well-defined and, as a topological space, $X$ has a mixed, plus, or minus decomposition as in Item \eqref{item:BB_decomposition_complex_manifold_mixed_plus_minus} of Definition \ref{maindefn:BB_decomposition_complex_manifold}. Moreover, their construction ensures that the maps $\pi_\alpha^\pm$ are projections that define continuous, complex vector bundles, $X_\alpha^\pm\to X_\alpha^0$, with zero sections $X_\alpha^0$ and this yields Item \eqref{item:Xalpha_0_section_Xalpha_pm} and part of Item \eqref{item:pi_Xalpha_pm_to_Xalpha_0_bundle}. In particular, this verifies the global properties of the Bia{\l}ynicki--Birula decomposition in Definition \ref{maindefn:BB_decomposition_complex_manifold}.

The remaining properties of the Bia{\l}ynicki--Birula decomposition in Definition \ref{maindefn:BB_decomposition_complex_manifold} are local. We can thus appeal to Theorem \ref{thm:Bialynicki-Birula_decomposition_compact_complex_Kaehler_manifold} and especially its proof by Fujiki, to justify these local properties of the Bia{\l}ynicki--Birula decomposition in Definition \ref{maindefn:BB_decomposition_complex_manifold}. As explained in Remark \ref{rmk:BB_decomposition_complex_manifold_C*_action_previous_versions}, Fujiki's proof does not require $X$ to be K\"ahler or compact under the hypotheses of Theorem \ref{mainthm:BB_decomposition_complex_manifold_C*_action} and Item \eqref{item:BB_complex_manifold_Az_holomorphic_BB_decomposition_exists}. This verifies Item \eqref{item:BB_complex_manifold_Az_holomorphic_BB_decomposition_exists} and completes the proof of Theorem \ref{mainthm:BB_decomposition_complex_manifold_C*_action}.
\end{proof}


\section{Functorial properties of Bia{\l}ynicki--Birula decompositions}
\label{sec:Functorial_properties_Bialynicki-Birula_decomposition}
In this section, we discuss certain functorial properties of Bia{\l}ynicki--Birula decompositions that allow one to easily construct Bia{\l}ynicki--Birula decompositions for new complex manifolds with little additional work, using such decompositions for given complex manifolds. Before proceeding, we recall some facts regarding the decomposition of complex analytic spaces or sets into connected components, closely following Grauert and Remmert \cite[Sections 9.2.1 and 9.2.2, pp. pp. 171--174]{Grauert_Remmert_coherent_analytic_sheaves}.

Every topological space $X$ is the set-theoretic union of its connected components $\{X_\alpha\}_{\alpha\in\sA}$, where each set $X_\alpha$ is a maximal connected subset of $X$. These sets are pairwise disjoint and closed in $X$, but not necessarily 
open in $X$. To see this, we recall from Munkres \cite[Theorem 23.4, p. 150]{Munkres_topology_second_edition} that if $A$ is a connected subspace of a topological space $X$ and $\bar A$ is its closure and $A\subseteq B\subseteq \bar A$, then $B$ is connected. Thus, each connected component $X_\alpha$ of a topological space $X$ is closed, since the closure $\bar X_\alpha$ of a connected component $X_\alpha$ is closed by \cite[Theorem 23.4, p. 150]{Munkres_topology_second_edition}. If a topological space $X$ has only finitely many connected components, then each connected component is also open in $X$, since its complement is a finite union of closed sets (see \cite[Section 25, p. 160]{Munkres_topology_second_edition}).

\begin{rmk}[Connected components of manifolds]
\label{rmk:Connected_components_manifolds}
If $X$ is a smooth or $\KK$-analytic manifold (for $\KK=\RR$ or $\CC$) with \emph{finitely} many connected components $\{X_\alpha\}_{\alpha\in\sA}$, then each connected component $X_\alpha$ is an open subset of $X$ by the preceding discussion and thus is an embedded smooth or $\KK$-analytic submanifold of $X$, respectively.
\end{rmk}

More generally, if $X$ is a $\KK$-analytic space, a discussion of the properties of its connected components becomes more involved. According to Grauert and Remmert \cite[Section 9.1.2, Theorem, p. 168]{Grauert_Remmert_coherent_analytic_sheaves}, a reduced complex analytic space $X$ can be defined to be \emph{(globally) irreducible} if $X\less X_\sing$ is connected, where $X_\sing \subset X$ is the singular set (as in Definition \ref{defn:Singular_locus_analytic_space}).

\begin{defn}[Decomposition of a complex analytic space into irreducible components]
\label{defn:Grauert_Remmert_9-2-2_definition}
(See Grauert and Remmert \cite[Section 9.2.2, Definition, p. 172]{Grauert_Remmert_coherent_analytic_sheaves}.)  
Let $(X,\sO_X)$ be a complex analytic space that is reduced (as in Definition \ref{defn:Reduced_locally_ringed_space}). A family $\{X_\alpha\}_{\alpha\in\sA}$ of irreducible complex analytic subsets of $X$ is called a \emph{decomposition of $X$ into irreducible components} if the following hold:
\begin{enumerate}  
\item $\{X_\alpha\}_{\alpha\in\sA}$ is a locally finite covering of $X$, and
\item $X_\alpha\not\subset X_\beta$ for all $\alpha\in\sA$ with $\alpha\neq\beta$. 
\end{enumerate}
\end{defn}

\begin{thm}[Global decomposition of a reduced complex analytic space into irreducible components]
\label{thm:Grauert_Remmert_9-2-2_theorem}
(See Grauert and Remmert \cite[Section 9.2.2, Theorem, p. 172]{Grauert_Remmert_coherent_analytic_sheaves}.)
Let $(X,\sO_X)$ be a reduced complex analytic space. Then there exists a unique decomposition of $X$ into irreducible components and every irreducible component of $X$ contains smooth points of $X$.
\end{thm}

Grauert and Remmert give an example \cite[Section 9.2.2, Theorem, p. 173]{Grauert_Remmert_coherent_analytic_sheaves} to illustrate their observation that irreducible components of $X$ are \emph{not} necessarily open in $X$. We now proceed with the

\begin{lem}[Bia{\l}ynicki--Birula decomposition for a $\CC^*$-invariant, topologically closed, complex analytic subspace: existence of subsets]
\label{lem:BB_decomposition_C*_invariant_complex_submanifold_subset_existence}
Let $X$ be a complex manifold that admits a holomorphic $\CC^*$ action, $\CC^*\times X\to X$, and let $(Y,\sO_Y)$ be a locally closed analytic subspace of $X$ in the sense of Definition \ref{defn:Closed_analytic_subspace_of_analytic_space} and such that $Y\subset X$ is a topologically closed subspace and $\CC^*$-invariant with at least one fixed point. If $X$ admits a mixed Bia{\l}ynicki--Birula decomposition in the sense of Definition \ref{maindefn:BB_decomposition_complex_manifold}, then the following hold:
\begin{subequations}
\label{eq:Y0_Ypm_identity}  
\begin{align}
  \label{eq:Y0_identity}
  Y^0 &= Y\cap X^0,
  \\
  \label{eq:Ypm_identities}
  Y^\pm &= Y\cap X^\pm,
\end{align}
\end{subequations}
where $Y^0$ is the set of fixed points of the $\CC^*$ action on $Y$ and
\begin{subequations}
  \label{eq:Xpm_Ypm}
  \begin{align}
  \label{eq:X+}
  X^+ &:= \left\{x\in X: \lim_{\lambda\to 0} \lambda\cdot x \in X^0 \right\},
    \\
    \label{eq:X-}
  X^- &:= \left\{x\in X: \lim_{\lambda\to \infty} \lambda\cdot x \in X^0 \right\},
    \\
    \label{eq:Y+}
  Y^+ &:= \left\{y\in Y: \lim_{\lambda\to 0} \lambda\cdot y \in Y^0 \right\},
    \\
    \label{eq:Y-}
  Y^- &:= \left\{y\in Y: \lim_{\lambda\to \infty} \lambda\cdot y \in Y^0 \right\}.
\end{align}
\end{subequations}
Furthermore, if $Y^0 = \sqcup_{\alpha,\beta}Y_{\alpha\beta}^0$, where the subsets
\begin{equation}
  \label{eq:Yab0}
  Y_{\alpha\beta}^0 := \left(Y\cap X_\alpha^0\right)_\beta
\end{equation}
denote the connected components of $Y^0$, and
\begin{subequations}
  \begin{align}
  \label{eq:Yab+} 
  Y_{\alpha\beta}^+ &:= \left\{y\in Y: \lim_{\lambda\to 0} \lambda\cdot y \in Y_{\alpha\beta}^0 \right\},
    \\
  \label{eq:Yab-}  
  Y_{\alpha\beta}^- &:= \left\{y\in Y: \lim_{\lambda\to \infty} \lambda\cdot y \in Y_{\alpha\beta}^0 \right\},
\end{align}
\end{subequations}
then 
\begin{equation}
  \label{eq:Yab_pm_identities}
  Y_{\alpha\beta}^\pm = \left(\pi^\pm\right)^{-1}\left(Y_{\alpha\beta}^0\right),
\end{equation}
with projections
\begin{equation}
  \label{eq:pi_ab_pm}
  \pi_{\alpha\beta}^\pm:Y_{\alpha\beta}^\pm \to Y_{\alpha\beta}^0
\end{equation}
given by the restrictions of the projections $\pi^\pm:X^\pm \to X^0$ to the subsets $Y_{\alpha\beta}^0 \subset X^0$ defined by
\begin{subequations}
  \begin{align}
  \label{eq:pi+}  
  \pi^+:X^+ \ni x \mapsto \lim_{\lambda\to 0} \lambda\cdot x \in X^0,
    \\
    \label{eq:pi-}  
  \pi^-:X^- \ni x \mapsto \lim_{\lambda\to \infty} \lambda\cdot x \in X^0.
\end{align}
\end{subequations}
\end{lem}

\begin{proof}
Since $Y$ is a $\CC^*$-invariant subset of $X$, Lemma \ref{lem:Fixed_points_group_actions_on_subsets_and_invariant_subsets} implies that $Y^0 = Y\cap X^0$. We label the connected components of $Y^0$ as indicated, where we allow for the fact that each subset $Y\cap X_\alpha^0$ need not be connected, even though (by definition) each subset $X_\alpha^0$ is connected. If $z \in Y\cap X^+$, then by definition \eqref{eq:X+} of $X^+$ we have
\[
  \lim_{\lambda\to 0} \lambda\cdot z \in X^0.
\]
Because $Y$ is $\CC^*$-invariant, we have $\lambda\cdot z \in Y$ for all $\lambda \in \CC^*$ since $z\in Y$ and $Y$ is $\CC^*$-invariant and thus, since our hypotheses assert that $Y$ is a topologically closed subset of $X$,
\[
  \lim_{\lambda\to 0} \lambda\cdot z \in Y.
\]
Hence, we obtain
\[
  \lim_{\lambda\to 0} \lambda\cdot z \in Y\cap X^0
\]
and therefore $z \in Y^+$ by definition \eqref{eq:Y+} of $Y^+$. We conclude that
\[
  Y\cap X^+ \subseteq Y^+.
\]
Conversely, if $z \in Y^+$, then by definition \eqref{eq:Y+} of $Y^+$ and the fact that $Y^0 = Y\cap X^0$, we have
\[
  \lim_{\lambda\to 0} \lambda\cdot z \in Y\cap X^0 \subset X^0
\]
and thus $z \in X^+$ by definition \eqref{eq:X+} of $X^+$. We conclude that
\[
   Y^+ \subseteq Y\cap X^+
\]
and hence that $Y^+ = Y\cap X^+$. An identical argument yields the equality $Y^- = Y\cap X^-$.

If $z \in Y_{\alpha\beta}^+$, then the definition \eqref{eq:Yab+} of $Y_{\alpha\beta}^+$ and the fact that $Y_{\alpha\beta}^0 \subset X^0$ implies that $z \in X^+$ by definition \eqref{eq:X+} of $X^+$. Moreover, the definition \eqref{eq:Yab+} of $Y_{\alpha\beta}^+$ and the definition \eqref{eq:pi+} of $\pi^+$ imply that $z \in (\pi^+)^{-1}(Y_{\alpha\beta}^0)$ and thus
\[
  Y_{\alpha\beta}^+ \subseteq \left(\pi^+\right)^{-1}\left(Y_{\alpha\beta}^0\right).
\]
Conversely, if $z \in (\pi^+)^{-1}(Y_{\alpha\beta}^0)$, then the definition \eqref{eq:pi+} of $\pi^+$ and the definition \eqref{eq:Yab+} of $Y_{\alpha\beta}^+$ imply that $z \in Y_{\alpha\beta}^+$. We conclude that
\[
  \left(\pi^+\right)^{-1}\left(Y_{\alpha\beta}^0\right) \subseteq Y_{\alpha\beta}^+
\]
and hence that $Y_{\alpha\beta}^+ = (\pi^+)^{-1}(Y_{\alpha\beta}^0)$. An identical argument yields the equality $Y_{\alpha\beta}^- = (\pi^-)^{-1}(Y_{\alpha\beta}^0)$ and this completes the proof.
\end{proof}

At this point, we can complete the

\begin{proof}[Proof of Theorem \ref{mainthm:BB_decomposition_C*_invariant_complex_analytic_subspace}]
Consider Item \eqref{item:BB_decomposition_C*_invariant_complex_analytic_subspace}. The existence of the subsets $Y^\pm \subset Y$ and $Y_p^\pm \subset Y$, for all $p \in Y^0$, and the equalities \eqref{eq:BB_decomposition_C*_invariant_complex_analytic_subspace} follow from Lemma \ref{lem:BB_decomposition_C*_invariant_complex_submanifold_subset_existence}. The fact that the subsets $Y^0$, $Y^\pm$, and $Y_p^\pm$, for all $p \in Y^0$, are locally closed complex analytic subspaces of $Y$ follows from their expressions as intersections in \eqref{eq:BB_decomposition_C*_invariant_complex_analytic_subspace}, from the fact that $X^0$, $X^\pm$, and $X_p^\pm$ are embedded complex submanifolds of $X$ by Definition \ref{maindefn:BB_decomposition_complex_manifold}, and from Definition \ref{defn:Analytic_subspace}. This proves Item \eqref{item:BB_decomposition_C*_invariant_complex_analytic_subspace}.

The assertion in Item \eqref{item:BB_decomposition_C*_invariant_complex_analytic_subspace_smooth_point} on the property of $\dim\sO_{Y,p}$ restates a conclusion of Theorem \ref{thm:Narasimhan_section_3-1_theorem_1_p_41}.

Consider Item \eqref{item:BB_dimensions_nullity_co-index_index_complex_analytic_subspace}. The first equality in each one of the dimension formulae \eqref{eq:BB_dimensions_nullity_co-index_index_complex_analytic_subspace} for the Bia{\l}ynicki--Birula nullity, co-index, and index of $p$ in $Y$ restates the corresponding equality in Definition \ref{maindefn:Stable_unstable_subspaces_BB_index_co-index_nullity}. The second equality in each one of the formulae \eqref{eq:BB_dimensions_nullity_co-index_index_complex_analytic_subspace} follows from Theorem \ref{thm:Narasimhan_section_3-1_theorem_1_p_41}. This verifies Item \eqref{item:BB_dimensions_nullity_co-index_index_complex_analytic_subspace}. The assertions in Items \eqref{item:BB_decomposition_C*_invariant_complex_analytic_subspace_nonempty_BB_subspaces} and \eqref{item:BB_index_coindex_at_p_positive_basic_implies_not_local_min_max} on implications of positivity of the integers $\beta_Y^0(p)$ or $\beta_Y^\pm(p)$ follow from Theorem \ref{thm:Narasimhan_section_3-1_theorem_1_p_41}.

Consider Item \eqref{item:BB_index_coindex_at_p_positive_basic_implies_not_local_min_max}. The Item \eqref{item:Frankel_almost_Hermitian_f_is_MB} in Frankel's Theorem \ref{mainthm:Frankel_almost_Hermitian} tells us that $f:X\to\RR$ is Morse--Bott at $p$ in the sense of Definition \ref{maindefn:Morse-Bott_function} and Item \eqref{item:Frankel_almost_Hermitian_WeightsAreEigenvalues} in Theorem \ref{mainthm:Frankel_almost_Hermitian} tells us that the eigenvalues of $\Hess_g f \in \End(T_pM)$ are given by the weights of the circle action on $T_pX$. That $S^1$ action on the tangent space $T_pX$ is induced by the $\CC^*$ action on $X$. Suppose that the Bia{\l}ynick--Birula index $\beta_Y^-(p)$ is positive. Since $\beta_Y^-(p) = \dim((Y_p^-)_\sm\cap U)$ by Item \eqref{item:BB_dimensions_nullity_co-index_index_complex_analytic_subspace}, we see that $(Y_p^-)_\sm\cap U$ is non-empty. By Item \eqref{item:BB_decomposition_C*_invariant_complex_analytic_subspace}, we know that $Y_p^- = Y\cap X_p^-$ and so $X_p^-$ must have positive dimension near $p$ as well. 
But $X_p^-$ is the smoothly embedded unstable manifold for $f$ through the point $p$ and because there existsa point $y\in Y\cap X_p^-$ with $y\neq p$, we necessarily have $f(y)<f(p)$. Thus, $p$ is not a local minimum of the restriction $f: Y \to \RR$. An almost identical argument shows that if Bia{\l}ynicki--Birula co-index $\beta_Y^+(p)$ is positive, then $p$ is not a local maximum of the restriction $f: Y \to \RR$. This completes the verification of Item \eqref{item:BB_index_coindex_at_p_positive_basic_implies_not_local_min_max} and hence the proof of Theorem \ref{mainthm:BB_decomposition_C*_invariant_complex_analytic_subspace}.
\end{proof}

We now consider a setting where we can expect the more refined properties of a Bia{\l}ynicki--Birula decomposition to hold as in Definition \ref{maindefn:BB_decomposition_complex_manifold}. Before proceeeding, we recall the

\begin{thm}[Remmert--Stein extension theorem for complex analytic sets]
\label{thm:Remmert-Stein_extension_theorem_analytic_sets}  
(See Grauert and Remmert \cite[Section 9.4.2, Extension Theorem for Analytic Sets, p. 181]{Grauert_Remmert_coherent_analytic_sheaves} or Remmert and Stein \cite{Remmert_Stein_1953}.)
Let $X$ be a complex analytic space, $d$ be an integer, and $T$ be a complex analytic set in $X$ with $\dim T < d$. If $A$ is a complex analytic set in $X\less T$ which has $\dim_pA \geq d$ for all points $p \in A$, then the topological closure $\bar A$ of $A$ in $X$ is a complex analytic subset of $X$ as in Definition \ref{defn:Analytic_set}  
\end{thm}

Bishop \cite{Bishop_1964} proved several generalizations of Theorem \ref{thm:Remmert-Stein_extension_theorem_analytic_sets} and Aguilar and Verjovsky \cite{Aguilar_Verjovsky_2021arxiv} provide a recent exposition of those results. See also Narasimhan \cite[Chapter III, Theorem 4, p. 47 and Corollary 3, p. 55]{Narasimhan_introduction_theory_analytic_spaces} for variants of Theorem \ref{thm:Remmert-Stein_extension_theorem_analytic_sets}.

\begin{lem}[Bia{\l}ynicki--Birula decomposition for a $\CC^*$-invariant, properly embedded complex submanifold: properties of subsets]
\label{lem:BB_decomposition_C*_invariant_complex_submanifold_subset_properties}
Continue the hypotheses of Lemma \ref{lem:BB_decomposition_C*_invariant_complex_submanifold_subset_existence} but assume now that $Y$ is a $\CC^*$-invariant, properly embedded complex submanifold. Then the subsets $Y^0$ and $\{Y_{\alpha\beta}^\pm\}$ comprise a Bia{\l}ynicki--Birula decomposition for $Y$ in the sense of Definition \ref{maindefn:BB_decomposition_complex_manifold}.
\end{lem}

\begin{proof}
Recall from Lee \cite[Chapter 5, p. 100]{Lee_john_smooth_manifolds} that an embedded submanifold $S$ of a smooth manifold $M$ is \emph{properly embedded} if the inclusion map $S \hookrightarrow M$ is a proper map. By Lee \cite[Proposition 5.5, p. 100]{Lee_john_smooth_manifolds}, if $S \hookrightarrow M$ is an embedded submanifold, then $S$ is properly embedded if and only if $S$ is a topologically closed subset of $M$. Thus, our hypotheses imply that $Y$ is a topologically closed subset of $X$.  

By hypothesis of Lemma \ref{lem:BB_decomposition_C*_invariant_complex_submanifold_subset_existence}, each map $\pi_\alpha^+:X_\alpha^+ \to X_\alpha^0$ is a $\CC^*$-equivariant, holomorphic, maximal-rank surjection
whose fibers are vector spaces. By Lemma \ref{lem:Fixed-point_set_complex_Lie_group_on_complex_manifold}, we see that $Y^0 = Y^{\CC^*}$ is an embedded complex submanifold of $Y$. Thus, Remark \ref{rmk:Connected_components_manifolds} implies that each connected component $Y_{\alpha\beta}^0$ of $Y^0$ is an embedded complex submanifold of $Y$. By hypothesis, $Y$ is an embedded complex submanifold of $X$, so we can conclude that $Y_{\alpha\beta}^0$ is an embedded complex submanifold of $X$. Moreover, $Y_{\alpha\beta}^0 = (Y\cap X_\alpha^0)_\beta \subset X_\alpha^0$ and $X_\alpha^0$ is an embedded complex submanifold of $X$, so $Y_{\alpha\beta}^0$ is an embedded complex submanifold of $X_\alpha^0$. The projection $\pi_{\alpha\beta}^+:Y_{\alpha\beta}^+ \to Y_{\alpha\beta}^0$ is the restriction of the projection $\pi_\alpha^+:X_\alpha^+ \to X_\alpha^0$ to the embedded complex submanifold $Y_{\alpha\beta}^0 \subset  X_\alpha^0$ and therefore $Y_{\alpha\beta}^+ = (\pi_\alpha^+)^{-1}(Y_{\alpha\beta}^0)$ is an embedded complex submanifold of $Y$. (See, for example, Lee \cite[Corollary 6.31, p. 144]{Lee_john_smooth_manifolds} in the case of smooth manifolds and note that the corresponding result for complex manifolds is obtained by replacing the role of the Implicit Mapping Theorem for smooth maps by the corresponding result for holomorphic maps.) Hence, each map $\pi_{\alpha\beta}^+:Y_{\alpha\beta}^+ \to Y_{\alpha\beta}^0$ is a $\CC^*$-equivariant, holomorphic, maximal-rank surjection whose fibers are vector spaces and $Y_{\alpha\beta}^0$ is a section of $Y_{\alpha\beta}^+$.

Let $N_{Y_{\alpha\beta}^0/Y_{\alpha\beta}^+}$ denote the the normal bundle of $Y_{\alpha\beta}^0$ in $Y_{\alpha\beta}^+$, and $N_{Y_{\alpha\beta}^0/Y_{\alpha\beta}^-}$ denote the the normal bundle of $Y_{\alpha\beta}^0$ in $Y_{\alpha\beta}^-$, and $N_{Y_{\alpha\beta}^0/Y}$ denote the the normal bundle of $Y_{\alpha\beta}^0$ in $Y$. Let
\[
  TY\restriction Y_{\alpha\beta}^0
  =
  TY_{\alpha\beta}^0 \oplus N_{Y_{\alpha\beta}^0/Y}^+ \oplus N_{Y_{\alpha\beta}^0/Y}^- 
\]
be the weight-sign decomposition of $TY\restriction Y_{\alpha\beta}^0$ defined by the $S^1$ action on $Y$ associated to the $\CC^*$ action. By our hypothesis of a Bia{\l}ynicki--Birula decomposition for $X$, we have
\[
  TX_\alpha^+\restriction X_\alpha^0
  =
  TX_\alpha^0 \oplus N_{X_\alpha^0/X}^+
\]
and by naturality, it follows that
\[
  TY_{\alpha\beta}^+\restriction Y_{\alpha\beta}^0
  =
  TY_{\alpha\beta}^0 \oplus N_{Y_{\alpha\beta}^0/Y}^+.
\]
Similarly, we have
\[
  TY_{\alpha\beta}^-\restriction Y_{\alpha\beta}^0
  =
  TY_{\alpha\beta}^0 \oplus N_{Y_{\alpha\beta}^0/Y}^-.
\]
Consequently,
\[
  N_{Y_{\alpha\beta}^0/Y_{\alpha\beta}^+} = N_{Y_{\alpha\beta}^0/Y}^+
  \quad\text{and}\quad
  N_{Y_{\alpha\beta}^0/Y_{\alpha\beta}^-} = N_{Y_{\alpha\beta}^0/Y}^-.
\]
Note that $Y_{\alpha\beta}^+ = (\pi_\alpha^+)^{-1}(Y_{\alpha\beta}^0)$ by an argument that is almost identical to the proof of the identities \eqref{eq:Yab_pm_identities}, where $\pi_\alpha^+: X_\alpha^+ \to X_\alpha^0$ is the natural projection, so $Y_{\alpha\beta}^+ \subset X_\alpha^+$.


Next, we verify that the sets $Y_{\alpha\beta}^\pm$ have the properties indicated by Items \eqref{item:Closure_Xalpha_pm_complex_analytic_subvariety_X} and \eqref{item:Xalpha_pm_Zariski_open_in_barXalpha_pm} in Definition \ref{maindefn:BB_decomposition_complex_manifold}. Suppose first that the sets $X^0$ and $Y^0=Y\cap X^0$ are connected, so we may omit the connected component labels $\alpha,\beta$, and consider $Y^+ = Y\cap X^+$. The closure of $Y^+$ in $Y$ is 
\[
  \bar Y^+ = \overline{Y\cap X^+} = Y\cap \bar X^+.
\]
Thus, $\bar Y^+$ is a complex analytic subvariety of $Y$ and $Y^+$ is Zariski-open in $\bar Y^+$, since $\bar X^+$ is a complex analytic subvariety of $X$ and $X^+$ is Zariski-open in $\bar X^+$. In general, we have
\[
  Y_{\alpha\beta}^+ = (Y\cap X_\alpha^+)_\beta,
\]
where the label $\beta$ indicates a restriction of the fiber bundle $X_\alpha^+\to X_\alpha^0$ to the component $Y_{\alpha\beta}^0 = (Y\cap X_\alpha^0)_\beta$. By repeating the preceding argument for each connected component $Y_{\alpha\beta}^0$ of $Y^0$, we see that the closure of $Y_{\alpha\beta}^+$ in $Y$ is equal to
\[
  \bar Y_{\alpha\beta}^+ = \overline{(Y\cap X_\alpha^+)_\beta} = (Y\cap \bar X_\alpha^+)_\beta.
\]
Thus, in this more general case, $\bar Y_{\alpha\beta}^+$ is a complex analytic subvariety of $Y$ and $Y_{\alpha\beta}^+$ is Zariski-open in $\bar Y_{\alpha\beta}^+$, since $\bar X_\alpha^+$ is a complex analytic subvariety of $X$ and $X_\alpha^+$ is Zariski-open in $\bar X_\alpha^+$.

Naturally, the same argument applies to $Y_{\alpha\beta}^-$ and $\bar Y_{\alpha\beta}^-$. This completes our verification of the properties of the subsets $Y_{\alpha\beta}^0$ and $Y_{\alpha\beta}^\pm$ required for existence of a Bia{\l}ynicki--Birula decomposition for $Y$.
\end{proof}

We next complete the

\begin{proof}[Proof of Theorem \ref{mainthm:BB_decomposition_C*_invariant_complex_submanifold}]
The conclusions are provided by Lemmas \ref{lem:BB_decomposition_C*_invariant_complex_submanifold_subset_existence} and \ref{lem:BB_decomposition_C*_invariant_complex_submanifold_subset_properties}. 
\end{proof}

We return to our development of functorial properties of Bia{\l}ynicki--Birula decompositions.

\begin{lem}[Bia{\l}ynicki--Birula decomposition for the product of two complex manifolds]
\label{lem:BB_decomposition_product_complex_manifolds}
Let $X, Y$ be finite-dimensional complex manifolds with holomorphic $\CC^*$ actions, $\CC^*\times X\to X$ and $\CC^*\times Y\to Y$, that each admit at least one fixed point. If $X$ and $Y$ admit Bia{\l}ynicki--Birula decompositions in the sense of Definition \ref{maindefn:BB_decomposition_complex_manifold}, then the product complex manifold $X\times Y$ admits a Bia{\l}ynicki--Birula decomposition for the induced holomorphic $\CC^*$ action, $\CC^*\times X\times Y \to X\times Y$.
\end{lem}

\begin{proof}
The subsets of $X\times Y$ that define its Bia{\l}ynicki--Birula decomposition are given by
\begin{align*}
  (X\times Y)^0 &= X^0\times Y^0,
  \\
  (X\times Y)^+ &= X^+\times Y^+,
  \\
  (X\times Y)^- &= X^-\times Y^-.
\end{align*}
Let $\{X_\alpha^0\}$ and $\{Y_\beta^0\}$ denote the connected components of $X$ and $Y$, respectively. According to Munkres \cite[Theorem 23.6, p. 150]{Munkres_topology_second_edition}, the sets $\{X_\alpha^0\times Y_\beta^0\}$ are connected since the sets $X_\alpha^0$ and $Y_\beta^0$ are connected. On the other hand, if $\alpha'\neq\alpha''$, then the union of the sets
\[
  X_{\alpha'}^0\times Y_\beta^0 \quad\text{and}\quad X_{\alpha''}^0\times Y_\beta^0
\]
is not connected for any $\beta$ and similarly, if $\beta'\neq\beta''$, then the union of the sets
\[
  X_\alpha^0\times Y_{\beta'}^0 \quad\text{and}\quad X_\alpha^0\times Y_{\beta''}^0
\]
is not connected for any $\alpha$. We conclude that the set $\{X_\alpha^0\times Y_\beta^0\}$ comprises the set of connected components of $X\times Y$.

The projections $\pi_\alpha^\pm:X_\alpha^\pm\to X_\alpha^0$ and $\pi_\beta^\pm:Y_\beta^\pm\to Y_\beta^0$ given by the Bia{\l}ynicki--Birula decompositions for $X$ and $Y$, respectively, yield projections
\[
  \pi_{\alpha\beta}^\pm: (X\times Y)_{\alpha\beta}^\pm \to (X\times Y)_{\alpha\beta}^0,
\]
where $\pi_{\alpha\beta}^\pm = \pi_\alpha^\pm\times \pi_\beta^\pm$ and
\begin{align*}
  (X\times Y)_{\alpha\beta}^0 &= X_\alpha^0\times Y_\beta^0,
  \\
  (X\times Y)_{\alpha\beta}^+ &= X_\alpha^+\times Y_\beta^+,
  \\
  (X\times Y)_{\alpha\beta}^- &= X_\alpha^-\times Y_\beta^-.
\end{align*}
The desired properties of the projections $\pi_{\alpha\beta}^\pm$ and sets $(X\times Y)_{\alpha\beta}^\pm$ are almost immediate consequences of the properties of the projections $\pi_\alpha^\pm$, $\pi_\beta^\pm$, and sets $X_\alpha^\pm$, $Y_\beta^\pm$. Hence, $X\times Y$ admits a Bia{\l}ynicki--Birula decomposition, as claimed.
\end{proof}

We note one further elementary functorial property of Bia{\l}ynicki--Birula decompositions.

\begin{lem}[Equivariant, holomorphic maps preserve Bia{\l}ynicki--Birula decompositions]
\label{lem:Equivariant_holomorphic_maps_preserve_BB_decompositions}  
Let $(X,\sO_X)$ and $(Y,\sO_Y)$ be complex analytic spaces as in Definition \ref{defn:Analytic_space} and assume that they have holomorphic $\CC^*$ actions, $\CC^*\times X\to X$ and $\CC^*\times Y\to Y$, that each admit at least one fixed point. If $X$ and $Y$ admit Bia{\l}ynicki--Birula decompositions in the sense of Definition \ref{maindefn:BB_decomposition_complex_analytic_space} and $F:X\to Y$ is a $\CC^*$-equivariant morphism of complex analytic spaces, then we have the following inclusions:
\[
  F(X^0) \subset Y^0, \quad F(X^\pm) \subset Y^\pm,
  \quad\text{and}\quad
  F(X_p^\pm) \subset Y_{F(p)}^\pm, \quad\text{for all } p \in X^0.
\]
\end{lem}

\begin{proof}
If $p \in X^0$, then $\lambda\cdot F(p) = F(\lambda\cdot p) = F(p)$, for all $\lambda \in \CC^*$, and so $F(p) \in Y^0$, which proves the first inclusion. If $x \in X_p^+$, so $\lim_{\lambda\to 0}\lambda\cdot x = p \in X^0$, then
\[
  \lim_{\lambda\to 0}\lambda\cdot F(x) = F\left(\lim_{\lambda\to 0}\lambda\cdot x\right) = F(p) \in Y^0,
\]
and so $F(x) \in Y_{F(p)}^+$. An almost identical argument shows that if $x \in X_p^-$, then $F(x) \in Y_{F(p)}^-$. This proves the second and third pairs of inclusions.   
\end{proof}

We can now proceed to complete the

\begin{proof}[Proof of Theorem \ref{mainthm:Upper_and_lower_bounds_Krull_dimensions_X0_Xpm_at_p}]
Since the problem is local, we may assume without loss of generality that $(X,\sO_X)$ is a local model space as in Definition \ref{defn:Analytic_model_space}, so $(X,\sO_X)$ is defined by a domain $D \subset \CC^n$ around the origin with $n=\dim T_pX$, an ideal $\sI \subset \sO_D$ with generators $f_1,\ldots,f_r$, support $X = \supp(\sO_D/\sI) \subset D$, and structure sheaf $\sO_X = (\sO_D/\sI)\restriction X$.
  
By Lemma \ref{lem:Blanchard_holomorphic_map_from_domain_around_fixed_point_into_vector_space}, since the origin $p\in X$ is a fixed point of the $S^1$ action, we may further assume without loss of generality that the induced action of $S^1 \subset \CC^*$ on $D \subset \CC^n$ is linear, that $D$ is $S^1$-invariant, and that the holomorphic map $F = (f_1,\ldots,f_r):D \to \Xi_p$ determines a linear (unitary) $S^1$ action on $\Xi_p = \CC^r$ such that $F$ is $S^1$-equivariant. The unitary action of $S^1$ on $\Xi_p$ uniquely determines a linear $\CC^*$ action on $\Xi_p$. According to Corollary \ref{cor:Circle-equivariant_complex_analytic_maps_preserve_weight-sign_decompositions} and identifying $T_pX = \CC^n$, the $S^1$-equivariant, holomorphic map $F:T_pX \supset D \to \Xi_p$ obeys the inclusions
\[
  F(D\cap T_p^0X) \subset \Xi^0 \quad\text{and}\quad F(D\cap T_p^\pm X) \subset \Xi_p^\pm.
\]
Since the complex vector space $T_pX$ admits a Bia{\l}ynicki--Birula decomposition, as described above, we can appeal to Lemma \ref{lem:Equivariant_holomorphic_maps_preserve_BB_decompositions} to see that the $\CC^*$-equivariant, complex analytic monomorphism $\iota:X \to D \subset T_pX$ obeys the inclusions
\[
  X^0 = \iota(X^0) \subset T_p^0X \quad\text{and}\quad X_p^\pm = \iota(X_p^\pm) \subset T_p^\pm X.
\]
By hypothesis of Theorem \ref{mainthm:Upper_and_lower_bounds_Krull_dimensions_X0_Xpm_at_p}, the complex analytic space $(X,\sO_X)$ admits a Bia{\l}ynicki--Birula decomposition as in Definition \ref{maindefn:BB_decomposition_complex_analytic_space} and, of course, so does the complex vector space $\Xi_p$, as described above. Hence, by composing these maps, wee that the $S^1$-equivariant, complex analytic morphism $F:X \to \Xi_p$ obeys the inclusions
\begin{equation}
  \label{eq:F_X0pm_subset_Xi0pm}
  F(X^0) \subset \Xi^0 \quad\text{and}\quad F(X_p^\pm) \subset \Xi_p^\pm,
\end{equation}
where $X^0$ and $X_p^\pm$ are the complex analytic subsets of $X$ described in Definition \ref{maindefn:BB_decomposition_complex_analytic_space}.

By denoting the canonical injections and projections, respectively, as
\begin{align*}
  \iota^0:T_p^0X \to T_pX \quad&\text{and}\quad \iota^\pm: T_p^\pm X \to T_pX,
  \\
  \pi^0:\Xi_p \to \Xi_p^0 \quad&\text{and}\quad \pi^\pm:\Xi_p \to \Xi_p^\pm,
\end{align*}
we may define
\[
  F^0 = \pi^0\circ F\circ\iota^0:D\cap T_p^0X \to \Xi_p^0
  \quad\text{and}\quad
  F^\pm = \pi^\pm\circ F\circ\iota^\pm:D\cap T_p^\pm X \to \Xi_p^\pm.
\]
Therefore, writing
\[
  D^0 := D\cap T_p^0X \quad\text{and}\quad D^\pm := D\cap T_p^\pm X
\]
and, setting $r^0 := \dim \Xi_p^0$ and $r^\pm := \dim \Xi_p^\pm$ and choosing complex bases for $\Xi_p^0$ and $\Xi_p^\pm$ to write 
\[
  F^0 = (f_1^0,\ldots,f_{r^0}^0):D^0 \to \Xi_p^0
  \quad\text{and}\quad
  F^\pm = (f_1^\pm,\ldots,f_{r^\pm}^\pm): D^\pm \to \Xi_p^\pm,
\]  
we see that the complex analytic subsets $X^0$ and $X_p^\pm$ of $X$ are defined, respectively, by ideals
\begin{align*}
  \sI_{X^0} &:= (f_1^0,\ldots,f_{r^0}^0) \subset \sO_{D^0},
  \\
  \sI_{X_p^+} &:= (f_1^+,\ldots,f_{r^+}^+) \subset \sO_{D^+},
  \\
  \sI_{X_p^-} &:= (f_1^-,\ldots,f_{r^-}^-) \subset \sO_{D^-}.
\end{align*}
Here, we use the conventions that $D^0 := \{0\}$ if $n^0 = 0$, and $D^+ := \{0\}$ if $n^+ = 0$, and $D^- := \{0\}$ if $n^- = 0$, while $\sI_{X^0} := (0)$ if $r^0 = 0$, and $\sI_{X_p^+} := (0)$ if $r^+ = 0$, and $\sI_{X_p^-} := (0)$ if $r^- = 0$. Therefore,
\begin{align*}
  X^0 &= \supp\left(\sO_{D^0}/\sI_{X^0}\right),
  \\
  X_p^+ &= \supp\left(\sO_{D^+}/\sI_{X_p^+}\right),
  \\
  X_p^- &= \supp\left(\sO_{D^-}/\sI_{X_p^-}\right),
\end{align*}
with the convention that $X^0 = \{p\}$ if $n^0 = 0$, and $X^+ = \{p\}$ if $n^+ = 0$, and $X^- = \{p\}$ if $n^- = 0$. Moreover, the local rings are given by
\begin{align*}
  \sO_{X^0,p} &= \sO_{D^0,p}/\sI_{X^0,p},
  \\
  \sO_{X_p^+,p} &= \sO_{D^+,p}/\sI_{X_p^+,p},
  \\
  \sO_{X_p^-,p} &= \sO_{D^-,p}/\sI_{X_p^-,p},
\end{align*}
and so the lower bounds in \eqref{eq:Upper_and_lower_bounds_Krull_dimensions_X0_Xpm_at_p} follow from Lemma \ref{lem:Algebraic_dimension_pX_geq_expdim_pX}.

We may assume without loss of generality that $r$ is the minimal number of generators of $\sI$, so that
$\dim\sO_{X,p} = \expdim_p X$, while the lower bounds in \eqref{eq:Upper_and_lower_bounds_Krull_dimensions_X0_Xpm_at_p} yield
\[
  \dim\sO_{X^0,p} \geq \expdim_p X^0, \quad \dim\sO_{X^+,p} \geq \expdim_p X^+, \quad\text{and}\quad \dim\sO_{X^-,p} \geq \expdim_p X^-.
\]
The equality of \eqref{eq:Expdim_X_at_p_equals_sum_virtual_BB_nullity_index_coindex_index_X_at_p} of expected dimensions thus yields the following inequality for the sum of the Krull dimensions,
\begin{equation}
  \label{eq:Krulldim_X_at_p_equals_sum_Krulldim_X0_and_X+_and_X-_at_p}
  \dim\sO_{X,p} \leq \dim\sO_{X^0,p} + \dim\sO_{X^+,p} + \dim\sO_{X^-,p},
\end{equation}
and hence the inequality \eqref{eq:Krulldim_X_at_p_equals_sum_BBnullity_and_BBcoindex_and_BBindex_at_p}, by Definition \ref{maindefn:Stable_unstable_subspaces_BB_index_co-index_nullity} of the Bia{\l}ynicky--Birula nullity, coindex, and index.
\end{proof}

\begin{rmk}[On the proof of Theorem \ref{mainthm:Upper_and_lower_bounds_Krull_dimensions_X0_Xpm_at_p}]
\label{rmk:Proof_of_mainthm_upper_and_lower_bounds_Krull_dimensions_X0_Xpm_at_p}  
Theorem \ref{thm:Cartan_linearization_theorem_equivariant_noncompact} implies that $F:D \to \Xi_p$ is \emph{$\CC^*$-equivariant in the sense of germs} as in Definition \ref{defn:Cairns_Ghys_1-3}: for each $\lambda \in \CC^*$, there is an open neighborhood $U_\lambda \subset D$ of the origin such that $F(\lambda\cdot z) = \lambda\cdot F(z)$, for all $z \in U_\lambda$. It would have been tempting in the proof of Theorem \ref{mainthm:Upper_and_lower_bounds_Krull_dimensions_X0_Xpm_at_p} to try to appeal to Lemma \ref{lem:Equivariant_holomorphic_maps_preserve_BB_decompositions} to show that the morphism $F:F \to \Xi_p$ obeys the key inclusions \eqref{eq:F_X0pm_subset_Xi0pm}, but $F$ is only $\CC^*$-equivariant in the sense of germs and not (fully) $\CC^*$-equivariant, as assumed by Lemma \ref{lem:Equivariant_holomorphic_maps_preserve_BB_decompositions}.  
\end{rmk}

\section{Hamiltonian functions for circle actions}
\label{sec:Hamiltonian_functions_circle_actions}
The examples in this section serve to confirm our sign conventions and help illustrate the relation described in Section \ref{sec:BB_decomposition_complex_manifold} between Morse theory for the Hamiltonian function for an $S^1$ action on a K\"ahler manifold and the Bia{\l}ynicki--Birula decomposition for a $\CC^*$ action that induces the $S^1$ action. We continue the notation of Section \ref{sec:Blowups_analytic_manifolds_along_embedded_analytic_submanifolds}, but assume $\KK=\CC$ and that we are given a faithful, unitary representation for the circle as in \eqref{eq:Circle_matrix_representation}. The homomorphism $\rho:S^1\to\U(n)$ is real analytic and defines a real analytic action $S^1 \times \CC^n \to \CC^n$.

\begin{exmp}[Hamiltonian function for a unitary $S^1$ action on a Hermitian vector space]
\label{exmp:Hamiltonian_function_linear_circle_action_vector_space}  
Define a real analytic function
\begin{equation}
  \label{eq:Hamiltonian_function_linear_circle_action} 
  f:\CC^n \ni (z_1,\ldots,z_n) \mapsto \frac{1}{2}\sum_{j=1}^n l_j|z_j|^2 \in \RR
\end{equation}
and observe that, writing $h(z,w) = \langle z,w \rangle = z\cdot\bar w = \sum_{j=1}^n z_j\bar w_j$ for all $z,w \in \CC^n$, 
\[
  df(z)(\xi)
  =
  \sum_{j=1}^n \Real\,\langle l_jz_j,\xi_j\rangle
  =
  \sum_{j=1}^n \Real l_jz_j\bar\xi_j, \quad\text{for all } \xi \in \CC^n.
\]
The definition \eqref{eq:Vector_field_generator_circle_action} of the vector field $\Theta$ generating the isometric circle action $S^1\times \CC^n\to\CC^n$ in \eqref{eq:Circle_matrix_representation} gives
\[
  \Theta_z = i(l_1z_1,\ldots,l_nz_n), \quad\text{for all } z = (z_1,\ldots,z_n) \in \CC^n.
\]
The standard expression (see Kobayashi \cite[Equation (7.6.7), p. 251]{Kobayashi_differential_geometry_complex_vector_bundles}, noting that we omit the factor of $2$) for the K\"ahler form $\omega$ on $\CC^n$,
\[
  \omega(z)(\eta,\xi)
  =
  \Imag\langle \eta,\xi\rangle
  =
  -\Real\,\langle i\eta,\xi\rangle
  =
  \Real\,\langle \eta,i\xi\rangle,
  \quad\text{for all } z \in \CC^n \text{ and } \eta,\xi \in \CC^n,
\]
and the standard identification of $\CC^n \ni z \mapsto x \in \RR^{2n}$ via $z_j = x_j+iy_j$ for $j=1,\ldots,n$ and $z = (z_1,\ldots,z_n) \mapsto  (x_1,y_1,\ldots,x_n,y_n) = x$ yields
\begin{multline*}
  \omega(\Theta_z,\xi)
  =
  -\Real\,\langle i\Theta_z,\xi\rangle
  =
  \Real\,\langle (l_1z_1,\ldots,l_nz_n),(\xi_1,\ldots,\xi_n)\rangle
  \\
  =
  \sum_{j=1}^n \Real \langle l_jz_j,\xi_j\rangle
=
\sum_{j=1}^n \Real l_jz_j\bar\xi_j,
  \quad\text{for all } z \in \CC^n \text{ and } \xi \in \CC^n.
\end{multline*}
Hence, by comparing the preceding expressions for $df(\xi)$ and $\omega(\Theta_z,\xi)$, respectively, we obtain
\[
  df(z)(\xi) = \omega(\Theta_z,\xi), \quad\text{for all } z \in \CC^n \text{ and } \xi \in \CC^n,
\]
just as in Section \ref{sec:BB_decomposition_complex_manifold}, so $f$ is a Hamiltonian function for the isometric $S^1$ action on $\CC^n$ with its standard symplectic (K\"ahler) form. The Hessian bilinear form on $\CC^n$ is given by
\[
  \hess f(z)(\eta,\xi)
  =
  \sum_{j=1}^n \Real l_j\eta_j\bar\xi_j
  =
  \Real\,\langle L\eta, \xi\rangle, \quad\text{for all } \eta,\xi \in \CC^n,
\]
where $L = \diag(l_1,\ldots,l_n) \in \Mat(n,\CC)$. If $g$ is the standard Euclidean metric on $\RR^{2n} \cong \CC^n$ given by $g(x,y) = \Real\,\langle z,w \rangle$, then
\[
  \omega(z)(\eta,\xi) = \Real\,\langle\eta,i\xi\rangle = g(\eta,i\xi),
  \quad\text{for all } z \in \CC^n \text{ and } \eta,\xi \in \CC^n,
\]
and so the gradient and Hessian operator, respectively, are given by
\begin{align*}
  \grad_g f(z) &= (l_1z_1,\ldots,l_nz_n) = -i\Theta_z \in \CC^n \cong \RR^{2n},
  \\
  \Hess_g f(z) &= \diag(l_1,\ldots,l_n) \in \Mat(n,\CC) \cong \Mat(2n,\RR),
                 \quad\text{for all } z \in \CC^n,
\end{align*}
in agreement with the conventions of Frankel \cite{Frankel_1959}. However, note that the Hamiltonian function $f$ in \eqref{eq:Hamiltonian_function_linear_circle_action} is \emph{not} proper, one of the hypotheses of Yang \cite[Theorem 4.12, p. 92]{Yang_2008} which provides a Bia{\l}ynicki--Birula decomposition for a complete K\"ahler manifold with a $\CC^*$ action.
\qed
\end{exmp}

\begin{exmp}[Hamiltonian function for an $S^1$ action on complex projective space]
\label{exmp:Hamiltonian_function_circle_action_complex_projective_space}
This example is a sequel to Example \ref{exmp:Hamiltonian_function_linear_circle_action_vector_space}. We assume that $n\geq 2$ and regard $[z_1,\ldots,z_n]$ as homogeneous coordinates for $\PP^{n-1} = \PP(\CC^n)$. We now define
\begin{equation}
  \label{eq:Hamiltonian_function_linear_circle_action_complex_projective_space} 
  f:\CC^n\less\{0\} \ni (z_1,\ldots,z_n) \mapsto \frac{1}{2|z|^2}\sum_{j=1}^n l_j|z_j|^2 \in \RR,
\end{equation}
where $z := (z_1,\ldots,z_n)$ and $|z|^2 := \sum_{j=1}^n |z_j|^2$. Observe that $f$ is well-defined on $\PP(\CC^n)$ (see also Yang \cite[Example 4.13, p. 93]{Yang_2008}) and that for all $\xi = (\xi_1,\ldots,\xi_n) \in \CC^n = T_z(\CC^n\less\{0\})$ and $z \in \CC^n\less\{0\}$, we have
\[
  df(z)(\xi)
  =
  \frac{1}{|z|^4}\Real\left(|z|^2\sum_{j=1}^n l_j z_j\bar\xi_j
    - \sum_{j=1}^n l_j|z_j|^2\sum_{k=1}^n z_k\bar \xi_k \right),
\]
that is,
\[
  df(z)
  =
  \frac{1}{|z|^4}\Real\left(|z|^2\sum_{j=1}^n l_j z_jd\bar z_j
    - \sum_{j=1}^n l_j|z_j|^2\sum_{k=1}^n z_kd\bar z_k \right).
\]
Recall (see McDuff and Salamon \cite[Example 4.3.3, p. 168]{McDuffSalamonSympTop3}) that
\[
  T_{[z]}\PP(\CC^n) = \Hom(\CC z, z^\perp),
\]
where $z^\perp \subset \CC^n$ is the complex hyperplane defined by the orthogonal complement of the line $\CC z$. In particular, if $\xi \in \CC^n$, then it defines a homomorphism, $\CC z \ni \lambda z \mapsto \lambda\xi^\perp \in z^\perp \cong \CC^{n-1}$, where $\xi^\perp := \xi - \pi_z\xi = \xi - \langle \xi,z\rangle/|z|$, and thus an element of $T_{[z]}\PP(\CC^n)$. Clearly, the differential $df([z]):T_{[z]}\PP(\CC^n) \to \RR$ is well-defined for all $[z] \in \PP(\CC^n)$ by the preceding expression for $df(z):T_z(\CC^n\less\{0\}) \to \RR$.

On the other hand, the K\"ahler two-form on $\PP(\CC^n)$ is given by (see Ballmann \cite[Example 4.10 (5), p. 43]{Ballmann_lectures_Kaehler_manifolds}, Griffiths and Harris \cite[Chapter 0, Section 2, Example 3, p. 30]{GriffithsHarris}, Huybrechts \cite[Example 3.1.9 (i), p. 117]{Huybrechts_2005}, or McDuff and Salamon \cite[Example 4.3.3, p. 168 and Example 4.3.4, p. 169]{McDuffSalamonSympTop3})
and keeping in mind that we adhere to the convention in Kobayashi \cite[Equation (7.6.7), p. 251]{Kobayashi_differential_geometry_complex_vector_bundles} (but omit the factor $2$)
\[
  \omega(z)(\eta,\xi)
  =
  -\Real h(z)(i\eta,\xi)
  =
  -\frac{1}{|z|^4}\Real\left(|z|^2\sum_{j=1}^n (i\eta_j)\bar\xi_j
  - \sum_{j=1}^n \bar z_jz_k(i\eta_j)\bar\xi_k\right),
\]
for all $\xi,\eta \in \CC^n = T_z(\CC^n\less\{0\})$ and $z \in \CC^n\less\{0\}$, that is,
\[
  \omega(z)
  =
  -\Real h(z)(i\cdot,\cdot)
  =
  -\frac{1}{|z|^4} \Real \left(|z|^2\sum_{j=1}^n idz_jd\bar z_j
  - \sum_{j,k=1}^n i\bar z_jz_kdz_jd\bar z_k\right),
\]
where $h$ denotes the Hermitian metric on $\PP(\CC^n)$ induced by the standard Hermitian metric on $\CC^n$. The vector field $\Theta = (il_1z_1,\ldots,il_nz_n)$ generating the isometric circle action $S^1 \times \CC^n\less\{0\} \to \CC^n\less\{0\}$ in \eqref{eq:Circle_matrix_representation} defines a vector field, also denoted $\Theta$, which generates the induced isometric $S^1$ action on $\PP(\CC^n)$. 
We see that
\begin{align*}
  \omega(z)(\Theta_z,\cdot)
  &=
  -\frac{1}{|z|^4}\Real\left(|z|^2\sum_{j=1}^n (i^2l_jz_j)d\bar z_j
    - \sum_{j,k=1}^n i\bar z_jz_k (il_jz_j)d\bar z_k\right)
  \\
  &=
  -\frac{1}{|z|^4}\Real\left(-|z|^2\sum_{j=1}^n l_jz_jd\bar z_j
    + \sum_{j,k=1}^n l_j|z_j|^2 z_kd\bar z_k\right)
    \\
  &=
  \frac{1}{|z|^4}\Real\left(|z|^2\sum_{j=1}^n l_jz_jd\bar z_j
    - \sum_{j=1}^n l_j|z_j|^2 \sum_{k=1}^n z_kd\bar z_k\right).
\end{align*}
In particular, we obtain
\[
  df(z) = \omega(z)(\Theta_z,\cdot), \quad\text{for all } z \in \CC^n\less\{0\},
\]
and the preceding identity holds on $\PP(\CC^n)$ as well, namely
\[
  df[z] = \omega[z](\Theta_{[z]},\cdot), \quad\text{for all } [z] \in \PP(\CC^n).
\]
Rather than derive the preceding identity using homogeneous coordinates on $\PP(\CC^n)$, we could also compute using affine local coordinates, using the more familiar expressions for the K\"ahler metric in those coordinates in the cited references.
\qed
\end{exmp}  

\section{Morse inequalities and Poincar\'e polynomial of a complex algebraic variety}
\label{sec:Morse_inequalities_Poincare_polynomial_smooth_complex_algebraic_variety}
We first briefly review the Morse inequalities and then results due to Bia{\l}ynicki--Birula for the Poincar\'e polynomial of a smooth, complex projective algebraic variety.

\subsection{Morse inequalities}
\label{subsec:Morse_inequalities}
Let $X$ be a topological space, $\FF$ be a field, $k$ be a non-negative integer,
\[
  b_k(X,\FF) := \dim_\FF H_k(X,\FF),
\]
be the $k$-th Betti number, given by the dimension over $\FF$ of the singular homology group $H_k(X,\FF)$ (see, for example, Munkres \cite[Chapter 4]{Munkres2ndEdition}), and define the \emph{Poincar\'e series} (see Atiyah and Bott \cite[Section 1, p. 529]{Atiyah_Bott_1983} or Nicolaescu \cite[Section 2.3, p. 57 and p. 59]{Nicolaescu_morse_theory}) by
\[
  P_{X,\FF}(t) := \sum_{k\geq 0} b_k(X,\FF) t^k.
\]
We shall suppress $\FF$ from our notation; typically, we take $\FF=\RR$ or $\QQ$. The \emph{Euler characteristic} of $X$ is
\[
  e(X) := \sum_{k\geq 0} (-1)^k b_k(X),
\]
which gives the equality
\[
  e(X) = P_X(-1).
\]
Suppose now that $X$ is a closed, finite-dimensional, smooth manifold and that $f:X\to\RR$ is a Morse function with \emph{Morse series} (see Atiyah and Bott \cite[Section 1, p. 529]{Atiyah_Bott_1983} or Nicolaescu \cite[Section 1.1, p. 7]{Nicolaescu_morse_theory})
\[
  P_f(t) := \sum_{p \in \Crit f} t^{\lambda_p^-(f)} = \sum_{k\geq 0} \mu_k(f)t^k,
\]
where $\lambda_p^-(f)$ is the \emph{Morse index} of $f$ at a critical point $p\in X$ and $\mu_k(f)$ is the number of critical points (\emph{Morse number}) of $f$ of index $k$. If $X$ has (real) dimension $d$, then a local maximum $p_{\max}$ will have Morse index $\lambda_{p_{\max}}^-(f)=d$ and a local minimum $p_{\min}$ will have Morse index $\lambda_{p_{\min}}^-(f)=0$.

One has the well-known identity (see Atiyah and Bott \cite[Section 1, p. 529]{Atiyah_Bott_1983} or Nicolaescu \cite[Corollary 2.15, p. 60]{Nicolaescu_morse_theory})
\[
  e(X) = P_f(-1).
\]
The Morse function $f$ satisfies the \emph{Morse inequalities} (see Atiyah and Bott \cite[Section 1, p. 529]{Atiyah_Bott_1983}): there exists a polynomial $R(t)$ with non-negative coefficients such that
\[
  P_f(t) - P_X(t) = (1 + t)R(t),
\]  
and so the coefficients of $P_f(t)$ dominate those of $P_X(t)$. One calls $f$ an \emph{$\FF$-perfect Morse function} if one has (see Atiyah and Bott \cite[Section 1, p. 529]{Atiyah_Bott_1983} or Nicolaescu \cite[Section 2.3, p. 61]{Nicolaescu_morse_theory})
\[
  P_f(t) = P_X(t), \quad\text{for all } t \in \FF,
\]
that is, if the Morse series is equal to the Poincar\'e series for $\FF$; if equality holds for all fields $\FF$, then $f$ is called a \emph{perfect Morse function}.

Atiyah and Bott \cite[Section 1, p. 529]{Atiyah_Bott_1983} provide two criteria for establishing when a Morse function $f$ is perfect. First, one has the \emph{Morse lacunary principle}: if the set $\{\lambda_p^-(f)\}$ of all Morse indices of $f$ contains no consecutive integers, then $f$ is perfect. In particular, if it can be shown that $f$ has only even indices at its critical points, then the lacunary principle immediately implies that $f$ is perfect. This method was used by Bott and Samelson \cite{Bott_Samelson_1958, Bott_Samelson_1961} to show that the energy function on the space of loops of a Lie group is perfect, with similar applications due to
\begin{inparaenum}[\itshape a\upshape)]
\item Hitchin \cite[Section 7, p. 96]{Hitchin_1987} for a Hamiltonian function associated to the circle action on the moduli space of Higgs pairs over a Riemann surface,
\item Kirwan \cite[Theorem 5.4, p. 64]{Kirwan_cohomology_quotients_symplectic_algebraic_geometry} for the square of a moment map associated to the action of a compact Lie group on a compact symplectic manifold,
\item Nakajima \cite[Sections 5.1 and 5.2]{Nakajima_lectures_hilbert_schemes_points_surfaces} for a Hamiltonian function associated to the circle action induced by a $\CC^*$ action on the Hilbert scheme of $n$ points on smooth, complex, projective algebraic surface.
\end{inparaenum}
Second, one has the \emph{completion principle}, as described Atiyah and Bott \cite[Section 1, p. 529]{Atiyah_Bott_1983} or Nicolaescu \cite[Proposition 2.17, p. 61]{Nicolaescu_morse_theory}, but this is considerably more difficult to apply.

The preceding discussion extends to the case where $f$ is Morse--Bott in the sense of Definition \ref{maindefn:Morse-Bott_function}. Following Atiyah and Bott \cite[Section 1, p. 530]{Atiyah_Bott_1983} and restricting to the case where $X$ is orientable, we let $\Crit f$ denote the set of connected components of critical submanifolds of $f$. We write $\lambda_C^-(f)$ for the Morse--Bott index of a connected critical submanifold $C$. The \emph{Morse--Bott series} for $f$ is defined to be
\[
  P_f(t) := \sum_{C \in \Crit f} t^{\lambda_C^-(f)}P_C(t), \quad\text{for all } t \in \FF,
\]
where $P_C(t)$ is the Poincar\'e series of $C$. The Morse--Bott function $f$ satisfies the \emph{Morse--Bott inequalities} just as in the case of a Morse function. One calls $f$ an \emph{$\FF$-perfect Morse--Bott function} if one has (see Atiyah and Bott \cite[Section 1, p. 531]{Atiyah_Bott_1983})
\[
  P_f(t) = P_X(t), \quad\text{for all } t \in \FF,
\]
that is, if the Morse--Bott series is equal to the Poincar\'e polynomial for $\FF$; if equality holds for all fields $\FF$, then $f$ is called a \emph{perfect Morse--Bott function}.

\begin{rmk}[Critical submanifolds with Morse--Bott index zero]
\label{rmk:Critical_submanifolds_Morse-Bott_index_zero}  
Recall that $b_0(X)$ is equal to the number of connected components of $X$ and if $X$ is connected, then $P_X(0) = b_0(X)=1$. In particular, $P_C(0) = b_0(C) = 1$. If $f$ is a perfect Morse--Bott function, then
\[
  \sum_{C \in \Crit f} 1 = \sum_{C \in \Crit f} P_C(0) = P_f(0) = P_X(0) = b_0(X) = 1.
\]
Hence, there is exactly one connected critical submanifold $C_{\min} \subset X$ with Morse--Bott index equal to zero, namely where $f$ attains its absolute minimum value.
\end{rmk}

\subsection{Poincar\'e polynomial of a smooth complex algebraic variety}
\label{subsec:Poincare_polynomial_smooth_complex_algebraic_variety}
Suppose that $X$ is a smooth, complex projective variety with a holomorphic action $\CC^*\times X\to X$ and fixed point set $X^0$ with decomposition $\sqcup_{\alpha=1}^r X_\alpha^0$ into connected (thus, irreducible) components. One has that $X^0$ and thus each connected component $X_\alpha^0$ is a smooth subvariety of $X$. Let $N_{X_\alpha^0/X}$ be the normal bundle of $X_\alpha^0$ in $X$ and write
\[
  N_{X_\alpha^0/X} = \bigoplus_{n\in\ZZ}N_{X_\alpha^0/X}(n),
\]
where $N_{X_\alpha^0/X}(n) \subset N_{X_\alpha^0/X}$ is the subbundle where $\CC^*$ acts on the fibers with weight $n$. Define
\[
  d_\alpha(n) := \rank_\CC N_{X_\alpha^0/X}(n),
  \quad
  d_\alpha^+ := \sum_{n>0} d_\alpha(n),
  \quad
  d_\alpha^- := \sum_{n<0} d_\alpha(n).
\]
With the preceding understood, one has the

\begin{thm}[Poincar\'e polynomial of a smooth, complex projective variety with a holomorphic $\CC^*$ action]
\label{thm:Bialynicki-Birula_1}
(See Bia{\l}ynicki--Birula \cite[Theorem 1, p. 1097]{Bialynicki_1974}.)  
Let $X$ be a smooth, complex projective variety. If $\CC^*\times X\to X$ is a holomorphic action, then its Poincar\'e polynomial is given by
\begin{equation}
\label{eq:Bialynicki-Birula_1}  
  P_X(t) = \sum_{\alpha=1}^r P_{X_\alpha^0}(t) t^{2d_\alpha^+}  = \sum_{\alpha=1}^r P_{X_\alpha^0}(t) t^{2d_\alpha^-}.
\end{equation}
\end{thm}

Note that if $f$ is a Hamiltonian function corresponding to the induced action on $X$ by the circle $S^1 \subset \CC^*$ and the fixed point set $X^0$ has dimension zero, then \eqref{eq:Bialynicki-Birula_1} yields the equalities,
\[
  P_X(t) = \sum_{\alpha=1}^r t^{2d_\alpha^+}  = \sum_{\alpha=1}^r t^{2d_\alpha^-} = P_f(t)
\]
and so, in particular, $f$ is a perfect Morse function.

\section[Bia{\l}ynicki--Birula decomposition for Hilbert schemes of points]{Bia{\l}ynicki--Birula decomposition for the Hilbert scheme of points in the complex projective plane}
\label{subsec:Biayynicki-Birula_decomposition_Hilbert_scheme_points_complex_projective_surface}
Let $X$ be a nonsingular, complex projective surface and $n\geq 1$ be an integer. A well-known result due to Fogarty \cite[Theorem 2.4, p. 517]{Fogarty_1968} asserts that the Hilbert scheme $\Hilb^n(X)$ (also denoted $X^{[n]}$) of $n$ points in $X$ (see Nakajima \cite[Definition 1.2, p. 6]{Nakajima_lectures_hilbert_schemes_points_surfaces}) is a non-singular scheme of dimension $2n$ and is connected if $X$ is connected by \cite[Proposition 2.3, p. 517]{Fogarty_1968}. The Hilbert--Chow morphism $\pi:\Hilb^n(X) \to \Sym^n(X)$ is birational by Fogarty \cite[Corollary 2.7, p. 518]{Fogarty_1968} or Mumford, Fogarty, and Kirwan \cite[Section 5.4]{Mumford_Fogarty_Kirwan_geometric_invariant_theory} (see also Nakajima \cite[Section 1.3, Theorem 1.15, p. 12]{Nakajima_lectures_hilbert_schemes_points_surfaces})). Here, $\Sym^n(X)$ is the $n$-th symmetric product of $X$ (see Nakajima \cite[Section 1.1, p. 6]{Nakajima_lectures_hilbert_schemes_points_surfaces}).

From Theorem \ref{thm:Hilbert_polynomial} it follows that the degree of the Hilbert polynomial of a subscheme is equal to the dimension of the subscheme, so if $Z$ is a closed subscheme of a quasi-projective scheme $X$ over $\CC$ with Hilbert polynomial $n$, then $Z$ is $0$-dimensional and so is supported at distinct closed points. In particular, the support of $Z$ contains at most $n$ closed points.

Nakajima constructs a K\"ahler metric on the Hilbert scheme $\Hilb^n(\CC^2)$ of $n$ points in $\CC^2$ (see Nakajima \cite[Chapter 3]{Nakajima_lectures_hilbert_schemes_points_surfaces}) and applies Morse theory to compute the Poincar\'e polynomial corresponding to a perfect Morse function $f$ that is a Hamiltonian function for the action of the compact torus $S^1\times S^1$ on $\Hilb^n(\CC^2)$ (see Nakajima \cite[Section 5.2]{Nakajima_lectures_hilbert_schemes_points_surfaces}). As he notes, his approach is equivalent to that of Ellingsrud and Str{\o}mme \cite{Ellingsrud_Stromme_1987}, who apply the Bia{\l}ynicki--Birula decomposition to $\Hilb^n(\CC^2)$ and $\Hilb^n(\PP^2)$. By regarding $\CC^2$ as the complement of a line $\ell_\infty \subset \PP^2$, Ellingsrud and Str{\o}mme view $\Hilb^n(\CC^2)$ as an open subscheme of $\Hilb^n(\PP^2)$ corresponding to subschemes of $\PP^2$ with support in $\CC^2$ (see \cite[Section 1, p. 343]{Ellingsrud_Stromme_1987}). To preserve consistency with Nakajima and Yoshioka \cite[Section 2, p. 319]{NakajimaInstCountI}, we choose $\ell_\infty = \{[0,z_1,z_2] \in \PP^2\}$ rather than adopt the choice $L = \{[z_0,z_1,0] \in \PP^2\}$ of Ellingsrud and Str{\o}mme \cite[Section 2, p. 346]{Ellingsrud_Stromme_1987}.

In Lemma \ref{lem:Fixed-point_sets_circle_actions_complex_projective_space} and Remark \ref{rmk:Fixed-point_sets_C*_actions_complex_projective_space}, we determined the fixed point sets corresponding to $\CC^*$ actions on $\PP^n$. In \cite[Section 2, p. 346]{Ellingsrud_Stromme_1987}, Ellingsrud and Str{\o}mme consider the action of the maximal torus $T_3 \subset\SL(3,\CC)$ comprising diagonal matrices and observe that induced action on $\PP^2$ has three fixed points, $P_0 = [1,0,0]$, $P_1 = [0,1,0]$, and $P_2 = [0,0,1]$. The action of $T_3$ on $\PP^2$ induces an action of $T_3$ on $\Hilb^n(\CC^2)$. If $Z \subset \PP^2$ is a subscheme with $\length(\sO_Z)=n$ that is a fixed point of this action, then the topological support of $Z$ is contained in the set $\{P_0,P_1,P_2\}$. Hence, one may write $Z = Z_0\cup Z_1\cup Z_2$, where $Z_i$ is supported in $P_i$ and corresponds to a fixed point in $\Hilb^{n_i}(\CC^2)$, where $\length(\sO_{Z_i}) = n_i$. According to \cite[Lemma 2.1, p. 346]{Ellingsrud_Stromme_1987} and its proof, the action of $T_3$ on $\Hilb^n(\PP^2)$ has only finitely many fixed points. A point $Z \in \Hilb^n(\PP^2)$ is a fixed point if and only if the corresponding homogeneous ideal $\sI_Z \subset \CC[z_0,z_1,z_2]$ is invariant under $T_3$, which is the case if and only if $\sI_Z$ is generated by monomials.

Let $F_0 := \{P_0\}$, $F_1 := \ell_\infty - \{P_0\}$, and $F_2 := \PP^2 - \ell_\infty$. Then $F_i \cong \CC^i$ for $i=0,1,2$ and these subsets define a cellular decomposition of $\PP^2$. The one-parameter subgroups $\lambda:\CC^* \to T_3$ inducing this cellular decomposition are given by $\lambda(t) = \diag(t^{w_0},t^{w_1},t^{w_2})$, where $w_0<w_1<w_2$ and $w_0 + w_1 + w_2 = 0$. For a sufficiently general choice of one-parameter subgroup $\lambda$, the fixed point subsets of $\lambda$ and $T_3$ on $\Hilb^n(\PP^2)$ coincide by \cite[Remark 1.7, p. 346]{Ellingsrud_Stromme_1987}.

For any triple $(n_0,n_1,n_2)$ of non-negative integers with $n=n_0+n_1+n_2$, one defines $W(n_0,n_1,n_2)$ to be the (locally closed) subset of $\Hilb^n(\PP^2)$ corresponding to subschemes $Z$ with $\length(Z_i)=n_i$ for $i=0, 1, 2$. Thus,
\[
  \Hilb^n(\PP^2) = \bigcup_{n_0+n_1+n_2=n}W(n_0,n_1,n_2).
\]
Let $\lambda$ be any one-parameter subgroup of $T_3$ that respects the cellular decomposition $(F_0,F_1,F_2)$ of $\PP^2$. Ellingsrud and Str{\o}mme observe that $\lambda$ induces a cellular decomposition of $\Hilb^n(\PP^2)$ and $W(n_0,n_1,n_2)$ is a union of cells from this decomposition. Indeed, let $Z \in W(n_0,n_1,n_2)$ and write $Z = Z_0\cup Z_1\cup Z_2$. If $t\to 0$, then $\lambda(t)(Z_i)$
approaches a subscheme supported in $P_i$. Thus, $W(n_0,n_1,n_2)$ has a cellular decomposition and 
\[
  W(n_0,n_1,n_2) \cong W(n_0,0,0) \times W(0,n_1,0) \times W(0,0,n_2).
\]
The cells contained in $W(n_0,0,0)$ (respectively, $W(0,n_1,0)$ or $W(0,0,n_2)$) are exactly those corresponding to fixed points supported in $P_0$ (respectively, $P_1$ or $P_2$). Therefore, one may focus on $T_3$-invariant subschemes of $\PP^2$ with one fixed point of $T_3$. Any such subscheme is contained in a $T_3$-invariant affine plane. Hence, it is enough to examine ideals of $\CC[x,y]$ of finite colength and that are invariant under the action of a two-dimensional torus $T_2 \subset \GL(2,\CC)$. If $\sI\subset\CC[x,y]$ is such an ideal, then it is generated by monomials in $x$ and $y$. The number $b_j := \inf\{k: x^jy^k \in \sI\}$ exists for each integer $j \geq 0$ and $b_j = 0$ for $j \gg 0$. Let $m$ be the least integer such that $b_m=0$. One sees that $b_0 \geq b_1 \geq \cdots \geq b_{m-1}\geq b_m=0$ and
\[
  \sum_{j=0}^m b_j = \length(\CC[x,y]/\sI) = n,
\]
where $\length(\CC[x,y]/\sI)$ is equal to the dimension of $\CC[x,y]/\sI$ as a complex vector space. Furthermore, $\{y^{b_0}, xy^{b_1}, \ldots, x^jy^{b_j},\ldots, x^m\}$ is a (not necessarily minimal) set of generators for $\sI$. One obtains a bijection between $T_2$-invariant ideals of colength $n$ in $\CC[x,y]$ and partitions of $n$. Ellingsrud and Str{\o}mme exploit the preceding facts and the associated Bia{\l}ynicki--Birula decomposition for $\Hilb^n(\PP^2)$ to prove their \cite[Theorem 1.1, p. 344]{Ellingsrud_Stromme_1987}, which states that the odd Betti numbers of $\Hilb^n(\PP^2)$ and $\Hilb^n(\CC^2)$ are zero and give explicit formulae for the even Betti numbers and the Euler characteristics.

If $n=1$, then $\dim_\CC(\CC[x,y]/\sI)=1$, so $\CC[x,y]/\sI \cong \CC$ and $\sI$ is a maximal ideal and thus $(y,x)$ is set of generators since $\sI$ is fixed by $T_2$, so $b_0=1, b_1=0$, and $m=1$ in the notation above. If $n=2$, then $\dim_\CC(\CC[x,y]/\sI)=2$ and, for example, $\sI$ may have a set of generators of the form $(y,xy,x^2) = (y,x^2)$ for $b_0=1,b_1=1,b_2=0$ and $m=2$. Another possibility for a set of generators is $(y^2,x,x^2) = (y^2,x)$ for $b_0=2,b_1=0$ and $m=1$.

To motivate the construction of $\Hilb^n(\CC^2)$, observe that if one has $n$ distinct points in $\CC^2$, that is, a point in the smooth locus of $Z \in \Sym^n(\CC^2)$, then there is a unique ideal $\sI\subset\CC[x,y]$ comprising all polynomial functions that vanish on $Z$. The coordinate ring $\CC[x,y]/\sI$ is a complex vector space of dimension $n$. This $Z$ is a length-$n$ subscheme of the affine plane $\CC^2$ and $\Hilb^n(\CC^2)$ is a compactification of the smooth locus obtained by considering all ideals of length $n$. See Qin \cite[Section 1.4]{Qin_Hilbert_schemes_points_infinite_dimensional_Lie_algebras} for further results concerning Hilbert schemes of points in complex projective schemes.


\chapter[Bia{\l}ynicki--Birula decompositions for blowups of complex manifolds]{Bia{\l}ynicki--Birula decompositions for blowups of complex manifolds: Algebraic approach}
\label{chap:BB_decomposition_blowup_complex_manifold_algebraic} 
In this chapter, we use the functorial property of blowups to produce holomorphic $\CC^*$ actions and Bia{\l}ynicki--Birula decompositions on the blowups of complex analytic spaces endowed with holomorphic $\CC^*$ actions and Bia{\l}ynicki--Birula decompositions using an algebraic approach. (In Chapter \ref{chap:BB_decomposition_blowup_complex_manifold_analytic}, we consider the corresponding problem using an analytic approach and also obtain more general versions of the main results in this chapter.) We begin in Section \ref{sec:Blowup_complex_linear_subspace_C*-action} by examining the model case of a linear $\CC^*$ action on a finite-dimensional complex vector space and the induced holomorphic $\CC^*$ action on the blowup of that vector space along a $\CC^*$-invariant complex linear subspace. Section \ref{sec:Bialynicki-Birula decomposition_blowup_complex_linear_subspace_C*-action} provides a description of the Bia{\l}ynicki--Birula decomposition for the induced holomorphic $\CC^*$ action on the blowup of a complex vector space along a $\CC^*$-invariant complex linear subspace, where the complex vector space is endowed with a linear $\CC^*$ action. In Section \ref{sec:C*_action_blowup_complex_analytic_model_space}, we discuss the holomorphic $\CC^*$ action and its fixed points on the blowup of a complex analytic model space along a $\CC^*$-invariant linear subspace of the ambient complex vector space containing the the domain $D$ defining the complex analytic model space, where the complex vector space is endowed with a linear $\CC^*$ action.  Section \ref{sec:Bialynicki-Birula decomposition_blowup_compact_complex_Kaehler_manifold_C*-action} establishes a Bia{\l}ynicki--Birula decomposition for the induced holomorphic $\CC^*$ action on the blowup of a compact complex, K\"ahler manifold along a $\CC^*$-invariant complex submanifold, where the given compact complex, K\"ahler manifold is endowed with a holomorphic $\CC^*$ action. (In Chapter \ref{chap:BB_decomposition_blowup_complex_manifold_analytic}, we remove the restrictions that the complex manifold be compact complex or K\"ahler.)

\section{Linear circle and $\CC^*$ actions and blowups along invariant complex linear subspaces}
\label{sec:Blowup_complex_linear_subspace_C*-action}
%
With the aid of Examples \ref{exmp:Hamiltonian_function_linear_circle_action_vector_space} and \ref{exmp:Hamiltonian_function_circle_action_complex_projective_space}, we can identify a Hamiltonian function for the induced $S^1$ action on the blowup of $\CC^n$ along a coordinate subspace.

\begin{lem}[Hamiltonian function for an $S^1$ action on the blowup of a complex vector space along an invariant linear subspace]
\label{lem:Hamiltonian_function_circle_action_blowup_complex_vector_space_along_invariant_subspace}
Let $Z$ be a complex linear subspace of a complex vector space $X$ with a Hermitian inner product and $\rho:S^1 \to \U(X)$ be a unitary representation that restricts to a unitary representation $\rho:S^1 \to \U(Z)$, that is, $Z$ is an $S^1$-invariant linear subspace. If the blowup $\Bl_Z(X)$ has its standard K\"ahler metric as in Proposition \ref{prop:Voisin_3-24_vector_space}, then the $S^1$ action $\tilde\rho:S^1\times\Bl_Z(X) \to \Bl_Z(X)$ induced by the $S^1$ action $\rho:S^1\times X \to X$ is isometric and Hamiltonian with respect to the corresponding Riemannian metric and symplectic form on $\Bl_Z(X)$, respectively.
\end{lem}

\begin{proof}
The $S^1$ action on the product complex K\"ahler manifold $X\times \PP(Z^\perp)$ is clearly Hamiltonian with respect to the Hamiltonian function and corresponding symplectic form on the product complex K\"ahler manifold $X\times \PP(Z^\perp)$ given by Examples \ref{exmp:Hamiltonian_function_linear_circle_action_vector_space} and \ref{exmp:Hamiltonian_function_circle_action_complex_projective_space}, where $Z^\perp \subset X$ is the orthogonal complement of $Z$ with respect to the Hermitian inner product on $X$. Hence, the induced $S^1$ action on the embedded complex submanifold $\Bl_Z(X) \hookrightarrow X\times \PP(Z^\perp)$ is Hamiltonian with respect to the restriction to $\Bl_Z(X)$ of the Hamiltonian function and symplectic form on $X\times \PP(Z^\perp)$.

Similarly, the $S^1$ action on the product complex K\"ahler manifold $X\times \PP(Z^\perp)$ is clearly isometric with respect to the corresponding Riemannian metric on the product complex K\"ahler manifold $X\times \PP(Z^\perp)$ given by Examples \ref{exmp:Hamiltonian_function_linear_circle_action_vector_space} and \ref{exmp:Hamiltonian_function_circle_action_complex_projective_space}. Hence, the induced $S^1$ action on the embedded complex submanifold $\Bl_Z(X) \hookrightarrow X\times \PP(Z^\perp)$ is isometric with respect to the restriction to $\Bl_Z(X)$ of the Riemannian metric on $X\times \PP(Z^\perp)$.
\end{proof}

\begin{rmk}[Application of Cartan's magic formula to existence of Hamiltonian functions]
\label{rmk:Application_Cartan_magic_formula_existence_Hamiltonian_functions}  
The fact that the $S^1$ actions in Examples \ref{exmp:Hamiltonian_function_linear_circle_action_vector_space} and \ref{exmp:Hamiltonian_function_circle_action_complex_projective_space} and Lemma \ref{lem:Hamiltonian_function_circle_action_blowup_complex_vector_space_along_invariant_subspace} are Hamiltonian with respect to \emph{some} smooth function $f:X\to\RR$ follows from Cartan's magic formula (see Lee \cite[Theorem 14.35, p. 372]{Lee_john_smooth_manifolds}). For a smooth manifold $X$ with a smooth $S^1$ action and generator $\xi \in C^\infty(TX)$ and symplectic form $\omega$ that is $S^1$-invariant in the sense of \eqref{eq:Circle_invariant_covariant_2-tensor}, we have $\sL_\xi\omega = 0$ (by Lee \cite[Theorem 12.37, p. 324]{Lee_john_smooth_manifolds}) and thus $d(\iota_\xi\omega) = 0$ since $d\omega = 0$ and Lee \cite[Equation (14.32), p. 372]{Lee_john_smooth_manifolds} gives
\[
  \sL_\xi\omega = \iota_\xi(d\omega) + d(\iota_\xi\omega).
\]
If in addition $H^1(X;\RR) = 0$, then $\iota_\xi\omega$ must be an exact one-form and thus $\iota_\xi\omega = df$, for some $f \in C^\infty(X,\RR)$. (See Frankel \cite[Section 3, Lemma 1, p. 2]{Frankel_1959} for the preceding argument on existence of a Hamiltonian function.)
According to Griffiths and Harris \cite[Chapter 4, Section 6, p. 605]{GriffithsHarris}, when $Z$ is an embedded complex manifold of a complex manifold $X$ and $\pi:\Bl_Z(X)\to X$ is the blowup with exceptional divisor $E$, the cohomology of $\Bl_Z(X)$ is given by
\begin{equation}
  \label{eq:Cohomology_blowup_complex_manifold_along_embedded_submanifold}
  H^\bullet(\Bl_Z(X)) = \pi^*H^\bullet(X)\oplus H^\bullet(E)/\pi^*H^\bullet(Z).
\end{equation}
Applying this formula with coefficients $\RR$ to compute de Rham cohomology for the case $X = \CC^n$ with $n\geq 2$ and coordinate subspace $Z=\CC^m$ with $0\leq m<n$, we note that $H^k(\CC^n;\RR) = (0)$ and $H^k(\CC^m;\RR) = (0)$ for $k>0$ while $H^k(\PP(\CC^{n-m});\RR) = (0)$ for odd $k$ and $H^k(\PP(\CC^{n-m});\RR) = \RR$ for even $k\leq 2(n-m)$ (see Hatcher \cite[Theorem 3.12, p. 212]{Hatcher}).
Because $E = \CC^m\times \PP(\CC^{n-m})$, it follows from the preceding formula that $H^1(\Bl_{\CC^m}(\CC^n);\RR) = (0)$, so Cartan's magic formula predicts the existence of a Hamiltonian function for the induced $S^1$ action on $\Bl_{\CC^m}(\CC^n)$ that we obtained explicitly in the proof of Lemma \ref{lem:Hamiltonian_function_circle_action_blowup_complex_vector_space_along_invariant_subspace}.
\end{rmk}

\begin{rmk}[Existence of Hamiltonian functions for circle actions on compact, complex K\"ahler manifolds]
\label{rmk:Existence_Hamiltonian_function_circle_action_compact_complex_Kaehler_manifold}
While Frankel \cite[Section 3, Lemma 1, p. 2]{Frankel_1959} provides the chomological criterion $H^1(X;\RR) = 0$ for existence of a Hamiltonian function for a circle action on a compact, complex K\"ahler manifold, Kobayashi gives a simpler one in \cite[Chapter III, Corollary 4.6, p. 95]{Kobayashi_transformation_groups_differential_geometry} (see also Fujiki \cite[Lemma 1.5, p. 801]{Fujiki_1979}): Let $(X,\omega)$ be a compact, complex K\"ahler manifold. If $S^1\times X\to X$ is a smooth circle action such that $S^1$ has at least one fixed point, $S^1$ acts by isometries, and the almost complex structure $J$ is circle-invariant, then there is a function $f\in C^\infty(X;\RR)$ such that $df = \omega(\xi,\cdot)$, where the vector field $\xi \in C^\infty(TX)$ is the infinitesimal generator of the circle action.
\end{rmk}  

\begin{rmk}[Examples of circle actions that are not Hamiltonian]
\label{rmk:Circle_actions_exist_not_Hamiltonian}  
It is important to recall that the not every circle action is Hamiltonian. If $X$ is the two-torus $\TT^2 = \RR/\ZZ\times\RR/\ZZ = S^1\times S^1$ with its standard symplectic form $d\theta_1\wedge d\theta_2$, then the vector fields $\xi_1 := \partial/\partial\theta_1$ and $\xi_2 := \partial/\partial\theta_2$ generating the two obvious circle actions are not Hamiltonian (see Cannas da Silva \cite[Section 18.1, Example, p. 128]{Cannas_da_Silva_lectures_on_symplectic_geometry}).
Of course, $H^1(X;\RR) = \RR\oplus\RR$ in this example.
\end{rmk}  

\begin{rmk}[Cohomology of the blowup of a complex manifold along an embedded complex submanifold]
\label{rmk:Cohomology_blowup_complex_manifold_along_embedded_submanifold}
For the purpose of this monograph, we only need the fact that $H^1(\Bl_Z(X);\RR) = 0$ when $X=\CC^n$ and $Z=\CC^m$, a coordinate subspace with $0\leq m<n$ and $n\geq 2$. However, we give a proof (explained to us by Thomas Leness) of a far more general result, where $X$ is any complex manifold with $H^1(X;\RR) = 0$ and $Z \hookrightarrow X$ is an embedded complex submanifold. We first observe that if $H^1(X;\RR)=0$, then the relation \eqref{eq:Cohomology_blowup_complex_manifold_along_embedded_submanifold} yields
\[
  H^1(\Bl_Z(X);\RR) = H^1(E;\RR)/\pi^*H^1(Z;\RR).
\]
Recall from Proposition \ref{prop:Huybrechts_2-5-3} and Definition \ref{defn:Huybrechts_2-5-4} that $E = \PP(N_{Z/X})$, where $N_{Z/X}$ is the normal bundle of $Z$ in $X$. The Leray--Hirsch Theorem (see Bott and Tu \cite[Chapter I, Theorem 5.11, p. 50]{BT}) implies that
\[
  H^\bullet(\PP(N_{Z/X});\RR) = H^\bullet(Z;\RR)\otimes_\RR H^\bullet(\PP(\CC^{n-m});\RR),
\]
if $X$ and $Z$ have complex dimensions $n$ and $m$, respectively, with $0\leq m<n$ and $n\geq 2$, noting that the bundle $\PP(N_{Z/X})$ has fiber $\PP(\CC^{n-m})$. Consequently, using the cohomology of $\PP(\CC^{n-m})$ given in Remark \ref{rmk:Application_Cartan_magic_formula_existence_Hamiltonian_functions}, the Leray--Hirsch Theorem yields
\[
  H^1(\PP(N_{Z/X});\RR) = H^1(Z;\RR).
\]
Alternatively, from the expression for the cohomology of the projectivization of a vector bundle in Bott and Tu \cite[Chapter IV, Equation (20.7), p. 270]{BT},
\[
  H^\bullet(\PP(N_{Z/X});\RR)
  =
  \left.H^\bullet(Z;\RR)[x]\right/\left(x^m + c_1(N_{Z/X})x^{m-1} + \cdots + c_m(N_{Z/X})\right),
\]  
and one again sees that $H^1(\PP(N_{Z/X});\RR) = H^1(Z;\RR)$. In other words, the bundle map $\pi:E \to Z$ induces an isomorphism $\pi^*:H^1(Z;\RR) \to H^1(E;\RR)$ and thus $H^1(\Bl_Z(X);\RR) = (0)$.
\end{rmk}

\section[Bia{\l}ynicki--Birula decomposition for the blowup of a complex vector space]{Bia{\l}ynicki--Birula decomposition for a $\CC^*$ action on the blowup of a complex vector space  along an invariant complex linear subspace}
\label{sec:Bialynicki-Birula decomposition_blowup_complex_linear_subspace_C*-action}
We continue the notation of Section \ref{sec:Blowup_complex_linear_subspace_C*-action}. For clarity in the following discussion, we write $X = \CC^n$ and define its weight-sign decomposition with respect to the unitary representation $\rho:S^1\to\U(X)$ by
\begin{equation}
  \label{eq:X_weight-sign_decomposition}
  X_p = X^0\oplus X_p^+\oplus X_p^-,
\end{equation}
where, letting $e_1,\ldots,e_n$ denote the standard basis for $\CC^n$ with $e_1=(1,0,\ldots,0),\ldots,e_n=(0,\ldots,0,1)$, 
\begin{multline}
  \label{eq:W0_W+_W-_for_linear_circle_representation}
  X^0 := \bigoplus_{\{k:l_k=0\}}\CC e_k,
  \quad\text{and}\quad
  X_p^+ := \{p\}\times \bigoplus_{\{k:l_k>0\}}\CC e_k \cong \bigoplus_{\{k:l_k>0\}}\CC e_k
  \\
  \text{and}\quad
  X_p^- := \{p\}\times \bigoplus_{\{k:l_k<0\}}\CC e_k \cong \bigoplus_{\{k:l_k<0\}}\CC e_k,
  \quad\text{for all } p \in X^0,
\end{multline}
and write
\begin{equation}
  \label{eq:W+_W-_bundles_over_W0_for_linear_circle_representation}
  X^+ := \bigoplus_{\{k:l_k\geq 0\}}\CC e_k \quad\text{and}\quad X^- := \bigoplus_{\{k:l_k\leq 0\}}\CC e_k
\end{equation}
for the product bundles over $X^0$:
\begin{align*}
  \bigoplus_{\{k:l_k>0\}}\CC e_k \times \bigoplus_{\{k:l_k=0\}}\CC e_k \to \bigoplus_{\{k:l_k=0\}}\CC e_k,
  \\
  \bigoplus_{\{k:l_k<0\}}\CC e_k \times \bigoplus_{\{k:l_k=0\}}\CC e_k \to \bigoplus_{\{k:l_k=0\}}\CC e_k.
\end{align*}
Clearly, $X^0$ is the set of fixed points of the circle action on $X$ defined by $\rho$. Note that the intersections $X^0\cap Z$, and $X_p^+\cap Z$, and $X_p^-\cap Z$ are nonempty (they all include the point $p \in X^0$) and clean in the sense of Definition \ref{defn:Clean_intersection}, being given by coordinate subspaces of $X$. Furthermore, the unions $X^0\cup Z$, and $X_p^+\cup Z$, and $X_p^-\cup Z$ have normal crossings (as in Section \ref{sec:Normal_crossing_divisors}), being given by unions of coordinate subspaces of $X$. \emph{However}, one cannot apply Corollary \ref{cor:Smoothness_strict_transform_submanifold} without knowing that the intersections $X^0\cap Z$, and $X_p^+\cap Z$, and $X_p^-\cap Z$ are regular. For example, one or more of $X^0$, $X_p^+$, or $X_p^-$ could be contained in or equal to $Z$.

The subspaces $X^0,X_p^+,X_p^- \subset X$ may be alternatively defined by analogy with the Bia{\l}ynicki--Birula decomposition for algebraic varieties over a field $\KK$ and algebraic $\KK^*$ actions given by (see Theorem \ref{thm:Milne_13-47} or Bia{\l}ynicki--Birula \cite{Bialynicki_1974}) or compact, complex K\"ahler manifolds with holomorphic $\CC^*$ actions (see Theorem \ref{thm:Bialynicki-Birula_decomposition_compact_complex_Kaehler_manifold} or Carrell and Sommese  \cite{Carrell_Sommese_1983, Carrell_Sommese_1979cmh}). Recall from Remark \ref{rmk:Classification_complex_representations_C*} that the unitary representation $\rho:S^1\to\U(n)$ in \eqref{eq:Circle_matrix_representation} uniquely extends
to a complex linear representation $\rho_\CC:\CC^*\to\GL(n,\CC)$ in \eqref{eq:C*_matrix_representation} with the same integer weights. The subspaces $X^0,X_p^+,X_p^- \subset X$ are thus characterized by the expressions below:
\begin{subequations}
  \label{eq:W0_W+_W-_for_linear_C*_representation}
  \begin{align}
    \label{eq:W0_for_linear_C*_representation}
    X^0 &= \{z\in X:\rho_\CC(\lambda)z = z, \text{for all }\lambda\in\CC^*\},
    \\
    \label{eq:W+_for_linear_C*_representation}
    X^+ &:= \left\{z\in X: \lim_{\lambda\to 0}\rho_\CC(\lambda)z \in X^0\right\}
            \quad\text{and}\quad X_p^+ = (\pi^+)^{-1}(p),
    \\
    \label{eq:W-_for_linear_C*_representation}
    X^- &:= \left\{z\in X: \lim_{\lambda\to \infty}\rho_\CC(\lambda)z \in X^0\right\}
            \quad\text{and}\quad X_p^- = (\pi^-)^{-1}(p),
\end{align}
\end{subequations}
where the associated product bundle projections are given by
\begin{equation}
  \label{eq:bWpm_to_W0_projections}
  \pi^\pm:X^\pm \to X^0.
\end{equation}
In particular, $X^0$ is the set of fixed points of the $\CC^*$ action \eqref{eq:C*_matrix_representation} on $X$. As discussed in Examples \ref{exmp:Bialynicki-Birula_decomposition_C^2} and \ref{exmp:Bialynicki-Birula_decomposition_CP^1}, Theorem \ref{thm:Bialynicki-Birula_decomposition_compact_complex_Kaehler_manifold} predicts the local structure of $X$ near a fixed point $p$, but not its global structure since $X$ is noncompact.

The Bia{\l}ynicki--Birula decomposition for a complex vector space $X$ with a $\CC^*$ action induced by a representation $\rho_\CC:\CC^*\to\GL(X)$ of the form \eqref{eq:C*_matrix_representation} is immediate from the definition of $\rho_\CC$. Rather than give an abstract proof of the existence of the Bia{\l}ynicki--Birula decomposition for the blowup $\Bl_Z(X)$ of $X$ along an invariant subspace $Z$ by extending the proof due to Carrell and Sommese \cite{Carrell_Sommese_1978ms} of Theorem \ref{thm:Bialynicki-Birula_decomposition_compact_complex_Kaehler_manifold} or the proof due to Yang of \cite[Theorem 4.12, p. 92]{Yang_2008}, we shall instead give an explicit proof based on the Bia{\l}ynicki--Birula decompositions for $X$ and $\PP(Z^\perp)$ and the fact that the Bia{\l}ynicki--Birula decomposition for the blowup $\Bl_Z(X)$ is equal to the restriction of the Bia{\l}ynicki--Birula decomposition for $X\times \PP(Z^\perp)$ to the embedded complex submanifold $\Bl_Z(X)$.

We now apply Lemmas \ref{lem:BB_decomposition_complex_projective_space}, \ref{lem:BB_decomposition_C*_invariant_complex_submanifold_subset_existence}, \ref{lem:BB_decomposition_C*_invariant_complex_submanifold_subset_properties}, and \ref{lem:BB_decomposition_product_complex_manifolds} to deduce the existence of a Bia{\l}ynicki--Birula decomposition for the blowup $\Bl_Z(X)$ of a complex vector space $X$ with a linear $\CC^*$ action along a $\CC^*$-invariant subspace $Z \subset X$. (In the forthcoming Theorem \ref{thm:BB_decomposition_C*_action_blowup_complex_manifold_along_invariant_submanifold}, we establish a more general version of Theorem \ref{thm:BB_decomposition_blowup_vector_space_along_linear_subspace_explicit} by using analytical rather than algebraic methods.)

\begin{thm}[Bia{\l}ynicki--Birula decomposition for the blowup of a complex vector space along a $\CC^*$-invariant complex linear subspace]
\label{thm:BB_decomposition_blowup_vector_space_along_linear_subspace_explicit}  
Let $X$ be a complex vector space of finite dimension two or more, $Z\subsetneq X$ be a complex linear subspace, and $\rho_\CC:\CC^* \to \GL(X)$ be a representation that restricts to a representation $\rho_\CC:\CC^* \to \GL(Z)$, so $Z$ is a $\CC^*$-invariant linear subspace, and let $\tilde\rho_\CC:\CC^* \to \Aut(\Bl_Z(X))$ be the unique holomorphic $\CC^*$ action on the blowup $\Bl_Z(X)$ induced\footnote{See Lemma \ref{lem:C*_equivariance_blowup_vector_space_along_linear_subspace} for this construction.} such that the blowup morphism $\pi:\Bl_Z(X)\to X$ is $\CC^*$-equivariant. Then $\Bl_Z(X)$ admits a Bia{\l}ynicki--Birula decomposition in the sense of Definition \ref{maindefn:BB_decomposition_complex_manifold}.
\end{thm}

\begin{rmk}[Bia{\l}ynicki--Birula decompositions need not commute with strict transforms]
\label{rmk:Submanifolds_blowup_strict_transforms}  
The fixed-point submanifolds $\Bl_Z(X)_\alpha^0$ or stable or unstable fibers $\Bl_Z(X)_{\alpha,\tilde p}^\pm$ in $\Bl_Z(X)$ provided by Theorem
\ref{thm:BB_decomposition_blowup_vector_space_along_linear_subspace_explicit} need \emph{not} coincide with the strict transforms of $X^0$ or $X_p^\pm$, where $\tilde p \in \pi^{-1}(p)$ and $\pi:\Bl_Z(X)\to X$ is the blowup morphism. For example, if the blowup center $Z$ contains $X^0$ or $X_p^\pm$, then the strict transforms of $X^0$ or $X_p^\pm$ are empty since the sets $X^0\less Z$ or $X_p^\pm\less Z$ are empty.
\end{rmk}  

\begin{proof}[Proof of Theorem
\ref{thm:BB_decomposition_blowup_vector_space_along_linear_subspace_explicit}]
From the proof of Lemma \ref{lem:HilbertSpaceClosedUnderC*Action}, we may assume without loss of generality that $X = \CC^n$ and $Z = \CC^m$ embedded as a coordinate subspace as in \eqref{eq:Coordinate_subspace_blowup_center} with $n\geq 2$ and $0\leq m<n$ and that $\rho_\CC$ is given by \eqref{eq:C*_matrix_representation}, so
\[
  \rho_\CC(\lambda)(z_1,\ldots,z_n)
  =
  \left(\lambda^{l_1}z_1,\ldots,\lambda^{l_n}z_n\right),
  \quad\text{for all } \lambda \in \CC^* \text{ and } (z_1,\ldots,z_n) \in \CC^n,
\]
for integers $l_1,\ldots,l_n \in \ZZ$. From the definitions of $X^0$ and $X^\pm$ in \eqref{eq:W0_W+_W-_for_linear_C*_representation}, we have
\begin{align*}
  X^0 &= \left\{(z_1,\ldots,z_n) \in \CC^n: z_j = 0 \text{ for all } j\notin J^0\right\},
  \\
  X^+ &= \left\{(z_1,\ldots,z_n) \in \CC^n: z_j = 0 \text{ for all } j\notin J^+\right\} \supseteq X^0,
  \\
  X^- &= \left\{(z_1,\ldots,z_n) \in \CC^n: z_j = 0 \text{ for all } j\notin J^-\right\} \supseteq X^0,      
\end{align*}
with linear projections $\pi^\pm: X^\pm \to X^0$, where
\begin{align*}
  J^0 &= \left\{j \in \ZZ: 1\leq j\leq n \text{ and } l_j = 0\right\},
  \\
  J^+ &= \left\{j \in \ZZ: 1\leq j\leq n \text{ and } l_j \geq 0\right\} \supseteq J^0,
  \\
  J^- &= \left\{j \in \ZZ: 1\leq j\leq n \text{ and } l_j \leq 0\right\} \supseteq J^0,
\end{align*}
and a subspace $X^0$, $X^+$, or $X^-$ is equal to the zero subspace $(0)$ if one of the corresponding subsets $J^0$, $J^+$, or $J^-$ is empty. In particular, $X$ clearly admits a Bia{\l}ynicki--Birula decomposition in the sense of Definition \ref{maindefn:BB_decomposition_complex_manifold}.

Moreover, because $0\leq m<n$ by hypothesis and $Z^\perp \cong \CC^{n-m}$, then $\PP(Z^\perp)$ is equal to either a point if $n-m-1 = 0$ or a complex projective space $\PP(\CC^{n-m})$ of dimension $n-m-1 > 0$ otherwise. In the second case, $\PP(\CC^{n-m})$ admits a Bia{\l}ynicki--Birula decomposition in the sense of Definition \ref{maindefn:BB_decomposition_complex_manifold} by Lemma \ref{lem:BB_decomposition_complex_projective_space}. We conclude that $X\times\PP(Z^\perp)$ admits a Bia{\l}ynicki--Birula decomposition in the sense of Definition \ref{maindefn:BB_decomposition_complex_manifold} by Lemma \ref{lem:BB_decomposition_product_complex_manifolds}, since both $X$ and $\PP(Z^\perp)$ admit such decompositions. Lastly, $\Bl_Z(X)$ admits a Bia{\l}ynicki--Birula decomposition in the sense of Definition \ref{maindefn:BB_decomposition_complex_manifold} by Lemmas \ref{lem:BB_decomposition_C*_invariant_complex_submanifold_subset_existence} and \ref{lem:BB_decomposition_C*_invariant_complex_submanifold_subset_properties} since $\Bl_Z(X) \hookrightarrow X\times\PP(Z^\perp)$ is an embedded complex submanifold by Huybrechts \cite[Example 2.5.2, p. 99]{Huybrechts_2005} and is $\CC^*$-invariant with respect to the induced $\CC^*$ action by \eqref{eq:C*_action_blowup}.
\end{proof}

We obtain the following useful consequence from Lemma \ref{lem:Equivariant_holomorphic_maps_preserve_BB_decompositions} and the statement of Theorem \ref{thm:BB_decomposition_blowup_vector_space_along_linear_subspace_explicit}, but we refer the reader to the forthcoming Corollary \ref{cor:Correspondence_between_Bialynicki-Birula_decompositions_manifolds} for a generalization of Corollary \ref{cor:C*-equivariance_blowup_map_preservation_Bialynicki-Birula_decomposition} that is proved by analytical methods.

\begin{cor}[$\CC^*$-equivariance of the blowup map and preservation of the Bia{\l}ynicki--Birula decomposition]
\label{cor:C*-equivariance_blowup_map_preservation_Bialynicki-Birula_decomposition}  
Continue the hypotheses and notation of Theorem \ref{thm:BB_decomposition_blowup_vector_space_along_linear_subspace_explicit}. If $\pi:\Bl_Z(X) \to X$ is the blowup morphism, then $\pi$ is $\CC^*$-equivariant with respect to the $\CC^*$ action $\rho_\CC$ on $X$ in \eqref{eq:C*_matrix_representation} and the induced $\CC^*$ action $\tilde\rho_\CC$ on $\Bl_Z(X)$ in \eqref{eq:C*_action_blowup}, and the following inclusions hold:
\begin{subequations}
  \label{eq:piBlZWalphapm0_in_Wpm0}
  \begin{align}
    \label{eq:piBlZWalpha0_in_W0}
  \pi\left(\Bl_Z(X)^0\right) &\subseteq X^0,
    \\
    \label{eq:piBlZWalpha+_in_W+}
  \pi\left(\Bl_Z(X)^+\right) &\subseteq X^+,
    \\
    \label{eq:piBlZWalpha-_in_W-}
  \pi\left(\Bl_Z(X)^-\right) &\subseteq X^-.
  \end{align}
\end{subequations}
Moreover, for all points $\tilde p \in \Bl_Z(X)^0$ and their images $p = \pi(\tilde p) \in X^0$, the following hold:
\begin{subequations}
  \label{eq:piBlZWalphatildep_pm_in_Wp_pm}
  \begin{align}
    \label{eq:piBlZWalphatildep+_in_Wp+}
  \pi\left(\Bl_Z(X)_{\tilde p}^+\right) \subseteq X_p^+,
    \\
    \label{eq:piBlZWalphatildep-_in_Wp-}
  \pi\left(\Bl_Z(X)_{\tilde p}^-\right) \subseteq X_p^-.
  \end{align}
\end{subequations}
\end{cor}

\begin{proof}
Equivariance of the blowup map $\pi$ with respect to the $\CC^*$ actions follows from Lemma \ref{lem:C*_equivariance_blowup_vector_space_along_linear_subspace} (see also Section \ref{sec:C*_action_blowup_complex_analytic_model_space} for an explicit description in local coordinates and Theorem \ref{thm:Equivariance_property_blowups_complex_manifolds} for a general statement of equivariance). Lemma \ref{lem:Equivariant_holomorphic_maps_preserve_BB_decompositions} now yields the conclusions.
\end{proof}

\begin{rmk}[Applications of the Bia{\l}ynicki--Birula Decomposition to gauge-theoretic moduli spaces of Higgs bundles]
\label{rmk:Applications_BB_decomposition_moduli_spaces_Higgs_bundles}  
For applications of the Bia{\l}ynicki--Birula Decomposition to gauge-theoretic moduli spaces of Higgs bundles, we refer to Gothen and Z\'{u}\~{n}iga--Rojas \cite{Gothen_Zuniga-Rojas_2017} and Hausel \cite{Hausel_1998}.
\end{rmk}

The inclusion \eqref{eq:piBlZWalpha0_in_W0} in Corollary \ref{cor:C*-equivariance_blowup_map_preservation_Bialynicki-Birula_decomposition} is strengthened in the following lemma to an equality between the fixed-point sets for the $\CC^*$ actions on $X$ and its blowup $\Bl_Z(X)$ along $Z$. Because of its importance, we also provide an indirect proof of the lemma in Remark \ref{rmk:Blowdown_C*_fixed-point_set_blowupW_equals_C*_fixed-point_set_W}. Again, we refer the reader to \eqref{eq:piBlZX0_is_X0} in the forthcoming Corollary  \ref{cor:Correspondence_between_Bialynicki-Birula_decompositions_manifolds} for a version of Lemma \ref{lem:Blowdown_C*_fixed-point_set_blowupW_equals_C*_fixed-point_set_W} that is proved by analytical methods.

\begin{lem}[Relation between the fixed-point sets of a linear $\CC^*$ action on a complex vector space and the induced $\CC^*$ action on the blowup along a $\CC^*$-invariant, linear subspace]
\label{lem:Blowdown_C*_fixed-point_set_blowupW_equals_C*_fixed-point_set_W}
Continue the hypotheses of Lemma \ref{lem:C*_equivariance_blowup_vector_space_along_linear_subspace}, so that $X$ is a finite-dimensional, complex vector space, $Z \subsetneq X$ is a complex linear subspace, and $\rho_\CC:\CC^* \to \GL(X)$ is a homomorphism such that $Z$ is invariant under the induced action of $\CC^*$ on $X$. If $\pi:\Bl_Z(X)\to X$ is the blowup of $X$ along $Z$ constructed in Section \ref{sec:Blowups_analytic_manifolds_along_embedded_analytic_submanifolds}, then the fixed-point subsets $\Bl_Z(X)^0 = \Bl_Z(X)^{\CC^*}$ and $X^0 = X^{\CC^*}$ are related by
\begin{equation}
  \label{eq:pi_fixed-point_set_blowup_W_equals_fixed-point_set_W}
  \pi\left(\Bl_Z(X)^0\right) = X^0.
\end{equation}
If $\Bl_Z(X)_\alpha^0$, for $\alpha = 1,\ldots,r$, denote the finitely many connected components of the fixed-point set $\Bl_Z(X)^0$ for the induced $\CC^*$ action on $\Bl_Z(X)$, then
\begin{equation}
  \label{eq:pi_fixed-point_set_blowup_Walpha_equals_fixed-point_set_W}
  \pi\left(\Bl_Z(X)_\alpha^0\right) = X^0, \quad\text{for } \alpha = 1, \ldots, r.
\end{equation}
\end{lem}  

\begin{proof}
As we saw in the proof of Lemma \ref{lem:C*_equivariance_blowup_vector_space_along_linear_subspace}, the representation $\rho_\CC:\CC^* \to \GL(X)$ is necessarily of the form \eqref{eq:C*_matrix_representation} with respect to suitable basis for $X$. If $\tilde p \in \Bl_Z(X)$ is a fixed point of the action \eqref{eq:C*_action_blowup} of $\CC^*$ on $\Bl_Z(X)$, then Lemma \ref{lem:Equivariant_holomorphic_maps_preserve_BB_decompositions} and the fact that $\pi:\Bl_Z(X)\to X$ is $\CC^*$-equivariant imply that $p = \pi(\tilde p)$  is a fixed point of the action \eqref{eq:C*_matrix_representation} of $\CC^*$ on $X$, so
\[
  \pi\left(\Bl_Z(X)^0\right) \subseteq X^0.
\]
To prove the reverse inclusion, given a fixed point $p \in X$ of the action \eqref{eq:C*_matrix_representation} of $\CC^*$ on $X$, we shall seek a fixed point $\tilde p \in \pi^{-1}(p) \subset \Bl_Z(X)$. Since the submanifold $\Bl_Z(X) \subset \PP(Z^\perp)\times X$ is invariant under the $\CC^*$ action on $\PP(Z^\perp)\times X$, Lemma \ref{lem:Fixed_points_group_actions_on_subsets_and_invariant_subsets} yields the following relation between sets of fixed points of the $\CC^*$ actions:
\begin{equation}
  \label{eq:Fixed-point_set_BlZW_equals_BlZW_cap_fixed-point_set_PZperp_times_W}
  \Bl_Z(X)^{\CC^*}
  =
  \Bl_Z(X)\cap \left(\PP(Z^\perp)\times X\right)^{\CC^*}.
\end{equation}
We need to prove that the fiber $\pi^{-1}(p) \subset \Bl_Z(X)$ contains a fixed point $\tilde p$ of the $\CC^*$ action on $\Bl_Z(X)$. Recall from \eqref{eq:Blowup_linear_subspace_projection_Kn} that $\pi$ is the projection onto the second factor in the product $\PP(Z^\perp)\times X$. By the incidence relation \eqref{eq:Blowup_linear_subspace_incidence_variety}, any point $\tilde p \in \pi^{-1}(p)$ has the form $(\ell,p)$, where $\ell \in \PP(Z^\perp) \subset Z^\perp \subseteq X$ is a complex line such that $p \in \langle Z,\ell\rangle$, the linear span over $\KK$ in $X$ of $Z$ and $\ell$. Let $\{e_1,\ldots,e_m\}$ denote the standard basis for $Z=\CC^m$ and $\pi_Z$ and $\pi_{Z^\perp}$ denote the orthogonal projections from $X=Z\oplus Z^\perp$ onto $Z$ and $Z^\perp = \CC^{n-m}$, respectively. If $\ell = \CC v$, for a necessarily non-zero vector $v \in Z^\perp$, we then have
\[
  p = a_0v + \sum_{i=1}^m a_ie_i \in \langle Z,\ell\rangle \subset Z^\perp \oplus Z,
\]
where $a_0v = \pi_{Z^\perp}p \in Z^\perp$, and $a_i \in \CC$ for $i=0,1,\ldots,m$, and
\[
  \lambda\cdot p = a_0(\lambda\cdot v) + \sum_{i=1}^m a_i\lambda^{l_i}e_i \in Z^\perp \oplus Z,
  \quad\text{for all } \lambda \in \CC^*,
\]
where the weights $l_1,\ldots,l_m$ are given by the action \eqref{eq:C*_matrix_representation} of $\CC^*$ on $X$. But $p$ is a fixed point of the $\CC^*$ action on $X$, so $\lambda\cdot p = p$ for all $\lambda \in \CC^*$. If $a_0$ is non-zero, then we must have $\lambda\cdot v = v$ for all $\lambda \in \CC^*$ and so $v$ is a fixed point of the $\CC^*$ action on $Z^\perp$. Consequently, to choose a line $\ell \in \PP(Z^\perp)$ such that $\tilde p = (\ell,p) \in \Bl_Z(X)$ \emph{and} $\tilde p$ is a fixed point of the $\CC^*$ action on $\PP(Z^\perp)\times X$ when $a_0$ is non-zero, we may choose $\ell = \CC v$, where $v$ is the induced fixed point of the $\CC^*$ action on $Z^\perp$. If $a_0$ is zero, then $p$ imposes no constraint on $\ell = \CC v$ through the incidence relation \eqref{eq:Blowup_linear_subspace_incidence_variety} and we may choose any $v\in Z^\perp$ such that $\ell \in \PP(Z^\perp)$ is a fixed point of the $\CC^*$ action on $\PP(Z^\perp)$. Lemma \ref{lem:Fixed-point_sets_circle_actions_complex_projective_space} and Remark \ref{rmk:Fixed-point_sets_C*_actions_complex_projective_space} classify the fixed-point subsets of the $\CC^*$ actions on $\PP(Z^\perp)$ induced by linear $\CC^*$ actions on $Z^\perp$. We have thus shown that the reverse inclusion,
\[
  X^0 \subseteq \pi\left(\Bl_Z(X)^0\right),
\]
also holds. The refinement \eqref{eq:pi_fixed-point_set_blowup_Walpha_equals_fixed-point_set_W} follows immediately from the proof of the identity \eqref{eq:pi_fixed-point_set_blowup_W_equals_fixed-point_set_W}. This completes the proof of the lemma.
\end{proof}

\begin{rmk}[Indirect proof of the identity \eqref{eq:pi_fixed-point_set_blowup_W_equals_fixed-point_set_W} from properties of blowup morphisms]
\label{rmk:Blowdown_C*_fixed-point_set_blowupW_equals_C*_fixed-point_set_W}   
The identity \eqref{eq:pi_fixed-point_set_blowup_W_equals_fixed-point_set_W} may also be inferred indirectly from properties of the blowup $\pi:\Bl_Z(X)\to X$ provided by Proposition \ref{prop:Huybrechts_2-5-3}. If $p \in X\less Z$ is a fixed point of the $\CC^*$ action, then the fact that the restriction of the blowup morphism, $\pi:\Bl_Z(X)\less E\to X\less Z$, to the complement of the exceptional divisor $E = \pi^{-1}(Z)$ is a $\CC^*$-equivariant, biholomorphic map implies that there is a unique fixed point $\tilde p \in \Bl_Z(X)$ such that $\pi(\tilde p) = p$. If $p \in Z$, then the fiber $\pi^{-1}(p) \subset E$ is identified with a copy of $\PP(T_pZ^\perp)$ and we may choose any fixed point $\tilde p \in \PP(T_pZ^\perp)$ of the $\CC^*$ action provided by Lemma \ref{lem:Fixed-point_sets_circle_actions_complex_projective_space} and Remark \ref{rmk:Fixed-point_sets_C*_actions_complex_projective_space}, where $T_pX = T_pZ \oplus T_pZ^\perp$ as an orthogonal direct sum.
\end{rmk}

The forthcoming Corollary \ref{cor:Correspondence_between_Bialynicki-Birula_decompositions_manifolds} provides a generalization of the following lemma from the case where $X$ is a complex vector space and $Z$ is a linear subspace to the case where $X$ is a complex manifold and $Z$ is an embedded complex submanifold. The following lemma sharpens Corollary \ref{cor:C*-equivariance_blowup_map_preservation_Bialynicki-Birula_decomposition} by establishing that the inclusions \eqref{eq:piBlZWalphatildep_pm_in_Wp_pm} in that result are actually equalities.

\begin{lem}[Preservation of the Bia{\l}ynicki--Birula decomposition in the blowup of a complex vector space along an invariant linear subspace]
\label{lem:Preservation_Bialynicki-Birula_decomposition_blowup_vector_space}  
Continue the hypotheses of Lemma \ref{lem:Blowdown_C*_fixed-point_set_blowupW_equals_C*_fixed-point_set_W}. Then
\begin{equation}
  \label{eq:Preservation_Bialynicki-Birula_decomposition_blowup_vector_space}
  \pi\left(\Bl_Z(X)^\pm\right) = X^\pm
\end{equation}
and for any $p \in X^0$ and $\tilde p \in \Bl_Z(X)^0$ such that $\pi(\tilde p) = p$,
\begin{equation}
  \label{eq:Preservation_Bialynicki-Birula_decomposition_blowup_vector_space_fibers}
  \pi\left(\Bl_Z(X)_{\tilde p}^\pm\right) = X_p^\pm.
\end{equation}
Moreover,
\begin{equation}
  \label{eq:Preservation_Bialynicki-Birula_decomposition_blowup_vector_space_fibers_alpha}
  \pi\left(\Bl_Z(X)_{\alpha,\tilde p}^\pm\right) = X_p^\pm, \quad\text{for } \alpha = 1,\ldots,r.
\end{equation}
\end{lem}

\begin{proof}
We first claim that  
\begin{equation}
  \label{eq:BlZWpm_equals_BlZW_cap_PZperp_times_Wpm}
  \Bl_Z(X)^\pm
  =
  \Bl_Z(X)\cap \left(\PP(Z^\perp)\times X\right)^\pm.
\end{equation}
We focus on the `$+$' version of the identity \eqref{eq:BlZWpm_equals_BlZW_cap_PZperp_times_Wpm}, as the proof of the `$-$' version is identical. If $\tilde q \in \Bl_Z(X)^+$, then there exists $\tilde p \in \Bl_Z(X)^0$ such that $\lambda\cdot \tilde q \to \tilde p$ as $\lambda \to 0$ and thus $\tilde q \in \Bl_Z(X)_{\tilde p}^+$. Since $\tilde p \in \Bl_Z(X)^0$, the identity \eqref{eq:Fixed-point_set_BlZW_equals_BlZW_cap_fixed-point_set_PZperp_times_W} yields
\[
  \tilde p \in \left(\PP(Z^\perp)\times X\right)^0,
\]
and so $\tilde q \in \left(\PP(Z^\perp)\times X\right)_{\tilde p}^+$. By assumption, we also have $\tilde q\in\Bl_Z(X)$, so this gives the inclusion
\[
  \Bl_Z(X)^+
  \subseteq
  \Bl_Z(X)\cap \left(\PP(Z^\perp)\times X\right)^+.
\]
To prove the reverse inclusion, observe that if $\tilde q \in \Bl_Z(X)\cap \left(\PP(Z^\perp)\times X\right)^+$, then there exists $\tilde p \in \left(\PP(Z^\perp)\times X\right)^0$ such that $\lambda\cdot \tilde q \to \tilde p \in \left(\PP(Z^\perp)\times X\right)^0$ as $\lambda \to 0$. But $\Bl_Z(X)$ is invariant under the action of $\CC^*$ and thus $\tilde p \in \Bl_Z(X)$, since $\Bl_Z(X)$ is a topologically closed subspace of $\PP(Z^\perp)\times X$. According to \eqref{eq:Fixed-point_set_BlZW_equals_BlZW_cap_fixed-point_set_PZperp_times_W}, we have
\[
  \tilde p \in \Bl_Z(X)^0 = \Bl_Z(X)\cap \left(\PP(Z^\perp)\times X\right)^0,
\]
and thus $\tilde q\in \Bl_Z(X)^+$. This yields the reverse inclusion
\[
  \Bl_Z(X)\cap \left(\PP(Z^\perp)\times X\right)^+
  \subseteq
  \Bl_Z(X)^+
\]
and proves the claim \eqref{eq:BlZWpm_equals_BlZW_cap_PZperp_times_Wpm}.

We now consider the assertion \eqref{eq:Preservation_Bialynicki-Birula_decomposition_blowup_vector_space} and again focus on the `$+$' version, as the proof of the `$-$' version is identical and can thus be omitted. The inclusion
\begin{equation}
\label{eq:Bialynicki-Birula_decomposition_blowup_vector_space_+_inclusion_W+}
  \pi\left(\Bl_Z(X)^+\right) \subseteq X^+
\end{equation}
is an immediate consequence of the inclusion \eqref{eq:piBlZWalpha+_in_W+} provided by Corollary \ref{cor:C*-equivariance_blowup_map_preservation_Bialynicki-Birula_decomposition}, so we focus on the reverse inclusion
\begin{equation}
\label{eq:W+_inclusion_Bialynicki-Birula_decomposition_blowup_vector_space_+}
  X^+ \subseteq \pi\left(\Bl_Z(X)^+\right).
\end{equation}
If $q\in X^+$, then there exists $p\in X^0$ such that $\lambda\cdot q \to p$ as $\lambda \to 0$ and thus $q \in X_p^+$. If $\tilde q \in \pi^{-1}(q)\cap \Bl_Z(X)$, then $q \in \langle \ell,Z \rangle$ and $\ell = \CC v$, for a necessarily non-zero vector $v \in Z^\perp$. To prove the inclusion \eqref{eq:W+_inclusion_Bialynicki-Birula_decomposition_blowup_vector_space_+}, we need to choose $v$ such that $q \in \langle v,Z \rangle$ and $\tilde q = (\ell,q) \in \Bl_Z(X)^+$. We may write
\[
  q = \pi_{Z^\perp}q + \pi_Z q.
\]
If $\pi_{Z^\perp}q = 0$, then choose any $v \in Z^\perp$ such that $\ell = \CC v\in \PP(Z^\perp)^+$ and observe that
\[
  \tilde q = (\ell,q) \in \left(\PP(Z^\perp)\times X\right)^+,
\]
where the inclusion follows from $(\ell,q) \in \PP(Z^\perp)^+\times X^+$ and the easily proved facts that
\begin{equation}
  \label{eq:Fixed-points_PZperp_times_W_equals_fixed-points_PZperp_times_fixed-points_W}
  \PP(Z^\perp)^0\times X^0
  =
  \left(\PP(Z^\perp)\times X\right)^0
\end{equation}
and thus
\begin{equation}
  \label{eq:(PZperp_times_W)pm_equals_PZperppm_times_Wpm}
  \PP(Z^\perp)^\pm\times X^\pm
  =
  \left(\PP(Z^\perp)\times X\right)^\pm.
\end{equation}
Since we also have $\tilde q = (\ell,q) \in \Bl_Z(X)$, then $\tilde q \in \Bl_Z(X)^+$ by \eqref{eq:BlZWpm_equals_BlZW_cap_PZperp_times_Wpm}. If $\pi_{Z^\perp}q \neq 0$, then choose $v = \pi_{Z^\perp}q$ and observe that $q \in \langle \ell,Z\rangle$ for $\ell = \CC v$, so $\tilde q \in \pi^{-1}(q) \cap \Bl_Z(X)$. Moreover
\[
  \lambda\cdot q \to p \quad\text{and}\quad
  \lambda\cdot\pi_{Z^\perp}q = \pi_{Z^\perp}(\lambda\cdot q) \to \pi_{Z^\perp}p \quad\text{and}\quad
  \lambda\cdot\pi_Zq = \pi_Z(\lambda\cdot q) \to \pi_Zp.
\]
Because $p \in X$ is a fixed point of the $\CC^*$ action on $X$, then $\pi_Zp \in Z$ and $\pi_{Z^\perp}p \in Z^\perp$ are fixed points of the $\CC^*$ action on $Z$ and $Z^\perp$, respectively. Because the $\CC^*$ action on $X$ induces a $\CC^*$ action on $\PP(Z^\perp)$, we must have $\pi_{Z^\perp}p \neq 0$ and $\ell_0 = \CC\pi_{Z^\perp}p \in \PP(Z^\perp)$ is a fixed point of the $\CC^*$ action on $\PP(Z^\perp)$. Thus
\[
  \lambda\cdot\tilde q = \lambda\cdot(\ell,q) =  (\lambda\cdot\ell,\lambda\cdot q)
  \to (\ell_0,p) \in \PP(Z^\perp)^0\times X^0 = \left(\PP(Z^\perp)\times X\right)^0
\]
and so $\tilde q \in \left(\PP(Z^\perp)\times X\right)^+$. But $(\ell_0,p) = \pi^{-1}(p) \in \Bl_Z(X)$, so $\tilde p = (\ell_0,p) \in \Bl_Z(X)^0$ and because $\tilde q \in \Bl_Z(X)$ as well, we have
\[
  \tilde q \in \Bl_Z(X)\cap \left(\PP(Z^\perp)\times X\right)^+
\]  
and thus $\tilde q \in \Bl_Z(X)^+$ by the identity \eqref{eq:BlZWpm_equals_BlZW_cap_PZperp_times_Wpm}. This proves the inclusion \eqref{eq:W+_inclusion_Bialynicki-Birula_decomposition_blowup_vector_space_+} and, together with the inclusion \eqref{eq:Bialynicki-Birula_decomposition_blowup_vector_space_+_inclusion_W+}, proves the `$+$' version of the identity \eqref{eq:Preservation_Bialynicki-Birula_decomposition_blowup_vector_space}.

Furthermore, by \eqref{eq:piBlZWalphatildep+_in_Wp+} in Corollary \ref{cor:C*-equivariance_blowup_map_preservation_Bialynicki-Birula_decomposition} we have the following inclusion of fibers,
\[
  \pi(\Bl_Z(X)_{\tilde p}^+) \subseteq X_p^+,
\]
for any $p \in X^0$ and $\tilde p \in \pi^{-1}(p)\cap \Bl_Z(X)^0$. The reverse inclusion,
\[
  X_p^+ \subseteq \pi(\Bl_Z(X)_{\tilde p}^+),
\]
for any $p \in X^0$ and $\tilde p \in \pi^{-1}(p)\cap \Bl_Z(X)^0$, is an immediate consequence of our proof of the inclusion \eqref{eq:W+_inclusion_Bialynicki-Birula_decomposition_blowup_vector_space_+}. This proves the `$+$' version of the identity \eqref{eq:Preservation_Bialynicki-Birula_decomposition_blowup_vector_space_fibers}, while the proof of the `$-$' version of the identity is identical and can be omitted. This completes the proof of Lemma \ref{lem:Preservation_Bialynicki-Birula_decomposition_blowup_vector_space}.
\end{proof}

\section[$\CC^*$ action and fixed points on the blowup of a complex analytic space]{$\CC^*$ action and fixed points on the blowup of a complex analytic space along a $\CC^*$-invariant subspace}
\label{sec:C*_action_blowup_complex_analytic_model_space}
The forthcoming Lemma \ref{lem:C*_action_blowup_complex_analytic_model_space} is a corollary of Theorem \ref{thm:Strict_transform_complex_analytic_model_space} and explicit descriptions in local coordinates provided below of the $\CC^*$ action on $\Bl_Z(X)$ induced by the the $\CC^*$ action on $X$ and the blowup $\pi:\Bl_Z(X)\to X$ of $X$ along a $\CC^*$-invariant linear subspace $Z \subset X$. The existence of fixed points asserted by Lemma \ref{lem:C*_action_blowup_complex_analytic_model_space} is proved by algebraic methods, whereas the corresponding result (for the blowup of a complex manifold rather than a complex analytic model space) is obtained by analytical methods in the proof of Theorem \ref{thm:BB_decomposition_C*_action_blowup_complex_manifold_along_invariant_submanifold}.

Suppose that $X=\CC^n$ and that $Z=\CC^m$ is the standard coordinate subspace in \eqref{eq:Coordinate_subspace_blowup_center} with $\KK=\CC$. Given a linear representation $\rho_\CC:\CC^*\to\GL(n,\CC)$ as in \eqref{eq:C*_matrix_representation}, we know from Lemma \ref{lem:C*_equivariance_blowup_vector_space_along_linear_subspace} that there exists a unique holomorphic action,
\begin{equation}
  \label{eq:C*_action_blowup}
  \tilde\rho_\CC:\CC^*\times \Bl_{\CC^m}(\CC^n) \to \Bl_{\CC^m}(\CC^n),
\end{equation}
such that the blowup map $\pi:\Bl_{\CC^m}(\CC^n) \to \CC^n$ in \eqref{eq:Blowup_linear_subspace_projection_Kn} is $\CC^*$-equivariant with respect to the actions $\tilde\rho_\CC$ on $\Bl_{\CC^m}(\CC^n)$ and $\rho_\CC$ on $\CC^n$. In the proof of Lemma \ref{lem:C*_equivariance_blowup_vector_space_along_linear_subspace}, we gave a global construction of $\tilde\rho_\CC$, but it will be useful to express $\tilde\rho_\CC$ in terms of $\CC^*$-invariant, local coordinate charts for $\Bl_{\CC^m}(\CC^n)$. To accomplish this, we will check that each coordinate domain $U_j \subset \Bl_{\CC^m}(\CC^n)$ in \eqref{eq:Blowup_linear_subspace_coordinate_patch} is
$\CC^*$-invariant with respect to the action $\tilde\rho_\CC$ on $\Bl_{\CC^m}(\CC^n)$ and define linear representations $\rho_j:\CC^*\to\GL(n,\CC)$, for $j=1,\ldots,n$, such that the following diagram commutes:
\begin{equation}
\label{eq:Blowup_local_coordinate_chart_circle_equivariance}
\begin{tikzcd}
  \Bl_{\CC^m}(\CC^n) \supset U_j \arrow[r, "\tilde\rho_\CC(\lambda)"] \arrow[d, "\varphi_j"']
  &U_j \subset \Bl_{\CC^m}(\CC^n) \arrow[d, "\varphi_j"]
    \\
    \CC^n \supset \varphi_j(U_j) \arrow[r, "\rho_j(\lambda)"] &\varphi_j(U_j) \subset \CC^n
  \end{tikzcd}
\end{equation}
We recall from the proof of Lemma \ref{lem:C*_equivariance_blowup_vector_space_along_linear_subspace} that the linear $\CC^*$ action $\rho_\CC:\CC^*\times\CC^n\to\CC^n$ restricts to a linear $\CC^*$ action on $\CC^{n-m} \hookrightarrow \CC^n$ and thus defines a $\CC^*$ action on $\PP(\CC^{n-m})$. The equations defining the blowup $\Bl_{\CC^m}(\CC^n) \subset \PP(\CC^{n-m})\times\CC^n$ in \eqref{eq:Blowup_linear_subspace_equations} are preserved by the induced $\CC^*$ action on $\PP(\CC^{n-m})\times \CC^n$ and so that action induces the global holomorphic $\CC^*$ action \eqref{eq:C*_action_blowup} on $\Bl_{\CC^m}(\CC^n)$. From this construction, we see that each coordinate domain $U_j \subset \Bl_{\CC^m}(\CC^n)$ in \eqref{eq:Blowup_linear_subspace_coordinate_patch} is
$\CC^*$-invariant with respect to the action $\tilde\rho_\CC$ on $\Bl_{\CC^m}(\CC^n)$, as claimed. Keeping in mind the definition of the charts $\varphi_j = (w(j)_1,\ldots,w(j)_n)$ defined coordinatewise by \eqref{eq:Blowup_linear_subspace_coordinates} for $j=1,\ldots,n$, we choose linear representations for $\CC^*$,
\begin{equation}
  \label{eq:C*_action_blowup_coordinate_patch_j}
  \rho_j:\CC^*\to\GL(n,\CC),
  \quad\text{where } \rho_j(\lambda)w_k
  :=
  \begin{cases}
    \lambda^{l_k}w_k, &\text{for } k = 1, \ldots,m,
    \\
    \lambda^{l_k-l_j}w_k, &\text{for } k = m+1, \ldots,\hat\jmath,\ldots,n,
    \\
    \lambda^{l_j}w_j, &\text{for } k = j,
  \end{cases}
\end{equation}
where $\hat\jmath$ indicates that the index $j$ is omitted from the index set $\{m+1,\ldots,n\}$, for $j=1,\ldots,n$. The coordinate charts $\varphi_j$ are thus clearly $\CC^*$-equivariant in the sense of the diagram \eqref{eq:Blowup_local_coordinate_chart_circle_equivariance}. Hence, each linear $\CC^*$ action \eqref{eq:C*_action_blowup_coordinate_patch_j} on the image $\varphi_j(U_j)\subset\CC^n$ of the $\CC^*$-invariant coordinate domain $U_j\subset\Bl_{\CC^m}(\CC^n)$ glues together, for $j=1,\ldots,n$, to give the global holomorphic $\CC^*$ action $\tilde\rho_\CC$ on $\Bl_{\CC^m}(\CC^n)$. We can now proceed to the statement and proof of the following generalization of Lemma \ref{lem:Blowdown_C*_fixed-point_set_blowupW_equals_C*_fixed-point_set_W}.

\begin{lem}[Existence of fixed points for the $\CC^*$ action on the blowup of a complex analytic model space along a $\CC^*$-invariant linear subspace]
\label{lem:C*_action_blowup_complex_analytic_model_space}
Continue the notation of Theorem \ref{thm:Strict_transform_complex_analytic_model_space}.
Let $\rho_\CC:\CC^* \to \GL(X)$ be a representation that restricts to a representation $\rho_\CC:\CC^* \to \GL(Z)$, that is, $Z$ is a $\CC^*$-invariant linear subspace. If the linear $\CC^*$ action on $X$ restricts to a $\CC^*$ action on the closed, complex analytic subspace $Y\subset X$, then there are unique holomorphic $\CC^*$ actions on the blowup $\Bl_Z(X)$ and on the strict transform $\widetilde Y \subset \Bl_Z(X)$ of $Y$ as in Definition \ref{defn:Strict_transform_analytic_space} such that the blowup morphism $\pi_X:\Bl_X(X)\to X$ is $\CC^*$-equivariant and
\[
  \pi_X\left((\widetilde Y)^0\right) = Y^0,
\]
where $Y^0 = Y^{\CC^*}$ is the subset of fixed points of the $\CC^*$ action on $Y$ and $(\widetilde Y)^0 = (\widetilde Y)^{\CC^*}$ is the subset of fixed points of the induced $\CC^*$ action on $\widetilde Y$.
\end{lem}

\begin{rmk}[Bia{\l}ynicki--Birula decompositions need not commute with strict transforms]
\label{rmk:BB_decompositions_blowup_complex_analytic_model_space_need_not_commute_with_strict_transforms}  
As explained in Remark \ref{rmk:Submanifolds_blowup_strict_transforms} in a simpler setting, the fixed-point subspace $(\widetilde Y)^0$ described in Lemma \ref{lem:C*_action_blowup_complex_analytic_model_space} need \emph{not} coincide with the strict transform of $Y^0$.
\end{rmk}

\begin{rmk}[Statement and proof of Lemma \ref{lem:C*_action_blowup_complex_analytic_model_space} with $\CC^*$ replaced by $S^1$]
\label{rmk:S1_action_blowup_complex_analytic_model_space}  
We note that $\CC^*$ may be replaced by $S^1$ in the statement Lemma \ref{lem:C*_action_blowup_complex_analytic_model_space} since the proof for the group $S^1$ is identical to that for $\CC^*$ except for elementary changes and simplifications.
\end{rmk}

\begin{proof}[Proof of Lemma \ref{lem:C*_action_blowup_complex_analytic_model_space}]
Theorem \ref{thm:Equivariance_property_blowups_complex_manifolds} implies that the holomorphic $\CC^*$ action on $X$ has a unique lift to a holomorphic $\CC^*$ action on the blowup $\Bl_Z(X)$ such that the blowup morphism $\pi_X:\Bl_Z(X) \to X$ is $\CC^*$-equivariant. Similarly, Theorem \ref{thm:Equivariance_property_blowups_complex_analytic_spaces} implies that the  holomorphic $\CC^*$ action on $Y$ has a unique lift to a holomorphic $\CC^*$ action on the blowup $\Bl_{Z\cap Y}(Y)$ such that the blowup morphism $\pi_Y:\Bl_{Z\cap Y}(Y) \to Y$ is $\CC^*$-equivariant, where $\pi_Y$ is the restriction of $\pi_X$ to the closed, complex analytic subspace $\Bl_{Z\cap Y}(Y) \subset \Bl_Z(X)$ provided by Corollary \ref{cor:Strict_transform_closed_subspace_under_blowup_analytic_space_along_subspace}.

According to Corollary \ref{cor:Strict_transform_closed_subspace_under_blowup_analytic_space_along_subspace}, the strict transform $\widetilde Y$ of $Y \subset X$ defined by the blowup $\pi_X:\Bl_Z(X)\to X$ is equal to $\Bl_{Z\cap Y}(Y)$. Lemma \ref{lem:Equivariant_holomorphic_maps_preserve_BB_decompositions} and the fact that the blowup morphism $\pi_Y:\Bl_{Z\cap Y}(Y) \to Y$ is $\CC^*$-equivariant ensures that 
\[
  \pi_Y\left(\Bl_{Z\cap Y}(Y)^{\CC^*}\right) \subset Y^{\CC^*}.
\]  
To prove the reverse inclusion, it suffices to show that if $p \in Y$ is a fixed point of that $\CC^*$ action, then there exists a fixed point $\tilde p \in \pi_X^{-1}(p) \cap \widetilde Y$ since that will yield the inclusion
\[
   Y^{\CC^*} \subset \pi\left(\Bl_{Z\cap Y}(Y)^{\CC^*}\right).
\]  
According to Lemma \ref{lem:HilbertSpaceClosedUnderC*Action}, we may choose a Hermitian inner product on $X$ such that the orthogonal complement $Z^\perp \subset X$ of $Z$ is also $\CC^*$-invariant. By Remark \ref{rmk:Classification_complex_representations_C*}, the restrictions of $\rho_\CC$ to $Z$ and $Z^\perp$ are each diagonal matrices of the form \eqref{eq:C*_matrix_representation} with respect to choices of orthonormal bases for $Z$ and $Z^\perp$. Hence, we may assume without loss of generality that $X=\CC^n$ and that $Z\subset X$ is the standard coordinate subspace $\CC^m\hookrightarrow\CC^n$ as in \eqref{eq:Coordinate_subspace_blowup_center} with $\KK=\CC$.

To prove the existence of fixed points $\tilde p \in \widetilde Y$, we may consider without loss of generality one coordinate domain $U_j \subset \Bl_Z(X)$ in our application of Theorem \ref{thm:Strict_transform_complex_analytic_model_space}, noting that $\Bl_Z(X)$ is covered by the coordinate domains $U_j$ for $j = m+1,\ldots,n$. With respect to the local coordinates $(w_1,\ldots,w_n)$ on $U_j\subset \Bl_Z(X)$, we know from \eqref{eq:C*_action_blowup_coordinate_patch_j} that $\lambda \in \CC^*$ acts on $U_j$ by
\begin{multline*}
  (w_1,\ldots,w_n)
  \\
  \mapsto (\lambda^{l_1}w_1,\ldots,\lambda^{l_m}w_m,\lambda^{l_{m+1}-l_j}w_{m+1},\ldots,\lambda^{l_{n-1}-l_j}w_{j-1},
  \lambda^{l_j}w_j,\lambda^{l_{n+1}-l_j}w_{j+1}\ldots,\lambda^{l_n-l_j}w_n),
\end{multline*}
where the local coordinates $(w_1,\ldots,w_n)$ on $U_j \subset \Bl_Z(X)$ are related via \eqref{eq:Blowup_linear_subspace_coordinates} to the coordinates $(z_1,\ldots,z_n)$ on $X \cong \CC^n$ by $w_k=z_k$ for $k=1,\ldots,m$ and $k=j$, while $w_k=z_k/z_j$ for $k=m+1,\ldots,j-1,j+1,\ldots,n$.

Because the point $p \in X$ is a fixed point of the linear $\CC^*$ action on $X$, then any of its coordinates $z_1(p),\ldots,z_n(p)$ with non-zero weights among $l_1,\ldots,l_n$ must be zero. Therefore, any of the coordinates $w_1(\tilde p) = z_1(p),\ldots,w_m(\tilde p)) = z_m(p)$ or $w_j(\tilde p) = z_j(p)$ of $\tilde p\in U_j$ with non-zero weights among $l_1,\ldots,l_m$ and $l_j$ must also be zero. Any of the coordinates
\[
  w_{m+1} = z_{m+1}/z_j, \ldots, w_{j-1} = z_{j-1}/z_j,
  w_{j+1} = z_{j+1}/z_j, \ldots, w_n = z_n/z_j,
\]
on $U_j$ of $\tilde p$ with non-zero weights among $l_{m+1}-l_j,\ldots,l_{j-1}-l_j, l_{j+1}-l_j,\ldots,l_n-l_j$ must also be zero. Consequently, we see from the equation \eqref{eq:Exceptional_divisor_support_equation_in_local_coordinates} defining the exceptional divisor $E\cap U_j$ via $w_j=0$ and the equations \eqref{eq:Strict_transform_support_equations_in_local_coordinates} defining the strict transform $\widetilde Y\cap U_j$ via
\[
  g_{k,j}(w_1,\ldots,w_{j-1},1,w_{j+1},\ldots,w_n) = 0, \quad\text{for } k = 1,\ldots,r,
\]
that there exist fixed points $\tilde p \in \pi_X^{-1}(p)$ of the induced $\CC^*$ action on $\Bl_Z(X)$ that belong to the intersection $\widetilde Y\cap U_j$. This observation for points $\tilde p$ holds because although the ratios $z_k/z_j$ are undefined for $z_j=0$, the holomorphic functions $g_{k,j}$ are well-defined even when $z_j=0$. Moreover, since $j$ was arbitrary, the preceding observations hold for each coordinate domain $U_j$, with $j=m+1,\ldots,n$. 
\end{proof}

The following useful result is of independent interest.

\begin{lem}[Fixed-point sets of the $\CC^*$ action and the induced $S^1$ action coincide for a holomorphic $\CC^*$ action on a complex analytic space]
\label{lem:S1_and_C*_fixed_points_coincide_on_complex_analytic_space}  
Let $(X,\sO_X)$ be a complex analytic space. If $\CC^*\times X\to X$ is a holomorphic $\CC^*$ action, then the fixed-point sets of the $\CC^*$ action on $X$ and the fixed-point sets of the induced $S^1$ action coincide, that is,
\[
  X^{\CC^*} = X^{S^1}.
\]  
\end{lem}

\begin{proof}
Clearly, $X^{\CC^*} \subset X^{S^1}$ and so it suffices to prove the reverse inclusion, $X^{S^1} \subset X^{\CC^*}$. Let $p \in X$ be a fixed point of the $S^1$ action on $X$, that is, $p \in X^{S^1}$. By the Kaup Linearization Theorem \ref{thm:Kaup_linearization_theorem}, there are an $S^1$-invariant open neighborhood $U \subset X$ of $p$, an $S^1$-invariant open neighborhood $V \subset T_pX$ of the origin, and a closed,  holomorphic embedding $\varphi:U \to V$ that is $S^1$-equivariant with respect to the action of $S^1$ on $U$ and the isotropy representation $\rho$ of $S^1$ on $T_pX \cong \CC^n$, where we write $n = \dim T_pX$.

From Lemma \ref{lem:Classification_complex_representations_circle_group}, every representation $\rho:S^1\to\GL(n,\CC)$ has the form \eqref{eq:Circle_matrix_representation} and is thus unitary, with integer weights $l_1,\ldots,l_n \in \ZZ$. By Remark \ref{rmk:Classification_complex_representations_C*}, the unitary representation $\rho:S^1\to\U(n)$ given by \eqref{eq:Circle_matrix_representation} extends uniquely to a complex representation $\rho_\CC:\CC^*\to\GL(n,\CC)$ given by \eqref{eq:C*_matrix_representation} with the same integer weights. From the preceding observation, we see that there is a unique extension of the isotropy representation $\rho$ of $S^1$ on $T_pX \cong \CC^n$ to a linear $\CC^*$ representation $\rho_{\CC}$ with the same integer weights.

The real analytic action $S^1\times X \to X$ defines a homomorphism $\Phi:S^1 \to \Aut(X)$ from $S^1$ to the group $\Aut(X)$ of real analytic isomorphisms of $X$ onto itself. Since $U$ is $S^1$-invariant, we obtain a homomorphism $\Phi:S^1 \to \Aut(U)$ to the group $\Aut(U)$ of real analytic isomorphisms of $U$ onto itself that obeys the following relation:
\[
  \Phi(e^{i\theta})(z) = \varphi^{-1}\left(\rho(e^{i\theta})\varphi(z)\right),
  \quad\text{for all } e^{i\theta} \in S^1 \text{ and } z \in U.
\]
Similarly, the holomorphic action $\CC^*\times X \to X$ defines a homomorphism $\Phi_\CC:\CC^* \to \Aut(X)$ from $\CC^*$ to the group of biholomorphic isomorphisms of $X$ onto itself. Although $U$ is not necessarily $\CC^*$-invariant and $X$ in Theorem \ref{thm:Cartan_linearization_theorem_equivariant_noncompact} is assumed to be a complex manifold rather than a complex analytic space as we more generally allow here, its method of proof extends without change to give a local group $\Phi_\CC:\CC^*\to\Aut_\loc(U)$ of biholomorphic maps in the sense of Definition \ref{defn:Local_transformation_group_analytic_space} that extends the relation obeyed by $\Phi$ and $\rho$:
\[
  \Phi_\CC(\lambda)(z) = \varphi^{-1}\left(\rho_{\CC}(\lambda)\varphi(z)\right),
  \quad\text{for all } \lambda \in \CC^* \text{ and } z \in U \text{ such that }
  \rho_{\CC}(\lambda)\varphi(z) \in \varphi(U).
\]
Because the integer weights of the representations $\rho$ and $\rho_\CC$ on $T_pX \cong \CC^n$ are the same, the fixed-point sets of the corresponding $S^1$ and $\CC^*$ actions on $\CC^n$ coincide. Hence, the fixed-point set of the local $\CC^*$ action $\Phi_\CC$ on $U$ coincides with the fixed-point set of the $S^1$ action $\Phi$ on $U$. In particular, $p \in X^{\CC^*}$ and thus $X^{S^1} \subset X^{\CC^*}$ and hence we obtain the stated equality.
\end{proof}

We proved a special case of the forthcoming Corollary \ref{cor:Blowdown_C*_fixed-point_set_blowup_Y_analytic_subspace_equals_C*_fixed-point_set_Y} by direct, algebraic calculation in Lemma \ref{lem:Blowdown_C*_fixed-point_set_blowupW_equals_C*_fixed-point_set_W}, for the blowup of a complex vector space along a $\CC^*$-invariant, linear subspace. We shall prove a more general version of Corollary \ref{cor:Blowdown_C*_fixed-point_set_blowup_Y_analytic_subspace_equals_C*_fixed-point_set_Y}
by analytical methods in the forthcoming Corollary \ref{cor:Correspondence_between_Bialynicki-Birula_decompositions_manifolds} (see equation \eqref{eq:piBlZX0_is_X0}) for the blowup of a complex manifold along a $\CC^*$-invariant, embedded complex submanifold). 

\begin{cor}[Relation between the fixed-point sets of a $\CC^*$ action on a closed, complex analytic subspace of a complex manifold and the induced $\CC^*$ action on the blowup along a $\CC^*$-invariant, embedded submanifold]
\label{cor:Blowdown_C*_fixed-point_set_blowup_Y_analytic_subspace_equals_C*_fixed-point_set_Y}  
Suppose $X$ is a complex manifold, $\CC^*\times X \to X$ is a holomorphic $\CC^*$ action, and $Z\subsetneq X$ is a $\CC^*$-invariant, embedded complex submanifold. If $Y \subseteq X$ is a $\CC^*$-invariant, closed, complex analytic subspace, then we have the following relation between the subsets of fixed points of the holomorphic $\CC^*$ action on $Y$ and the holomorphic $\CC^*$ action on $\Bl_{Z\cap Y}(Y)$ induced by Theorem \ref{thm:Equivariance_property_blowups_complex_analytic_spaces}:
\[
  \pi_Y\left(\Bl_{Z\cap Y}(Y)^0\right) = Y^0.
\]
\end{cor}

\begin{rmk}[Bia{\l}ynicki--Birula decompositions need not commute with strict transforms]
\label{rmk:BB_decompositions_blowup_complex_analytic_space_need_not_commute_with_strict_transforms}  
As explained in Remark \ref{rmk:Submanifolds_blowup_strict_transforms} in a simpler setting, the fixed-point subspace $\Bl_{Z\cap Y}(Y)^0$ described in Corollary \ref{cor:Blowdown_C*_fixed-point_set_blowup_Y_analytic_subspace_equals_C*_fixed-point_set_Y} need \emph{not} coincide with the strict transform $\Bl_{Z\cap Y}(Y^0)$ of $Y^0$.
\end{rmk}

\begin{proof}[Proof of Corollary \ref{cor:Blowdown_C*_fixed-point_set_blowup_Y_analytic_subspace_equals_C*_fixed-point_set_Y}]
According to Corollary \ref{cor:Strict_transform_closed_subspace_under_blowup_analytic_space_along_subspace}, the strict transform $\widetilde Y$ of $Y \subset X$ defined by the blowup $\pi_X:\Bl_Z(X)\to X$ is equal to $\Bl_{Z\cap Y}(Y)$. Lemma \ref{lem:Equivariant_holomorphic_maps_preserve_BB_decompositions} and the fact that the blowup morphism $\pi_X:\Bl_Z(X) \to X$ is $\CC^*$-equivariant ensures that 
\[
  \pi_Y\left(\Bl_{Z\cap Y}(Y)^{\CC^*}\right) \subset Y^{\CC^*},
\]
noting that $\pi_Y$ is equal to the restriction of $\pi_X$ to $\Bl_{Z\cap Y}(Y) \subset \Bl_Z(X)$. We shall give two overlapping proofs of the reverse inclusion. For the first proof, we observe that it suffices to prove that
\[
   Y^{S^1} \subset \pi_Y\left(\Bl_{Z\cap Y}(Y)^{S^1}\right),
\]
since $Y^{\CC^*} = Y^{S^1}$ and $\Bl_{Z\cap Y}(Y)^{\CC^*} = \Bl_{Z\cap Y}(Y)^{S^1}$ by virtue of Lemma \ref{lem:S1_and_C*_fixed_points_coincide_on_complex_analytic_space}. Let $p \in Y^{S^1}$. By Corollary \ref{cor:Cartan_linearization_theorem_equivariant}, there an $S^1$-invariant, open neighborhood $U \subset X$ of $p$, an open neighborhood $V \subset T_pX$ of the origin that is invariant under the induced isotropy action of $S^1$ on $T_pX$, and an $S^1$-equivariant, biholomorphic map $\varphi:U\to V$ such that $\varphi(Z\cap U) = T_pZ\cap V$ and $\varphi(p)=0$. In particular, $\varphi(Y\cap U)$ is an $S^1$-invariant, closed, complex analytic subspace of $V$ and we are now in the setting of Lemma \ref{lem:C*_action_blowup_complex_analytic_model_space}, just with $S^1$ replacing the role of $\CC^*$. Hence, by Remark \ref{rmk:S1_action_blowup_complex_analytic_model_space}, there exists $\tilde p \in \Bl_{Z\cap Y}(Y)^{S^1}$ such that $\pi_X(\tilde p) = p$. This proves the reverse inclusion since $p$ was arbitrary.

For the second proof of the reverse inclusion, we bypass the simplification afforded by Lemma \ref{lem:S1_and_C*_fixed_points_coincide_on_complex_analytic_space} and instead replace our appeal to Corollary \ref{cor:Cartan_linearization_theorem_equivariant}, which provides local linearization for the action of $S^1$ near a fixed point, by an appeal to Theorem \ref{thm:Cartan_linearization_theorem_equivariant_noncompact}, which provides local linearization in the sense of germs for the action of $\CC^*$ near a fixed point. The reverse inclusion now follows just as in the first proof by applying Lemma \ref{lem:C*_action_blowup_complex_analytic_model_space} for a $\CC^*$ action on a complex analytic model space.
\end{proof}

\section[Bia{\l}ynicki--Birula decomposition for blowup of compact complex manifold]{Bia{\l}ynicki--Birula decomposition for a $\CC^*$ action on the blowup of a compact complex, K\"ahler manifold along an invariant complex submanifold}
\label{sec:Bialynicki-Birula decomposition_blowup_compact_complex_Kaehler_manifold_C*-action}
In Theorem \ref{thm:BB_decomposition_blowup_vector_space_along_linear_subspace_explicit}, we established the existence of a Bia{\l}ynicki--Birula decomposition for a $\CC^*$ action on the blowup of a complex vector space along a  $\CC^*$-invariant, complex linear subspace. The forthcoming Theorem \ref{thm:BB_decomposition_C*_action_blowup_compact_Kaehler_manifold_invariant_submanifold}   establishes the corresponding result for the blowup of a compact, complex K\"ahler manifold along a $\CC^*$-invariant, embedded complex submanifold. In the forthcoming Theorem \ref{thm:BB_decomposition_C*_action_blowup_complex_manifold_along_invariant_submanifold}, we shall remove the hypotheses that the complex manifold be compact or K\"ahler.

\begin{thm}[Bia{\l}ynicki--Birula decomposition for $\CC^*$ action on the blowup of a compact, complex K\"ahler manifold along an embedded, invariant, complex submanifold]
\label{thm:BB_decomposition_C*_action_blowup_compact_Kaehler_manifold_invariant_submanifold}  
Let $X$ be a compact, complex, finite-dimensional K\"ahler manifold and $Z \subset X$ be an embedded, complex submanifold. If $\Phi:\CC^* \to \Aut(X)$ is a homomorphism from the group $\CC^*$ onto a subgroup of the group $\Aut(X)$ of biholomorphic automorphisms of $X$ such that the corresponding action $\CC^*\times X\to X$ is holomorphic with at least one fixed point in $X$ and leaves $Z$ invariant, then the following hold for the blowup $\Bl_Z(X)$ of $X$ along $Z$:
\begin{enumerate}
\item\label{item:BlZX_compact_complex_Kaehler_manifold}
$\Bl_Z(X)$ is a compact, complex, K\"ahler manifold of dimension equal to that of $X$.
\item\label{item:BlZX_C*_action_lifts_to_BlZX_from_X_compact_Kaehler_manifold}
The homomorphism $\Phi:\CC^* \to \Aut(X)$ lifts uniquely to a homomorphism $\Bl_Z(\Phi):\CC^* \to \Aut(\Bl_Z(X))$ such that the blowup morphism $\pi:\Bl_Z(X)\to X$ is $\CC^*$-equivariant, the following diagram commutes (for $\pi_{\Aut}$ defined by Proposition \ref{prop:Functorial_property_blowup_complex_manifold}), and the corresponding action $\CC^*\times \Bl_Z(X)\to \Bl_Z(X)$ is holomorphic:
\[
  \begin{tikzcd}
    &\Aut(\Bl_Z(X)) \arrow[d, "\pi_{\Aut}"]
    \\
    \CC^* \arrow[r, "\Phi"] \arrow[ru, "\Bl_Z(\Phi)"] &\Aut(X)
  \end{tikzcd}
\]    
\item\label{item:BlZX0_non-empty_compact_Kaehler_manifold}
The subset $\Bl_Z(X)^0$ of fixed points of the $\CC^*$ action on $\Bl_Z(X)$ is related to the subset $X^0$ of fixed points of the $\CC^*$ action on $X$ by the identity
\[
  \pi\left(\Bl_Z(X)^0\right) = X^0.
\]  
\item\label{item:BB_decomposition_BlZX_compact_Kaehler_manifold}
The blowup $\Bl_Z(X)$ inherits plus and minus Bia{\l}ynicki--Birula decompositions as in Definition \ref{maindefn:BB_decomposition_complex_manifold} from the plus and minus Bia{\l}ynicki--Birula decompositions of $X$ provided by Theorem \ref{thm:Bialynicki-Birula_decomposition_compact_complex_Kaehler_manifold}.
\end{enumerate}
\end{thm}

\begin{proof}
Consider Item \eqref{item:BlZX_compact_complex_Kaehler_manifold}. The assertion that $\Bl_Z(X)$ is a complex manifold is given by Proposition \ref{prop:Huybrechts_2-5-3} and the facts that $\Bl_Z(X)$ is compact and K\"ahler are given by Proposition \ref{prop:Voisin_3-24}. Item \eqref{item:BlZX_C*_action_lifts_to_BlZX_from_X_compact_Kaehler_manifold} is an immediate consequence of Theorem \ref{thm:Equivariance_property_blowups_complex_manifolds}. Item \eqref{item:BlZX0_non-empty_compact_Kaehler_manifold} follows from Corollary \ref{cor:Blowdown_C*_fixed-point_set_blowup_Y_analytic_subspace_equals_C*_fixed-point_set_Y} by choosing $Y = X$.

Item \eqref{item:BB_decomposition_BlZX_compact_Kaehler_manifold} thus follows from Items \eqref{item:BlZX_compact_complex_Kaehler_manifold}, \eqref{item:BlZX_C*_action_lifts_to_BlZX_from_X_compact_Kaehler_manifold}, \eqref{item:BlZX0_non-empty_compact_Kaehler_manifold}, and Theorem \ref{thm:Bialynicki-Birula_decomposition_compact_complex_Kaehler_manifold} on the existence of a Bia{\l}ynicki--Birula decomposition for a $\CC^*$ action on a compact, complex, K\"ahler manifold with a non-empty set $\Bl_Z(X)^0$ of fixed points.  This completes the proof of Theorem \ref{thm:BB_decomposition_C*_action_blowup_compact_Kaehler_manifold_invariant_submanifold}.
\end{proof}

\begin{rmk}[Alternative proof of Item \eqref{item:BlZX0_non-empty_compact_Kaehler_manifold} in Theorem \ref{thm:BB_decomposition_C*_action_blowup_compact_Kaehler_manifold_invariant_submanifold}]
\label{item:BlZX0_non-empty_compact_Kaehler_manifold_direct_proof}  
If $p \in X$ is a fixed point of the $\CC^*$ action, so $p \in X^0$, we may choose an $S^1$-invariant, open neighborhood $U$ of $p$ in $X$ and a holomorphic coordinate chart $\varphi:X\supset U \to T_pX$ that is $S^1$-equivariant with respect to the induced action of $S^1 \subset \CC^*$ on $X$ and the (linear) isotropy action of $S^1$ on $T_pX$ and such that the action of $\CC^*$ is holomorphically linearizable at $p$ in the sense of Definition \ref{defn:Cairns_Ghys_1-3} (that is, in the sense of germs). Assume that $\varphi(p) = 0 \in T_pX$ and let $V = \varphi(U) \subset T_pX$ denote the open neighborhood of the origin in $T_pX$ defined by $(\varphi,U)$. Assume further that $X$ and $Z$ have complex dimensions $n$ and $0\leq m<n$, respectively, so $T_pX \cong \CC^n$ and $T_pZ \cong \CC^m$. Theorem \ref{thm:Cartan_linearization_theorem_equivariant_noncompact} ensures that $\varphi(U\cap Z) = V\cap T_pZ \cong V\cap\CC^m$. We can now appeal to the method of proof of  Proposition \ref{prop:Huybrechts_2-5-3} (see Huybrechts \cite[Section 2.5, pp. 99--100]{Huybrechts_2005}) and recall that, with respect to the coordinate chart $\varphi:X\supset U \to \CC^n$, the restriction of the blowup $\pi:\Bl_Z(X) \to X$ to $\pi:\Bl_Z(X)|_{\pi^{-1}(U)} \to U$ is defined by the restriction of the blowup $\pi:\Bl_{\CC^m}(\CC^n) \to \CC^n$ to $\pi:\Bl_{\CC^m}(\CC^n)|_{\pi^{-1}(V)} \to V$, where $V := \varphi(U) \subset \CC^n$. We now apply Lemma \ref{lem:Blowdown_C*_fixed-point_set_blowupW_equals_C*_fixed-point_set_W} to find $\tilde p \in \Bl_{\CC^m}(\CC^n)$ and thus $\tilde p \in \Bl_Z(X)^0$ such that $\pi(\tilde p) = p$. This proves that $X^0 \subset \pi(\Bl_Z(X)^0)$ and the simpler reverse inclusion $\pi(\Bl_Z(X)^0)\subset X^0$ follows, as usual, from Lemma \ref{lem:Equivariant_holomorphic_maps_preserve_BB_decompositions}.
\end{rmk}

\chapter[Bia{\l}ynicki--Birula decompositions for blowups of complex manifolds]{Bia{\l}ynicki--Birula decompositions for blowups of complex manifolds: Analytic approach}
\label{chap:BB_decomposition_blowup_complex_manifold_analytic}
In Section \ref{sec:Bialynicki-Birula decomposition_blowup_complex_linear_subspace_nonlinear_C*-action}, we develop the Bia{\l}ynicki--Birula decomposition for a holomorphic $\CC^*$ action on the blowup of a noncompact complex manifold along a $\CC^*$-invariant complex submanifold, where the given noncompact complex manifold is endowed with a holomorphic $\CC^*$ action. Section \ref{sec:Correspondence_between_Bialynicki-Birula_decompositions_manifold_and_blowup} provides precise details of the correspondence between Bia{\l}ynicki--Birula decompositions of a complex manifold and its blowup along an invariant, complex submanifold, together with equalities of dimensions.

\section[Bia{\l}ynicki--Birula decomposition for blowup of non-compact complex manifold]{Bia{\l}ynicki--Birula decomposition for a $\CC^*$ action on the blowup of a noncompact complex manifold along an invariant complex submanifold}
\label{sec:Bialynicki-Birula decomposition_blowup_complex_linear_subspace_nonlinear_C*-action}
In Theorem \ref{thm:BB_decomposition_blowup_vector_space_along_linear_subspace_explicit}, we established the existence of a Bia{\l}ynick--Birula decomposition for the blowup of a complex vector space along a $\CC^*$-invariant complex linear subspace. While Theorem \ref{thm:BB_decomposition_C*_action_blowup_compact_Kaehler_manifold_invariant_submanifold} provides a Bia{\l}ynicki--Birula decomposition for a $\CC^*$ action on the blowup of a complex manifold along a $\CC^*$-invariant complex submanifold, its hypotheses that the complex manifold be compact and K\"ahler (which allowed us to appeal to Theorem \ref{thm:Bialynicki-Birula_decomposition_compact_complex_Kaehler_manifold} due to Carrell and Sommese) are inconvenient. With these motivations, we shall apply Theorem \ref{thm:Extension_C*_actions_Hermitian_complex_manifolds} to prove the following generalization of Theorems \ref{thm:BB_decomposition_blowup_vector_space_along_linear_subspace_explicit} and \ref{thm:BB_decomposition_C*_action_blowup_compact_Kaehler_manifold_invariant_submanifold} that relaxes their restrictive hypotheses.

\begin{thm}[Bia{\l}ynicki--Birula decomposition for $\CC^*$ action on the blowup of a complex manifold along an invariant, complex submanifold]
\label{thm:BB_decomposition_C*_action_blowup_complex_manifold_along_invariant_submanifold}
Continue the hypotheses of Theorem \ref{thm:BB_decomposition_C*_action_blowup_compact_Kaehler_manifold_invariant_submanifold}, but omit the  hypotheses that $X$ is compact or K\"ahler. Then $\Bl_Z(X)$ is a complex manifold of dimension equal to that of $X$ and, in particular, Items \eqref{item:BlZX_compact_complex_Kaehler_manifold} and \eqref{item:BlZX_C*_action_lifts_to_BlZX_from_X_compact_Kaehler_manifold} continue to hold. If $X$ admits a plus (respectively, minus or mixed) Bia{\l}ynicki--Birula decomposition in the sense of Definition \ref{maindefn:BB_decomposition_complex_manifold}, then $\Bl_Z(X)$ also admits such a plus (respectively, minus or mixed) Bia{\l}ynicki--Birula decomposition and, in particular, Items \eqref{item:BlZX0_non-empty_compact_Kaehler_manifold} and \eqref{item:BB_decomposition_BlZX_compact_Kaehler_manifold} continue to hold. 
\end{thm}

\begin{proof}
The first conclusion that $\Bl_Z(X)$ is a complex manifold of dimension equal to that of $X$ is provided by Proposition \ref{prop:Huybrechts_2-5-3}, just as in the proof of Item \eqref{item:BlZX_compact_complex_Kaehler_manifold} in Theorem \ref{thm:BB_decomposition_C*_action_blowup_compact_Kaehler_manifold_invariant_submanifold}. Item \eqref{item:BlZX_C*_action_lifts_to_BlZX_from_X_compact_Kaehler_manifold} in Theorem \ref{thm:BB_decomposition_C*_action_blowup_compact_Kaehler_manifold_invariant_submanifold} --- which asserts that the holomorphic $\CC^*$ action on $X$ lifts uniquely to a holomorphic $\CC^*$ action on $\Bl_Z(X)$ so that the blowup morphism $\pi:\Bl_Z(X) \to X$ is $\CC^*$-equivariant --- continues to hold.

Let $\tilde z \in \Bl_Z(X)$. We claim that $\tilde z \in \Bl_Z(X)^+$ or $\Bl_Z(X)^-$, where $\Bl_Z(X)^\pm$ are the subsets defined in Definition \ref{maindefn:BB_decomposition_complex_manifold}. By hypothesis, the subset $X^0\subset X$ of fixed points of the $\CC^*$ action on $X$ is non-empty and $X$ has a (possibly mixed) Bia{\l}ynicki--Birula decomposition and so $z := \pi(\tilde z) \in X^+$ or $X^-$, where $X^\pm$ are as in Definition \ref{maindefn:BB_decomposition_complex_manifold}. We may assume without loss of generality that $z\in X^+$, in which case $\lambda\cdot z \to z_0 \in X^0$ for some $z_0 \in X^0$, as $\lambda \to 0$. (The analysis for the case $z\in X^-$ proceeds \mutatis where $\lambda\cdot z \to z_0 \in X^0$ for some $z_0 \in X^0$, as $\lambda \to \infty$.) Let 
\[
  A_z:\CC^* \ni \lambda \mapsto \lambda\cdot z \in X
\]
be the $\CC^*$-equivariant, holomorphic map \eqref{eq:Orbit_point_under_C*_action_complex_manifold} provided by our hypothesis of a holomorphic action $\CC^*\times X \to X$. Let $D$ be an open disk with closure $\bar D$ centered at the origin in $\CC$. By shrinking the radius of $D$ if necessary, we may assume that $A_z(\bar D)$ is contained in the domain $U$ of a coordinate chart $(U,\varphi)$ on $X$ with $\varphi(z_0) = 0 \in \CC^n$. Since $\bar D$ is compact and $A_z$ is continuous, the image $A_z(\bar D) \subset X$ is compact and so the preimage $\pi^{-1}(A_z(\bar D)) \subset \Bl_Z(X)$ is compact since the blowup morphism $\pi:\Bl_Z(X)\to X$ is proper by the proof of Proposition \ref{prop:Huybrechts_2-5-3}. By replacing $\tilde z$ by $\lambda\cdot\tilde z$ and $z$ by $\lambda\cdot z$ with sufficiently small $|\lambda| \in (0, 1]$ if necessary, we may assume without loss of generality that $z = \pi(\tilde z) \in A_z(D)$, since $\lambda\cdot z \to z_0 = A_z(0)$ as $\lambda\to 0$ (by our hypothesis that $X$ has a Bia{\l}ynicki--Birula decomposition in the sense of Definition \ref{maindefn:BB_decomposition_complex_manifold}), and hence that $\tilde z \in \pi^{-1}(A_z(D))$.

Let $\CC^* \times \Bl_Z(X) \to \Bl_Z(X)$ denote the holomorphic action induced by the given holomorphic action $\CC^* \times X \to X$ and Theorem \ref{thm:Equivariance_property_blowups_complex_manifolds}. The orbit
\begin{equation}
  \label{eq:Orbit_point_under_C*_action_blowup_complex_manifold}
  \Bl_Z(A_z):\CC^* \ni \lambda \mapsto \lambda\cdot \tilde z \in \Bl_Z(X)
\end{equation}
has the property that $\Bl_Z(A_z)(D)$ is contained in a compact subset of $\Bl_Z(X)$ since $\CC^*$-equivariance of the blowup morphism $\pi:\Bl_Z(X)\to X$ provided by Theorem \ref{thm:Equivariance_property_blowups_complex_manifolds} ensures that $\pi(\Bl_Z(A_z)(D)) = A_z(D)$, while the map $A_z:D\to X$ in \eqref{eq:Orbit_point_under_C*_action_complex_manifold} is continuous (and in fact holomorphic) by our hypothesis that $X$ has a Bia{\l}ynicki--Birula decomposition in the sense of Definition \ref{maindefn:BB_decomposition_complex_manifold}. In particular, $A_z(\bar D)\subset X$ is compact and so $\Bl_Z(A_z)(\bar D) = \pi^{-1}(A_z(\bar D)) \subset \Bl_Z(X)$ is compact.

The restriction $\pi:\Bl_{Z\cap U}(U) \to U$ of the blowup morphism is isomorphic to the restriction to $\varphi(U) \subset \CC^n$ of the model $\pi:\Bl_{\CC^m}(\CC^n)\to\CC^n$ described in Section \ref{sec:Blowups_analytic_manifolds_along_embedded_analytic_submanifolds}, where $Z\cap U$ is identified with an open neighborhood of the origin in a coordinate subspace $\CC^m \hookrightarrow \CC^n$ via a submanifold coordinate chart $(Z\cap U,\varphi\restriction Z)$. Furthermore, by Corollary \ref{cor:Cartan_linearization_theorem_equivariant}, we may assume that $U$ is $S^1$-invariant, that the induced $S^1$ action on $U$ is linear and unitary with respect to the standard inner product on $\CC^n$, and that $\CC^m$ is $S^1$-invariant. The blowup $\Bl_{Z\cap U}(U)$ has a K\"ahler metric given by the restriction of the K\"ahler metric on $\Bl_{\CC^m}(\CC^n)$ provided by Proposition \ref{prop:Voisin_3-24_vector_space}. According to Lemma \ref{lem:Hamiltonian_function_circle_action_blowup_complex_vector_space_along_invariant_subspace}, the induced (real analytic) $S^1$ action on $\Bl_{Z\cap U}(U)$ is Hamiltonian with respect to the symplectic form $\omega$ defined by the K\"ahler metric. Hence, Theorem \ref{thm:Extension_C*_actions_Hermitian_complex_manifolds} implies that $\Bl_Z(A_z)$ in \eqref{eq:Orbit_point_under_C*_action_blowup_complex_manifold} extends to a $\CC^*$-equivariant, holomorphic map,
\[
  \Bl_Z(A_z):\CC \to \Bl_Z(X),
\]
and so the limit $\tilde z_0 := \lim_{\lambda\to 0}\lambda\cdot\tilde z$ exists and necessarily belongs to $\Bl_Z(X)^0$, so $\tilde z \in \Bl_Z(X)^+$, as claimed.

A slight modification of the preceding argument also proves that Item \eqref{item:BlZX0_non-empty_compact_Kaehler_manifold} holds in Theorem \ref{thm:BB_decomposition_C*_action_blowup_complex_manifold_along_invariant_submanifold}. To see this, let $z_0\in X^0$. Because $X$ has a (mixed) Bia{\l}ynicki--Birula decomposition, we may write $z_0 = \lim_{\lambda\to 0}\lambda\cdot z$ for some $z\in X^+$ or $z_0 = \lim_{\lambda\to \infty}\lambda\cdot z$ for some $z\in X^-$. We may assume the former limit holds since the argument for the latter limit is almost identical. Choose $\tilde z \in \pi^{-1}(z)$ and observe that our proof that $\tilde z \in \Bl_Z(X)^+$ gave $\tilde z_0 := \lim_{\lambda\to 0}\lambda\cdot\tilde z \in \Bl_Z(X)^0$. The definitions of $z_0$ and $\tilde z_0$ and continuity and $\CC^*$-equivariance of the morphism $\pi$ yield
\[
  \pi(\tilde z_0)
  =
  \pi\left(\lim_{\lambda\to 0}\lambda\cdot\tilde z\right)
  =
  \lim_{\lambda\to 0}\pi(\lambda\cdot\tilde z)
  =
  \lim_{\lambda\to 0}\lambda\cdot\pi(\tilde z)
  =
  \lim_{\lambda\to 0}\lambda\cdot z
  = z_0,
\]
and so we obtain
\begin{equation}
\label{eq:X0_in_piBlZX0}
  X^0 \subseteq \pi\left(\Bl_Z(X)^0\right),
\end{equation}
The simpler reverse inequality
\begin{equation}
\label{eq:piBlZX0_in_X0}
  \pi\left(\Bl_Z(X)^0\right) \subseteq X^0,
\end{equation}
follows from Lemma \ref{lem:Equivariant_holomorphic_maps_preserve_BB_decompositions}. By combining the inclusions \eqref{eq:X0_in_piBlZX0} and \eqref{eq:piBlZX0_in_X0}, we obtain the equality
\[
 \pi\left(\Bl_Z(X)^0\right) = X^0
\]
and this proves that Item \eqref{item:BlZX0_non-empty_compact_Kaehler_manifold} holds in Theorem \ref{thm:BB_decomposition_C*_action_blowup_complex_manifold_along_invariant_submanifold}.

In particular, the subset $\Bl_Z(X)^0 = \Bl_Z(X)^{\CC^*}$ of fixed points of the holomorphic $\CC^*$ action on $\Bl_Z(X)$ is non-empty. The remaining properties of the Bia{\l}ynicki--Birula decomposition in Definition \ref{maindefn:BB_decomposition_complex_manifold} in Theorem \ref{thm:BB_decomposition_C*_action_blowup_compact_Kaehler_manifold_invariant_submanifold} are local and their verification follows exactly as in the proof of Theorem \ref{thm:BB_decomposition_C*_action_blowup_compact_Kaehler_manifold_invariant_submanifold}, where $X$ and thus $\Bl_Z(X)$ are compact. This verifies Item \eqref{item:BB_decomposition_BlZX_compact_Kaehler_manifold} in Theorem \ref{thm:BB_decomposition_C*_action_blowup_complex_manifold_along_invariant_submanifold} and completes the proof.
\end{proof}

\section[Bia{\l}ynicki--Birula decompositions of a complex manifold and its blowup]{Equality of dimensions and correspondence between Bia{\l}ynicki--Birula decompositions of a complex manifold and its blowup along an invariant, complex submanifold}
\label{sec:Correspondence_between_Bialynicki-Birula_decompositions_manifold_and_blowup}  
We have the following generalization of Corollary \ref{cor:C*-equivariance_blowup_map_preservation_Bialynicki-Birula_decomposition} and Lemma \ref{lem:Blowdown_C*_fixed-point_set_blowupW_equals_C*_fixed-point_set_W}.

\begin{cor}[Correspondence between Bia{\l}ynicki--Birula decompositions of a complex manifold and its blowup along an invariant, complex submanifold]
\label{cor:Correspondence_between_Bialynicki-Birula_decompositions_manifolds}  
Continue the hypotheses and notation of Theorem \ref{thm:BB_decomposition_C*_action_blowup_complex_manifold_along_invariant_submanifold}. Then the following identities hold:
\begin{subequations}
  \label{eq:piBlZXpm_in_Xpm}
  \begin{align}
    \label{eq:piBlZX0_is_X0}
    \pi\left(\Bl_Z(X)^0\right) &= X^0,
    \\
    \label{eq:piBlZX+_is_X+}
    \pi\left(\Bl_Z(X)^+\right) &= X^+,
    \\
    \label{eq:piBlZX-_is_X-}
  \pi\left(\Bl_Z(X)^-\right) &= X^-.
  \end{align}
\end{subequations}
Moreover, for all points $\tilde p \in \Bl_Z(X)^0$ and their images $p = \pi(\tilde p) \in X^0$, the following hold:
\begin{subequations}
  \label{eq:piBlZXtildep_pm_is_Xp_pm}
  \begin{align}
    \label{eq:piBlZXtildep+_is_Xp+}
   \pi\left(\Bl_Z(X)_{\tilde p}^+\right) &= X_p^+,
    \\
    \label{eq:piBlZXtildep-_is_Xp-}
   \pi\left(\Bl_Z(X)_{\tilde p}^-\right) &= X_p^-.
  \end{align}
\end{subequations}
\end{cor}

\begin{rmk}[Bialynicki--Birula decompositions need not commute with strict transforms]
\label{rmk:Bialynicki-Birula_decompositions_and_strict_transforms}
For the reasons that we described in Remark \ref{rmk:Submanifolds_blowup_strict_transforms} in the simpler case of a blowup of complex vector space along linear subspace, it need \emph{not} be the case that any one of $\Bl_Z(X)^0$, $\Bl_Z(X)^\pm$, or $\Bl_Z(X)_{\tilde p}^\pm$ is equal to the strict transform (as in Definition \ref{defn:Strict_transform_analytic_space}) $\widetilde{X^0}$, $\widetilde{X^\pm}$, or $\widetilde{X_p^\pm}$ of $X^0$, $X^\pm$, or $X_p^\pm$, respectively.
\end{rmk}  

\begin{proof}[Proof of Corollary \ref{cor:Correspondence_between_Bialynicki-Birula_decompositions_manifolds}]
The equality \eqref{eq:piBlZX0_is_X0} restates that Item \eqref{item:BlZX0_non-empty_compact_Kaehler_manifold} holds in Theorem \ref{thm:BB_decomposition_C*_action_blowup_complex_manifold_along_invariant_submanifold}.

We next observe that the proofs of the inclusions \eqref{eq:piBlZWalphapm0_in_Wpm0} in Corollary \ref{cor:C*-equivariance_blowup_map_preservation_Bialynicki-Birula_decomposition} extend without change to yield the inclusions:
\begin{subequations}
  \label{eq:piBlZX0pm_in_X0pm}
  \begin{align}
    \label{eq:piBlZX+_in_X+}
    \pi\left(\Bl_Z(X)^+\right) &\subseteq X^+,
    \\
    \label{eq:piBlZX-_in_X-}
    \pi\left(\Bl_Z(X)^-\right) &\subseteq X^-.
  \end{align}
\end{subequations}
Similarly, the proofs of the inclusions \eqref{eq:piBlZWalphatildep_pm_in_Wp_pm} extend without change to yield the inclusions:
\begin{subequations}
  \label{eq:piBlZXtildep_pm_in_Xp_pm}
  \begin{align}
    \label{eq:piBlZXtildep+_in_Xp+}
   \pi\left(\Bl_Z(X)_{\tilde p}^+\right) &\subseteq X_p^+,
    \\
    \label{eq:piBlZXtildep-_in_Xp-}
   \pi\left(\Bl_Z(X)_{\tilde p}^-\right) &\subseteq X_p^-,
  \end{align}
\end{subequations}
for all $\tilde p\in\Bl_Z(X)$ with $\pi(\tilde p) = p$.

To prove the reverse inclusion
\begin{equation}
  \label{eq:X+_in_piBlZX+}
  X^+ \subseteq \pi\left(\Bl_Z(X)^+\right),
\end{equation}
we may choose $z \in X^+$ and observe that the proof of \eqref{eq:X0_in_piBlZX0} gave a point $\tilde z \in \Bl_Z(X)^+$ such that $\pi(\tilde z) = z$ and this yields the reverse inclusion \eqref{eq:X+_in_piBlZX+}. The proof of the reverse inclusion
\begin{equation}
  \label{eq:X-_in_piBlZX-}
  X^- \subseteq \pi\left(\Bl_Z(X)^-\right),
\end{equation}
follows \mutatis that of the proof of \eqref{eq:X+_in_piBlZX+}. Combining the inclusions \eqref{eq:piBlZX+_in_X+} and \eqref{eq:X+_in_piBlZX+} yields the identity \eqref{eq:piBlZX+_is_X+} and the identity \eqref{eq:piBlZX-_is_X-} similarly follows from the inclusions \eqref{eq:piBlZX-_in_X-} and \eqref{eq:X-_in_piBlZX-}.

The proofs of the reverse inclusions
\begin{subequations}
  \label{eq:Xp_pm_in_piBlZXtildep_pm}
  \begin{align}
    \label{eq:Xp+_in_piBlZXtildep+}
    X_p^+ &\subseteq \pi\left(\Bl_Z(X)_{\tilde p}^+\right),
    \\
    \label{eq:Xp-_in_piBlZXtildep-}
    X_p^- &\subseteq \pi\left(\Bl_Z(X)_{\tilde p}^-\right),
  \end{align}
\end{subequations}
are immediate consequences of the proofs of the inclusions \eqref{eq:X+_in_piBlZX+} and \eqref{eq:X-_in_piBlZX-}.  Combining the inclusions \eqref{eq:piBlZXtildep+_in_Xp+} and \eqref{eq:Xp+_in_piBlZXtildep+} yields the identity \eqref{eq:piBlZXtildep+_is_Xp+} and the identity \eqref{eq:piBlZXtildep-_is_Xp-} similarly follows from the inclusions \eqref{eq:piBlZXtildep-_in_Xp-} and \eqref{eq:Xp-_in_piBlZXtildep-}. This completes the proof of Corollary \ref{cor:Correspondence_between_Bialynicki-Birula_decompositions_manifolds}.
\end{proof}

We shall need the following corollaries of Sard's Theorem (see Lee \cite[Theorem 6.10, p. 129]{Lee_john_smooth_manifolds}). 

\begin{prop}[Measure of the image of a smooth map]
\label{prop:Measure_image_smooth_map}
(See Lee \cite[Corollary 6.11, p. 131]{Lee_john_smooth_manifolds}.)  
Suppose $M$ and $N$ are finite-dimensional, smooth ($C^\infty$) manifolds with or without boundary, and $F:M\to N$ is a $C^\infty$ map. If $\dim M < \dim N$, then $F(M)$ has measure zero in $N$. 
\end{prop}

See Lee \cite[Problem 6.1, p. 147]{Lee_john_smooth_manifolds} for the outline of a simple proof of Proposition \ref{prop:Measure_image_smooth_map} that does not depend on the full strength of Sard's Theorem. Proposition \ref{prop:Measure_image_smooth_map} immediately yields the

\begin{cor}[Dimension of the codomain of a smooth surjective map]
\label{cor:Dimension_codomain_smooth_surjective_map}
Continue the hypotheses of Proposition \ref{prop:Measure_image_smooth_map}. If $F:M\to N$ is surjective, then $\dim M \geq \dim N$.
\end{cor}  

\begin{lem}[Dimensions of the fixed-point, stable, and unstable submanifolds in Bia{\l}ynicki--Birula decomposition of blowup of a complex vector space along an invariant linear subspace]
	\label{lem:Dimensions_fixed-point_stable_unstable_submanifolds_BB_decomposition_blowup_vector_space}
	Continue the hypotheses of Theorem \ref{thm:BB_decomposition_C*_action_blowup_complex_manifold_along_invariant_submanifold}. If $\alpha \in \sA$ and $p \in X^0$ and $\tilde p \in \Bl_Z(X)_\alpha^0$ obey $\pi(\tilde p) = p$, then 
	\begin{subequations}
		\label{eq:Dimensions_fixed-point_stable_unstable_submanifolds_BB_blowup_vector_space_equalities_alpha}
		\begin{align}
			\label{eq:Dimension_fixed-point_submanifold_BB_blowup_vector_space_alpha}
			\dim_\CC \Bl_Z(X)_\alpha^0 &= \dim_\CC T_pX^0,
			\\
			\label{eq:Dimension_stable_submanifold_BB_blowup_vector_space_alpha}
			\dim_\CC \Bl_Z(X)_{\alpha,\tilde p}^+ &= \dim_\CC X_p^+, 
			\\
			\label{eq:Dimension_unstable_submanifold_BB_blowup_vector_space_alpha}
			\dim_\CC \Bl_Z(X)_{\alpha,\tilde p}^- &= \dim_\CC X_p^-.
		\end{align}  
	\end{subequations}
\end{lem}

\begin{proof}
	The blowup morphism $\pi:\Bl_Z(X) \to X$ constructed in Section \ref{sec:Blowups_analytic_manifolds_along_embedded_analytic_submanifolds} is a smooth, surjective map. By Theorem \ref{thm:BB_decomposition_C*_action_blowup_complex_manifold_along_invariant_submanifold} (and Definition \ref{maindefn:BB_decomposition_complex_manifold}), each connected component $\Bl_Z(X)_\alpha^0$ of the fixed-point set $\Bl_Z(X)^0$ is an embedded, complex submanifold of $\Bl_Z(X)$ and by \eqref{eq:piBlZX0_is_X0} in Corollary \ref{cor:Correspondence_between_Bialynicki-Birula_decompositions_manifolds}, the restriction $\pi:\Bl_Z(X)_\alpha^0\to X^0$ is surjective, as well as smooth. Hence, Corollary \ref{cor:Dimension_codomain_smooth_surjective_map} implies that
	\begin{equation}
		\label{eq:Dimension_fixed-point_submanifold_BB_blowup_vector_space_alpha_geq}
		\dim_\CC \Bl_Z(X)_\alpha^0 \geq \dim_\CC X^0.
	\end{equation}
	By Theorem \ref{thm:BB_decomposition_C*_action_blowup_complex_manifold_along_invariant_submanifold} (and Definition \ref{maindefn:BB_decomposition_complex_manifold}), the fiber $\Bl_Z(X)_{\alpha,\tilde p}^+$ of the projection $\pi_\alpha^+:\Bl_Z(X)_\alpha^+\to \Bl_Z(X)_\alpha^0$ is an embedded, complex submanifold of $\Bl_Z(X)$ and by \eqref{eq:piBlZXtildep+_is_Xp+} Corollary \ref{cor:Correspondence_between_Bialynicki-Birula_decompositions_manifolds}, the restriction $\pi:\Bl_Z(X)_{\alpha,\tilde p}^+\to X_p^+$ is surjective, as well as smooth. Hence, Corollary \ref{cor:Dimension_codomain_smooth_surjective_map} implies that
	\begin{equation}
		\label{eq:Dimension_stable_submanifold_BB_blowup_vector_space_alpha_geq}
		\dim_\CC \Bl_Z(X)_{\alpha,\tilde p}^+ \geq \dim_\CC X_p^+.
	\end{equation}
	The same reasoning yields the inequality
	\begin{equation}
		\label{eq:Dimension_unstable_submanifold_BB_blowup_vector_space_alpha_geq}
		\dim_\CC \Bl_Z(X)_{\alpha,\tilde p}^- \geq \dim_\CC X_p^-.
	\end{equation}
	But the Bia{\l}ynicki--Birula decomposition in Theorem
	\ref{thm:BB_decomposition_C*_action_blowup_complex_manifold_along_invariant_submanifold} also yields the equality
	\[
	\dim_\CC \Bl_Z(X)_\alpha^0 + \dim_\CC \Bl_Z(X)_{\alpha,\tilde p}^+ + \dim_\CC \Bl_Z(X)_{\alpha,\tilde p}^-
	=
	\dim_\CC \Bl_Z(X),
	\]
	and the equality below is immediate from the definition of the subspaces $X^0$ and $X_p^\pm$ of $X$:
	\[
	\dim_\CC X^0 + \dim_\CC X_p^+ + \dim_\CC X_p^- = \dim_\CC X.
	\]
	Lastly, $\dim_\CC\Bl_Z(X) = \dim_\CC X$ by the construction of the blowup morphism $\pi:\Bl_Z(X) \to X$. Combining these observations yields the equality
	\[
	\dim_\CC \Bl_Z(X)_\alpha^0 + \dim_\CC \Bl_Z(X)_{\alpha,\tilde p}^+ + \dim_\CC \Bl_Z(X)_{\alpha,\tilde p}^-
	=
	\dim_\CC X^0 + \dim_\CC X_p^+ + \dim_\CC X_p^-.
	\]
	Combining the preceding equality with the inequalities \eqref{eq:Dimension_fixed-point_submanifold_BB_blowup_vector_space_alpha_geq}, \eqref{eq:Dimension_stable_submanifold_BB_blowup_vector_space_alpha_geq}, and \eqref{eq:Dimension_unstable_submanifold_BB_blowup_vector_space_alpha_geq} yields the desired equalities \eqref{eq:Dimension_fixed-point_submanifold_BB_blowup_vector_space_alpha} and \eqref{eq:Dimensions_fixed-point_stable_unstable_submanifolds_BB_blowup_vector_space_equalities_alpha}.
\end{proof}

We can now complete the proofs of the remainder of the principal results stated in Chapter \ref{chap:Introduction}.

\begin{proof}[Proof of Theorem \ref{mainthm:BB_decomposition_blowup_complex_manifold_C*_action_along_submanifold}]
Item \eqref{item:BlZX_complex_manifold} follows from Proposition \ref{prop:Huybrechts_2-5-3}. Item \eqref{item:BlZX_C*_action_lifts_to_BlZX_from_X_complex_manifold} follows from Proposition \ref{prop:Functorial_property_blowup_complex_manifold} and Theorem \ref{thm:Equivariance_property_blowups_complex_manifolds}. Items \eqref{item:BlZX0_non-empty_complex_manifold} and \eqref{item:BB_decomposition_BlZX_complex_manifold} follow from Theorem \ref{thm:BB_decomposition_C*_action_blowup_complex_manifold_along_invariant_submanifold}. Item \eqref{item:piBlZXpm0_is_Xpm0_complex_manifold} follows from Corollary \ref{cor:Correspondence_between_Bialynicki-Birula_decompositions_manifolds}. The equalities \eqref{eq:Complex_dimension_BlZX0_fibers_pm_complex_manifold} in Item \eqref{item:Dimension_BlZXalphapm0_complex_manifold} follow from Lemma \ref{lem:Dimensions_fixed-point_stable_unstable_submanifolds_BB_decomposition_blowup_vector_space} and
the equalities \eqref{eq:BBsignature_preserved_by_blowup} follow from the equalities \eqref{eq:Complex_dimension_BlZX0_fibers_pm_complex_manifold} and the Definition \ref{maindefn:Stable_unstable_submanifolds_BB_index_co-index_nullity} of the Bia{\l}ynicki--Birula signature.
\end{proof}

\begin{proof}[Proof of Theorem \ref{mainthm:BB_decomposition_strict_transform_complex_analytic_subspace}]
The conclusions in Item \eqref{item:X'_and_Y'_complex_manifolds} follow from Embedded Resolution of Singularities for complex analytic spaces (see Theorems \ref{thm:Embedded_resolution_of_singularities_analytic_space} and \ref{thm:Embedded_resolution_of_singularities_algebraic_scheme}). The conclusions in Items \eqref{item:C*_action_lifts_to_X'_from_X}, \eqref{item:X'0_non-empty_complex_manifold}, \eqref{item:BB_decomposition_X'_complex_manifold}, \eqref{item:PiX'0pm_in_Xpm_and_PiX'p'_pm_is_Xp_pm_complex_manifold}, and \eqref{item:Dimension_X'alphapm_X'alpha0_complex_manifold} follow from repeated application of Theorem \ref{mainthm:BB_decomposition_blowup_complex_manifold_C*_action_along_submanifold} for the behavior of the Bia{\l}ynicki--Birula decomposition for a $\CC^*$ action on a complex manifold under blowup along a $\CC^*$-invariant, embedded complex submanifold, noting that the resolution morphism $\Pi:X'\to X$ can be expressed as a composition of $\CC^*$-equivariant blowup morphisms $\pi_i:X_i\to X_{i-1}$ along centers $Z_{i-1}\subset X_{i-1}$ that are $\CC^*$-invariant, embedded complex submanifolds.
\end{proof}

\begin{proof}[Proof of Corollary \ref{maincor:BB_decomposition_strict_transform_complex_analytic_subspace}]
We begin by observing that if the resolution morphism $\Pi:X'\to X$ comprises a single blowup morphism $\pi:\Bl_Z(X)\to X$, then Item \eqref{item:Y'0_non-empty_strict_transform_complex_analytic_subspace} follows from our hypothesis that $Y^0$ is non-empty, the identity
\[
  \pi\left(\Bl_{Z\cap Y}(Y)^0\right) = Y^0
\]
provided by Corollary \ref{cor:Blowdown_C*_fixed-point_set_blowup_Y_analytic_subspace_equals_C*_fixed-point_set_Y}, and the fact that $\Bl_{Z\cap Y}(Y)$ is equal to the strict transform $Y' = \widetilde Y$ of $Y$ by Corollary \ref{cor:Strict_transform_closed_subspace_under_blowup_analytic_space_along_subspace} since the the combination of the preceding facts show that $Y^{\prime,0}$ is non-empty. In general, $\Pi:X'\to X$ comprises a composition of finitely many blowup morphisms as in Theorems \ref{thm:Embedded_resolution_of_singularities_algebraic_scheme} and  \ref{thm:Embedded_resolution_of_singularities_analytic_space}. In this general setting, Item  \eqref{item:Y'0_non-empty_strict_transform_complex_analytic_subspace} follows by repeated application of the result for a single blowup just described, noting that each intermediate strict transform $Y_i$ is a closed, complex analytic subspace of the intermediate strict transform $X_i$, which is a complex manifold with a holomorphic action $\CC^*\times X_i\to X_i$ that leaves $Y_i$ invariant with at least one fixed point and which is equipped with a $\CC^*$-equivariant blowup morphism $\pi_i:X_i\to X_{i-1}$ along a center $Z_{i-1}\subset X_{i-1}$ that is a $\CC^*$-invariant, embedded complex submanifold.

The conclusions in Item \eqref{item:BB_decomposition_strict_transform_complex_analytic_subspace} follow from the equalities \eqref{eq:BB_decomposition_C*_invariant_complex_submanifold} in Theorem \ref{mainthm:BB_decomposition_C*_invariant_complex_submanifold} since the strict transform $Y'$ is an embedded, complex submanifold of the strict transform $X'$, which is a complex manifold with a holomorphic action $\CC^*\times X'\to X'$ that leaves $Y'$ invariant with at least one fixed point. Item \eqref{item:PiY'0pm_in_Ypm_strict_transform_complex_analytic_subspace} follows by repeated application of Corollary \ref{cor:Correspondence_between_Bialynicki-Birula_decompositions_manifolds}.

Consider Item \eqref{item:BB_dimensions_nullity_co-index_index_strict_transform_complex_analytic_subspace}. The equality \eqref{eq:dim_p_Yprime_equals_BB_nullity_plus_coindex_plus_index} is an immediate consequence of \eqref{eq:Dim_Tp_X_equals_BB_nullity_plus_coindex_plus_index}. According to Theorem \ref{thm:Narasimhan_section_3-1_theorem_1_p_41}, there is an open neighborhood $U \subset Y$ of $p$ such that if $\dim\sO_{Y^0,p} > 0$, then $\dim\sO_{Y^0,p} = \dim (Y^0)_\sm\cap U$. Similarly, if $\dim\sO_{Y^+,p} > 0$, then $\dim\sO_{Y^+,p} = \dim (Y^+)_\sm\cap U$, while if $\dim\sO_{Y^-,p} > 0$, then $\dim\sO_{Y^-,p} = \dim (Y^-)_\sm\cap U$. According to \eqref{eq:PiY'0pm_in_Ypm_strict_transform_complex_analytic_subspace} and \eqref{eq:PiY'p'_pm_is_Yp_pm_strict_transform_complex_analytic_subspace} in Item \eqref{item:PiY'0pm_in_Ypm_strict_transform_complex_analytic_subspace} of Corollary \ref{maincor:BB_decomposition_strict_transform_complex_analytic_subspace}, the resolution morphism $\Pi:X'\to X$ yields the following surjections:
\[
  \Pi: Y^{\prime,0} \to \dim Y^0, \quad \Pi:Y_{p'}^{\prime,+} \to Y_p^+, \quad\text{and}\quad \Pi:Y_{p'}^{\prime,-} \to Y_p^-.
\]
Hence, Corollary \ref{cor:Dimension_codomain_smooth_surjective_map} yields the following inequalities:
\[
  \dim_pY^{\prime,0} \geq \dim (Y^0)_\sm\cap U,
  \quad \dim_pY^{\prime,+} \geq \dim (Y^+)_\sm\cap U,
  \quad\text{and}\quad \dim_pY^{\prime,-} \geq \dim (Y^-)_\sm\cap U.
\]
Therefore, we obtain
\[
  \dim_pY^{\prime,0} \geq \dim\sO_{Y^0,p},
  \quad \dim_pY^{\prime,+} \geq \dim\sO_{Y^+,p},
  \quad\text{and}\quad \dim_pY^{\prime,-} \geq \dim\sO_{Y^-,p}.
\]
By combining the preceding inequality with the Definitions \ref{maindefn:Stable_unstable_submanifolds_BB_index_co-index_nullity} and \ref{maindefn:Stable_unstable_subspaces_BB_index_co-index_nullity} of the Bia{\l}ynicki--Birula signature, we obtain the inequalities \eqref{eq:BB_nullity_coindex_index_preserved_resolution_singularities}. Moreover, because $\dim_pY^{\prime,0} = \dim Y_\sm \cap U = \dim_pY$ (via Theorems \ref{thm:Narasimhan_section_3-1_theorem_1_p_41}, \ref{thm:Embedded_resolution_of_singularities_algebraic_scheme}, and \ref{thm:Embedded_resolution_of_singularities_analytic_space}), by combining the preceding inequalities with \eqref{eq:dim_p_Yprime_equals_BB_nullity_plus_coindex_plus_index}, we obtain the inequality \eqref{eq:dim_p_Y_equals_BB_nullity_plus_coindex_plus_index}. This verifies Item  \eqref{item:BB_dimensions_nullity_co-index_index_strict_transform_complex_analytic_subspace}.

Consider Item \eqref{eq:BB_index_coindex_at_p_positive_implies_not_local_min_max}. This assertion simply repeats that of Item \eqref{item:BB_index_coindex_at_p_positive_basic_implies_not_local_min_max} in Theorem \ref{mainthm:BB_decomposition_C*_invariant_complex_analytic_subspace} and so its proof is an immediate consequence. This completes the verification of Item \eqref{eq:BB_index_coindex_at_p_positive_implies_not_local_min_max} and the proof of Corollary \ref{maincor:BB_decomposition_strict_transform_complex_analytic_subspace}.
\end{proof}

\appendix

\chapter{Technical Results and Definitions}
\label{chap:Technical_results_definitions}
We begin in Section \ref{sec:Open_cone_neighborhoods_closed_linear_subspaces_of_Hilbert_space} by giving definitions of open cone neighborhoods of closed linear subspaces of a Hilbert space. Section \ref{sec:Transversality_normal_crossing_divisors} provides definitions of transversal, clean, and normal crossing intersections. We conclude in Section \ref{sec:Adapted_analytic_coordinates} with a discussion of adapted analytic coordinate charts, extending the classical definition of submanifold coordinate charts.

\section{Open cone neighborhoods of closed linear subspaces of a Hilbert space}
\label{sec:Open_cone_neighborhoods_closed_linear_subspaces_of_Hilbert_space}
The local properties of a Morse--Bott function near a critical point is facilitated via the concept of an open cone neighborhood of a linear subspace of an inner product space, a concept that we introduce here and motivated by the usual definition of a right circular, double cone around an axis in $\RR^3$ with aperture angle $2\theta$, for $\theta\in[0,\pi/2)$, and apex at the origin.

\begin{defn}[Open cone neighborhood of a closed linear subspace]
\label{defn:Cone_neighborhood_subspace}  
Let $\KK=\RR$ or $\CC$ and $E$ be a Hilbert space over $\KK$. If $P \subseteqq E$ is a closed, linear, non-zero subspace and $\theta\in(0,\pi/2]$ is a constant, then
\begin{equation}
  \label{eq:Cone_neighborhood_subspace}  
  \Cone_\theta P
  := \left\{ v \in E: \|\pi_P^\perp v\|_E < \|\pi_Pv\|_E\tan\theta \text{ or } v = 0\right\},
\end{equation}
is an \emph{open cone neighborhood with axis $P \subset E$ and aperture $2\theta$}, where $\pi_P:E\to P$ is orthogonal projection and $\pi_P^\perp := \id_E - \pi_P$. 
\end{defn}

Note that $\Cone_\theta P\less\{0\}$ is an open subset of $E$, with $\Cone_{\pi/2} P = E$ and $\Cone_\theta E = E$. In our applications, we shall primarily be interested in open cone neighborhoods of linear subspaces $P\subsetneqq E$ with aperture $\pi/2$, in which case $\theta=\pi/4$ and $\tan\theta=1$.

Our definition \eqref{eq:Cone_neighborhood_subspace} of $\Cone_\theta P$ may be written in a more insightful equivalent form. We recall that the cosine of the angle between two non-zero vectors in an inner product space over $\KK=\RR$ or $\CC$ is defined as usual by (see, for example, Axler \cite[Exercise 6.A.14]{Axler_linear_algebra_done_right} when $\KK=\RR$)
\begin{equation}
  \label{eq:Inner_product_definition_cosine_angle_between_two_vectors}
  \cos \angle (v, w)
  :=
  \frac{\Real\langle v, w\rangle_E}{\|v\|_E \|w\|_E} \in [-1,1], \quad\text{for } v, w \in E\less\{0\}.
\end{equation}
We now make the

\begin{defn}[Angle between a non-zero vector and a closed linear subspace]
\label{defn:Angle_between_vector_and_subspace}
Continue the assumptions of Definition \ref{defn:Cone_neighborhood_subspace}. We define the angle between a non-zero vector $v \in E$ and the closed linear subspace $P \subseteqq E$ by
\begin{equation}
  \label{eq:Angle_between_vector_and_subspace}
  \angle (v, P)
  :=
  \angle (v, \pi_Pv),
\end{equation}
where $\cos\angle (v, \pi_Pv)$ is defined by \eqref{eq:Inner_product_definition_cosine_angle_between_two_vectors}.
\end{defn}

We then have the

\begin{lem}[Alternative definition of an open cone neighborhood of a closed linear subspace]
\label{lem:Properties_open_cones}
Continue the assumptions of Definition \ref{defn:Cone_neighborhood_subspace}. Then
\begin{equation}
  \label{eq:Cone_neighborhood_subspace_angle_definition}
  \Cone_\theta P = \left\{ v \in E: 0 \leq \angle(v,P) < \theta \text{ or } v = 0 \right\}.
\end{equation}
\end{lem}

\begin{proof}
We first write the expression on the right-hand side of \eqref{eq:Cone_neighborhood_subspace_angle_definition} in the equivalent form
\[
  C_\theta P := \left\{ v \in E: \cos\angle(v,\pi_Pv) > \cos \theta, \text{ if } v\neq 0 \right\}.
\]
Observe that for any non-zero vector $v\in E$ we have
\begin{equation}
  \label{eq:Plane_geometry_definition_cosine_angle_between_two_vectors}
  \cos\angle(v,\pi_Pv) = \frac{\|\pi_Pv\|_E}{\|v\|_E}, \quad\text{for } v \in E\less\{0\}.
\end{equation}
Indeed, by definition \eqref{eq:Inner_product_definition_cosine_angle_between_two_vectors} we see that
\[
  \cos \angle (v, \pi_Pv)
  =
  \frac{\Real\langle v, \pi_Pv\rangle_E}{\|v\|_E \|\pi_Pv\|_E}
  =
  \frac{\langle \pi_Pv, \pi_Pv\rangle_E}{\|v\|_E \|\pi_Pv\|_E}
  =
  \frac{\|\pi_Pv\|_E}{\|v\|_E}, 
\]
and thus \eqref{eq:Plane_geometry_definition_cosine_angle_between_two_vectors} follows, which agrees with the formula from plane geometry which also gives
\begin{equation}
  \label{eq:Plane_geometry_definition_sine_angle_between_two_vectors}
  \sin\angle(v,\pi_Pv) = \frac{\|\pi_P^\perp v\|_E}{\|v\|_E}, \quad\text{for } v \in E\less\{0\}.
\end{equation}
(Alternatively, to verify \eqref{eq:Plane_geometry_definition_sine_angle_between_two_vectors} we may use the identity $\|v\|_E^2 = \|\pi_Pv\|_E^2 + \|\pi_P^\perp v\|_E^2$ --- see, for example, Rudin \cite[Theorem 4.11 (d)]{RudinRealComplex} --- and observe that
\[
  \sin^2\angle(v,\pi_Pv) = 1 - \cos^2\angle(v, \pi_Pv) = 1 - \frac{\|\pi_Pv\|_E^2}{\|v\|_E^2}
  = \frac{\|v\|_E^2-\|\pi_Pv\|_E^2}{\|v\|_E^2} = \frac{\|\pi_P^\perp v\|_E^2}{\|v\|_E^2},
\]
and apply the fact that $\angle(v,\pi_Pv)$ is nonnegative by definition \eqref{eq:Cone_neighborhood_subspace_angle_definition} when taking square roots.) We thus obtain
\begin{equation}
  \label{eq:Tan_angle_between_two_vectors}
  \tan\angle(v,\pi_Pv) = \frac{\|\pi_P^\perp v\|_E}{\|\pi_Pv\|_E}, \quad\text{for } v \in E\less\{0\},
\end{equation}
by combining \eqref{eq:Plane_geometry_definition_cosine_angle_between_two_vectors} and \eqref{eq:Plane_geometry_definition_sine_angle_between_two_vectors}. Therefore, if $v\neq 0$,
\begin{multline*}
  v\in C_\theta P \iff 0\leq \angle(v,\pi_Pv) < \theta \iff 0\leq \tan\angle(v,\pi_Pv) < \tan\theta
  \\
  \iff\frac{\|\pi_P^\perp v\|_E}{\|\pi_Pv\|_E} < \tan\theta \iff v \in \Cone_\theta P,
\end{multline*}
and thus $C_\theta P = \Cone_\theta P$, so definitions \eqref{eq:Cone_neighborhood_subspace} and \eqref{eq:Cone_neighborhood_subspace_angle_definition} are equivalent.
\end{proof}

The next lemma provides a useful application of cone neighborhoods.

\begin{lem}[Orthogonal projection of a closed linear subspace of an open cone]
\label{lem:Orthogonal_projection_subspace_cone_neighborhood_subspace}
Continue the assumptions of Definition \ref{defn:Cone_neighborhood_subspace}. If $Q \subseteqq E$ is a closed, linear subspace such that $Q \subset \Cone_\theta P$, then orthogonal projection $\pi_P:Q \to P$ is a monomorphism and if in addition $Q$ and $P$ are finite-dimensional with $\dim Q = \dim P$, then $\pi_P:Q \to P$ is an isomorphism.
\end{lem}

\begin{proof}
Suppose $v \in \Cone_\theta P$. Lemma \ref{lem:Properties_open_cones} implies that either $v=0$ or
\[
  \cos\angle (v,P) > \cos\theta.
\]
By hypothesis, $\theta \in [0,\pi/2)$ and so $\cos\theta > 0$. For any $w\in E\less\{0\}$, Definition \ref{defn:Angle_between_vector_and_subspace} gives $\cos\angle (w,P) = \cos\angle (w,\pi_Pw)$. Hence, either $v=0$ or
\[
  \cos\angle (v,\pi_Pv) > \cos\theta.
\]
The identity \eqref{eq:Plane_geometry_definition_cosine_angle_between_two_vectors} asserts that for any any $w\in E\less\{0\}$,
\[
  \cos\angle (w,\pi_Pw) = \frac{\|\pi_Pw\|_E}{\|w\|_E}.
\]
Therefore, either $v=0$ or
\[
  \frac{\|\pi_Pv\|_E}{\|v\|_E} > \cos\theta.
\]  
Hence, if $\pi_Pv=0$, then we must have $v=0$ and so $\pi_P:Q\to E$ is injective since $v\in\Cone_\theta P$ was arbitrary. If $\dim Q = \dim P < \infty$, then $\pi_P:Q\to E$ is bijective and thus an isomorphism of finite-dimensional vector spaces.    
\end{proof}

\begin{rmk}[Angles between subspaces of an inner product space]
\label{sec:Angles_between_subspaces_finite-dimensional_inner_product_space}
We shall appeal to Definition \ref{defn:Cone_neighborhood_subspace} when seeking to perturb non-zero linear subspace $P\subset E$ to a linear subspace $\widehat P \subset E$ of the same dimension that is transverse, relative to $E$, to another non-zero linear subspace $V \subsetneqq E$ and also contained in $\Cone_\theta P\cup\{0\}$, so that in this sense the ``angle between subspaces'' $P$ and $\widehat P$ is less than $\theta$.  

The definition of \emph{principal angles} between between linear subspaces of Euclidean space was originated by Jordan \cite{Jordan_1875} and we outline his idea here, following the exposition by Miao and Ben--Israel \cite[Section 1]{Miao_Ben-Israel_1992}; see also Afriat \cite{Afriat_1957}, Gal\'{a}ntai and Heged\H{u}s \cite[Definition 2]{Galantai_Hegedus_2006}, Hotelling \cite{Hotelling_1936}, and Jiang \cite[Section 5]{Jiang_1996}. Let $L, M$ be linear subspaces of $\RR^n$ with $\dim L = l \leq \dim M = m$. Then the \emph{principal angles} between $L$ and $M$,
\begin{equation}
  \label{eq:Miao_Ben-Israel_1-1}
  0 \leq \theta_1 \leq \theta_2 \leq \cdots \leq \theta_l \leq \frac{\pi}{2}
\end{equation}
are defined by
\begin{multline}
  \label{eq:Miao_Ben-Israel_1-2}
  \cos\theta_i := \frac{\langle v_i, w_i\rangle}{\|v_i\| \|w_i\|}
  \\
  = \max\left\{\frac{\langle v, w\rangle}{\|v\| \|w\|}: v\in L, w \in M, v\perp v_k, w\perp w_k, k = 1\ldots, i-1\right\},
\end{multline}
where
\begin{equation}
  \label{eq:Miao_Ben-Israel_1-3}
  (v_i, w_i) \in L\times M, \quad i = 1,\ldots,l,
\end{equation}
are the corresponding $l$ pairs of \emph{principal vectors}. Miao and Ben--Israel note that
\begin{equation}
  \label{eq:Miao_Ben-Israel_1-4}
  \theta_1 = \cdots = \theta_k = 0 < \theta_{k+1} \iff \dim L\cap M = k, 
\end{equation}
and that if $\dim L = \dim M = 1$, then $\theta_1$ is the (nonobtuse) angle between the lines $L$ and $M$. Jiang \cite[Definition, p. 116]{Jiang_1996} defines a \emph{higher dimensional angle} $\theta$ between $L$ and $M$ with the property that \cite[Theorem 5]{Jiang_1996}
\[
  \cos\theta = \cos\theta_1\cdot\cos\theta_2\cdots\cos\theta_l.
\]
See Gal\'{a}ntai and Heged\H{u}s \cite{Galantai_Hegedus_2006}, Gunawan, Neswan, and Setya-Budhi \cite{Gunawan_Neswan_Setya-Budhi_2005}, Halmos \cite{Halmos_1969}, and Knyazev, Jujunashvili, and Argentati \cite{Knyazev_Jujunashvili_Argentati_2010} for related results, including extensions to complex inner product spaces by Gal\'{a}ntai and Heged\H{u}s \cite{Galantai_Hegedus_2006} and infinite-dimensional inner product spaces by Knyazev, Jujunashvili, and Argentati \cite{Knyazev_Jujunashvili_Argentati_2010}.
\end{rmk}

\section{Transversal, clean, and normal crossing intersections}
\label{sec:Transversality_normal_crossing_divisors}
We generalize the usual definition of transversal intersection in differential topology.

\begin{defn}[Transversal intersection of smooth submanifolds of an ambient smooth manifold]
\label{defn:Transversal_intersection}
(See Lee \cite[Chapter 6, p. 143]{Lee_john_smooth_manifolds} for the case $k=2$.)  
If $k\geq 2$ is an integer and $S_1,\ldots,S_k$ are embedded smooth submanifolds of a smooth manifold $M$ and $p\in S_1\cap \cdots \cap S_k$ is a point, then $S_1,\ldots,S_k$ \emph{intersect transversely at $p$ relative to $M$} if their tangent spaces at $p$ obey
\begin{equation}
  \label{eq:Transverse_intersection_submanifolds}
  T_pS_1 + \cdots + T_pS_k = T_pM
\end{equation}
and write $S_1\transv_p \cdots \transv_pS_k \subset M$ (or $\rel M$). One says that $S_1,\cdots, S_k$ \emph{intersect transversely relative to $M$} and write $S_1\transv\cdots\transv S_k \subset M$ (or $\rel M$) if \eqref{eq:Transverse_intersection_submanifolds} holds for all $p\in S_1\cap \cdots \cap S_k$. 
\end{defn}

In Definition \ref{defn:Transversal_intersection}, we emphasize the dependence of transversal intersection on the ambient manifold $M$. However, one also has $S_1\transv\cdots\transv S_k$ if the subset $S_1\cap \cdots \cap S_k$ is empty and this trivial case is independent of the embedding. For our applications, we shall need a generalization of this concept of transversal intersection from standard differential topology to a concept of clean intersection discussed, in various forms, by Bott \cite[Section 5, p. 194]{Bott_1956}, Faber and Hauser \cite[p. 379 and pp. 390--391]{Faber_Hauser_2010}, Hauser \cite[Chapter 2, Section 10, p. 380]{Hauser_2003} and \cite[Definition 3.19, p. 15]{Hauser_2014}, Hu \cite[Section 1, p. 4737]{Hu_2003}, and Li \cite[Section 5.1, p. 553]{Li_2009}.

\begin{defn}[Clean intersection of smooth submanifolds of an ambient smooth manifold]
\label{defn:Clean_intersection}
(Compare Bott \cite[Definition 5.1, p. 195]{Bott_1956}, Faber and Hauser \cite[p. 390]{Faber_Hauser_2010}, and Li \cite[Section 5.1.1, p. 554]{Li_2009}.)
Continue the assumptions of Definition \ref{defn:Transversal_intersection}. The submanifolds $S_1,\cdots, S_k$ \emph{intersect cleanly at $p$} if there is an open neighborhood $U\subset M$ of $p$ such that $S_1\cap \cdots \cap S_k\cap U$ is an embedded smooth submanifold of $U$ and their tangent spaces at $p$ obey
\begin{equation}
  \label{eq:Clean_intersection_tangent_spaces}
  T_p(S_1\cap \cdots \cap S_k) = T_pS_1 \cap \cdots \cap T_pS_k.
\end{equation}
If \eqref{eq:Clean_intersection_tangent_spaces} holds at every point of $S_1\cap \cdots \cap S_k$, then $S_1,\cdots, S_k$ \emph{intersect cleanly}.
\end{defn}

Faber and Hauser \cite[p. 390]{Faber_Hauser_2010} and Li \cite[Section 5.1.1, p. 554]{Li_2009} provide variants of Definition \ref{defn:Clean_intersection} in the categories of analytic varieties, algebraic varieties, and schemes, whereas Bott's definition is in the category of smooth manifolds. When $k=2$, we recall that if $S_1$ and $S_2$ intersect transversely at $p\in S_1\cap S_2$ relative to $M$, then \eqref{eq:Clean_intersection_tangent_spaces} necessarily holds (see Lee \cite[Exercise 6-10, p. 148]{Lee_john_smooth_manifolds}) and thus $S_1$ and $S_2$ intersect cleanly at $p$ in the sense of Definition \ref{defn:Clean_intersection}. In algebraic geometry, the following variant of Definition \ref{defn:Clean_intersection} is frequently used in the context of resolution of singularities.

\begin{defn}[Normal crossings intersection of subvarieties of an ambient variety]
\label{defn:Normal_crossings_intersection}
(Compare Bodn\'ar \cite[Section 2, p. 4, Definition 2]{Bodnar_2004}, Bruschek and Wagner \cite[Section 3, p. 137]{Bruschek_Wagner_2011}, Faber and Hauser \cite[p. 379 and pp. 390--391]{Faber_Hauser_2010}, Hauser \cite[Chapter 2, Section 10, p. 380]{Hauser_2003} and \cite[Section 3, Definition 3.15, p. 14 and Proposition 3.17, p. 15]{Hauser_2014}.)
Let $k\geq 2$ be an integer and $W_1,\ldots,W_k$ be (algebraic or analytic) subvarieties of an ambient (algebraic or analytic) variety $X$. If $p\in W_1\cap\cdots\cap W_k$ is a regular point of $X$, then $W_1,\ldots, W_k$ have \emph{normal crossings at $p$} if there is a system of local coordinates $x_1,\ldots,x_n$ for $X$ around $p$ such that each subvariety $W_k$ is defined locally as the zero locus of a subset of those coordinates. The subvarieties $W_1,\ldots, W_k$ have \emph{normal crossings} if they have normal crossings at every point of $W_1\cap\cdots\cap W_k$.

A variety $X$ has \emph{normal crossings}, or is a \emph{normal crossings
variety}, if its irreducible components $X_1, \ldots, X_k$ have normal crossings.
\end{defn}

According to Hauser \cite[Proposition 3.17, p. 15]{Hauser_2014}, the condition in Definition \ref{defn:Normal_crossings_intersection} that $W_k$ is defined locally as the zero locus of a subset of the coordinates $x_1,\ldots,x_n$ for $X$ around $p$ is equivalent to the assertion that each germ $(W_j,p)$ in $(X,p)$ is defined by a radical monomial ideal for $j=1,\ldots,k$. (In the case of schemes, Hauser \cite[Remark 3.18, p. 15]{Hauser_2014} notes that the monomial ideal need not be radical.) If $W$ is a closed subscheme of a regular ambient scheme $X$, then $W$ is a \emph{normal crossings scheme} if it can be defined locally by a monomial ideal (see Hauser \cite[Appendix D, p. 394]{Hauser_2003}).

If the subvarieties $W_1,\ldots, W_k$ have normal crossings at $p$ as in Definition \ref{defn:Normal_crossings_intersection}, then (see Faber and Hauser \cite[p. 391]{Faber_Hauser_2010}) all possible intersections
\[
  \bigcap_{j\in J}W_j, \quad\text{with } J \subset \{1,\ldots,k\},
\]
are smooth (in the sense of schemes) and the subvarieties $W_1,\ldots, W_k$ intersect cleanly at $p$.

\section{Adapted analytic coordinates}
\label{sec:Adapted_analytic_coordinates}
The forthcoming Lemma \ref{lem:Adapted_analytic_coordinates} is is implicit in applications of the concept of normal crossing intersections in the category of analytic spaces, but we include a proof since we were unable to identify a reference. We shall assume for simplicity in the hypotheses of Lemma \ref{lem:Adapted_analytic_coordinates} that the intersection $S_1\cap \cdots \cap S_A$ is a point, but one should be able to extend the result to allow $S_1\cap \cdots \cap S_A$ to be an embedded, positive-dimensional, $\KK$-analytic submanifold such that $T_p(S_1\cap \cdots \cap S_A) = T_pS_1\cap \cdots \cap T_pS_A$, for all $p \in S_1\cap \cdots \cap S_A$, and thus again a clean intersection in the sense of Definition \ref{defn:Clean_intersection}.

\begin{lem}[Adapted analytic coordinates]
\label{lem:Adapted_analytic_coordinates}
Let $\KK=\RR$ or $\CC$ and $A\geq 2$ be an integer and $S_1, \ldots, S_A$ be embedded, positive-dimensional, $\KK$-analytic submanifolds of a $\KK$-analytic, finite-dimensional manifold $X$. Assume that
\begin{equation}
  \label{eq:Sum_tangent_spaces_equals_direct_sum_tangent_spaces}
  T_pS_1 + \cdots + T_pS_A = T_pS_1 \oplus \cdots \oplus T_pS_A.
\end{equation}
If $s_a := \dim S_a$ for $a=1,\ldots,A$ and $\dim X = m$, then around each point $p \in R$ there exists a $\KK$-analytic local coordinate chart $(\varphi,U)$ such that $\varphi(p)=0\in\KK^m$ and $\varphi(U\cap S_a) = \varphi(U)\cap\KK^{s_a}$, for $a=1,\ldots,A$, where the coordinate subspaces are defined by
\[
  \KK^{s_a} := \{x\in\KK^m: x_j = 0, \text{ for } j \notin J_a\} \hookrightarrow \KK^m,
\]
and $J_a := \{s_1+\cdots+s_{a-1}+1,\ldots,s_1+\cdots+s_{a-1}+s_a\}$, for $a=1,\ldots,A$, with the convention that $s_1+s_0$ or $s_1+s_1$ are replaced by $0$ or $s_1$ when $a=1$ or $2$, respectively.
\end{lem}

We note that the hypothesis \eqref{eq:Sum_tangent_spaces_equals_direct_sum_tangent_spaces} is implied by an apparently weaker assumption
\[
  \dim  T_pS_1 + \cdots + \dim T_pS_A = \dim\left(T_pS_1 + \cdots + T_pS_A\right),
\]
according to Hoffman and Kunze \cite[Section 6.6, p. 213, Exercise 2]{Hoffman_Kunze_linear_algebra}.

\begin{proof}[Proof of Lemma \ref{lem:Adapted_analytic_coordinates}]
It suffices to consider the case $A=2$, since the general case $A\geq 2$ differs only in notational complexity. The coordinate subspaces of $\KK^m$ corresponding to $S_1$ and $S_2$ ar given by
\begin{align*}
  \KK^{s_1} &:= \{x\in\KK^m: x_j = 0, \text{ for } j=s_1+1,\ldots,m\},
  \\
  \KK^{s_2} &:= \{x\in\KK^m: x_j = 0, \text{ for } j=1,\ldots,s_1 \text{ and } j=s_1+s_2+1,\ldots,m\}.
\end{align*}
Because $S_1, S_2 \subset X$ are embedded $\KK$-analytic submanifolds of $X$, there exist an open neighborhood $V\subset\KK^m$ of the origin and $\KK$-analytic embeddings $\phi_1:\KK^{s_1}\cap V \to X$ and $\phi_2:\KK^{s_2}\cap V \to X$ such that $\phi_1(\KK^{s_1}\cap V) = S_1\cap U_1$ and $\phi_2(\KK^{s_2}\cap V) = S_2\cap U_2$, where $U_1,U_2 \subseteq U$ are open neighborhoods of $p$ in $X$. The existence of the embeddings $\phi_1$ and $\phi_2$, given by $\phi_1^{-1} = \varphi_1\restriction (S_1\cap U_1)$ and $\phi_2^{-1} = \varphi_2\restriction (S_2\cap U_2)$ for submanifold coordinate charts $(U_1,\varphi_1)$ and $(U_2,\varphi_2)$, follows from Lee \cite[Theorem 4.12, p. 81, and Theorem 5.8, p. 101]{Lee_john_smooth_manifolds} or Guillemin and Pollack \cite[Section 1.3, Local Immersion Theorem, p. 15]{Guillemin_Pollack}. While the latter results are only stated for smooth manifolds, one obtains the analogous results for $\KK$-analytic manifolds by replacing in their proofs the role of the Inverse Mapping Theorem for smooth maps of smooth manifolds by the role of the Inverse Mapping Theorem for $\KK$-analytic maps of $\KK$-analytic manifolds.

Let $(\psi,U)$ be a $\KK$-analytic local coordinate chart around $p$ such that $\psi(p)=0\in\KK^m$ and, writing $\KK^{s_1+s_2} = \KK^{s_1}\times \KK^{s_2}$, define a $\KK$-analytic map $\phi:\KK^{s_1+s_2}\cap V \to X$ by
\[
  \phi
  :=
  \psi^{-1}\circ\left(\psi\circ\phi_1 \times \psi\circ\phi_2\right).
\]
Observe that
\begin{subequations}
\label{eq:Adapted_analytic_parameterization_pair_submanifolds}  
\begin{align}
  \label{eq:Adapted_analytic_parameterization_pair_submanifold_1}
  \phi(\KK^{s_1} \cap V) &= \phi_1(\KK^{s_1} \cap V) = S_1\cap U_1,
  \\
  \label{eq:Adapted_analytic_parameterization_pair_submanifold_2}
  \phi(\KK^{s_2} \cap V) &= \phi_2(\KK^{s_2} \cap V) = S_2\cap U_2.
\end{align}  
\end{subequations}
The differential $\d\phi(0):\KK^{s_1+s_2} \to T_pX$ is injective since, if $v_1\in \KK^{s_1}$ and $v_2\in \KK^s$, then
\[
  \d\phi(0)(v_1,v_2) = d\phi_1(0)v_1 + d\phi_2(0)v_2 \in T_pS_1 + T_pS_2 = T_pS_1 \oplus T_pS_2,
\]
where the final equality follows from the hypothesis \eqref{eq:Sum_tangent_spaces_equals_direct_sum_tangent_spaces}. Therefore, $\d\phi(0)(v_1,v_2) = 0 \implies d\phi_1(0)v_1 = 0 \in T_pS_1$ and $d\phi_2(0)v_2 = 0 \in T_pS_2$. Hence, $v_1 = 0 \in \KK^{s_1}$ and $v_2 = 0 \in \KK^{s_2}$ because the maps $\phi_1$ and $\phi_2$ are embeddings and therefore the differentials $d\phi_1(0)$ and $d\phi_2(0)$ are injective. Thus, after possibly shrinking $V$, the map $\phi$ is a $\KK$-analytic embedding of $\KK^{s_1+s_2} \cap V$ onto an embedded, open, $\KK$-analytic submanifold $S = \phi(\KK^{s_1+s_2}\cap V) \subset X$ by Lee \cite[Proposition 5.22, p. 110]{Lee_john_smooth_manifolds} (although stated for smooth manifolds, that result holds for $\KK$-analytic manifolds by reasoning similar to that used earlier in this proof). Consequently, after possibly shrinking $U$ and again applying Lee \cite[Theorem 5.8, p. 101]{Lee_john_smooth_manifolds}, there is a $\KK$-analytic coordinate chart $(\varphi,U)$ such that $\varphi(U\cap S) = \varphi(U)\cap \KK^{s_1+s_2}$ and $\phi^{-1} = \varphi \restriction (U\cap S)$. In particular, after possibly shrinking $U$, the identities \eqref{eq:Adapted_analytic_parameterization_pair_submanifolds} yield
\begin{equation}
  \label{eq:Adapted_analytic_coordinate_chart_pair_submanifolds}
  \varphi(U\cap S_1) = \varphi(U)\cap \KK^{s_1}
  \quad\text{and}\quad
  \varphi(U\cap S_2) = \varphi(U)\cap \KK^{s_2}.
\end{equation}
and so the chart $(\varphi,U)$ has the desired properties.
\end{proof}

\section{Local transformation group on an analytic space}
\label{sec:Local_transformation_group_analytic_space}
We give the following analogue of the standard definition of an additive, local one-real-parameter group of diffeomorphisms of a smooth manifold in Lee \cite[Chapter 9, pp. 211--212]{Lee_john_smooth_manifolds}.

\begin{defn}[Local transformation group on an analytic space]
\label{defn:Local_transformation_group_analytic_space}  
(See Akhiezer \cite[Section 1.2, p. 9]{Akhiezer_lie_group_actions_complex_analysis} for the case $\KK=\CC$.)  
Let $G$ be a topological group, $\KK=\RR$ or $\CC$, and $(X,\sO_X)$ a $\KK$-analytic space in the sense of Definition \ref{defn:Analytic_space}. Let $\Pi_X$ denote the collection of all pairs $\pi = (U_\pi,V_\pi)$, where $U_\pi$ and $V_\pi$ are open subsets of $X$ such that $U_\pi \Subset V_\pi$. Suppose that for each $\pi \in \Pi_X$ there is an open neighborhood $G_\pi$ of the identity $\id_G \in G$ and a map $\Phi_\pi:G_\pi \to \An(U_\pi,V_\pi)$, the Fr\'echet space of $\KK$-analytic maps from $U_\pi$ into $V_\pi$. The system $\{\Phi_\pi\}$ defines a \emph{local (continuous) $G$-action} on $X$ and $(G, \{\Phi_\pi\})$ is a \emph{local (topological) transformation group} of $X$ if the following conditions are satisfied:
\begin{enumerate}
\item For all $g, h \in G_\pi$ such that $k := gh \in G_\pi$, one has
  \[
    \Phi_\pi(g)\circ\Phi_\pi(h)\restriction U_{\pi,h} = \Phi_\pi(k)\restriction U_{\pi,k},
  \]
  where $U_{\pi,h} := \{x \in U_\pi: \Phi_\pi(h)x \in U_\pi\}$;
\item $\Phi_\pi(\id_G) = \id_X$;
\item For all $\pi,\rho \in \Pi_X$ and $g \in G_\pi \cap G_\rho$, one has
  \[
    \Phi_\pi(g)\restriction U_\pi\cap U_\rho = \Phi_\rho(g)\restriction U_\pi\cap U_\rho,
  \]
so that $gx := \Phi_\pi(g)x$ is independent of the choice of $\pi$ with $x \in U_\pi$ and $g \in G_\pi$;
\item For any two open sets $U \Subset U_\pi$ and $V \subset V_\pi$, the set
  \[
    W := W_{C,V} := \left\{g \in G_\pi | g\cdot C \subset V \right\}, 
  \]
  is open in $G_\pi$, where $C := \bar U$, and the map \eqref{eq:Action_topological_group_structure_sheaf} is continuous for all $f \in \sO_X(V)$.
\end{enumerate}
\end{defn}

We refer to Akhiezer \cite[Section 1.2, pp. 9--10]{Akhiezer_lie_group_actions_complex_analysis} for the concept of equivalence of two local actions $\{\Phi_\pi\}$, $\{\Psi_\pi\}$ of the same group $G$. For convenience, we shall write $\Phi:G\to\Aut_\loc(X)$ for a local action $\{\Phi_\pi\}$ of a topological group $G$ on an analytic space $X$.

\section{Sheaves of ideals for strict transforms in blowups}
\label{subsec:Sheaves_ideals_strict_transforms_resolution_blowups}
It is useful to explicitly identify the ideals that define the strict transform $\widetilde Y \subset \widetilde X$ of $Y \subset X$ in Definition \ref{defn:Strict_transform}. For this discussion, which extends that of Remark \ref{rmk:Strict_transform_colon_ideals}, we shall rely on Bravo and Villamayor \cite[Section 11, p. 43]{Bravo_Villamayor_2012claynotes}, Cutkosky \cite[Section 4.1, p. 38]{Cutkosky_resolution_singularities}, and Hauser \cite[Section 3, p. 16, following Exercise 8]{Hauser_2006}, \cite[Definition 6.2, p. 30]{Hauser_2014},  \cite{Hauser_encyc_math_blowup_algebra}, who identify these ideals in the categories of affine varieties and schemes.

\subsection{Strict transforms of ideals in the categories of affine algebraic and analytic varieties}
\label{subsec:Strict_transforms_ideals_category_affine_algebraic_varieties}
We first give a brief description of this identification in the category of affine algebraic and analytic varieties. As in Hauser \cite[Section 3, p. 16, following Exercise 8]{Hauser_2006}, let $\sI$ be an ideal in $\KK[x_1, \ldots , x_n]$, and let $f_1, \ldots, f_q$ be a (local) \emph{Macaulay basis} of $\sI$, namely, a generator system whose \emph{initial homogeneous forms} (of minimal degree) generate the ideal of all initial forms of elements of $\sI$ (see the forthcoming Definition \ref{defn:Initial_form}). Then the strict transform $\tilde\sI$ of $\sI = (f_1, \ldots, f_q)$ under the blowup of $\KK^n$ in a regular center is generated by the strict transforms $\tilde f_1, \ldots, \tilde f_q$ of the generators $f_1, \ldots, f_q$ of $\sI$:
\[
  \tilde\sI = (\tilde f_1, \ldots, \tilde f_q),
\]
with
\[
  \tilde f_i = m_E^{-\ord_P f_i}\cdot f_i^*, \quad\text{for } i=1,\ldots,q,
\]
where $f_1^*, \ldots, f_q^*$ denote the total transforms (pullbacks by the blowup map $\pi$) of $f_1, \ldots, f_q$ and $\ord_P f_i$ is the order of vanishing of $f_i$ at $P\in Y$ (see \cite[Section 3, p. 15, following Exercise 3]{Hauser_2006}), while $m_E=0$ is the reduced equation for the exceptional divisor $E \subset \widetilde X$ and, in local affine charts, $m_E$ is a monomial in one of the variables. See Hironaka \cite[Chapter III]{Hironaka_1964-I-II}, where the statement is proved in the local setting for the formal power series, $\KK[[x_1, \ldots , x_n]]$.

\subsection{Ideal quotients, affine varieties, and polynomial ideals}
\label{secsub:Ideal_quotients}
Before discussing the identification of strict transforms of ideals in the category of schemes in more detail, we need to recall the

\begin{defn}[Ideal quotients]
\label{defn:Ideal_quotients}
(See Atiyah and MacDonald \cite[Chapter 1, p. 8]{Atiyah_Macdonald_introduction_commutative_algebra}, Cox, Little, and O'Shea \cite[Chapter 4, Section 4, Definition 5, p. 200]{Cox_Little_OShea_ideals_varieties_algorithms}, Ene and Herzog \cite[Section 1.2.1, p. 6]{Ene_Herzog_grobner_bases_commutative_algebra}, and Greul and Pfister \cite[Definition 1.3.14, p. 27]{Greuel_Pfister_singular_introduction_commutative_algebra}.  
If $I, J$ are ideals in a ring $R$, then the \emph{ideal quotient} (or \emph{colon ideal}) of $I$ by $J$ is
\[
  (I:J) = \{x\in R: xJ \subset I\}.
\]
(One may also denote this ideal by $I:J$ or $[I:J]$, or by $I:_RJ$ to emphasize its dependence on the ambient ring $R$ and avoid ambiguity.) Moreover,
\[
  IJ^\infty = \{x\in R: \exists\, n \in \NN\cup\{0\} \text{ such that } xJ^n \subset I\}
\]
is called the \emph{saturation of $I$ with respect to $J$}.
\end{defn}

For geometric interpretations of ideal quotients in the category of affine varieties, we refer to Cox, Little, and O'Shea \cite[Chapter 4, Section 4, Definition 5, p. 200]{Cox_Little_OShea_ideals_varieties_algorithms}. Thus, if $\KK$ is any field and $\sI$ and $\sJ$ are ideals in $\KK[x_1,\ldots,x_n]$, then the Zariski closure of $\VV(\sI)\less\VV(\sJ)$ is given by \cite[Chapter 4, Section 4, Theorem 10, p. 203]{Cox_Little_OShea_ideals_varieties_algorithms}
\[
  \overline{\VV(\sI)\less\VV(\sJ)} \subset \VV(\sI:\sJ^\infty),
\]
where $\VV(\sI)$ is the affine variety associated to $\sI$ (see \cite[Chapter 1, Section 2, Definition 1, p. 5 and Section 4, Proposition 4, p. 31]{Cox_Little_OShea_ideals_varieties_algorithms}). If $\KK$ is algebraically closed, then \cite[Chapter 4, Section 4, Theorem 10, p. 203]{Cox_Little_OShea_ideals_varieties_algorithms}
\[
  \overline{\VV(\sI)\less\VV(\sJ)} = \VV(\sI:\sJ^\infty).
\]
Furthermore, if $\KK$ is algebraically closed and $\sI$ is radical, then \cite[Chapter 4, Section 4, Corollary 11, p. 204]{Cox_Little_OShea_ideals_varieties_algorithms}
\[
  \overline{\VV(\sI)\less\VV(\sJ)} = \VV(\sI:\sJ).
\]
According to \cite[Chapter 4, Section 4, Proposition 9, p. 202]{Cox_Little_OShea_ideals_varieties_algorithms}, one has
\[
  \sI \subset (\sI:\sJ) \subset (\sI:\sJ^\infty)
\]
and $(\sI:\sJ^\infty) = (\sI:\sJ^N)$ for all sufficiently large integers $N$. For further geometric properties of ideal quotients, we refer to Greul and Pfister \cite[pp. 81--83]{Greuel_Pfister_singular_introduction_commutative_algebra}.


\subsection{Bases for ideals in polynomial and formal power series rings}
\label{subsec:Bases_ideals_polynomial_and_formal_power_series_rings}
We begin with the

\begin{defn}[Initial form, initial ideal, and Macaulay basis]
\label{defn:Initial_form}  
(See Hauser \cite[Definition 6.5, p. 31]{Hauser_2014}.)
Let $\KK[x_1,\ldots,x_n]$ be the polynomial ring over $K$, considered with the natural grading given by the degree. Denote by $\inn(g)$ the homogeneous form of lowest degree of a non-zero polynomial $g$ of $\KK[x_1,\ldots,x_n]$, called the \emph{initial form} of $g$. Set $\inn(0) = 0$. For a non-zero ideal $I$, denote by $\inn(I)$ the ideal generated by all initial forms $\inn(g)$ of elements $g$ of $I$, called the \emph{initial ideal} of $I$. Elements $g_1, \ldots, g_k$ of an ideal $I$ of $\KK[x_1,\ldots,x_n]$ are a \emph{Macaulay basis} of $I$ if their initial
forms $\inn(g_1),\ldots,\inn(g_k)$ generate $\inn(I)$. 
\end{defn}

In Hironaka \cite[Section III.1, p. 208, Definition 3]{Hironaka_1964-I-II}, a Macaulay basis
was called a standard basis, but that is now used for a slightly more specific concept, as described in the Remark \ref{rmk:Standard_bases} below. Because $\KK[x_1,\ldots,x_n]$ is a Noetherian ring, any ideal possesses a Macaulay basis (see Hauser \cite[p. 31]{Hauser_2014}). 

\begin{prop}[Generators of the strict transform of an ideal]
\label{prop:Strict_transform_ideal}  
(See Hauser \cite[Proposition 6.6, p. 31]{Hauser_2014} or Hironaka \cite[Section III.2, p. 216, Lemma 6,  and Section III.6, p. 238, Theorem 5]{Hironaka_1964-I-II}.)
The strict transform of an ideal under blowup in a regular center is generated by the strict transforms of the elements of a Macaulay basis of the ideal.
\end{prop}  

The proof of Proposition \ref{prop:Strict_transform_ideal} relies on the Grauert--Hironaka--Galligo Division Theorem for the ring $\KK[[x_1,\ldots,x_n]]$ of formal power series; see Alonso, Castro--Jim\'enez, and Hauser \cite[Theorem 4.1]{Alonso_Castro-Jimenez_Hauser_2018} or Rond \cite[Theorem 10.1]{Rond_2018} for its statement, further references, discussions of the result and its applications, and generalizations.

\begin{rmk}[Concepts of standard bases]
\label{rmk:Standard_bases}  
(See Hauser \cite[Remark 6.7, p. 31]{Hauser_2014}.)
Continue the notation of Definitions \ref{defn:Strict_transform} and \ref{defn:Initial_form}.)
The strict transform of a Macaulay basis at a point $\tilde p$ of $\widetilde X$ need not be a Macaulay basis. This is however the case if the Macaulay basis is \emph{reduced} and the sequence of its orders has remained constant at $\tilde p$ (see  Hironaka \cite[Section III.8, p. 254, Lemma 20]{Hironaka_1964-I-II}. Instead of the grading of $\KK[x_1,\ldots,x_n]$ by degree, consider a grading such that all homogeneous elements are one-dimensional and generated by
monomials, that is, a grading induced by a monomial order on $\NN^n$. Then the initial form of
a polynomial and the initial ideal are both monomial. In this case, Macaulay bases
are called \emph{standard bases}. A standard basis $g_1, \ldots , g_k$ is \emph{reduced} if no monomial of the tails $g_i - \inn(g_i)$ belongs to $\inn(I)$. If the monomial order is degree compatible, that is, the induced grading a refinement of the natural grading of $\KK[x_1,\ldots,x_n]$ by degree, then the strict transforms of the elements of a standard basis of $I$ generate the strict
transform of the ideal. 
\end{rmk}

\begin{rmk}[Gr\"obner, Macaulay, $H$, and standard bases]
\label{rmk:Grobner_Macaulay_H_Standard_bases}
The analogues of \emph{Gr\"obner bases} for ideals in power series rings are called \emph{standard bases}. For introductions to Gr\"obner bases, we refer to Adams and Loustaunau \cite{Adams_Loustaunau_introduction_grobner_bases}, Becker and Weispfenning \cite{Becker_Weispfenning_grobner_bases}, Ene and Herzog \cite{Ene_Herzog_grobner_bases_commutative_algebra}, and Hibi \cite[Chapter 1]{Hibi_grobner_bases}. For a treatment of standard bases for ideals $I$ in the ring $\CC[\langle x_1,\ldots,x_n\rangle]$ of convergent power series, we refer to de Jong and Pfister \cite[Chapter 7]{DeJong_Pfister_local_analytic_geometry}. For treatments of both Gr\"obner bases for ideals in polynomial rings and standard bases for ideals in power series rings, we refer to Broer, Hoveijn, Lunter, and Vegter \cite[Chapter 6]{Broer_Hoveijn_Lunter_Vegter_bifurcations_hamiltonian_systems} and Cox, Little, and O'Shea \cite{Cox_Little_OShea_using_algebraic_geometry}. For treatments of Gr\"obner and Macaulay bases for ideals in polynomial rings and standard bases for ideals in power series rings, we refer to Kreuzer and Robbiano \cite{Kreuzer_Robbiano_computational_commutative_algebra_2}. For the relationship between Macaulay (also called \emph{$H$-bases} by Sauer) and Gr\"obner bases, we refer to Sauer \cite[Section 4, p. 2299]{Sauer_2001} and M\"oller \cite{Moller_1987}.
\end{rmk}

\subsection{Affine schemes and strict transform of an ideal via ideal quotients}
\label{subsec:Formula_strict_transform_ideal_via_ideal_quotients}
The abstract formula \eqref{eq:Strict_transform_ideal_as_union_colon_ideals} for the strict transform of an ideal is discussed further by Bravo and Villamayor \cite[Section 11, p. 43]{Bravo_Villamayor_2012claynotes}, Bravo, Encinas, and Villamayor \cite[Section 7, p. 372]{Bravo_Encinas_Villamayor_2005}, Cutkosky \cite[Section 4.1, p. 38]{Cutkosky_resolution_singularities}, and Schwede \cite{Schwede_2012mathoverflow_ideal_strict_transform}. We closely follow Bravo and Villamayor \cite[Section 11, p. 43]{Bravo_Villamayor_2012claynotes} in the following discussion. Let $A$ be a ring and $I \subset A$ be an ideal. If $I = (f_1, \ldots, f_r)$, then the blow-up $\Bl_I(A)$ of $A$ at $I$, namely
\[
  \Bl_I(A) \to \Spec(A),
\]
is defined by patching the affine morphisms
\[
  \Spec(A_i) \to \Spec(A),
\]
where $A_i := A[f_1/f_i,\ldots,f_r/f_i]$, a subring of $A_{f_i}$.

If $J \subset A$ is another ideal, then the \emph{total transform} $J^*$ of $J$ in $\Bl_I(A)$ is the coherent $\Bl_I(A)$-ideal defined by patching the extended ideals $JA_i$ in $A_i$. The ideal $J^*$ can also be defined by the homogeneous ideal $JR$ (see \cite[Remark 8.16 (1), p. 30]{Bravo_Villamayor_2012claynotes}). Indeed, the inclusion of graded modules,
\[
  0 \xrightarrow{} JR \xrightarrow{} R,
\]
defines an ideal on $\Bl_I(A)$, which, in each affine chart $\Spec(A_i)$, is also the ideal $JA_i$.

Let $\bar f_1, \ldots, f_r$ denote the image of $f_1, \ldots, f_r$ in $B := A/J$, so $\bar f_1, \ldots, f_r$ are generators of $\bar I := IB$. Consider
\[
  \bar R := B \oplus\bigoplus_{j=1}^\infty\bar I^j
\]
and observe that one obtains an exact sequence
\begin{equation}
  \label{eq:Bravo_Villamayor_2012_11-2-1}
  0 \xrightarrow{} H \xrightarrow{} R  \xrightarrow{} \bar R   \xrightarrow{} 0
\end{equation}
for some homogeneous ideal $H$. One can check that (see \cite[Section 11, p. 43]{Bravo_Villamayor_2012claynotes})
\[
  [H]_i = I^i\cap J, \quad\text{for all } i \geq 0.
\]  
The homogeneous ideal $H$ defines a $\Bl_I(A)$-ideal which is the \emph{strict transform} of $J$ in $\Bl_I(A)$. The restriction of the strict transform to the affine chart $\Spec(A_i)$ is given by an ideal, say $J_i$, in $A_i$. There is an inclusion of homogeneous ideals,
\[
  JR \subset H,
\]
so $JA_i$ is contained in $J_i$, as ideals in $A_i$, for each index $i\geq 0$. The restriction to the affine chart $\Spec(A_i)$ of the exact sequence of $\Bl_I(A)$-modules obtained from the short exact sequence \eqref{eq:Bravo_Villamayor_2012_11-2-1} is given by
\begin{equation}
  \label{eq:Bravo_Villamayor_2012_11-2-2}
  0 \xrightarrow{} J_i \xrightarrow{} A_i  \xrightarrow{} B_i \xrightarrow{} 0
\end{equation}
where
\[
B_i := B[\bar f_1/\bar f_i,\ldots,\bar f_r/\bar f_i],
\]
and which is a subring of $B_{\bar f_i}$.

%
%

\bibliography{/Users/pfeehan/Dropbox/LATEX/Bibinputs/master,/Users/pfeehan/Dropbox/LATEX/Bibinputs/mfpde}
\bibliographystyle{amsplain-nodash}

\end{document}